\documentclass[reqno]{report}

\usepackage{mathtools} 
\usepackage{amsmath, amssymb, amsthm}
\usepackage{stmaryrd} 
\usepackage{tikz}
\usepackage{tikz-cd}
\usepackage{xcolor}
\usepackage{etoolbox} 
\usepackage{hyperref}
\usepackage[alphabetic]{amsrefs} 
\usepackage[T1]{fontenc}
\usepackage{mathrsfs} 
\usepackage{newpxtext}
\usepackage[euler-digits, euler-hat-accent]{eulervm}
\usepackage[top=1in,bottom=1in,left=1in,right=1in]{geometry}  
\usepackage[nameinlink,capitalize]{cleveref}
\usepackage{autonum}
\usepackage{titlesec}
\titleformat{\subsection}[runin]
    {\normalfont\bfseries}
    {\thesubsection}
    {0.5em}
    {}
    []
\usepackage{makeidx}
\usepackage[refpage]{nomencl}
\usepackage{enumitem}

\makeindex
\makenomenclature

\renewcommand{\nompreamble}{%
This is a list of symbols that are either persistent during a large portion of the book,
important concepts in a given chapter, or simply fundamental in nature. The page
number denotes either where the notation is defined, or its first appearance.
Symbols that are only used temporarily or sporadically are not included and we
will try to provide cross-references in the text instead.

Due to the technicality involved in this book, the
notations may have to be changed locally, either due to conflicts with other
local notations, or to make certain part more readable. We will provide
warnings when potential confusion may arise in the text.

We also try to make font choices more indicative for concepts that have multiple
analogous versions:
\begin{itemize}[noitemsep]
    \item Bold symbols such as \(\bM\) are used for invariant theoretic concepts
        in the split or absolute setting;
    \item Fraktur or roman symbols such as \(\FRM, G\) are used for invariant
        theoretic concepts in the non-split setting;
    \item Calligraphic symbols such as \(\cM,\cA\) are used for arithmetic
        concepts in the local or global setting;
    \item Sans serif symbols are also used for some concepts
        used in the local setting, albeit less commonly; most notably for some arc-
        or loop-related constructions such as \(\Gr\) or \(\Hk\).
\end{itemize}
Many invariant theoretic constructions have both split and non-split
versions. We put both versions into one entry whenever possible.
}

\makeatletter
\renewcommand\@pnumwidth{2em}
\makeatother

\hypersetup{
    colorlinks,
    linkcolor={blue!50!black},
    citecolor={red!50!black},
    urlcolor={green!80!black}
}

\newtheorem{theorem}[subsection]{Theorem}
\newtheorem{proposition}[subsection]{Proposition}
\newtheorem{lemma}[subsection]{Lemma}
\newtheorem{corollary}[subsection]{Corollary}

\theoremstyle{definition}
\newtheorem{definition}[subsection]{Definition}

\newtheorem{example}[subsection]{Example}

\theoremstyle{remark}
\newtheorem{remark}[subsection]{Remark}

\numberwithin{equation}{section}

\Crefname{equation}{}{}
\crefname{equation}{}{}
\Crefname{section}{\S}{\S\S}
\crefname{section}{\S}{\S\S}
\Crefname{subappendix}{\S}{\S\S}
\crefname{subappendix}{\S}{\S\S}

\makeatletter
\def\@secnumfont{\bfseries}
\makeatother

\AtBeginDocument{
    \tikzset{AIIkeys/.cd,x/.initial=0,y/.initial=0}
    \tikzdeclarecoordinatesystem{AII}{%
        \tikzset{AIIkeys/.cd,#1}
        \pgfpointxy{2*\pgfkeysvalueof{/tikz/AIIkeys/x}-2*0.5*\pgfkeysvalueof{/tikz/AIIkeys/y}}{2*0.866*\pgfkeysvalueof{/tikz/AIIkeys/y}}
    }

\DeclareFontEncoding{LS1}{}{}
\DeclareFontSubstitution{LS1}{stix2}{m}{n}
\DeclareSymbolFont{CormaExOperators}{LS1}{stix2}{m}{n}
\DeclareSymbolFont{CormaMathbb}{LS1}{stix2bb}{m}{n}
\DeclareSymbolFontAlphabet{\mathbb}{CormaMathbb}
\DeclareMathSymbol{\CormaQED}{\mathord}{CormaExOperators}{"D1}

\newcommand{\mathup}{\mathrm}
\newcommand{\symup}{\mathrm}
\newcommand{\symbf}{\mathbf}
\newcommand{\symbb}{\mathbb}
\newcommand{\symsf}{\mathsf}
\newcommand{\symcal}{\mathcal}
\newcommand{\symscr}{\mathscr}
\newcommand{\symfrak}{\mathfrak}

\newcommand{\hy}{{\operatorname{-}}}

\renewcommand{\emptyset}{\varnothing}

\newcommand{\bA}{{\symbf{A}}}
\newcommand{\bB}{{\symbf{B}}}
\newcommand{\bC}{{\symbf{C}}}
\newcommand{\bD}{{\symbf{D}}}
\newcommand{\bE}{{\symbf{E}}}

\newcommand{\bG}{{\symbf{G}}}
\newcommand{\bH}{{\symbf{H}}}
\newcommand{\bI}{{\symbf{I}}}
\newcommand{\bJ}{{\symbf{J}}}

\newcommand{\bL}{{\symbf{L}}}
\newcommand{\bM}{{\symbf{M}}}
\newcommand{\bN}{{\symbf{N}}}

\newcommand{\bP}{{\symbf{P}}}

\newcommand{\bR}{{\symbf{R}}}
\newcommand{\bS}{{\symbf{S}}}
\newcommand{\bT}{{\symbf{T}}}
\newcommand{\bU}{{\symbf{U}}}

\newcommand{\bW}{{\symbf{W}}}
\newcommand{\bX}{{\symbf{X}}}

\newcommand{\bZ}{{\symbf{Z}}}

\newcommand{\bFa}{{\symbf{a}}}
\newcommand{\bFb}{{\symbf{b}}}
\newcommand{\bFc}{{\symbf{c}}}

\newcommand{\bFe}{{\symbf{e}}}
\newcommand{\bFf}{{\symbf{f}}}

\newcommand{\bFm}{{\symbf{m}}}
\newcommand{\bFn}{{\symbf{n}}}
\newcommand{\bFo}{{\symbf{o}}}

\newcommand{\bFs}{{\symbf{s}}}

\newcommand{\bFx}{{\symbf{x}}}

\newcommand{\bbA}{{\symbb{A}}}
\newcommand{\bbB}{{\symbb{B}}}
\newcommand{\bbC}{{\symbb{C}}}

\newcommand{\bbF}{{\symbb{F}}}

\newcommand{\bbM}{{\symbb{M}}}
\newcommand{\bbN}{{\symbb{N}}}

\newcommand{\bbP}{{\symbb{P}}}
\newcommand{\bbQ}{{\symbb{Q}}}
\newcommand{\bbR}{{\symbb{R}}}

\newcommand{\bbZ}{{\symbb{Z}}}

\newcommand{\cA}{{\symcal{A}}}
\newcommand{\cB}{{\symcal{B}}}
\newcommand{\cC}{{\symcal{C}}}
\newcommand{\cD}{{\symcal{D}}}
\newcommand{\cE}{{\symcal{E}}}
\newcommand{\cF}{{\symcal{F}}}
\newcommand{\cG}{{\symcal{G}}}
\newcommand{\cH}{{\symcal{H}}}

\newcommand{\cJ}{{\symcal{J}}}
\newcommand{\cK}{{\symcal{K}}}
\newcommand{\cL}{{\symcal{L}}}
\newcommand{\cM}{{\symcal{M}}}

\newcommand{\cO}{{\symcal{O}}}
\newcommand{\cP}{{\symcal{P}}}

\newcommand{\cR}{{\symcal{R}}}
\newcommand{\cS}{{\symcal{S}}}
\newcommand{\cT}{{\symcal{T}}}
\newcommand{\cU}{{\symcal{U}}}
\newcommand{\cV}{{\symcal{V}}}

\newcommand{\sA}{{\symscr{A}}}
\newcommand{\sB}{{\symscr{B}}}
\newcommand{\sC}{{\symscr{C}}}

\newcommand{\sG}{{\symscr{G}}}
\newcommand{\sH}{{\symscr{H}}}

\newcommand{\sQ}{{\symscr{Q}}}

\newcommand{\sX}{{\symscr{X}}}
\newcommand{\sY}{{\symscr{Y}}}

\newcommand{\FRc}{{\symfrak{c}}}

\newcommand{\FRg}{{\symfrak{g}}}

\newcommand{\FRm}{{\symfrak{m}}}

\newcommand{\FRp}{{\symfrak{p}}}

\newcommand{\FRx}{{\symfrak{x}}}
\newcommand{\FRy}{{\symfrak{y}}}

\newcommand{\FRA}{{\symfrak{A}}}

\newcommand{\FRC}{{\symfrak{C}}}
\newcommand{\FRD}{{\symfrak{D}}}
\newcommand{\FRE}{{\symfrak{E}}}

\newcommand{\FRJ}{{\symfrak{J}}}

\newcommand{\FRM}{{\symfrak{M}}}

\newcommand{\FRR}{{\symfrak{R}}}
\newcommand{\FRS}{{\symfrak{S}}}
\newcommand{\FRT}{{\symfrak{T}}}

\newcommand{\FRV}{{\symfrak{V}}}

\newcommand{\SFM}{{\symsf{M}}}

\newcommand{\SFQ}{{\symsf{Q}}}

\newcommand{\widebar}[1]{\mkern 2.0mu\overline{\mkern-2.0mu#1\mkern-2.0mu}\mkern 2.0mu}
\renewcommand{\bar}{\widebar}

\newcommand{\x}{\times}

\newcommand{\longot}{\longleftarrow}
\newcommand{\longto}{\longrightarrow}

\newcommand{\ar}{\rightarrow}

\newcommand{\hookto}{\hookrightarrow}

\newcommand{\olto}{\overleftarrow}
\newcommand{\orto}{\overrightarrow}
\newcommand{\olrto}{\overleftrightarrow}
\newcommand{\defeq}{\coloneqq}
\let\eqdef\undefined
\newcommand{\eqdef}{\eqqcolon}
\newcommand{\ul}{\underline}

\newcommand{\notion}{\emph}

\providecommand{\given}{}
\newcommand\SetSymbol[1][]{%
    \nonscript\:#1\vert
    \allowbreak
    \nonscript\:
    \mathopen{}}
\DeclarePairedDelimiterX{\Set}[1]{\{}{\}}{%
    \renewcommand\given{\SetSymbol[\delimsize]}
    #1
}
\DeclarePairedDelimiterX{\restr}[1]{}{\vert}{#1}
\DeclarePairedDelimiterX{\Stack}[1]{[}{]}{%
    #1
}

\newcommand{\OneHalf}{{\frac{1}{2}}}

\newcommand{\dd}{\symup{d}} 
\newcommand{\DD}{\symup{D}} 

\DeclarePairedDelimiterX{\abs}[1]{\lvert}{\rvert}{%
    \ifblank{#1}{\:\cdot\:}{#1}
}
\DeclarePairedDelimiterX{\norm}[1]{\lVert}{\rVert}{%
    \ifblank{#1}{\:\cdot\:}{#1}
}
\DeclarePairedDelimiterX{\lrangle}[1]{\langle}{\rangle}{%
    \ifblank{#1}{\:\cdot\:}{#1}
}
\DeclarePairedDelimiterX{\powser}[1]{[\![}{]\!]}{%
    \ifblank{#1}{\:\cdot\:}{#1}
}
\DeclarePairedDelimiterX{\lauser}[1]{(\!(}{)\!)}{%
    \ifblank{#1}{\:\cdot\:}{#1}
}

\DeclareMathOperator{\vol}{vol} 

\DeclareMathOperator{\ord}{ord} 
\DeclareMathOperator{\Cnt}{\#} 
\DeclareMathOperator{\Hom}{Hom} 
\DeclareMathOperator{\Ext}{Ext} 
\DeclareMathOperator{\Aut}{Aut} 
\DeclareMathOperator{\Out}{Out} 
\DeclareMathOperator{\rk}{rk} 
\DeclareMathOperator{\Tr}{Tr} 
\DeclareMathOperator{\Nm}{Nm} 
\DeclarePairedDelimiter{\Ggen}{\langle}{\rangle}%
\newcommand{\id}{\symup{id}}
\newcommand{\Id}{\symup{id}}
\newcommand{\ev}{\symup{ev}} 
\newcommand{\dbq}[3]{{#1\backslash #2/#3}} 
\DeclareMathOperator{\coker}{coker} 
\DeclareMathOperator{\img}{Im} 

\DeclareMathOperator{\Sym}{Sym}
\newcommand{\Disc}{\symup{Disc}} 

\DeclareMathOperator{\RH}{H} 
\newcommand{\RHc}[1][]{%
    \ifblank{#1}{%
        {\symup{H}}_{\symup{c}}
    }{
        {\symup{H}}_{\symup{c},#1}
    }
} 
\DeclareMathOperator{\PH}{\prescript{\FRp}{}{H}} 
\newcommand{\ESP}{\symup{E}} 
\newcommand{\RDF}{\symbf{R}} 
\newcommand{\LDF}{\symbf{L}} 
\newcommand{\Ptr}{\prescript{\FRp}{}{\tau}} 

\DeclareMathOperator{\Cox}{Cox} 

\DeclareMathOperator{\Spec}{Spec}

\DeclareMathOperator{\RSpec}{\ul{Spec}}

\DeclareMathOperator{\Res}{Res} 
\DeclareMathOperator{\codim}{codim} 
\DeclareMathOperator{\TanB}{T} 
\newcommand{\CoTB}{\Omega} 
\DeclareMathOperator{\PicS}{\cP ic} 
\DeclareMathOperator{\supp}{supp} 
\newcommand{\Sing}{{\symup{sing}}} 
\newcommand{\Sg}{{\symup{sing}}} 
\newcommand{\Sm}{{\symup{sm}}} 
\newcommand{\Gm}[1][]{%
    \ifblank{#1}{%
        {\symbb{G}}_{\symup{m}}
    }{
        {\symbb{G}}_{\symup{m},#1}
    }
}
\newcommand{\Ga}[1][]{%
    \ifblank{#1}{%
        {\symbb{G}}_{\symup{a}}
    }{
        {\symbb{G}}_{\symup{a},#1}
    }
}
\DeclareMathOperator{\BG}{\bbB\!} 
\DeclareMathOperator{\Bun}{Bun} 
\DeclareMathOperator{\Tot}{Tot} 
\newcommand{\git}{\mathbin{/\mkern-6mu/}} 
\newcommand{\pr}{{\symup{pr}}} 
\newcommand{\Gr}{\symsf{Gr}} 
\newcommand{\Hk}{\symsf{Hk}} 
\newcommand{\Cartan}{\symsf{C}} 
\newcommand{\ShSec}{\Gamma} 
\newcommand{\tfree}{{\symup{tf}}} 
\DeclareMathOperator{\Hilb}{Hilb} 
\DeclareMathOperator{\IHom}{\underline{Hom}} 
\DeclareMathOperator{\RIHom}{\symbf{R}\underline{Hom}} 
\DeclareMathOperator{\IExt}{\underline{Ext}} 
\DeclareMathOperator{\IAut}{\underline{Aut}} 
\DeclareMathOperator{\Irr}{Irr} 
\newcommand{\Red}{{\symup{red}}} 
\newcommand{\Loop}{\symbb{L}} 
\newcommand{\Arc}[1][]{
    \ifblank{#1}{%
        \Loop^+
    }{
        \Loop^{+ #1}
    }
}
\newcommand{\IC}{{\symup{IC}}} 
\DeclareMathOperator{\TateM}{T} 
\DeclareMathOperator{\DCat}{D} 
\DeclareMathOperator{\Dbc}{D^b_c} 
\DeclareMathOperator{\Ob}{Ob} 

\DeclareMathOperator{\val}{val} 
\DeclareMathOperator{\Char}{char} 
\DeclareMathOperator{\Gal}{Gal} 
\newcommand{\Frob}{\symup{Frob}} 

\newcommand{\Qlb}{\bar{\bbQ}_\ell}

\newcommand{\sep}{{\symup{s}}} 

\newcommand{\TypeA}{\symup{A}} 
\newcommand{\TypeB}{\symup{B}}
\newcommand{\TypeC}{\symup{C}}
\newcommand{\TypeD}{\symup{D}}
\newcommand{\TypeE}{\symup{E}}
\newcommand{\TypeF}{\symup{F}}
\newcommand{\TypeG}{\symup{G}}

\DeclareMathOperator{\Lie}{Lie} 
\DeclareMathOperator{\Ad}{Ad} 
\DeclareMathOperator{\ad}{ad} 
\DeclareMathOperator{\Cent}{C} 
\DeclareMathOperator{\Norm}{N} 
\newcommand{\dual}{\check} 
\newcommand{\LD}[1]{\prescript{L}{}{#1}} 
\DeclareMathOperator{\GL}{GL} 
\DeclareMathOperator{\PGL}{PGL} 
\DeclareMathOperator{\SL}{SL} 
\DeclareMathOperator{\SYP}{Sp} 
\DeclareMathOperator{\SO}{SO} 
\DeclareMathOperator{\Mat}{Mat} 
\newcommand{\La}[1]{\symfrak{#1}} 
\DeclarePairedDelimiterX\Pair[2]{\langle}{\rangle}{#1,#2}

\newcommand{\Der}{{\symup{der}}} 
\newcommand{\AB}{{\symup{ab}}} 
\newcommand{\CharG}{{\symbb{X}}} 
\newcommand{\CoCharG}{\check{\symbb{X}}} 
\newcommand{\Roots}{\Phi} 
\newcommand{\CoRoots}{\check{\Roots}} 
\newcommand{\PosRts}{{\Roots_+}} 
\newcommand{\PosCoRts}{{\CoRoots_+}} 
\newcommand{\SimRts}{\Delta} 
\newcommand{\SimCoRts}{\check{\SimRts}} 
\newcommand{\Rt}{\alpha} 
\newcommand{\CoRt}{\check{\alpha}} 
\newcommand{\Wt}{\varpi} 
\newcommand{\CoWt}{\check{\varpi}} 
\newcommand{\SC}{{\symup{sc}}} 
\newcommand{\AD}{{\symup{ad}}} 
\newcommand{\reg}{{\symup{reg}}} 
\newcommand{\ssim}{{\symup{ss}}} 
\newcommand{\rss}{{\symup{rs}}} 
\newcommand{\srs}{{\symup{srs}}} 
\DeclareMathOperator{\Env}{Env} 
\DeclareMathOperator{\inv}{inv} 
\DeclareMathOperator{\Inv}{Inv} 

\newcommand{\Sat}{\symup{Sat}} 
\newcommand{\Rep}{\symup{Rep}} 
\newcommand{\Alg}{{\symup{alg}}} 

\newcommand{\bSPL}{\boldsymbol{\wp}}
\newcommand{\dbSPL}{\dual{\boldsymbol{\wp}}}
\newcommand{\SPL}{\wp}

\newcommand{\bOmega}{\boldsymbol{\Omega}}

\DeclareMathOperator{\FM}{\cF\cM}
\newcommand{\BD}[1][]{ 
  \ifblank{#1}{
      \cB^!
  }{ 
      \cB^{! #1}
  }
}
\newcommand{\KV}{\symsf{M}} 
\newcommand{\KVU}{\symsf{Y}} 
\newcommand{\SP}{\cM} 
\renewcommand{\Frob}{\sigma} 
\newcommand{\OGT}{\vartheta} 
\newcommand{\hST}{{\symup{st}}} 
\newcommand{\bcQ}{\sQ}
\newcommand{\PSP}{\Sigma} 

\DeclareMathOperator{\occ}{occ}
\DeclareMathOperator{\cl}{cl}
\DeclareMathOperator{\OI}{\symbf{O}} 
\DeclareMathOperator{\SOI}{\symbf{SO}} 
\newcommand{\vC}{\check{C}} 
\newcommand{\iH}{\ul{\sH}} 
\newcommand{\GASch}{\symsf{Q}} 
\newcommand{\OGASch}[1][]{ 
  \ifblank{#1}{
      \GASch^!
  }{ 
      \GASch^{! #1}
  }
}
\newcommand{\disj}{{\symup{dj}}} 
\newcommand{\ACT}{{\symup{act}}} 
\newcommand{\WP}{\Xi} 
\newcommand{\TD}[1]{\prescript{\Theta}{}{#1}} 
\newcommand{\fS}{f} 
\newcommand{\CRWT}{\symup{wt}} 
\newcommand{\NC}{\symcal{NC}} 
\newcommand{\ANI}{{\symup{ani}}} 
\newcommand{\ELL}{{\symup{ell}}} 
\renewcommand{\notion}[1]{\emph{#1}}
\newcommand{\inotion}[1]{\notion{#1}\index{#1}}
}

\title{Multiplicative Hitchin Fibrations\\and the\\Fundamental Lemma}
\author{X. Griffin Wang}
\date{June 9, 2025}

\begin{document}

\maketitle
\setcounter{tocdepth}{1}

\begin{abstract}
    Let \(k\) be a finite field and let \(G\) be a reductive group over
    \(k\powser{\pi}\). Suppose \(\Char(k)\) is larger than twice the Coxeter
    number of \(G\), we prove the standard endoscopic fundamental lemma for the
    spherical Hecke algebra of \(G\) using multiplicative Hitchin fibrations.
\end{abstract}

\tableofcontents

\chapter{Introduction}%
\label{chap:introduction}

In this book, we prove the standard
endoscopic fundamental lemma for the spherical Hecke algebra of an unramified
reductive group and its unramified endoscopic groups.
A more precise version of the following result may
be found in \Cref{thm:FL_main}.

\begin{theorem}
    \label[theorem]{thm:FL_intro}
    Let \(k=\bbF_q\) be a finite field with \(q\) element and characteristic
    \(p\), \(\cO=k\powser{\pi}\) the ring of power series of one variable over \(k\)
    and \(F=k\lauser{\pi}\) its field of fractions. Let \(G\) be a reductive
    group scheme over \(\cO\) whose Coxeter number is less than
    \(p/2\). Let \((\kappa,\OGT_\kappa,\xi)\) be an endoscopic datum of
    \(G\) over \(\cO\) and \(H\) the corresponding endoscopic group. Then we
    have equality in orbital integrals
    \begin{align}
        \Delta_0(\gamma_H,\gamma)\OI_{a}^\kappa\bigl(f^{V},\dd t\bigr)
        =\SOI_{a_H}\Bigl(f_{H,\xi}^{V},\dd t\Bigr),
    \end{align}
    where \(a\) and \(a_H\) are matching strongly regular semisimple conjugacy
    classes represented by \(\gamma\in G(F)\) and \(\gamma_H\in H(F)\)
    respectively, \(f^V\) is the
    Satake function associated with an \(\LD{G}\)-representation
    \(V\) and \(f_{H,\xi}^{V}\) is the corresponding Satake function of \(H\)
    by transferring through \(L\)-embedding
    \(\xi\), and \(\Delta_0(\gamma_H,\gamma)\) is the (appropriately defined
    version of the) transfer factor.
\end{theorem}

The Lie algebra analogue of \Cref{thm:FL_intro} was proved by B.C.~Ng\^o in his
monumental work \cite{Ng10}. Based on Ng\^o's result, people were able to deduce
several analogous results. First, when combined with works of Waldspurger
\cite{Wal06}, one may deduce the Lie algebra fundamental lemma over local number
fields. Waldspurger's another work \cite{Wal97} then allows us to deduce the
analogue of \Cref{thm:FL_intro} for the unit element of the spherical Hecke
algebra over local number fields. Later, Hales \cite{Ha95} used some simple
trace formula to further obtain the result for the entire spherical Hecke
algebra, again over local number fields. Finally, Casselman, Cely and Hales in
\cite{CaCeHa19} used motivic integration to deduce a slightly weaker version of
\Cref{thm:FL_intro} from \cite{Ha95}, where the characteristic of \(k\) is
assumed to be ``large'' rather than having an explicit lower bound.

For many applications in number theory, the historical results recounted above
may be viewed as ``complete''. However, there are still deficiencies from a
development point of view. For example, the chain of logic across several similar
yet different situations, namely number fields versus function
fields, and Lie algebra versus group, seems unnecessarily convoluted. It would
already be very interesting if a direct proof can be found for the group case.
For another, the proof of \cite{Ng10} is entirely geometric and provides
great insight into understanding the trace formula for Lie algebras through
geometric lens, almost all of which is not reflected in later works mentioned
above. Moreover, the trace formula for Lie algebra is only a watered-down
version of the \textit{bona fide} trace formula, and many interesting phenomena
are lost. Therefore, a direct, geometric proof of \Cref{thm:FL_intro} in the
style of \cite{Ng10} is very much desirable.

For these reasons, the organizational structure of \cite{Ng10} is used as a base template for
this book, and comparison with the Lie algebra case will be a recurring theme.
The central object of this book is the multiplicative analogue of Hitchin
fibration in the case of the adjoint action of a reductive group. We will call
it the \notion{multiplicative Hitchin fibration}, or \notion{mH-fibration} for
short. Such generalization in arithmetic setting was first envisioned in
\cite{FN11}, and some important preliminary explorations were done notably in
\cites{Bo15,Bo17,BoCh18,Ch22}.

In recent years, the so-called relative Langlands program has been of great
interest, where a much more general class of trace formulae involving spherical
varieties are considered. Therefore, the current book also serves as a proof of
concept for future development.
To this end, in this introduction we will try to
give an overview of this book with an emphasis on the big picture, new features
and potential development points.

The geometry of mH-fibration is incredibly rich and beautiful and has
profound connection with representation theory, therefore it is important to
give enough spotlight to all the key features in this introduction so that the
reader can have a grasp of what is to come. However, to keep the introduction
concise and clear enough, we will also try to avoid using too many notations
due to the technical complexity of this subject.

\section{On Hitchin-type Fibrations} 
\label{sec:on_hitchin_type_fibrations}

One of the motivations behind introducing Hitchin fibration to arithmetic
problems like fundamental lemma is that global objects tend to behave much
better than local ones, even if they are less computationally accessible at
times. Therefore, it would be beneficial to understand the general principles on
these Hitchin-type global constructions. We will use mH-fibrations
as the primary example for this section.

\subsection{}
There have been many attempts to prove the fundamental lemma for reductive groups or Lie
algebras in the past focusing on the local picture.
Over Archimedean local fields Shelstad was able to prove the general
smooth transfer directly at group level (see \cite{She79}), thus bypassing the
need of a fundamental lemma. Over non-Archimedean
fields, however, the problem seems much harder, and most proofs before
\cite{Ng10} were on one specific group at a time.
The reader can see the introduction of \cite{Ng10} for more historical
details.

One of the most successful general results on the local front is
perhaps the conditional proof for unramified conjugacy classes by Goresky,
Kottwitz and MacPherson in \cite{GKM04}. Although quite an impressive framework
in itself, it also shows the limitation of local geometric method. For one, as
Ng\^o already pointed out in \cite{Ng10}, it depends on a purity conjecture of
the cohomology of related affine Springer fibers, which is only partially known
due to affine Springer fibers usually being highly singular. Another serious
obstacle pointed out by \cite{Ng10} is that it crucially depends on the
conjugacy class being unramified so that there is a large torus acting on the
affine Springer fibers.

\subsection{}
Roughly speaking, the local geometric method is based on the fact that orbital
integrals may be interpreted as certain kind of point-counting on related affine
Springer fibers, which in turn reduces the fundamental lemma to some cohomological statement using
a variant of Grothendieck--Lefschetz trace formula. The affine Springer fiber in
question, can be roughly understood as certain subset of maps from a formal disc
\(X_v=\Spec{\cO}\) to the quotient stack \(\Stack{\La{g}/G}\), together with
some other naturally attached data. To move from the local picture to the global
one, one only needs to replace the formal disc by a smooth projective curve
\(X\). Because \(\La{g}\) is affine and \(X\) is projective, one needs to add
some auxiliary twist, otherwise the global object we end up getting will be
trivial. In Lie algebra case, such twist is provided by the natural
\(\Gm\)-action on \(\La{g}\) viewed as a vector space.

More formally, the classical Hitchin fibration can be formulated as follows: by
Chevallay restriction theorem, the GIT quotient \(\FRc=\La{g}\git G\simeq
\La{t}\git W\) is an affine space. The \(\Gm\)-scaling on \(\La{g}\) commutes
with the adjoint action, so the map \(\La{g}\to \FRc\) is \(\Gm\)-equivariant.
Consider the two-stage map
\begin{align}
    \Stack*{\La{g}/G\x\Gm}\longto\Stack*{\FRc/\Gm}\longto \BG{\Gm},
\end{align}
where \(\BG{\Gm}\) is the classifying stack of \(\Gm\). Apply the
mapping stack functor from \(X\) to the first map, we obtain the
\notion{universal Hitchin fibration} associated with \(\La{g}\) and \(X\),
lying over the Picard stack of \(X\). Taking
the fiber over the canonical sheaf of \(X\), or more generally any line bundle,
we obtain the Hitchin fibration in the usual sense.

\subsection{}
Without going further into the global geometry, we must first ask what is the
analogue for the above setup in the group case. Clearly
we are now interested in the quotient stack \(\Stack{G/G}\) (where \(G\) acts on
itself by adjoint action) in place of \(\Stack{\La{g}/G}\), but we also need a
natural twisting similar to the \(\Gm\)-action on \(\La{g}\).
If \(G=\GL_n\), for example, we have the scaling action of
\(\Gm\) by viewing \(G\) as a subspace of \(\Mat_n\), and for a group \(G\) with
a non-trivial central torus \(Z_0\), we may use the action of \(Z_0\).

There are several problems using \(Z_0\), however.
For one, this action heavily depends on the
situation, and does not exist if \(G\) is
semisimple. Even if \(Z_0\) is not trivial, it still does not provide
enough twisting if it is anisotropic. There is also another difficulty:
even for \(G=\GL_1=\Gm\), if we twist \(G\) by a \(\Gm\)-torsor \(\cL\), then
the induced bundle \(G_\cL\) is just
\(\cL\) itself, which has no global sections unless \(\cL\)
is trivial. Comparing to the Lie algebra case when \(\La{g}=\bbA^1\):
\(\La{g}_\cL\) is the line bundle associated with the \(\Gm\)-torsor which does
have a lot of sections provided that the degree is large enough.

\subsection{}
There is a relatively straightforward way to solve this lack-of-twist
problem: since there is no
non-constant map from \(X\) to \(\Stack{G/G}\), we may instead consider rational
maps. This is the same idea behind global Hecke stacks, and in order to obtain a
reasonable geometric object, we also need a systematic way to control the poles
of the rational maps, similar to truncated Hecke stacks. This is in fact the
formulation used by
some prototype of mH-fibrations in the field of mathematical physics
(see, for example, \cite{HM02}), as well as its first appearance in arithmetic
setting in \cite{FN11}.

One downside (among many others) of this ``Hecke-style'' formulation is that
some technical aspects become very messy and unwieldy if \(G^\SC\neq G^\AD\) or
\(G\) is not split. On the other hand, the mapping stack approach used in Lie
algebra case is very clean when dealing with non-split groups, so pursuing that
route seems more promising.

Since we cannot guarantee that \(G\) contains a central torus, we can
certainly add one to it by augmenting the group. The question is of course what
torus to add. Moreover, the \(G=\Gm\) versus \(\La{g}=\bbA^1\) example above
shows we need to add some ``boundaries'' to \(G\) so that the mapping stack
contains more than just constant maps. This is where reductive monoids enter the
picture.

\subsection{}
A primal example of reductive monoid is that of \(n\x n\) matrices
\(\Mat_n\). The group \(\GL_n\) embeds in
\(\Mat_n\) as an open subset and is its group of units. Therefore,
\(\Mat_n-\GL_n\) is the ``boundary'' we could add to \(\GL_n\).
Coincidentally, \(\Mat_n\cong\La{gl}_n\), so our first example of mH-fibration
is the same Hitchin fibration associated with \(\La{gl}_n\).

The \(\Mat_n\) example is somewhat misleading, because it is essentially the
only case where the monoid is smooth.
In general, the Lie algebra is replaced by a \notion{very flat} reductive monoid
\(\FRM\) (see \Cref{sec:review_of_very_flat_reductive_monoids}), and the role of
\(\Gm\) is replaced by the center \(Z_\FRM\) of the
unit group \(\FRM^\x\). There is a maximal toric variety \(\FRT_\FRM\) playing
the role of Cartan subalgebra, and it is known that \(\FRC_\FRM\defeq\FRM\git
G\simeq\FRT_\FRM\git W\). However, there is an extra step for mH-fibrations.

\subsection{}
The monoid \(\FRM\) comes with an \notion{abelianization} map \(\FRM\to\FRA_\FRM\)
introduced by Vinberg in \cite{Vi95}  by
taking the GIT quotient by \(G^\SC\x G^\SC\)-multiplication on the left
and right. We will be primarily interested in the case where the commutative
monoid \(\FRA_\FRM\) is an affine space (the general case can be dealt with
after some technical work but does not add anything new). It can be shown that
\begin{align}
    \FRC_\FRM\simeq\FRA_\FRM\x\FRC,
\end{align}
where \(\FRC\defeq G^\SC\git G^\SC\) is an affine space.
The upshot is that we have a three-stage map
\begin{align}
    \Stack*{\FRM/G\x
    Z_\FRM}\longto\Stack*{\FRC_\FRM/Z_\FRM}\longto\Stack*{\FRA_\FRM/Z_\FRM}\longto\BG{Z_\FRM},
\end{align}
and we may apply the mapping stack functor to these maps. The resulting
stacks are denoted by
\begin{align}
    \cM_X^+\longto\cA_X^+\longto\cB_X^+\longto\Bun_{Z_\FRM}.
\end{align}

\subsection{}
There is a maximal open Deligne--Mumford substack \(\cB_X\subset\cB_X^+\) which
is essentially a moduli
stack of ``divisors'' on \(X\) with values in the dominant cocharacters of \(G\)
(see \Cref{sec:boundary_divisors,sec:Constructions}). Such divisors are
usually called \notion{colored divisors} in some literature, and are called
\notion{boundary divisors} in this book. Taking inverse images, we have stacks
\begin{align}
    \cM_X\longto\cA_X\longto\cB_X\longto \Bun_{Z_\FRM},
\end{align}
the first step of which is the \notion{universal mH-fibration} associated with
\(G\), \(\FRM\) and \(X\). Unlike Lie algebra case where the \(\Gm\)-torsor is
fixed, we will \emph{not} fix a
\(Z_\FRM\)-torsor.

In earlier papers like \cites{Bo15,Bo17,Ch22}, \(\FRM\) is restricted to the
so-called \notion{Vinberg monoid} or \notion{universal monoid} of
\(G^\SC\) (denoted by \(\Env(G^\SC)\)), while in \cite{FN11} the monoid is what
is commonly called an
\notion{\(L\)-monoid} associated with a single fixed dominant cocharacter
\(\lambda\) of \(G\) (denoted by \(\FRM(\lambda)\)). In this book we will
consider an arbitrary very flat monoid
\(\FRM\), thus unifying different formulations into one framework. The main
reason for doing this is not merely technical and is only apparent when we start
to consider endoscopic groups.

\subsection{}
\label{sub:intro_future_generalizations}
We will continue discussing mH-fibrations in the next section,
and close this one with a few comments on potential generalizations.
Instead of limiting ourselves to just \(G\) acting on \(G\) of \(\FRM\), we
may consider a \(G\)-space \(M\) and an equivariant partial compactification
\(\FRM\). Suppose we also have some group \(Z_\FRM\) of multiplicative type acting on
\(\FRM\) commuting with the \(G\)-action, then we may also consider the
Hitchin-type fibration associated with \(\Stack*{\FRM/G\x Z_\FRM}\to
\Stack*{(\FRM\git G)/Z_\FRM}\).

In \cite{SV17} and later \cite{BZSV24}, the authors proposed a unified framework for
relative trace formulae using spherical varieties. In the group case, \(G\)
is viewed as a \(G\x G\)-spherical variety, and \(\Env(G^\SC)\) is a
horospherical contraction of \(G^\SC\) as a spherical variety and also the
spectrum of the Cox ring of the wonderful compactification of \(G^\AD\) (one
may consult \cite{Ti11} for terminologies). Therefore, a particularly promising
generalization is the following: suppose we have a spherical homogeneous space
\(G/H\) which admits a wonderful compactification \(\bar{G/H}\), we may consider
the spectrum \(\FRM\) of the Cox ring of \(\bar{G/H}\). In this case, there
is a canonical torus \(Z_\FRM\) acting on \(\FRM\) commuting with the action of
\(G\), and so it descends to the abelianization \(\FRA_\FRM\defeq\FRM\git G\)
similar to the monoid case. The abelianization \(\FRA_\FRM\) is known to be
an affine space. Thus, we may consider the Hitchin-type fibration associated
with any subgroup (for example, isotropy group \(H\) itself) of \(G\) acting on
\(\FRM\).

\section{Local Model of Singularity}
\label{sec:intro_local_model_sing}

After we define mH-fibrations using reductive monoids, a lot of their geometric
properties can be laid out in a fashion parallel to \cite{Ng10}, often
without too much modification to the argument. Doing so requires detailed
study of the invariant theory of the \(G\)-action on \(\FRM\), which is largely
done in \cites{Bo15,Ch22}, and we only record the results and
bridge some remaining gaps. One main exception to all this ``house-keeping''
is studying the deformation of multiplicative Higgs (mHiggs) bundles, which has not been
done so far in a satisfactory way.

\subsection{}
In \cite{Ng10}*{\S~4.14}, Ng\^o proved that the total stack of Hitchin
fibration is smooth under the assumption that the twisting line bundle is very
ample. His result is based on the tangent-obstruction theory of Higgs
bundles due to Biswas and Ramanan in \cite{BiRa94}. The group case is more
complicated and the total stack of mH-fibration is usually not smooth. As a
result, we must replace constant sheaf by the intersection complex in our later
study of cohomologies. The question is then what the singularities of the
total stack look like, and in \cite{FN11}*{Conjecture~4.2} the authors
conjectured that the intersection complex is exactly the globalized version of
some Satake sheaf. We will prove this conjecture essentially in full in this book by
establishing a tangent-obstruction theory for mHiggs bundles. Below is a rather
vaguely stated corollary of
\Cref{thm:local_singularity_model_weak,thm:local_singularity_model_main}, and
proves \cite{FN11}*{Conjecture~4.2}:
\begin{theorem}
    \label[theorem]{thm:simplified_local_singularity_model_in_intro}
    Let \(h_X\colon\cM_X\to\cA_X\) be the mH-fibration associated with \(G\),
    \(\FRM\), and curve \(X\), and let \(\cB_X\) be its moduli of boundary
    divisors. The pullback complex \(\ev^*\IC_{\Stack{\GASch_X}}\) is, up to
    shifts and Tate twists, isomorphic to
    the intersection complex on \(\cM_X\) after restricting to a ``very
    large'' open subset of \(\cA_X\). Here \(\Stack{\GASch_X}\to\cB_X\) is certain
    truncated global Hecke stack induced by monoid \(\FRM\) and parametrized by
    boundary divisors in \(\cB_X\)
    (see \Cref{sec:global_affine_schubert_scheme}), and
    \(\ev\colon\cM_X\to\Stack{\GASch_X}\) is an appropriately defined evaluation
    map over \(\cB_X\).
\end{theorem}

We avoid explaining more notations in
\Cref{thm:simplified_local_singularity_model_in_intro} because they are
technical and will drag this introduction longer than necessary. However, an
explanation of the phrase ``very large'' is due in order to understand the novelty
in this result. In \cite{Bo17}, a weaker version of
\Cref{thm:simplified_local_singularity_model_in_intro} was proved using
\textit{ad hoc} method. It did not try to establish a tangent-obstruction theory
in that paper, and thus had to impose a much stronger technical condition that
is even more difficult to state.

The gist is that Bouthier's version of ``very
large open subset'' is not large enough and \emph{never} includes any stratum coming
from endoscopy (outside some very rare cases, essentially only when
\(G=\SL_2\)). It makes the old result
unusable in tackling the fundamental lemma. In our new results, the size of
the open subset is on par with the analogous result in \cite{Ng10}.

\subsection{}
The most important implication of
\Cref{thm:simplified_local_singularity_model_in_intro} is that the
intersection complex on \(\cM_X\) is essentially a product of local
\(\IC\)-sheaves on affine Schubert varieties. More generally, if we have at each
place \(v\) of \(X\) a local Satake sheaf \(\bcQ_v\), we can assemble them into
a pure perverse sheaf \(\bcQ\) on \(\cM_X\). This is essentially a
geometrization of representations of global \(L\)-groups. Once we start to consider
endoscopic groups (and their mH-fibrations), it will be very clear how to transfer
\(\bcQ\) to endoscopic side, and we will be able to geometrize the stabilization
process as well, similar to \cite{Ng10}*{\S~6.4}.

\subsection{}
We briefly go back to the hypothesized generalization involving a
wonderful variety \(\bar{G/H}\) at the end of
\Cref{sec:on_hitchin_type_fibrations}. Since in this case we can still define
the moduli stack of boundary divisors using \(\Stack*{\FRA_\FRM/Z_\FRM}\), we
should be able to define \(\Stack*{\GASch_X}\) analogously (see
\Cref{sec:global_affine_schubert_scheme}). Naturally, we expect to have a
statement of geometric Satake isomorphism and a local model of singularity as
well, but to do so we must first have a good framework dealing with
\(\ell\)-adic sheaves on the quotient of an arc space by an arc group,
which does not seem to exist in published form at the time of writing.


\section{Geometric Transfer} 
\label{sec:dual_groups_and_geometric_transfer_map}

Utilizing universal monoid \(\FRM=\Env(G^\SC)\) and the associated mH-fibration,
A.~Bouthier and
J.~Chi in \cites{Bo15,BoCh18}, and later J.~Chi in \cite{Ch22}
were able to prove the dimension
formula for the multiplicative affine Springer fibers of \(G\).
As we have discussed before, A.~Bouthier in \cite{Bo17} also obtained global
results by studying the singularities of the total stack of
mH-fibrations. However, once one
attempts to connect mH-fibrations of \(G\) to those of endoscopic group \(H\),
things quickly fall apart if one only sticks with \(\Env(G^\SC)\) and
\(\Env(H^\SC)\).

\subsection{}
In Lie algebra case, if \(H\) happens to be a subgroup of \(G\), then we have the
natural map \(\Stack{\La{h}/H}\to\Stack{\La{g}/G}\), and since both \(\La{h}\) and
\(\La{g}\) are reductive Lie algebras, the general results for Hitchin
fibrations can be applied to both sides. In general, \(H\) is not a subgroup of
\(G\), nevertheless there is still a canonical
map from \(\FRc_H=\La{h}\git H\) to \(\FRc\), which induces a
closed embedding from the Hitchin base of \(H\) to that of \(G\). This map is
all we need for establishing the geometric (spherical) transfer in Lie algebra case.

In group case, however, there is usually no way to directly relate
\(\Env(H^\SC)\) to \(\Env(G^\SC)\) even if \(H\) is a subgroup of \(G\), and
more generally there is no map from \(\FRC_{\Env(H^\SC)}=\Env(H^\SC)\git H\) to
\(\FRC_{\Env(G^\SC)}=\Env(G^\SC)\git G\) either. Using the general theory of
(not necessarily very flat) reductive monoids, one can still produce some
canonical monoid \(\FRM_H'\) for \(H\) such
that \(\FRM_H'\git H\) does map to \(\FRC_{\Env(G^\SC)}\), and if \(H\) happens
to be a subgroup of \(G\), \(\FRM_H'\) is just the closure of the group
\(H'\subset\FRM^\x\) where \(H'\) has the same semisimple type as \(H\).
However, such monoid \(\FRM_H'\) in general is not a very flat monoid, so
its corresponding mH-fibration will not have very good geometric properties.

\subsection{}
This, of course, is not a coincidence, and (in my opinion) is where
the story becomes fascinating. The picture is better understood once we
move to the dual group side, and the solution to this difficulty will manifest
itself once we do so.

The fundamental lemma is usually stated as an equality between the evaluation of
orbital integrals at a spherical function. In Lie algebra case, the
function is already given, namely the characteristic function on \(\La{g}(\cO)\)
and that on \(\La{h}(\cO)\),
so all we need is to figure out how to match conjugacy classes.
In group case, however, we have to geometrize both the matching of conjugacy
classes and the matching of functions. The latter
has close connection to the representations of \(L\)-groups thanks to
geometric Satake isomorphism.

\subsection{}
Temporarily going back to the local setting, we can see that the effect of the
\(\Gm\)-twisting in the Lie algebra case is that we are not really considering
the characteristic function on \(\La{g}(\cO)\), but rather those on scaled subsets
\(\pi^{-d}\La{g}(\cO)\subset \La{g}(F)\) for arbitrary \(d\ge 0\). Since these
subsets do not look very different from 
\(\La{g}(\cO)\), we are still effectively only considering a single function. In
group case, however, things get complicated because the analogue
would be
Cartan decomposition and Satake functions on the closure of double cosets
\(G(\cO)\pi^\lambda G(\cO)\), and these functions certainly do not look alike.

\subsection{}
Recall that to each mH-fibration we also have the moduli \(\cB_X\) of boundary
divisors, which are essentially divisors valued in dominant cocharacters of
\(G\). The cocharacters \(\lambda\) that can appear depends on the choice of monoid
\(\FRM\), essentially (but not strictly) in a bijective way.
\Cref{thm:simplified_local_singularity_model_in_intro} implies that
we can encode Satake sheaves that are supported on
affine Schubert variety \(\Gr_G^{\le\lambda}\) into \(\cM_X\).
Using geometric Satake isomorphism, by carefully choosing \(\FRM\), we can
encode any representation of \(L\)-group \(\LD{G}\) into mH-fibrations.

\subsection{}
If we have an admissible embedding, otherwise known as \(L\)-embedding, of \(L\)-groups
\(\LD{H}\to\LD{G}\), we can ``transfer'' a representation of \(\LD{G}\) to
\(\LD{H}\), namely by restriction. Since we have a (roughly bijective)
correspondence between monoid \(\FRM\) and representations of \(L\)-groups,
it gives hints to how to relate mH-fibrations for \(G\) to those
for \(H\).

With this consideration about \(L\)-group in mind, we will be able to construct for any
given monoid \(\FRM\) of \(G\) a
canonically associated monoid \(\FRM_H\) (see \Cref{eqn:def_endoscopic_monoid})
of \(H\) such that there exists a canonical map of quotient space
\begin{align}
    \FRC_{\FRM,H}=\FRM_H\git H\longto\FRC_\FRM.
\end{align}

\subsection{}
Finding the correct monoid, however, is not the end of the story.
The reason is very simple: mH-base \(\cA_X\) is constructed from
stack \(\Stack{\FRC_\FRM/Z_\FRM}\), but the monoids \(\FRM\) and \(\FRM_H\) have
different centers, and there is usually no map from \(Z_{\FRM,H}\) (the center of
\(\FRM_H^\x\)) to \(Z_\FRM\) whatsoever. For example, suppose \(G=\SL_2\) and
\(\FRM=\Mat_2\), and \(H\) is a \(1\)-dimensional torus, then \(\FRM_H=\bbA^2\) is the
monoid of diagonal matrices in \(\FRM\). In this case \(Z_{\FRM,H}\) is the
group of invertible diagonal matrices while \(Z_\FRM\) is the scalar
matrices, so there is no map from \(\Stack*{\FRC_{\FRM,H}/Z_{\FRM,H}}\) to
\(\Stack*{\FRC_\FRM/Z_\FRM}\). In general, \(Z_{\FRM,H}\) will be much
larger than \(Z_\FRM\).

The solution to this problem is very simple: we always have a canonical
subgroup \(Z_\FRM^\kappa\subset Z_{\FRM,H}\) that maps canonically onto
\(Z_\FRM\).  Essentially, \(Z_\FRM^\kappa\) is the ``preimage'' of
\(Z_\FRM\) in \(\FRM_H\), which is necessarily contained in \(Z_{\FRM,H}\).
In our example above where \(\FRM=\Mat_2\) and \(\FRM_H=\bbA^2\),
\(Z_\FRM^\kappa\) is equal to \(Z_\FRM\). In general,
the group \(Z_\FRM^\kappa\) is again much larger
than \(Z_\FRM\). We now have a map
\begin{align}
    \label{eqn:intro_geom_transfer_map}
    \Stack{\FRC_{\FRM,H}/Z_\FRM^\kappa}\longto \Stack{\FRC_\FRM/Z_\FRM},
\end{align}
compatible with the map between classifying stacks
\(\BG{Z_\FRM^\kappa}\to\BG{Z_\FRM}\).

The map \eqref{eqn:intro_geom_transfer_map} is analogous to the map \(\Stack{\FRc_H/\Gm}\to\Stack{\FRc/\Gm}\) in Lie
algebra case.
Since for both \(\La{g}\) and \(\La{h}\) the twisting is given by
\(\Gm\), we may pick a single \(\Gm\)-torsor of large degree, as is done by most
literature. In group case, we cannot afford to fix
any torsor because \(Z_\FRM^\kappa\) and \(Z_\FRM\) are different, and all we
can do is to apply the mapping stack functor to the whole map
\eqref{eqn:intro_geom_transfer_map}. It turns out to be the correct formulation.
In short, we obtain a finite unramified map
\begin{align}
    \nu_\cA\colon\cA_{H,X}^\kappa\longto \cA_X,
\end{align}
which, unlike the Lie algebra case, is usually \emph{not} a
closed embedding. Strictly speaking, here we need to replace \(\cA_X\) by an open
subset (the reduced locus \(\cA_X^\heartsuit\)), but we ignore such issues in this introduction
for simplicity.

\subsection{}
Using map \eqref{eqn:intro_geom_transfer_map}, we will be able to
match conjugacy classes of \(G\) and \(H\). The transfer of functions is
more complicated because it involves global \(L\)-embeddings. For simplicity, suppose we have
an irreducible representation of the dual group \(\dual{\bG}\) with
highest-weight \(\lambda\), it decomposes into a direct sum of irreducible
representations of \(\dual{\bH}\). Note that for  a \(\dual{\bH}\)-highest
weight \(\lambda_H\) therein, it may have multiplicity greater than \(1\).

As it turns out, the choice of
monoids \(\FRM\) and \(\FRM_H\), aside from conjugacy classes, will also match
\(\lambda\) with the set of \(\lambda_H\) in the decomposition. We can then restore the multiplicity of
\(\lambda_H\) by tensoring with an appropriate multiplicity space, viewed as
a constant sheaf on \(\Gr_H^{\le\lambda_H}\). We will refrain from talking
about the technical details involving global \(L\)-embeddings, however. See
\Cref{sec:transfer_of_Satake_sheaves} for more details.


\section{Support Theorem} 
\label{sec:support_theorem_and_beyond}

Continue following the same general strategy in \cite{Ng10},
the next and perhaps most crucial piece is the support theorem.
We have already discussed several difficulties and features unique to the group
case, and this part is no exception.

\subsection{}
Ng\^o's support theorem in \cite{Ng10} has three main ingredients aside from the
geometric transfer map already established. The first is the symmetry of Hitchin
fibers, namely a Picard stack acting on the Hitchin fibration induced by the
regular centralizer. Then from this symmetry one can stratify the Hitchin base
in several ways. These two ingredients imply that Hitchin fibrations satisfy
certain properties making them so-called \notion{\(\delta\)-regular abelian
fibrations}.
Ng\^o then proves a support theorem for a general \(\delta\)-regular abelian
fibration with a smooth total stack. The support theorem at this level of generality is a very powerful
tool, but on its own it is not quite enough for the fundamental lemma, notably
because it only establishes an upper bound for the supports.

The reason that a general support theorem cannot pin down the exact set of
supports is because the argument is ultimately based on dimension
considerations. It is essentially a vastly optimized Goresky--MacPherson
inequality (see \cite{Ng10}*{\S~7.3}). As such, they cannot
deal with other contributions that are discrete in nature. In the case of
\(\delta\)-regular abelian fibrations, this extra ``discrete contribution''
comes from the irreducible components of the fibers.

Ng\^o is able to
determine the exact amount of such discrete contribution for Hitchin
fibrations using a concrete description of the irreducible components, and shows
that the only such contribution comes from endoscopy.

\subsection{}
In group case, the first two ingredients remain largely unchanged. First, the
regular centralizer can be defined similarly, albeit with more technical
efforts, so that there is a Picard action. Second, the stratifications on the
mH-base can be done by a parallel fashion, showing that mH-fibrations are
\(\delta\)-regular abelian fibrations.

There is a new ingredient necessary for mH-fibrations: the local model of
singularity. Unlike its connection to \(L\)-groups as we already discussed, its
role here is rather simple to explain: suppose \(f\colon X\to Y\) is a map of
varieties over \(\bar{k}\) of pure relative dimension \(d\), then the
cohomological degree of \(\RDF{f}_{!}\Qlb\) is bounded above by \(2d\), and at
degree \(2d\) its stalks are spanned by irreducible components of the fibers of
\(f\). However, if we replace \(\Qlb\) by the intersection complex \(\IC_X\),
then \(\RDF{f}_{!}\IC_X\) is not necessarily cohomologically bounded above by
\(2d-\dim{X}\). Using local model of singularity and a simple spectral sequence,
we will show that such upper bound does hold for mH-fibrations, and at top
degree the contribution can be described using irreducible components of the
fibers. The idea of using spectral sequence was originally due to A.~Bouthier in
an unpublished note.

\subsection{}
With these input, we are able to prove a generalization of Ng\^o's support
theorem for mH-fibrations, which reduces the problem to describe irreducible
components of mH-fibers. Unlike Lie algebra case, the description of the
irreducible components is much more complicated, and there are additional
contributions other than endoscopy. In order to understand it, we use a simple
example to illustrate the idea.

Suppose we have two effective Cartier divisors
\(D'<D\) on \(X\), we will have a natural inclusion of section spaces
\begin{align}
    \RH^0(X,\cO_X(D'))\subset\RH^0(X,\cO_X(D)).
\end{align}
Similarly, for (ordinary) Hitchin fibrations, we can embed a Hitchin base into another Hitchin
base with larger degree:
\begin{align}
    \RH^0(X,\FRc_{D'})\subset\RH^0(X,\FRc_D).
\end{align}
We call the image of this embedding an \notion{inductive stratum} in the
latter space. Of course, such strata are not very important in Lie algebra
case, but the story is different for mH-fibrations.

There is a similar partial order on boundary divisors in
\(\cB_X\) induced by the dominance order on cocharacters. For a fixed boundary
\(b\in\cB_X\), let \(\cA_b\) be the fiber of \(\cA_X\) over \(b\), then we have
a natural inclusion \(\cA_{b'}\subset\cA_b\) for \(b'<b\), defined in a similar
way. The union of these linear subspaces from boundary
subdivisors induce a stratification on \(\cA_b\), and we can let \(b\)
vary in \(\cB_X\) to obtain a stratification on whole \(\cA_X\). These
strata are the inductive strata in mH-bases.

Consider the endoscopic locus
\(\cA^{\kappa}\subset\cA_X\) that is the union of the images of mH-bases
\(\cA_{H,X}^\kappa\)
for all endoscopic groups \(H\) associated with endoscopic parameter \(\kappa\).
Let \(S_H^\kappa\) be the set of images in \(\cA_X\) of inductive strata in
\(\cA_{H,X}^\kappa\), so it is a set of closed subsets of \(\cA_X\).
Now we can state a rough form of support theorem for
mH-fibrations (see
\Cref{thm:stable_support_of_mH_main,thm:support_endoscopic_main} for precise
statements):
\begin{theorem}
    \label[theorem]{thm:intro_supprot_theorem}
    The complex \(h_{X,*}\IC_{\cM_X}\) is pure perverse and decomposes into a direct sum
    of \(\kappa\)-isotypic components \((h_{X,*}\IC_{\cM_X})_\kappa\) where \(\kappa\)
    ranges over possible endoscopic parameters, and the set of supports of
    perverse sheaves \(\PH^\bullet(h_{X,*}\IC_{\cM_X})_\kappa\) is a subset of
    the union of \(S_H^\kappa\), where \(H\) ranges over all possible endoscopic
    groups associated with \(\kappa\).
\end{theorem}

\subsection{}
As it turns out, the contribution of inductive strata can be described locally
using \emph{unramified} multiplicative affine Springer fibers, whose irreducible
components have an explicit connection with MV-cycles. This leads to the final
component of the proof of \Cref{thm:FL_intro} which we discuss in the next
section.

\section{Transfer Factor and Kashiwara Crystals}
\label{sec:intro_TF_and_crystals}

In \cite{Ng10}, any problem involving the transfer factor is avoided because it
has already been studied in \cite{Ko99}, and can be conveniently skipped over if
one utilizes the Kostant section. In group case, there does not seem to be any
corresponding result. T.~Hales gives a simplified definition of transfer factors
for unramified groups in \cite{Ha93}, but it is still far more complicated than
the Lie algebra case in \cite{Ko99}. In case of unramified conjugacy classes,
however, things can be more explicitly described by studying the Frobenius
action on the set of irreducible components of multiplicative affine Springer
fibers.

\subsection{}
Recall from earlier we mentioned the proof of fundamental lemma for
unramified elements by \cite{GKM04} that is conditional on the purity conjecture
of affine Springer fibers. Even though it is still a conjecture at the time of
writing, it is widely believed to be true, and it is reasonable to expect the
same is true for their multiplicative analogues. Consequently, it is easy to
believe that fundamental lemma can be more or less deduced from its asymptotic
form by throwing away ``lower order terms'' and only considering the
contribution from irreducible components. Of course, the latter is no longer a
mere belief thanks to support theorem, even though we do not have local purity
results.

In \cite{Ch22}, Chi described the geometric irreducible components of unramified
multiplicative affine Springer fibers by relating them to MV-cycles, and the latter
give canonical bases of representations of the dual group. What is left
is to describe the Frobenius action on those components and the
\(\kappa\)-twisting. As MV-cycles form canonical bases in
representations, a naive guess would be that the components on \(G\)-side
are related to those on \(H\)-side via the restriction of representations from
\(\LD{G}\) to \(\LD{H}\). This is, of course, too naive because it seems
impossible to incorporate \(\kappa\)-twisting with the restriction
functor.

\subsection{}
The correct answer is to not regard MV-cycles as basis elements in
representations, but in the associated Kashiwara crystals, or crystals for
short. We will not review the theory of crystals here, but heuristically they
capture the combinatorial essence of highest-weight representations, and
any highest-weight representation has a canonically associated (normal)
crystal. One can sometimes think of a crystal basis as a set of particularly
nicely chosen
vector basis in the corresponding representation (although this point of
view is not very accurate conceptually).

An important advantage of crystal bases is that they afford a canonical Weyl
group action, while in contrast there is no Weyl group action on
representations. Moreover, the restriction to standard Levi
subgroups becomes very straightforward to describe. The tradeoff, however, is
that restriction functor in general becomes much more complicated: given a crystal and a
(non-Levi) subgroup, there is no straightforward way to obtain a crystal of the
subgroup by ``restricting'', and extra work needs to be done depending on the
situation. This plays to our favor, because for endoscopic \(L\)-embeddings, this extra
work is precisely the \(\kappa\)-twisted matching of irreducible components of
multiplicative affine Springer fibers of \(G\) and of \(H\).

\subsection{}
The transfer factor defined in \cite{LS87} is a product of five terms:
\begin{align}
    \Delta=\Delta_{\symup{I}}\Delta_{\symup{II}}\Delta_{\symup{III}_1}\Delta_{\symup{III}_2}\Delta_{\symup{IV}}.
\end{align}
The last term \(\Delta_{\symup{IV}}\) is just cohomological twist, and
\(\Delta_{\symup{III}_1}\) is the coefficient in defining \(\kappa\)-orbital
integrals and also relatively straightforward. The remaining terms, however, are
rather complicated, and individually they depend on choices of certain auxiliary
data. Fortunately, for unramified conjugacy classes, there are choices that are
more or less canonical, and by using those choices, each individual
term can be visualized in multiplicative affine Springer fibers.

As a result,
we obtain an \notion{asymptotic fundamental lemma} described using
crystals. We will refrain from stating the precise result due to the complexity
in notations (see
\Cref{thm:asymptotic_FL}).
Even though this asymptotic version of the fundamental lemma looks a lot simpler than
the original, it still seems quite difficult to prove. Nevertheless, there are
special cases that are simple enough to prove by hand, and those cases are
abundant enough for us to prove the fundamental lemma by induction. By reversing
the argument, we can prove the asymptotic fundamental lemma in full as
well.

\subsection{}
It is interesting to note that in this book, the asymptotic fundamental lemma
is deduced from the full fundamental lemma, even though the former looks much
simpler. It would be interesting to find a purely local proof of the asymptotic
fundamental lemma.

As a final note, we expect our framework of mH-fibrations can be generalized
to a large class of spherical varieties hence used to tackle problems in
relative trace formulae in general. For each case of such generalizations, we
also expect that crystals will play a similar role in describing inner
structures of \(L\)-packets, including but not limited to endoscopy.


\section{Structure of this Book}
Here is a summary of the content in each chapter. Many chapters are arranged
very much like the way in \cite{Ng10} not only because the logical structure
therein is well-organized, but it is also more convenient to compare with Lie
algebra case.

In \Cref{chap:reductive_monoids_and_invariant_theory} we introduce some basic
notations and setups in the absolute setting. We will review some standard facts
about reductive monoids and its invariant theory under the adjoint action. The
highlight of this chapter is the construction of the endoscopic monoid and the
transfer map at invariant-theoretic level. We will also give a formal statement
of the fundamental lemma at the end.

\Cref{chap:multiplicative_valuation_strata} studies the multiplicative
valuation strata. It is a purely technical tool for this book and is mostly
self-contained and can be read by itself. \Cref{chap:global_constructions} is
also auxiliary. It collects many global constructions needed for our
discussions, including the moduli of boundary divisors and a construction of
global affine Schubert schemes using reductive monoids. The latter is somewhat
interesting by itself as by using monoids we are able to define affine Schubert
schemes without referring to affine Grassmannians.

\Cref{chap:generalized_affine_springer_fibers} reviews constructions and
properties of multiplicative affine Springer fibers for groups.
Most of the chapter is very
similar to what is done in the Lie algebra case, with two most notable new
additions: the connection with MV-cycles and the connection with Kashiwara crystals.
We also discuss how to translate fundamental lemma into the language of monoids.

\Cref{chap:multiplicative_hitchin_fibrations} studies the constructions of
mH-fibrations and its basic properties. It is more or less a global parallel to
\Cref{chap:generalized_affine_springer_fibers} and also similar to
the Lie algebra case. The local model of singularity is also stated here, but
the proof is postponed to \Cref{chap:deformation}, where we study the
deformations of mHiggs bundles. In \Cref{chap:deformation}, we
first review some general facts about deformations of mapping stacks, and then
use it to study mH-bases. Using certain duality between the mH-base and
mH-fibers, we study the deformation of mHiggs bundles by modifying the
argument for mH-bases.

In \Cref{chap:stratifications}, we discuss various
useful stratifications on the mH-base. The particularly new feature is the
notion of inductive subsets, which will be important in our support theorem.
In \Cref{chap:cohomologies}, we study cohomologies based on those
stratifications. The global endoscopic transfer, both at the level of conjugacy
classes and functions, in the form of global Satake sheaves, will be explained
here. We will also state the geometric stabilization theorem.
Another interesting object we introduce here is a new kind of Hecke-type stack,
using which we can upgrade the traditional product formula over a point into a
family. The latter will be particularly useful when we study the top ordinary
cohomologies, as studying them beyond just the rank of stalks will be necessary.

\Cref{chap:support_theorem} contains the proof of the support theorem. The proof
itself is surprisingly close to the case of Lie algebra with perhaps the only
important new ingredient being the local model of singularity proved in
\Cref{chap:deformation}. However, due to mH-fibers having
more complicated description of irreducible components, the implication of the
support theorem differs from the Lie algebra case, and this is where inductive
subsets and the connection with MV-cycles and crystals come into play.

In \Cref{chap:counting_points}, we review the point-counting framework done by
the last chapter of \cite{Ng10}, and extend it slightly to make it more
convenient. After that, in \Cref{chap:proof_of_fundamental_lemma}, we will
finally be able to prove the fundamental lemma after several reductions and
inductions.

There are two appendices. In \Cref{chap:Review_on_Transfer_Factors} we review
the definition of standard transfer factor defined in \cite{LS87}. The language
is slightly altered to be coherent with the main part of this book. In
\Cref{chapA:Review_of_Kashiwara_Crystals} we give some basic facts about
Kashiwara crystals for the reader's convenience.

\section{Acknowledgement}
This project was initiated during my time as a PhD student at The University of
Chicago, and a substantial portion was developed in my PhD thesis.
I am extremely grateful to my PhD advisor Ng\^o Bao Ch\^au for suggesting the entire
project and his constant support. I thank Jingren Chi for helping me during the
early stage of this project by explaining many details of his thesis and
answering numerous questions. I thank Minh-Tam Trinh for his constant
encouragement, interest, and many inspiring discussions. I
also thank Tsao-Hsien Chen, Thomas Hameister, Benedict Morrissey, Yiannis
Sakellaridis, Sug Woo Shin, Zhiwei Yun, and Xinwen Zhu for their interest in
this project and delightful discussions. Furthermore, I am especially thankful to Akshay
Venkatesh for helpful suggestions on structuring this book. Last but not
least, I thank the anonymous referees for their careful reading of the
manuscript, and many great pieces of advice that improved parts of this work.

I would also like
to express my gratitude towards the Department of Mathematics at The
University of Chicago, without whose generous support during my PhD study this
book would certainly not have been possible. The
preparation of this manuscript was finished during my stay at the Institute for
Advanced Study as a member, and I am grateful for their hospitality and providing such
stimulating academic atmosphere.
This material is partially based upon work supported by the National
Science Foundation under Grant No.~DMS-1926686.

\chapter{Reductive Monoids and Invariant Theory}%
\label{chap:reductive_monoids_and_invariant_theory}

We start by reviewing some facts about the invariant theories of the
adjoint action of a reductive group \(G\) on the simply-connected cover
\(G^\SC\) of its derived subgroup and on a very flat reductive
monoid \(\FRM\) associated with \(G^\SC\). The new results in this chapter are
mostly contained in \Cref{sec:endoscopic_groups_inv_theory}, where the
canonical endoscopic monoid \(\FRM_H\) associated with any given monoid \(\FRM\)
and endoscopic group \(H\) is defined. The endoscopic monoid plays a key role
throughout this book. Proofs are mostly omitted for
well-known facts, but references will be given when possible.

\section{Quasi-split Forms}
\label{sec:Quasi-split Forms}

In this section, we set up some basic notations used throughout this book. In
particular, we introduce the notion of quasi-split forms.

\subsection{}
Let \(k\)
    \nomenclature[\(k \)]{\(k\)}{a finite field}
be a finite field with \(q\)
    \nomenclature[\(q \)]{\(q\)}{the cardinality of \(k\)}
elements and \(\bar{k}\)
    \nomenclature[\(k_bar \)]{\(\bar{k}\)}{a fixed algebraic closure of \(k\)}
be a fixed algebraic closure of \(k\). Let \(\bG\)
    \nomenclature[\(G"bold \)]{\(\bG\)}{a split reductive group}
be a split connected reductive group defined over \(k\) of rank \(n_G\)
    \nomenclature[\(n_G \)]{\(n_G\)}{the rank of \(\bG\) or \(G\)}
and semisimple rank \(r\).
    \nomenclature[\(r \)]{\(r\)}{the semisimple rank of \(\bG\) or \(G\)}
We assume \(p=\Char(k)\)
    \nomenclature[\(p \)]{\(p\)}{the characteristic of \(k\)}
is larger than twice the Coxeter number of \(\bG\). We fix once and for all a
split maximal \(k\)-torus \(\bT\)
    \nomenclature[\(T"bold \)]{\(\bT, T\)}{a fixed maximal torus of \(\bG\) or \(G\)}
of \(\bG\) and a Borel subgroup \(\bB\)
    \nomenclature[\(B"bold \)]{\(\bB, B\)}{a fixed Borel subgroup of \(\bG\) (resp.~\(G\)) containing \(\bT\) (resp.~\(T\))}
containing \(\bT\). Let \(\bU\)
    \nomenclature[\(U"bold \)]{\(\bU,U\)}{the unipotent radical of \(\bB\) or \(B\)}
be the unipotent radical of \(\bB\). Let \((\CharG,\Roots, \CoCharG,\CoRoots)\)
    \nomenclature[\(X"bbold' \)]{\(\CharG\)}{the character group of \(\bT\)}
    \nomenclature[\(X"bbold'_T \)]{\(\CharG(T)\)}{the character group of a torus \(T\)}
    \nomenclature[\(X"bbold'check \)]{\(\CoCharG\)}{the cocharacter group of \(\bT\)}
    \nomenclature[\(X"bbold'check_T \)]{\(\CoCharG(T)\)}{the cocharacter group of a torus \(T\)}
    \nomenclature[\(Phi' \)]{\(\Roots\)}{the set of roots}
    \nomenclature[\(Phi'check \)]{\(\CoRoots\)}{the set of coroots}
be the root datum associated with \((\bG,\bT)\), and
\(\SimRts=\{\Rt_1,\ldots,\Rt_r\}\)
    \nomenclature[\(Delta \)]{\(\SimRts\)}{the set of simple roots}
    \nomenclature[\(alpha_i \)]{\(\Rt_i\)}{the simple roots of \(\bG\)}
(resp.~\(\PosRts\))
    \nomenclature[\(Phi'_+ \)]{\(\PosRts\)}{the set of positive roots}
the set of simple (resp.~positive) roots determined by \(\bB\), and let
\(\SimCoRts=\{\CoRt_1,\ldots,\CoRt_r\}\)
    \nomenclature[\(Delta'check \)]{\(\SimCoRts\)}{the set of simple coroots}
    \nomenclature[\(alpha'check_i \)]{\(\CoRt_i\)}{the simple coroots of \(\bG\)}
(resp.~\(\PosCoRts\))
    \nomenclature[\(Phi'check_+ \)]{\(\PosCoRts\)}{the set of positive coroots}
the simple (resp.~positive) coroots. Let \(\CharG_+\)
    \nomenclature[\(X"bbold'_+ \)]{\(\CharG_+\)}{the set of dominant characters of \(\bT\)}
(resp.~\(\CoCharG_+\))
    \nomenclature[\(X"bbold'check_+ \)]{\(\CoCharG_+\)}{the set of dominant cocharacters of \(\bT\)}
be the set of dominant characters (resp.~cocharacters). Let \(\bW=W(\bG,\bT)\)
    \nomenclature[\(W \)]{\(\bW, W\)}{the Weyl group of \(\bT\) (resp.~\(T\)) in \(\bG\) (resp.~\(G\))}
    \nomenclature[\(W_G_T \)]{\(W(G,T)\)}{the Weyl group of a torus \(T\) in a group \(G\)}
be the Weyl group, and let \(w_0\in\bW\)
    \nomenclature[\(w_0 \)]{\(w_0\)}{the longest Weyl element}
be the longest element with respect to \(\bB\).

\subsection{} 
Let \(\bG^\Der\)
    \nomenclature[\(G^der \)]{\(G^\Der\)}{the derived subgroup of a group \(G\)}
be the derived subgroup of  \(\bG\), and let \(\bG^{\SC}\)
    \nomenclature[\(G^sc \)]{\(G^\SC\)}{the simply-connected group of \(G^\Der\)}
    \nomenclature[\(.sc \)]{\((\cdot)^\SC, (\cdot)_\SC\)}{the (pre)image in \(G^\SC\)}
(resp.~\(\bG^{\AD}\))
    \nomenclature[\(G^ad \)]{\(G^\AD\)}{the adjoint group of a reductive group \(G\)}
    \nomenclature[\(.ad \)]{\((\cdot)^\AD, (\cdot)_\AD\)}{the (pre)image in \(G^\AD\)}
be the universal cover (resp.~adjoint quotient) of \(\bG^\Der\). We also use
\(H^\SC\) (resp.~\(H^\AD\)) to denote the preimage (resp.~image) of any subset
\(H\subset \bG\) in \(\bG^\SC\) (resp.~\(\bG^\AD\)). Let
\(\{\Wt_1,\ldots,\Wt_r\}\)
    \nomenclature[\(pi_var_i \)]{\(\Wt_i\)}{the fundamental weights of \(\bG^\SC\)}
be the set of fundamental weights of \(\bG^{\SC}\), and \((\rho_i,V_{\Wt_i})\)
    \nomenclature[\(rho_i \)]{\(\rho_i\)}{the \(i\)-th fundamental representation (Weyl module) of \(\bG^\SC\)}
    \nomenclature[\(V_pi_var_i \)]{\(V_{\Wt_i}\)}{the \(k\)-space of \(\rho_i\)}
be the corresponding Weyl module with highest-weight \(\Wt_i\). Let
\(\{\CoWt_1,\ldots,\CoWt_r\}\)
    \nomenclature[\(pi'check'var_i \)]{\(\CoWt_i\)}{the fundamental coweights of \(\bG^\SC\)}
be the set of fundamental coweights. Let \(\rho\)
    \nomenclature[\(rho \)]{\(\rho\)}{the half-sum of positive roots}
(resp.~\(\dual{\rho}\))
    \nomenclature[\(rho'check \)]{\(\dual{\rho}\)}{the half-sum of positive coroots}
be the half-sum of all positive roots (resp.~coroots).

\subsection{}
We fix a \(k\)-pinning \(\bSPL=(\bT,\bB,\bFx_+\defeq \Set{\bU_\alpha}_{\alpha\in\SimRts})\) 
    \nomenclature[\(P_ \)]{\(\bSPL,\SPL\)}{a fixed pinning of \(\bG\) or \(G\)}
of \(\bG\), where \(\bU_\alpha\colon \Ga\to \bG\)
    \nomenclature[\(U_alpha"bold \)]{\(\bU_\Rt\)}{the root vector associated with simple root \(\Rt\)}
is a one-parameter unipotent group associated with simple root \(\alpha\). The
associated groups \(\bG^\Der\), \(\bG^\SC\) and \(\bG^\AD\) all come with
pinnings induced by \((\bT,\bB,\bFx_+)\). Using this \(k\)-pinning, we may
identify the group \(\Out(\bG)\)
    \nomenclature[\(Out_G \)]{\(\Out(G)\)}{the group of outer automorphisms of a split group \(G\)}
of outer automorphisms of \(G\) with a subgroup of \(\Aut_k(\bG)\).
This is a discrete group, possibly infinite. Given any \(k\)-scheme \(X\), we
can consider \'etale \(\Out(\bG)\)-torsors over \(X\). 

\subsection{}
Given such a \(\Out(\bG)\)-torsor \(\OGT_G\),
    \nomenclature[\(theta_G_var \)]{\(\OGT_G\)}{the \(\Out(\bG)\)-torsor inducing \(G\)}
we may obtain a quasi-split twisted form
\begin{align}
    G=\bG\x^{\Out(\bG)}\OGT_G
    \nomenclature[\(G \)]{\(G\)}{an \(\Out(\bG)\)-twisted form of \(\bG\)}
\end{align}
of \(\bG\) on \(X\). This is a reductive group scheme over \(X\) together with a
pinning \(\SPL=(T,B,x_+)\)
relative to \(X\). The torsor \(\OGT_G\) also induces a
\(\Out(\bG^\AD)=\Out(\bG^\SC)\)-torsor, still denoted by \(\OGT_G\). It induces
quasi-split forms \(G^\AD\) and \(G^\SC\) over \(X\). The Weyl group \(\bW\)
induces a Weyl group scheme \(W\)
over \(X\), which is an \'etale group scheme.

\subsection{}
If we fix a geometric point \(\infty\in X\), then any \'etale \(\Out(\bG)\)-torsor
can be given by a continuous homomorphism
\begin{align}
    \OGT^\bullet_G\colon\pi_1(X,\infty)\longto \Out(\bG).
    \nomenclature[\(theta_G_var_bullet \)]{\(\OGT^\bullet_G\)}{the pointed version of \(\OGT_G\)}
\end{align}
In this way the group \(G\) comes with a canonical geometric point \(\infty_G\)
    \nomenclature[\(infinity_G \)]{\(\infty_G\)}{the canonical geometric point induced by \(\OGT_G^\bullet\)}
over \(\infty\), and we use \((G,\infty_G)\) for this pointed twisted form.

\subsection{}
If \(G\) is \emph{any} reductive group  over an algebraically closed field with a
fixed maximal torus \(T\) and Weyl group \(W\), we define \(\WP_G(\lambda)\)
    \nomenclature[\(Xi_G_lambda \)]{\(\WP_G(\lambda)\)}{the set of (co)characters \(\mu\)
    contained in the convex hull of the Weyl orbit of a (co)character
    \(\lambda\), such that \(\lambda-\mu\) is in the root lattice}
    \nomenclature[\(Xi_lambda \)]{\(\WP(\lambda)\)}{same as \(\WP_G(\lambda)\) when \(G\) is clear from the context}
for any \(\lambda\in\CharG(T)\) to be the set of weights \(\mu\in\CharG(T)\)
in the convex hull of the orbit \(W\lambda\) such that \(\lambda-\mu\) is
contained in the root lattice of \(T\) in \(G\). If \(G\) is clear from the
context, we will simply use \(\WP(\lambda)\).

\section{Invariant Theory of the Group}%
\label{sec:invariant_theory_of_the_group}
In this section, we review the invariant theory of the adjoint
action of \(\bG\) on \(\bG^\SC\). Later in
\Cref{sec:invariant_theory_of_reductive_monoids} we will enhance most results in
this section to reductive monoids.

\subsection{}
The group \(\bG\) acts on \(\bG^\SC\) by adjoint action. 
Let \(\chi\colon \bG^\SC\to \bG^\SC\git \bG\)
    \nomenclature[\(chi \)]{\(\chi\)}{the GIT quotient map \(\bG^\SC\to\bC\) or \(G^\SC\to \FRC\)}
be the GIT-quotient map.

\begin{theorem}[\cite{St65}]
    The inclusion \(\bT^\SC\hookto \bG^\SC\) induces an isomorphism
    \begin{align}
        \bC\defeq \bT^\SC\git\bW\stackrel{\sim}{\longto} \bG^\SC\git
        \bG^\SC=\bG^\SC\git\bG.
        \nomenclature[\(C"bold \)]{\(\bC,\FRC\)}{the GIT quotient \(\bG^\SC\git \bG\) or \(G^\SC\git G\)}
    \end{align}
    In addition, both schemes are isomorphic to affine space \(\bbA^r\) whose
    coordinates are given by the traces \(\chi_i\)
    \nomenclature[\(chi_i \)]{\(\chi_i\)}{the trace of the fundamental representation \(\rho_i\)}
    of the fundamental representations \((\rho_i,V_{\Wt_i})\) of \(\bG^\SC\).
\end{theorem}

\begin{definition}
    A \(\bG\)-orbit \(\bG(\gamma)\) (\(\gamma\in \bG^\SC\)) is called
    \notion{regular}\index{orbit!regular} if its stabilizers have minimal dimension. It is called
    \notion{semisimple}\index{orbit!semisimple} if it contains an element in
    \(\bT^\SC\). It is called \notion{regular semisimple}\index{orbit!regular semisimple} if it is both regular
    and semisimple. The union of regular (resp.~semisimple, resp.~regular
    semisimple) orbits is denoted by \(\bG^{\SC}_\reg\)
    (resp.~\(\bG^\SC_\ssim\), resp.~\(\bG^\SC_\rss\)).
\end{definition}

\subsection{}
The regular locus \(\bG^\SC_\reg\)
    \nomenclature[\(.reg \)]{\((\cdot)^\reg, (\cdot)_\reg\)}{the (restriction to the) regular locus}
is open and the restriction of \(\chi\) to
\(\bG^\SC_\reg\) is smooth and surjective. The semisimple locus is 
dense but not open in \(\bG^\SC\). Their intersection \(\bG^\SC_\rss\),
    \nomenclature[\(.rs \)]{\((\cdot)^\rss, (\cdot)_\rss\)}{the (restriction to the) regular semisimple locus}
however,
remains open dense in \(\bG^\SC\), and its image in \(\bC\) is denoted by
\(\bC^\rss\)\index{locus!regular semisimple}.
Consider the \notion{discriminant function}\index{function!non-extended discriminant}
\index{discriminant!function, non-extended} on \(\bT\):
\begin{align}
    {\Disc}\defeq \prod_{\alpha\in\Roots}(1-e^{\alpha}).
    \nomenclature[\(D{}isc \)]{\(\Disc\)}{the non-extended discriminant function}
\end{align}
This function is \(\bW\)-invariant, hence descends to a function on \(\bC\),
still denoted by \({\Disc}\). It defines an effective principal divisor \(\bD\)
    \nomenclature[\(D"bold \)]{\(\bD,\FRD\)}{the non-extended discriminant divisor on \(\bC\) or \(\FRC\)}
on \(\bC\), called the \notion{discriminant divisor}\index{divisor!non-extended
discriminant}\index{discriminant!divisor, non-extended}. 
In fact, \(\bD\) is a reduced divisor, so it makes sense to talk about its
singular locus \(\bD^\Sing\).
    \nomenclature[\(.sing \)]{\((\cdot)^\Sing, (\cdot)_\Sing\)}{the singular locus}

\begin{proposition}[\cite{St65}]
    The regular semisimple locus \(\bC^\rss\)
    is equal to the complement of the
    divisor \(\bD\) and thus is an open subset of \(\bC\). Moreover, it is
    exactly the locus over which the fiber of \(\chi\) consists of a single
    \(\bG\)-orbit.
\end{proposition}

\subsection{}%
\label{sub:steinberg_quasi_section}

The quotient map \(\chi\) admits many sections, just like in the Lie algebra
case.  The difference is that the section we want to use lacks an explicit
formula, unlike for example the Kostant section for Lie algebras. Instead, what
Steinberg explicitly constructed is a quasi-section. A morphism \(f\) in a
category is called a
\notion{quasi-section}\index{quasi-section} to morphism \(g\) if \(g\circ f\) is an \emph{automorphism} (not
necessarily identity).

\begin{definition}\label[definition]{def:Coxeter_element}
    Fix our choice of simple roots \(\SimRts\) earlier.
    A \notion{Coxeter datum}\index{Coxeter!datum}  is a pair \((\xi, \dot{S})\), where
    \begin{enumerate}
        \item \(\xi\colon \Set{1,\ldots,r}\to \SimRts\) is a bijection (i.e. a
            total ordering on the set of simple roots),
        \item \(\dot{S}\) is a set of representatives \(\dot{s}_{\alpha}\) in
            \(\Norm_{\bG}(\bT)\)
            \nomenclature[\(N_G_T \)]{\(\Norm_{G}(T)\)}{the normalizer of a subgroup \(T\) in a group \(G\)}
            of simple reflections determined by \(\SimRts\).
    \end{enumerate}
    A \notion{Coxeter element}\index{Coxeter!element} (after fixing a set of simple reflections)
    of \(\bW\) is one that can be written as
    \(w=w_{\xi}=s_{\xi(1)}\cdots s_{\xi(r)}\) for some \(\xi\).  Denote by
    \(\Cox(\bW,\SimRts)\)
    \nomenclature[\(C{}ox_W_Delta \)]{\(\Cox(W,\SimRts)\)}{the set of Coxeter
    elements of a Weyl group \(W\) with respect to a set of simple reflections \(\SimRts\)}
    the set of all Coxeter elements of \(\bW\).
\end{definition}

Fix a Coxeter datum \((\xi,\dot{S})\).  Let \(\beta_j=\xi(j)\) be the simple
root corresponding to \(s_{\xi(j)}\). Recall we have one-parameter groups
\(\bU_\alpha\) in the pinning \(\bSPL\), such that
\(\Ad_t(\bU_\alpha(x))=\bU_\alpha(\alpha(t)x)\) for all \(t\in \bT\) and \(x\in
\Ga\). Let
\begin{align}
    \epsilon^{(\xi,\dot{S})}\colon \bC\cong\bbA^r &\longto \bG^\SC\\
    (x_1,\ldots,x_r)&\longmapsto
    \prod_{j=1}^r\bU_{\beta_{j}}(x_{j})\dot{s}_{\xi(j)},
    \nomenclature[\(epsilon_xi_S \)]{\(\epsilon^{(\xi,\dot{S})}\)}{the non-extended Steinberg quasi-section induced by a Coxeter datum \((\xi,\dot{S})\)}
\end{align}
where the product on the right-hand side is considered taken in the specified
order. This is the Steinberg quasi-section associated with Coxeter datum
\((\xi,\dot{S})\).  We also call the image \(\img{\epsilon^{(\xi,\dot{S})}}\) a
\notion{Steinberg cross-section}\index{Steinberg cross-section}. We summarize the results in the following
theorem.

\begin{theorem}
    [\cite{St65}]
    \label[theorem]{thm:Steinberg_section_group}
    For each pair \((\xi,\dot{S})\), the map \(\epsilon^{(\xi,\dot{S})}\) is a
    quasi-section of \(\chi\).  Moreover, \(\img{\epsilon^{(\xi,\dot{S})}}\) is
    contained in the regular locus, and has transversal intersection with each
    regular orbit.
\end{theorem}

Once \(\bSPL\) is fixed, the construction of Steinberg quasi-section relies on
two choices: a total ordering of the simple roots, and representatives of simple
reflections. The influence of the choices is summarized below.

\begin{proposition}
    [\cite{St65}*{Lemma~7.5 and Proposition~7.8}]
    \label[proposition]{prop:st_quasisection_influence_of_choices}
    For any Coxeter data \((\xi,\dot{S})\) and \((\xi',\dot{S}')\):
    \begin{enumerate}
        \item If \(\xi=\xi'\), then there exists \(t,t'\in \bT\) such that
            \(\img{\epsilon^{(\xi,\dot{s}')}}=t'\img{\epsilon^{(\xi,\dot{S})}}=t\img{\epsilon^{(\xi,\dot{S})}}t^{-1}\).
        \item If \(\dot{S}=\dot{S}'\), then for any \(x,x'\in\bbA^r\) such that
            \(x_{\xi(j)}=x'_{\xi'(j)}\) for \(1\le j\le r\),
            \(\epsilon^{(\xi,\dot{S})}(x)\) and
            \(\epsilon^{(\xi',\dot{S})}(x')\) are \(\bG\)-conjugate. In fact, such
            conjugation can be made functorially for any \(k\)-algebra \(R\).
            In other words, the transporter from
            \(\epsilon^{(\xi,\dot{S})}(\bC)\)
            to \(\epsilon^{(\xi',\dot{S})}(\bC)\) is a trivial torsor under the
            centralizer group scheme over \(\epsilon^{(\xi,\dot{S})}(\bC)\).
    \end{enumerate}
\end{proposition}

\subsection{}%
\label{sub:regular_centralizer}

The universal centralizer group scheme \(\bI\to \bG^\SC\)
    \nomenclature[\(I \)]{\(\bI,I\)}{the universal stabilizer group scheme of \(\bG\)
    (resp.~\(G\)) acting on \(\bG\) or \(\bG^\SC\) (resp.~\(G\) or \(G^\SC\)) depending on the context}
restricts to a smooth group
scheme over the regular locus \(\bI^\reg\to \bG^\SC_\reg\).
Since generically the
centralizer is a torus, \(\bI^\reg\to \bG^\SC_\reg\) is a commutative group scheme.
Therefore, one can utilize the descent argument in \cite{Ng10}*{Lemme~2.1.1} to
obtain the following result.

\begin{proposition}
    There is a unique smooth commutative group scheme \(\bJ\to \bC\)
    \nomenclature[\(J \)]{\(\bJ,\FRJ\)}{the regular centralizer group scheme over \(\bC\) or \(\FRC\)}
    with a \(\bG\)-equivariant isomorphism
    \begin{align}
        \chi^*\bJ|_{\bG^\SC_{\reg}}\stackrel{\sim}{\longto} \bI^{\reg},
    \end{align}
    which can be extended to a homomorphism \(\chi^*\bJ\to \bI\).
\end{proposition}

\begin{definition}
    The group scheme \(\bJ\to\bC\) is called the \notion{regular
    centralizer}\index{regular centralizer}.
\end{definition}

There is another description of \(\bJ\) using the cameral cover \(\pi\colon
\bT^\SC\to\bC\), similar to the Lie algebra case as in \cite{DoGa02} and
\cite{Ng10}.  Consider the trivial family of torus \(\pr_2\colon \bT\x \bT^\SC\to
\bT^\SC\)
(\(\pr_2\) means the second projection), and
we have Weil restriction
\begin{align}
    \Pi_{\bG}\defeq \pi_*(\bT\x \bT^\SC),
\end{align}
on which \(\bW\) acts diagonally. Since \(\pi\) is finite flat, \(\Pi_{\bG}\) is
representable.  Moreover, since \(\bT\x \bT^\SC\) is smooth, so is \(\Pi_{G}\). Over
\(\bC^\rss\), \(\pi\) is Galois \'etale with Galois group \(\bW\), hence
\(\Pi_{\bG}^\rss\) is a torus. Since \(\Char(k)\) does not divide the order
of \(\bW\), we have a smooth group scheme \(\bJ^1\) over \(\bC\)
\begin{align}
    \bJ^1\defeq \Pi_{\bG}^{\bW}.
    \nomenclature[\(J^1 \)]{\(\bJ^1,\FRJ^1\)}{the slightly larger Galois-theoretic analogue of \(\bJ\) or \(\FRJ\)}
\end{align}
Let \(\bJ^0\subset\bJ^1\)
    \nomenclature[\(J^0"bold \)]{\(\bJ^0,\FRJ^0\)}{the fiberwise neutral component of \(\bJ\) or \(\FRJ\)}
be the open subgroup scheme of fiberwise neutral component.

Similar to \cite{Ng10}*{D\'efinition~2.4.5}, we consider this subfunctor
\(\bJ'\) of \(\bJ^1\):
for a \(\bC\)-scheme \(S\), \(\bJ'(S)\) consists of points 
\begin{align}
    f\colon S\x_{\bC}\bT^\SC\to \bT
\end{align}
such that for every geometric point \(x\in S\x_{\bC}\bT^\SC\), if
\(s_\Rt(x)=x\) for a root \(\Rt\), then \(\Rt(f(x))\neq -1\). Notice that on
\(\bT^\SC\), the condition \(s_\Rt(x)=x\) is the same as \(x\in\ker{\Rt}\).

\begin{lemma}
    The subfunctor \(\bJ'\) is representable by an open subgroup scheme of
    \(\bJ^1\) containing \(\bJ^0\).
\end{lemma}
\begin{proof}
    The proof is entirely parallel to \cite{Ng10}*{Lemme~2.4.6}.
    Indeed, it suffices
    to prove this claim after finite flat base change to \(\bT^\SC\).

    On \(\bT^\SC\), the discriminant divisor is the
    union (with multiplicities)
    of subgroups \(\ker{\Rt}\subset\bT^\SC\) for all roots \(\Rt\in \Roots\). By
    adjunction,
    we have a map \(\bJ^1\x_{\bC}\bT^\SC\to\bT\x \bT^\SC\), whose restriction
    to \(\ker{\Rt}\) factors through fixed-point subgroup
    \(\bT^{s_\alpha}\x\ker{\Rt}\). Therefore, we have a map
    \begin{align}
        \bJ^1\x_{\bC}\ker{\Rt}\stackrel{\Rt}{\longto}\Set{\pm 1}.
    \end{align}
    The inverse image of \(-1\) is an open and closed subset of
    \(\bJ^1\x_{\bC}\ker{\Rt}\), hence a Cartier divisor on
    \(\bJ^1\x_{\bC}\bT^\SC\). The subfunctor \(\bJ'\x_{\bC}\bT^\SC\) is the
    complement of these Cartier divisors, hence an open subgroup scheme. In
    addition, \(\bJ'\) contains \(\bJ^0\).
\end{proof}

\begin{proposition}
    There exists a canonical open embedding of group schemes \(\bJ\to \bJ^1\)
    that identifies \(\bJ\) with subgroup scheme \(\bJ'\).
\end{proposition}
\begin{proof}
    The claim about open embedding is proved in \cite{Ch22}*{Proposition~2.3.2}.
    In fact, the argument in \textit{loc. cit.}
    shows that if \(\bG=\bG^\SC\), then \(\bJ\to\bJ^1\) is an
    isomorphism. Moreover, since on \(\bT^\SC\), \(\ker{\Rt}\) and
    \(\bT^{\SC,s_\Rt}\) coincide for any root \(\Rt\), we also see that
    \(\bJ'=\bJ^1\).

    In general, we have a canonical map \(\bJ^\SC\to\bJ\) and
    \(\bJ\) is generated by the image of \(\bJ^\SC\)
    and \(\bZ_{\bG}\). On the other hand, the map \(\bT^\SC\to\bT\) induces
    canonical map \(\bJ^{\SC,1}\to\bJ^1\) compatible with the map
    \(\bJ^\SC\to\bJ\). By definition of \(\bJ'\) we also have a third map 
    \(\bJ^{\SC,\prime}\to\bJ'\). This means that the image of \(\bJ^{\SC,1}\) is
    contained in \(\bJ'\). Since on \(\bZ_{\bG}\) all roots are trivial, we know
    that the image of \(\bJ\) is contained in \(\bJ'\).

    It remains to show that \(\bJ\to\bJ'\) is an isomorphism. Then we can
    repeat the ``codimension \(2\)'' argument 
    of \cite{Ch22}*{Proposition~2.3.2} (see also
    \cite{Ng10}*{Proposition~2.4.7}) to reduce the problem to the case where
    \(\bG=\SL_2\), \(\GL_2\), or \(\PGL_2\). The direct calculation is omitted.
\end{proof}

\begin{remark}
    \label[remark]{rmk:connectedness_of_FRJ_vs_center_of_G}
    The proof we present here is subtly different from
    \cite{Ng10}*{Proposition~2.4.7}, where there is a surjection
    \(\pi_0(\bZ_{\bG})\to\pi_0(\bJ_{\La{g}})=\bJ_{\La{g}}/\bJ_{\La{g}}^0\) (here
    \(\bJ_{\La{g}}\) stands for the regular centralizer for the Lie algebra of
    \(\bG\)). It is not so in group case. 
    For example, suppose \(\bG\) is simple of type \(\TypeG_2\). Then
    \(\bG=\bG^\AD=\bG^\SC\). Let \(x\in\bT\) be such that
    \(\Cent_{\bG}(x)\cong\SL_3\). Let \(u\in U_3\defeq\bU\cap\SL_3\) be a
    regular unipotent element in \(\SL_3\). Assuming \(\Char(k)\) is large
    enough, then \(\Cent_{\bG}(xu)\) contains \(\Cent_{\SL_3}(u)=Z_{\SL_3}\x
    \Cent_{U_3}(u)\) with finite index. So \(\Cent_{\bG}(xu)\) is disconnected
    but \(\bZ_{\bG}\) is trivial.
\end{remark}

\subsection{}

Let \(G\) be a quasi-split form of \(\bG\) over \(k\)-scheme \(X\), induced by
an \(\Out(\bG)\)-torsor \(\OGT_G\). We can twist almost every construction in
the invariant theory of \(\bG\) by \(\OGT_G\) and obtain a twisted form over
\(X\). First we still have adjoint action of \(G\)
on \(G^\SC\), with invariant quotient \(\chi\colon G^\SC\to \FRC=G^\SC\git G\).
We also have the natural isomorphism
\begin{align}
    T\git W\stackrel{\sim}{\longto} G^\SC\git G.
\end{align}
The discriminant is invariant under \(\Out(\bG)\), so we have a divisor \(\FRD\subset\FRC\)
relative to \(X\). The open loci of regular and regular
semisimple orbits are stable under \(\Out(\bG)\), and so we still have the
twisted form of regular centralizer \(\FRJ\to\FRC\)
over \(X\), as well as its Galois construction \(\FRJ^1\).

\subsection{}
An important difference here compared to Lie algebra case is that Steinberg
quasi-section does not necessarily exist, since it may not be stable under
\(\Out(\bG)\). Nevertheless, one may obtain a weaker result as follows.
The Steinberg quasi-sections depends on a choice of representatives of simple
reflections in \(\bG^\SC\). Using the pinning \(\bSPL\), one can
make such a choice that is stable under \(\Out(\bG)\). Indeed, the root
vectors in \(\bFx_+\) can each be extended to a unique \(\La{sl}_2\)-triple, hence
an opposite pinning with root vectors denoted by \(\bFx_-\). For each simple
root \(\alpha_i\in\SimRts\) (under \(\bB\)), we let
\begin{align}
    \dot{s}_i=\bU_{\alpha_i}(1)\bU_{-\alpha_i}(-1)\bU_{\alpha_i}(1),
    \nomenclature[\(s_dot_i \)]{\(\dot{s}_i\)}{the image of the Tits section of the \(i\)-th simple root}
\end{align}
where \(\bU_{\alpha_i}\) (resp.~\(\bU_{-\alpha_i}\)) are the one-parameter unipotent
subgroups determined by \(\bFx_+\) (resp.~\(\bFx_-\)).

The Steinberg quasi-section also depends on a choice of ordering on \(\SimRts\).
When \(\bG^\SC\) does not contain any simple factor of type \(\TypeA_{2m}\), any
two simple roots conjugate under \(\Out(\bG)\) are not linked in the Dynkin
diagram. Therefore, by grouping together \(\Out(\bG)\) orbit when making the
ordering, we can make such Steinberg quasi-section equivariant under \(\Out(\bG)\).
Thus, a section to \(\chi\) always exists as long as \(\bG^\SC\)
does not have any simple factor of type \(\TypeA_{2m}\).

In fact, we have a slightly stronger result: 
let \(\bG'\) be the direct factor of \(\bG^\SC\) consisting of all its simple
factors of types \(\TypeA_{2m}\) (for various \(m\)), which is preserved by
\(\Out(\bG)\). Then a Steinberg quasi-section exists for \(G^\SC\) as long as
the twist \(G'\) of \(\bG'\) is split.

In the remaining cases (i.e., \(G'\) is a non-trivial outer twist),
we cannot hope to construct a section whose image lies in the regular locus.
Rather, we have a weaker result by Steinberg, which is still very helpful later.
\begin{theorem}
    [\cite{St65}*{Theorem~9.8}]
    Let \(G\) be a quasi-split semisimple and simply-connected group over
    a perfect field \(K\), then \(G(K)\to \FRC(K)\) is surjective.
\end{theorem}
The statement in \textit{loc.~cit.} requires \(K\) to be perfect, but the proof
in fact works for arbitrary field provided that \(\Char(K)\) is not too small
for \(G\) (e.g., larger than twice the Coxeter number of
\(G\)). If \(K\) is perfect, then the result refines
to that \(G^{\ssim}(K)\to\FRC(K)\) is surjective.
We will postpone the details until \Cref{thm:A_2m_rationality_fix} 
since we need to extend the above result to reductive monoids.

\section{Review of Very Flat Reductive Monoids}%
\label{sec:review_of_very_flat_reductive_monoids}

The Lie algebra of \(\bG\) carries a natural \(\Gm\)-action from its vector
space structure, which is useful in global constructions. The group \(\bG^\SC\),
however, has no such symmetry built in. Therefore, we must embed \(\bG^\SC\)
into a \(\bG\)-space where a similar ``scaling'' is possible. The quintessential
example of this is the embedding of \(\SL_n\) into \(\Mat_n\). In general, we
use the theory of very flat reductive monoids.

\subsection{}%

An algebraic semigroup over \(k\) is just a \(k\)-scheme of
finite type \(M\) together with a multiplication morphism \(m\colon M\x M\to
M\), such that the usual commutative diagram of associativity holds. If there
exists a multiplicative identity \(e\colon \Spec{k}\to M\), then \(M\) is an
algebraic monoid over \(k\). We will only consider monoids that are affine,
integral and normal as schemes. If the subgroup \(M^\x\)
    \nomenclature[\(.times \)]{\((\cdot)^\x\)}{the group of invertible elements in a monoid}
of invertible elements
of \(M\) is a reductive group, then we call \(M\) a \notion{reductive
monoid}\index{monoid!reductive}.
In this case, we denote the derived subgroup of \(M^\x\) by \(M^\Der\),
    \nomenclature[\(M_der \)]{\(M^\Der\)}{the derived subgroup of \(M^\x\), where \(M\) is a monoid}
and call it the \notion{derived subgroup}\index{subgroup!derived} of \(M\).

\begin{example}
    A (normal) toric variety \(M\) is a reductive monoid, whose unit group
    \(M^\x\) is the torus acting on it, and whose derived subgroup is the
    trivial group. The variety \(\Mat_n\) with matrix multiplication is also a
    reductive monoid, whose unit group is \(\GL_n\) and derived subgroup is
    \(\SL_n\).
\end{example}

The category of normal reductive monoids is classified by Renner (see e.g.,
\cite{Re05})
with the help of their unit groups. It is well known that over an algebraically
closed field \(K\), the category of normal affine toric
varieties \(A\) of a fixed torus \(T\) is equivalent to that of strictly
convex and saturated cones \(\sC\subset\CoCharG(T)\).
    \nomenclature[\(C"scr \)]{\(\sC\)}{the cocharacter cone associated with an affine toric variety}
If \(T\subset G\) is a maximal torus in a reductive group over \(K\) with Weyl
group \(W\), and suppose
\(\sC\) is stable under \(W\)-action, then the \(W\)-action extends over \(A\).
Renner's classification theorem states that reductive monoids \(M\) with unit
group \(G\) is classified by such cones \(\sC\). 
More precisely, we have the following result.

\begin{theorem}[\cite{Re05}*{Theorems~5.2 and 5.4}]
    \label[theorem]{thm:monoid_classification_over_k_bar}
    Let \(G\) be a reductive group over an algebraically closed field with
    maximal torus \(T\) and Weyl group \(W=W(G,T)\).
    \begin{enumerate}
        \item Let \(T\subset A\) be a normal affine toric variety such that the
            \(W\)-action on \(T\) extends over \(A\). Then there exists a normal
            monoid \(M\) with unit group \(G\) such that \(\bar{T}\)
            \nomenclature[\(T_bar \)]{\(\bar{T}\)}{the closure of a subgroup \(T\) in a monoid \(M\)}
            (the Zariski closure of \(T\) in \(M\)) is isomorphic to \(A\).
        \item Let \(M\) be any normal reductive monoid with unit group \(G\), then
            the submonoid \(\bar{T}\) is normal.
        \item If \(M_1\) and \(M_2\) are such that we have a commutative
            diagram of monoids
            \begin{equation}
                \begin{tikzcd}
                    \bar{T}_1 \ar[d] & T_1\ar[l] \ar[r]\ar[d] & G_1\ar[d]\\
                    \bar{T}_2  & T_2\ar[l]\ar[r] & G_2
                \end{tikzcd}
            \end{equation}
            Then this diagram extends to a unique homomorphism \(M_1\to M_2\).
            Moreover, if the vertical arrows are isomorphisms, then \(M_1\cong
            M_2\).
    \end{enumerate}
\end{theorem}

\begin{remark}
    \label[remark]{rmk:monoid_classification_over_base_field}
    The result is proved over an algebraically closed field, but it is not hard
    to see that if \(G\) is split over a perfect field \(k\), 
    \(T\) is a split maximal torus,
    and all the maps in the commutative diagram 
    are defined over \(k\), then \(M_1\to M_2\) is also defined
    over \(k\) by looking at \(\Gal(\bar{k}/k)\)-action. The existence of
    \(M\) given \(G\) and \(A\) is proved using \(k\)-rational representations
    which are well-defined since \(G\) is split. See \textit{loc. cit.} for
    details.
\end{remark}

\subsection{}%
We first give an explicit description of a very important monoid \(\bM\) with
\(\bM^\Der\simeq \bG^\SC\). Consider the group
\begin{align}
    \bG_+\defeq (\bT^\SC\x \bG^\SC)/\bZ^\SC,
    \nomenclature[\(G_+ \)]{\(\bG_+,G_+\)}{the extended group \((\bT^\SC\x \bG^\SC)/\bZ^\SC\) or \((T^\SC\x G^\SC)/Z^\SC\)}
\end{align}
where the center \(\bZ^\SC\) of \(\bG^\SC\) acts on \(\bT^\SC\x \bG^\SC\)
anti-diagonally. There is a
maximal torus \(\bT_+=(\bT^\SC\x \bT^\SC)/\bZ^\SC\)
    \nomenclature[\(T_+ \)]{\(\bT_+,T_+\)}{the maximal torus of \(\bG_+\) or \(G_+\)}
in \(\bG_+\), and we denote
\(\bZ_+=\bZ_{\bG_+}\cong\bT^\SC\).
    \nomenclature[\(Z_+ \)]{\(\bZ_+,Z_+\)}{the center of \(\bG_+\) or \(G_+\)}
The character lattice of \(\bT_+\) is identified as
\begin{align}
    \CharG(\bT_+)=\Set*{(\lambda,\mu)\in \CharG(\bT^\SC)\x\CharG(\bT^\SC)\given
                \lambda-\mu\in \CharG(\bT^\AD)},
\end{align}
and its cocharacter lattice is identified as
\begin{align}
    \CoCharG(\bT_+)=\Set*{(\dual{\lambda},\dual{\mu})\in
    \CoCharG(\bT^\AD)\x\CoCharG(\bT^\AD)\given \dual{\lambda}+\dual{\mu}\in
        \CoCharG(\bT^\SC)}.
\end{align}
So we have characters \((\alpha,0)\) (\(\alpha\in\SimRts\)) of \(\bT_+\) which are
also one-dimensional representations of \(\bG_+\). On the other hand, any \(k\)-rational
representation \(\rho\) of highest-weight \(\lambda\) of \(\bG^\SC\) extends to
\(\bG_+\) by the formula
\begin{align}
    \rho_{+}(z,g)=\lambda(z)\rho(g),
\end{align}
which is easily seen well-defined. Consider the representation of \(\bG_+\)
\begin{align}
    \rho_\star\defeq \bigoplus_{i=1}^r(\alpha_i,0)\oplus
    \bigoplus_{i=1}^r \rho_{i+}.
    \nomenclature[\(rho_star \)]{\(\rho_\star\)}{the representation of \(\bG_+\) defining the universal monoid}
\end{align}
We take the normalization of the closure of \(\bG_+\) in \(\rho_\star\), 
and denote it by \(\Env(\bG^\SC)\).
    \nomenclature[\(Env_G \)]{\(\Env(G)\)}{the universal monoid of a semisimple group \(G\)}
This is the \notion{universal monoid}\index{universal!monoid}\index{monoid!universal}
associated with \(\bG^\SC\), constructed by Vinberg (\cite{Vi95}) in characteristic \(0\) and
Rittatore (\cite{Ri01}) in all characteristics.

\subsection{}%
Given any reductive monoid \(\bM\)
    \nomenclature[\(M"bold \)]{\(\bM\)}{a split reductive monoid, usually assumed to be very flat}
with derived subgroup \(\bG^\SC\), we have an action
of \(\bG^\SC\) by multiplication on the left, and another one by
inverse-multiplication on the right. These two actions commute, hence
combine to a \(\bG^\SC\x \bG^\SC\)-action on \(\bM\). The GIT quotient
\begin{align}
    \alpha_{\bM}\colon \bM\longto \bA_{\bM}\defeq \bM\git (\bG^\SC\x \bG^\SC)
    \nomenclature[\(A_M \)]{\(\bA_\bM,\FRA_\FRM\)}{the abelianization of monoid \(\bM\) or \(\FRM\)}
    \nomenclature[\(alpha_M \)]{\(\alpha_\bM,\alpha_\FRM\)}{the abelianization map of monoid \(\bM\) or \(\FRM\)}
\end{align}
is called the \notion{abelianization}\index{map!abelianization}\index{abelianization map} of \(\bM\), in
the sense that it is the largest quotient monoid of \(\bM\) that is commutative.

\begin{theorem}
    [\cites{Vi95,Ri01}]
    \label[theorem]{thm:abelianization_first_properties}
    Let \(\bZ_{\bM}\) be the center of \(\bM^\x\), \(\bZ_{\bM,0}\) its neutral
    component, and \(\bZ_0=\bZ_{\bM,0}\cap
    \bG^\SC\), then we have that
    \begin{enumerate}
        \item \(\bA_{\bM}\) is a normal commutative monoid with unit group
            \(\alpha_{\bM}(\bZ_{\bM})\simeq \bZ_{\bM,0}/\bZ_0\),
        \item \(\alpha_{\bM}^{-1}(1)=\bG^\SC\) and
            \(\alpha_{\bM}^{-1}(\alpha_{\bM}(\bZ_{\bM,0}))=\bG_+\),
        \item \(\alpha_{\bM}(\bar{\bZ_{\bM,0}})=\bA_{\bM}\), and \(\bA_{\bM}\simeq
            \bar{\bZ_{\bM,0}}\git\bZ_0\), where \(\bar{\bZ_{\bM,0}}\) denotes the Zariski closure of
            \(\bZ_{\bM,0}\) in \(\bM\).
    \end{enumerate}
\end{theorem}

\subsection{}
Given a homomorphism of two reductive monoids \(\phi\colon \bM'\to \bM\), by the
universal property of GIT quotient it induces a homomorphism \(\phi_{\bA}\colon
\bA_{\bM'}\to\bA_{\bM}\) that fits into the commutative square
\begin{equation}\label{eqn:abelianization_square}
    \begin{tikzcd}
        \bM'\ar[r, "\phi"]\ar[d, "\alpha_{\bM'}",swap] & \bM\ar[d, "\alpha_{\bM}"]\\
        \bA_{\bM'}\ar[r, "\phi_{\bA}"] & \bA_{\bM}
    \end{tikzcd}.
\end{equation}

\begin{definition}
    The homomorphism \(\phi\colon \bM'\to \bM\) is called
    \notion{excellent}\index{excellent homomorphism} if \eqref{eqn:abelianization_square} is
    Cartesian.
\end{definition}

\begin{definition}
    \label[definition]{def:very_flat_monoid}
    A reductive monoid \(\bM\) is called \notion{very flat}\index{monoid!very
    flat} if \(\alpha_{\bM}\colon \bM\to \bA_{\bM}\) is flat with integral fibers. Let
    \(\FM(\bG^\SC)\)
    \nomenclature[\(F{}M"cal_G \)]{\(\FM(G)\)}{the category of very flat monoids associated with a semisimple group \(G\)}
    be the category in which:
    \begin{enumerate}
        \item objects are very flat reductive monoids \(\bM\) with \(\bM^{\Der}\cong
            \bG^\SC\);
        \item morphisms from \(\bM'\) to \(\bM\) are excellent
            homomorphisms \(\phi\colon \bM'\to \bM\).
    \end{enumerate}
    Let \(\FM_0(\bG^\SC)\)
    \nomenclature[\(F{}M"cal_0_G \)]{\(\FM_0(G)\)}{the full subcategory of \(\FM(G)\) consisting of monoids with \(0\)}
    be the full subcategory of \(\FM(\bG^\SC)\) whose objects are
    very flat monoids with element \(0\) (one such that \(0x=x0=0\) for all \(x\in \bM\)).
\end{definition}

\begin{theorem}
    [\cites{Vi95,Ri01}]
    \label[theorem]{thm:env_semigroup}
    If \(\bM\in \FM(\bG^\SC)\), then
    \begin{align}
        \Hom_{\FM(\bG^\SC)}(\bM,\Env(\bG^\SC))\neq\emptyset, 
    \end{align}
    and is a singleton if \(\bM\in \FM_0(\bG^\SC)\). In other words,
    \(\Env(\bG^\SC)\) is a versal (resp.~universal) object in \(\FM(\bG^\SC)\)
    (resp.~\(\FM_0(\bG^\SC)\)).
\end{theorem}

\begin{remark}
    \Cref{thm:env_semigroup} implies that a homomorphism
    \(\bM\to\Env(\bG^\SC)\) will induce a homomorphism \(\phi_{\bZ_{\bM}}\colon
    \bZ_{\bM}\to \bT^\SC\)
    \nomenclature[\(phi_Z_M \)]{\(\phi_{\bZ_\bM},\phi_{Z_\FRM}\)}{the map
    \(\bZ_{\bM}\to \bT^\SC\) (resp.~\(Z_{\FRM}\to T^\SC\)) induced by a map
    \(\bM\to\Env(\bG^\SC)\) (resp.~\(\FRM\to\Env(G^\SC)\)) for very flat monoid \(\bM\) or \(\FRM\)}
    whose restriction on \(\bZ^\SC\) is the identity. In case
    that \(\bM=\Env(\bG^\SC)\), \(\phi_{\bZ_{\bM}}=\Id_{\bT^\SC}\).
\end{remark}

\subsection{}%
\label{sub:different_monoids}
The abelianization of the universal monoid \(\Env(\bG^\SC)\) is an affine space of
dimension \(r\), whose coordinate functions can be exactly given by 
\(e^{(\alpha,0)}\) where \(\alpha\in\SimRts\) are simple roots. In the language
of toric varieties, \(\bA_{\Env(\bG^\SC)}\) is the \(\bT^\AD\)-toric variety associated
with cone generated by fundamental coweights \(\CoWt_i\in\CoCharG(\bT^\AD)\). 

If \(\bA\) is any affine normal toric variety, and \(\phi\colon\bA\to
\bA_{\Env(\bG^\SC)}\) is a homomorphism of algebraic monoids, then the pullback
of \(\Env(\bG^\SC)\) through \(\phi\) gives an object \(\bM\in\FM(\bG^\SC)\)
with abelianization \(\bA_{\bM}=\bA\), and by \Cref{thm:env_semigroup},
every \(\bM\in\FM(\bG^\SC)\) arises in this way.

An important class of monoids
consists of those corresponding to the map \(\bbA^1\to \bA_{\Env(\bG^\SC)}\)
given by a \emph{dominant} cocharacter \(\lambda\in\CoCharG(\bT^\AD)\).
We will denote such monoid by \(\bM(\lambda)\). More generally, we may
also consider the monoid formed using a tuple of dominant cocharacters
\(\ul{\lambda}=(\lambda_1,\ldots,\lambda_m)\) (allowing
repetitions), or equivalently, a multiplicative map \(\bbA^m\to
\bA_{\Env(\bG^\SC)}\). We shall denote this monoid by
\(\bM(\ul{\lambda})\).
    \nomenclature[\(M_lambda \)]{\(\bM(\ul{\lambda}),\FRM(\ul{\lambda})\)}{the
    standard very flat monoid associated with a tuple of cocharacters
\(\ul{\lambda}\) (that are stable under monodromy in non-split case)}

Moreover, if \(\ul{\lambda}\) is a
tuple of dominant cocharacters in \(\CoCharG(\bT)\), it induces a tuple
\(\ul{\lambda}_\AD\) of dominant cocharacters in \(\CoCharG(\bT^\AD)\).
The induced map \(\bbA^m\to\bA_{\Env(G^\SC)}\) also produces a monoid, still
denoted by \(\bM(\ul{\lambda})\). When \(\ul{\lambda}\) is a
singleton, \(\bM(\ul{\lambda})\) is often called an
\notion{\(L\)-monoid}\index{monoid!\(L\)-}
in the literature.

\begin{example}
    When \(\bG^\SC=\SL_n\) and \(\bM=\Mat_n\), the abelianization map is the
    determinant, and the excellent map
    \(\bM\to\Env(\SL_n)\) is the pullback of
    \begin{align}
        \bA_{\bM}\cong\bbA^1\longto \bA_{\Env(\SL_n)}\cong\bbA^{n-1}
    \end{align}
    corresponding to cocharacter \(\CoWt_{n-1}\), and the map \(\phi_{\bZ_{\bM}}\)
    is the cocharacter
    \(\CoRt_1+2\CoRt_2+\cdots+(n-1)\CoRt_{n-1}=n\CoWt_{n-1}\). In other words,
    \(\Mat_n=\bM(\CoWt_{n-1})\).
\end{example}

\subsection{}%
It is sometimes convenient to define a \notion{numerical boundary
divisor}\index{divisor!numerical boundary} \(\bE_{\bM}\)
    \nomenclature[\(E"bold_M \)]{\(\bE_\bM,\FRE_\FRM\)}{the numerical boundary divisor of monoid \(\bM\) or \(\FRM\)}
on \(\bA_{\bM}\) as the complement of \(\bA_{\bM}^\x\), with
reduced scheme structure. It is a Weil divisor and
when \(\bA_{\bM}\) is factorial, it is an effective Cartier divisor. In this
book we will mostly be interested in those \(\bM\) such that \(\bA_{\bM}\) is
isomorphic to an affine space \(\bbA^m\), in which case \(\bE_{\bM}\) is cut out
by the product of the \(m\) coordinate functions, hence principal.  Note that
\(\bE_{\bM}\) is in general \emph{not} the pullback of \(\bE_{\Env(\bG^\SC)}\),
but it always contains the latter after passing to the underlying topological
spaces. 

\subsection{}%
\label{sub:section_of_abelianization}
The abelianization map admits a section as follows: for \(\bM=\Env(\bG^\SC)\), let
\(\bT^\AD\to \bT_+\) be the map \(t\mapsto (t,t^{-1})\). It is well-defined and
extends to a map \(\delta_{\bM}\colon \bA_{\bM}\to \bM\),
    \nomenclature[\(delta_M \)]{\(\delta_\bM,\delta_\FRM\)}{the anti-diagonal section of the abelianization map of a very flat monoid \(\bM\) or \(\FRM\)}
which is in fact a section of
\(\alpha_\bM\). For a general monoid in \(\FM(\bG^\SC)\), 
the formula is \(z\mapsto (z,\phi_{\bZ_{\bM}}(z)^{-1})\).

\subsection{}%
\label{sub:cones_of_maximal_toric_variety}
Let \(\bT_{\bM}\)
    \nomenclature[\(T_M \)]{\(\bT_\bM,T_\FRM\)}{the maximal torus of \(\bM^\x\) or \(\FRM^\x\)}
be the maximal torus of \(\bM^\x\) containing \(\bT^\SC\), and \(\bar{\bT}_{\bM}\)
    \nomenclature[\(T_M_bar \)]{\(\bar{\bT}_\bM,\FRT_\FRM\)}{the maximal toric variety of \(\bM\) or \(\FRM\)}
its closure in \(\bM\). It is a normal affine toric variety under \(\bT_{\bM}\),
and when \(\bM=\Env(\bG^\SC)\), its character cone \(\sC^*_{\bar{\bT}_{\bM}}\)
and cocharacter cone \(\sC_{\bar{\bT}_{\bM}}\) have the following description:
\begin{align}
    \sC^*_{\bar{\bT}_{\bM}}&=\Set*{(\lambda,\mu)\in\CharG(\bT_+)\given \lambda\ge
        w(\mu), w\in \bW}\\
                &=\bbN\left(\Set*{(\alpha,0)\given
                \alpha\in\SimRts}\cup\Set*{(\Wt_i,w(\Wt_i))\given 1\le i\le
                r\text{ and }
        w\in \bW}\right),\\
    \sC_{\bar{\bT}_{\bM}}&=\Set*{(\dual{\lambda},\dual{\mu})\in\CoCharG(\bT_+)\given
    \dual{\lambda}\in\CoCharG(\bT^\AD)_+\text{ and }\dual{\lambda}\ge
    -w(\mu), \forall w\in \bW}.
    \nomenclature[\(C"scr^* \)]{\(\sC^*\)}{the character cone associated with an affine toric variety}
\end{align}
It is known that \(\bar{\bT}_{\bM}\) is Cohen--Macaulay.

\subsection{}
It is also possible to talk about the category \(\FM(\bG)\) for any semisimple
group \(\bG\). In this case, let \(\bG^\SC\to \bG\) be the natural isogeny and
let \(\bZ'\) be the kernel, then the universal monoid \(\Env(\bG)\) for
\(\bG\) is none other than the GIT quotient \(\Env(\bG^\SC)\git\bZ'\), and any
\(\bM\in\FM(\bG)\) is of the form \(\bM'\git\bZ'\) for some
\(\bM'\in\FM(\bG^\SC)\). The details can be found in \cite{Ri01} for example.
We may also talk about abelianization map of
\(\bM\in\FM(\bG)\), and since \(\bM\) is affine and \(\bG\) is reductive, it is
easy to see that \(\bA_\bM=\bA_{\bM'}\).

We will rarely talk about monoids in \(\FM(\bG)\) for non-simply-connected
groups \(\bG\), for two reasons: the first is that monoids in \(\FM(\bG^\SC)\)
are already sufficient for studying fundamental lemma and many other aspects of
representation theory, and the second is that the invariant quotient
\(\bG\git\bG\) is in general ill-behaved (and so is \(\bM\git \bG\) by
extension) and not suitable for developing a geometric theory. However, one can
see \Cref{sec:global_affine_schubert_scheme,sec:factorizations} for a brief use
of monoids in \(\FM(\bG^\AD)\).

\section{Invariant Theory of Reductive Monoids}%
\label{sec:invariant_theory_of_reductive_monoids}

In this section we discuss the invariant theory of reductive monoids. Most
results are built upon ones in \Cref{sec:invariant_theory_of_the_group}.

\subsection{}
Let \(\bM\in\FM(\bG^\SC)\) be any monoid.  We have the adjoint action of
\(\bG^\SC\)
on \(\bM\). It can be viewed as the restriction of the \(\bG^\SC\x
\bG^\SC\)-action to the diagonal embedding of \(\bG^\SC\). The center
\(\bZ^\SC\) acts trivially, so the action factors through \(\bG^\AD\), hence lifts to a
\(\bG\)-action on \(\bM\). The GIT quotient space \(\bM\git
\bG\) maps to \(\bA_{\bM}\). We also have the \notion{cameral cover}\index{cameral!cover}
\(\pi_{\bM}\colon\bar{\bT}_{\bM}\to \bar{\bT}_{\bM}\git\bW\).
    \nomenclature[\(pi_M \)]{\(\pi_\bM,\pi_\FRM\)}{the cameral cover associated with very flat monoid \(\bM\) or \(\FRM\)}

On the other hand, the fundamental
representation \(\rho_i\) extends to \(\Env(\bG^\SC)\) by its definition, 
still denoted by \(\rho_{i+}\).
    \nomenclature[\(rho_i+ \)]{\(\rho_{i+}\)}{the extension of \(\rho_i\) to \(\bM\in\FM(\bG^\SC)\)}
Therefore, the character function \(\chi_{i+}\)
    \nomenclature[\(chi_i+ \)]{\(\chi_{i+}\)}{the trace of \(\rho_{i+}\)}
makes sense for arbitrary \(\bM\in\FM(\bG^\SC)\) after realizing it as a
pullback of \(\Env(\bG^\SC)\). Thus, we have a map \(\bM\to
\bA_{\bM}\x \bC\), which factors through the GIT quotient map \(\chi_{\bM}\colon \bM\to
\bC_{\bM}\defeq \bM\git\bG\).
    \nomenclature[\(C"bold_M \)]{\(\bC_\bM,\FRC_\FRM\)}{the GIT quotient \(\bM\git\bG\) or \(\FRM\git G\)}
    \nomenclature[\(chi_M \)]{\(\chi_\bM,\chi_\FRM\)}{the GIT quotient map \(\bM\to\bC_\bM\) or \(\FRM\to\FRC_\FRM\)}

\begin{theorem}
    The maps \(\bar{\bT}_{\bM}\to\bM\to\bA_{\bM}\x\bC\) induce canonical isomorphisms
    \begin{align}
        \bar{\bT}_{\bM}\git\bW\simeq\bC_{\bM}\simeq\bA_{\bM}\x\bC.
    \end{align}
    In fact, the first isomorphism holds for \emph{any} (not necessarily very
    flat) monoid with \(\bM^\Der\cong\bG^\SC\).
\end{theorem}
\begin{proof}
    The first isomorphism is proved in \cite{Re88}, but the reader can also find
    it in \cite{Re05} for a more modern reference. We do not know where the
    second isomorphism was first proved, but it was at least proved in
    \cite{Bo15}. However, the reader should also see \cites{BoCh18,Ch22} since
    \cite{Bo15} contains some error unrelated to the current theorem.
\end{proof}

\begin{corollary}
    \label[corollary]{cor:cameral_cover_is_Cohen_Macaulay}
    The cameral cover \(\pi_{\bM}\) is a Cohen--Macaulay morphism, in other words, it is flat
    and has Cohen--Macaulay fibers.
\end{corollary}
\begin{proof}
    Since being a Cohen--Macaulay morphism is stable under  base change, it
    suffices to prove the claim for \(\bM=\Env(\bG^\SC)\), in which case
    \(\bC_{\bM}\) is regular. Since \(\bar{\bT}_{\bM}\) is Cohen--Macaulay, and
    \(\pi_{\bM}\) is finite, it is a flat morphism. Then it is a general result
    that flat morphism between locally Noetherian schemes with a Cohen--Macaulay
    source must be Cohen--Macaulay. See \cite{StacksP}*{\href{https://stacks.math.columbia.edu/tag/0C0X}{Tag 0C0X}}.
\end{proof}

\begin{proposition}
    Given a Steinberg quasi-section \(\epsilon^{(\xi,\dot{S})}\) of the group
    \(\bG^\SC\), the map
    \begin{align}
        \epsilon_{\bM}^{(\xi,\dot{S})}\colon \bA_{\bM}\x\bC&\longto \bM\\
        (a,x)&\longmapsto \delta_{\bM}(a)\epsilon^{(\xi,\dot{S})}(x)
        \nomenclature[\(epsilon_M_xi_S \)]{\(\epsilon_{\bM}^{(\xi,\dot{S})}\)}{the extended Steinberg quasi-section of \(\bM\in\FM(\bG^\SC)\)}
    \end{align}
    defines a quasi-section of \(\chi_{\bM}\) whose image lies in the regular locus
    \(\bM^\reg\). Moreover, the union of
    \(\img\epsilon_{\bM}^{(\xi,\dot{S})}\)
    over all Coxeter data \((\xi,\dot{S})\) generates \(\bM^\reg\) under
    \(\bG^\SC\)-action.
\end{proposition}
\begin{proof}
    This is a corollary of \cite{Ch22}*{Lemma~2.2.8, Propositions~2.2.9, 2.2.19}
    for \(\bM=\Env(\bG^\SC)\), and the general
    case is deduced by pulling back.
\end{proof}

\subsection{}%
The Weyl group \(\bW\) acts on \(\bar{\bT}_{\bM}\). Similar to the
group case, a \(\bG\)-orbit in \(\bM\) is called
\notion{regular}\index{orbit!regular} if the stabilizer
achieves minimal dimension, \notion{semisimple}\index{orbit!semisimple} if it
contains an element in \(\bar{\bT}\), and \notion{regular
semisimple}\index{orbit!regular semisimple} if it is
both regular and semisimple. 

\subsection{}%
The regular semisimple locus \(\bM^\rss\) is open and smooth, and can be
characterized by a discriminant function extending \({\Disc}\) from the group
case. Indeed, we only need to define it for
\(\bM=\Env(\bG^\SC)\) and pull it back to general \(\bM\). For \(\bM=\Env(\bG^\SC)\), let
\begin{align}
    {\Disc_+}\defeq e^{(2\rho,0)}\prod_{\alpha\in\Roots}(1-e^{(0,\alpha)}).
    \nomenclature[\(D{}isc_+ \)]{\(\Disc_+\)}{the extended discriminant function}
\end{align}
This function extends to \(\bar{\bT}_{\bM}\) and is \(\bW\)-Invariant, hence descends to a
function on \(\bC_{\bM}\). It is called the \notion{extended discriminant
function}\index{function!extended discriminant}\index{discriminant!function,
extended}, and defines the \notion{extended discriminant
divisor}\index{divisor!extended discriminant}\index{discriminant!divisor, extended} \(\bD_{\bM}\subset\bC_{\bM}\).
    \nomenclature[\(D"bold_M \)]{\(\bD_\bM,\FRD_\FRM\)}{the extended discriminant divisor of very flat monoid \(\bM\) or \(\FRM\)}
The complement of \(\bD_{\bM}\) is the regular semisimple locus
\(\bC_{\bM}^\rss\), whose preimage under \(\chi_{\bM}\) is \(\bM^\rss\).
For convenience, we will also regard \(\bE_\bM\) as a divisor on \(\bC_\bM\) by
pulling back (and it will not cause any confusion), and denote the
\notion{invertible locus}\index{locus!invertible} \(\bC_\bM^\x\defeq\bA_\bM^\x\x\bC\).

\begin{lemma}
    \label[lemma]{lem:disc_and_num_bd_intersects_properly}
    The divisor \(\bD_{\bM}\) intersects \(\bC_{\bM}-\bC_{\bM}^\x\) properly. In other
    words, the codimension of
    \(\bD_{\bM}^\x\defeq\bD_{\bM}\cap(\bC_{\bM}-\bC_{\bM}^\x)\) in \(\bC_{\bM}\)
    is at least \(2\). In particular, set-theoretically \(\bD_{\bM}\) is the
    closure of \(\bD_{\bM}^\x\).
\end{lemma}
\begin{proof}
    The proof is essentially in \cite{Ch22}*{Lemma~2.4.2}, and we reproduce it
    here. Without loss of generality, we may assume that \(\bM\) has \(0\)
    (since we can always find a very flat reductive monoid with \(0\) containing
    \(\bM\) by enlarging the cocharacter cone).
    When \(\bM=\Env(\bG^\SC)\), consider the idempotent
    \(e_{\emptyset,\SimRts}\) (see \Cref{sub:desc_of_big_cell_idempotents}). It
    can be explicitly computed that \(e_{\emptyset,\SimRts}\) is regular semisimple. Its image in
    \(\bA_{\bM}\) is \(0\) which is contained in every irreducible component of
    \(\bA_{\bM}-\bA_{\bM}^\x\). Using the universal property of \(\Env(\bG^\SC)\),
    we see that for any \(\bM\) the image of the regular semisimple locus is dense in every
    irreducible component of \(\bA_{\bM}-\bA_{\bM}^\x\), and we are done.
\end{proof}

\begin{lemma}
    \label[lemma]{lem:extended_discriminant_divisor_is_reduced}
    The divisor \(\bD_{\bM}\) is a reduced divisor.
\end{lemma}
\begin{proof}
    Since \(\bD_{\bM}\) is a principal divisor in a Cohen--Macaulay scheme, it is
    itself Cohen--Macaulay. By \Cref{lem:disc_and_num_bd_intersects_properly}, we only need to prove
    that it is reduced over the invertible locus \(\bD_{\bM}^\x\). But then it
    suffices to consider the divisor in \(\bT_{\bM}\) cut out by
    \((1-e^{(0,\Rt)})(1-e^{(0,-\Rt)})\) for a single root \(\Rt\), and it
    reduces to groups with derived subgroup \(\SL_2\),
    in which case we can directly compute.
    The argument is parallel to for example \cite{Ng10}*{Lemme~1.10.1}. 
    We leave the details to the reader.
\end{proof}

\subsection{}%
Contrary to the Lie algebra case, there are in general more than one open orbits
in the fibers of \(\chi_{\bM}\), in other words, the fibers of
\(\chi_{\bM}^\reg\colon \bM^\reg\to \bC_{\bM}\) are no longer homogeneous
\(\bG\)-spaces.

\subsection{}%
We would like to define the regular centralizer group scheme \(\bJ_{\bM}\to\bC_{\bM}\)
similarly to the group case, but the original descent argument needs some
adaptation, due to the fact that a
fiber of \(\chi_{\bM}\) may contain multiple regular orbits. 

The key observation used by \cite{Ch22} is that the numerical boundary divisor
\(\bE_{\bM}\)
and the discriminant divisor \(\bD_{\bM}\) intersect properly
(\Cref{lem:disc_and_num_bd_intersects_properly}). The descent
argument works over \(\bM^\x\cup \bM^\rss\) which is an open subset whose
complement has codimension at least \(2\). We leave the details to
\cite{Ch22}*{Lemma~2.4.2}.

\begin{proposition}
    There is a unique smooth commutative group scheme \(\bJ_{\bM}\to \bC_{\bM}\)
    \nomenclature[\(J_M \)]{\(\bJ_\bM,\FRJ_\FRM\)}{the regular centralizer group scheme over \(\bC_\bM\) or \(\FRC_\FRM\)}
    with a \(\bG\)-equivariant isomorphism
    \begin{align}
        \chi_{\bM}^*\bJ_{\bM}|_{\bM^{\reg}}\stackrel{\sim}{\longto} \bI_{\bM}^{\reg},
    \end{align}
    which can be extended to a homomorphism \(\chi_{\bM}^*\bJ_{\bM}\to \bI_{\bM}\).
\end{proposition}

There is also a description of the regular centralizer using cameral cover
\(\pi_{\bM}\colon\bar{\bT}_{\bM}\to \bC_{\bM}\). As in the group case, let
\begin{align}
    \Pi_{\bM}=\pi_{\bM*}(\bT\x \bar{\bT}_{\bM}),
\end{align}
and
\begin{align}
    \bJ^1_{\bM}=\Pi_{\bM}^{\bW}.
    \nomenclature[\(J_M_1 \)]{\(\bJ_\bM^1,\FRJ_\FRM^1\)}{the analogue of \(\bJ^1\) (resp.~\(\FRJ^1\)) for \(\bJ_\bM\) (resp.~\(\FRJ_\FRM\))}
\end{align}
Similarly, we have subfunctor \(\bJ_{\bM}'\):
for a \(\bC_{\bM}\)-scheme \(S\), \(\bJ_{\bM}'(S)\) consists of points 
\begin{align}
    f\colon S\x_{\bC_{\bM}}\bar{\bT}_{\bM}\to \bT
\end{align}
such that for every geometric point \(x\in S\x_{\bC_{\bM}}\bar{\bT}_{\bM}\), if
\(s_\Rt(x)=x\) for a root \(\Rt\), then \(\Rt(f(x))\neq -1\). With the same
proof as in the group case, this is an open subgroup scheme of \(\bJ^1_{\bM}\)
containing the fiberwise neutral component \(\bJ^0_{\bM}\).
    \nomenclature[\(J_M_0 \)]{\(\bJ_\bM^0,\FRJ_\FRM^0\)}{the analogue of \(\bJ^0\) (resp.~\(\FRJ^0\)) for \(\bJ_\bM\) (resp.~\(\FRJ_\FRM\))}

\begin{proposition}\label[proposition]{prop:Galois_reg_cent_monoid}
    There is a canonical open embedding
    \begin{align}
        \bJ_{\bM}\longto \bJ^1_{\bM}
    \end{align}
    that identifies \(\bJ_{\bM}\) with \(\bJ_{\bM}'\).
\end{proposition}
\begin{proof}
    The open embedding claim is proved in \cite{Ch22}*{Proposition~2.4.7}. We
    even have \(\bJ_{\bM}=\bJ_{\bM}^1=\bJ_{\bM}'\) if \(\bG=\bG^\SC\). The
    point is that the complement of \(\bC_{\bM}^\x\cup\bC_{\bM}^\rss\) has
    codimension \(2\), thus we only need to prove the claim over this open
    locus, but then it is a consequence of the group case and regular semisimple
    case, either of which is easy.
\end{proof}

\begin{remark}
    In \cite{Ch22}*{Proposition~2.4.7}, the open embedding
    \(\bJ_\bM\to\bJ^1_\bM\) relies on the Grothendieck's simultaneous resolution
    for monoids, which is analogues to the well-known group case:
    \begin{equation}
        \begin{tikzcd}
            \widetilde{\bM}\defeq\bG^\SC\x^{\bB^\SC}\bar{\bB}_\bM \ar[r, "\tilde{\pi}_{\bM}"]\ar[d] & \bM \ar[d]\\
            \bar{\bT}_{\bM} \ar[r, "\pi_{\bM}"] & \bC_{\bM}
        \end{tikzcd},
    \end{equation}
    where \(\bar{\bB}_\bM\) is the closure of the Borel \(\bB_\bM\)
    \nomenclature[\(B"bold_M \)]{\(\bB_\bM,B_\FRM\)}{the Borel subgroup of \(\bM^\x\) (resp.~\(\FRM^\x\)) induced by \(\bB\) (resp.~\(B\))}
    of \(\bM\) induced by
    \(\bB\). The crucial part is to define the vertical map on the left by
    defining a canonical projection map
    \(\bar{\bB}_\bM\to\bar{\bT}_\bM\). However, the definition of such
    projection in \textit{loc.~cit.} (due to \cite{Bo17}) is wrong. The correct
    definition is given in \Cref{prop:def_of_monoidal_Borel_projection} below.
    Moreover, \cite{Bo17} wrongly claims that the square becomes Cartesian after
    restricting to \(\bM^\reg\). In fact, it is only Cartesian over
    \(\bM^\rss\cup\bM^{\x,\reg}\), which is still sufficient for defining
    \(\bJ_\bM\to\bJ^1_\bM\).
\end{remark}

\begin{corollary}
    \label[corollary]{cor:M_reg_is_union_of_gerbs}
    The map \(\Stack{\bM^\reg/\bG}\to\bC_{\bM}\) is a finite union of
    \(\bJ\)-gerbes.
\end{corollary}

\subsection{}%
There is an open subset  \(\bM^\circ\subset\bM\)
    \nomenclature[\(M^circ \)]{\(\bM^\circ,\FRM^\circ\)}{the big-cell locus of very flat monoid \(\bM\) or \(\FRM\)}
containing \(\bM^\reg\)
that has significant representation-theoretic meaning. We will call it
the \notion{big-cell locus}\index{locus!big-cell}.
We will define it for \(\Env(\bG^\SC)\), and for general \(\bM\in\FM(\bG^\SC)\)
one simply takes the preimage under any morphism \(\bM\to\Env(\bG^\SC)\) in the
same category.

For \(\Env(\bG^\SC)\), the big-cell locus is the set
\begin{align}\label{eqn:big_cell_locus_def}
    \Env(\bG^\SC)^\circ=\Set*{\,x\in\Env(\bG^\SC)\given \rho_{i+}(x)\neq 0, 1\le
    i\le r\,}.
\end{align}
This is readily seen an open locus and can be shown to contain the regular
locus. The central torus \(\bZ_+\) acts on \(\Env(\bG^\SC)^\circ\), in
fact freely, and the quotient is isomorphic to the wonderful compactification
\(\bar{\bG^\AD}\) of
\(\bG^\AD\) (cf.,~\cite{Bo15}*{Proposition~2.2}). Therefore, \(\Env(\bG^\SC)^\circ\) is
smooth. Therefore, \(\Env(\bG^\SC)^\circ\) is the total space of a
\(\bZ_+\)-torsor over the projective variety \(\bar{\bG^\AD}\), which is quasi-affine,
while \(\Env(\bG^\SC)\) is its affinization, which may be seen as some sort of affine cone.

\subsection{}%
\label{sub:desc_of_big_cell_idempotents}
There is another description of the big-cell locus using the idempotents.
Idempotents are very important in studying reductive monoids because they
represent \(\bM^\x \x \bM^\x\)-orbits in \(\bM\). 
Idempotents allows a more conceptually pleasing
description of \(\Env(\bG^\SC)^\circ\), while \eqref{eqn:big_cell_locus_def}, although concise, is not
very revealing.

\begin{definition}
    A pair \((I,J)\) of subsets of \(\SimRts\) is called
    \notion{essential}\index{subset(s) of simple roots!essential pair of} if
    no connected components of \(\SimRts-J\) in the Dynkin diagram
    is entirely contained in \(I\)
\end{definition}
Pairs \((I,J)\) of subsets of \(\SimRts\) can be partially ordered by
inclusion condition.
There is a finite partially ordered set of idempotents \(e_{I,J}\in\Env(\bG^\SC)\)
    \nomenclature[\(e_I_J \)]{\(e_{I,J}\)}{the idempotent of \(\Env(\bG^\SC)\) associated with an essential pair \((I,J)\)}
labeled by essential pairs \((I,J)\). The \(\bG_+\x \bG_+\)-orbits in
\(\Env(\bG^\SC)\) are in bijection
with \(e_{I,J}\) in an order-preserving way. In other words,
\(\bG_+e_{I,J}\bG_+\) is contained in the closure of \(\bG_+e_{I',J'}\bG_+\) if
and only if \(I\subset I'\) and \(J\subset J'\).
Similarly, the \(\bT^\AD\)-orbits in
\(\bA_{\Env(\bG^\SC)}\) is in bijection with idempotents \(e_I\)
    \nomenclature[\(e_I \)]{\(e_{I}\)}{the idempotent of \(\bA_{\Env(\bG^\SC)}\) associated with \(I\subset\SimRts\)}
for \(I\subset \SimRts\), and readers can easily figure out what they are. We have that
\(\alpha_{\Env(\bG^\SC)}(e_{I,J})=e_I\). In
particular, \(e_{\SimRts,\SimRts}=1\) and \(e_{\emptyset,\emptyset}=0\). The
big-cell locus is the union
\begin{align}
    \Env(\bG^\SC)^\circ=\bigcup_{I\subset\SimRts}\bG_+ e_{I,\SimRts}\bG_+.
\end{align}
Then it is readily seen that \(\Env(\bG^\SC)^\circ\) is smooth because the restriction
of \(\alpha_{\Env(\bG^\SC)}\) to \(\Env(\bG^\SC)^\circ\) is flat with smooth fibers (homogeneous spaces 
of \(\bG^\SC\x \bG^\SC\)), and the base \(\bA_{\Env(\bG^\SC)}\) is itself
smooth.

\subsection{}
More properties of the idempotents \(e_{I,J}\) are described by Vinberg in
\cite{Vi95}, and we recall it here for future reference. 

For any subset \(I\subset \SimRts\), let \(\bP_I\subset\bG^\SC\) be the standard
    \nomenclature[\(P_I"bold \)]{\(\bP_I\)}{the standard parabolic subgroup associated with \(I\subset\SimRts\)}
parabolic subgroup containing \(\bB^\SC\)
corresponding to \(I\), and \(\bP_I^-\)
    \nomenclature[\(.- \)]{\((\cdot)^-\)}{the subgroup opposite to a parabolic subgroup or its unipotent radical}
its opposite.
Let \(\bU_I\)
    \nomenclature[\(U_I"bold \)]{\(\bU_I\)}{the unipotent radical of \(\bP_I\)}
(resp. \(\bU_I^-\))
the unipotent radical of \(\bP_I\) (resp. \(\bP_I^-\)), and
\(\bL_I=\bP_I\cap \bP_I^-\)
    \nomenclature[\(L_I"bold \)]{\(\bL_I\)}{the standard Levi subgroup contained in \(\bP_I\)}
the standard Levi subgroup.
Let \(\pr_I\) (resp. \(\pr_I^-\)) the projection from
\(\bP_I\) (resp. \(\bP_I^-\)) to \(\bL_I\).
Let \(\bP_{I+}, \bU_{I+}, \bL_{I+}\), etc. be
the same constructions in \(\bG_+\).

For each \(I,J\subset \SimRts\), define \(I^c\) to be the set \(\SimRts - I\),
\(I^\circ\) to be the \notion{interior}\index{subset(s) of simple roots!interior
of} of \(I\), that is, those in \(I\) that
is not joint with \(I^c\) by an edge in the Dynkin diagram, and
\(\Sigma_{I,J}=(I\cap J^\circ)\cup J^c\). Let
\(D_I\) be the abelian monoid in \(\CharG(\bT^\SC)\) generated by \(I\), and
\(C_J\) the one generated by such \(\Wt_j\) that \(\Rt_j\in J\). Using
identification \(\CharG(\bT_+)\subset
\CharG(\bZ_+)\x\CharG(\bT^\SC)=\CharG(\bT^\SC)\x\CharG(\bT^\SC)\), we let
\begin{align}
    F_{I,J}&=\Set{(\lambda_1,\lambda_2)\in \CharG(\bT_+)\mid \lambda_1-\lambda_2\in D_I, \lambda_2\in C_J},\\
    \bT_{I,J}&=\Set{t\in \bT_+\mid \lambda(t)=1\text{ for all }\lambda\in F_{I,J}}.
\end{align}
The stabilizer of \(e_{I,J}\) in \(\bG_+\x\bG_+\) for an essential pair
\((I,J)\) is the subgroups of \(\bP_{\Sigma_{I,J}+}\x\bP_{\Sigma_{I,J}+}^-\)
consisting of elements \((g,g^-)\) such that
\begin{align}
    \label{eqn:idempotent_stab_info}
    \pr_{\Sigma_{I,J}+}(g)\equiv
    \pr_{\Sigma_{I,J}+}^-(g^-)\bmod\bL_{J^c}^\Der \bT_{I,J}.
\end{align}
The idempotent \(e_{I,J}\) itself is characterized by
\begin{align}
    (\Rt_i,0)(e_{I,J})&=
    \begin{cases}
        1 & \Rt_i\in I\\
        0 & \Rt_i\not\in I
    \end{cases},\\
    (\Wt_j,\Wt_j)(e_{I,J})&= 
    \begin{cases}
        1 & \Rt_j\in J\\
        0 & \Rt_j\not\in J
    \end{cases},\\
    (\Wt_j,w(\Wt_j))(e_{I,J})&=
    \begin{cases}
        1 & \Rt_j\in J, \Wt_j-w(\Wt_j)\in D_{I\cap J^\circ}\\
        0 & \text{otherwise}
    \end{cases}.
\end{align}

\begin{proposition}
    \label[proposition]{prop:def_of_monoidal_Borel_projection}
    There exists a canonical projection map
    \begin{align}
        \bar{\bB}_\bM\longto\bar{\bT}_\bM
    \end{align}
    whose restriction to \(\bB_\bM\) is the quotient map
    \(\bB_\bM\to\bT_\bM\simeq\bB_\bM/\bU_\bM\).
\end{proposition}
\begin{proof}
    Let \(\lambda\in\CoCharG^\SC\) be a strictly dominant cocharacter and
    \(t\in\lambda(\Gm)\subset\bT^\SC\) be a regular element.
    By \cite{Re05}*{Theorem~13.5}, the submonoid \(\bM_t\defeq\bM\) commuting
    with \(t\) is generated by the centralizer of \(t\) in \(\bB_\bM\) and
    elements of shape \(\dot{w}_1e_{I,J}\dot{w}_2\) (where
    \(\dot{w}_1,\dot{w}_2\in\Norm_{\bG^\SC}(\bT^\SC)\)) that commute with \(t\).

    By \eqref{eqn:idempotent_stab_info} and a straightforward computation
    similar to \cite{Ch22}*{Lemma~2.2.14}, it is easy to see that we must have
    \begin{align}
        (w_1^{-1}-w_2)(\lambda)\in \bL_{J^c}^\Der.
    \end{align}
    By a classical result of Borel-Tits
    (\cite{BoTi65}*{Proposition~12.17}), this means that
    \begin{align}
        (w_1^{-1}-w_2)(\lambda)=(ww_2-w_1)(\lambda)
    \end{align}
    for some \(w\in\bW_{J^c}\), or in other words, \(w_1ww_2\) fixes
    \(\lambda\). Note that
    \(e_{I,J}\dot{w}e_{I,J}=e_{I,J}e_{I,J}=e_{I,J}\), and so
    \begin{align}
        \dot{w}_1e_{I,J}\dot{w}_2=\dot{w}_1e_{I,J}\dot{w}_2(\dot{w}_2^{-1}\dot{w}^{-1}\dot{w}_1^{-1})\dot{w}_1e_{I,J}\dot{w}_2.
    \end{align}
    This is the condition required by \cite{Re05}*{Theorem~13.8} for \(\bM_t\) to
    be irreducible. Since such locus clearly contains \(\bT_\bM\) as a dense
    open subset (because \(t\) is regular), we then have
    \(\bM_t=\bar{\bT}_\bM\).

    Clearly, \(\bB_\bM\) is contained in the attractor locus of the adjoint
    action of \(\lambda(\Gm)\), thus so is \(\bar{\bB}_\bM\) because everything
    is affine. Thus the projection map is defined to be the contraction map
    \begin{align}
        \bar{\bB}_\bM\longto \bar{\bB}_\bM\git\lambda(\Gm)
    \end{align}
    by taking the limit at parameter \(0\). To see that this map does not depend
    on the choice of \(\lambda\), we note that
    \begin{align}
        \bigl(\bar{\bB}_\bM\git\lambda(\Gm)\bigr)\git\bT^\SC\simeq
        \bar{\bB}_\bM\git\bT^\SC,
    \end{align}
    and we are done.
\end{proof}

\subsection{}%
The central group \(\bZ_{\bM}\)
    \nomenclature[\(Z_M \)]{\(\bZ_\bM,Z_\FRM\)}{the full central subgroup of a very flat monoid \(\bM\) or \(\FRM\)}
acts on \(\bM\) by translation.
This action commutes with the adjoint action of \(\bG\), hence descends to an
action on \(\bC_{\bM}\), making \(\chi_{\bM}\) a \(\bZ_{\bM}\)-equivariant map. The
\(\bZ_{\bM}\)-action on \(\bC_{\bM}\) has a simple description as follows: on
\(\bA_{\bM}\) it is simply the translation by torus
\(\alpha_{\bM}(\bZ_{\bM})=\bZ_{\bM}/\bZ^\SC\), while on \(\bC\) it is the translation
action given by weights \(\Wt_i\circ\phi_{\bZ_{\bM}}\) if we use the coordinates
\(\chi_{i+}\) on \(\bC\) (recall that \(\phi_{\bZ_{\bM}}\) is the map \(\bZ_{\bM}\to
\bT^\SC\) induced by a choice of morphism \(\bM\to\Env(\bG^\SC)\) in
\(\FM(\bG^\SC)\)). This action can be lifted to an action of \(\bZ_{\bM}\) on
\(\bJ_{\bM}\) compatible with the group scheme structure by looking at the
construction of \(\bJ_{\bM}\).

On the other hand, any choice of a Steinberg quasi-section is far from being
\(\bZ_{\bM}\)-equivariant. However, one can rectify this with some technical
modification. This modification first appears in \cite{Bo15} and later used by
\cite{Ch22}. However, the proof in \cite{Bo15} contains an elementary but
serious mistake, so we include a corrected proof here. None of the results in
\cite{Ch22} is affected by this error.

\begin{proposition}\label[proposition]{prop:St_section_Z_equivariance}
    For each Coxeter datum \((\xi,\dot{S})\), one can define an action
    \(\tau_{\bM}^{(\xi,\dot{S})}\) of \(\bZ_{\bM}\) on \(\bM\) such that
    \begin{align}\label{eqn:St_section_equivariance}
        \epsilon_{\bM}^{(\xi,\dot{S})}\circ\tau_{\bC_{\bM}}(z^c)
        =\tau_{\bM}^{(\xi,\dot{S})}(z)\circ\epsilon_{\bM}^{(\xi,\dot{S})},
    \end{align}
    where \(\tau_{\bC_{\bM}}\) is the natural action of \(\bZ_{\bM}\) on \(\bC_{\bM}\), and
    \(c=\abs{\bZ^\SC}\). Moreover, for a fixed \(z\),
    \(\tau_{\bM}^{(\xi,\dot{S})}(z)\) is a composition of translation by  \(z^c\)
    and conjugation by an element in \(\bT^\SC\) determined by a homomorphism
    \(\bZ_{\bM}\to \bT^\SC\) independent of \(z\).
\end{proposition}
\begin{proof}
    We re-label \(\Set{\Wt_i}\) in such a way that \(\Wt_j\) is the weight
    corresponding to root \(\beta_j=\alpha_{\xi(j)}\).  To simplify notations,
    for \(z\in \bZ_{\bM}\), we will denote \(\phi_{\bZ_{\bM}}(z)\) simply by \(z_T\).

    Let \((a,x) \in \chi_{\bM}(\bM^\x)\subset\bC_{\bM}\). Fix a \(\Wt_i\) (\(1\le i\le
    r\)), and a weight vector 
    \(0\neq v\in V_{\Wt_i}[\mu]\), where \(\mu\le \Wt_i\) is a weight
    of \(\rho_{i}\) such that \(\mu=\sum_{j=1}^r m_j\Wt_j\).  We have (to
    simplify notations we omit \(\rho_{i+}\) in the computations):
    \begin{align}
        \left[ \epsilon^{(\xi,\dot{S})}(x) \right](v) 
          &= \left( \prod_{j=1}^r\bU_{\beta_{j}}(x_{j})\dot{s}_{\beta_j} \right)(v)\\
          &= \sum_{\substack{k_j\ge -m_j\\1\le j\le r}}\left(
          \prod_{l=1}^rx_l^{k_l+m_l} \right)v_{\ul{k}},
    \end{align}
    where \(\ul{k}=\left(k_1,\ldots,k_r\right)\) is a multi-index,
    \(v_{\ul{k}}\) is some vector, independent of any \(x_j\), of weight
    \begin{align}
        \mu_{\ul{k}}\defeq \mu+0_{\ul{k}}, 
    \end{align}
    and
    \begin{align}
        0_{\ul{k}}\defeq \sum_{d=1}^{r} \sum_{1\le
                l_1<\cdots< l_d\le r} \left[k_{l_d} \prod_{e=1}^{d-1}\left(
            -\Pair{\beta_{l_{e+1}}}{\beta_{l_{e}}^\vee} \right)\beta_{l_1}\right].
    \end{align}
    Note that \((k_1,\ldots,k_r)\mapsto 0_{\ul{k}}\) is a group isomorphism from
    \(\bbZ^r\) to root lattice \(\bbZ\Roots\).
  
    Thus, we have that
    \begin{align}
        \left[ \epsilon_{\bM}^{(\xi,\dot{S})}(a,x) \right](v)
        &= \delta_{\bM}(a)\left[ \epsilon^{(\xi,\dot{S})}(x) \right](v)\\
        &=(z_a,z_{a,T}^{-1})\sum_{\substack{k_j\ge -m_j\\1\le j\le r}}\left(
            \prod_{l=1}^r x_l^{k_l+m_l}\right)v_{\ul{k}}\\
        &=\Wt_i(z_{a,T})\sum_{\substack{k_j\ge -m_j\\1\le j\le r}}\left(
            \prod_{l=1}^rx_l^{k_l+m_l} \right)\mu_{k}(z_{a,T})^{-1}v_{\ul{k}},
    \end{align}
    where \(z_a\in \bZ_{\bM}\) is some element that \(\alpha_{\bM}(z_a)=a\); and that
    \begin{align}\label{eqn:tau_FRCM_action}
        \left[\epsilon_{\bM}^{(\xi,\dot{S})}
            \left(\tau_{\bC_{\bM}}(z)(a,x)\right)\right](v)
        &=(zz_{a},(z_Tz_{a,T})^{-1})\epsilon^{(\xi,\dot{S})}(\Wt_j(z_T)x_j)(v)\\
        &=\Wt_i(z_Tz_{a,T})\sum_{\substack{k_j\ge -m_j\\1\le j\le r}}\left(
            \prod_{l=1}^r(\Wt_l(z_T)x_1)^{k_l+m_l}
            \right)\mu_{k}(z_Tz_{a,T})^{-1}v_{\ul{k}},
    \end{align}
    for any \(z\in \bZ_{\bM}\).

    On the other hand, for \(t\in \bT^\SC\),
    \begin{align}\label{eqn:adjoint_action_and_weight_vector}
        \left[\Ad_{t}\left(\epsilon_{\bM}^{(\xi,\dot{S})}(a,x)\right)\right](v)
        &=t\left[\epsilon_{\bM}^{(\xi,\dot{S})}(a,x)\right]t^{-1}(v)\\
        &=\Wt_i(z_{a,T})\mu(t)^{-1}\sum_{\substack{k_j\ge -m_j\\1\le j\le
            r}}\left( \prod_{l=1}^rx_l^{k_l+m_l} \right)\mu_{k}(tz_{a,T}^{-1})v_{\ul{k}}\\
        &=\Wt_i(z_{a,T})\sum_{\substack{k_j\ge -m_j\\1\le j\le
            r}}\left( \prod_{l=1}^rx_l^{k_l+m_l} \right)
            0_{\ul{k}}(t)\mu_{k}(z_{a,T})^{-1}v_{\ul{k}}.
    \end{align}

    Consider commutative diagram
    \begin{equation}\label{eqn:Cox_modified_T_iso}
        \begin{tikzcd}
            \bZ_{\bM} \ar[r, "\Wt_\bullet"] & \Gm^r\ar[d, "\sim", swap] & 
                \bT^\SC \ar[l, "0_{(\bullet)}", swap]\ar[ld, "\beta_\bullet", swap]
                    \ar[d, "z\mapsto z^c"]\\
                  & \beta_\bullet(\bT^\SC) \ar[r]& \bT^\SC
        \end{tikzcd},
    \end{equation}
    where \(\beta_\bullet\) is the map
    \begin{align}
        t\longmapsto (\beta_1(t),\ldots,\beta_r(t)),
    \end{align}
    \(\Wt_\bullet\) is the map
    \begin{align}
        z\longmapsto (\Wt_1(z_T),\ldots,\Wt_r(z_T)),
    \end{align}
    and \(0_{(\bullet)}\) is the one
    \begin{align}
        t\longmapsto (0_{(1,0,\ldots, 0)}(t), \ldots, 
            0_{(0,\ldots, 1,\ldots, 0)}(t), \ldots, 
            0_{(0,\ldots, 0,1)}(t)).
    \end{align}
    Denote by \(\psi\) the map \(\bZ_{\bM}\to \bT^\SC\) from the upper-left to
    the lower-right in \eqref{eqn:Cox_modified_T_iso}, then we have that
    \(\Wt_\bullet(z^c)=0_{(\bullet)}(\psi(z))\).

    Now For \(z\in \bZ_{\bM}\), define \(\tau_{\bM}^{(\xi,\dot{S})}(z)\) to be the
    composition of translation by \((z^c,1)\) and conjugation by
    \(\psi(z)z_T^{-c}\), in other words, for \((t, g)\in \bM^\x\),
    \begin{align}
        \tau_{\bM}^{(\xi,\dot{S})}(z)\colon (t,g)\mapsto (z^ct,\Ad_{\psi(z)z_T^{-c}}(g)).
    \end{align}

    Then one sees from \eqref{eqn:tau_FRCM_action} and
    \eqref{eqn:adjoint_action_and_weight_vector} that
    \begin{align}
        \left[\tau_{\bM}^{(\xi,\dot{S})}(z)
            \left(\epsilon_{\bM}^{(\xi,\dot{S})}(a,x) \right)\right](v)
        &=(z^c,1)\Wt_i(z_{a,T})\sum_{\substack{k_j\ge -m_j\\1\le j\le r}}
            \left(\prod_{l=1}^rx_l^{k_l+m_l} \right)0_{\ul{k}}(\psi(z)z_T^{-c})
            \mu_{k}(z_{a,T})^{-1}v_{\ul{k}}\\
        &=\Wt_i(z_T^cz_{a,T})\sum_{\substack{k_j\ge -m_j\\1\le j\le r}}
            \left( \prod_{l=1}^rx_l^{k_l+m_l} \right)
            0_{\ul{k}}(\psi(z)z_T^{-c})\mu_{k}(z_{a,T})^{-1}v_{\ul{k}}\\
        &=\Wt_i(z_T^cz_{a,T})\sum_{\substack{k_j\ge -m_j\\1\le j\le r}}
            \left( \prod_{l=1}^r(\Wt_l(z_T)^cx_l)^{k_l+m_l} \right)
            \mu_{k}(z_T^cz_{a,T})^{-1}v_{\ul{k}}\\
        &=\left[ \epsilon_{\bM}^{(\xi,\dot{S})}
           \left(\tau_{\bC_{\bM}}(z^c)(a,x)\right)\right](v).
  \end{align}

    Finally, clearly the images under \(\alpha_{\bM}\) of
    \(\tau_{\bM}^{(\xi,\dot{S})}(z)\left(\epsilon_{\bM}^{(\xi,\dot{S})}(a,x) \right)\) and
    \(\epsilon_{\bM}^{(\xi,\dot{S})}\left(\tau_{\bC_{\bM}}^{(\xi,\dot{S})}(z^c)(a,x)\right)\)
    are the same, being \(\alpha_{\bM}(z^c,1)a\).

    Therefore, since \(v\), \(\mu\), and \(\Wt_i\) are arbitrary, we know that
    \(\tau_{\bM}^{(\xi,\dot{S})}(z)\circ\epsilon_{\bM}^{(\xi,\dot{S})}
    =\epsilon_{\bM}^{(\xi,\dot{S})}\circ\tau_{\bC_{\bM}}(z^c)\)
    when restricted to \(\chi_{\bM}(\bM^\x)\). Since \(\chi_{\bM}(\bM^\x)\) is dense in
    \(\bC_{\bM}\), we are done.
\end{proof}

\subsection{}%
\label{sub:equivariant_ST_sec}
The GIT quotient map \(\chi_{\bM}\) induces map
\begin{align}
    \bar{\chi_{\bM}}\colon \Stack{\bM/\bG}\longto \bC_{\bM},
\end{align}
which is \(\bZ_{\bM}\)-equivariant. We have the further induced map
\begin{align}
    \Stack{\chi_{\bM}}\colon \Stack{\bM/(\bG\x \bZ_{\bM})}\longto \Stack{\bC_{\bM}/\bZ_{\bM}}.
\end{align}
Choose a Coxeter datum \((\xi,\dot{S})\), the quasi-section
\(\epsilon_{\bM}^{(\xi,\dot{S})}\) induces a quasi-section of \(\bar{\chi_{\bM}}\), but
not of \(\Stack{\chi_{\bM}}\) unless \(c=\abs{\bZ^\SC}=1\). 

To fix this, recall we have \(\psi\colon \bZ_{\bM}\to \bT^\SC\) in the proof of
\Cref{prop:St_section_Z_equivariance}, also viewed as a morphism into
\(\bT\) by abuse of notations. Let
\begin{align}
    \Psi\colon \bZ_{\bM}&\longto \bT\\
    z&\longmapsto \psi(z)z_T^{-c}.
\end{align}
We define stack \(\Stack{\bM/(\bG\x \bZ_{\bM})}_{[c]}\) using pullback Cartesian diagram
\begin{equation}
    \begin{tikzcd}
        \Stack{\bM/(\bG\x \bZ_{\bM})}_{[c]}\ar[r]\ar[d] & \Stack{\bM/(\bG\x \bZ_{\bM})}\ar[d]\\
        \BG{\bZ_{\bM}}\ar[r, "\mu\,\mapsto\,\mu^{\otimes c}"] & \BG{\bZ_{\bM}}
    \end{tikzcd},
\end{equation}
where we use \(\BG{G}\)
    \nomenclature[\(B"bbold. \)]{\(\BG{(\cdot)}\)}{the classifying stack of a group}
to denote the classifying stack of a group \(G\);
and similarly we define \(\Stack*{\bC_{\bM}/\bZ_{\bM}}_{[c]}\).
By \Cref{prop:St_section_Z_equivariance}, we obtain a quasi-section
\begin{align}
    \Stack*{\epsilon_{\bM}^{(\xi,\dot{S})}}_{[c]}\colon 
    \Stack*{\bC_{\bM}/\bZ_{\bM}}_{[c]}\longto \Stack{\bM/(\Psi\x\Id)(\bZ_{\bM}))}_{[c]},
\end{align}
which, by composition, induces quasi-section
\begin{align}
    \Stack*{\epsilon_{\bM}^{(\xi,\dot{S})}}_{[c]}\colon 
    \Stack*{\bC_{\bM}/\bZ_{\bM}}_{[c]}\longto \Stack{\bM/(\bG\x \bZ_{\bM})}_{[c]}.
\end{align}
Finally, note that the images of all these quasi-sections lie inside the regular
locus (that is, the images of \(\bM^\reg\) in the respective quotient stacks).

\subsection{}

The \(\Out(\bG^\SC)\) action on \(\bG^\SC\)
induces an action on \(\Env(\bG^\SC)\), hence we have a quasi-split universal
monoid \(\Env(G^\SC)\) on \(X\).
The actions of \(\Out(\bG^\SC)\) on maximal toric variety
\(\bar{\bT}_{\Env(\bG^\SC)}\) and on Weyl group \(\bW\) can be combined
into an action of \(\bW\rtimes \Out(\bG^\SC)\) on \(\bar{\bT}_{\Env(\bG^\SC)}\) and
a compatible action of \(\Out(\bG^\SC)\) on \(\bC_{\Env(\bG^\SC)}\). Let
\(\FRT_{\Env(G^\SC)}\) and \(\FRC_{\Env(G^\SC)}\) be the respective induced
twisted forms. After twisting by \(\OGT_G\), we obtain invariant map
\begin{align}
    \chi_{\Env(G^\SC)}\colon \Env(G^\SC)\longto \FRC_{\Env(G^\SC)},
\end{align}
and the quotient space is naturally isomorphic to \(\FRT_{\Env(G^\SC)}\git W\).
The regular centralizer \(\FRJ_{\Env(G^\SC)}\) is
also well-defined over \(X\), and so is its Galois
description \(\FRJ_{\Env(G^\SC)}^1\).

The Steinberg quasi-section, however, is not necessarily defined unless
\(G^\SC\) either has no non-split simple factor of type \(\TypeA_{2m}\), or all of its such
factors are split. The various
choices are carefully made in constructing the quasi-section as in the group
case. On the other hand,
if such conditions are satisfied and a Steinberg quasi-section
\(\epsilon_{\Env(G^\SC)}^{(\xi,\dot{S})}\) is defined, then the
equivariant version \(\Stack{\epsilon_{\Env(G^\SC)}^{(\xi,\dot{S})}}_{[c]}\) can also
be defined since the relevant constructions are \(\Out(\bG)\)-equivariant.

\subsection{}
If \(\bM\in\FM(\bG^\SC)\) is a very flat monoid such that the
\(\Out(\bG)\)-action on \(\bG^\SC\) extends over \(\bA_{\bM}\) compatible with map
\(\bA_{\bM}\to\bA_{\Env(\bG^\SC)}\), then we have twisted forms \(\FRM\)
    \nomenclature[\(M"frak \)]{\(\FRM\)}{a very flat reductive monoid associated with a quasi-split semisimple group}
(resp.~\(\FRT_{\FRM}\), resp.~\(\FRC_{\FRM}\), etc.)
of \(\bM\)
(resp.~\(\bar{\bT}_{\bM}\), resp.~\(\bC_{\bM}\), etc.) over \(X\). More generally,
if \(\bM\) is such that \(\bA_{\bM}\) is stable under the monodromy determined
by \(\vartheta_G^\bullet\) (but not necessarily stable under
\(\Out(\bG)\)), then we also have the pointed twisted form \(\FRM\). We also
have non-pointed version by using appropriate \'etale coverings of \(X\)
instead. The category of such monoids is denoted by \(\FM(G^\SC)\), and
\(\FM_0(G^\SC)\) is the full subcategory of monoids with \(0\).

For future convenience, we make the following definition:
\begin{definition}
    \label[definition]{def:toric_standard_type}
    Let \(A\) be a torus over \(X\) acting on a flat affine \(X\)-scheme \(\FRA\)
    such that each fiber is a \(A\)-toric variety.
    We say \(\FRA\) is \notion{of standard type}\index{monoid!toric, of standard type} if
    \begin{align}
        A=p_*(\Gm\x X')\subset \FRA=p_*(\bbA^1\x X')
    \end{align}
    where \(p\colon X'\to X\) is some finite \'etale cover.
    We call \(\FRM\in\FM(G^\SC)\) \notion{of standard
    type}\index{monoid!reductive, of standard type} if its abelianization \(\FRA_\FRM\)
    is (equivalently, \(\bA_\bM\) is isomorphic to an affine space).
\end{definition}

\subsection{}
In case \(G^\SC\) has a non-split factor whose type is a product of types
\(\TypeA_{2m}\) for various \(m\), we want to show the following result
extending Steinberg's:
\begin{theorem}
    \label[theorem]{thm:A_2m_rationality_fix}
    Let \(\FRM\in\FM(G^\SC)\) for a quasi-split group \(G\) over \(X=\Spec{K}\)
    of some field \(K\) whose characteristic is larger than twice the Coxeter
    number of \(G\). Then \(\FRM(K)\to\FRC_{\FRM}(K)\) is surjective.
\end{theorem}
\begin{proof}
    Since \(\FRM\) is the fiber product of \(\FRA_{\FRM}\) with \(\Env(G^\SC)\)
    over \(\FRA_{\Env(G^\SC)}\), it suffices to prove the result for
    \(\FRM=\Env(G^\SC)\). The result of Steinberg is the same as setting
    \(\FRM=G^\SC\). Moreover, it suffices to assume \(\bG^\SC\) consists solely
    of types \(\TypeA_{2m}\). Our proof is a modified version of Steinberg's proof.

    Indeed, we start with simple group \(\bG=\SL_{2m+1}\). We label the simple roots
    (and fundamental weights, etc.) from one end of the Dynkin diagram to the other
    end, so that the middle two simple roots are \(\Rt_m\) and \(\Rt_{m+1}\).
    Let \(\Rt=\Rt_m+\Rt_{m+1}\). One can verify that the set
    \begin{align}
        \SimRts'=\Set{\Rt_1,\ldots,\Rt_{m-1},\Rt,\Rt_{m+2},\ldots,\Rt_{2m}}
    \end{align}
    generates a root subsystem of type \(\TypeA_{2m-1}\), and it induces a
    subgroup \(\bG'\subset\bG\) isomorphic to \(\SL_{2m}\) which is stable under
    \(\Out(\bG)\). This also identifies \(\Out(\bG)\) with \(\Out(\bG')\) using
    the pinning \(\bSPL\). The fundamental representations \((\rho_i,V_i)\)
    of \(\bG\) are just \(i\)-th exterior products of the standard
    representation of \(\bG\),
    and simple linear algebra shows that the restriction of \((\rho_i,V_i)\) to
    \(\bG'\) decomposes into two irreducible representations, one with
    highest-weight \(\Wt_i\), and the other with highest-weight which we shall
    denote by
    \(\Wt_i'\in\CharG(\bT)\). The weight \(\Wt_i'\) is a weight in \(V_i\),
    and the difference \(\Wt_i-\Wt_i'\) is a linear
    combination of simple roots of \(\bG\) with coefficients in \(\bbN\).

    Consider map
    \begin{align}
        \epsilon''\colon\bbA^{2m-1}&\longto \bG'\subset\bG\\
        x=(x_1,\ldots,x_{m-1},x_\Rt,x_{m+2},\ldots,x_{2m})&\longmapsto
        \bU_\Rt(x_\Rt)\dot{s}_\Rt\prod_{i=1}^{m-1}(\bU_i(x_i)\dot{s}_i\bU_{2m-i+1}(x_{2m-i+1})\dot{s}_{2m-i+1}),
    \end{align}
    where the representatives of reflections \(\dot{s}_i\) and
    \(\dot{s}_\Rt=\dot{s}_{m+1}\dot{s}_m\dot{s}_{m+1}\)
    are so chosen that they are stable under \(\Out(\bG)\). Steinberg shows that
    it is a closed embedding, and
    for fixed \(1\le i\le 2m\) and a weight \(\mu\) in \(V_i\), the trace of
    \(\pi_\mu \rho_i(\epsilon''(x))\pi_\mu\) is zero unless \(\mu=\Wt_i\) or
    \(\mu=\Wt_i'\), where \(\pi_\mu\colon V_i\to V_i\) is the projection to the
    weight space of weight \(\mu\), in which cases the weight multiplicities are
    both \(1\). Furthermore, letting \(x_0=x_{2m+1}=1\), we have
    \begin{align}
        \Tr(\pi_{\Wt_i}\rho_i(\epsilon''(x))\pi_{\Wt_i})&=\begin{cases}
            x_i & i\neq m-1, m\\
            x_\alpha & \text{otherwise}
        \end{cases}\\
        \Tr(\pi_{\Wt_i'}\rho_i(\epsilon''(x))\pi_{\Wt_i'})
        &=\begin{cases}
            x_{i-1} & 1\le i\le m\\
            x_{i+1} & m+1\le i\le 2m
        \end{cases}
    \end{align}

    The abelianization \(\bA_{\bM}\) is isomorphic to \(\bbA^{2m}\) with
    coordinates given by simple roots. The subtorus \(\bA_{\bM}^\x\) is identified
    with \(\bT^\AD\). The section \(\delta_{\bM}\) of
    abelianization map is induced by the anti-diagonal map
    \begin{align}
        \bT^\AD&\longto \bT_+\\
        a&\longmapsto (a,a^{-1})
    \end{align}
    Now let \(\tilde{\epsilon}''\colon \bA_{\bM}\x\bbA^{2m-1}\to\bM\) be the
    product \(\delta_{\bM}\epsilon''\). The image of \(\delta_{\bM}\) lies in
    \(\bar{\bT}_+\), thus preserves each weight spaces in \(V_i\), so
    \(\pi_\mu \rho_i(\tilde{\epsilon}''(a,x))\pi_\mu\) is still zero unless
    \(\mu=\Wt_i\) or \(\Wt_i'\), and
    \begin{align}
        \Tr(\pi_{\Wt_i}\rho_i(\tilde{\epsilon}''(a,x))\pi_{\Wt_i})&=\begin{cases}
            x_i & i\neq m-1, m\\
            x_\alpha & \text{otherwise}
        \end{cases}\\
        \Tr(\pi_{\Wt_i'}\rho_i(\tilde{\epsilon}''(a,x))\pi_{\Wt_i'})
        &=\begin{dcases}
            x_{i-1}\prod_{j=1}^{2m}a_j^{c_{ij}} & 1\le i\le m\\
            x_{i+1}\prod_{j=1}^{2m}a_j^{c_{ij}} & m+1\le i\le 2m
        \end{dcases}
    \end{align}
    where \(c_{ij}\) is such that \(\Wt_i-\Wt_i'=\sum_j c_{ij}\Rt_j\).
    Summarizing this part, we have that
    \begin{align}
        \chi_i(\tilde{\epsilon}''(a,x))=\begin{dcases}
            x_i+x_{i-1}\prod_{j=1}^{2m}a_j^{c_{ij}} & 1\le i\le m-1\\
            x_i+x_{i+1}\prod_{j=1}^{2m}a_j^{c_{ij}} & m+2\le i\le 2m\\
            x_\Rt+x_{m-1}\prod_{j=1}^{2m}a_j^{c_{ij}} & i=m\\
            x_\Rt+x_{m+2}\prod_{j=1}^{2m}a_j^{c_{ij}} & i=m+1
        \end{dcases}
    \end{align}

    On the other hand, consider another map
    \begin{align}
        \epsilon'''\colon \bbA^{2m-1}\x\Gm &\longto \bG\\
        (x,t)&\longmapsto
        u_mu_{m+1}\bU_\Rt(x_\Rt)\dot{s}_\Rt\dual{\Rt}(t)
        \prod_{i=1}^{m-1}(\bU_i(x_i)\dot{s}_i\bU_{2m-i+1}(x_{2m-i+1})\dot{s}_{2m-i+1}),
    \end{align}
    where \(1\neq u_m\in\bU_m(K)\) and \(1\neq u_{m+1}\in\bU_{m+1}(K)\) are two
    arbitrarily chosen elements. 
    Steinberg shows that this map is also a
    closed embedding. He also shows that if the commutator
    \([u_{m+1},u_m]\) is \(\bU_\Rt(1)\), \(\pi_\mu
    \rho_i(\epsilon''(x))\pi_\mu\) has trace \(0\) unless \(\mu=\Wt_i\) or
    \(\mu=\Wt_i'\), as well as
    \begin{align}
        \Tr(\pi_{\Wt_i}\rho_i(\epsilon'''(x,t))\pi_{\Wt_i})&=\begin{cases}
            x_i & i\neq m-1, m\\
            t x_\alpha & i=m\\
            t x_\Rt+t & i=m+1
        \end{cases}\\
        \Tr(\pi_{\Wt_i'}\rho_i(\epsilon'''(x,t))\pi_{\Wt_i'})
        &=\begin{cases}
            x_{i-1} & 1\le i\le m\\
            x_{i+1} & m+1\le i\le 2m
        \end{cases}
    \end{align}
    Let \(\tilde{\epsilon}'''=\delta_{\bM}\epsilon'''\), then we similarly have
    \begin{align}
        \Tr(\pi_{\Wt_i}\rho_i(\tilde{\epsilon}''(a,x,t))\pi_{\Wt_i})&=\begin{cases}
            t x_i & i\neq m-1, m\\
            t x_\alpha & i=m\\
            t x_\Rt+t & i=m+1
        \end{cases}\\
        \Tr(\pi_{\Wt_i'}\rho_i(\tilde{\epsilon}''(a,x,t))\pi_{\Wt_i'})
        &=\begin{dcases}
            x_{i-1}\prod_{j=1}^{2m}a_j^{c_{ij}} & 1\le i\le m\\
            x_{i+1}\prod_{j=1}^{2m}a_j^{c_{ij}} & m+1\le i\le 2m
        \end{dcases}
    \end{align}
    hence
    \begin{align}
        \chi_i(\tilde{\epsilon}'''(a,x,t))=\begin{dcases}
            t x_i+x_{i-1}\prod_{j=1}^{2m}a_j^{c_{ij}} & 1\le i\le m-1\\
            t x_i+x_{i+1}\prod_{j=1}^{2m}a_j^{c_{ij}} & m+2\le i\le 2m\\
            t x_\Rt+x_{m-1}\prod_{j=1}^{2m}a_j^{c_{ij}} & i=m\\
            t x_\Rt+t+x_{m+2}\prod_{j=1}^{2m}a_j^{c_{ij}} & i=m+1
        \end{dcases}
    \end{align}

    Let \(\tilde{\epsilon}'\) be the disjoint union of \(\tilde{\epsilon}''\)
    and \(\tilde{\epsilon}'''\), then it is not hard to see that
    \(\tilde{\epsilon}'\) is a bijection on points valued in \(K\)-fields between
    \(\bA_{\bM}\x(\bbA^{2m-1}\cup(\bbA^{2m-1}\x\Gm))\) and \(\bC_{\bM}\). Moreover,
    since \([u_{m},u_{m+1}]\in\bU_\Rt\), \(u_{m+1}u_m\bU_\Rt\)
    is stable under \(\Out(\bG)\). This means that the image of
    \(\epsilon'''\) is stable under \(\Out(\bG)\). Thus, there is an induced
    action of \(\Out(\bG)\) on \(\bA_{\bM}\x(\bbA^{2m-1}\cup(\bbA^{2m-1}\x\Gm))\)
    making \(\chi_{\bM}\circ\tilde{\epsilon}'\) an \(\Out(\bG)\)-equivariant map.
    Therefore, \(\tilde{\epsilon}'\) also induces bijection on \(K\)-rational
    points after any \(\Out(\bG)\)-twisting.

    If \(\bG\) is a product of groups \(\SL_{2m+1}\) for various \(m\) (allowing
    repetitions), then we simply take the direct product of
    \(\tilde{\epsilon}'\) of each simple factor in such way that if \(m_1=m_2\),
    then \(u_{m_1}=u_{m_2}\) and \(u_{m_1+1}=u_{m_2+1}\). In this way the image
    of \(\tilde{\epsilon}'\) will be stable under \(\Out(\bG)\). Hence, we are
    done.
\end{proof}

\section{Endoscopic Groups and Endoscopic Monoids}
\label{sec:endoscopic_groups_inv_theory}

In this section, we define one of the most important notions of this book being
the endoscopic monoid. We will also establish the transfer map at
invariant-theoretic level. As a byproduct of our effort, we can also define Levi
monoids in a similar way.

\subsection{}
Given the pinning
\(\bSPL=(\bT,\bB,\bFx_+)\), there is a pinning
\(\dbSPL=(\dual{\bT},\dual{\bB},\dual{\bFx}_+)\)
    \nomenclature[\(P_'check \)]{\(\dbSPL\)}{the pinning dual to \(\bSPL\)}
    \nomenclature[\(T"bold'check \)]{\(\dual{\bT}\)}{the dual torus in \(\dbSPL\)}
    \nomenclature[\(B"bold'check \)]{\(\dual{\bB}\)}{the Borel subgroup of \(\dual{\bG}\) in \(\dbSPL\)}
on the dual group \(\dual{\bG}\).
    \nomenclature[\(G"bold'check \)]{\(\dual{\bG}\)}{the dual group of \(G\) over \(\Qlb\) or \(\bbC\)}
Let \(\kappa\in\dual{\bT}\),
    \nomenclature[\(kappa \)]{\(\kappa\)}{the semisimple element in the dual group that is part of an endoscopic datum}
and \(\dual{\bH}\)
    \nomenclature[\(H"bold'check \)]{\(\dual{\bH}\)}{the connected centralizer of \(\kappa\) in \(\dual{\bG}\)}
be the \emph{connected} centralizer
of \(\kappa\) in \(\dual{\bG}\). Then \(\dual{\bH}\) has a maximal torus
\(\dual{\bT}\) and a Borel subgroup
\(\bB_{\dual{\bH}}=\dual{\bH}\cap\dual{\bB}\).
    \nomenclature[\(B"bold'check_H \)]{\(\bB_{\dual{\bH}}\)}{the Borel subgroup of \(\dual{\bH}\) induced by \(\dual{\bB}\)}
Taking the dual root datum of
\(\dual{\bH}\) determined by \((\dual{\bT},\bB_{\dual{\bH}})\), one obtains a
split \(k\)-group \(\bH\).
    \nomenclature[\(H"bold \)]{\(\bH\)}{the split dual group of \(\dual{\bH}\) over \(k\)}
Note that \(\Out(\bG)=\Out(\dual{\bG})\) and
similarly for \(\bH\) and \(\dual{\bH}\). We also have the root datum, the Weyl
group, and so on associated with \(\bH\), which we will denote by putting
\(\bH\) in the subscript.
    \nomenclature[\(.H,H"bold \)]{\((\cdot)_\bH, (\cdot)_H\)}{the analogue for \(\bH\) or \(H\) of constructions (such as root datum) related to \(\bG\) or \(G\)}

\subsection{}
Following \cite{Ng10}*{\S\S~1.8--1.9}, the centralizer
\((\dual{\bG}\rtimes\Out(\bG))_\kappa\) of \(\kappa\) in
\(\dual{\bG}\rtimes\Out(\bG)\) fits into a short exact sequence
\begin{align}
    1\longto \dual{\bH}\longto (\dual{\bG}\rtimes\Out(\bG))_\kappa\longto
    \pi_0(\kappa)\longto 1,
    \nomenclature[\(pi_0_kappa \)]{\(\pi_0(\kappa)\)}{the component group of the centralizer of \(\kappa\) in \(\dual{\bG}\rtimes\Out(\bG)\)}
\end{align}
where \(\pi_0(\kappa)\) is the component group of
\((\dual{\bG}\rtimes\Out(\bG))_\kappa\). Because \(\Out(\bG)\) is discrete, 
the projection \(\dual{\bG}\rtimes\Out(\bG)\to\Out(\bG)\) induces a canonical map
\begin{align}
    \bFo_{\bG}(\kappa)\colon \pi_0(\kappa)\longto \Out(\bG).
    \nomenclature[\(o_G_kappa \)]{\(\bFo_{\bG}(\kappa)\)}{the canonical map \(\pi_0(\kappa)\longto \Out(\bG)\)}
\end{align}
The action of \((\dual{\bG}\rtimes\Out(\bG))_\kappa\) on its normal subgroup
\(\dual{\bH}\) gives a homomorphism into \(\Aut(\dual{\bH})\), which further
maps into \(\Out(\dual{\bH})\). Hence, we have another canonical map
\begin{align}
    \bFo_{\bH}(\kappa)\colon \pi_0(\kappa)\longto \Out(\bH).
    \nomenclature[\(o_H_kappa \)]{\(\bFo_{\bH}(\kappa)\)}{the canonical map \(\pi_0(\kappa)\longto \Out(\bH)\)}
\end{align}

\begin{definition}
    Let \(G\) be a reductive group on \(X\) given by an \(\Out(\bG)\)-torsor
    \(\OGT_G\). An \notion{endoscopic datum}\index{endoscopic!datum} \((\kappa,\OGT_\kappa)\)
    \nomenclature[\(theta_kappa_var \)]{\(\OGT_\kappa\)}{the \(\pi_0(\kappa)\)-torsor in an endoscopic datum}
    of \(G\)
    is a pair where \(\kappa\in\dual{\bT}\), and \(\OGT_\kappa\) is a
    \(\pi_0(\kappa)\)-torsor such that the \(\Out(\bG)\)-torsor induced by it
    through \(\bFo_{\bG}(\kappa)\) is isomorphic to \(\OGT_G\). The
    \notion{endoscopic group}\index{endoscopic!group} \(H\)
    \nomenclature[\(H \)]{\(H\)}{an endoscopic group of \(G\)}
    associated with \((\kappa,\OGT_\kappa)\) is the
    twisted form of \(\bH\) induced by \(\Out(\bH)\)-torsor \(\OGT_H\), itself
    induced by \(\OGT_\kappa\) through \(\bFo_{\bH}(\kappa)\).
\end{definition}

\subsection{}
There is also a pointed variant using representations of \(\pi_1(X,\infty)\) after
fixing a geometric point \(\infty\in X\).
\begin{definition}
    A \notion{pointed endoscopic datum}\index{endoscopic!datum, pointed} of
    \((G,\infty_G)\) is a pair
    \((\kappa,\OGT_\kappa^\bullet)\), where \(\kappa\in\dual{\bT}\) and
    \(\OGT_\kappa^\bullet\)
    \nomenclature[\(theta_kappa_var_bullet \)]{\(\OGT_\kappa^\bullet\)}{the pointed version of \(\OGT_\kappa\)}
    is a continuous homomorphism
    \(\pi_1(X,\infty)\to\pi_0(\kappa)\) lying over \(\OGT_G^\bullet\).
    The \notion{pointed endoscopic group}\index{endoscopic!group, pointed}
    \((H,\infty_H)\)
    \nomenclature[\(infty_H \)]{\(\infty_H\)}{the canonical geometric point induced by \(\OGT_\kappa^\bullet\)}
    is the pointed twisted
    form induced by \(\OGT_\kappa^\bullet\) through \(\bFo_{\bH}(\kappa)\).
\end{definition}

\subsection{}
Experienced readers may have noticed that we missed a crucial ingredient in the
definition of endoscopic data, namely the admissible embeddings (or
\(L\)-embeddings) of relevant \(L\)-groups. Such
information is not present in \cite{Ng10} because it turns out to be irrelevant
in the Lie algebra case
if one exploits the properties of Kostant sections (see \cite{Ko99} for
example).

In group case, we are not as fortunate and
\(L\)-embeddings (even the existence of such) is one of the main reasons
why endoscopic theories are so complicated. However, since in our considerations
both \(G\) and \(H\) are quasi-split, the relevance of \(L\)-embeddings
does not appear until much later, and for the majority of this book it is not
necessary. For this reason we will only bring up
\(L\)-embeddings when it becomes relevant. The reader can refer to
\Cref{sec:Appendix_endoscopic_groups} for a more complete definition of
endoscopic data.

\subsection{}
Through \(\bFo_{\bG}(\kappa)\) and \(\bFo_{\bH}(\kappa)\) respectively, we can
define actions of \(\bW\rtimes\pi_0(\kappa)\) and
\(\bW_{\bH}\rtimes\pi_0(\kappa)\) on \(\bT\). The following lemma is crucial to
settling the compatibility questions regarding these two actions. 
\begin{lemma}[\cite{Ng10}*{Lemme~1.9.1}]
    \label[lemma]{lem:kappa_G_H_compatibility}
    There is a canonical homomorphism
    \begin{align}
        \bW_{\bH}\rtimes \pi_0(\kappa)\longto \bW\rtimes \pi_0(\kappa)
    \end{align}
    whose restriction to \(\bW_{\bH}\) is the inclusion \(\bW_{\bH}\subset\bW\),
    and it induces identity on quotient group \(\pi_0(\kappa)\). Moreover, such
    homomorphism is compatible with the actions of the two groups on \(\bT\).
\end{lemma}

\begin{remark}
    Note that on contrary the \(\pi_0(\kappa)\)-actions on \(\bT\) induced
    respectively by \(\bFo_{\bG}(\kappa)\) and \(\bFo_{\bH}(\kappa)\) are not
    compatible in general.
\end{remark}

\subsection{}
Let \(\bM\in\FM(\bG^\SC)\) be a very flat monoid with maximal toric variety
\(\bar{\bT}_{\bM}\). Then the root datum of \(\bH\) induces a group \(\bH_{\bM}\)
with maximal torus \(\bT_{\bM}\) and \(\bH_{\bM}^\SC=\bH^\SC\).
According to \Cref{thm:monoid_classification_over_k_bar} (also see
\Cref{rmk:monoid_classification_over_base_field}), the diagram
\begin{align}
    \bar{\bT}_{\bM}\supset\bT_{\bM}\subset \bH_{\bM}
\end{align}
defines up to isomorphism a unique monoid \(\bM_{\bH}'\)
    \nomenclature[\(M_H"bold^prime \)]{\(\bM_\bH'\)}{the (usually non-flat) monoid of \(\bH\) with maximal toric variety \(\bar{\bT}_\bM\)}
with unit group
\(\bH_{\bM}\). However, \(\bM_{\bH}'\) is usually not very flat (by looking at
the cone of its maximal toric variety; cf.~\Cref{sub:cones_of_maximal_toric_variety}), and
\(\bH_{\bM}^\Der\) is not necessarily simply-connected. 
Therefore, we need to find a very flat monoid in \(\FM(\bH^\SC)\) that remedies
this problem. We shall see
that it has to do with the fact that an irreducible representation of
\(\dual{\bG}\) is no longer irreducible when restricted to \(\dual{\bH}\) but
decomposes into a direct sum of irreducible ones. The
correct monoid for \(\bH\) is exactly guided by this decomposition.

\subsection{}
Suppose \(G\) is obtained from \(\OGT_G^\bullet\). To simplify notations we
let \(\Theta=\pi_1(X,\infty)\) and \(\TD{G}=\dual{\bG}\rtimes\Theta\).
Let \(V\) be a finite-dimensional representation of \(\TD{G}\)
over \(\Qlb\). It canonically decomposes into irreducible
\(\dual{\bG}\)-representations as follows:
\begin{align}
    V\simeq\bigoplus_{\theta}V_\theta\otimes\Hom_{\dual{\bG}}(V_\theta,V),
\end{align}
where \(\theta\) ranges over \(\dual{\bG}\)-highest weights and \(V_\theta\) is
the corresponding irreducible representation. Let
\begin{align}
    -w_0(\theta_1),\ldots,-w_0(\theta_m)
\end{align}
be those \(\theta\) with non-zero multiplicities in \(V\). Then the monodromy
\(\Theta\) acts on the set of \(\theta_i\) through \(\Out(\bG)\). The
multiplicity spaces \(\Hom_{\dual{\bG}}(V_\theta,V)\) will play an important
role later, but we ignore them for now.

Let \(\bA_{\bM}\) be the affine space whose cone is freely
generated by cocharacters \(\theta_1,\ldots,\theta_m\). The map \(G\to G^\AD\)
induces canonical morphism of toric varieties
\(\bA_{\bM}\to\bA_{\Env(\bG^\SC)}\) by sending each \(\theta_i\) to its image
in the adjoint cocharacter lattice, and let \(\bM\) be the monoid whose
abelianization is \(\bA_\bM\).
The elements in the cone of \(\bar{\bT}_{\bM}\) are then of the form
\begin{align}
    (c_1\theta_1,\ldots,c_m\theta_m,\mu_\AD)
\end{align}
for \(c_i\in\bbN\) and a unique \(\mu\in \WP_{\dual{\bG}}[ -w_0(\sum_{i=1}^m
c_i\theta_i) ]\) whose image in the adjoint group is \(\mu_\AD\).

\subsection{}
Suppose furthermore that \(\OGT_G^\bullet\) is induced by
\(\OGT_\kappa^\bullet\) and \(H\) is the endoscopic group. We denote similarly
\(\TD{H}=\dual{\bH}\rtimes\Theta\).
The restriction of representation \(V_{-w_0(\theta_i)}\) to \(\dual{\bH}\)
decomposes as follows:
\begin{align}
    V_{-w_0(\theta_i)}\simeq
    \bigoplus_{j=1}^{e_i}V_{-w_{\bH,0}(\lambda_{ij})}^H
    \otimes\Hom_{\dual{\bH}}(V_{-w_{\bH,0}(\lambda_{ij})}^H,V_{-w_0(\theta_i)}),
\end{align}
where \(w_{\bH,0}\) is the longest element of \(\bW_{\bH}\) with respect to
\(\bB_{\bH}\), \(-w_{\bH,0}(\lambda_{ij})\) are distinct \(\dual{\bH}\)-highest
weights and \(V_{-w_{\bH,0}(\lambda_{ij})}^H\) are the corresponding irreducible
\(\dual{\bH}\)-representations with non-zero multiplicities.
Although there is no canonical \(\TD{H}\)-action on
\(V\), we still have an action on the set of all \(\lambda_{ij}\) and this
action factors through \(\Theta\). Moreover,
the assignment \(\lambda_{ij}\mapsto\theta_i\) is \(\Theta\)-equivariant
thanks to \Cref{lem:kappa_G_H_compatibility}.
Let \(\bM_{\bH}\in\FM_0(\bH^\SC)\)
    \nomenclature[\(M_H"bold \)]{\(\bM_\bH\)}{the split form of the endoscopic monoid associated with \(\FRM\in\FM(G^\SC)\)}
be the monoid whose
abelianization \(\bA_{\bM,\bH}\)
    \nomenclature[\(.M_H"bold \)]{\((\cdot)_{\bM,\bH},(\cdot)_{\FRM,H}\)}{constructions related to \(\bM_\bH\) or \(\FRM_H\)}
is the affine space with cone freely generated by all
\(\lambda_{ij}\). The cone of the maximal toric variety \(\bar{\bT}_{\bM,\bH}\)
consists of elements of the following form:
\begin{align}
    (a_{ij}\lambda_{ij},\mu_{H,\AD})
\end{align}
where \(a_{ij}\in\bbN\), \(\mu_H\in\WP_{\dual{\bH}}[ -w_{\bH,0}(\sum_{i,j}a_{ij}\lambda_{ij}) ]\) and
\(\mu_{H,\AD}\) its image under map \(H\to H^\AD\).

For technical reasons, we also define another closely related monoid \(\bM_\bH^\star\)
    \nomenclature[\(M_H"bold^star \)]{\(\bM_\bH^\star\)}{the split form of the coarse endoscopic monoid}
as follows: for each individual \(i\), we replace the set of
\(\lambda_{ij}\) with its \(\Theta\)-stable subset
consisting of elements that are maximal with respect to the \(\bH\)-dominance
order. In other words, if \(j\neq j'\), then neither
\(\lambda_{ij}\le_\bH\lambda_{ij'}\) nor \(\lambda_{ij'}\le_\bH\lambda_{ij}\)
is true. The monoid \(\bM_\bH^\star\) is the monoid corresponding to this subset
of \(\lambda_{ij}\). Clearly we have a canonical closed embedding \(\bM_\bH^\star\to
\bM_\bH\) that is the pullback of map of abelianizations
\(\bA_{\bM,\bH}^\star\to\bA_{\bM,\bH}\).
    \nomenclature[\(.^star \)]{\((\cdot)^\star\)}{constructions related to coarse endoscopic or Levi monoids}

\subsection{}
We now construct a \(\bW_{\bH}\)-equivariant homomorphism
\begin{align}
    \tilde{\nu}_{\bH}\colon\bar{\bT}_{\bM,\bH}\longto\bar{\bT}_{\bM},
    \nomenclature[\(nu'tilde_H \)]{\(\tilde{\nu}_{\bH},\tilde{\nu}_H\)}{the
    canonical map \(\bar{\bT}_{\bM,\bH}\to\bar{\bT}_{\bM}\) or \(\FRT_{\FRM,H}\to\FRT_{\FRM}\)}
\end{align}
and as a consequence we will 
obtain by \Cref{rmk:monoid_classification_over_base_field} 
a homomorphism of monoids \(\bM_{\bH}\to\bM_{\bH}'\).
Indeed, the weight \((a_{ij}\lambda_{ij},\mu_{H,\AD})\) can be uniquely
written as
\begin{align}
    \left(a_{ij}\lambda_{ij},-w_{\bH,0}\Bigl(\sum_{i,j}a_{ij}\lambda_{ij}\Bigr)
    -\sum_{\CoRt\in\SimCoRts_{\bH}}h_{\CoRt}^\mu\CoRt\right),
\end{align}
where \(h_{\CoRt}^\mu\in\bbN\). We send this element to
\begin{align}
    \left(\sum_{j=1}^{e_1} a_{1j}\theta_1,\ldots,\sum_{j=1}^{e_m}a_{mj}\theta_m,
        -w_{\bH,0}\Bigl(\sum_{i,j}a_{ij}\lambda_{ij}\Bigr)
    -\sum_{\CoRt\in\SimCoRts_{\bH}}h_{\CoRt}^\mu\CoRt\right).
\end{align}
One may check that it is indeed a \(\bW_{\bH}\)-equivariant homomorphism.
By composition, we also have a \(\bW_{\bH}\)-equivariant homomorphism
\begin{align}
    \tilde{\nu}_{\bH}^\star\colon\bar{\bT}_{\bM,\bH}^\star\longto\bar{\bT}_{\bM}.
\end{align}

\begin{lemma}
    \label[lemma]{lem:surjectivity_of_cones_from_TMH_to_TM}
    The map of cones of toric varieties \(\bar{\bT}_{\bM,\bH}\) and
    \(\bar{\bT}_\bM\) induced by \(\tilde{\nu}_\bH\) is surjective.
    Consequently, the homomorphism \(\bT_{\bM,\bH}\to\bT_\bM\) is surjective
    with connected kernel. The same is true if \(\bar{\bT}_{\bM,\bH}\)
    (resp.~\(\bT_{\bM,\bH}\)) is replaced by \(\bar{\bT}_{\bM,\bH}^\star\)
    (resp.~\(\bT_{\bM,\bH}^\star\)).
\end{lemma}
\begin{proof}
    This is because for each \(i\) the weights inside
    \(V_{-w_{\bH,0}(\lambda_{ij})}^H\) (\(1\le j\le e_i\)) account for all
    weights in \(V_{-w_0(\theta_i)}\), and the same remains true if we only consider
    \(\bH\)-maximal \(\lambda_{ij}\)'s.
\end{proof}

\subsection{}
\label{sub:Z_kappa_def}
We have \(\bW_{\bH}\)-homomorphisms
\begin{align}
    \bT_{\bM,\bH}\longto\bT_{\bM}\longto\bT^\AD\longto \bT_{\bH^\AD}.
\end{align}
It implies that the diagonalizable group (but not necessarily a torus)
\begin{align}
    \bZ_{\bM}^\kappa\defeq 
    \tilde{\nu}_{\bH}^{-1}(\bZ_{\bM})
    \quad\text{(resp. }\bZ_{\bM}^{\kappa\star}\defeq
    \bZ_{\bM}^\kappa\cap\bM_\bH^\star\text{)}
    \nomenclature[\(Z_M_kappa \)]{\(\bZ_{\bM}^\kappa,Z_\FRM^\kappa\)}{the preimage of \(\bZ_\bM\) (resp.~\(Z_\FRM\)) in \(\bT_{\bM,\bH}\) (resp.~\(T_{\FRM,H}\))}
\end{align}
is contained in the center of \(\bM_{\bH}\) (resp.~\(\bM_\bH^\star\)). As a
result, the maps
\begin{align}
    \Stack*{\bar{\bT}_{\bM,\bH}^\star/\bZ_{\bM}^{\kappa\star}}\longto
    \Stack*{\bar{\bT}_{\bM,\bH}/\bZ_{\bM}^\kappa}\longto\Stack*{\bar{\bT}_{\bM}/\bZ_{\bM}}
\end{align}
are generically isomorphisms. In fact, we can improve it to the following
statement:
\begin{lemma}
    \label[lemma]{lem:arc_lifting_to_endoscopic_monoid}
    Let \(\cO=\bar{k}\powser{\pi}\) and \(F=\bar{k}\lauser{\pi}\). Then the map
    \begin{align}
        \bar{\bT}_{\bM,\bH}(\cO)\cap\bT_{\bM,\bH}(F)\longto\bar{\bT}_{\bM}(\cO)\cap\bT_{\bM}(F)
    \end{align}
    is surjective. Moreover, suppose
    \(t\in\bT_{\bM,\bH}/\bZ_{\bM}^\kappa(F)\) extends to a point
    \(t_\cO\in\Stack*{\bar{\bT}_{\bM}/\bZ_{\bM}}(\cO)\), then there exists at least
    one and finitely many ways to extend \(t\) to a point
    in \(\Stack*{\bar{\bT}_{\bM,\bH}/\bZ_{\bM}^\kappa}(\cO)\) lying over
    \(t_\cO\). The same remain true if we replace everything by its respective
    \(\star\)-version.
\end{lemma}
\begin{proof}
    The first claim follows from
    \Cref{lem:surjectivity_of_cones_from_TMH_to_TM}.
    For the second claim, since \(\cO\) has algebraically closed residue
    field, \(\bZ_{\bM}\)-torsors over \(\cO\) are trivial, and so we can lift
    \(t_\cO\) to a point \(t'_\cO\in\bar{\bT}_{\bM}(\cO)\cap\bT_{\bM}(F)\).
    Using the first claim just proved, let
    \(\tilde{t}\) be the lift of \(t'_\cO\) in
    \(\bar{\bT}_{\bM,\bH}(\cO)\cap\bT_{\bM,\bH}(F)\), then its image in
    \(\Stack*{\bar{\bT}_{\bM,\bH}/\bZ_{\bM}^\kappa}(\cO)\) is the desired
    extension of \(t\).

    Since \(\bZ_{\bM}^\kappa\)-torsors over \(\cO\) are also trivial, the
    isomorphism class of any lift of \(t_\cO\) is uniquely determined by the
    lift of cocharacter \((c_i\theta_i,\mu)\) to
    \((a_{ij}\lambda_{ij},\mu')\).
    For a fixed \(t'\), \(c_i\) and \(\mu\) are fixed, so the set of possible
    \(a_{ij}\) is finite, and \(\mu'\) is uniquely
    determined by \(\mu\). This proves finiteness of possible lifts. The same
    proof works verbatim for the \(\star\)-version.
\end{proof}

\subsection{}
Since \(\bZ_\bM^\kappa\) is usually not the whole center \(\bZ_{\bM,\bH}\), it
is important to check that it is still big enough so that its action on
\(\bA_{\bM,\bH}\) has good properties. For this purpose, let \(a_{ij}\) be the
cardinality of the \(\bW_\bH\)-orbit of \(\lambda_{ij}\), then
it is easy to check that after projecting to \(\bT^\AD\), we have
\begin{align}
    \label{eqn:a_ij_for_Z_M_H}
    0\le_{\bH}
    -w_{\bH,0}\Bigl(\sum_{ij}a_{ij}\lambda_{ij}\Bigr)_\AD\in\CoCharG(\bT^\AD).
\end{align}
This shows that \((a_{ij}\lambda_{ij},0)\in\CoCharG(\bT_{\bM,\bH})\) is a
cocharacter of \(\bZ_\bM^\kappa\), and its projection to each \(\bbA^1\)-factor
of \(\bA_{\bM,\bH}\)
corresponding to \(\lambda_{ij}\) is dominant. It implies in particular that the
\(\bZ_\bM^\kappa\)-action on \(\bA_{\bM,\bH}\) has a unique closed orbit being
the origin. The same is again true for \(\bZ_\bM^{\kappa\star}\) acting on
\(\bA_{\bM,\bH}^\star\).

\begin{lemma}
    The map \(\tilde{\nu}_{\bH}\) induces maps of abelianizations of \(\bM_{\bH}\),
    \(\bM_{\bH}'\) and \(\bM\):
    \begin{align}
        \bA_{\bM,\bH}\longto \bA_{\bM,\bH}'\longto\bA_{\bM}.
    \end{align}
    The preimage of invertible locus \(\bA_{\bM}^\x\) in \(\bA_{\bM,\bH}\) is
    precisely \(\bA_{\bM,\bH}^{\x}\). Similarly, its preimage in
    \(\bA_{\bM,\bH}^{\star}\) is \(\bA_{\bM,\bH}^{\star\x}\).
\end{lemma}
\begin{proof}
    The first map is induced by the universal property of taking
    respective GIT quotients of \(\bM_\bH\) and \(\bM_\bH'\) by
    \(\bH^\SC\x\bH^\SC\) (recall we have a map \(\bM_\bH\to\bM_\bH'\) above).
    The second one can be seen as follows:
    the ring of regular functions \(k[\bA_{\bM,\bH}']\) is generated by
    characters in \(k[\bar{\bT}_{\bM}]\) perpendicular to coroots in \(\bH\),
    while \(k[\bA_{\bM}]\) is one generated by characters perpendicular to coroots
    in \(\bG\), and \(\CoRoots_{\bH}\subset\CoRoots\).
    The composition \(\bA_{\bM,\bH}\to\bA_{\bM}\) can be described at
    cocharacter level by \((a_{ij}\lambda_{ij})\mapsto (c_i\theta_i)\), where
    \(c_i=\sum_{j=1}^{e_i}a_{ij}\). So it is readily seen that the preimage of
    \(\bA_{\bM}^\x\) is \(\bA_{\bM,\bH}^{\x}\). The result is then trivial for
    \(\bA_{\bM,\bH}^\star\).
\end{proof}

\subsection{}
We have \(\bW\rtimes\Theta\) acting on \(\bar{\bT}_{\bM}\) by combining the
actions of \(\bW\) and \(\Theta\), hence also an action of
\(\bW_{\bH}\rtimes\Theta\) on \(\bar{\bT}_{\bM}\) induced by the canonical
map in \Cref{lem:kappa_G_H_compatibility}.

\begin{lemma}
    \label[lemma]{lem:kappa_compatibility_on_monoid}
    There is a canonical action of \(\bW_{\bH}\rtimes\Theta\) on
    \(\bar{\bT}_{\bM,\bH}\) (resp.~\(\bar{\bT}_{\bM,\bH}^\star\)) making
    \(\tilde{\nu}_{\bH}\) (resp.~\(\tilde{\nu}_{\bH}^\star\)) a
    \(\bW_{\bH}\rtimes\Theta\)-equivariant map.
\end{lemma}
\begin{proof}
    The same proof works for both \(\tilde{\nu}_{\bH}\) and
    \(\tilde{\nu}_{\bH}^\star\) so we only consider the former.
    To avoid confusion, we write \(\bW\rtimes\Theta\) as
    \(\bW\rtimes_{\bG}\Theta\) and \(\bW_{\bH}\rtimes\Theta\) as
    \(\bW_{\bH}\rtimes_{\bH}\Theta\).
    It suffices to prove compatibility on subgroup
    \(1\rtimes_{\bH}\Theta\). Let \(\sigma\in 1\rtimes_{\bH}\Theta\)
    be an element.

    The action of \(\sigma\) on \(\bT\) through
    \(\bW\rtimes_{\bG}\Theta\) is compatible with its natural action on
    \(\bT_{\bH^\SC}\) through \(\bW_{\bH}\rtimes\Out(\bH)\). In particular,
    \(w_{\bH,0}\) is fixed by \(\sigma\). The
    action of \(\sigma\) on \(\CoCharG(\bT)\)
    preserves subset \(\Set{\lambda_{ij}}\) because as an element of
    \(1\rtimes_\bH\Theta\) it must send one \(\dual{\bH}\)-highest weight to
    another. By the same reason \(\sigma\) also stabilizes the simple coroots of
    \(\bH\). 
    Then the result follows
    from the definition of map \(\tilde{\nu}_{\bH}\).
\end{proof}

Combining \Cref{lem:kappa_compatibility_on_monoid}
with the action of \(\Out(\bH)\) on \(\Env(\bH^\SC)\), we
obtain a canonical action of \(1\rtimes_{\bH}\Theta\) on \(\bM_{\bH}\).
As a result, we obtain morphisms of monoids over \(X\) by taking
\(\OGT_\kappa^\bullet\)-twists:
\begin{align}
    \FRM_H^\star\longto\FRM_H\longto \Env(H^\SC).
\end{align}
Note, however, that there is no action of \(\Out(\bH)\) on
\(\bar{\bT}_{\bM,\bH}\) or \(\bM_{\bH}\).
\begin{definition}
    \label{eqn:def_endoscopic_monoid}
    The monoid \(\FRM_H\)
    \nomenclature[\(M"frak_H \)]{\(\FRM_H\)}{the endoscopic monoid associated with \(\FRM\in\FM(G^\SC)\) and endoscopic group \(H\)}
    (resp.~\(\FRM_H^\star\)) is called the
    \notion{endoscopic monoid}\index{monoid!endoscopic}\index{endoscopic!monoid}
    (resp.~\notion{coarse endoscopic monoid}\index{monoid!coarse
    endoscopic}\index{endoscopic!monoid, coarse}) associated
    with monoid \(\FRM\) and endoscopic group \(H\).
\end{definition}

The group \(\Theta\) acts canonically on \(\bZ_{\bM}^\kappa\) and
\(\bZ_{\bM,\bH}\) because
\(\bW_{\bH}\) and \(\bW\) act trivially on central tori, and the map
in \Cref{lem:kappa_G_H_compatibility} induces identity on
\(\pi_0(\kappa)\) as quotient groups.
The map \(\bZ_{\bM}^\kappa\to \bZ_{\bM}\) is \(\Theta\)-equivariant.

\subsection{}
We have a commutative diagram of invariant quotients
\begin{equation}
    \begin{tikzcd}
        & \bar{\bT}_{\bM,\bH}\git\bW_{\bH} \ar[r, "\sim"]\ar[d] & \bM_{\bH}\git\bH \ar[d]\\
        \bar{\bT}_{\bM}\git\bW & \bar{\bT}_{\bM}/\bW_{\bH} \ar[l]\ar[r, "\sim"] & \bM_{\bH}'\git\bH
    \end{tikzcd}
\end{equation}
which induces canonical map
\begin{align}
    \label{eqn:absolute_transfer_map_def}
    \nu_{\bH}\colon\bC_{\bM,\bH}\defeq\bM_{\bH}\git\bH\longto\bC_{\bM}.
    \nomenclature[\(nu_H \)]{\(\nu_{\bH},\nu_H\)}{the canonical map \(\bC_{\bM,\bH}\to\bC_{\bM}\) or \(\FRC_{\FRM,H}\to\FRC_{\FRM}\)}
\end{align}
Since the preimage of \(\bA_{\bM}^\x\) in \(\bA_{\bM,\bH}\) is
\(\bA_{\bM,\bH}^\x\), we have \(\nu_{\bH}^{-1}(\bC_{\bM}^\x)=\bC_{\bM,\bH}^\x\).
We also have the coarse version \(\bC_{\bM,\bH}^\star\).
For convenience, we let \(\bC_{\bM,\bH}'\) be the quotient
\(\bar{\bT}_{\bM}\git\bW_{\bH}\).
\begin{lemma}
    \label[lemma]{lem:arc_lifting_to_endoscopic_quotient}
    The restriction of maps
    \begin{align}
        \Stack*{\bC_{\bM,\bH}/\bZ_{\bM}^\kappa}\longto \Stack*{\bC_{\bM,\bH}'/\bZ_{\bM}}
        \longto \Stack*{\bC_{\bM}/\bZ_{\bM}}
    \end{align}
    to \(\Stack*{\bC_{\bM}^{\x,\rss}/\bZ_{\bM}}\)
    are \'etale, and the first one is even an isomorphism over this subset.
    Moreover, let \(\cO=\bar{k}\powser{\pi}\) and \(F=\bar{k}\lauser{\pi}\), and suppose
    \(a_{\bH}\in\Stack*{\bC_{\bM,\bH}^\x/\bZ_{\bM}^\kappa}(F)\) extends to a point
    \(a\in\Stack*{\bC_{\bM}/\bZ_{\bM}}(\cO)\), then there exists at least one and 
    finitely many ways to extend \(a_{\bH}\) to a point
    in \(\Stack*{\bC_{\bM,\bH}/\bZ_{\bM}^\kappa}(\cO)\) lying over \(a\). The
    same is true if we replace \(\Stack*{\bC_{\bM,\bH}/\bZ_{\bM}^\kappa}\) by
    \(\Stack*{\bC_{\bM,\bH}^\star/\bZ_{\bM}^{\kappa\star}}\).
\end{lemma}
\begin{proof}
    The first claim is clear from definition. The second claim is proved using
    \Cref{lem:arc_lifting_to_endoscopic_monoid}. Indeed, since
    \(\bar{\bT}_{\bM,\bH}\to\bC_{\bM,\bH}\) is finite, we may lift \(a_{\bH}\)
    to a point \(t_{\bH}\in\Stack*{\bar{\bT}_{\bM}/\bZ_{\bM}}(F_l)\) for some finite
    tamely ramified extension \(F_l\) of \(F\), and there are only finitely many
    non-isomorphic ways to do so if we require that \(F_l\) is chosen to be as
    small as possible.

    Since \(\bar{\bT}_{\bM}\to\bC_{\bM}\) is also finite, the image of \(t_{\bH}\) in
    \(\Stack*{\bar{\bT}_{\bM}/\bZ_{\bM}}\) extends to a unique \(\cO_l'\)-point \(t\)
    lying over \(a\) by valuative criteria for properness, where \(\cO_l'\) is
    a finite extension of \(\cO_l\). By our assumption on \(\Char(k)\) relative
    to the rank of \(\bG\), any \(\bZ_{\bM}\)-torsor over \(F_l\) is
    trivializable over a tamely ramified extension, so \(\cO_l'\) can be chosen
    to be tamely ramified over \(\cO_l\), hence over \(\cO\). Thus, we may
    replace \(\cO_l\) by \(\cO_l'\) and then \(t_{\bH}\) extends to a point in
    \(\Stack*{\bar{\bT}_{\bM}/\bZ_{\bM}}(\cO_l)\).

    \Cref{lem:arc_lifting_to_endoscopic_monoid}
    shows that there is at least one and at most finitely many ways to
    extend \(t_{\bH}\) to an \(\cO_l\)-point in
    \(\Stack*{\bar{\bT}_{\bM,\bH}/\bZ_{\bM}^\kappa}\). The image of such
    extension in \(\Stack*{\bC_{\bM,\bH}/\bZ_{\bM}^\kappa}\) is an
    \(\cO_l\)-point over \(a\) extending \(a_{\bH}\). But it is also an
    \(F\)-point, thus it must be an \(\cO\)-point. Moreover, any extension of
    \(a_{\bH}\) to an \(\cO\)-point can be obtained in this way: indeed, any such
    \(\cO\)-point lifts to some \(\cO_l\)-point in
    \(\Stack*{\bar{\bT}_{\bM,\bH}/\bZ_{\bM}^\kappa}\) (using the fact that
    \(\bZ_{\bM}^\kappa\)-torsors over \(\cO\) are trivial and valuative criteria
    for finite map of varieties \(\bar{\bT}_{\bM,\bH}\to\bC_{\bM,\bH}\)).

    Therefore, the set of extensions of
    \(a_{\bH}\) to \(\cO\)-points is necessarily finite because each step above
    yields finitely many possibilities and every such extension of
    \(a_{\bH}\) can be obtained in this way. The proof for the
    \(\star\)-version is the same.
\end{proof}

\subsection{}
Let \(\bC_{\bM,\bH}^{\bG\hy\rss}\)
    \nomenclature[\(.G_rs \)]{\((\cdot)^{\bG\hy\rss},(\cdot)^{G\hy\rss}\)}{the (restriction to the) \(\bG\)- or \(G\)-regular semisimple locus}
be the preimage of \(\bC_{\bM}^\rss\) under
\(\nu_{\bH}\). Recall the extended discriminant functions:
\begin{align}
    {\Disc_+}&= e^{(2\rho,0)}\prod_{\alpha\in\Roots}(1-e^{(0,\alpha)})\\
    {\Disc_{\bH,+}}&=e^{(2\rho_{\bH},0)}\prod_{\alpha\in\Roots_{\bH}}(1-e^{(0,\alpha)})
\end{align}
and note that \(\Roots_{\bH}\subset\Roots\) and \(\Roots_{\bH,+}\subset\PosRts\).
Over the invertible locus, \(\Disc_+/\Disc_{\bH,+}\) is a regular function,
hence \(\nu_{\bH}(\bD_{\bM,\bH}^\x)\subset\bD_{\bM}^\x\). Since \(\bD_{\bM,\bH}^\x\)
is dense in \(\bD_{\bM,\bH}\) and \(\bD_{\bM,\bH}\) is reduced, we see that
\(\nu_{\bH}(\bD_{\bM,\bH})\subset\bD_{\bM}\) in the scheme-theoretic sense.
As a result, \(\bC_{\bM,\bH}^{\bG\hy\rss}\subset\bC_{\bM,\bH}^\rss\) and
\begin{align}
    \Disc_+/\Disc_{\bH,+}=e^{(2\rho-2\rho_{\bH},0)}\prod_{\Rt\in\Roots-\Roots_{\bH}}(1-e^{(0,\Rt)})
\end{align}
defines a \(\bW_{\bH}\)-invariant function on \(\bar{\bT}_{\bM,\bH}\).
The following is an analogue of \cite{Ng10}*{Lemme~1.10.2}:
\begin{lemma}
    Let \(\Psi\subset \Roots-\Roots_{\bH}\) be such that for any pair of
    roots \(\Set{\pm\Rt}\subset \Roots-\Roots_\bH\), the intersection
    \(\Set{\pm\Rt}\cap\Psi\) contains exactly one element. Let \(\rho_\Psi\) be
    the half-sum of elements in \(\Psi\). Then
    the formula
    \begin{align}
        r_\bH^\bG\defeq e^{(\rho-\rho_G,-\rho_\Psi)}\prod_{\Rt\in\Psi}(1-e^{(0,\Rt)})
        \nomenclature[\(r_H^G \)]{\(r_\bH^\bG,r_H^G\)}{the extended resultant function for \(\bG\) (resp.~\(G\)) and \(\bH\) (resp.~\(H\))}
    \end{align}
    defines a \(\bW_\bH\)-invariant function on \(\bar{\bT}_{\bM,\bH}\) such
    that
    \begin{align}
        \label{eqn:equality_between_resultant_and_two_discriminants}
        (r_\bH^\bG)^2(-1)^{\abs*{\Roots_{\bH,+}}}\Disc_{\bH,+}=(-1)^{\abs*{\PosRts}}\Disc_+.
    \end{align}
\end{lemma}
\begin{proof}
    As a rational function on \(\bT_{\bM,\bH}\), \(r_\bH^\bG\) clearly satisfies
    \eqref{eqn:equality_between_resultant_and_two_discriminants} by definition.
    Since \(\bar{\bT}_{\bM,\bH}\) is normal and \(\Disc_+/\Disc_{\bH,+}\) is a
    regular function, \(r_\bH^\bG\) must also be regular on
    \(\bar{\bT}_{\bM,\bH}\) because it is a rational function that satisfies a
    monic polynomial over the integrally closed ring \(k[\bar{\bT}_{\bM,\bH}]\).
    It remains to prove \(\bW_\bH\)-invariance of \(r_\bH^\bG\). The proof is
    similar to that of \cite{Ng10}*{Lemme~1.10.2}.

    Choose a subset \(\Psi_\bH\subset \Roots_\bH\) such that for any
    pair of roots \(\Set{\pm\Rt}\subset \Roots_\bH\) the intersection
    \(\Psi_\bH\cap \Set{\pm\Rt}\) contains exactly one element, and let
    \(\Psi_\bG=\Psi_\bH\cup \Psi\). Let \(\rho_{\Psi,\bH}\)
    (resp.~\(\rho_{\Psi,\bG}\)) be the half-sum of elements in \(\Psi_\bH\)
    (resp.~\(\Psi_\bG\)), then \(\rho_{\Psi,\bG}=\rho_\Psi+\rho_{\Psi,\bH}\).
    For any \(w\in\bW\), it is straightforward to see that
    \begin{align}
        w\Bigl(e^{(\rho,-\rho_{\Psi,\bG})}\prod_{\Rt\in\Psi_\bG}(1-e^{(0,\Rt)})\Bigr)
        =(-1)^{\abs*{w(\Psi_\bG)-\Psi_\bG}}e^{(\rho,-\rho_{\Psi,\bG})}\prod_{\Rt\in\Psi_\bG}(1-e^{(0,\Rt)}).
    \end{align}
    We claim that the sign
    \begin{align}
        (-1)^{\abs*{w(\Psi_\bG)-\Psi_\bG}}
    \end{align}
    is independent of \(\Psi_\bG\). Indeed, suppose \(\Psi_\bG'\subset \Roots\)
    is obtained from \(\Psi_\bG\) by replacing a single element
    \(\Rt\in\Psi_\bG\) with \(-\Rt\). If \(w(\Rt)\in\Set{\pm\Rt}\), then
    \(w(\Psi_\bG)-\Psi_\bG\) and
    \(w(\Psi_\bG')-\Psi_\bG'\) have the same number of elements, 
    otherwise the difference in numbers of
    elements is either \(0\) or \(2\). Therefore, we always have
    \begin{align}
        (-1)^{\abs*{w(\Psi_\bG)-\Psi_\bG}}=(-1)^{\abs*{w(\Psi_\bG')-\Psi_\bG'}},
    \end{align}
    and the claim follows by induction. A particular choice of \(\Psi_\bG\) is
    \(\PosRts\), and in this case we have
    \begin{align}
        \abs*{w(\PosRts)-\PosRts}=\ell_\bG(w)=\dim(\bB w\bB/\bB),
    \end{align}
    where \(\ell_\bG(w)\) is the length of \(w\) in \(\bW\) with respect to
    \(\bB\). Putting together, we have that for any \(\Psi_\bG\),
    \begin{align}
        w\Bigl(e^{(\rho,-\rho_{\Psi,\bG})}\prod_{\Rt\in\Psi_\bG}(1-e^{(0,\Rt)})\Bigr)
        =(-1)^{\ell_\bG(w)}e^{(\rho,-\rho_{\Psi,\bG})}\prod_{\Rt\in\Psi_\bG}(1-e^{(0,\Rt)}).
    \end{align}
    Similarly, for \(w\in\bW_\bH\),
    \begin{align}
        w\Bigl(e^{(\rho,-\rho_{\Psi,\bH})}\prod_{\Rt\in\Psi_\bH}(1-e^{(0,\Rt)})\Bigr)
        =(-1)^{\ell_\bH(w)}e^{(\rho,-\rho_{\Psi,\bH})}\prod_{\Rt\in\Psi_\bH}(1-e^{(0,\Rt)}).
    \end{align}
    Since \((-1)^{\ell_\bG(w)}\) is also the determinant of the reflection
    matrix of \(w\) acting on \(\CharG(\bT)_\bbQ\), and the actions of
    \(\bW\) and \(\bW_\bH\) are compatible with inclusion
    \(\bW_\bH\subset\bW\), we have that for \(w\in\bW_\bH\),
    \begin{align}
        (-1)^{\ell_\bG(w)}=(-1)^{\ell_\bH(w)}.
    \end{align}
    As a result, \(r_\bH^\bG\) is \(\bW_\bH\)-invariant, and we are done.
\end{proof}

\begin{definition}
    The \notion{resultant divisor}\index{divisor!resultant} \(\bR_{\bH}^\bG\)
        \nomenclature[\(R_H^G \)]{\(\bR_\bH^\bG,\FRR_H^G\)}{the extended resultant divisor on \(\bC_{\bM,\bH}\) or \(\FRC_{\FRM,H}\)}
    (resp.~\(\bR_{\bH}^{\bG\star}\)) is the principal divisor on
    \(\bC_{\bM,\bH}\) (resp.~\(\bC_{\bM,\bH}^\star\)) defined by \(r_\bH^\bG\).
\end{definition}

\subsection{}
By previous discussion, we have equality of Cartier divisors on \(\bC_{\bM,\bH}\):
\begin{align}
    \nu_{\bH}^*\bD_{\bM}=\bD_{\bM,\bH}+2\bR_{\bH}^\bG.
\end{align}
However, unlike the Lie algebra case, the resultant divisor is not reduced in
general as seen in the following counterexample:
\begin{example}
    Let \(\bG=\SL_2\) and \(\bM\) be the monoid corresponding to the cocharacter
    \(4\CoWt\) (\(\CoWt\) being the fundamental coweight). The universal
    monoid of \(\bG\) is just \(\Mat_2\) and its abelianization is \(\bbA^1\)
    given by the usual determinant. The map from \(\bA_\bM\cong\bbA^1\) to
    \(\bA_{\Mat_2}\) is raising to the \(4\)-th power. Let \(a,b\) be the
    coordinates of the maximal toric variety in \(\Mat_2\) (which is the
    diagonal matrices), then the maximal toric variety \(\bar{\bT}_\bM\) is
    given by the ring
    \begin{align}
        k[a,b,s]/(s^4-ab).
    \end{align}
    Let \(\dual{\bH}\) be the diagonal torus in \(\dual{\bG}=\PGL_2\), then
    one can compute the maximal toric variety of \(\bM_\bH\) to be just
    \(\bbA^5\), with coordinates denoted by \(A,B,C,D,E\). The map
    \(\bar{\bT}_{\bM,\bH}\to \bar{\bT}_\bM\) is given  by the ring homomorphism
    \begin{align}
        k[a,b,s]/(s^4-ab)&\longto k[A,B,C,D,E]\\
        a&\longmapsto A^4B^3C^2D\\
        b&\longmapsto BC^2D^3E^4\\
        s&\longmapsto ABCDE
    \end{align}
    The group \(\bH\) has no root, so the discriminant divisor of \(\bM_\bH\) is
    empty, and \(\nu_\bH^*\bD_\bM=2\bR_\bH^\bG\). It is also straightforward to
    see that the divisor \(\nu_\bH^*\bD_\bM\) is cut out by function
    \begin{align}
        (a-b)^2=\bigl(BC^2D(A^4B^2-D^2E^4)\bigr)^2.
    \end{align}
    The divisor \(\bR_\bH^\bG\) is thus cut out by
    \begin{align}
        BC^2D(A^4B^2-D^2E^4),
    \end{align}
    which contains a non-reduced irreducible component cut out by \(C^2\).
\end{example}

\subsection{}
We also note that in general \(\bR_\bH^\bG\) strictly contains the divisor
\begin{align}
    (\bR_\bH^\bG)'=\bar{\bR_\bH^\bG\cap\bC_{\bM,\bH}^\x},
        \nomenclature[\(R_H^G_prime \)]{\((\bR_\bH^\bG)',(\FRR_H^G)'\)}{the non-boundary part of \(\bR_\bH^\bG\) or \(\FRR_H^G\)}
\end{align}
where we take the scheme-theoretic closure on the right-hand side. Since
\(\bC_{\bM,\bH}\) is factorial, \((\bR_\bH^\bG)'\) is also a principal divisor.
In this example above, \((\bR_\bH^\bG)'\) is
cut out by \(A^4B^2-D^2E^4\), and it misses out the components cut out by
\(BC^2D\). It is easy to see that \((\bR_\bH^\bG)'\) is always reduced because it is
generically so (by reducing to groups of semisimple rank \(1\)) and is
Cohen--Macaulay. As a result, \(\bR_\bH^\bG\) is the sum of reduced divisor
\((\bR_\bH^\bG)'\) and multiples of some components in
\(\bE_{\bM,\bH}=\bC_{\bM,\bH}-\bC_{\bM,\bH}^\x\).

\subsection{}
The non-reducedness of \(\bR_\bH^\bG\) is sensitive to the choice of monoid
\(\bM\), or equivalently, the representations of \(\dual{\bG}\) we are trying to
encode. Roughly speaking (see
\Cref{prop:characterizing_non_reducedness_of_resultant}), it only happens when a
\(\dual{\bG}\)-representation with highest-weight \(-w_0(\theta_i)\) contains a
\(\dual{\bH}\)-representation with highest-weight \(-w_{\bH,0}(\lambda_{ij})\),
such that the difference between \(-w_0(\theta_i)\) and the
\(\dual{\bG}\)-dominant element in \(-\bW\lambda_{ij}\) has at least height
\(2\).

To illustrate, in the example
above, the non-reduced component is cut out by \(C^2\), which secretly
corresponds to the \(0\)-weight space in the representation \(\Sym^4\Qlb^2\).
In this case, the difference between \(4\CoWt=2\CoRt\) and
\(0\) is \(2\CoRt\), which has height \(2\). On the other hand,
\(D\) corresponds to weight \(-\CoRt\), whose \(\dual{\bG}\)-dominant element in
its \(\bW\)-orbit is \(\CoRt\), which is height-\(1\) away from highest-weight
\(2\CoRt\), and correspondingly the component cut out by \(D\) is reduced.
A more precise statement is the following:

\begin{proposition}
    \label[proposition]{prop:characterizing_non_reducedness_of_resultant}
    Let \(\bE_{\bM,\bH,ij}\subset\bE_{\bM,\bH}\) be the irreducible component
    cut out by the coordinate dual to \(\lambda_{ij}\) in \(\bC_{\bM,\bH}\). Let
    \(\lambda_{ij}'\) be the unique \(\bG\)-dominant element in the Weyl group orbit
    \(\bW(-w_{\bH,0}(\lambda_{ij}))=-\bW\lambda_{ij}\) and
    \begin{align}
        h_{ij}=\Pair{\rho}{-w_0(\theta_i)-\lambda_{ij}'}.
    \end{align}
    Then we have equality of principal divisors
    \(\bR_\bH^\bG=(\bR_\bH^\bG)'+(\bR_\bH^\bG)''\) where
    \begin{align}
        (\bR_\bH^\bG)''=\sum_{ij}h_{ij}\bE_{\bM,\bH,ij}.
        \nomenclature[\(R_H^G_prime_prime \)]{\((\bR_\bH^\bG)'',(\FRR_H^G)''\)}{the boundary part of \(\bR_\bH^\bG\) or \(\FRR_H^G\)}
    \end{align}
    Similarly, we have the coarse version by taking intersection with
    \(\bC_{\bM,\bH}^\star\).
\end{proposition}
\begin{proof}
    The divisors \(\bE_{\bM,\bH,ij}\) are affine subspaces in \(\bC_{\bM,\bH}\),
    and we only need to compute their multiplicities in \(\bR_\bH^\bG\). To do
    so, it suffices to take a \(\cO\)-valued point (where
    \(\cO=\bar{k}\powser{\pi}\)) \(a_\bH\in\bC_{\bM,\bH}(\cO)\) such that the
    valuation of divisor \(a_\bH^*\bE_{\bM,\bH,ij}\) is \(1\) for a unique pair
    of \(i\) and \(j\) and \(0\) otherwise, and compute the valuation of
    \(a_\bH^*\bR_\bH^\bG\). Since \(\bD_{\bM,\bH}\) intersects properly with
    \(\bE_{\bM,\bH}\), we only need to consider those \(a_\bH\) contained in
    \(\bC_{\bM,\bH}^\rss\), and in such case the valuation of
    \(a_\bH^*\bR_\bH^\bG\) is half the valuation of \(a_\bH^*\bD_\bM\). When
    such valuation reaches minimum at some \(a_{\bH,0}\), the valuation of
    \(a_{\bH,0}^*(\bR_\bH^\bG)'\) must be \(0\) because \((\bR_\bH^\bG)'\) also
    intersects with \(\bE_{\bM,\bH,ij}\) properly, and so this minimal value is
    equal to the multiplicity of \(\bE_{\bM,\bH,ij}\).

    Given such an \(a_\bH\), we may find some ramified extension
    \(\cO_l=\bar{k}\powser{\pi^{1/l}}\) of \(\cO\) such that \(a_\bH\) lifts to a
    point \(t_\bH\in\bar{\bT}_{\bM,\bH}(\cO_l)\). Let
    \(F_l=\bar{k}\lauser{\pi^{1/l}}\), then we may also view \(t_\bH\) as a point in
    \(\bT_{\bM,\bH}(F_l)\). That the arc \(a_\bH\) intersects
    \(\bE_{\bM,\bH,ij}\) transversally for a unique pair of \(i\) and \(j\)
    implies that
    \begin{align}
        t_\bH=\pi^{(l\lambda_{ij},\mu)/l}t_{\bH,0},
    \end{align}
    where \(\mu\le_\bH -w_{\bH,0}(l\lambda_{ij})\) and
    \(t_{\bH,0}\in\bT_{\bM,\bH}(\cO_l)\). Without loss of generality, we may
    assume that \(\mu\) is \(\bH\)-dominant. The valuation of \(\Disc_{\bH,+}\)
    on \(t_\bH\) may be computed as follows:
    \begin{align}
        \val_\cO\bigl(\Disc_{\bH,+}(t_\bH)\bigr)
        =\Pair{2\rho_\bH}{-w_{\bH,0}(\lambda_{ij})-\mu/l}
        +\sum_{\substack{\Rt\in\Roots_\bH\\\Pair{\Rt}{\mu}=0}}\val_\cO(1-\Rt(t_{\bH,0})).
    \end{align}
    Since we assume \(a_\bH\in\bC_{\bM,\bH}^\rss\), this valuation is \(0\),
    and so \(\mu=-w_{\bH,0}(l\lambda_{ij})\).

    The image of \(t_\bH\) in \(\bT_\bM(F_l)\) can be written as
    \begin{align}
        t=\pi^{(\theta_i,-w_{\bH,0}(\lambda_{ij}))}t_0,
    \end{align}
    where \(t_0\in\bT_\bM(\cO_l)\) and \(t_0\) is the image of \(t_{\bH,0}\)
    under the map \(\bT_{\bM,\bH}\to\bT_\bM\). We then similarly have
    \begin{align}
        \val_\cO\bigl(\Disc_{+}(t)\bigr)
        =\Pair{2\rho}{-w_{0}(\theta_i)-\lambda_{ij}'}
        +\sum_{\substack{\Rt\in\Roots\\\Pair{\Rt}{-w_{\bH,0}(\lambda_{ij})}=0}}\val_\cO(1-\Rt(t_{0})).
    \end{align}
    This valuation reaches minimum when \(t_0\) is sufficiently general so that
    \(\val_\cO(1-\Rt(t_0))=0\) for all \(\Rt\in\Roots\) perpendicular to
    \(-w_{\bH,0}(\lambda_{ij})\). Note that this minimum can be achieved for
    some \(t_{\bH,0}\in \bT_{\bM,\bH}(\cO)\), and for those \(a_\bH\), we have
    \begin{align}
        \val_{\cO}(a_\bH^*\bD_\bM)=2\val_{\cO}(a_\bH^*\bR_\bG^\bH)
        =\Pair{2\rho}{-w_{0}(\theta_i)-\lambda_{ij}'}=2h_{ij}.
    \end{align}
    This finishes the proof.
\end{proof}

\subsection{}
We have regular centralizer \(\bJ_{\bM,\bH}\) of the adjoint
\(\bH\)-action on \(\bM_{\bH}\) as well as the group scheme \(\bJ_{\bM,\bH}^1\) on
\(\bC_{\bM,\bH}\) using Galois description of regular centralizer. 
Similar to Lie algebra case, we have the following results.

\begin{lemma}
    \label[lemma]{lem:reg_cent_transfer_split}
    There is a canonical commutative diagram over \(\bC_{\bM,\bH}\)
    \begin{equation}
        \begin{tikzcd}
              \nu_{\bH}^*\bJ_{\bM}\ar[r]  \ar[d]& \bJ_{\bM,\bH}  \ar[d]\\
              \nu_{\bH}^*\bJ_{\bM} \ar[r]^1     & \bJ_{\bM,\bH}^1    
        \end{tikzcd}
    \end{equation}
    such that all arrows are isomorphisms over \(\bC_{\bM,\bH}^{\bG\hy\rss}\).
\end{lemma}
\begin{proof}
    We have canonical homomorphism by adjunction
    \begin{align}
        \nu_{\bH}^*\pi_{\bM*}\bT\longto \pi_{\bM,\bH*}\tilde{\nu}_{\bH}^*\bT.
    \end{align}
    Taking invariant under either \(\bW\) or \(\bW_{\bH}\), we have
    \begin{align}
        \nu_{\bH}^*\bJ_{\bM}^1
        =\nu_{\bH}^*(\pi_{\bM*}\bT)^\bW\longto \nu_{\bH}^*(\pi_{\bM*}\bT)^{\bW_{\bH}}\longto
        (\pi_{\bM,\bH*}\tilde{\nu}_{\bH}^*\bT)^{\bW_{\bH}}=\bJ_{\bM,\bH}^1.
    \end{align}
    This map is clearly an isomorphism over the \(\bG\)-regular semisimple locus
    (which is contained in \(\bH\)-regular semisimple locus).
    Using the identification of \(\bJ_{\bM}\) with \(\bJ_{\bM}'\), we see that if
    the condition defining \(\bJ_{\bM}'\) is satisfied, then the analogous
    condition for \(\bJ_{\bM,\bH}'\) is also satisfied since
    \(\Roots_{\bH}\subset\Roots\). Therefore, the map
    \(\nu_{\bH}^*\bJ_{\bM}\to\bJ_{\bM,\bH}^1\) factors through
    \(\bJ_{\bM,\bH}\), and we are done.
\end{proof}

\subsection{}
\label{sub:absolute_resultant_divisor}
To establish the transfer between regular semisimple orbit of \(\FRM\) and
\(G\)-regular semisimple orbit in \(\FRM_H\), we need to construct a canonical map
\begin{align}
    \nu_H\colon \FRC_{\FRM,H}\longto \FRC_{\FRM}.
\end{align}
Through \(\bFo_{\bG}\), we have twisted form \(\FRT_{\FRM}\) of \(\bar{\bT}_{\bM}\),
and the canonical action of \(1\rtimes_{\bH}\Theta\) induces twisted form
\(\FRT_{\FRM,H}\) of \(\bar{\bT}_{\bM,\bH}\). However, \(\tilde{\nu}_{\bH}\)
is not \(\Theta\)-equivariant since the image of
\(1\rtimes_{\bH}\Theta\) is not \(1\rtimes_{\bG}\Theta\), therefore
there is no natural map between \(\FRT_{\FRM}\) and \(\FRT_{\FRM,H}\).
Nevertheless, we still have the canonical map \(\nu_H\) because the action of
\(\bW_{\bH}\rtimes\Theta\) on \(\bC_{\bM,\bH}\) factors through quotient
\(\Theta\) and same is true for \(\bW\rtimes \Theta\) acting on
\(\bC_{\bM}\).

Similarly, we have relation between twisted discriminant divisors
\begin{align}
    \nu_H^*\FRD_{\FRM}=\FRD_{\FRM,H}+2\FRR_H^G,
\end{align}
where \(\FRR_H^G\) is a principal divisor on \(\FRC_{\FRM,H}\). We also have a
decomposition of principal divisors
\begin{align}
    \FRR_H^G=(\FRR_H^G)'+(\FRR_H^G)'',
\end{align}
where divisor \((\FRR_H^G)'\) is the maximal subdivisor that intersects
with \(\FRE_{\FRM,H}\) properly and is reduced, while \((\FRR_H^G)''\) is
characterized by \Cref{prop:characterizing_non_reducedness_of_resultant}.
The map between abelianizations also have natural twisted form
\begin{align}
    \underline{\nu}_H\colon\FRA_{\FRM,H}\longto\FRA_{\FRM}.
\end{align}
For the same reason, the central diagonalizable group 
\(\bZ_{\bM}^\kappa\) has twisted form
\(Z_{\FRM}^\kappa\), and it maps naturally to \(\FRA_{\FRM,H}^\x\) and
surjectively onto \(Z_{\FRM}\).

\subsection{}
As a continued trend, everything so far also has a coarse version, so we have maps
\(\nu_H^\star\colon\FRC_{\FRM,H}^\star\to\FRC_\FRM\) and
\(\underline{\nu}_H^\star\colon\FRA_{\FRM,H}^\star\to\FRA_{\FRM}\),
divisors \(\FRD_{\FRM,H}^\star\), \(\FRR_{H}^{G\star}\), and so on.

\subsection{}
We establish a canonical map
\(\nu_H\colon\FRJ_{\FRM}\to\FRJ_{\FRM,H}\) as follows: regarding \(\OGT_\kappa\) as
a finite \'etale map \(\OGT_\kappa\colon X_\kappa\to X\), we have a finite flat
map
\begin{align}
    \pi_{\kappa}\colon X_\kappa\x_X\FRT_{\FRM}\simeq X_\kappa\x \bar{\bT}_{\bM}\longto \FRC_{\FRM}
\end{align}
such that over any geometric point \(\bar{v}\in X\) it is generically a
\(\bW\rtimes\Theta\)-torsor.
Then we may alternatively describe \(\FRJ_{\FRM}^1\) as the fixed-point scheme
\begin{align}
    \FRJ_{\FRM}^1=\pi_{\kappa*}(X_\kappa\x
    \bar{\bT}_{\bM}\x\bT)^{\bW\rtimes\Theta},
\end{align}
and \(\FRJ_{\FRM}\) can be identified with the subfunctor whose \(S\)-points for a
\(\FRC_{\FRM}\)-scheme \(S\) consists of maps
\begin{align}
    f\colon S\x_{\FRC_{\FRM}}(X_\kappa\x\bar{\bT}_{\bM})\longto \bT
\end{align}
such that for any geometric point \(x\in
S\x_{\FRC_{\FRM}}(X_\kappa\x\bar{\bT}_{\bM})\), if \(s_\Rt(x)=x\) for a root
\(\Rt\), then \(\Rt(f(x))\neq -1\).
We also have 
\begin{align}
    \FRJ_{\FRM,H}^1=\pi_{\kappa*}(X_\kappa\x
    \bar{\bT}_{\bM,\bH}\x\bT)^{\bW_{\bH}\rtimes\Theta},
\end{align}
and a similar description for \(\FRJ_{\FRM,H}\) using roots
\(\Rt\in\Roots_{\bH}\).
Using commutative diagram
\begin{equation}
    \begin{tikzcd}
        X_\kappa\x\bar{\bT}_{\bM,\bH} \ar[r, "\Id\x\tilde{\nu}_{\bH}"]\ar[d] 
                & X_\kappa\x\bar{\bT}_{\bM} \ar[d]\\
        \FRC_{\FRM,H} \ar[r, "\nu_H"] & \FRC_{\FRM}
    \end{tikzcd}
\end{equation}
and the same argument in \Cref{lem:reg_cent_transfer_split}, we have:
\begin{lemma}
    \label[lemma]{lem:reg_cent_transfer_general}
    There is a canonical commutative diagram over \(\FRC_{\FRM,H}\)
    \begin{equation}
        \begin{tikzcd}
            \nu_H^*\FRJ_{\FRM} \ar[r]\ar[d] & \FRJ_{\FRM,H}\ar[d]\\
            \nu_H^*\FRJ_{\FRM}^1\ar[r] & \FRJ_{\FRM,H}^1
        \end{tikzcd}
    \end{equation}
    such that all arrows are isomorphisms over \(\FRC_{\FRM,H}^{G\hy\rss}\).
\end{lemma}

\subsection{}
We have induced maps of quotient stacks
\begin{align}
    \Stack{\FRC_{\FRM,H}^\star/Z_{\FRM}^{\kappa\star}}\longto
    \Stack{\FRC_{\FRM,H}/Z_{\FRM}^\kappa}\longto\Stack{\FRC_{\FRM,H}'/Z_{\FRM}}\longto \Stack{\FRC_{\FRM}/Z_{\FRM}},
\end{align}
where the first two maps are generically isomorphisms over the invertible locus. The discriminant
divisor, the regular centralizer, etc. all descend to these quotients. However,
even if one has a Steinberg quasi-section over \(\FRC_{\FRM,H}\), it may not
descend to \(\Stack{\FRC_{\FRM,H}/Z_{\FRM}^\kappa}\).

\subsection{}
\label{sub:invariant_theory_for_levi}
Notice that a large part of our results for the endoscopic monoid only
depends on the fact that the dual group of \(\bH\) embeds into that of \(\bG\)
as a subgroup. This means that we can try to do the same for
Levi subgroups as well.

Let \(\bS\subset\bT\) be a subtorus, then the centralizer \(\bL\) of \(\bS\) in
\(\bG\) is a Levi subgroup containing maximal torus \(\bT\). Using the
pinning \(\bSPL\), its dual \(\dual{\bL}\) can be identified with a Levi
subgroup in \(\dual{\bG}\). Choosing a sufficiently general element \(\gamma\) 
in the center of \(\dual{\bL}\), we can realize \(\dual{\bL}\) as the centralizer
of \(\gamma\) in \(\dual{\bG}\). Therefore, similar to endoscopic group \(\bH\),
we have for any monoid \(\bM\in\FM(\bG^\SC)\) a canonical associated monoid
\(\bM_{\bL}\in\FM(\bL^\SC)\)
    \nomenclature[\(M_L \)]{\(\bM_\bL,\FRM_L\)}{the Levi monoid associated with Levi subgroup \(\bL\subset\bG\) or \(L\subset G\)}
using the same construction. We also have the
non-flat monoid \(\bM_{\bL}'\), and canonical homomorphism of reductive monoids
\begin{align}
    \bM_{\bL}\longto\bM_{\bL}'\longto\bM,
\end{align}
where the second arrow is induced  by inclusion \(\bL\subset\bG\) (which is not
necessarily present for endoscopic groups). Let \(\bC_{\bM,\bL}\)
(resp.~\(\bC_{\bM,\bL}'\)) be the GIT
quotient of \(\bM_{\bL}\) (resp.~\(\bM_{\bL}'\)) by \(\bL\), and
let \(\bZ_{\bM}^\bL\) be the preimage
of \(\bZ_{\bM}\) in \(\bM_{\bL}\), then we have canonical maps
\begin{align}
    \Stack*{\bC_{\bM,\bL}/\bZ_{\bM}^\bL}\longto\Stack*{\bC_{\bM,\bL}'/\bZ_{\bM}}\longto
    \Stack*{\bC_{\bM}/\bZ_{\bM}},
\end{align}
where the second map is finite, and the first map is generically an isomorphism,
and \Cref{lem:arc_lifting_to_endoscopic_quotient} holds (after replacing
\(\bH\) by \(\bL\) and \(\bZ_{\bM}^\kappa\) by \(\bZ_{\bM}^\bL\)).

Let \(\tilde{\bL}\) be the centralizer of \(\bS\) in \(\bG\rtimes\Out(\bG)\),
then its neutral component is \(\bL\), and its quotient group of connected
component \(\pi_0(\tilde{\bL})\) is a subgroup of \(\Out(\bG)\). Suppose that
the image of \(\Theta\) is contained in \(\pi_0(\tilde{\bL})\),
then we have induced maps similar to endoscopic case
\begin{align}
    \Stack*{\FRC_{\FRM,L}/Z_{\FRM}^L}\longto\Stack*{\FRC_{\FRM,L}'/Z_{\FRM}}\longto
    \Stack*{\FRC_{\FRM}/Z_{\FRM}}.
\end{align}
Note that in this case \(\bS\) is a split subtorus in the center of \(L\).
Of course, we can also construct the \emph{coarse} Levi monoid \(\bM_\bL^\star\) and all
other related notions, but they are not so relevant to us in this book, so we
will ignore them.

\section{Fundamental Lemma for Spherical Hecke Algebras} 
\label{sec:fundamental_lemma_for_spherical_hecke_algebras}

In this section, we give the statement of fundamental lemmas for spherical Hecke
algebras. We vaguely follow the same outline as in the first chapter of
\cite{Ng10}, but there are also many new ingredients necessary.

\subsection{}
Let \(F_v=k\lauser{\pi_v}\)
    \nomenclature[\(F_v \)]{\(F_v\)}{the local field \(k_v\lauser{\pi_v}\) at place \(v\)}
be a field of Laurent series over \(k\) and
\(\cO_v=k\powser{\pi_v}\)
    \nomenclature[\(O"cal_v \)]{\(\cO_v\)}{the complete local ring \(k_v\powser{\pi_v}\) at place \(v\)}
its ring of integers. We fix a separable closure
\(F_v^\sep\)
    \nomenclature[\(F_sep \)]{\(F^\sep\)}{a separable closure of a field \(F\)}
of \(F_v\), and let \(\Gamma_v=\Gal(F_v^\sep/F_v)\).
    \nomenclature[\(Gamma_v \)]{\(\Gamma_v\)}{the local Galois group \(\Gal(F_v^\sep/F_v)\)}
As before, let
\(G\) be a reductive group over \(\cO_v\) obtained from \(\bG\) and
\(\Out(\bG)\)-torsor \(\OGT_G\). Let \(\gamma\in G(F_v)\) be a semisimple
element. Since the derived
subgroup of \(G\) is not necessarily simply-connected, the centralizer
\(I_\gamma\)
    \nomenclature[\(I_gamma \)]{\(I_\gamma\)}{the centralizer of \(\gamma\) in \(G\)}
of \(\gamma\) in \(G\) may not be a torus even if \(\gamma\) is
regular. For example, the
diagonal matrix with entries \(1\) and \(-1\), viewed as an element in
\(\PGL_2\) is regular, but its centralizer is disconnected. To remedy this, we
need a slightly stronger variant of regular semisimplicity:
\begin{definition}
    A semisimple element (or orbit) \(\gamma\in G(F_v)\) is called \notion{strongly
    regular semisimple}\index{orbit!strongly regular semisimple} if its
    centralizer \(I_\gamma\) is a torus. The locus
    of strongly regular semisimple in \(G\git G\) is denoted by
    \(\FRC_G^\srs\)\index{locus!strongly regular semisimple}.
    \nomenclature[\(.srs \)]{\((\cdot)^\srs,(\cdot)_\srs\)}{the (restriction to the) strongly regular semisimple locus}
\end{definition}

\begin{lemma}
    The subset \(\FRC_G^\srs\) is open dense in \(G\git G\),
    and the natural map \(T\git W\to G\git G\) is an isomorphism over
    \(\FRC_G^\srs\). Moreover, \(\FRC_G^\srs\) is precisely the locus over which
    the cameral cover \(T\to T\git W\) is \'etale.
\end{lemma}
\begin{proof}
    The question is geometric, so we may replace \(F_v\) by its separable
    closure, and so any semisimple element is conjugate to an element in \(T\).
    Clearly, the cameral cover \(T\to T\git W\) is a \(W\)-torsor over any
    strongly regular semisimple element.

    Conversely, it is known (see \cite{Hu95}*{\S~2.2}) that the component group
    of the centralizer of an
    element in \(T\) is generated by elements in \(\Norm_G(T)\). This
    implies that any point in \(T\) at which the cameral cover \(T\to T\git W\)
    is \'etale is strongly regular semisimple. That \(T\git W\to G\git G\) is an
    isomorphism over \(\FRC_G^\srs\) can be easily deduced from the case where
    \(G^\Der=G^\SC\).
\end{proof}

\subsection{}
We will assume from now on that \(\gamma\) is strongly regular semisimple,
and let \(a\in\FRC_G^\srs(F_v)\) be the image of \(\gamma\).
If \(\gamma'\in G(F_v)\) is another \(F_v\)-point over \(a\), then (with our
assumption on \(\Char(k)\)) there exists
some \(g\in G(F_v^\sep)\) such that \(\gamma'=g\gamma g^{-1}\). For
any \(\sigma\in\Gamma_v\), we have \(g\sigma(g)^{-1}\in
I_{\gamma}(F_v^\sep)\), and so \(\sigma\mapsto g\sigma(g)^{-1}\) defines a
Galois cocycle in \(I_\gamma\), whose cohomology class is denoted by
\begin{align}
    \inv(\gamma,\gamma')\in \RH^1(F_v,I_{\gamma}).
    \nomenclature[\(inv_gamma_gamma' \)]{\(\inv(\gamma,\gamma')\)}{the relative
    position between two conjugacy classes in the same geometric conjugacy class}
\end{align}
This class depends only on the \(G(F_v)\)-conjugacy class of \(\gamma'\) but not
on \(g\). The image of \(\inv(\gamma,\gamma')\) in \(\RH^1(F_v,G)\) is trivial,
and in this way we have a bijection between the kernel of the map
\begin{align}
    \RH^1(F_v,I_\gamma)\longto \RH^1(F_v,G),
\end{align}
and the set of \(G(F_v)\)-conjugacy classes inside the
\(G(F_v^\sep)\)-conjugacy class of \(\gamma\).

\subsection{}
Fixing an \(F_v^\sep\)-point on \(\OGT_G\), we obtain a homomorphism
\(\OGT_G^\bullet\colon\Gamma_v\to \Out(\bG)\). By Tate--Nakayama duality, the
cohomological set \(\RH^1(F_v,G)\) has a natural abelian group structure
characterized by its Pontryagin dual:
\begin{align}
    \RH^1(F_v,G)^*\simeq\pi_0\bigl((\bZ_{\dual{\bG}})^{\OGT_G^\bullet(\Gamma_v)}\bigr).
    \nomenclature[\(.* \)]{\((\cdot)^*\)}{the Pontryagin dual, or the dual vector
    space/bundle, depending on the context}
\end{align}

Choose a geometric point \(x_a\in \bT(F_v^\sep)\)
    \nomenclature[\(x_a \)]{\(x_a\)}{a point in \(\bT(F_v^\sep)\) lying over a geometric conjugacy class \(a\) to rigidify the identification \(\FRJ_a(F_v^\sep)\cong \bT(F_v^\sep)\)}
over \(a\), then it
corresponds to a lift of \(\OGT_G^\bullet\) to a homomorphism
\begin{align}
    \pi_a^\bullet\colon\Gamma_v\longto \bW\rtimes\Out(\bG).
    \nomenclature[\(pi_a_bullet \)]{\(\pi_a^\bullet\)}{the monodromy
    \(\Gamma_v\longto \bW\rtimes\Out(\bG)\) lifting \(\OGT_G^\bullet\) induced by \(x_a\)}
\end{align}
Although \(G\) is not an object in \(\FM(G^\SC)\), the regular centralizer
\(\FRJ_a\) can still be defined for strongly regular semisimple elements and its
Galois description still holds using the same proof because
\(T^\srs\to\FRC_G^\srs\) is \'etale. Thus, we have isomorphisms over \(F_v\)
\begin{align}
    I_\gamma\simeq \FRJ_a\simeq
    \Spec{F_v^\sep}\wedge^{\Gamma_v,\pi_a^\bullet}\bT.
\end{align}
Using Tate--Nakayama duality again, we have
\begin{align}
    \RH^1(F_v,\FRJ_a)^*\simeq\pi_0\bigl(\dual{\bT}^{\pi_a^\bullet(\Gamma_v)}\bigr).
\end{align}
Note that here the isomorphism depends on the choice of \(x_a\). The
inclusion \(\iota\colon \dual{\bT}\to\dual{\bG}\) is \(\Gamma_v\)-equivariant up
to conjugacy in the following sense: for any \(t\in\dual{\bT}\) and any
\(\sigma\in\Gamma_v\), we always have that \(\OGT_G^\bullet(\sigma)(\iota(t))\)
and \(\iota(\pi_a^\bullet(\sigma)(t))\) are \(\dual{\bG}\)-conjugate. Therefore,
we have an induced canonical map
\begin{align}
    \pi_0\bigl((\bZ_{\dual{\bG}})^{\OGT_G^\bullet(\Gamma_v)}\bigr)\longto
    \pi_0\bigl(\dual{\bT}^{\pi_a^\bullet(\Gamma_v)}\bigr),
\end{align}
whose Pontryagin dual is the map \(\RH^1(F_v,I_\gamma)\longto \RH^1(F_v,G)\).

\subsection{}
We fix a Haar measure \(\dd g_v\) on \(G(F_v)\) such that \(G(\cO_v)\) has
volume \(1\), as well as a non-zero Haar measure \(\dd t_v\) on \(\FRJ_a(F_v)\).
Using the canonical isomorphism \(I_\gamma\simeq \FRJ_a\), we have an induced Haar
measure on \(I_\gamma(F_v)\). With these choices, for any locally constant and
compactly supported function \(f\) on \(G(F_v)\), we may define orbital integral
\begin{align}
    \OI_\gamma(f,\dd t_v)=\int_{I_\gamma(F_v)\backslash G(F_v)}f(g_v^{-1}\gamma
    g_v)\frac{\dd g_v}{\dd t_v}.
    \nomenclature[\(OI_gamma \)]{\(\OI_\gamma\)}{the orbital integral induced by \(\gamma\)}
\end{align}

\begin{definition}
    Let \(\kappa\in \dual{\bT}^{\pi_a^\bullet(\Gamma_v)}\). We define the
    \notion{\(\kappa\)-orbital integral}\index{orbital integral!\(\kappa\)-} of \(a\) as the sum
    \begin{align}
        \OI_a^\kappa(f,\dd
        t_v)=\sum_{\gamma'}\Pair{\inv(\gamma,\gamma')}{\kappa}^{-1}\OI_\gamma(f,\dd
        t_v),
        \nomenclature[\(OI_kappa_gamma \)]{\(\OI_\gamma^\kappa\)}{the \(\kappa\)-orbital integral induced by \(\gamma\)}
    \end{align}
    where \(\gamma'\) ranges over all \(G(F_v)\)-conjugacy classes over \(a\),
    and \(\gamma\) is a fixed choice. When \(\kappa=1\), we denote
    \(\OI_a^\kappa\) by \(\SOI_a\),
    \nomenclature[\(S{}OI_gamma \)]{\(\SOI_\gamma\)}{the stable orbital integral induced by \(\gamma\)}
    called the \notion{stable orbital integral}\index{orbital integral!stable}.
\end{definition}
Note that \(\kappa\)-orbital integral is sensitive to the choice of base point
\(\gamma\) if \(\kappa\neq 1\), while the stable orbital integral is not. The
former also depends on a choice of the geometric point \(x_a\).

\subsection{}
\label{sub:endo_transfer_map_in_group_form_srs}
Let \((\kappa,\OGT_\kappa)\) be an endoscopic datum of \(G\) and \(H\) is the
endoscopic group.  We may view the \(\pi_0(\kappa)\)-torsor as an \'etale cover
\(\OGT_\kappa\colon X_\kappa\to X\). The canonical homomorphism in
\Cref{lem:kappa_G_H_compatibility} induces map
\begin{align}
    \FRC_H^{G\hy\srs}\defeq X_\kappa\x\bT^\srs \git \bW_\bH\rtimes\pi_0(\kappa)
    \longto X_\kappa\x\bT^\srs\git \bW\rtimes\pi_0(\kappa)\simeq \FRC_G^\srs,
\end{align}
where \(\FRC_H^{G\hy\srs}\) is called the \notion{strongly \(G\)-regular
semisimple locus}\index{locus!strongly \(G\)-regular
semisimple}\index{orbit!strongly \(G\)-regular semisimple}. It is easy to see
that it is an open subset of \(\FRC_H^\srs\).
We say a conjugacy class \(a\in\FRC_G^\srs(F_v)\) and
\(a_H\in\FRC_H^{G\hy\srs}(F_v)\) \notion{match each other}\index{orbit!matching} if \(a\) is the image of
\(a_H\). In this case, we have a canonical isomorphism
\(\FRJ_a\simeq\FRJ_{H,a_H}\). When the endoscopic datum is pointed, we shall
always choose the same point in
\(\bT(F_v^\sep)\) for both \(x_a\) and \(x_{a_H}\) lying over the base point of
the endoscopic datum, but the choice
of \(G(F_v)\)-conjugacy class \(\gamma\) over \(a\) remains arbitrary. The choice of \(x_{a_H}\)
does not affect the stable orbital integral of \(a_H\), but it is already
implicit in the definition of transfer factor that \(x_{a_H}=x_a\).

\subsection{}
Now we make a digression and discuss the Satake functions. There will be more
related details in
\Cref{chap:generalized_affine_springer_fibers,chap:global_constructions}, but
a comprehensive reference is \cite{Zh17}.

For split group \(\bG\), we may consider its affine Grassmannian defined as
\(\Gr_{\bG}=\Loop{\bG}/\Arc{\bG}\),
    \nomenclature[\(L{}Y"bbold \)]{\(\Loop{Y},\Loop_v{Y}\)}{the loop space of a \(k\)- or \(k_v\)-scheme \(Y\)}
    \nomenclature[\(L{}Y"bbold^+ \)]{\(\Arc{Y},\Arc_v{Y}\)}{the arc space of a \(k\)- or \(k_v\)-scheme \(Y\)}
    \nomenclature[\(G{}r_G \)]{\(\Gr_{G},\Gr_{G,v}\)}{the local affine Grassmannian (at place \(v\)) of a group \(G\)}
whose \(k\)-points is the
quotient set \(\bG(F_v)/\bG(\cO_v)\). It is an ind-projective ind-scheme over
\(k\). The arc group \(\Arc{\bG}\) acts
naturally on \(\Gr_{\bG}\), and its orbits are given by the Cartan decomposition
\begin{align}
    \bG(F_v)=\coprod_{[\lambda]}\bG(\cO_v)\pi_v^\lambda\bG(\cO_v),
\end{align}
where \([\lambda]\) ranges over the \(\bW\)-orbits in \(\CoCharG(\bT)\), and we
represent each \(\bW\)-orbit \([\lambda]\) by the unique \(\bB\)-dominant
coweight \(\lambda\) within. The
\(\Arc{\bG}\)-orbits, denoted by \(\Gr_{\bG}^\lambda\)
    \nomenclature[\(G{}r_G_lambda \)]{\(\Gr_{G}^\lambda,\Gr_{G,v}^\lambda\)}{the affine Schubert cell (at place \(v\)) of a group \(G\) and type \(\lambda\)}
for each \(\lambda\), form a
stratification of \(\Gr_{\bG}\) such that \(\Gr_{\bG}^\mu\) is contained in the
closure of \(\Gr_{\bG}^\lambda\) if and only if \(\mu\le \lambda\).
Let \(\Gr_{\bG}^{\le\lambda}\)
    \nomenclature[\(G{}r_G_lambda_closure \)]{\(\Gr_{G}^{\le\lambda},\Gr_{G,v}^{\le\lambda}\)}{the affine Schubert variety (at place \(v\)) of a group \(G\) and type \(\lambda\)}
be the closure of \(\Gr_{\bG}^\lambda\).

We may consider the standard intersection complex \(\IC^\lambda\)
    \nomenclature[\(I{}C^lambda \)]{\(\IC^\lambda\)}{the intersection complex of \(\Gr_G^{\le\lambda}\)}
with \(\Qlb\)-coefficients (\(\ell\) is coprime to \(p\)) on each
\(\Gr_{\bG}^{\le\lambda}\). By geometric Satake isomorphism, these complexes are
exactly the simple objects (up to isomorphism) in the category of
\(\Arc{\bG}\)-equivariant perverse sheaves on \(\Gr_{\bG}\), and
the latter is equivalent to the Tannakian category of finite dimensional
representations of \(\dual{\bG}\). Under such equivalence, we may identify
\(V_\lambda\),
    \nomenclature[\(V_lambda \)]{\(V_\lambda\)}{the irreducible representation of \(\dual{\bG}\) with highest weight \(\lambda\)}
the irreducible representation  of \(\dual{\bG}\) of highest
weight \(\lambda\), with cohomology space of \(\IC^\lambda\)
\begin{align}
    V_\lambda\simeq\RH^\bullet(\Gr_{\bG},\IC^\lambda).
\end{align}

On the other hand, by Grothendieck's sheaf-function dictionary, each
\(\IC^\lambda\) induces a function \(\fS^\lambda\) on the set of \(k\)-points of
\(\Gr_{\bG}\), hence also a \(\bG(\cO_v)\)-bi-invariant function on
\(\bG(F_v)\), still denoted by \(\fS^\lambda\). By induction, we can easily see
that the collection of \(\fS^\lambda\) form a basis of the spherical
algebra different from the characteristic functions on each
\(\bG(\cO_v)\pi_v^\lambda\bG(\cO_v)\). In other words, we have
\begin{align}
    \cH_{\bG,0}=\bigoplus_{\lambda}\Qlb \fS^\lambda,
    \nomenclature[\(H"cal_0 \)]{\(\cH_{G,0}\)}{the spherical Hecke algebra of an unramified group \(G\) over a local field}
\end{align}
where \(\cH_{\bG,0}\) denotes the spherical Hecke algebra of \(\bG\), or
\(\bG(\cO_v)\)-biinvariant locally constant and compactly supported functions on
\(\bG(F_v)\).

\subsection{}
When \(G\) is a quasi-split twist of \(\bG\) over \(\cO_v\), the action of
\(\Gamma_v\) on \(\CoCharG(T)\) factors through \(\Gal(\bar{k}/k)\). Let
\(\Frob_v\in\Gal(\bar{k}/k)\) be the geometric Frobenius element and also choose
a lifting to \(\Gamma_v\). In this case, the geometric Satake equivalence is
slightly more complicated. Following \cite{Zh17}*{\S~5.5}, the category of
sheaves to consider is \(\Arc{G}\)-equivariant perverse sheaves on \(\Gr_G\),
which is equivalent to continuous representations of
\(\LD{G}^{\symup{geo}}\). Here the group \(\LD{G}^{\symup{geo}}\) is the geometric
version of \(L\)-group of \(G\) induced by the geometric \(\Gamma_v\)-action on
\(\Gr_G\). In general, \(\LD{G}^{\symup{geo}}\) is not isomorphic to
\(\LD{G}\)
    \nomenclature[\(L_G \)]{\(\LD{G}\)}{the \(L\)-group of \(G\)}
due to the existence of Tate twists, but in the case where \(G\) is
defined over \(\cO_v\), there is a canonical isomorphism
\begin{align}
    \LD{G}^{\symup{geo}}\simeq\LD{G}
\end{align}
of topological groups that restricts to identity on \(\dual{\bG}\), after
choosing once and for all a half Tate twist \(\Qlb(1/2)\).

For arithmetic
purposes, we need to consider a smaller category: the category of
sheaves that becomes constant along \(\Arc{G}_{k'}\)-orbits for some finite
extension \(k'/k\), which is then equivalent to the category
\(\Rep_{\LD{G}}^\Alg\) of algebraic
representations of \(\LD{G}\), in other words, representations of \(\LD{G}\)
that are pullbacks of algebraic representations of
\(\dual{\bG}\rtimes\Gal(k'/k)\) for some finite extension \(k'/k\).

\begin{definition}
    For any representation \(V\in\Rep_{\LD{G}}^\Alg\), let \(\IC^V\)
    \nomenclature[\(I{}C^V \)]{\(\IC^V\)}{the Satake sheaf corresponding to the \(\LD{G}\)-representation \(V\)}
    be the
    corresponding \notion{Satake sheaf}\index{sheaf!Satake}\index{Satake!sheaf} and \(\fS^V\in\cH_{G,0}\)
    \nomenclature[\(f^V \)]{\(\fS^V\)}{the Satake function corresponding to the \(\LD{G}\)-representation \(V\)}
    be the \notion{Satake function}\index{function!Satake}\index{Satake!function} induced by \(\IC^V\).
\end{definition}

\subsection{}
A more concrete description of \(\IC^V\) and \(\fS^V\) is as follows:
\(V\) decomposes after restricting to \(\dual{\bG}\) as
\begin{align}
    \Res_{\dual{\bG}}^{\LD{G}}V\simeq \bigoplus_\lambda
    V_\lambda\otimes\Hom_{\dual{\bG}}(V_\lambda,V),
\end{align}
where \(\lambda\) ranges over all \(\dual{\bG}\)-highest weights. The
\(\dual{\bG}\)-action on the multiplicity space \(\Hom_{\dual{\bG}}(V_\lambda,V)\)
is trivial, so for any \(\sigma\in\Gal(\bar{k}/k)\) and any \(\lambda\), there
is a canonical map
\begin{align}
    \sigma\colon \Hom_{\dual{\bG}}(V_\lambda,V)\longto
    \Hom_{\dual{\bG}}(V_{\sigma(\lambda)},V).
\end{align}
If \(\lambda=\sigma(\lambda)\), then for any
\(\phi\in\Hom_{\dual{\bG}}(V_{\sigma(\lambda)},V)\), then the image of
\(V_\lambda\otimes\phi\) under \(\sigma\) is another copy of \(V_\lambda\) in
\(V\), hence of the form \(V_\lambda\otimes\sigma(\phi)\). By algebraicity and
since \(\Gal(\bar{k}/k)\) is topologically cyclic, \(\Hom_{\dual{\bG}}(V_\lambda,V)\)
decomposes into \(\Gal(\bar{k}/k)\)-eigenspaces. As a result, we have a direct
sum decomposition of \(\LD{G}\)-representations
\begin{align}
    V_\lambda\otimes\Hom_{\dual{\bG}}(V_\lambda,V)=\bigoplus_\phi
    V_\lambda\otimes\phi,
\end{align}
where \(\phi\) ranges over characters of \(\Gal(\bar{k}/k)\) in
\(\Hom_{\dual{\bG}}(V_\lambda,V)\). On the other hand, suppose we have two
algebraic \(\LD{G}\)-representations \(V_1\) and \(V_2\) whose restrictions to
\(\dual{\bG}\) are isomorphic, and we denote the \(\sigma\)-actions by
\(\sigma_1\) and \(\sigma_2\) respectively. Choose a \(\dual{\bG}\)-equivariant
isomorphism \(\iota\colon V_1\to V_2\), then
\(\sigma_1\iota\sigma_2^{-1}\iota^{-1}\) is a
\(\dual{\bG}\)-equivariant automorphism of \(V_1\), hence must be a scalar
\(\chi_{V_1,V_2}(\sigma)\). One can also verify that \(\chi_{V_1,V_2}\) is a
character of \(\Gal(\bar{k}/k)\) and does not depend on the choice of \(\iota\).
This implies that \(V_2\) is isomorphic to \(V_1\otimes \chi_{V_1,V_2}\) as
\(\LD{G}\)-representations. With this, we can normalize the
\(V_\lambda\)-isotypic subspace as
\begin{align}
    V_\lambda\otimes\Hom_{\dual{\bG}}(V_\lambda,V)=\bigoplus_\phi
    \RH^\bullet(\Gr_G,\IC^\lambda)\otimes\phi=\RH^\bullet(\Gr_G,\IC^\lambda)\otimes\biggl(\bigoplus_\phi\phi\biggr),
\end{align}
where \(\LD{G}\) acts on the cohomology space via \(\LD{G}^{\symup{geo}}\) and
\(\phi\) are \(\Gal(\bar{k}/k)\)-characters of finite order.

More generally, suppose \(\sigma\) is a topological generator of
\(\Gal(\bar{k}/k)\), and let \(e\ge 1\) be the minimal number so that
\(\sigma^e(\lambda)=\lambda\), then
\begin{align}
    \bigoplus_{i=0}^{e-1}V_{\sigma^i(\lambda)}\otimes\Hom_{\dual{\bG}}(V_{\sigma^i(\lambda)},V)
\end{align}
is induced from \(V_\lambda\otimes\Hom_{\dual{\bG}}(V_\lambda,V)\) by a cyclic
permutation, and so is completely determined by the \(\sigma^e\)-action on
the same. This allows us to reduce the general case to that where \(\lambda=\sigma(\lambda)\).

View \(\Hom_{\dual{\bG}}\bigl(\RH^\bullet(\Gr_{G,\bar{k}},\IC^\lambda),V\bigr)\)
as a constant sheaf on \(\Gr_{G,\bar{k}}^{\le\lambda}\), then the direct sum
\begin{align}
    \IC^V=\bigoplus_\lambda \IC^\lambda\otimes
    \Hom_{\dual{\bG}}\bigl(\RH^\bullet(\Gr_{G,\bar{k}},\IC^\lambda),V\bigr)
\end{align}
descends to \(\Gr_G\) whose Satake transform is precisely \(V\), and its induced
function is then \(\fS^V\). At this point, it is safe to simply identify and use
\(V_\lambda\) and \(\RH^\bullet(\Gr_{G,\bar{k}},\IC^\lambda)\) interchangeably,
and we shall always do so when necessary. We make the following definition for
convenience:
\begin{definition}
    The set of \(\lambda\) such that \(\Hom_{\dual{\bG}}(V_\lambda,V)\neq 0\) is
    called the \notion{support of \(V\)}\index{support!of representations of
    \(L\)-groups}.
\end{definition}

\subsection{}
Let \(H\) be an endoscopic group of \(G\). There is no canonical way to obtain
a representation of \(\LD{H}\) by restricting \(V\) because we have to choose an
admissible embedding \(\xi\). Nevertheless, after choosing \(\xi\) we obtain
the restricted representation
\begin{align}
    V_{H,\xi}'\defeq\Res_{\xi(\LD{H})}^{\LD{G}}V,
\end{align}
as well as Satake function \(\fS_{H}^{V_{H,\xi}'}\in\cH_{H,0}\). There is a minor
technical adjustment necessary due to cohomological shifts: for each
\(\dual{\bG}\)-highest weight \(\lambda\) in the support of \(V\) and each
\(\lambda_H\) in the \(\dual{\bH}\)-support of \(V_\lambda\), we modify the
Frobenius action on the corresponding
\(\IC_H^{\lambda_H}\) by
\begin{align}
    (-1)^{\Pair{2\rho}{\lambda}-\Pair{2\rho_H}{\lambda_H}},
\end{align}
and denote the resulting modification of \(\fS_{H}^{V_{H,\xi}'}\) by
\(\fS_{H,\xi}^V\). Similarly, we let the Satake sheaf (resp.~representation)
corresponding to \(\fS_{H,\xi}^V\) by \(\IC_{H,\xi}^V\) (resp.~\(V_{H,\xi}\)).
\begin{definition}
    \label[definition]{def:local_transfer_of_Satake_stuff}
    The function \(\fS_{H,\xi}^V\)
    \nomenclature[\(f^V_H_xi \)]{\(\fS_{H,\xi}^V\)}{the transfer of \(\fS^V\) to \(H\) via \(L\)-embedding \(\xi\)}
    is called the
    \notion{transfer}\index{transfer!of Satake functions} of \(\fS^V\) via
    \(\xi\). Similarly, we call \(\IC_{H,\xi}^V\) (resp.~\(V_{H,\xi}\))
    \nomenclature[\(I{}C^V_H_xi \)]{\(\IC_{H,\xi}^V\)}{the transfer of \(\IC^V\) to \(H\) via \(L\)-embedding \(\xi\)}
    \nomenclature[\(V_H_xi \)]{\(V_{H,\xi}\)}{the restriction of \(\LD{G}\)-representation \(V\) to \(\LD{H}\) via \(L\)-embedding \(\xi\)}
    the transfer of \(\IC^V\) (resp.~\(V\)).
\end{definition}

\subsection{}
Now we are ready to state the fundamental lemma, the main theorem of this book.
We fix an arbitrary admissible embedding \(\xi\colon \LD{H}\to\LD{G}\).

\begin{theorem}
    [Fundamental Lemma]
    \label[theorem]{thm:FL_main}
    For any \(V\in\Rep_{\LD{G}}^\Alg\), we have equality
    \begin{align}
        \Delta_0(\gamma_H,\gamma)\OI_a^\kappa(\fS^V,\dd
        t_v)=\SOI_{a_H}(\fS_{H,\xi}^V,\dd t_v),
    \end{align}
    where \(\Delta_0(\gamma_H,\gamma)\)
    is the transfer factor defined in
    \Cref{sec:Appendix_Definition_of_Transfer_Factors}.
\end{theorem}

\subsection{}
As a final comment of this chapter, \Cref{thm:FL_main} is
stated purely in terms of the groups \(G\) and \(H\), and has nothing to do with
reductive monoids so far. The necessary connection will be made in
\Cref{sec:matching_orbits}. The basic idea is that the monoid \(\FRM\) will
simultaneously encode information about conjugacy class \(\gamma\) and
representation \(V\), and the endoscopic monoid \(\FRM_H\) will do the same for
\(H\). This is already manifested in the construction of \(\FRM_H\) itself,
and will be a recurring theme throughout this book. We believe this is not
special to standard endoscopy but should lead to a general framework for
studying functoriality geometrically.

\chapter{Multiplicative Valuation Strata}%
\label{chap:multiplicative_valuation_strata}

In this chapter we study the \notion{root valuation strata}
\index{root valuation!stratum} in the multiplicative
setting analogous to those in \cite{GKM09} in the Lie algebra setting. The
main result of this section is a codimension formula for the valuation strata,
which will later become a key ingredient for the application of support theorem
to multiplicative Hitchin fibrations. The notations in this chapter are
slightly different from others in order to reduce cluttering.

\section{Arc Spaces of Tori and Congruence Subgroups}%
\label{sec:arc_spaces_of_tori_and_congruence_subgroups}

\subsection{}
Let \(F=\bar{k}\lauser{\pi}\) and \(\cO=\bar{k}\powser{\pi}\) be the rings of Laurent
series and power series over \(\bar{k}\) respectively. Let \(\val_F\) be the
normalized valuation on \(F\) such that \(\val_F(\pi)=1\).  Let \(F_\infty\)
    \nomenclature[\(F_infinity  \)]{\(F_\infty\)}{the maximal tamely ramified
    extension of a geometric local field \(F\) (i.e., the residue field of \(F\) is algebraically closed)}
be the maximal tamely ramified extension of \(F\) inside a fixed algebraic
closure \(\bar{F}\). For each \(l\ge 1\) not divisible by \(p\), we choose
\(\pi^{1/l}\in F_\infty\) and a primitive root of unity \(\zeta_l\)
    \nomenclature[\(zeta_l  \)]{\(\zeta_l\)}{a compatible system of primitive roots of unity of various orders \(l\)}
    \nomenclature[\(pi_l  \)]{\(\pi_l\)}{a compatible system of tamely ramified
    uniformizers \(\pi^{1/l}\) for a given uniformizer \(\pi\) in a local field}
in such a compatible way that \((\pi^{1/(lm)})^m=\pi^{1/l}\) and
\(\zeta_{lm}^m=\zeta_l\). This way one obtains a unique element
\(\tau_\infty\in\Gal(F_\infty/F)\)
    \nomenclature[\(tau_infinity  \)]{\(\tau_\infty\)}{the unique element in \(\Gal(F_\infty/F)\) such that \(\tau_\infty(\pi^{1/l})=\zeta_l\pi^{1/l}\)}
such that \(\tau_\infty(\pi^{1/l})=\zeta_l\pi^{1/l}\). Clearly \(\tau_\infty\) is a
topological generator of \(\Gal(F_\infty/F)\).  Let \(F_l=F[\pi^{1/l}]\),
    \nomenclature[\(F_l \)]{\(F_l\)}{the purely tamely ramified extension of degree \(l\) of a local field \(F\)}
\(\cO_l=\cO[\pi^{1/l}]\),
    \nomenclature[\(O"cal_l \)]{\(\cO_l\)}{the valuation ring of \(F_l\)}
and \(\tau_l\) be the image of \(\tau_\infty\) in \(\Gal(F_l/F)\).

\subsection{}
For \(m\in\bbN\), we let \(\cJ_m=\bar{k}[\pi]/(\pi^{m+1})\) be the ring of
\(m\)-jets over \(\bar{k}\). 
For a scheme \(Y\) over \(F\), let \(\Loop{Y}\) be the \inotion{loop
space}\index{space!loop} of \(Y\). In other words, for any \(\bar{k}\)-algebra \(R\),
\(\Loop{Y}(R)=Y(R\otimes_{\bar{k}} F)\). If \(Y\) is defined over \(\cO\), let
\(\Arc_m{Y}\)
    \nomenclature[\(L{}"bbold_n \)]{\(\Arc_n{Y}\)}{the \(n\)-jet space (\(n\in\bbN\)) of a scheme \(Y\)}
be the \notion{\(m\)-jet space}\index{jet space}\index{space!jet} of \(Y\) so that
\(\Arc_m{Y}(R)=Y(R\otimes_{\bar{k}} \cJ_m)\), and \(\Arc{Y}\defeq\varprojlim\Arc_m{Y}\)
be the \inotion{arc space}\index{space!arc} of \(Y\). If \(Y\) is a \(\bar{k}\)-scheme, we let
\(\Loop{Y}\defeq\Loop{Y_{F}}\), \(\Arc_m{Y}\defeq\Arc_m{Y_{\cO}}\),
and \(\Arc{Y}\defeq\Arc{Y_{\cO}}\).

\subsection{}%
Let \(T\) be a torus over \(\bar{k}\) and \(\La{t}\) its Lie algebra.  Let
\(T_n\) be the congruent subgroup
\begin{align}
    T_n\defeq \ker\left(\Arc{T}\to \Arc_n{T}\right).
\end{align}
Since we are primarily interested in the topological properties, 
we use \(T_n\) for both this group scheme and its \(\bar{k}\)-points.

We have canonical isomorphisms \(T(\cO)\simeq T(\bar{k})\x T_0\), 
and \(T_n/T_{n+1}\simeq \ker\left[T(\cJ_{n+1})\to T(\cJ_{n})\right]\), hence
\begin{align}
    T_n/T_{n+1}\simeq \Set*{\phi\in\Hom_{\bar{k}\hy\mathup{Alg}}(\bar{k}[T],\cJ_{n+1}))\given
    \phi(\lambda)\equiv 1\bmod{\pi^{n+1}},\forall \lambda\in\CharG(T)}.
\end{align}
Therefore, for any \(\phi\in T_n/T_{n+1}\) and \(\lambda\in\CharG(T)\),
\(\phi(\lambda)=1+a(\lambda)\pi^{n+2}\), where
\(a\in\Hom_\bbZ(\CharG(T),\bar{k})\simeq \La{t}\). 
Thus, \(T_n/T_{n+1}\simeq\La{t}\) for
all \(n\ge 0\). Under this isomorphism, for any \(\lambda\in\CharG(T)\), \(t\in
T_n\) and its image \(\bar{t}\in \La{t}\), we have \(\lambda(t)-1\equiv
\dd\lambda(\bar{t})\pi^{n+1}\bmod\pi^{n+2}\).

\subsection{}%
Let \(l\in\bbZ_+\) be coprime to \(p\), and suppose
\(A\) is a cyclic group of order \(l\) acting on \(\Arc{T}\) compatible with
filtrations \(T_n\) and canonical map \(T\to\Arc{T}\).
Note then \(T_n/T_{n+1}\simeq \La{t}\) is a \(\bar{k}\)-linear
representation of \(A\). The following lemma is necessary when we
later consider some twisted forms of \(T(\cO)\).
\begin{lemma}\label[lemma]{lem:cyclic_action_surjective}
    The induced maps on fixed points \(T_n^A\to (T_n/T_{n+1})^A\)
    and \(T(\cO)^A\to T(\bar{k})^A\) are surjective.
\end{lemma}
\begin{proof}
    The claim about \(T(\cO)\to T(\bar{k})\) is trivial as the projection splits.
    View \(T_n/T_{n+1}\) as \(\bar{k}\)-vector space \(\La{t}\).
    Let \(\sigma\in A\) be a generator. Since \(l\) is invertible in \(\bar{k}\),
    we know the \(\bar{k}\)-linear map 
    \begin{align}
        T_n/T_{n+1}&\longto (T_n/T_{n+1})^A\\
        x&\longmapsto \Nm_{\sigma}(x)\defeq x+\sigma(x)+\cdots+\sigma^{l-1}(x)
    \end{align}
     is surjective. Lift \(x\) to \(T_n\), denoted by \(t\), 
    then \(\Nm_\sigma(t)\defeq t\sigma(t)\cdots\sigma^{l-1}(t)\) is
    \(A\)-invariant and maps to \(\Nm_\sigma(x)\).
\end{proof}

\section{Root Valuation Functions and Filtrations}%
\label{sec:root_valuation_functions_and_filtrations}

\subsection{}%
Let \(t\in T(F_\infty)\), and \(\lambda\in\CharG(T)\). Define
\begin{align}
    r_t(\lambda)=\val_F(1-\lambda(t)).
    \nomenclature[\(r_t \)]{\(r_t\)}{the root valuation function induced by \(t\in T(F_\infty)\)}
\end{align}
Then if \(t\in T(\cO_\infty)\), \(r_t(\lambda+\mu)\ge\min\{r_t(\lambda),r_t(\mu)\}\)
and reaches equality if \(r_t(\lambda)\neq r_t(\mu)\). Note that \(r_t\) can
take \(\infty\) as value. 

Using the fixed system of uniformizers \(\pi^{1/l}\), we can decompose
\(T(F_\infty)\) as
\begin{align}
    T(F_\infty)\cong T(\cO_\infty)\x \CoCharG(T)_{\bbZ_{(p)}}\subset
    T(\cO_\infty)\x\CoCharG(T)_\bbQ.
\end{align}
Then one can
uniquely write \(t\in T(F_\infty)\) as a product \(t_0\pi^{\lambda/l}\) for some
\(t_0\in T(\cO_\infty)\), \(\lambda\in\CoCharG(T)\), and positive integer \(l\)
coprime to \(p\).

\subsection{}%
Now suppose \(G\) is a connected reductive group over \(\bar{k}\) and \(T\) is a
maximal torus of \(G\). Let \(t\in T(F_\infty)\), then \(r_t\) induces a function on
the set of roots \(\Roots\) by restriction, still denoted by \(r_t\), called the
\notion{root valuation function}\index{root valuation!function}\index{function!root valuation} induced by \(t\).
We can also define the
\notion{discriminant valuation}\index{discriminant!valuation} of \(t\) by
\begin{align}
    d(t)=\sum_{\alpha\in\Roots}r_t(\alpha)\in\bbQ\cup\Set{\infty},
    \nomenclature[\(d_t \)]{\(d(t)\)}{the non-extended discriminant valuation of \(t\in T(F_\infty)\)}
\end{align}
which is finite if and only if \(t\) is regular semisimple. Note that \(d(t)\)
is none other than the \(F\)-valuation of the discriminant \(\Disc(t)\)
introduced in \Cref{sec:invariant_theory_of_the_group}.

Write \(t=t_0\pi^{\lambda/l}\), then the rational cocharacter
\(\lambda_\AD/l\in\CoCharG(T^\AD)_\bbQ\) can be recovered from \(r_t\) as
follows. Observe that \(r_t(\alpha)=r_t(-\alpha)\) if and only if
\(\Pair{\alpha}{\lambda}=0\), and otherwise
\(\Set{r_t(\alpha),r_t(-\alpha)}=\Set{0,s}\) for some \(s<0\).
Define a new function
\begin{align}
    r_t^-(\alpha)= \begin{cases}
        r_t(\alpha) & r_t(\alpha)<0,\\
        0 & r_t(\alpha)=r_t(-\alpha),\\
        -r_t(-\alpha) & r_t(-\alpha)<0.
    \end{cases}
    \nomenclature[\(r_t_-_star \)]{\(r_t^-,r_t^\star\)}{two modified versions of \(r_t\)}
\end{align}
Then \(r_t^-\) extends to a homomorphism \(\CharG(T)_{\AD,\bbQ}\to \bbQ\). Since 
any set of simple roots form a basis of \(\CharG(T)_{\AD,\bbQ}\), 
this extension to a homomorphism is unique. Thus, \(\lambda_\AD/l\) is
recovered from \(r_t^-\). As a result, let \(\bar{\lambda}\) be the image of
\(\lambda\) in \(\CoCharG(G/G^\Der)\), then the pairs \((r_t,\lambda/l)\) and
\((r_t,\bar{\lambda}/l)\) carry the same amount of information about \(t\).

For later convenience, we also define a modified form \(r_t^\star\) of \(r_t\)
by
\begin{align}
    r_t^\star(\alpha)=\min\Set{r_t(\alpha),r_t(-\alpha)}.
\end{align}
Note that one can recover \(r_t\) from \(r_t^\star\) and \(r_t^-\).
In fact, we can define \(r^-\) and \(r^\star\) 
for any function \(r\) on \(\Roots\) such that
either \(0\le r(\alpha)=r(-\alpha)\le \infty\) or 
\(\Set{r(\alpha),r(-\alpha)}=\Set{0,s}\) for some \(s<0\).

\subsection{}%
To study the pair \((r_t,\bar{\lambda}/l)\) incurred by \(t\in T(F_\infty)\) in general,
we first study those \(t\in T(\cO)\), in which case \(r_t\) takes values in
\(\bbN\cup\Set{\infty}\) and \(\lambda\)-part is trivial.
For any root system \(\Roots\) and 
any functions \(r\colon\Roots\to\bbN\cup\Set{\infty}\), 
we may define filtration (of subsets) of \(\Roots\)
\begin{align}
    \Roots_n(r)\defeq\Set{\alpha\in\Roots\given r(\alpha)\ge n},
    \nomenclature[\(Phi_n_r \)]{\(\Phi_n(r)\)}{the subset
    \(\Set{\alpha\in\Roots\given r(\alpha)\ge n/l}\) (\(n\in\bbZ\)) of roots induced by a function \(r\) valued in \(\bbZ/l\)}
\end{align}
where \(n\in \bbN\cup\Set{\infty}\), and we will simply use \(\Roots_n\) if the
function \(r\) is clear from the context. Clearly this filtration stabilizes
after finite steps and
\begin{align}
    \bigcap_{0\le n<\infty}\Roots_n=\Roots_\infty.
\end{align}
To study this filtration for those \(r=r_t\), 
we first record two theorems concerning root subsystems of \(\Roots\).
\begin{theorem}[Slodowy's criterion, \cite{Sl80}*{Corollary~3.5}]
    A root subsystem \(\Roots'\) 
    of \(\Roots\) comes from a Levi subgroup if and
    only if \(\Roots'\) is \(\bbQ\)-closed in \(\Roots\), in other words,
    \(\bbQ\Roots'\cap\Roots=\Roots'\).
\end{theorem}

\begin{theorem}[Deriziotis's criterion, \cite{Hu95}*{\S~2.15}]
    \label[theorem]{thm:Deriziotis_criterion}
    A root subsystem \(\Roots'\) 
    of \(\Roots\) comes from a pseudo-Levi subgroup (i.e. the connected
    centralizer of a semisimple element of \(G\)) if and
    only if after conjugation by \(W\), 
    \(\Roots'\) has a basis which  is a \emph{proper}
    subset of the simple \emph{affine} roots
    of \(\Roots\), in other words, simple roots together with negative of the
    highest root.
\end{theorem}

\begin{lemma}\label[lemma]{lem:Qclosed_for_t_T0}
    Let \(t\in T_0\), then \(\Roots_n(r_t)\) is \(\bbQ\)-closed.
\end{lemma}
\begin{proof}
    Clearly \(\Roots_n\) is \(\bbZ\)-closed for each \(n\). Following the proof
    in \cite{GKM09}*{14.1.1}, if \(\alpha\in\Roots\) is a \(\bbQ\)-combination of roots
    in \(\Roots_n\), then \(d\alpha\in\bbZ\Roots_n\) for some \(d\) coprime
    to \(p\) (here only the fact that \(\Roots_n\) is \(\bbZ\)-closed
    is needed). Thus, \(r_t(d\alpha)\ge n\). But since \(t\in T_0\) and
    \(d\) is coprime to \(p\), we must have \(r_t(\alpha)\ge n\) as well,
    hence the lemma.
\end{proof}

The following corollary is not used in subsequent parts of this book, but
highlights the subtle difference between the multiplicative and additive cases
(cf.~\cite{GKM09}*{3.4.1}).
\begin{corollary}
    There is a Zariski-dense open subset of \(T(\bar{k})\) such that the conclusion 
    in \Cref{lem:Qclosed_for_t_T0} holds for its preimage in
    \(T(\cO)\).
\end{corollary}
\begin{proof}
    The collection of possible \(d\) (chosen to be as small as possible each
    time) appearing in the proof of
    \Cref{lem:Qclosed_for_t_T0} is a finite set \(S\). Let \(U\subset
    T(\bar{k})\) be the complement of those \(t\in T(\bar{k})\) such that
    \(1\neq \alpha(t)\in\mu_d\) (\(\mu_d\) is the set of \(d\)-th roots of
    unity) for some \(\alpha\in\Roots\) and some \(d\in
    S\). Then it is clear that the same proof goes through for elements in the
    preimage of \(U\).
\end{proof}

\begin{remark}
    It is not true even in characteristic \(0\)
    that \Cref{lem:Qclosed_for_t_T0} holds for all \(t\in
    T(\cO)\). For example, suppose \(G\) is of type \(\TypeB_2\) or 
    \(\TypeC_2\), and let \(t\) be an element such that \(\alpha(t)=-1+\pi\) for
    both positive short root \(\alpha\), then \(\Roots_1\) for \(r_t\) is
    exactly the set of all long roots, which is not \(\bbQ\)-closed.
\end{remark}

\begin{proposition}\label[proposition]{prop:characterize_rt_by_Roots_n}
    Let \(t\in T(\cO)\). Then \(\Roots_1(r_t)\) is a root subsystem of 
    pseudo-Levi type, and \(\Roots_n\) is \(\bbQ\)-closed in \(\Roots_1\)
    for all \(n\ge 1\).
\end{proposition}
\begin{proof}
    Write \(t=xt_0\), where \(x\in T(\bar{k})\) and \(t_0\in T_0\). Then
    \begin{align}
        \Roots_1(r_t)=\Set*{\alpha\in\Roots\given\alpha(x)=1},
    \end{align}
    which is the same as the roots of the connected centralizer of \(x\) in
    \(G\), hence the first claim.

    For the second claim, note that for \(n\ge 1\), we have that
    \begin{align}
        \Roots_n(r_t)=\Roots_n(r_{t_0})\cap
        \Roots_1(r_t)\subset\Roots_n(r_{t_0}).
    \end{align}
    Since \(\Roots_n(r_{t_0})\) is \(\bbQ\)-closed in \(\Roots\) by
    \Cref{lem:Qclosed_for_t_T0}, we have that
    \begin{align}
        \Roots_n(r_t)
        &\subset \bbQ\Roots_n(r_t)\cap\Roots_1(r_t)\\
        &\subset \bbQ\Roots_n(r_{t_0})\cap\Roots_1(r_t)\\
        &=\bbQ\Roots_n(r_{t_0})\cap[\Roots\cap\Roots_1(r_t)]\\
        &=[\bbQ\Roots_n(r_{t_0})\cap\Roots]\cap\Roots_1(r_t)\\
        &=\Roots_n(r_{t_0})\cap\Roots_1(r_t)\\
        &=\Roots_n(r_t).
    \end{align}
    Therefore, every inclusion is an equality, and in particular
    \(\Roots_n(r_t)=\bbQ\Roots_n(r_t)\cap\Roots_1(r_t)\), as claimed.
\end{proof}

\subsection{}%
Let \(L\subset G\) be a Levi subgroup containing \(T\), and
\(\Roots_L\subset \Roots\) be its root subsystem. Then 
\(\Roots_L\) is \(\bbQ\)-closed by Slodowy's criterion.

\begin{lemma}\label[lemma]{lem:find_t_for_r_two_step}
    For any \(0\le n< \infty\), we can find \(t_{n}\in T_n\) such that
    \begin{align}
        r_{t_{n}}(\alpha)= \begin{cases}
            n+1 & \alpha\in\Roots\setminus \Roots_L,\\
            \infty & \alpha\in \Roots_L.
        \end{cases}
    \end{align}
\end{lemma}
\begin{proof}
    Let \(Z_L\) be the \emph{connected} center of \(L\).
    Let \(X\in \La{z}_L\defeq\Lie(Z_L)\simeq
    (Z_L)_{n}/(Z_L)_{n+1}\) be an element that \(\alpha(X)\neq 0\) for all
    \(\alpha\in\Roots\setminus\Roots_L\). This \(X\) exists since \(\La{z}_L\)
    is exactly the intersection of kernels of \(\alpha\in \Roots_L\) and no
    other root vanishes identically on \(\La{z}_L\) (cf.~\cite{GKM09}*{14.2}).
    Lift \(X\) to \(t_{n}\in (Z_L)_{n}\), and we are done.
\end{proof}

\begin{lemma}\label[lemma]{lem:find_t_for_r_full}
    Suppose we have a function \(r\colon\Roots\to \bbZ_+\cup\Set{\infty}\) and
    the associated filtration \(\Roots_n\).
    If \(\Roots_n\) is a \(\bbQ\)-closed root subsystem in \(\Roots\) 
    for each \(1\le n\le \infty\),
    then we can find \(t\in T_0\) such that \(r_t=r\).
\end{lemma}
\begin{proof}
    For each (finite) \(i\ge 0\), let \(L_{i+1}\) 
    be the Levi type subgroup with root system
    \(\Roots_{i+i}\), and let \(t_{i}\) be as in
    \Cref{lem:find_t_for_r_two_step} for pair \(L_{i+2}\subset L_{i+1}\) and \(n=i\).
    Then
    \begin{align}
        r_{t_{i}}(\alpha)= \begin{cases}
            \text{something}\ge i+1 &\alpha\in\Roots\setminus\Roots_{i+1},\\
            i+1 & \alpha\in\Roots_{i+1}\setminus\Roots_{i+2},\\
            \infty & \alpha\in \Roots_{i+2}.
        \end{cases}
    \end{align}
    Thus, for a fixed \(i\ge 0\), 
    and \(\alpha\in \Roots_{i+1}\setminus\Roots_{i+2}\), 
    we have that
    \begin{align}
        r_{t_{j}}(\alpha)= \begin{cases}
            \text{something}\ge j+1 & j>i,\\
            j+1 & j= i,\\
            \infty & j<i.
        \end{cases}
    \end{align}
    Let \(t=\prod_{i=0}^\infty t_i\), which necessarily converges as \(t_i\in
    T_i\).  Since for
    \(\alpha\in\Roots_{i+1}\setminus \Roots_{i+2}\), \(r_{t_j}(\alpha)\) reaches
    unique minimum at \(j=i\), 
    we have that \(r_t(\alpha)=r_{t_{i}}(\alpha)=i+1\), and for
    \(\alpha\in\Roots_\infty\), it is clear that \(r_t(\alpha)=\infty\),
    as desired.
\end{proof}

\begin{theorem}\label[theorem]{thm:characterize_rt_by_filtration}
    Suppose we have a function \(r\colon\Roots\to \bbN\cup\Set{\infty}\) and
    the associated filtration \(\Roots_n\). Then \(r=r_t\) for some \(t\in
    T(\cO)\) if and only if 
    \(\Roots_1\subset\Roots\) is a root subsystem of pseudo-Levi type and
    \(\Roots_n\) is a \(\bbQ\)-closed subsystem in \(\Roots_1\)
    for each \(1\le n\le \infty\). Moreover, \(t\in T(\cO)\) is such that
    \(r_t=r\) if and only if we can write \(t=x\prod_{n=0}^\infty t_n\)
    where \(x\in T(\bar{k})\) and \(t_n\in T_n\) are such that \(\Roots_1\) is the set
    of roots taking trivial value on \(x\), and
    \begin{align}\label{eqn:decompose_t_property}
        r_{t_n}(\alpha)= \begin{cases}
            \text{something}\ge n+1 & \alpha\in\Roots_1\setminus\Roots_{n+1},\\
            n+1 & \alpha\in\Roots_{n+1}\setminus \Roots_{n+2},\\
            \infty & \alpha\in\Roots_{n+2}.
        \end{cases}
    \end{align}
\end{theorem}
\begin{proof}
    For the first claim,
    the ``only if'' part is proved in
    \Cref{prop:characterize_rt_by_Roots_n}. Now for the ``if'' part.
    Let \(x\in T(\bar{k})\) be such that \(\Roots_1\) is the root system of its
    connected centralizer \(L\), which is a connected reductive group
    with maximal torus \(T\). Apply \Cref{lem:find_t_for_r_full} to \(L\)
    and restricted function \(r|_{\Roots_1}\), then one can find \(t_0\in T_0\)
    such that \(r_{t_0}|_{\Roots_1}=r|_{\Roots_1}\). This means that \(t=xt_0\)
    is the desired element.

    For the second claim, the only nontrivial part is the ``only if'' part.
    The existence of \(x\) is also clear, so we may replace \(t\) by \(x^{-1}t\)
    and \(G\) by the connected centralizer of \(x\),
    and in turn we may assume \(t\) is already in \(T_0\). 
    Let \(L_n\) be the Levi subgroup of \(G\) corresponding to
    \(\Roots_n\), \(Z_{L_n}\) its connected center, and
    \(\La{z}_{L_n}=\Lie{Z_{L_n}}\).
    Let \(X_0\) be the image of \(t\) in \(\La{t}\), then \(\alpha(X_0)=0\) for
    all \(\alpha\in\Roots_2\), so we can lift it to \(t_0\in Z_{L_2,0}\) (again,
    the subscript \(0\) here means the congruent subgroup, not neutral
    component), then
    \(t_0\) satisfies \eqref{eqn:decompose_t_property}. Moreover, let
    \(\Roots_n'=\Roots_n(r_{tt_0^{-1}})\), then we have that
    \(\Roots_2'=\Roots\), \(\Roots_n'\supset\Roots_n\) for all \(n\), and 
    \(\Roots_n'\setminus\Roots_{n+1}'\supset\Roots_n\setminus\Roots_{n+1}\)
    for \(n\ge 2\). Doing the same argument
    for \(tt_0^{-1}\) but with every index
    increased by \(1\), we have some \(t_1\in T_1\) satisfying
    \eqref{eqn:decompose_t_property} (\textit{a priori} only for \(\Roots_n'\),
    but as
    \(\Roots_n'\setminus\Roots_{n+1}'\supset\Roots_n\setminus\Roots_{n+1}\),
    same is true for \(\Roots_n\)). Continuing this process, we see that the
    infinite product
    \begin{align}
        t\prod_{n=0}^\infty t_n^{-1}
    \end{align}
    converges to some \(\cO\)-point \(z\) in the center of \(G\). 
    Absorb \(z\) into \(t_0\) and we are done.
\end{proof}

\subsection{}%
Next we consider \(t\in T(F)\). Then \(t=t_0\pi^\lambda\) for some \(t_0\in
T(\cO)\) and \(\lambda\in\CoCharG(T)\). We form the functions \(r_t^-\) and
\(r_t^\star\) as in
\Cref{sec:root_valuation_functions_and_filtrations}.
In this case, let \(\Roots_n\defeq\Roots_n(r_t^\star)\) for \(0\le n\le \infty\),
then we have filtration
\begin{align}
    \Roots\supset \Roots_0\supset\cdots\supset
    \Roots_n\supset\cdots\supset\Roots_\infty\supset \emptyset.
\end{align}
Since \(\Roots_0\) is the set of roots that vanish on \(\lambda\), it is in
particular \(\bbQ\)-closed in \(\Roots\), and thus is a root subsystem of Levi
type.  Then for each \(0\le n\le\infty\), it is clear that
\begin{align}\label{eqn:Roots_n_and_Roots_n_t0_for_T_F}
    \Roots_n=\Roots_0\cap \Roots_n(r_{t_0}).
\end{align}
So \(\Roots_1\) is a root subsystem of pseudo-Levi type in \(\Roots_0\), and
for each \(1\le n\le \infty\),
\(\Roots_n\) is \(\bbQ\)-closed in \(\Roots_1\). To summarize it, we have
\begin{align}\label{eqn:filtration_for_T_F}
    \Roots\stackrel{\text{Levi}}{\supset}\Roots_0
    \stackrel{\text{pseudo-Levi}}{\supset}\Roots_1
    \stackrel{\text{Levi}}{\supset}\Roots_{1\le n\le \infty}.
\end{align}

On the other hand, given a function \(r\colon\Roots\to\bbZ\cup\Set{\infty}\), a
necessary condition for \(r\) to be equal to some \(r_t\) is that 
\(r^-\) and \(r^\star\) as in
\Cref{sec:root_valuation_functions_and_filtrations} are well-defined.
Suppose \(r\) is such a function.
\begin{theorem}\label[theorem]{thm:characterize_rt_by_filtration_T_F}
    There exists \(t\in T(F)\) such that \(r=r_t\) if and only if \(r^-\)
    extends to a homomorphism \(\CharG(T)\to\bbZ\), and the filtration
    \(\Roots_n\) satisfies \eqref{eqn:filtration_for_T_F}.
\end{theorem}
\begin{proof}
    The ``only if'' part is proved by discussions above and in
    \Cref{sec:root_valuation_functions_and_filtrations}. We now prove the ``if''
    part. Let \(\lambda\in\CoCharG(T)\) be one of extensions of \(r^-\) (which
    is not unique unless \(G\) is semisimple). By
    \Cref{thm:characterize_rt_by_filtration}, we can find \(t_0\in
    T(\cO)\) such that \eqref{eqn:Roots_n_and_Roots_n_t0_for_T_F} holds.
    Let \(t=\pi^\lambda t_0\) and we are done.
\end{proof}

\subsection{}%

As discussed before, we have associated pair \((\bar{\lambda}_t,r_t)\) for any \(t\in
T(F)\) where \(\bar{\lambda}_t\in\CoCharG(G/G^\Der)\). 
Let \(S_G^1\) be the set of all
pairs \((\bar{\lambda},r)\) where \(\bar{\lambda}\in\CoCharG(G/G^\Der)\) and
\(r\colon\Roots\to \bbQ\cup\Set{\infty}\). 
Then we have a partition of \(T(F)\) by
\begin{align}
    T(F)=\coprod_{(\bar{\lambda},r)\in S_G^1}T(F)_{(\bar{\lambda},r)},
\end{align}
where \(T(F)_{(\bar{\lambda},r)}\) is the set of all \(t\in T(F)\) such that
\((\bar{\lambda}_t,r_t)=(\bar{\lambda},r)\). Then by
\Cref{thm:characterize_rt_by_filtration_T_F}, \(T(F)_{(\bar{\lambda},r)}\)
is non-empty if and only if the following conditions are satisfied:
\begin{enumerate}
    \item \(r\) takes values in \(\bbZ\cup\Set{\infty}\);
    \item \(r^-\) and \(r^\star\) are defined;
    \item \(r^-\)  extends to (necessarily unique)
        \(\lambda_\AD\in\CoCharG(T^\AD)\) such that
        \begin{align}
            (\bar{\lambda},\lambda_\AD)\in\CoCharG(G/G^\Der)\x\CoCharG(T^\AD)
        \end{align}
        lies in the image of \(\CoCharG(T)\);
    \item the filtration \(\Roots_n(r^\star)\) satisfies
            \eqref{eqn:filtration_for_T_F}.
\end{enumerate}
For convenience, we will call \(\bar{\lambda}\) (resp.~\(\lambda_\AD\)) 
the central (resp.~adjoint) component of \(\lambda\).

\subsection{}%
We now study \(T(F_\infty)\) in general. Let \(t\in T(F_\infty)\), then \(t\in
T(F_l)\) for some \(l\). We
know that the image of \(\tau_\infty\) in the cyclic group \(\Gal(F_l/F)\) is a
generator. Suppose there is some \(w\in W\) with order \(l\). 
Then we know that \(r_t\) takes values in \(\bbZ/l\cup\Set{\infty}\), and 
we are mostly interested in the case when \(\tau_\infty w(t)=t\) (note that
the Galois action commutes with the Weyl group action). To that end we define
\begin{align}
    T_w(F)\defeq\Set*{t\in T(F_l)\given \tau_\infty w(t)=t}.
    \nomenclature[\(T_w_F \)]{\(T_w(F)\)}{the subset of \(T(F_\infty)\) (\(T\) a
    maximal torus of a reductive group \(G\)) fixed by \(\tau_\infty w\) with \(w\) a Weyl element}
\end{align}
Then we immediately have that if \(t=t'\pi^{\lambda_t/l}\),
then \(\lambda_t\) is fixed by \(w\), and \(t'\in T(\cO_l)\) is such that
\begin{align}
    \tau_\infty w(t')\zeta_l^{\lambda_t}=t'.
\end{align}
Let \(T_{F_l,n}\) (\(n\in\bbN\)) be the congruent subgroup of \(T(\cO_l)\) with
\(F\) replaced by \(F_l\), then we can write canonically \(t'=xt_0\) for \(x\in
T(\bar{k})\) and \(t_0\in T_{F_l,0}\). Since \(\zeta_l\in \bar{k}\), we must have 
\begin{align}\label{eqn:Tw_eigen_condition}
    \tau_\infty w(t_0)=t_0\text{ and } (1-w)(x)=\zeta_l^{\lambda_t}.
\end{align}
An immediate observation is that \(\bar{\lambda}_t/l\in\CoCharG(G/G^\Der)\).

As in the \(F\)-rational case, we can define \(r_t^-\) and \(r_t^\star\), as
well as filtration
\begin{align}
    \Roots_n=\Roots_n(r_t^\star)\defeq
    \Set*{\alpha\in\Roots\given r_t^\star(\alpha)\ge \frac{n}{l}}, (0\le
    n\le\infty).
\end{align}
By the same argument with \(F\) replaced by \(F_l\), we know that \(\Roots_n\)
satisfies \eqref{eqn:Roots_n_and_Roots_n_t0_for_T_F}.
We also have that \(r_t^-\) and \(\lambda_{t,\AD}/l\) mutually determines each
other. Therefore, \(t\) induces pair \((\bar{\lambda}_t/l,r_t)\) where the first 
component is an integral cocharacter of \(G/G^\Der\).

\subsection{}%
To fully characterize what pairs \((\bar{\lambda}/l,r)\) arises in this way,
again we first assume \(\lambda_t=0\), in other words, \(t\in T(\cO_l)\). In
fact, it is better to start with \(t\in T_{F_l,0}\),
in which case \(\Roots=\Roots_0=\Roots_1\).

Given \(t\in T_{F_l,0}\), by \Cref{thm:characterize_rt_by_filtration}, we
may write \(t\) as the infinite product of \(t_n\in T_{F_l,n}\) such that
\begin{align}\label{eqn:decompose_t_property_twisted}
    r_{t_n}(\alpha)= \begin{cases}
        \text{something}\ge \frac{n+1}{l} & \alpha\in\Roots_1\setminus\Roots_{n+1},\\
        \frac{n+1}{l} & \alpha\in\Roots_{n+1}\setminus\Roots_{n+2},\\
        \infty & \alpha\in\Roots_{n+2}.
    \end{cases}
\end{align}
Using the same notations in the proof of
\Cref{thm:characterize_rt_by_filtration},
since \(\tau_\infty w(t)=t\), \(\tau_\infty w\) preserves each \(\Roots_n\), and
we see that \(X_0\in \La{z}_{L_2}^{\tau_\infty w}\) (here the action of
\(\tau_\infty\) on \(\La{t}\), when viewed as quotient \(T_n/T_{n+1}\),
is multiplication by \(\zeta_l^{n+1}\), and the action of \(w\) is the usual one).
By \Cref{lem:cyclic_action_surjective}, we can require
the lift \(t_0\) to be contained in \(Z_{L_2,0}^{\tau_\infty w}\), hence
\(tt_0^{-1}\) is also fixed by \(\tau_\infty w\).
Inductively doing so in the construction of each \(t_n\),
we can make each \(t_n\) fixed by \(\tau_\infty w\).

\begin{definition}
    Suppose \(w\in W\) is of order \(l\) and 
    acting on a \(\bar{k}\)-vector space \(V\).
    For any integer \(i\), define \(V(w,i)\subset V\) 
    to be the maximal subspace where \(w\) acts as \(\zeta_l^{-i}\).
\end{definition}

\begin{lemma}
    Suppose we have a function \(r\colon\Roots\to \bbZ_+/l\cup\Set{\infty}\) 
    and the
    associated filtration \(\Roots_n\). Then we may find \(t\in T_{w,0}\defeq
    T_w(\cO)\cap T_{F_l,0}\) such
    that \(r_t=r\) if and only if \(\Roots_n\) is \(\bbQ\)-closed in
    \(\Roots\) for all \(1\le n\le\infty\), and the set
    \begin{align}\label{eqn:eigencondition_twisted_strata}
        \Set*{X_n\in\La{t}(w,n+1)\given 
            \begin{array}{l}
            \alpha(X_n)=0,\forall \alpha\in \Roots_{n+2},\\
            \alpha(X_n)\neq 0,\forall \alpha\in\Roots_{n+1}\setminus\Roots_{n+2}
            \end{array}
        }
    \end{align}
    is non-empty for all \(0\le n<\infty\).
\end{lemma}
\begin{proof}
    For the ``only if'' direction, simply choose \(X_n\) to be the image of
    \(t_n\) in the discussion above. For the ``if'' direction, by
    \Cref{lem:cyclic_action_surjective}, we can lift \(X_n\) to \(t_n\in
    T_{w,n}\defeq T_w(\cO)\cap T_{F_l,n}\), and then \(t=\prod_{n=0}^\infty t_n\)
    would be as desired.
\end{proof}
\begin{remark}
    Note that the condition \eqref{eqn:eigencondition_twisted_strata} being
    non-empty for all \(n\) automatically implies that each \(\Roots_n\) (hence
    also \(r\)) is
    preserved by \(w\), hence also by \(\tau_\infty w\).
\end{remark}

\begin{corollary}\label[corollary]{cor:characterizationn_of_Twr}
    Suppose we have a function \(r\colon\Roots\to \bbN/l\cup\Set{\infty}\) 
    and the
    associated filtration \(\Roots_n\). Then we may find \(t\in T_{w}(\cO)\)
    such that \(r_t=r\) if and only if the following conditions are satisfied
    \begin{enumerate}
        \item \(\Roots_1\subset\Roots\) is a root subsystem of pseudo-Levi type
            with corresponding reductive subgroup \(H\), such that the set
            \begin{align}
                \Set*{x\in T(\bar{k})^w\given H\text{ is the connected centralizer of }x}
            \end{align}
            is non-empty (in particular \(\Roots_1\) is preserved by \(w\));
        \item \(\Roots_n\) is \(\bbQ\)-closed in
            \(\Roots_1\) for all \(1\le n\le\infty\);
        \item the set \eqref{eqn:eigencondition_twisted_strata}
            is non-empty for all \(0\le n<\infty\).
    \end{enumerate}
\end{corollary}

\subsection{}%
Now we consider general pairs \((\bar{\lambda}/l,r)\), where \(\bar{\lambda}/l\)
is an (integral) element of \(\CoCharG(G/G^\Der)\), and
\(r\colon\Roots\to\bbQ\cup\Set{\infty}\). In other words,
\((\bar{\lambda}/l,r)\in S_G^1\).
Similar to \(F\)-rational case,
we stratify \(T_w(F)\) by such pairs, denoted by
\(T_w(F)_{(\bar{\lambda}/l,r)}\). 
    \nomenclature[\(T_w_F_lambda_r \)]{\(T_w(F)_{(\bar{\lambda}/l,r)}\)}{the
    subset of \(T_w(F)\) with boundary cocharacter \(\bar{\lambda}/l\) and root valuation function \(r\)}

From previous discussions, we see that 
\(T_w(F)_{(\bar{\lambda}/l,r)}\) is
non-empty if and only if the following conditions are satisfied:
\begin{enumerate}
    \item \(r\) takes values in \(\bbZ/l\cup\Set{\infty}\);
    \item \(r^-\) and \(r^\star\) are defined, and \(r^-\) extends to an
        element \(\lambda_\AD/l\in\CoCharG(T^\AD)^w/l\), such that
        \(\bar{\lambda}\) and \(\lambda_\AD\) are the central and
        adjoint components of some \(\lambda\in\CoCharG(T)\) respectively;
    \item let \(L_\lambda\) be the Levi subgroup of \(G\) determined by
        \(\lambda\) (or equivalently, \(\Roots_0\)),
        then \(\Roots_1\) is of pseudo-Levi type in \(\Roots_0\) corresponding
        to subgroup \(H<L_\lambda\), such that the set 
        \begin{align}\label{eqn:eigncentral_condition_twisted_strata}
            \Set*{x\in
            (1-w)^{-1}(\zeta_l^\lambda)\given H\text{ is the connected
            centralizer of }x\text{ in }L_\lambda} 
        \end{align}
        is non-empty;
    \item \(\Roots_n\) is \(\bbQ\)-closed in \(\Roots_1\) for all \(1\le
        n\le\infty\), and the set \eqref{eqn:eigencondition_twisted_strata}
            is non-empty for all \(0\le n<\infty\).
\end{enumerate}
Note that these conditions automatically imply that the filtration \(\Roots_n\)
is preserved by \(w\).
Moreover, given \(t\in T_w(F)\), we have that 
\(t\in T_w(F)_{(\bar{\lambda}/l,r)}\), if and only if we can write
\begin{align}
    t=\pi^{\lambda/l}x\prod_{n=0}^\infty t_n,
\end{align}
 where \(x\) is as in
 \eqref{eqn:eigncentral_condition_twisted_strata}, and \(t_n\in T_{F_l,n}\) is fixed
 by \(\tau_\infty w\) such that \(r_{t_n}\) satisfies
\eqref{eqn:decompose_t_property_twisted} for all \(0\le n<\infty\). 

\section{Cylinders in Reductive Monoids}%
\label{sec:cylinders_in_reductive_monoids}

\subsection{}%
We briefly review some facts about the arc spaces of a smooth affine scheme
\(X\) over \(\cO\).  Since \(X\) is \(\cO\)-smooth, every jet scheme of \(X\) is
smooth over \(\bar{k}\), and each consecutive map \(\Arc_{n+1}{X}\to\Arc_n{X}\)
is an affine space bundle of relative dimension \(\dim_{\cO}{X}\). 

\begin{definition}
    A subset \(Z\) of \(\Arc{X}\) is called
    \notion{\(n\)-admissible}\index{subset!admissible, of arcs} or an
    \notion{\(n\)-cylinder}\index{cylinder!of arcs} if it is the preimage of a
    constructible subset of \(\Arc_n{X}\) for some \(n\ge 0\). A cylinder is
    called \notion{open}\index{cylinder!open}
    (resp.~\notion{closed}\index{cylinder!closed},
    \notion{locally closed}\index{cylinder!locally closed}) if it is the preimage of some
    open (resp.~closed, locally closed) subset of \(\Arc_n{X}\).
\end{definition}

\begin{definition}
    Let \(Z\subset \Arc{X}\) be an \(n\)-cylinder. Then the
    \notion{codimension}\index{codimension!of cylinder}
    of \(Z\) in \(\Arc{X}\) is defined as the codimension of \(Z_m\) in
    \(\Arc_m{X}\) for any (and every) \(m\ge n\).
\end{definition}

\begin{definition}
    Suppose \(Y\) is a locally closed \(n\)-cylinder of \(\Arc{X}\).
    Then we call \(Y\) \notion{non-singular}\index{cylinder!non-singular}
    if for any (and every) \(m\ge n\) the reduced subscheme \(Y_m\in
    \Arc_m{X}\) is non-singular. We call a map \(g\colon Y\to Z\) of
    any locally closed \(n\)-cylinders \notion{smooth}\index{cylinder!smooth map} if
    \(g_m\) is smooth for all \(m\ge n\).
\end{definition}

\subsection{}
Recall the universal monoid \(\FRM=\Env(G)\) of a semisimple and simply-connected
group \(G=\bG^\SC\), and the closure \(\FRT=\FRT_\FRM\) of the extended maximal torus \(T_+\).
We know \(\bar{k}[\FRT]\) is spanned by a saturated, strictly convex 
cone \(\sC^*_{\FRT}\subset
\CharG(T_+)\), whose dual \(\sC_{\FRT}\subset\CoCharG(T_+)\) determines a
stratification of non-degenerate arcs
\begin{align}
    \label{eqn:Arc-FRT-strat}
    \FRT(\cO)\cap
    T_+(F)=\coprod_{\lambda\in\sC_{\FRT}}\Arc{\FRT}^\lambda(\bar{k}),
\end{align}
where \(\Arc{\FRT}^\lambda=\pi^\lambda\Arc{T_+}\). Note that
\(\Arc{\FRT}^\lambda\) is isomorphic to \(\Arc{T_+}\) as abstract
\(\bar{k}\)-schemes.

For \(w\in W\) be of order \(l\) and \(\lambda\) fixed by \(w\), we also have
the \(w\)-twisted form of \(\Arc{\FRT}^\lambda\) defined in an
obvious way:
\begin{align}
    \FRT_w^{\lambda/l}\defeq
    \left(\Res_{\cO_l/\cO}\pi^{\lambda/l}T_{+,\cO_l}\right)^{\tau_\infty
    w},
    \nomenclature[\(T"frak_w_lambda \)]{\(\FRT_w^{\lambda/l}\)}{the subset of
    \(\FRT_\FRM(\cO_l)\) with boundary \(\lambda/l\) and fixed by \(\tau_\infty w\), viewed as a \(\cO\)-space}
\end{align}
and the \(\bar{k}\)-points of \(\Arc{\FRT_w^{\lambda/l}}\) are simply those \(t\in
\pi^{\lambda/l}T_+(\cO_l)\) that are fixed by \(\tau_\infty w\).
Since \(l\) is coprime to \(\Char(k)\), the scheme
\(\FRT_w^{\lambda/l}\) is \(\cO\)-smooth, hence notions like cylinders and their
codimensions make sense for its arc space.

It is clear that any non-empty stratum \(T_{+,w}(F)_{(\bar{\lambda}/l,r)}\) as
in \Cref{sec:root_valuation_functions_and_filtrations} is
entirely contained in a unique \(\Arc{\FRT}_w^{\lambda/l}(\bar{k})\) if
\(\lambda\) is contained in the cone \(\sC_{\FRT}\). In such case, we
denote the stratum by
\begin{align}
    \FRT_w(\cO)_{(\bar{\lambda}/l,r)}\defeq T_{+,w}(F)_{(\bar{\lambda}/l,r)}.
    \nomenclature[\(T"frak_w_lambda_r \)]{\(\FRT_w(\cO)_{(\bar{\lambda}/l,r)}\)}{the
    subset of \(\FRT_w^{\lambda/l}(\cO)\) with root valuation function \(r\)}
\end{align}

\subsection{}%
Choose a geometric point \(x_+\in T_+(F^s)^\rss\) lying over some
\(a\in\FRC_{\FRM}(\cO)\), then it induces a homomorphism
\begin{align}
    \rho_a\colon \Gamma_F=\Gal(F^s/F)\longto W,
\end{align}
the conjugacy class of whose image depends only on \(a\). Since \(a\) is tamely
ramified (by our assumption on \(\Char(k)\)), one sees that the image \(W_a\) is
cyclic and generated by the image of \(\tau_\infty\), denoted by \(w_a\).
 Let \(\gamma\in\FRM^\x(F)\) be a point over \(a\), and we define the
 \inotion{ramification index} \(c\) of \(\gamma\) or \(a\) to be
\begin{align}
    \label{eqn:def_ramification_index}
    c(\gamma)=c(a)\defeq \dim{T_+}-\dim{T_+^{w_a}}.
    \nomenclature[\(c_a \)]{\(c(a), c_v(a)\)}{the ramification index (at place \(v\)) of a local conjugacy class \(a\in\FRC_\FRM(\cO_v)\)}
\end{align}
It is in fact equal to the difference
\(\rk_{\bar{F}}(G_{\gamma})-\rk_F(G_{\gamma})\), where \(\rk_F\) means the
\(F\)-split rank of the centralizer \(G_{\gamma}\). So \(c(\gamma)\) 
can be defined for arbitrary reductive group \(G\) (not necessarily semisimple
and simply-connected) and \(\gamma\in G(F)^\rss\).

Recall we have extended discriminant function \({\Disc_+}\in k[\FRC_{\FRM}]\).
The extended discriminant valuation of \(a\) is defined by
\(d_+(a)\defeq\val_F(\Disc_+(a))\).
    \nomenclature[\(d_+_a \)]{\(d_+(a), d_{v+}(a)\)}{the extended discriminant
    valuation (at place \(v\)) of a local conjugacy class \(a\in\FRC_\FRM(\cO_v)\)}
We define the \notion{local
\(\delta\)-invariant}\index{\(\delta\)-!invariant, local}\index{local!\(\delta\)-invariant}
of \(a\) to be
\begin{align}
    \delta(a)\defeq\frac{d_+(a)-c(a)}{2},
    \nomenclature[\(delta_v_a \)]{\(\delta(a), \delta_{v}(a)\)}{the local
        \(\delta\)-invariant (at place \(v\)) of a local conjugacy class \(a\in\FRC_\FRM^\rss(\cO_v)\)}
\end{align}
which is in fact necessarily in \(\bbN\) (due to the fact that it is the
dimension of certain multiplicative affine Springer fiber; see
\Cref{thm:GASF_dimension_formula}).

\subsection{}%
Let us clarify some relations between
Levi subgroups and pseudo-Levi subgroups in an arbitrary connected reductive
group \(G\) containing maximal torus \(T\). 
Suppose \(L_1\) is a pseudo-Levi subgroup containing 
\(T\) with root subsystem \(\Roots_1\), and let \(\Roots'\) be its \(\bbQ\)-closure in
\(\Roots\), which gives a Levi subgroup \(L\) of \(G\) containing \(L_1\). 
Let \(Z(L_1)\) be the
center of \(L_1\), and \(Z(L_1)_0\) the neutral component. Then the centralizer
of \(Z(L_1)_0\) in \(G\) is \(L\). 
Let \(x\in Z(L_1)\) be such that the connected centralizer
is exactly \(L_1\) and \(t\in Z(L_1)_0\). Then \(\alpha(xt)\neq 1\) for all
\(\alpha\not\in \Roots_1\) if \(t\) is general enough. This means that the
closure of the set
\begin{align}
    \Set*{x\in T\given (G_x)_0=L_1}
\end{align}
in \(T\) is the union of some connected components of \(Z(L_1)\). In particular,
it has the same dimension as \(Z(L_1)\). We denote the set of such components by
\(\pi_0^\circ(Z(L_1))\). 

\begin{proposition}\label[proposition]{prop:codim_Tplus_r}
    Suppose \(r\) takes values in \(\bbN\) and \(T_+(\cO)_r\neq\emptyset\). 
    Let \(L_n\) be the (pseudo-)Levi subgroup
    of \(G_+\) determined by \(\Roots_n(r)\) and \(\La{z}_n=\Lie(Z(L_n))\).
    Then we have the following:
    \begin{enumerate}
        \item The closure of \(T_+(\cO)_r\) in \(\Arc{T_+}\) is the union of some
            connected components of the subgroup
            \begin{align}
                \label{eqn:closure_of_Tplus_r}
                \Set*{t_+\in T_+(\cO)\given r_{t_+}(\alpha)\ge r(\alpha)}.
            \end{align}
        \item \(T_+(\cO)_r\) is non-singular and \(\pi_0(T_+(\cO)_r)\) is in
            bijection with \(\pi_0^\circ(Z(L_1))\). 
        \item The codimension of \(T_+(\cO)_r\) in \(\Arc{T_+}\) is
            \begin{align}
                \label{eqn:codim_Tplus_r}
                \sum_{n=0}^\infty
                n\dim_{\bar{k}}(\La{z}_{n+1}/\La{z}_{n})=\sum_{n=1}^\infty\dim_{\bar{k}}(\La{t}_+/\La{z}_n).
            \end{align}
    \end{enumerate}
\end{proposition}
\begin{proof}
    Let \(X\) be the set in \eqref{eqn:closure_of_Tplus_r}, which is easily seen
    a subgroup hence non-singular. Clearly 
    \(X\) is closed in \(\Arc{T_+}\) and contains 
    \(T_+(\cO)_r\) as an open subset. So \(T_+(\cO)_r\) is also non-singular.

    Similar to the proof of
    \Cref{thm:characterize_rt_by_filtration}, let \(L_n\) be the
    connected reductive subgroups of \(G_+\) 
    given by \(\Roots_n(r)\), and \(\La{z}_n\subset \La{t}_+\) be the Lie algebra of
    the center of \(L_n\). Let \(S_1=Z(L_1)\) and \(S_n\subset
    Z(L_n)_{n-2}\subset \Arc{Z(L_n)_0}\) be the lift of \(\La{z}_n\) for \(n\ge
    2\).  Then \(S_n=\Set{1}\) if \(n>r(\alpha)\) for all \(\alpha\).  Then it is
    clear by the factorization in
    \Cref{thm:characterize_rt_by_filtration} that the multiplication map
    \begin{align}
        \label{eqn:proj_map_closure_of_Tplus_r}
        \phi\colon\prod_{n=1}^\infty S_n\to X
    \end{align}
    is surjective. Therefore, \(\pi_0(X)\) is bounded by
    \(\pi_0(Z(L_1))\). We also have a projection \(X\to Z(L_1)\), which implies
    the bijection \(\pi_0(X)\cong \pi_0(Z(L_1))\).
    The bijection \(\pi_0(T_+(\cO)_r)\cong\pi_0^\circ(Z(L_1))\) is then clear as the
    projection of \(\phi^{-1}(T_+(\cO)_r)\) in \(S_n\) (\(n\ge 2\)) is dense.

    To prove the codimension formula, note 
    that \(\codim_{\Arc{T_+}}(T_+(\cO)_r)=\codim_{\Arc{T_+}}(X)\). So we only
    need to show the formula for \(X\), which is easily deduced from
    \eqref{eqn:proj_map_closure_of_Tplus_r}.
\end{proof}

\begin{corollary}\label[corollary]{cor:codim_Tplus_w_r}
    Fix \(w\in W\) with \(\ord(w)=l\). 
    Suppose \(r\) takes values in \(\bbN/l\) and \(T_{+,w}(\cO)_r\neq\emptyset\).
    Let \(L_n\) be the (pseudo-)Levi subgroup
    of \(G_+\) determined by \(\Roots_n(r)\) and \(\La{z}_n=\Lie(Z(L_n))\).
    Then we have the following:
    \begin{enumerate}
        \item The closure of \(T_{+,w}(\cO)_r\) in \(\Arc{T_{+,w}}\) 
            is the union of some connected components of the subgroup
            \begin{align}
                \label{eqn:closure_of_Tplus_w_r}
                \Set*{t_+\in T_{+,w}(\cO)\given r_{t_+}(\alpha)\ge r(\alpha)}.
            \end{align}
        \item \(T_{+,w}(\cO)_r\) is non-singular and \(\pi_0(T_{+,w}(\cO)_r)\)
            is in bijection with \(\pi_0^\circ(Z(L_1)^w)\), the latter defined
            as the preimage of \(\pi_0^\circ(Z(L_1))\) in \(\pi_0(Z(L_1)^w)\). 
        \item The codimension of \(T_{+,w}(\cO)_r\) in \(\Arc{T_{+,w}}\) is
            \begin{align}
                \label{eqn:codim_Tplus_w_r}
                \sum_{n=0}^\infty\dim_{\bar{k}}(\La{t}_+/\La{z}_{n+1})(w,n).
            \end{align}
    \end{enumerate}
\end{corollary}
\begin{proof}
    The first and second claim is a result of \eqref{eqn:Tw_eigen_condition} and 
    \Cref{prop:codim_Tplus_r}. The last claim is deduced from
    \Cref{cor:characterizationn_of_Twr} and
    \Cref{prop:codim_Tplus_r}, and the fact that the action of \(w\)
    on any \(\bar{k}\)-vector space is semisimple due to our assumption on
    \(\Char(k)\).
\end{proof}

\begin{corollary}\label[corollary]{cor:codim_FRT_w_r}
    Fix \(w\in W\) with \(\ord(w)=l\). 
    Suppose \(r\) takes values in \(\bbZ/l\) 
    and \(\FRT_w(\cO)_{(\bar{\lambda}/l,r)}\neq\emptyset\).
    Let \(L_n\) be the (pseudo-)Levi subgroup
    of \(G_+\) determined by \(\Roots_n(r^\star)\) and \(\La{z}_n=\Lie(Z(L_n))\).
    Then we have the following:
    \begin{enumerate}
        \item The closure of \(\FRT_{w}(\cO)_{(\bar{\lambda}/l,r)}\) in \(\Arc{\FRT_{w}^{\lambda/l}}\) 
            is the union of some connected components of
            \begin{align}
                \label{eqn:closure_of_FRT_w_r}
                \Set*{t_+\in \Arc{\FRT_w^{\lambda/l}}(\bar{k})\given r_{t_+}(\alpha)\ge r(\alpha)},
            \end{align}
            which is itself a union of torsors under a subgroup scheme of
            \(\Arc{T_{+,w}}\).
        \item \(\FRT_{w}(\cO)_{(\bar{\lambda}/l,r)}\) is non-singular and
            \(\pi_0(\FRT_{w}(\cO)_{(\bar{\lambda}/l,r)})\)
            is in bijection with \(\pi_0^\circ(Z(L_1)^w)\), the latter defined
            as the preimage of \(\pi_0^\circ(Z(L_1))\) in \(\pi_0(Z(L_1)^w)\). 
        \item The codimension of \(\FRT_{w}(\cO)_{(\bar{\lambda}/l,r)}\) in 
            \(\Arc{\FRT_{w}^{\lambda/l}}\) is
            \begin{align}
                \label{eqn:codim_FRT_w_r}
                \sum_{n=0}^\infty\dim_{\bar{k}}(\La{t}_+/\La{z}_{n+1})(w,n),
            \end{align}
    \end{enumerate}
\end{corollary}
\begin{proof}
    This is a straightforward result of the discussion at the end of
    \Cref{sec:root_valuation_functions_and_filtrations}, \eqref{eqn:Tw_eigen_condition} and
    \Cref{cor:codim_Tplus_w_r} (the same proof applies).
\end{proof}

\subsection{}%
Recall we have the set \(S_G^1\) of pairs \((\bar{\lambda},r)\) where
\(\bar{\lambda}\in\CoCharG(G/G^\Der)\) and \(r\colon\Roots\to
\bbQ\cup\Set{\infty}\) for any connected reductive group \(G\). 
The Weyl group \(W\) acts on \(W\x S_G^1\) by conjugation on the first factor
and the action on \((\bar{\lambda},r)\) is the most obvious one:
\begin{align}
    w(\bar{\lambda},r)=(\bar{\lambda},wr\defeq \alpha\mapsto r(w^{-1}\alpha)).
\end{align}
Let \(\cS_G=W\backslash(W\x S_G^1)\), whose elements will be denoted by
\([w,\bar{\lambda}/l,r]\) (where \(l=\ord(w)\), and \(\bar{\lambda}/l\) is
integral). In case where the group is \(G_+\), we
will simply use \(\cS=\cS_{G_+}\). Let 
\(\cS_G^\rss\subset\cS_G\) be the subsets where \(r\) takes finite
values. Clearly, if \(t\in T_w(F)_{(\bar{\lambda}/l,r)}\), then \(u(t)\in
T_{uwu^{-1}}(F)_{(\bar{\lambda}/l,ur)}\) for any \(u\in W\). 
This justifies the following definition.

\begin{definition}
    For \(s\in \cS\), define \(\FRC_{\FRM}(\cO)_s\) be the image
    of \(\FRT_{w}(\cO)_{(\bar{\lambda}/l,r)}\) in \(\FRC_{\FRM}(\cO)\) 
    for any (and every)
    representative \((w,\bar{\lambda}/l,r)\) of \(s\).
\end{definition}

\begin{lemma}
    The set of \(\bar{k}\)-points of \(\Loop^\flat{\FRC}_\FRM\defeq
    \Arc{\FRC_{\FRM}}-\Arc{\FRE_{\FRM}}-\Arc{\FRD_{\FRM}}\) is the \emph{disjoint}
    union of strata \(\FRC_\FRM(\cO)_s\) for all \(s\in\cS^\rss\), where
    \(\FRE_{\FRM}\) (resp.~\(\FRD_{\FRM}\)) is the numerical boundary divisor
    (resp.~discriminant divisor) as in
    \Cref{sec:review_of_very_flat_reductive_monoids,sec:invariant_theory_of_reductive_monoids}.
\end{lemma}
\begin{proof}
    Straightforward; cf.~\cite{GKM09}*{7.3}.
\end{proof}

\subsection{}%
Consider \(\cO_l\)-morphism
\begin{align}
    \chi_{\lambda/l}\colon T_{+,\cO_l}&\longto \FRC_{\FRM,\cO_l}\\
    x &\longmapsto \chi_{\FRM}(\pi^{\lambda/l}x),
\end{align}
which induces \(\cO\)-morphism
\begin{align}
    \chi_{(w,\lambda/l)}
    \colon \left(\pi^{\lambda/l}\Res_{\cO_l/\cO}T_{+,\cO_l}\right)^{\tau_\infty w}&\longto
    \FRC_{\FRM,\cO}.
\end{align}
This further
induces map \(\Arc\FRT_w^{\lambda/l}\to \Arc{\FRC_{\FRM}}\), which is precisely the
restriction of \(\Arc{\chi_{\FRM}}\) to \(\Arc\FRT_w^{\lambda/l}\).
Fix any \(x\in\Arc\FRT_w^{\lambda/l}(\bar{k})\), let
\(a=\chi(x)\in\Arc{\FRC_{\FRM}}(\bar{k})\), 
then \(\chi_{\FRM}\) induces \(\tau_\infty w\)-equivariant 
isomorphism of \(F_l\)-tangent spaces
\begin{align}
    (\dd \chi_{\FRM,F_l})_x\colon \TanB_x(\pi^{\lambda/l}T_{+,F_l})\longto
    \TanB_{a}\FRC_{\FRM,F_l}\simeq \FRC_{\FRM,F_l},
\end{align}
which restricts to an injection of \(\cO_l\)-modules
\begin{align}
    (\dd \chi_{\FRM,F_l})_x\colon \La{t}_{\cO_l}\longto \FRC_{\FRM,\cO_l},
\end{align}
Taking \(\tau_\infty w\)-invariants, we have \(\cO\)-linear injection
\begin{align}
    (\dd \chi_w)_x\colon Q_x\defeq \Res_{\cO_l/\cO}\La{t}_{\cO_l}^{\tau_\infty
    w}\longto \FRC_{\FRM,\cO}.
\end{align}
Choosing an \(\cO\)-basis on both sides, then \((\dd \chi_w)_x\) is 
represented by a matrix in \(\GL_{2r}(F)\). Changing uniformizer \(\pi^{1/l}\)
or the \(\cO\)-bases does not change the
valuation of the determinant of the matrix, so \(\val_F\det(\dd \chi_w)_x\) is
well-defined for \(x\). 

\begin{proposition}\label[proposition]{prop:det_computation}
    Suppose \(x\in \FRT_w(\cO)_{(\bar{\lambda}/l,r)}\) is generically
    regular semisimple, then
    \begin{align}
        \label{eqn:val_det_d_chi}
        \val_F\det(\dd
        \chi_w)_x&=\sum_{i=1}^r\Pair{\alpha_i}{\bar{\lambda}/l}+\frac{d_+(a)+c(a)}{2}\\
            &=\sum_{i=1}^r\Pair{\alpha_i}{\bar{\lambda}/l}+\delta(a)+c(a).
    \end{align}
\end{proposition}
\begin{proof}
    The proof is completely parallel to the argument in
    \cite{GKM09}, but the
    technical counterparts will occupy several pages.
    We extend \((\dd\chi_w)_x\) to \(\cO_l\):
    \begin{align}
        \Id_{F_l}\otimes (\dd\chi_w)_x\colon \cO_l\otimes_\cO Q_x\longto \FRC_{\cO_l},
    \end{align}
    which is the restriction of map
    \begin{align}
        (\dd\chi_{\FRM,F_l})_x\colon \La{t}_{\cO_l}\longto \FRC_{\cO_l}
    \end{align}
    to submodule \(\cO_l\otimes_\cO Q_x\).  Therefore,
    \begin{align}
        \val_F\det(\dd\chi_w)_x =\val_F\det(\dd\chi_{\FRM,F_l})_x
        +\frac{1}{l}\dim_{\bar{k}}\frac{\La{t}_{\cO_l}}{\cO_l\otimes_\cO Q_x}.
    \end{align}
    So we reduce to the \(w=1\) case (replacing \(\cO\) with \(\cO_l\) and
    \(\val_F\) with \(\val_{F_l}\)) and the claim that
    \begin{align}\label{eqn:c_as_dim_of_quotient}
        \dim_{\bar{k}}\frac{\La{t}_{\cO_l}}{\cO_l\otimes_\cO Q_x}=\frac{lc(a)}{2}.
    \end{align}

    The equation \eqref{eqn:c_as_dim_of_quotient} is proved using the same
    argument in \cite{Be96}. By definition, \(c(a)\) is the
    dimension of the largest subspace of \(\La{t}\) (over \(\bar{k}\)) that does not
    contain a trivial representation of \(w\).
    On the other hand, we have
    \begin{align}
        Q_x=\left[\bigoplus_{i=0}^{l-1}\La{t}(i)\pi^{i/l}\right]\otimes_{\bar{k}}\cO,
    \end{align}
    where \(\La{t}(i)\) is the eigenspace of \(w\) 
    in \(\La{t}\) with eigenvalue
    \(\zeta_l^{-i}\). Therefore, we can describe the quotient space as follows:
    \begin{align}
        \frac{\La{t}_{\cO_l}}{\cO_l\otimes_\cO Q_x} 
        &\simeq \bigoplus_{i=0}^{l-1}\left[\,\bigoplus_{j=i+1}^{l-1}\La{t}(j)\,\right],
    \end{align}
    which implies \eqref{eqn:c_as_dim_of_quotient}.

    The argument for \(w=1\) (and \(l=1\)) case is a slight 
    generalization of the argument in
    \cite{St74}*{pp.~125--127}. For completeness, we include a detailed
    argument here. Without loss of generality, we may assume \(\lambda_\AD\) is
    \emph{anti-dominant}. 
    Let \(\SimRts_\lambda\) be the simple roots vanishing on \(\lambda\),
    \(\Roots_\lambda\) the induced root subsystem, and
    \(W_\lambda\subset W\) the subgroup 
    generated by the reflections corresponding to roots in
    \(\SimRts_\lambda\). It is harmless to assume that \(\SimRts_\lambda\) contains
    first \(s\) simple roots.

    Choose generator basis (in the given
    order)
    \begin{align}
        \Set*{e^{(\alpha_1,0)},\ldots,e^{(\alpha_n,0)}, e^{(\Wt_1,\Wt_1)},\ldots,
        \ldots,e^{(\Wt_n,\Wt_n)}}
    \end{align}
    of \(\bar{k}[T_+]\), and 
    \begin{align}
        \Set*{e^{(\alpha_1,0)},\ldots,e^{(\alpha_n,0)}, \chi_{1,+},
        \ldots,\chi_{r,+}}
    \end{align}
    of \(\bar{k}[\FRC_{\FRM}]\), where \(\chi_{i,+}(z,t)=\Wt_i(z)\chi_i(t)\). Let \(e_i=\dd
    e^{(\alpha_i,0)}\), \(f_i=e^{-(\Wt_i,\Wt_i)}\dd e^{(\Wt_i,\Wt_i)}\), and
    \(g_i=\dd \chi_{i,+}\). We need to compute
    the valuation of linear map
    \begin{align}
        \dd (\chi_\lambda)_x\colon 
        \cO e_1\wedge\cdots\wedge e_r\wedge g_1\wedge \cdots\wedge g_n
        \longto 
        \cO e_1\wedge\cdots\wedge e_r\wedge f_1\wedge \cdots\wedge f_r,
    \end{align}
    identified with \(A\in \Mat_1(\cO)\simeq\cO\).  We claim
    that \(A\) is \(W_\lambda\)-skew symmetric, in other words, \(w(A)=\det(w)A\)
    for all \(w\in W_\lambda\). Indeed, fix any \(w\in W_\lambda\), and suppose
    that
    \begin{align}
        w(\Wt_i)=\sum_{i=1}^r n_{ij}\Wt_j,
    \end{align}
    then
    \begin{align}
        w(f_i) &= e^{-(\Wt_i,w(\Wt_i))}\dd e^{(\Wt_i,w(\Wt_i))}\\
               &= e^{(w(\Wt_i)-\Wt_i,0)}e^{-(w(\Wt_i),w(\Wt_i))}
               \dd \left[e^{(w(\Wt_i)-\Wt_i,0)}e^{(w(\Wt_i),w(\Wt_i))}\right]\\
               &= \left(\prod_{j=1}^r e^{-(n_{ij}\Wt_j,n_{ij}\Wt_j)}\right)
               \dd\left(\prod_{j=1}^r e^{(n_{ij}\Wt_j,n_{ij}\Wt_j)}\right)\\
               &\qquad +\text{terms involving at least one }e_j\\
               &= \sum_{j=1}^r n_{ij}f_j+\text{terms involving at least one }e_j.
    \end{align}
    Therefore, since \(e_i\) and \(g_i\) are fixed by \(w\),
    \begin{align}
        w(A e_1\wedge\cdots\wedge e_r\wedge f_1\wedge \cdots\wedge
        f_r)&=w(A)e_1\wedge\cdots\wedge e_r\wedge w(f_1)\wedge \cdots\wedge
        w(f_r)\\
            &= w(A)\det(n_{ij})e_1\wedge\cdots\wedge e_r\wedge
            f_1\wedge\cdots\wedge f_r\\
            &= w(A)\det(w)e_1\wedge\cdots\wedge e_r\wedge
            f_1\wedge\cdots\wedge f_r.
    \end{align}
    This means that \(w(A)=\det(w)^{-1}A=\det(w)A\) as claimed.

    Let \(\WP_i\) (resp.~\(\WP_\nu\)) be the set of weights in the Weyl module of
    \(G\) of highest weight \(\Wt_i\) (resp.~any dominant \(\nu\)), 
    and \(m_{i,\mu}\) the multiplicity of \(\mu\in\WP_i\).
    We can expand \(g_i\) into linear combinations of \(e_j\) and \(f_j\) as
    follows:
    \begin{align}
        g_i&=\pi^{\Pair{(\Wt_i,\Wt_i)}{\lambda}}\dd
        e^{(\Wt_i,\Wt_i)}+\sum_{\Wt_i\neq \mu\in\WP_i}
        m_{i,\mu}\pi^{\Pair{(\Wt_i,\mu)}{\lambda}}\dd e^{(\Wt_i,\mu)}\\
           &=\pi^{\Pair{(\Wt_i,\Wt_i)}{\lambda}}e^{(\Wt_i,\Wt_i)}f_i
               +\sum_{\Wt_i\neq \mu\in\WP_i}
               m_{i,\mu}\pi^{\Pair{(\Wt_i,\mu)}{\lambda}}e^{(\Wt_i,\mu)}e^{-(\mu,\mu)}\dd
                   e^{(\mu,\mu)}\\
           &\qquad +\text{terms involving at least one }e_j\\
           &=\pi^{\Pair{(\Wt_i,\Wt_i)}{\lambda}}e^{(\Wt_i,\Wt_i)}f_i
               +\sum_{\Wt_i\neq \mu\in\WP_i}
               \left[m_{i,\mu}\pi^{\Pair{(\Wt_i,\mu)}{\lambda}}e^{(\Wt_i,\mu)}\sum_{j=1}^r
               n_{\mu,j}f_j\right]\\
           &\qquad +\text{terms involving at least one }e_j.
    \end{align}
    This implies that
    \begin{align}
        B\defeq\prod_{i=1}^r\pi^{-\Pair{\alpha_i}{\bar{\lambda}}}A&=
        \pi^{\Pair{(\rho,\rho)}{\lambda}}e^{(\rho,\rho)}+
        \sum_{\mu\in\WP_\rho}C_\mu
        \pi^{\Pair{(\rho,\mu)}{\lambda}}e^{(\rho,\mu)}
    \end{align}
    for some integers \(C_\mu\). As \(B\) is \(W_\lambda\)-skew symmetric
    (because \(A\) is), we
    have that
    \begin{align}
        B&=\sum_{w\in W_\lambda}\det(w)\pi^{\Pair{(\rho,w(\rho))}{\lambda}}e^{(\rho,w(\rho))}+
        \sum_{\substack{\mu\in\WP_\rho\\\SimRts_\lambda\text{-dominant}}}C_\mu
        \left(\sum_{w\in
        W_\lambda}\det(w)\pi^{\Pair{(\rho,w(\mu))}{\lambda}}e^{(\rho,w(\mu))}\right)\\
        &=\pi^{\Pair{(\rho,\rho)}{\lambda}}\sum_{w\in W_\lambda}\det(w)e^{(\rho,w(\rho))}+
        \sum_{\substack{\mu\in\WP_\rho\\\SimRts_\lambda\text{-dominant}}}C_\mu
        \pi^{\Pair{(\rho,\mu)}{\lambda}}\left(\sum_{w\in
        W_\lambda}\det(w)e^{(\rho,w(\mu))}\right).
    \end{align}
    Observe that if \(\mu\) is not \emph{strictly} \(\SimRts_\lambda\)-dominant, then
    the summation
    \begin{align}
        \sum_{w\in
        W_\lambda}\det(w)\pi^{\Pair{(\rho,w(\mu))}{\lambda}}e^{(\rho,w(\mu))}
    \end{align}
    equals \(0\) because its summands cancel pairwise.

    We then claim that for strictly \(\SimRts_\lambda\)-dominant 
    \(\mu\), \(\Pair{(\rho,\mu)}{\lambda}\) reaches minimum if and only if
    \(\mu=\rho\). Indeed, suppose 
    \begin{align}
        \mu&=\rho - \sum_{i=1}^s p_i\alpha_i -
        \sum_{i=s+1}^r p_i\alpha_i =\rho +\sum_{i=1}^r q_i\Wt_i
    \end{align}
    for some \(p_i\in\bbN\) and \(q_i\in\bbZ\). If \(p_i>0\) for any \(i>s\),
    we are done since \(\lambda_\AD\) is anti-dominant and
    \(\Pair{\alpha_i}{\lambda_\AD}<0\) if \(i>s\). Now suppose
    \(p_i=0\) for all \(i>s\). Since \(\mu\) is strictly
    \(\SimRts_\lambda\)-dominant and \(\Pair{\rho}{\dual{\alpha}_i}=1\),
    we must have \(q_i\ge 0\) for \(1\le i\le s\). 
    Let \(\dual{\rho}_\lambda\) be the half-sum of 
    the positive coroots in \(\Roots_\lambda\).
    Then we have
    \begin{align}
        \Pair{\rho-\mu}{\dual{\rho}_\lambda}=\sum_{i=1}^s p_i\ge 0,
    \end{align}
    but on the other hand
    \begin{align}
        \Pair{\rho-\mu}{\dual{\rho}_\lambda}=-\sum_{i=1}^s\Pair{\Wt_i}{\dual{\rho}_\lambda}q_i\le 0,
    \end{align}
    which means \(p_i=q_i=0\) for all \(1\le i\le s\) and \(\mu=\rho\). Hence,
    the claim is proved. 

    For strictly \(\SimRts_\lambda\)-dominant \(\mu\),
    a well-known fact due to Weyl gives that
    \begin{align}
        \sum_{w\in
        W_\lambda}\det(w)e^{(\rho,w(\mu))}=e^{(\rho,\mu)}\prod_{\alpha\in\Roots_\lambda^+}\left(1-e^{-(0,\alpha)}\right).
    \end{align}
    Thus, we have that
    \begin{align}
        \val_F(B)&=\val_F\left[\pi^{\Pair{(\rho,\rho)}{\lambda}}\sum_{w\in
        W_\lambda}\det(w)e^{(\rho,w(\rho))}\right.\\
            &\qquad
            \left.+\sum_{\substack{\mu\in\WP_\rho\\\text{strictly }\SimRts_\lambda\text{-dominant}}}C_\mu
        \pi^{\Pair{(\rho,\mu)}{\lambda}}\left(\sum_{w\in
    W_\lambda}\det(w)e^{(\rho,w(\mu))}\right)\right]\\
            &=\val_F\left(\prod_{\alpha\in\Roots_\lambda^+}\left(1-e^{-(0,\alpha)}\right)\right)
            +\val_F\left[\pi^{\Pair{(\rho,\rho)}{\lambda}}
            +\sum_{\substack{\mu\in\WP_\rho\\\text{strictly
        }\SimRts_\lambda\text{-dominant}}}
    \pi^{\Pair{(\rho,\mu)}{\lambda}}\right]\\
    &=\val_F\left(\prod_{\alpha\in\Roots_\lambda^+}\left(1-e^{-(0,\alpha)}\right)\right)+\Pair{(\rho,\rho)}{\lambda}\\
    &=\frac{d_+(a)}{2}.
    \end{align}
    So 
    \begin{align}
        \val_F(A)&=\sum_{i=1}^r\Pair{\alpha_i}{\bar{\lambda}}+\frac{d_+(a)}{2},
    \end{align}
    as desired. This concludes the proof.
\end{proof}

\begin{corollary}\label[corollary]{cor:det_computation_general}
    Let \(\FRM\) be an arbitrary very flat reductive monoid associated with \(G\)
    such that its abelianization is an affine space \(\FRA_\FRM=\bbA^m\) with
    coordinates \(e^{\alpha_1'},\ldots,e^{\alpha_m'}\).
    Let \(\FRT_\FRM\) be the maximal toric variety in \(\FRM\). Let \(x\in
    \FRT_{\FRM,w}(\cO)_{(\lambda/l,r)}\) be generically regular semisimple whose
    image is \(a\in\FRC_\FRM(\cO)\). Then
    \begin{align}
        \label{eqn:val_det_d_chi_general_monoid}
        \val_F\det(\dd
        \chi_{\FRM,w})_x&=b(a)+\frac{d_+(a)+c(a)}{2}\\
                     &=b(a)+\delta(a)+c(a),
    \end{align}
    where \(b(a)\) is the boundary valuation
    \(\sum_{i=1}^m\Pair{\alpha_i'}{\bar{\lambda}/l}\).
\end{corollary}
\begin{proof}
    Note that the notion of 
    \(d_+(a)\) makes sense as it is the valuation of the pullback of
    the discriminant function (from \(\Env(G^\SC)\) to \(\FRM\)), and that
    the notion of \(c(a)\) also makes sense because it only depends on \(w\).

    The character lattice of the maximal torus in \(\FRM\) is freely spanned by
    the pullbacks of functions \(e^{\alpha_i'}\) from \(\FRA_\FRM\), and
    \(e^{(\Wt_i,\Wt_i)}\) from \(\FRT_{\Env(G^\SC)}\). The
    proof then proceeds essentially  the same in the proposition above; notably
    the computation for \(B\) is repeated word-by-word.
\end{proof}

\section{Codimension Formula of Valuation Strata}%
\label{sec:codimension_formula_val_strata}

 We can now use \eqref{eqn:val_det_d_chi_general_monoid} to derive a codimension
formula for the valuation strata for arbitrary very flat reductive
monoid \(\FRM\) of \(G\) whose abelianization is an affine space.

\begin{lemma}\label[lemma]{lem:cylinder_to_cylinder_one_point}
    Let \(f\colon \Gm[\cO]^r\to\bbA_\cO^r\) be an \(\cO\)-morphism, viewed as an
    \(r\)-vector of rational functions in \(r\) variables and
    \(\cO\)-coefficients. Let \(x\in \Gm^r(\cO)\) such that \(d=\val_F(\det \DD_x
    f)\)
    is finite. Then for any \(m>d\), we have that
    \begin{align}
        f(x+\pi^m)=f(x)+\DD_x f(\pi^m),
    \end{align}
    where \(\pi^m\) means \(\pi^m\) times the \(\cO\)-tangent space at \(x\).
\end{lemma}
\begin{proof}
    This is just a tiny modification of the claim in \cite{GKM09}*{Lemma~10.3.1},
    but since they did not provide a general statement, we include the
    adaptation here.

    We may regard
    \(\Gm^r(\cO)\) as a subset of \(\bbA^r(\cO)\) and 
    identify the \(\cO\)-tangent space at \(x\) 
    with \(r\)-copies of \(\cO\), denoted by \(\Lambda\). 
    Take the Taylor expansion of \(f\) at \(x\) over \(\cO\),
    whose linear terms gives the matrix \(A=\DD_x f\). Since \(m>d\), we have
    \(m\ge 1\), so in turn we have for any \(h\in \pi^m \Lambda\),
    \begin{align}
        f(x+h)\equiv f(x)+(\DD_x f)(h)\bmod \pi^{2m}\Lambda.
    \end{align}
    The proof then proceeds as in \cite{GKM09}*{Lemma~10.3.1}.
\end{proof}

\begin{lemma}
    The valuation strata \(\FRC_{\FRM}(\cO)_{[w,\bar{\lambda}/l,r]}\) are
    cylinders.
\end{lemma}
\begin{proof}
    It is clear from the fact that \(\FRT_w(\cO)_{(\bar{\lambda}/l,r)}\) are
    cylinders, \Cref{cor:det_computation_general}, and
    \Cref{lem:cylinder_to_cylinder_one_point} (applied to \(\pi^\lambda
    T_\FRM(\cO)\cong T_\FRM(\cO)\to \FRC_\FRM(\cO)\)).
\end{proof}

\begin{theorem}
    \label[theorem]{thm:codim_valuation_strata}
    The codimension of \(\FRC_\FRM(\cO)_{[w,\bar{\lambda}/l,r]}\) in \(\FRC_\FRM(\cO)\)
    is given by
    \begin{align}
        b+\frac{d_+ + c}{2}+e=b+\delta+c+e,
    \end{align}
    where \(b\), \(d_+\), \(c\), \(\delta\) are as in
    \Cref{cor:det_computation_general} (which are necessarily
    constant over \(\FRC_\FRM(\cO)_{[w,\bar{\lambda}/l,r]}\)), and \(e\) is the
    codimension of \(\FRT_w(\cO)_{(\bar{\lambda}/l,r)}\) in
    \(\Arc{\FRT^{\lambda/l}_w}\).
\end{theorem}
\begin{proof}
    Let \(U=\FRT_w(\cO)_{(\bar{\lambda}/l,r)}\) and
    \(V=\FRC_\FRM(\cO)_{[w,\bar{\lambda}/l,r]}\). Since they are both cylinders, we
    may pass to some jet space modulo \(\pi^N\) for some big \(N\), after which
    \(U\) becomes smooth and \(V\) is constructible, hence the notion of 
    codimension makes sense for \(V\). Since \(U\)
    is smooth (\Cref{cor:codim_FRT_w_r})
    and the original \(\cO\)-map from \(U\) to \(V\) have constant valuation
    (\Cref{cor:det_computation_general}), we may find an open dense
    subset of \(V\) modulo \(\pi^N\) over which \(U\to V\) is smooth of
    relative dimension \(b+\delta+c\). Then the claim follows from standard
    dimension counting (cf.~\cite{GKM09}*{\S~5}).
\end{proof}

\begin{remark}
    In \cite{GKM09}, the authors proved more properties about valuation strata, such
    as their non-singularity and the smoothness of the respective maps between the strata in
    \(\La{t}\) and those in \(\La{t}\git W\). Their argument is by induction
    on Levi subgroups, and it can be generalized to the group case using Levi
    monoids defined in \Cref{sec:endoscopic_groups_inv_theory}. We leave out the
    details since we will not use this result.
\end{remark}

\subsection{}
We make the following definition for future convenience:
\begin{definition}
    \label[definition]{def:newton_rss}
    We say \(a\in\FRC_\FRM(\cO)\) is \inotion{\(\nu\)-regular semisimple}
    if \(c=e=0\), or equivalently, \(a\) is unramified and 
    any \(t\in \pi^\lambda T_\FRM(\cO)\subset\FRT_\FRM(\cO)\)
    lying over \(a\) satisfies
    \begin{align}
        \val_F(1-\Rt(t))=0
    \end{align}
    for all root \(\Rt\) such that \(\Pair{\Rt}{\lambda}=0\).
\end{definition}

\chapter{Multiplicative Affine Springer Fibers}%
\label{chap:generalized_affine_springer_fibers}

In this chapter, we review the multiplicative affine Springer fibers, also
known as Kottwitz--Viehmann (KV) varieties. Analogous to the Lie
algebra case, the point-count of these
schemes over \(k\) encode information of orbital integrals of \(G\).

\section{Definition and Generalities}%
\label{sec:definition_and_generalities}

\subsection{}
Consider \(\cO_v=k_v\powser{\pi}\) and \(F_v=k_v\lauser{\pi}\) for
a finite extension \(k_v/k\).
    \nomenclature[\(k_v \)]{\(k_v\)}{the residue field of the local field \(F_v\) at a place \(v\)}
Let \(X_v=\Spec{\cO_v}\)
    \nomenclature[\(X_v \)]{\(X_v\)}{the formal disc \(\Spec{\cO_v}\)}
be the corresponding
formal disc, and \(X_v^\bullet=\Spec{F_v}\)
    \nomenclature[\(X_v_bullet \)]{\(X_v^\bullet\)}{the punctured formal disc \(\Spec{F_v}\)}
the punctured disc. Let \(v\) be the
closed point of \(X_v\), and \(\eta_v\)
    \nomenclature[\(eta_v \)]{\(\eta_v\)}{the generic point of \(X_v\)}
the generic point.
Let \(\breve{X}_v=\Spec{\cO_v\hat{\otimes}_k\bar{k}}\)
    \nomenclature[\(X'breve_v \)]{\(\breve{X}_v\)}{the base change \(\Spec{\cO_v\hat{\otimes}_k\bar{k}}\) to \(\bar{k}\) over \(k\) of \(X_v\)}
then we have isomorphism
\begin{align}
    \breve{X}_v=\coprod_{\bar{v}\colon k_v\to\bar{k}} \breve{X}_{\bar{v}},
\end{align}
where \(\bar{v}\) ranges over all \(k\)-field embeddings of \(k_v\) into
\(\bar{k}\), and \(\breve{X}_{\bar{v}}\cong\Spec{\bar{k}\powser{\pi}}\).
    \nomenclature[\(X'breve_v_bar \)]{\(\breve{X}_{\bar{v}}\)}{the geometric
    formal disc over \(X_v\) corresponding to a \(k\)-embedding \(\bar{v}\colon k_v\to \bar{k}\)}
If we choose a geometric point \(\bar{\eta}_{\bar{v}}\)
    \nomenclature[\(eta_v_bar_bar \)]{\(\bar{\eta}_{\bar{v}}\)}{a geometric generic point of \(\breve{X}_{\bar{v}}\)}
on \(\breve{X}_{\bar{v}}\) over
its generic point \(\eta_{\bar{v}}\),
    \nomenclature[\(eta_v_bar \)]{\(\eta_{\bar{v}}\)}{the generic point of \(\breve{X}_{\bar{v}}\)}
then we have a short exact sequence
\begin{align}
    1\longto I_v\longto \Gamma_v\longto \Gal(\bar{k}/\bar{v}(k_v))\longto 1,
\end{align}
where \(\Gamma_v=\pi_1(\eta_{v},\bar{\eta}_{\bar{v}})\) is the Galois group of
\(F_v\), and \(I_v=\pi_1(\breve{X}_{\bar{v}}^\bullet,\bar{\eta}_{\bar{v}})\)
    \nomenclature[\(I_v \)]{\(I_v\)}{the inertia group in \(\Gamma_v\)}
is the inertia group. We shall use \(\breve{\cO}_{\bar{v}}\)
    \nomenclature[\(O"cal'breve_v_bar \)]{\(\breve{\cO}_{\bar{v}}\)}{the ring of regular functions on \(\breve{X}_{\bar{v}}\)}
and \(\breve{F}_{\bar{v}}\)
    \nomenclature[\(F'breve_v_bar \)]{\(\breve{F}_{\bar{v}}\)}{the ring of regular functions on \(\breve{X}_{\bar{v}}^\bullet\)}
to denote the
ring of functions on \(\breve{X}_{\bar{v}}\) and \(\breve{X}_{\bar{v}}^\bullet\),
    \nomenclature[\(X'breve_v_bar^bullet \)]{\(\breve{X}_{\bar{v}}^\bullet\)}{the geometric
    punctured disc in \(\breve{X}_{\bar{v}}\)}
respectively.
We let \(\Frob_v\)
    \nomenclature[\(sigma_v \)]{\(\Frob_v\)}{the geometric Frobenius element of
    \(\bar{k}/k_v\), or a fixed lifting of such to \(\Gamma_v\)}
to be the \emph{geometric} Frobenius element of \(k_v\)
inside \(\bar{k}\), which is independent of the embedding \(\bar{v}\), and we
also view it as a Frobenius element of \(\Gamma_v\) by fixing a lifting in the
latter.

Let \(G\) be a quasi-split reductive group over \(\cO_v\) associated with a
\(\Out(\bG)\)-torsor \(\OGT_G\) on \(X_v\) and similarly a pointed
version \(\OGT_G^\bullet\) if we fix a geometric point \(\bar{v}\) of \(X_v\).
Equivalently, if we fix a geometric generic point \(\bar{\eta}_{\bar{v}}\), we
have a homomorphism \(\OGT_G^\bullet\colon\Gamma_v\to\Out(\bG)\) which factors
through \(\Gal(\bar{k}/\bar{v}(k_v))\), because \(\OGT_G\) is always trivial over
\(\breve{X}_v\).

\subsection{}
Let \(\lambda\in\CoCharG(T)_+\)
be a dominant \(F_v\)-rational cocharacter.
We denote by \(\Cartan_{G,v}^\lambda\)
    \nomenclature[\(C"sf_G_v_lambda \)]{\(\Cartan_G^\lambda,\Cartan_{G,v}^\lambda\)}{the Cartan double coset in \(G(F_v)\) with type \(\lambda\)}
the double
coset \(G(\cO_v)\pi^{\lambda} G(\cO_v)\) in Cartan decomposition of
\(G(F_v)\). We also denote by \(\Cartan_{G,v}^{\le\lambda}\)
    \nomenclature[\(C"sf_G_v_lambda_le \)]{\(\Cartan_G^{\le\lambda},\Cartan_{G,v}^{\le\lambda}\)}{the closure of \(\Cartan_{G,v}^\lambda\)}
the union of all
\(\Cartan_{G,v}^\mu\) such that \(\mu\le \lambda\) (i.e.,
\(\lambda-\mu\in\bbN\SimCoRts\)).  Let \(\gamma\in G(F_v)\) be a regular
semisimple element. We are interested in the sets
\begin{align}\label{eqn:KV_def_as_sets}
    \begin{split}
        \KV_{G,v}^\lambda(\gamma)&\defeq\Set*{g\in G(F_v)/G(\cO_v)\given
        \Ad_g^{-1}(\gamma)\in \Cartan_{G,v}^{\lambda}},\\
        \KV_{G,v}^{\le\lambda}(\gamma)&\defeq\Set*{g\in G(F_v)/G(\cO_v)\given
            \Ad_g^{-1}(\gamma)\in \Cartan_{G,v}^{\le \lambda}},
        \nomenclature[\(M{}gamma"sf_G_lambda \)]{\(\KV_{G,v}^\lambda(\gamma)\)}{the set-theoretic version of \(\cM_{G,v}^{\lambda}(\gamma)\)}
        \nomenclature[\(M{}gamma"sf_G_lambda_le \)]{\(\KV_{G,v}^{\le\lambda}(\gamma)\)}{the set-theoretic version of \(\cM_{G,v}^{\le\lambda}(\gamma)\)}
    \end{split}
\end{align}
on which we will later impose structures of \(k\)-varieties.

In order to best establish the connection to reductive monoids, 
it is best to generalize the definition
\eqref{eqn:KV_def_as_sets} as follows. Let \(G^\AB=G/G^\Der\), then the
homomorphism \(T\to G^\AB\x T^\AD\) is \'etale. We have the induced monomorphism
of lattices \(\CoCharG(T)\to\CoCharG(G^\AB)\x \CoCharG(T^\AD)\) with finite
index, and the same holds after taking \(\Gamma_v\)-invariants. 
Denote the image of \(\lambda\) under this map by
\((\lambda_\AB, \lambda_\AD)\). It is clear that for any \(\gamma\in G(F_v)\), the
condition \(\gamma\in \Cartan_G^\lambda\) is the same as
\(\gamma_\AD\in\Cartan_{G^\AD}^{\lambda_\AD}\) and \(\gamma_\AB\in
\pi^{\lambda_\AB}G^\AB(\cO_v)\). Therefore, if we define, for 
\(\lambda\in\CoCharG(T^\AD)\) and \(\gamma\in G^\AD(F_v)\),
\begin{align}\label{eqn:general_KV_def_as_sets}
    \begin{split}
        \KV_{G,v}^{\lambda}(\gamma)&\defeq\Set*{g\in G(F_v)/G(\cO_v)\given
        \Ad_g^{-1}(\gamma)\in \Cartan_{G^\AD,v}^{\lambda}},\\
        \KV_{G,v}^{\le\lambda}(\gamma)&\defeq\Set*{g\in G(F_v)/G(\cO_v)\given
            \Ad_g^{-1}(\gamma)\in \Cartan_{G^\AD,v}^{\le \lambda}},
    \end{split}
\end{align}
then the sets in \eqref{eqn:KV_def_as_sets}, if non-empty, are respectively
isomorphic to the sets in \eqref{eqn:general_KV_def_as_sets}, provided that
\(\gamma\) lifts to a point in \(G(F_v)\). By replacing \(X_v\)
with \(\breve{X}_v\), we have analogously defined sets, which we will call the
\(\bar{k}\)-points of the corresponding sets, and they are exactly the set of
\(\bar{k}\)-points of the corresponding schemes once we define them.

\subsection{}%
\label{sub:cartan_decomposition}

Let \(\FRM\in\FM(G^\SC)\). Recall we have the abelianization \(\FRA_{\FRM}\) and its
subtorus \(\FRA_{\FRM}^\x\). We call an \(F_v\)-rational cocharacter
\(\lambda\in\CoCharG(\FRA_{\FRM}^\x)\)
\notion{dominant}\index{dominant!cocharacter}\index{cocharacter!dominant} if it
is contained in the cocharacter cone
\(\sC(\FRA_{\FRM})\) or equivalently, \(\pi^{\lambda}\in\FRA_{\FRM}(\cO_v)\). For
each dominant \(\lambda\), we define a reductive monoid \(\cO_v\)-scheme
\(\FRM^{\lambda}\) by the pullback diagram
\begin{equation}
    \begin{tikzcd}
        \FRM^{\lambda} \ar[r] \ar[d] & \FRM\x \FRA_{\FRM}^\x \ar[d, "{(x,z)\mapsto
        \alpha_{\FRM}(x)z}"]\\
        \Spec{\cO_v}\ar[r, "\pi^{\lambda}"] & \FRA_{\FRM}
    \end{tikzcd},
    \nomenclature[\(M"frak_lambda \)]{\(\FRM^\lambda\)}{the \(\cO_v\)-subscheme
    of \(\FRM_{\cO_v}\) of points with boundary cocharacter \(\lambda\)}
\end{equation}
and by replacing \(\FRM\) with the big-cell locus \(\FRM^\circ\), we obtain an open
subscheme \(\FRM^{\circ\lambda}\) of \(\FRM^{\lambda}\).
\begin{remark}
    Note the difference in notation compared to \cite{Ch22}: if \(G=G^\SC\) and
    \(\FRM=\Env(G)\), then our \(\FRM^\lambda\) is the same as
    \(\mathup{Vin}_G^{-w_0(\lambda)}\) (not \(\mathup{Vin}_G^{\lambda}\)) in
    \cite{Ch22}.  The reason is that in \cite{Ch22}, they only
    consider the universal monoid, and taking the involution by
    \(-w_0\) saves a lot of notations later on.  Here more general framework is
    considered where there may not be a natural involution by \(-w_0\) on
    \(\sC(\FRA_{\FRM})\). 
\end{remark}

Recall for \(\FRM\in\FM(G^\SC)\), we can choose an excellent morphism \(\FRM\to
\Env(G^\SC)\) so that we have an induced map
\(\CoCharG(\FRA_{\FRM}^\x)\to\CoCharG(T^\AD)\). Let \(\lambda_\AD\) be the image
of \(\lambda\in \CoCharG(\FRA_{\FRM}^\x)\) under this map. Using the same argument 
in \cite{Ch22}*{Lemma~2.5.1}, one can easily show the following
results.
\begin{lemma}\label[lemma]{lem:arc_space_of_monoid_and_Cartan_cells}
    We have a disjoint union of \(\FRM^\x(\cO_v)\)-stable subsets
    \begin{align}
        \FRM(\cO_v)\cap \FRM^\x(F_v)=\bigcup_{\lambda\in
        \sC(\FRA_{\FRM})}\FRM^{\lambda}(\cO_v).
    \end{align}
    Moreover, let \(g\in \FRM^\x(F_v)\), then \(g\in \FRM^{\lambda}(\cO_v)\)
    (resp.~\(\FRM^{\circ\lambda}(\cO_v)\)) if and only if
    \(\alpha_{\FRM}(g)\in \pi^{\lambda} \FRA_{\FRM}^\x(\cO_v)\), and the image of \(g\) in
    \(G^{\AD}(F_v)\) belongs to \(\Cartan_{G^{\AD}}^{\le -w_0(\lambda_\AD)}\)
    (resp.~\(\Cartan_{G^{\AD}}^{-w_0(\lambda_\AD)}\)).
\end{lemma}

\subsection{}%
\label{sub:hodge_newton_decomposition_and_kottwitz_map}

For \(\gamma\in G(F_v)\), one can define a dominant element
\(\nu_\gamma\in\CoCharG(T)_\bbQ\),
    \nomenclature[\(nu_gamma \)]{\(\nu_\gamma\)}{the Newton point of an \(F_v\)-point \(\gamma\) in a reductive group}
called the \notion{Newton point}\index{Newton!point} of
\(\gamma\), that captures the \(F_v\)-valuations (called the (Newton)
\notion{slopes}\index{Newton!slope}) of eigenvalues of \(\gamma\) in \(G\)-representations. See, for
example, \cite{KoVi12}*{\S~2}. In \cite{Ko97}, Kottwitz defines (after fixing an
algebraic closure \(F_v\to\bar{F}_v\), or equivalently, a geometric point
\(\bar{\eta}_{\bar{v}}\) over \(\eta_v\))
 a canonical group homomorphism
\begin{align}
    \kappa_G\colon G(F_v)\longto \CharG(Z(\dual{G})^{I_v})^{\Frob_v}.
    \nomenclature[\(kappa_G \)]{\(\kappa_G\)}{the Kottwitz map \(G(F_v)\to \CharG(Z(\dual{G})^{I_v})^{\Frob_v}\)}
\end{align}
where \(\Frob_v\) is any Frobenius element of \(\Gamma_v\).
Here our group is unramified, so it
simplifies to a homomorphism
\begin{align}
    \kappa_G\colon G(F_v)\longto (\CoCharG(T)/\bbZ\CoRoots)_{F_v}=\pi_1(G)_{F_v}.
\end{align}
This is the \inotion{Kottwitz map}\index{map!Kottwitz}. One can also see \cite{KoVi12}*{\S~3.1} for a
description in this simplified situation.

Let \(p_G\)
    \nomenclature[\(p_G \)]{\(p_G\)}{the canonical map \(\CoCharG(T)\to\pi_1(G)\)}
be the natural quotient \(\CoCharG(T)\to\pi_1(G)\). 
The following observation is crucial:
\begin{lemma}\label[lemma]{lem:find_gamma_lambda_for_gamma}
    Let  \(\gamma\in G^\AD(F_v)\) and \(\lambda\in\CoCharG(T^\AD)^{\Gamma_v}\).  Suppose
    \(\kappa_{G}(\gamma)=p_{G}(\lambda)\). Then there exists an element
    \(\gamma_\lambda\in \Env(G^\SC)^\x(F_v)\) such that
    \begin{enumerate}
        \item the image of \(\gamma_\lambda\) in \(G^\AD(F_v)\) is \(\gamma\)
        \item
            \(\alpha_{\Env(G^\SC)}(\gamma_\lambda)\in\pi^{-w_0(\lambda)}\FRA_{\Env(G^\SC)}^\x(\cO_v)\).
    \end{enumerate}
\end{lemma}
\begin{proof}
    This is basically \cite{Ch22}*{Lemma~3.1.5}, except for the fact that
    residue
    field \(k_v\) is no longer algebraically closed (or \(G\) no longer split).
    The key in Chi's proof is the
    surjectivity of map \(\bG_+(F_v)\to\bG^\AD(F_v)\). After twisting, it is still
    true that \(G_+(F_v)\to G^\AD(F_v)\) is surjective, because the kernel is the
    unramified torus \(T^\SC\) and the Frobenius acts by permuting a basis of
    \(\CharG(T^\SC)\). In other words, \(T^\SC\) is an induced torus, 
    thus \(\RH^1(F_v, T)=0\). The remaining part of the proof
    then proceeds easily.
\end{proof}

\subsection{}%
\label{sub:non_emptyness}

The proposition below settles
the question when the sets \eqref{eqn:KV_def_as_sets} and
\eqref{eqn:general_KV_def_as_sets} are non-empty:
\begin{proposition}[\cite{Ch22}*{Proposition~3.1.6}]\label[proposition]{prop:KV_nonempty_criteria}
    Suppose \(\gamma\in G(F_v)\) (resp.~\(G^\AD(F_v)\))
    is regular semisimple and \(\lambda\in\CoCharG(T)^{\Gamma_v}\)
    (resp.~\(\CoCharG(T^\AD)^{\Gamma_v}\)), then the followings are
    equivalent:
    \begin{enumerate}
        \item \label{item:KV_nonempty_criterion_raw}
            The set of \(\bar{k}\)-points in \(\KV_{G,v}^\lambda(\gamma)\) is non-empty;
        \item \label{item:KV_nonempty_criterion_less_equal}
            The set of \(\bar{k}\)-points in \(\KV_{G,v}^{\le\lambda}(\gamma)\) is non-empty;
        \item \label{item:KV_nonempty_criterion_Hodge_Newton_Kottwitz}
            \(\kappa_{G}(\gamma)=p_{G}(\lambda)\), and
            \(\nu_{\gamma}\le_\bbQ \lambda\), the latter meaning
            \(\lambda-\nu_\gamma\) is a non-negative \(\bbQ\)-combination of
            positive roots;
        \item \label{item:KV_nonempty_criterion_Kottwitz_monoid}
            \(\kappa_{G}(\gamma)=p_{G}(\lambda)\), and
            \(\chi_{\Env(G^\SC)}(\gamma_\lambda)\in \FRC_{\Env(G^\SC)}(\cO_v)\),
            where \(\gamma_\lambda\) is as defined in
            \Cref{lem:find_gamma_lambda_for_gamma}.
    \end{enumerate}
\end{proposition}

\subsection{}
Recall that the \notion{affine
Grassmannian}\index{Grassmannian!affine}\index{affine!Grassmannian} \(\Gr_{G,v}\) 
is the functor sending a \(k\)-scheme \(S\) to the set
of pairs \((E,\phi)\) where \(E\) is a \(G\)-torsor on \(X_{v,S}=X_v\hat{\x}S\) and
\(\phi\) is a trivialization of \(E\) over
\(X_{v,S}^\bullet=X_v^\bullet\hat{\x}S\). Here \(X_{v,S}\) is the completion of
\(X_v\x_k S\) at \({v}\x_k S\), and \(X_{v,S}^\bullet\) is the (open) complement
of \(\Set{v}\x_k S\). \textit{A priori}, \(X_{v,S}\) is only a formal scheme,
not a scheme, but in this simple case it is easy to see it is representable by
a scheme. It is
well-known that \(\Gr_{G,v}\) is represented by an ind-projective ind-scheme of
ind-finite-type over \(k\). If we fix a \(k\)-embedding \(k_v\to\bar{k}\), we have an
analogous definition of affine Grassmannian \(\Gr_{G,\bar{v}}\) over
\(\bar{k}\). We have the isomorphism
\begin{align}
    \Gr_{G,v}\x_{\Spec{k}}\Spec{\bar{k}}\simeq \prod_{\bar{v}\colon k_v\to\bar{k}}\Gr_{G,\bar{v}}.
\end{align}
Since \(G\) has connected
fibers, we may also regard \(\Gr_{G,v}\) as the quotient sheaf
\(\Loop_v{G}/\Arc_v{G}\),
where \(\Arc_v{G}\) is the arc space functor of \(G\) sending \(S\) to the set
\(G(X_{v,S})\), and \(\Loop_v{G}\) is the loop space functor sending \(S\) to 
\(G(X_{v,S}^\bullet)\).
The reader may consult \cite{Zh17} for a comprehensive review or
\Cref{sec:review_BD_affine_grassmannian} for a basic summary.

\begin{remark}
    The construction and representability of
    affine Grassmannian \(\Gr_{G,v}\) requires only \(G\) to be a smooth affine
    group scheme, not necessarily reductive. The arc space and loop space
    functors both make sense for any \(k\)-scheme, and in some sense 
    behave ``well enough'' for normal varieties.
\end{remark}

\subsection{}%
We now impose algebraic structures on \eqref{eqn:KV_def_as_sets} promised at
the beginning. Following \cite{Ch22}, we introduce two approaches to this, one
of which relates to reductive monoids.

We regard the Cartan double coset \(\Cartan_{G,v}^\lambda\) as a subsheaf
\(\Arc_v{G}\pi^\lambda\Arc_v{G}\subset \Loop_v{G}\). We define a sheaf
\(\cM_{G,v}^\lambda(\gamma)\)
    \nomenclature[\(M{}gamma"cal_G_v_lambda \)]{\(\cM_{G,v}^\lambda(\gamma)\)}{the
    big-cell locus of \(\cM_{G,v}^{\le\lambda}(\gamma)\)}
over \(k\) to be the sheaf associated with presheaf
\begin{align}
    \Spec{R}\longmapsto \Set*{g\in\Gr_{G,v}(R)\given g^{-1}\gamma g\in
    \Cartan_{G,v}^\lambda(R)},
\end{align}
which has an ind-scheme structure. Similarly, we define
\(\cM_{G,v}^{\le\lambda}(\gamma)\).
    \nomenclature[\(M{}gamma"cal_G_v_lambda_le \)]{\(\cM_{G,v}^{\le\lambda}(\gamma)\)}{the
    multiplicative affine Springer fiber over \(k\) associated with \(\gamma\in G(F_v)\) and type \(\lambda\)}
The algebraic structure on
\(\KV_{G,v}^\lambda(\gamma)\) is then the reduced ind-subscheme structure
induced by \(\cM_{G,v}^\lambda(\gamma)\). Similarly, we have an algebraic
structure on \(\KV_{G,v}^{\le \lambda}(\gamma)\).

Another approach to the algebraic structure on \(\KV_{G,v}^\lambda(\gamma)\) is to
use reductive monoids.
Let \(\FRM\in\FM(G^\SC)\) and \(\gamma_{\FRM}\in
\FRM^\x(F_v)^\rss\) such that \(a=\chi_{\FRM}(\gamma_{\FRM})\in\FRC_{\FRM}(\cO_v)\).
\begin{definition}
    We make multiple closely related definitions:
    \begin{enumerate}
        \item The \notion{multiplicative affine Springer fiber}\index{multiplicative!affine Springer fiber}
            \(\cM_{G,v}(\gamma_{\FRM})\)
            \nomenclature[\(M{}gamma"cal_M_G_v \)]{\(\cM_{G,v}(\gamma_\FRM)\)}{the
            multiplicative affine Springer fiber over \(k\) associated with \(G\) and \(\gamma_\FRM\in \FRM(\cO_v)\cap\FRM^\x(F_v)\)}
            associated with \(\gamma_{\FRM}\) is
            the functor that associates to \(k\)-scheme \(S\) the
            isomorphism classes of pairs \((h,\iota)\) where \(h\) is an
            \(X_{v,S}\)-point of \(\Stack*{\FRM/G}\) over \(a\):
            \begin{equation}
                \begin{tikzcd}
                    X_{v,S}\ar[r, "h"] \ar[d] &
                    \Stack*{\FRM/G}\ar[d, "\chi_{\FRM}"]\\ X_v\ar[r,
                    "a"] & \FRC_{\FRM}
                \end{tikzcd},
            \end{equation}
            and \(\iota\) is an isomorphism between the restriction of \(h\) to
            \(X_{v,S}^\bullet\) and the \(X_{v,S}^\bullet\)-point of
            \(\Stack*{\FRM/G}\) induced by \(\gamma_{\FRM}\).
        \item The (open) subfunctors \(\SP^\circ_{G,v}(\gamma_{\FRM})\)
            \nomenclature[\(M{}gamma"cal_M_G_v_circ \)]{\(\cM_{G,v}^{\circ}(\gamma_\FRM)\)}{the
            big-cell locus of \(\cM_{G,v}(\gamma_\FRM)\)}
            (resp.~\(\SP^\reg_{G,v}(\gamma_{\FRM})\)) is defined similarly by
            replacing \(\FRM\) with \(\FRM^\circ\) (resp.~\(\FRM^\reg\)). 
        \item For \(a\in\FRC_{\FRM}(\cO_v)\cap \FRC_{\FRM}^\x(F_v)\), we denote by
            \(\SP_{G,v}(a)\)
            \nomenclature[\(M{}a"cal_G_v \)]{\(\cM_{G,v}(a),\cM_v(a)\)}{the
            multiplicative affine Springer fiber over \(k\) associated with \(G\) and \(a\in\FRC_\FRM(\cO_v)\cap\FRC_\FRM^\x(F_v)\)}
            (resp.~\(\SP_{G,v}^\circ(a)\), etc.)
            \nomenclature[\(M{}a"cal_G_v_circ \)]{\(\cM_{G,v}^{\circ}(a)\)}{the big-cell locus of \(\cM_{G,v}(a)\)}
            a fixed
            choice of \(\SP_{G,v}(\gamma_{\FRM})\)
            (resp.~\(\SP_{G,v}^\circ(\gamma_{\FRM})\), etc.) where
            \(\chi_{\FRM}(\gamma_{\FRM})=a\) (which always exists thanks to
            \Cref{thm:A_2m_rationality_fix}). The dependence on such choice will
            be emphasized if not clear from context.
        \item When the group \(G\) is clear from context, we drop it from
            subscripts.
    \end{enumerate}
\end{definition}

Let \(\gamma\in G^\AD(F_v)^\rss\) and \(\lambda\in \sC(\FRA_{\FRM})\) a dominant
\(F_v\)-cocharacter. Recall we have \(\lambda_\AD\in\CoCharG(T^\AD)_+\). Suppose
\(\KV_{G,v}^{\lambda_\AD}(\gamma)\) is non-empty, then by
\Cref{prop:KV_nonempty_criteria}, we have the element
\(\gamma_\lambda\in \Env(G^\SC)^\x(F_v)^\rss\) defined by
\Cref{lem:find_gamma_lambda_for_gamma}. We may lift \(\gamma_\lambda\) to
an element in \(\FRM^\x(F_v)\), still denoted by \(\gamma_\lambda\).  Let
\(a=\chi_{\FRM}(\gamma_\lambda)\), then it is easy to see that we have
homeomorphisms
\begin{align}
    \label{eqn:identify_KV_with_SP}
    \cM_{G,v}^{\le\lambda_\AD}(\gamma)\cong \SP_{G,v}(a)\text{ and }
    \cM_{G,v}^{\lambda_\AD}(\gamma)\cong \SP_{G,v}^{\circ}(a).
\end{align}

Conversely, for any \(a\in\FRC_{\FRM}(\cO_v)\cap \FRC_{\FRM}^\x(F_v)^\rss\) whose image in
\(\FRA_{\Env(G^\SC)}\) is contained in \(\pi^{-w_0(\lambda_\AD)} T^\AD(\cO_v)\).
Suppose \(\gamma_{\FRM}\) maps to \(a\) under \(\chi_{\FRM}\) (such element always
exists by \Cref{thm:A_2m_rationality_fix}). Suppose \(\gamma\in G^\AD(F_v)\)
is the image of \(\gamma_{\FRM}\) in \(G^\AD(F_v)\), then
\eqref{eqn:identify_KV_with_SP} still holds. Note that a \(k\)-point 
may not exist for such functors. 

\subsection{}
Like the affine Grassmannian, we may also replace \(X_v\) by
\(\breve{X}_{\bar{v}}\), so we have \(\cM_{G,\bar{v}}^{\le\lambda}(\gamma)\),
\(\cM_{G,\bar{v}}^\lambda(\gamma)\), and so
on. If we base change \(\cM_{G,v}^\lambda(\gamma)\) to \(\bar{k}\), we obtain
an isomorphism
\begin{align}
    \cM_{G,v}^{\le\lambda}(\gamma)\x_{\Spec{k}}\Spec{\bar{k}}\simeq \prod_{\bar{v}\colon k_v\to
    \bar{k}}\cM_{G,\bar{v}}^{\le\lambda}(\gamma),
\end{align}
and similarly for \(\SP_{G,v}(a)\), and so on.

\subsection{}
\label{sub:base_point_for_MASF_general}
For our purposes in this book, we also need a variant of multiplicative affine
Springer fibers by allowing a central twist to the monoid.
Let \(\cL\) be a \(Z_{\FRM}\)-torsor over \(\cO_v\), viewed as a
\(\cO_v\)-point of \(\BG{Z_\FRM}\). Given a point
\begin{align}
    a\in\Stack*{\FRC_\FRM/Z_\FRM}(\cO_v)\cap\Stack*{\FRC_{\FRM}^{\x,\rss}/Z_\FRM}(F_v)
\end{align}
lying over \(\cL\), it induces a point
\(a_0\in\Stack*{\FRC_{\Env(G^\SC)}/Z_{\Env(G^\SC)}}(\cO_v)\). Since
\(Z_{\Env(G^\SC)}\) is an induced torus, its torsors are trivial over \(\cO_v\),
so \(a_0\) lifts to a point in \(\FRC_{\Env(G^\SC)}(\cO_v)\), which in turn lifts to an
\(F_v\)-point in \(\Env(G^\SC)\). Using the Cartesian diagram
\begin{equation}
    \begin{tikzcd}
        \Stack*{\FRM/Z_\FRM} \ar[r]\ar[d] & \Stack*{\Env(G^\SC)/Z_{\Env(G^\SC)}} \ar[d]\\
        \Stack*{\FRC_\FRM/Z_\FRM} \ar[r] &
        \Stack*{\FRC_{\Env(G^\SC)}/Z_{\Env(G^\SC)}}
    \end{tikzcd}
\end{equation}
we obtain a point \(\gamma_v\in\Stack*{\FRM/Z_\FRM}(F_v)\) lying over \(a\).
By replacing \(\FRM\) with \(\FRM_\cL\), we may also define \(\SP_{G,v}(a)\) for
\(a\) using \(\gamma_v\).

\section{Dimension Formula}%
\label{sec:dimension_formula}

In the following few sections we consider some geometric properties of
multiplicative affine Springer fibers, and thus it is harmless to base change to
\(\bar{k}\) for our discussion.

Suppose we have a non-empty (as a sheaf)
multiplicative affine Springer fiber \(\SP_{G,\bar{v}}^\lambda(\gamma)\), then
\(\gamma_\lambda\) as in \Cref{lem:find_gamma_lambda_for_gamma} exists, and
let \(a=\chi_{\Env(G^\SC)}(\gamma_\lambda)\). Let
\(d_{\bar{v}+}(a)=\deg_{\bar{v}}(\FRD_a)\) be the degree of extended discriminant
divisor at \(\bar{v}\), which is equal to
\begin{align}
    \Pair{2\rho}{\lambda}+d_{\bar{v}}(\gamma),
\end{align}
where \(d_{\bar{v}}(\gamma)\) is the (non-extended) discriminant valuation
\begin{align}
    d_{\bar{v}}(\gamma)=d_{\bar{v}}(a)\defeq\sum_{\Rt\in\Roots}\val_{F_{v}}(1-\Rt(x_a^\bullet))
\end{align}
for any \(x_a^\bullet\in T_\FRM(F_v^\sep)\) lying over \(a\).
For a fixed choice of \(x_a^\bullet\), it determines a trivialization of the
\(W\)-torsor \(\breve{X}_{\bar{v}}\x_{a,\FRC_\FRM}\FRT_\FRM\), or
in other words, a homomorphism
\begin{align}
    \pi_a^\bullet\colon I_v\longto W.
\end{align}
By our assumption of
\(\Char(k)\), this map factors through the tame inertia group, hence has cyclic
image in \(W\). Define
\begin{align}
    c_{\bar{v}}(a)\defeq \dim\La{t}-\dim\La{t}^{\pi_a^\bullet(I_v)},
\end{align}
and the \notion{local \(\delta\)-invariant}\index{\(\delta\)-!invariant, local}\index{local!\(\delta\)-invariant}
\begin{align}
    \delta_{\bar{v}}(a)\defeq\frac{d_{\bar{v}+}(a)-c_{\bar{v}}(a)}{2}.
\end{align}
Note that \(c_{\bar{v}}\) (resp.~\(\delta_{\bar{v}}\)) is the same as \(c\)
(resp.~\(\delta\)) defined in \Cref{sec:cylinders_in_reductive_monoids}.

\begin{theorem}[\cite{Ch22}*{Theorem~1.2.1}]
    \label[theorem]{thm:GASF_dimension_formula}
    The functors \(\SP_{G,\bar{v}}^{\lambda,\Red}(\gamma)\) and
    \(\SP_{G,\bar{v}}^{\le\lambda,\Red}(\gamma)\)
    are represented by an equidimensional \(\bar{k}\)-scheme locally of finite type, 
    with dimension
    \begin{align}
        \dim\SP_{G,\bar{v}}^\lambda(\gamma)=\dim\SP_{G,\bar{v}}^{\le\lambda}(\gamma)
        &=\Pair{\rho}{\lambda}+\frac{d_{\bar{v}}(\gamma)-c_{\bar{v}}(\gamma)}{2}\\
        &=\frac{d_{\bar{v}+}(a)-c_{\bar{v}}(a)}{2}\\
        &=\delta_{\bar{v}}(a).
    \end{align}
\end{theorem}

\section{Symmetry and Irreducible Components}%
\label{sec:Symmetry and Irreducible Components}

Given \(a\in\FRC_{\FRM}(\cO_v)\cap \FRC_{\FRM}^\x(F_v)^\rss\), one can pull back the
regular
centralizer group 
scheme \(\FRJ_{\FRM}\) to \(a\), viewed as an \(\cO_v\)-scheme denoted by
\(\FRJ_a\).
    \nomenclature[\(J"frak_a \)]{\(\FRJ_a\)}{the regular centralizer scheme induced by either a local point \(a\in\FRC_\FRM(\cO_v)\) or a global point \(a\in\cA_X\)}
Define the \notion{local Picard}\index{Picard!functor, local} functor as the affine Grassmannian
\begin{align}
    \cP_v(a)\defeq \Gr_{\FRJ_a,v}.
    \nomenclature[\(P"cal_v_a \)]{\(\cP_v(a)\)}{the local Picard functor \(\Gr_{\FRJ_a,v}\)}
\end{align}
This group functor naturally acts on \(\SP_v(a)\) as follows:
any \(S\)-point of \(\cM_v(a)\) is a tuple \((E,\phi,\iota)\) where \(E\) is a
\(G\)-bundle on \(X_{v,S}\), \(\phi\) is a \(G\)-equivariant map \(E\to\FRM\),
and \(\iota\) is an isomorphism \((E_0,\gamma_\FRM)\to (E,\phi)\) over the
punctured disc \(X_{v,S}^\bullet\), where \(E_0\) is the trivial \(G\)-torsor.
On the other hand, a point of \(\cP_v(a)\) is a \(\FRJ_a\)-torsor \(E_{\FRJ}\)
on \(X_{v,S}\) with a trivialization \(\tau\) over \(X_{v,S}^\bullet\). It sends
\((E,\phi,\iota)\) to the following tuple \((E',\phi',\iota')\): the pair
\((E',\phi')\) is defined as
\begin{align}
    \phi'\colon E'\defeq E\x_{\phi,\FRM}^{\FRJ_a}E_{\FRJ}\longto \FRM,
\end{align}
where \(\FRJ_a\) acts on the fibers of \(\phi\) through the canonical map
\(\chi_{\FRM}^*\FRJ_{\FRM}\to I_{\FRM}\). The trivialization \(\tau\) induces an
isomorphism
\begin{align}
    E|_{X_{v,S}^\bullet}\longto E'|_{X_{v,S}^\bullet},
\end{align}
whose composition with \(\iota\) is \(\iota'\).

\subsection{}
In contrast to the Lie algebra case,
there may be more than one open \(\cP_v(a)\)-orbit in
\(\SP_v(a)\). Even
worse, the union of all open orbits is not dense in \(\SP_v(a)\), meaning there are
some irreducible components of \(\SP_v(a)\) that are stratified into infinitely many
orbits.

Fortunately, the free action of a sublattice of \(\cP_v(a)\) is still
present (see \Cref{sec:neron_models_and_connected_components}), and so the
action of \(\cP_v(a)\) on the set
\(\Irr(\SP_v(a))\) of irreducible components is still quite nice.
We shall see in
\Cref{sec:connection_with_mv_cycles,sec:Connection_with_Kashiwara_Crystals} that
the latter has deep connection with both representations of \(L\)-groups and the
transfer factor.

\subsection{}
We may also base change to \(\bar{k}\), then we have
\begin{align}
    \cP_v(a)\x_{\Spec{k}}\Spec{\bar{k}}\simeq\prod_{\bar{v}\colon k_v\to \bar{k}}\cP_{\bar{v}}(a).
\end{align}
The action of \(\cP_v(a)\) on \(\SP_v(a)\) is compatible with each direct factor.
So starting from this point we will base change everything to \(\bar{k}\) and 
consider \(\SP_{\bar{v}}(a)\) for a fixed \(\bar{v}\) (equivalently, we replace
\(X_v\) with \(\breve{X}_{\bar{v}}\)).

\subsection{}
Suppose \(\gamma\in G(\breve{F}_{\bar{v}})\) is such that the Newton point
\(\nu=\nu_\gamma\) is
an integral cocharacter in \(\CoCharG(T)\). We may consider more generally
\(\gamma\in G^\AD(\breve{F}_{\bar{v}})\). This happens when \(\gamma\) is
unramified, but not exclusively.  Let \(m_{\lambda\nu}\)
    \nomenclature[\(m_lambda_mu \)]{\(m_{\lambda\nu}\)}{the weight multiplicity of weight \(\nu\) in \(V_\lambda\)}
be the weight
multiplicity of the weight space \(\nu\) in the irreducible representation of
\(\dual{\bG}\) with highest-weight \(\lambda\).  If \(\nu\in W\lambda\), 
then \(m_{\lambda\nu}=1\). When \(\lambda\) is a central cocharacter, in other
words, \(\Pair{\alpha}{\lambda}=0\) for all \(\alpha\in\Roots\), or
equivalently, \(\lambda_\AD=0\), and
\(\SP_{G,\bar{v}}^\lambda(\gamma)\neq\emptyset\), then we always have \(\nu=\lambda\) and
so \(m_{\lambda\nu}=1\).

Using \eqref{eqn:identify_KV_with_SP}, the condition of \(\lambda\) being
central corresponds to that the image of \(a\) in \(\FRC_{\Env(G^\SC)}\) is
contained in \(\FRC_{\Env(G^\SC)}^\x\). The condition \(\gamma\)
being unramified corresponds to \(c_{\bar{v}}(a)=0\), because \(c_{\bar{v}}(a)\) is the
difference between the absolute rank and \(\breve{F}_{\bar{v}}\)-split rank of the
centralizer of \(\gamma_\lambda\).

\begin{definition}
    With the notations above, we call
    \(a\in\FRC_{\FRM}(\cO_{v})\cap\FRC_{\FRM}^\x(F_{v})^\rss\)
    \notion{unramified}\index{orbit!unramified}\index{unramified!orbit} if \(c_{\bar{v}}(a)=0\) for one (equivalently, all)
    \(\bar{v}\), \notion{invertible}\index{invertible!orbit}\index{orbit!invertible} if it is
    contained in \(\FRC_{\FRM}^\x(\cO_v)\), and
    \notion{central}\index{central!orbit}\index{orbit!central} if its image in
    \(\FRC_{\Env(G^\SC)}\) is contained in \(\FRC_{\Env(G^\SC)}^\x\).
\end{definition}

\begin{theorem}[\cite{Ch22}*{Corollaries~3.5.3 and 3.8.2}]
    \label[theorem]{thm:local_irr_components_weight_mult}
    Suppose \(a\in\FRC_{\FRM}(\breve{\cO}_{\bar{v}})\cap\FRC_{\FRM}^\x(\breve{F}_{\bar{v}})\)
    is either unramified or central, 
    then the number of irreducible components in
    \(\SP_{\bar{v}}(a)\) modulo the action of \(\cP_{\bar{v}}(a)\) is
    \(m_{\lambda\nu}\). Moreover, when \(a\) is central, there is a
    unique open dense \(\cP_{\bar{v}}(a)\)-orbit in \(\SP_{\bar{v}}(a)\) being the
    regular locus \(\SP_{\bar{v}}^\reg(a)\).
\end{theorem}

\begin{corollary}[\cite{Ch22}*{Corollary~3.8.3}]
    \label[corollary]{cor:central_lambda_and_d_1_implies_reg}
    Suppose \(a\) is central and \(d_+(a)\le 1\), then
    \(\SP_{\bar{v}}(a)=\SP_{\bar{v}}(a)^\reg\) and is itself a
    \(\cP_{\bar{v}}(a)\)-torsor.
\end{corollary}

The weight multiplicity \(m_{\lambda\nu}\) in
\Cref{thm:local_irr_components_weight_mult} is obtained by connecting
\(\SP_{\bar{v}}(a)\) with MV-cycles, which we will elaborate in
\Cref{sec:connection_with_mv_cycles}.

\section{N\'eron Models and Connected Components}%
\label{sec:neron_models_and_connected_components}

Following \cite{Ng10} and \cite{Ch22}, we have a more precise description of the
\(\cP_{\bar{v}}(a)\) using N\'eron models. By definition, the N\'eron model
of \(\FRJ_a\) is a unique (up to a unique isomorphism) smooth
group scheme \(\FRJ_a^\flat\)
    \nomenclature[\(J"frak_a^flat \)]{\(\FRJ_a^\flat\)}{the N\'eron model of \(\FRJ_a\)}
over \(\breve{\cO}_{\bar{v}}\) together with a
homomorphism \(\FRJ_a\to\FRJ_a^\flat\) that is an isomorphism over the generic
point and satisfies the following universal property: for any smooth 
group scheme \(J\) over \(\breve{\cO}_{\bar{v}}\) with an \(\breve{F}_{\bar{v}}\)-isomorphism
\(J_{\breve{F}_{\bar{v}}}\to\FRJ_{a,\breve{F}_{\bar{v}}}^\flat\), there is a
canonical lift to a \(\breve{\cO}_{\bar{v}}\)-homomorphism \(J\to\FRJ_a^\flat\).
N\'eron model necessarily exists for smooth commutative \(\breve{\cO}_{\bar{v}}\)-group
schemes, and can be explicitly constructed using cameral covers.

\subsection{}
Recall we have the cameral cover \(\pi_{\FRM}\colon\FRT_{\FRM}\to\FRC_{\FRM}\). The group
scheme \(\FRJ_{\FRM}\) is canonically isomorphic to an open subgroup of the fixed
point scheme
\begin{align}
    \FRJ_{\FRM}^1=\left(\pi_{\FRM*}(T\x\FRT_{\FRM})\right)^W.
\end{align}
The scheme \(\FRJ_a\) is the pullback of \(\FRJ_{\FRM}\) via \(a\).
Let local cameral cover \(\pi_a\colon\tilde{X}_a\to \breve{X}_{\bar{v}}\) 
    \nomenclature[\(pi_a \)]{\(\pi_a\)}{the pullback of the cameral cover \(\pi_\FRM\) to \(a\in\FRC_\FRM(\cO_v)\) or \(a\in\cA_X\)}
    \nomenclature[\(X'tilde_a \)]{\(\tilde{X}_a\)}{the cameral curve induced by  \(a\in\FRC_\FRM(\cO_v)\) or \(a\in\cA_X\)}
be the pullback of \(\pi_{\FRM}\) through \(a\). Then by proper base change, one can
also define \(\FRJ_a\) and \(\FRJ_a^1\) using the same construction applied to
\(\pi_a\).

The cameral cover \(\pi_{\FRM}\) is flat with a Cohen--Macaulay source 
and a regular target. Thus,
\(\tilde{X}_a\) is Cohen--Macaulay as well. Since \(a\) is generically
regular semisimple, \(\tilde{X}_a^\bullet\) is a regular scheme, being an
\'etale \(W\)-cover of \(\breve{X}_{\bar{v}}^\bullet\). So
\(\tilde{X}_a\) is a reduced scheme.
Since \(\tilde{X}_a\) is one-dimensional, its
normalization \(\tilde{X}_a^\flat\) is regular. The N\'eron model can be
shown to be the group scheme
\begin{align}
    \FRJ_a^\flat=\pi_{a*}^\flat(T\x\tilde{X}_a^\flat)^W,
\end{align}
where \(\pi_a^\flat\)
    \nomenclature[\(pi_a^flat \)]{\(\pi_a^\flat\)}{the normalization of the
    local or global cameral cover \(\pi_a\)}
    \nomenclature[\(X'tilde_a^flat \)]{\(\tilde{X}_a^\flat\)}{the normalization of \(\tilde{X}_a\)}
is the natural map \(\tilde{X}_a^\flat\to
\breve{X}_{\bar{v}}\). The same proof in \cite{Ng10}*{\S~3.8.2} applies to the case here.

\begin{lemma}
    \label[lemma]{lem:local_delta_formula_by_Neron}
    We have formula
    \begin{align}
        \delta_{\bar{v}}(a)=\dim{\cP_{\bar{v}}(a)} =
        \dim_{\bar{k}}\bigl(\La{t}_{\breve{\cO}_{\bar{v}}}\otimes_{\breve{\cO}_{\bar{v}}}(\tilde{\cO}_{a}^\flat/\tilde{\cO}_{a})\bigr)^W,
    \end{align}
    where \(\tilde{\cO}_{a}\) (resp.~\(\tilde{\cO}_{a}^\flat\)) is
    the ring of functions of \(\tilde{X}_a\) (resp.~\(\tilde{X}_a^\flat\)).
\end{lemma}
\begin{proof}
    This is essentially \cite{Ng10}*{Corollaire~3.8.3} and
    \cite{Ch22}*{Corollary~3.3.4}. The key is that \(\FRJ_a\) and
    \(\FRJ_a^\flat\) are smooth, so we only need to compute the dimension of the 
    tangent space at
    the identity of \(\cP_{\bar{v}}\), and then it is clear from the Galois
    descriptions of \(\FRJ_a\) and \(\FRJ_a^\flat\).
\end{proof}

\subsection{}
The morphism \(\FRJ_a\to\FRJ_a^\flat\) induces morphism of group ind-schemes
\begin{align}
    p_{\bar{v}}\colon\cP_{\bar{v}}(a)\longto \cP_{\bar{v}}^\flat(a)\defeq
    \Gr_{\FRJ_a^\flat,\bar{v}}.
    \nomenclature[\(P"cal_v_a_flat \)]{\(\cP_v^\flat(a)\)}{the local Picard functor \(\Gr_{\FRJ_a^\flat,v}\) of the N\'eron model}
\end{align}

With exact same proof, we have:
\begin{lemma}[\cite{Ng10}*{Lemme~3.8.1}]
    \label[lemma]{lem:kernel_of_local_to_Neron}
    The group \(\cP_{\bar{v}}^\flat(a)\) is homeomorphic to a finitely generated
    free abelian group (viewed as a discrete \(\bar{k}\)-scheme). The map
    \(p_{\bar{v}}\) is surjective, and its kernel \(\cR_{\bar{v}}(a)\)
    \nomenclature[\(R"cal_v_a \)]{\(\cR_{v}(a)\)}{the kernel of the map \(\cP_v(a)\to\cP_v^\flat(a)\)}
    is an
    affine group scheme of finite type over \(\bar{k}\).
\end{lemma}

\begin{corollary}
    The dimension of \(\cR_{\bar{v}}(a)\) is exactly the local
    \(\delta\)-invariant \(\delta_{\bar{v}}(a)\).
\end{corollary}

\subsection{}
Since \(\pi_0(\cP_{\bar{v}}(a))^\Red\) is a finitely generated abelian group and
\(\cR_{\bar{v}}(a)\)
is affine of finite type, \(\pi_0(\cP_{\bar{v}}^\flat(a))\) must be homeomorphic
to the largest free quotient of \(\pi_0(\cP_{\bar{v}}(a))\). Call this lattice
\(\Lambda_a\). Since \(\Lambda_a\) is free, we can choose a lifting of it to
\(\cP_{\bar{v}}(a)\) to define an action of \(\Lambda_a\) on
\(\SP_{\bar{v}}(a)\).
\begin{proposition}
    \label[proposition]{prop:projective_quotient_of_GASF}
    The action of \(\Lambda_a\) on \(\SP_{\bar{v}}^\Red(a)\) is free, and the
    quotient \(\SP_{\bar{v}}^\Red(a)/\Lambda_a\) is a projective
    \(\bar{k}\)-scheme.
\end{proposition}
\begin{proof}
    This is essentially \cite{Ch22}*{Theorem~3.6.2}, in which the freeness is
    proved for an explicitly constructed lifting
    \(\Lambda_a\to\SP_{\bar{v}}(a)\), and the quotient
    \(\SP_{\bar{v}}^\Red(a)/\Lambda_a\) is shown to be a proper algebraic space.
    So we only need to slightly strengthen these statements.

    First, freeness holds for an arbitrary lift of \(\Lambda_a\), because using
    \cite{Ch22}*{Theorem~3.6.2}, the stabilizer of any point in
    \(\SP_{\bar{v}}(a)\) is a subgroup of \(\Lambda_a\) of finite type, hence
    trivial because \(\Lambda_a\) is a lattice.

    Moreover, in \textit{loc. cit.}, Chi showed that \(\SP_{\bar{v}}^\Red(a)\)
    is the union of \(\Lambda_a\)-translations of a fundamental domain which is
    an algebraic variety. This union is locally finite in the sense that the
    fundamental domain only has non-trivial intersection with finitely many
    translated copies. The last claim is straightforward by applying
    Kazhdan--Lusztig's argument \cite{KaLu88}*{Propositions~2.1 and 3.1} to
    Chi's description of \(\SP_{\bar{v}}^\Red(a)\) (which is itself an
    adaptation of the Lie algebra analogue in \cite{KaLu88}).

    As a
    result, the quotient space \(\SP_{\bar{v}}^\Red(a)/\Lambda_a\) is a scheme by
    \cite{DeGa70}*{Th\'eor\`em~3.2}, because \(\Lambda_a\) is a
    \(\bar{k}\)-group scheme locally of finite type, acting freely and
    \emph{locally algebraically} (meaning the orbit of any point has only
    finite intersection with any affine open subset of
    \(\SP_{\bar{v}}^\Red(a)\)) on \(\SP_{\bar{v}}^\Red(a)\).
\end{proof}

\begin{corollary}
    The stabilizers of the \(\cP_{\bar{v}}^\Red(a)\)-action
    are affine and contained in \(\cR_{\bar{v}}(a)\).
\end{corollary}

\subsection{}
Following \cite{Ng10}*{\S~3.9}, we give a precise description of connected
components of \(\cP_{\bar{v}}(a)\).
Recall we have open subscheme \(\FRJ_a^0\subset\FRJ_a\) of fiberwise neutral
component. The quotient \(\FRJ_a/\FRJ_a^0\) is supported over the closed point
of \(\breve{X}_{\bar{v}}\), with fiber
\begin{align}
    \pi_0(\FRJ_a)=\FRJ_a(\breve{\cO}_{\bar{v}})/\FRJ_a^0(\breve{\cO}_{\bar{v}}).
\end{align}

We have exact sequence
\begin{align}
    1\longto \pi_0(\FRJ_a)\longto
    \FRJ_a(\breve{F}_{\bar{v}})/\FRJ_a^0(\breve{\cO}_{\bar{v}})\longto
    \FRJ_a(\breve{F}_{\bar{v}})/\FRJ_a(\breve{\cO}_{\bar{v}})\longto 1.
\end{align}
Since \(\bar{k}\) is algebraically closed, this is the same as exact sequence
\begin{align}
    1\longto \pi_0(\FRJ_a)\longto
    \cP_{\bar{v}}^0(a)(\bar{k})\longto
    \cP_{\bar{v}}(a)(\bar{k})\longto 1,
\end{align}
where \(\cP_{\bar{v}}^0(a)=\Gr_{\FRJ_a^0,\bar{v}}\).
In other words, the morphism \(\cP_{\bar{v}}^0(a)\to\cP_{\bar{v}}(a)\) is
surjective. Thus, we have exact sequence
\begin{align}
    \label{eqn:pi_0_cP_a_rough}
    \pi_0(\FRJ_a)\longto \pi_0(\cP_{\bar{v}}^0(a))\longto
    \pi_0(\cP_{\bar{v}}(a))\longto 1.
\end{align}

For any finitely generated abelian group \(\Lambda\), we let \(\Lambda^\odot\)
    \nomenclature[\(.odot \)]{\((\cdot)^\odot\)}{Cartier duality between finitely
    generated abelian groups and diagonalizable goups over \(\Qlb\)}
be its \(\Qlb\)-Cartier dual
\begin{align}
    \Lambda^{\odot}=\Spec{\Qlb[\Lambda]},
\end{align}
and conversely, for any diagonalizable group \(A\) over \(\Qlb\), we let
\(A^\odot=\CharG(A)\) be its character group.

Now we fix a trivialization of \(\OGT_G\) over \(\breve{X}_{\bar{v}}\), which
identifies \(G\) with \(\bG\) together with associated pinnings. Over
\(\breve{X}_{\bar{v}}^\bullet\), the cameral cover is a \(\bW\)-\'etale cover. If
we fix a geometric point in \(\tilde{X}_{a}\) over the geometric generic
point \(\bar{\eta}_{\bar{v}}\in \breve{X}_{\bar{v}}\), we have a homomorphism
\(\pi_a^\bullet\colon I_v\to\bW\) (recall that \(I_v\) is the inertia group of
\(F_v\)).

\begin{proposition}
    \label[proposition]{prop:local_pi_0_description}
    After fixing a geometric point of \(\tilde{X}_a\) lying over
    \(\bar{\eta}_{\bar{v}}\) as above, we have canonical isomorphisms of
    diagonalizable groups
    \begin{align}
        \pi_0(\cP_{\bar{v}}^0)^\odot&\simeq \dual{\bT}^{\pi_a^\bullet(I_v)},\\
        \pi_0(\cP_{\bar{v}})^\odot&\simeq \dual{\bT}(\pi_a^\bullet(I_v)),
    \end{align}
    where \(\dual{\bT}(\pi_a^\bullet(I_v))\) is a subgroup of
    \(\dual{\bT}^{\pi_a^\bullet(I_v)}\) consisting of elements \(\kappa\)
    such that \(\pi_a^\bullet(I_v)\) is contained in the Weyl group 
    \(\bW_{\bH}\) of the neutral component \(\dual{\bH}\) of the
    centralizer of \(\kappa\) in \(\dual{\bG}\).
\end{proposition}

The proof is similar to \cite{Ng10}*{Proposition~3.9.2}, but there is one
part that must be modified, due to the following subtle difference: in Lie algebra case,
Ng\^o uses the fact that if \(Z_G\) is connected, then the special fiber of
\(\FRJ_a\) is connected, while this is not true in multiplicative case (see
\Cref{rmk:connectedness_of_FRJ_vs_center_of_G}). Therefore, we will give a modified
proof. First we have a lemma concerning
\(\cP_{\bar{v}}^{\flat,0}(a)\) associated with the neutral component of N\'eron
model \(\FRJ_a^{\flat,0}\).

\begin{lemma}[\cite{Ng10}*{Lemme~3.9.3}]
    The homomorphism \(\cP_{\bar{v}}^0(a)\to\cP_{\bar{v}}^{\flat,0}(a)\) induces
    an isomorphism
    \begin{align}
        \pi_0(\cP_{\bar{v}}^0(a))\longto
        \FRJ_a(\breve{F}_{\bar{v}})/\FRJ_a^{\flat,0}(\breve{\cO}_{\bar{v}}).
    \end{align}
\end{lemma}
\begin{proof}
    Both \(\FRJ_a^0\) and \(\FRJ_a^{\flat,0}\) have connected fibers, thus
    following the exact same steps for obtaining \eqref{eqn:pi_0_cP_a_rough}, 
    we have exact sequence
    \begin{align}
        1=\pi_0(\FRJ_a^{\flat,0})
        \longto\pi_0(\cP_{\bar{v}}^0(a))
        \longto \pi_0(\cP_{\bar{v}}^{\flat,0}(a))
        =\FRJ_a(\breve{F}_{\bar{v}})/\FRJ_a^{\flat,0}(\breve{\cO}_{\bar{v}})
        \longto 1
    \end{align}
    as desired.
\end{proof}

\begin{lemma}
    \label[lemma]{lem:pi_0_local_Picard_with_coinv_cochar}
    We have isomorphism
    \begin{align}
        \FRJ_a(\breve{F}_{\bar{v}})/\FRJ_a^{\flat,0}(\breve{\cO}_{\bar{v}})\simeq
        \CoCharG(\bT)_{\pi_a^\bullet(I_v)}.
    \end{align}
\end{lemma}
\begin{proof}
    This is \cite{Ng10}*{Lemme~3.9.4} and is a special case of
    \cite{Ko85}*{Lemma~2.2}. It only uses certain properties of the
    functor \(A\mapsto A(\breve{F}_{\bar{v}})/A^{\flat,0}(\breve{\cO}_{\bar{v}})\) on the
    category of \(\breve{F}_{\bar{v}}\)-tori, as well as the fact that
    \(\FRJ_a|_{\breve{F}_{\bar{v}}}\) is a \(\breve{F}_{\bar{v}}\)-torus.
\end{proof}

\begin{proof}[Proof of \Cref{prop:local_pi_0_description}]
    By duality, we have
    \(\CoCharG(\bT)_{\pi_a^\bullet(I_v)}=\CharG(\dual{\bT})_{\pi_a^\bullet(I_v)}\),
    so by combining the two lemmas above, we obtain
    \begin{align}
        \pi_0(\cP_{\bar{v}}^0(a))\simeq \CoCharG(\bT)_{\pi_a^\bullet(I_v)},
    \end{align}
    and by taking Cartier dual we have
    \begin{align}
        \pi_0(\cP_{\bar{v}}^0(a))^\odot\simeq \dual{\bT}^{\pi_a^\bullet(I_v)}.
    \end{align}
    This is the first claim.

    For the second claim, we use \(z\)-extensions. We have exact sequence of
    reductive group schemes
    \begin{align}
        1\longto G\longto G_1\longto C\longto 1,
    \end{align}
    where \(G_1\) has connected center. It is obtained from its split
    counterparts twisted by an \(\Out(\bG)\)-torsor. We also have exact sequence
    of tori
    \begin{align}
        1\longto T\longto T_1\longto C\longto 1
    \end{align}
    and its split counterpart. Both \(G\) and \(G_1\) act on
    \(\FRM\in\FM(G^\SC)\), and they are compatible. In fact, we have
    \(\FRC_{\FRM}=\FRM\git G\simeq \FRM\git G_1\). Let
    \(\FRJ_1\to\FRC_\FRM\)
    be the regular centralizer scheme associated with \(G_1\)-action. Note that
    contrary to the situation in \cite{Ng10}*{Proposition~3.9.2}, it
    may have disconnected special fiber. We still have exact sequence
    \begin{align}
        1\longto \FRJ\longto \FRJ_1\longto C\longto 1,
    \end{align}
    and the fiber over \(a\)
    \begin{align}
        1\longto \FRJ_a\longto \FRJ_{1,a}\longto C_a\longto 1.
    \end{align}
    Since \(\FRJ_a\to\FRJ_{1,a}\) is proper,
    \(\FRJ_a(\breve{\cO}_{\bar{v}})=\FRJ_a(\breve{F}_{\bar{v}})\cap\FRJ_{1,a}(\breve{\cO}_{\bar{v}})\).
    This implies that the map
    \begin{align}
        \FRJ_a(\breve{F}_{\bar{v}})/\FRJ_a(\breve{\cO}_{\bar{v}})\longto
        \FRJ_{1,a}(\breve{F}_{\bar{v}})/\FRJ_{1,a}(\breve{\cO}_{\bar{v}})
    \end{align}
    is injective. Since \(C_a(\breve{F}_{\bar{v}})/C_a(\breve{\cO}_{\bar{v}})\) is a finitely
    generated free abelian group, the neutral components in \(\cP_{\bar{v}}(a)\)
    and \(\cP_{\bar{v},1}(a)\)  are homeomorphic. Thus,
    \(\pi_0(\cP_{\bar{v}}(a))\to \pi_0(\cP_{\bar{v},1}(a))\) is injective. This
    means that \(\pi_0(\cP_{\bar{v}}(a))\) can be identified with the image
    of \(\pi_0(\cP_{\bar{v}}^0(a))\) in \(\pi_0(\cP_{\bar{v},1}(a))\).

    Using the same argument as in the first claim, we have
    \begin{align}
        \pi_0(\cP_{\bar{v}}^0(a))&\simeq \CoCharG(\bT)_{\pi_a^\bullet(I_v)},\\
        \pi_0(\cP_{\bar{v},1}^0(a))&\simeq \CoCharG(\bT_1)_{\pi_a^\bullet(I_v)}.
    \end{align}
    Thus, \(\pi_0(\cP_{\bar{v}}(a))^\odot\) may be identified with a subgroup of the image of
    the homomorphism
    \begin{align}
        \dual{\bT}_1^{\pi_a^\bullet(I_v)}\longto
        \dual{\bT}^{\pi_a^\bullet(I_v)}.
    \end{align}
    Let \(\kappa_1\in\dual{\bT}_1^{\pi_a^\bullet(I_v)}\) with image
    \(\kappa\in\dual{\bT}\). Since \(\dual{\bG}_1\) has a simply-connected
    derived subgroup, the centralizer \(\dual{\bH}_1\) of \(\kappa_1\) (a
    semisimple element) in \(\dual{\bG}_1\) is connected, and its image in
    \(\dual{\bG}\) is \(\dual{\bH}\). Thus, \(\pi_a^\bullet(I_v)\) is contained
    in \(\bW_{\bH}\). This shows that
    \begin{align}
        \pi_0(\cP_{\bar{v}})^\odot\subset\dual{\bT}(\pi_a^\bullet(I_v)).
    \end{align}

    To show the opposite inclusion, we only need to show that we in fact have
    \begin{align}
        \pi_0(\cP_{\bar{v},1}^0(a))=\pi_0(\cP_{\bar{v},1}(a)),
    \end{align}
    or in other words, any \(\kappa_1\in\dual{\bT}_1^{\pi_a^\bullet(I_v)}\),
    viewed as a \(\Qlb\)-valued character of
    \(\FRJ_{1,a}(\breve{F}_{\bar{v}})\), is trivial on
    \(\FRJ_{1,a}(\breve{\cO}_{\bar{v}})\). Indeed, note that we can always lift
    \(\kappa_1\) to an endoscopic datum \((\kappa_1,\OGT_{\kappa_1})\) of \(G_1\)
    over \(\cO_v\hat{\otimes}_kk'\) for some sufficiently large extension
    \(k'/k\), and we let \(\FRM_{H}\) be the base-change
    of the associated endoscopic monoid to \(\breve{\cO}_{\bar{v}}\). Then \(a\) must be
    the image of some \(a_{H}\in \FRC_{\FRM,H}(\breve{\cO}_{\bar{v}})\), and we
    have canonical injection
    \begin{align}
        \FRJ_{1,a}(\breve{\cO}_{\bar{v}})\subset\FRJ_{H_1,a_H}(\breve{\cO}_{\bar{v}}).
    \end{align}
    Replacing \(G_1\) by endoscopic group \(H_1\), it suffices to prove that if
    \(\kappa_1\) lies in the center of 
    \(\dual{\bG}_1\), then its restriction to \(\FRJ_{1,a}(\breve{\cO}_{\bar{v}})\) is
    trivial.

    Recall that we have surjective map (note that we no longer assume that \(Z_{G_1}\) is
    connected)
    \begin{align}
        Z_{G_1,0}(\breve{\cO}_{\bar{v}})\x\FRJ_{a}^\SC(\breve{\cO}_{\bar{v}})\longto\FRJ_{1,a}(\breve{\cO}_{\bar{v}}),
    \end{align}
    because the residue field is algebraically closed.
    The restriction of \(\kappa_1\) to \(\FRJ_{a}^\SC(\breve{\cO}_{\bar{v}})\)
    is trivial because the image of \(\kappa_1\) in
    \(\dual{\bG}^\AD\) (the dual group of \(G^\SC\)) is the identity, and the
    restriction of \(\kappa_1\) to \(Z_{G_1,0}(\breve{\cO}_{\bar{v}})\) is also
    trivial, because the \emph{connected} center \(Z_{G_1,0}\) always lies in
    \(\FRJ_{a,1}^0\). This finishes the proof.
\end{proof}

\section{Approximations}
\label{sec:approximation_of_multiplicative_affine_springer_fibers}

Later on we will need to connect local objects with global ones. Therefore, it would be
useful to have some approximation results for multiplicative affine
Springer fibers. More precisely, we will prove the following proposition:

\begin{proposition}
    \label[proposition]{prop:approximation_GASF_local}
    For a fixed \(a\in\FRC_\FRM(\cO_v)\cap\FRC_\FRM^{\x,\rss}(F_v)\), let
    \(\cM_v(a)\) be defined using an arbitrarily fixed \(\gamma\in\FRM(F_v)\)
    lying over \(a\). Then there
    exists some integer \(N\) such that for any \(a'\in\FRC_\FRM(\cO_v)\) with
    \(a\equiv a'\bmod \pi_v^N\), there exists a choice of
    \(\gamma_0'\in\FRM(F_v)\) lying over \(a'\) with which we may define
    \(\cM_v(a')\) so that we have \(k\)-isomorphisms
    \(\cM_v(a)\cong\cM_v(a')\) and \(\cP_v(a)\cong\cP_v(a')\) compatible with the action of
    \(\cP_v(a)\) (resp.~\(\cP_v(a')\)) on \(\cM_v(a)\) (resp.~\(\cM_v(a')\)).
    Moreover, we also have isomorphism of groupoids
    \begin{align}
        \Stack*{\cM_v(a)/\cP_v(a)}(k)\cong
        \Stack*{\cM_v(a')/\cP_v(a')}(k)
    \end{align}
    regardless of which \(\gamma\) and \(\gamma_0'\) we choose.
\end{proposition}
The proof is a slight enhancement of
\cite{Ch22}*{Theorem~5.1.1}, which is in turn an adaptation of
\cite{Ng10}*{Proposition~3.5.1}. First we need a few lemmas.

\begin{lemma}
    \label[lemma]{lem:MASF_determined_by_reg_cent}
    Let \(g\in G(F_v)\) and \(\gamma_0\in\FRM(F_v)\) such that
    \(a=\chi_\FRM(\gamma_0)\in\FRC_\FRM^{\x,\rss}(F_v)\cap\FRC_\FRM(\cO_v)\).
    Let \(I_{\gamma_0}^\FRM\) be the centralizer of \(\gamma_0\) in \(\FRM^\x\).
    Then we have \(\Ad_g^{-1}(\gamma_0)\in\FRM(\cO_v)\) if and only if
    \begin{align}
        \Ad_g^{-1}(\gamma_0\FRJ_a^\FRM(\cO_v))\subset \FRM(\cO_v),
    \end{align}
    where \(\FRJ_a^\FRM(\cO_v)\) is treated as a subset of
    \(I_{\gamma_0}^\FRM(F_v)\) via 
    the canonical isomorphism \(\FRJ_a^\FRM(F_v)\simeq I_{\gamma_0}^\FRM(F_v)\).
\end{lemma}
\begin{proof}
    The ``if'' direction is trivial because
    \(\gamma_0\in\gamma_0\FRJ_a^\FRM(\cO_v)\) and we now prove the ``only
    if'' direction.
    Let \(\gamma=\Ad_g^{-1}(\gamma_0)\in\FRM(\cO_v)\). The centralizer
    \(I_{\gamma}^\FRM\) is a group scheme over \(\cO_v\) there is a canonical map
    \(\FRJ_a^\FRM\to I_\gamma^\FRM\). Over \(F_v\), the
    composition of canonical isomorphisms
    \begin{align}
        I_{\gamma_0,F_v}^\FRM\longto \FRJ_{a,F_v}^\FRM\longto
        I_{\gamma,F_v}^\FRM
    \end{align}
    is given by \(\Ad_g^{-1}\). Moreover, the
    image of \(\FRJ_a^\FRM(\cO_v)\) in \(I_\gamma^\FRM(F_v)\) is
    contained in \(I_\gamma^\FRM(\cO_v)\), so we have
    \begin{align}
        \Ad_g^{-1}(\FRJ_a^\FRM(\cO_v))\subset I_{\gamma}^\FRM(\cO_v)\subset
        \FRM^\x(\cO_v).
    \end{align}
    This implies that
    \begin{align}
        \Ad_g^{-1}(\gamma_0\FRJ_a^\FRM(\cO_v))\subset\gamma I_\gamma^\FRM(\cO_v)\subset\FRM(\cO_v),
    \end{align}
    as desired.
\end{proof}

\begin{remark}
    If \(\gamma_0\) happens to be a point in \(\FRM^\reg(\cO_v)\), then we may
    also replace \(\FRJ_a^\FRM(\cO_v)\) in the statement by
    \(I_{\gamma_0}^\FRM(\cO_v)\).
\end{remark}

\begin{lemma}
    \label[lemma]{lem:approximation_of_reg_cent}
    Let \(\gamma_0,\gamma_0'\in\FRM^\reg(\cO_v)\cap\FRM^{\x,\rss}(F_v)\) with
    \(\gamma_0\equiv\gamma_0'\bmod \pi_v\). Let
    \(a=\chi_\FRM(\gamma_0)\) and \(a'=\chi_\FRM(\gamma_0')\). Suppose there
    exists a \(W\)-equivariant isomorphism between cameral
    covers \(\iota\colon\tilde{X}_a\cong\tilde{X}_{a'}\) lifting
    the identity modulo \(\pi_v\). Then there exists some \(g\in G^\SC(\cO_v)\)
    such that
    \begin{align}
        \Ad_g^{-1}(I_{\gamma_0}^\FRM)=I_{\gamma_0'}^\FRM.
    \end{align}
\end{lemma}
\begin{proof}
    Let \(R_a\) (resp.~\(R_{a'}\)) be the ring of functions on \(\tilde{X}_a\)
    (resp.~\(\tilde{X}_{a'}\)) and \(F_a=R_a\otimes_{\cO_v}F_v\)
    (resp.~\(F_{a'}=R_{a'}\otimes_{\cO_v}F_v\)). Let \(R_a^\flat\)
    (resp.~\(R_{a'}^\flat\)) be the normalization of \(R_a\) (resp.~\(R_{a'}\))
    in \(F_a\) (resp.~\(F_{a'}\)).

    Since \(\gamma_0\) (resp.~\(\gamma_0'\)) is regular at special point \(v\),
    we have canonical isomorphism between \(\FRJ_a^\FRM\)
    (resp.~\(\FRJ_{a'}^\FRM\)) and \(I_{\gamma_0}^\FRM\) (resp.~\(I_{\gamma_0'}^\FRM\)).
    Since the regular centralizer is canonically determined by the cameral cover
    (using the Galois description \Cref{prop:Galois_reg_cent_monoid}), so is
    \(I_{\gamma_0}^\FRM\) (resp.~\(I_{\gamma_0'}^\FRM\)). Therefore, the
    \(W\)-equivariant isomorphism
    \(\tilde{X}_a\cong\tilde{X}_{a'}\) induces isomorphism of \(\cO_v\)-group
    schemes
    \begin{align}
        \iota^\FRM\colon I_{\gamma_0}^\FRM\stackrel{\sim}{\longto}
        I_{\gamma_0'}^\FRM.
    \end{align}

    We may describe \(\iota^\FRM\) more explicitly as follows: the canonical inclusion
    \begin{align}
        \FRJ_a^\FRM(\cO_v)\subset T_\FRM(R_a)^W
    \end{align}
    extends to an isomorphism
    \begin{align}
        I_{\gamma_0}^\FRM(F_v)\stackrel{\sim}{\longto}T_\FRM(F_a)^W,
    \end{align}
    which is given by \(\Ad_h^{-1}\) for some \(h\in \FRM^\x(F_a)\). The image of
    \(h\) in flag variety \(\FRM^\x/B_{\FRM^\x}\) extends to an
    \(R_a^\flat\)-point by properness and is canonically defined by the proof of
    \Cref{prop:Galois_reg_cent_monoid}. Similarly, we have \(h'\) for
    \(\gamma_0'\). Since
    \(\gamma_0\equiv\gamma_0'\bmod\pi_v\) and \(\iota\equiv \Id\bmod \pi_v\), we
    may choose \(h\) and \(h'\) such that \(h\iota^*(h')^{-1}\in
    \FRM^\x(R_a^\flat)\) and \(h\equiv \iota^*(h')\bmod\pi_v\).

    Since \(\iota^\FRM\) is equal to \(\Ad_{h\iota^*(h')^{-1}}^{-1}\) on
    \(F_v\)-points, we have
    \begin{align}
        \iota^\FRM(\gamma_0)\in I_{\gamma_0'}^\FRM(F_v)\cap
        \FRM(R_a^\flat)=I_{\gamma_0'}^\FRM(F_v)\cap \FRM(\cO_v),
    \end{align}
    and \(\chi_\FRM(\iota^\FRM(\gamma_0))=a\). We also have that
    \begin{align}
        \iota^\FRM(\gamma_0)\equiv\gamma_0\equiv\gamma_0'\bmod \pi_v,
    \end{align}
    meaning \(\iota^\FRM(\gamma_0)\) is regular at special point \(v\).
    As a consequence, we have
    \(I_{\iota^\FRM(\gamma_0)}^\FRM=I_{\gamma_0'}^\FRM\) as subgroups of
    \(\FRM^\x\) over \(\cO_v\). Finally, since the map
    \begin{align}
        G^\SC\x_X \FRM^\reg\longto \FRM^\reg\x_{\FRC_\FRM}\FRM^\reg
    \end{align}
    is smooth, we can find \(g\in G^\SC(\cO_v)\) with \(g\equiv 1\bmod \pi_v\)
    such that \(\Ad_g^{-1}(\gamma_0)=\iota^\FRM(\gamma_0)\), and so
    \begin{align}
        \Ad_g^{-1}(I_{\gamma_0}^\FRM)=I_{\iota^\FRM(\gamma_0)}^\FRM=I_{\gamma_0'}^\FRM.
    \end{align}
    This finishes the proof.
\end{proof}

\begin{proof}
    [Proof of \Cref{prop:approximation_GASF_local}]
    Without loss of generality, we may assume that \(k_v=k\). Similar to
    \cite{Ng10}*{Lemme~3.5.2}, we can find some \(N>0\) such that for any
    \(a'\in\FRC_\FRM(\cO_v)\) with \(a\equiv a'\bmod \pi_v^N\) there exists a
    \(W\)-equivariant isomorphism \(\tilde{X}_a\cong \tilde{X}_{a'}\) over
    \(\cO_v\) lifting the identity modulo \(\pi_v^N\). Let \(\gamma_0\) and
    \(\gamma_0'\) be the respective (arbitrarily chosen) points over \(a\) and
    \(a'\) defining \(\cM_v(a)\) and \(\cM_v(a')\). Let \(\gamma_1\)
    (resp.~\(\gamma_1'\)) be the point in \(\FRM^\reg(\breve{\cO}_{\bar{v}})\)
    induced by a fixed Steinberg quasi-section over \(\bar{k}\). By
    \Cref{lem:approximation_of_reg_cent}, we may find \(g_1\in
    G^\SC(\breve{\cO}_{\bar{v}})\) such that
    \begin{align}
        \Ad_{g_1}^{-1}(I_{\gamma_1}^\FRM)=I_{\gamma_1'}^\FRM.
    \end{align}
    Let \(h,h'\in G^\SC(\breve{F}_{\bar{v}})\) be such that
    \(\Ad_h(\gamma_0)=\gamma_1\) and
    \(\Ad_{h'}(\gamma_0')=\gamma_1'\). Then we have that
    \begin{align}
        \Ad_g^{-1}(I_{\gamma_0}^\FRM)=I_{\gamma_0'}^\FRM,
    \end{align}
    where \(g=h^{-1}g_1 h'\).
    By \Cref{lem:approximation_in_arc_monoid}, a technical result that we will
    prove later, and enlarging \(N\) if necessary, we can ensure that
    \begin{align}
        (\gamma_1')^{-1}\iota^\FRM(\gamma_1)\in
        G^\SC(\breve{\cO}_{\bar{v}})\cap
        I_{\gamma_1'}^\FRM(\breve{F}_{\bar{v}})\subset
        I_{\gamma_1'}^\FRM(\breve{\cO}_{\bar{v}}),
    \end{align}
    so that we also have
    \begin{align}
        \Ad_g^{-1}\bigl(\gamma_0\FRJ_{a}^\FRM(\breve{\cO}_{\bar{v}})\bigr)
        &=\Ad_{h'}^{-1}\bigl(\Ad_{g_1}^{-1}(\gamma_1)\bigr)\FRJ_{a'}^\FRM(\breve{\cO}_{\bar{v}})\\
        &=\Ad_{h'}^{-1}\bigl(\iota^\FRM(\gamma_1)\bigr)\FRJ_{a'}^\FRM(\breve{\cO}_{\bar{v}})\\
        &\subset\Ad_{h'}^{-1}\bigl(\gamma_1'I_{\gamma_1'}^\FRM(\breve{\cO}_{\bar{v}})\bigr)\FRJ_{a'}^\FRM(\breve{\cO}_{\bar{v}})\\
        &=\gamma_0'\FRJ_{a'}^\FRM(\breve{\cO}_{\bar{v}}).
    \end{align}
    If the Steinberg quasi-section exists over \(\cO_v\), then we may
    let \(\gamma_0=\gamma_1\), \(\gamma_0'=\gamma_1'\) and \(g\in
    G^\SC(\cO_v)\). In this case the proposition is an immediate consequence of
    \Cref{lem:MASF_determined_by_reg_cent}: the isomorphism between
    \(\cM_{\bar{v}}(a)\) and \(\cM_{\bar{v}}(a')\) is given by
    \begin{align}
        \cM_{\bar{v}}(a)&\stackrel{\sim}{\longto} \cM_{\bar{v}}(a')\\
        m&\longmapsto g^{-1}m,
    \end{align}
    and that between \(\cP_{\bar{v}}(a)\) and \(\cP_{\bar{v}}(a')\) is just
    \(\Ad_g^{-1}\).

    Even if the Steinberg quasi-section may not exist over \(\cO_v\), we still
    have that
    \(\Ad_{g}^{-1}\) maps \(\gamma_0\FRJ_a^\FRM(\breve{\cO}_{\bar{v}})\) to
    \(\gamma_0'\FRJ_{a'}^\FRM(\breve{\cO}_{\bar{v}})\) and
    \(\FRJ_a^\FRM(\cO_v)\subset I_{\gamma_0}^\FRM(F_v)\) to \(\FRJ_{a'}^\FRM(\cO_v)\subset
    I_{\gamma_0'}^\FRM(F_v)\). Moreover, since \(\Ad_g^{-1}\) equals \(\iota^\FRM\)
    as isomorphisms between \(F_v\)-group schemes \(\FRJ_a^\FRM\) and
    \(\FRJ_{a'}^\FRM\), it maps \(F_v\)-points to \(F_v\)-points, and so
    we have
    \begin{align}
        \Ad_g^{-1}(\gamma_0)\in
        \gamma_0'\FRJ_{a'}^\FRM(\breve{\cO}_{\bar{v}})\cap
        I_{\gamma_0'}^\FRM(F_v)=\gamma_0'\FRJ_{a'}^\FRM(\cO_v).
    \end{align}
    This implies that \(g\Frob_v(g)^{-1}\in I_{\gamma_0}^\FRM(F_{\bar{v}})\), which
    in turn means that we have isomorphism of \(k\)-points
    \begin{align}
        \Stack*{\cM_v(a)/\cP_v(a)}(k)\longto\Stack*{\cM_v(a')/\cP_v(a')}(k).
    \end{align}
    Finally, replacing \(\gamma_0\) by \(\Ad_g^{-1}(\gamma_0)\) and
    swapping roles of \(a\) (resp.~\(\gamma_0\), etc.) with \(a'\) (resp.~\(\gamma_0'\),
    etc.), we see that \(\gamma_0'\) can be so chosen that we have compatible
    \(k\)-isomorphisms (the first one is in fact even an equality)
    \begin{align}
        \cM_v(a)&\cong\cM_v(a'),\\
        \cP_v(a)&\cong\cP_v(a').
    \end{align}
    This finishes the proof of \Cref{prop:approximation_GASF_local}.
\end{proof}


\section{Connection with MV-cycles and Transfer Factors}%
\label{sec:connection_with_mv_cycles}

In this section we study the special case of unramified conjugacy classes. In
this case, \cite{Ch22} established a bijection between the irreducible
components of \(\cM_{G,\bar{v}}(a)\) modulo \(\cP_{\bar{v}}(a)\)-action
and that of certain Mirkovi\'c--Vilonen (MV) cycles. This was done at
geometric level, and here we expand their result to include Frobenius actions.
Seemingly a straightforward task, it turns out to be quite deep
because what we will observe is a concrete example of local transfer factor at
work.

\subsection{}
First, let us consider the case where \(G=\bG\) is split, and \(\gamma\in
\bT(\breve{F}_{\bar{v}})^\rss\). Let \(\lambda\) be a dominant cocharacter of
\(\bT\), and suppose \(\gamma_\lambda\in \bM\) as in
\Cref{lem:find_gamma_lambda_for_gamma} exists for some \(\bM\in\FM(\bG^\SC)\).
Since \(\gamma\in\bT\), we have \(\gamma_\lambda\in\bar{\bT}_\bM\). Let
\(a=\chi_\FRM(\gamma_\lambda)\) and suppose it lies in
\(\bC_\bM(\breve{\cO}_{\bar{v}})\). Then the \(\breve{F}_{\bar{v}}\)-torus
\(\bJ_{a,\breve{F}_{\bar{v}}}\) is canonically isomorphic to the maximal torus
\(\bT_{\breve{F}_{\bar{v}}}\) itself. Thus, its N\'eron model is just
\(\bT_{\breve{\cO}_{\bar{v}}}\). Recall we have the lattice \(\Lambda_a\) being
the largest free quotient of \(\pi_0(\cP_{\bar{v}}(a))^\Red\). In this case
\(\Lambda_a\) is isomorphic to \(\CoCharG(\bT)\), which, after choosing a
uniformizer of \(\breve{F}_{\bar{v}}\), can be regarded as a subgroup of
\(\cP_{\bar{v}}(a)\).

\subsection{}
Following \cite{Ch22}*{\S~3.5.1}, let
\begin{align}
    \KVU_{\bG,\bar{v}}^{\le\lambda}(\gamma)=\Set*{u\in
        \Loop_{\bar{v}}{\bU}/\Arc_{\bar{v}}{\bU}\given
        \Ad_u^{-1}(\gamma) \in
    \Cartan_{\bG}^{\le\lambda}},
    \nomenclature[\(Y"sf_G_v_lambda_gamma \)]{\(\KVU_{G,\bar{v}}^\lambda(\gamma),\KVU_{G,\bar{v}}^{\le\lambda}(\gamma)\)}{the
    analogue of \(\cM_{G,\bar{v}}^{\lambda}(\gamma)\) or \(\cM_{G,\bar{v}}^{\le\lambda}(\gamma)\) for an unramified \(\gamma\) by replacing \(G\) with \(U\)}
\end{align}
and \(\tilde{\KVU}_{\bG,\bar{v}}^{\le\lambda}(\gamma)\) be its preimage in
\(\Loop_{\bar{v}}{\bU}\), both with reduced ind-scheme structure.
Let \(S_{\mu,\bar{v}}\) be the semi-infinite orbit
\begin{align}
    S_{\mu,\bar{v}}=\Loop_{\bar{v}}{U}\pi^{\mu}\Arc_{\bar{v}}{G}/\Arc_{\bar{v}}{G}\subset
    \Gr_{G,\bar{v}}
\end{align}
for any \(\mu\in\CoCharG(T)\), and let
\(S_{G,\mu}^{\lambda}=S_{\mu}^{\lambda}\)
    \nomenclature[\(S_G_mu_lambda \)]{\(S_{G,\mu}^\lambda,S_{G,\mu}^{\le\lambda}\)}{the big-cell locus of \(S_{G,\mu}^{\le\lambda}\)}
be the
intersection \(S_{\mu,\bar{v}}\cap\Cartan_{\bG,\bar{v}}^{\lambda}\). We similarly have
\(S_{\mu}^{\le\lambda}\).
    \nomenclature[\(S_G_mu_lambda_le \)]{\(S_{G,\mu}^{\le\lambda}\)}{the union of MV-cycles of type \(\lambda\) and weight \(\mu\)}

We also denote by \(\Hk_{\bG,\bar{v}}^{\le\lambda}\)
    \nomenclature[\(H{}k_G_mu_lambda_le \)]{\(\Hk_{G,v}^{\le\lambda}\)}{the finite-type local Hecke stack at place \(v\) of \(G\) with type \(\lambda\)}
the local Hecke stack truncated at \(\lambda\), in other words, formally it is
the prestack defined by quotient
\begin{align}
    \Hk_{\bG,\bar{v}}^{\le\lambda}\defeq\Stack*{\Arc_{\bar{v}}{\bG}\backslash\Gr_{\bG,\bar{v}}^{\le\lambda}}.
\end{align}
We will skip a more detailed definition because we will only use its
\(\bar{k}\)-points in this section, but see \Cref{def:equivariant_global_affine_Schubert}.

We have commutative diagram
\begin{equation}
    \label{eqn:diagram_connecting_GASF_to_MV}
    \begin{tikzcd}
        \tilde{\KVU}_{\bG,\bar{v}}^{\le\lambda}(\gamma)\ar[r, "{u\mapsto
        \Ad_u^{-1}(\gamma)}"]\ar[d, "{/\Arc_{\bar{v}}{\bU}}", swap] &
        \Cartan_{\bG,\bar{v}}^{\le\lambda}\cap
        \Loop_{\bar{v}}{\bU}\gamma\ar[d, "{/\Arc_{\bar{v}}{\bU}}"]\\ 
        \KVU_{\bG,\bar{v}}^{\le\lambda}(\gamma)\ar[d]\ar[rd] &
        S_{\nu_\gamma}^{\le\lambda}\ar[d]\\
        \CoCharG(\bT)\x\KVU_{\bG,\bar{v}}^{\le\lambda}(\gamma)\ar[d,
        "{(\mu,u)\mapsto \pi^\mu u}", swap] &
        \Hk_{\bG,\bar{v}}^{\le\lambda}\\
        \SP_{\bG,\bar{v}}^{\le\lambda}(\gamma)\ar[ru, "{g\mapsto \Ad_g^{-1}(\gamma)}", swap]
    \end{tikzcd}.
\end{equation}
In this diagram, the arrows marked with ``\(/\Arc_{\bar{v}}{\bU}\)'' are
\(\Arc_{\bar{v}}{\bU}\)-torsors, and the top horizontal
arrow, after taking
quotient by a sufficiently small congruent subgroup of \(\Arc_{\bar{v}}{\bU}\),
is a ``homotopy equivalence'' in the following sense:
\begin{definition}
    Let \(Y_1\), \(Y_2\) be \(k\)-schemes of finite type. We say \(Y_1\) and
    \(Y_2\) are \inotion{\(\bbA\)-equivalent} if they can be connected
    by a chain of roofs between \(k\)-schemes
    \begin{align}
        Y\longot \tilde{Y}\longto Y',
    \end{align}
    where the arrows are locally trivial fibration of affine spaces. In this
    case, we denote \(Y_1\leftrightsquigarrow Y_2\).
\end{definition}

We continue describing \eqref{eqn:diagram_connecting_GASF_to_MV}.
The lower-left vertical arrow is a bijection on \(\bar{k}\)-points and induces a
stratification compatible with the ind-scheme topology on
\(\SP_{\bG,\bar{v}}^\lambda(\gamma)\) induced by that on \(\Gr_{\bG,\bar{v}}\). 
So the closure of \(\KVU_{\bG,\bar{v}}^\lambda(\gamma)\) in
\(\SP_{\bG,\bar{v}}^\lambda(\gamma)\) is a fundamental domain of the free
\(\CoCharG(\bT)\)-action, and the irreducible components of
\(\Stack{\cM_{\bG,\bar{v}}(a)/\cP_{\bar{v}}(a)}\) are identified with
those of \(\KVU_{\bG,\bar{v}}^\lambda(\gamma)\).

\subsection{}
In order to study non-split elements, it is first necessary to have some
compatibility results in the
split case with respect to standard Levi subgroups \(\bL_\Rt\) generated by a simple root
\(\Rt\). Let \(\bP_\Rt\) be the standard parabolic containing \(\bB\) and
\(\bN_\Rt\) the unipotent radical of \(\bP_\Rt\). Let \(\bB_\Rt\) be the
standard Borel of \(\bL_\Rt\) containing \(\Rt\), and \(\bU_\Rt\) its unipotent
radical. Let \(s_\Rt\in\bW\) be the simple reflection determined by \(\Rt\). For
simplicity of notations, for any subgroup \(H\subset \bG\), we let
\begin{align}
    \SFM_{H}
    &=\Set*{h\in \Gr_{H,\bar{v}}\given
    \Ad_h^{-1}(\gamma) \in \Cartan_{\bG,\bar{v}}^{\le\lambda}},\\
    \tilde{\SFM}_{H}
    &=\Set*{h\in \Loop_{\bar{v}}{H}\given
    \Ad_h^{-1}(\gamma) \in \Cartan_{\bG,\bar{v}}^{\le\lambda}}.
\end{align}
Note that when \(H=\bL_\Rt\), \(\SFM_H\) is the union of
\(\cM_{\bL_\Rt,\bar{v}}^{\le_{\bL_\Rt}\mu}\) with \(\mu\le_\bG \lambda\).

\newsavebox{\KVDIMa}
\newsavebox{\KVDIMb}
\newsavebox{\KVDIMc}
\newsavebox{\KVDIMd}
\sbox{\KVDIMa}{%
    \fbox{
    \begin{tikzcd}[ampersand replacement = \&, column sep=small]
        \&\& \tilde{\SFM}_{\bB} \ar[ld]\ar[rd] \&\& \\
        \& \tilde{\SFM}_{\bB_\Rt} \ar[ld]\ar[rd]\&  \& \tilde{\SFM}_{\bP_\Rt} \ar[ld]\ar[rd]\& \\
        \tilde{\SFM}_{\bT} \&\& \tilde{\SFM}_{\bL_\Rt}  \&\& \tilde{\SFM}_{\bG}
    \end{tikzcd}%
    }
}
\sbox{\KVDIMb}{%
    \fbox{
    \begin{tikzcd}[ampersand replacement = \&, column sep=small]
        \&\& \Gr_{\bB}^{\le_\bG\lambda} \ar[ld]\ar[rd] \&\& \\
        \& \Gr_{\bB_\Rt}^{\le_\bG\lambda} \ar[ld]\ar[rd]\&  \& \Gr_{\bP_\Rt}^{\le_\bG\lambda} \ar[ld]\ar[rd]\& \\
        \Gr_{\bT}^{\le_\bG\lambda} \&\&  \Gr_{\bL_\Rt}^{\le_\bG\lambda} \&\& \Gr_{\bG}^{\le_\bG\lambda}
    \end{tikzcd}%
    }
}
\sbox{\KVDIMc}{%
    \fbox{
    \begin{tikzcd}[ampersand replacement = \&, column sep=small]
        \&\& \SFM_{\bB} \ar[ld]\ar[rd] \&\& \\
        \& \SFM_{\bB_\Rt} \ar[ld]\ar[rd]\&  \& \SFM_{\bP_\Rt} \ar[ld]\ar[rd]\& \\
        \SFM_{\bT} \&\& \SFM_{\bL_\Rt}  \&\& \SFM_{\bG}
    \end{tikzcd}%
    }
}
\sbox{\KVDIMd}{%
    \fbox{
    \begin{tikzcd}[ampersand replacement = \&, column sep=small]
        \&\& \Hk_{\bB}^{\le_\bG\lambda} \ar[ld]\ar[rd] \&\& \\
        \& \Hk_{\bB_\Rt}^{\le_\bG\lambda} \ar[ld]\ar[rd]\&  \& \Hk_{\bP_\Rt}^{\le_\bG\lambda} \ar[ld]\ar[rd]\& \\
        \Hk_{\bT}^{\le_\bG\lambda} \&\&  \Hk_{\bL_\Rt}^{\le_\bG\lambda} \&\& \Hk_{\bG}^{\le_\bG\lambda}
    \end{tikzcd}%
    }
}
\begin{figure*}[t]
\begin{equation}
    \label{eqn:Levi_induction_of_MV_big_diagram}
    \begin{tikzcd}[column sep=large]
        \usebox{\KVDIMa} \ar[r, "\widetilde\ACT", "g\mapsto
        \Ad_g^{-1}(\gamma)"']\ar[d, "\tilde{p}"'] & \usebox{\KVDIMb} \ar[d, "p"]\\
        \usebox{\KVDIMc} \ar[r, "\ACT", "g\mapsto \Ad_g^{-1}(\gamma)"'] & \usebox{\KVDIMd}
    \end{tikzcd}
\end{equation}
\end{figure*}
We have diagram \eqref{eqn:Levi_induction_of_MV_big_diagram} where the diamond
inside each framed subdiagram is Cartesian.
To see that the arrows from \(\bP_\Rt\)-related objects to \(\bL_\Rt\)-related
ones are well-defined: suppose \(p=mn\in\Loop_{\bar{v}}\bP_\Rt\) where
\(m\in\Loop_{\bar{v}}\bL_\Rt\) and \(n\in\Loop_{\bar{v}}\bN_\Rt\) is such
that \(n^{-1}m^{-1}\gamma mn\in \Cartan_{\bG,\bar{v}}^{\le\lambda}\). Suppose
\(m^{-1}\gamma m\in\Cartan_{\bL_\Rt,\bar{v}}^\mu\), then we can find \(m_0,m_0'\in
\Arc_{\bar{v}}\bL_\Rt\), such that \(m^{-1}\gamma m=m_0^{-1}\pi^\mu m_0'\). This
implies that
\begin{align}
    m_0n^{-1}m^{-1}\gamma mnm_0^{\prime-1}
    &= m_0n^{-1}m_0^{-1}\pi^\mu m_0'nm_0^{\prime-1}\\
    &=(m_0nm_0^{-1})^{-1}\big(\pi^\mu (m_0'nm_0^{\prime-1})\pi^{-\mu}\big)\pi^\mu\\
    &\in \Loop_{\bar{v}}\bU\pi^\mu\Arc_{\bar{v}}\bG\cap
    \Cartan_{\bG,\bar{v}}^{\le\lambda}.
\end{align}
This implies that \(\mu\le_\bG\lambda\), hence the maps
\(\tilde{\SFM}_{\bP_\Rt}\to\tilde{\SFM}_{\bL_\Rt}\) and
\(\SFM_{\bP_\Rt}\to\SFM_{\bL_\Rt}\) are well-defined. The well-definedness of
the other two maps is easy.
That the four diamonds are Cartesian is easy to check using the fact that
\(\bP_\Rt/\bB=\bL_\Rt/\bB_\Rt\).

\subsection{}
We make another convention of notations: for any \(\ell\)-adic complex \(\cF\)
on any node inside the lower-left or the upper-right diagram, we use \(\cF^*\)
to denote the pullback along a lower-right-pointing arrow, and \(\cF_!\) the
shriek pushforward along a lower-left-pointing arrow, if such arrows exist.
Note that \((\cF^*)_!\simeq (\cF_!)^*\) by proper base change, so \(\cF_!^*\)
is unambiguously defined, whenever the relevant arrows exist.
For example, if \(\cF=\IC^\lambda\) on \(\Gr_\bG^{\le\lambda}\), then
\(\cF_{!\,!}^{**}\) is the hyperbolic localization functor corresponding to
restriction functor \(\Res^{\dual{\bG}}_{\dual{\bT}}\) on the representation
side, and \(\cF_!^*\) corresponds to \(\Res^{\dual{\bG}}_{\dual{\bL}_{\Rt}}\),
and so on.

\subsection{}
Let \(\bB_{\Rt}'\subset\bL_\Rt\) be the Borel subgroup containing \(\bT\)
opposite to \(\bB_\Rt\), and let \(\bB'\subset \bP_\Rt\subset \bG\) be the
corresponding Borel subgroup in \(\bG\). Then we have schemes \(\SFM_{\bB'}\)
and \(\SFM_{\bB_{\Rt}'}\). Suppose \(Z\subset \SFM_\bB\) and
\(Z'\subset\SFM_{\bB'}\) are two respective irreducible components such that
their images in \(\SFM_\bG=\cM_{\bG,\bar{v}}^{\le\lambda}(\gamma)\) coincide.
Then they have the same image in \(\SFM_{\bP_\Rt}\) as well, hence also in
\(\SFM_{\bL_\Rt}\).

Suppose further that \(S=\widetilde{\ACT}\circ\tilde{p}^{-1}(Z)\) is the
corresponding MV-cycle in \(S_{\nu_\gamma}^{\le\lambda}\), and similarly for
\(Z'\) we have \(S'\), an
MV-cycle formed using \(\bU'\) instead of \(\bU\). The cohomology class
\([S]\in\RH^\bullet(\Gr_\bG,\IC^\lambda)\) may be regarded as a vector in weight
space \(V_\lambda[\nu_\gamma]\). Similarly, \([S']\) is regarded as a vector in
\(V_{s_\Rt(\lambda)}'[\nu_\gamma]\), where \(V_{s_\Rt(\lambda)}'\) means the
\(\dual{\bG}\)-representation with \(\bU'\)-highest weight \(s_\Rt(\lambda)\).
Through the canonical isomorphisms
\begin{align}
    V_\lambda\simeq \RHc^\bullet(\Gr_\bG,\IC^\lambda)\simeq V_{s_\Rt(\lambda)}',
\end{align}
we may canonically identify \(V_{s_\Rt(\lambda)}'[\nu_\gamma]\) with
\(V_{\lambda}[s_\Rt( \nu_\gamma )]\). Clearly, if \((Z_1,Z_1')\) is another pair
like \((Z,Z')\) such that \(\widetilde{\ACT}\circ\tilde{p}^{-1}(Z_1)=S\), then
we also have \(\widetilde{\ACT}\circ\tilde{p}^{-1}(Z_1')=S'\). Thus, the
assignment \([S]\mapsto [S']\) induces a well-defined isomorphism
\begin{align}
    \bFs_\Rt^\lambda\colon V_\lambda[\nu_\gamma]\longto V_\lambda[s_\Rt(\nu_\gamma)].
\end{align}
The following lemma is crucial for relating transfer factor to
representations of the dual group:
\begin{lemma}
    \label[lemma]{lem:MV_transfer_factor_key}
    Let \(\dual{s}_\Rt\in\dual{\bG}\) be any lifting of \(s_\Rt\). Then for any
    \([S]\) as above, there exists some \(c_S\in\Qlb^\x\) such that
    \(\bFs_\Rt^\lambda([S])=c_S\dual{s}_\Rt([S])\).
\end{lemma}

\subsection{}
In order to prove \Cref{lem:MV_transfer_factor_key}, we need to reduce it to an
analogous statement for the Levi subgroup \(\bL_\Rt\), for which direct
computation becomes feasible.
By functoriality of geometric Satake, the image of \([S]\) in
\(\RHc^\bullet(\Gr_{\bL_\Rt},(\IC^\lambda)_!^*)\) may be regarded as a vector in
\(\Res_{\dual{\bL}_\Rt}^{\dual{\bG}}V_\lambda[\nu_\gamma]\).
Decompose \(V_\lambda\) into \(\dual{\bL}_\Rt\)-irreducible representations:
\begin{align}
    \Res_{\dual{\bL}_\Rt}^{\dual{\bG}}V_\lambda
    \simeq\bigoplus_{\mu}V_{\Rt,\mu}\otimes\Hom_{\dual{\bL}_\Rt}(V_{\Rt,\mu},V_\lambda),
\end{align}
where \(\mu\) ranges over \(\bB_\Rt\)-dominant coweights of \(\bL_\Rt\).
For each \(\mu\), we also have a similarly defined map:
\begin{align}
    \bFs_{\bL_\Rt,\Rt}^\mu\colon V_{\Rt,\mu}[\nu_\gamma]\longto V_{\Rt,\mu}[s_\Rt(\nu_\gamma)].
\end{align}
In order to characterize the multiplicity space
\(\Hom_{\bL_\Rt}(V_{\Rt,\mu},V_\lambda)\), we examine the cohomology of
\((\IC^\lambda)_!^*\) more closely using more subtle properties of MV-cycles.

\subsection{}
Given irreducible MV-cycle \([S]\), one may view it as a subset of
\(\Gr_{\bP_\Rt}\), so the closure of its image in
\(\Gr_{\bL_\Rt}\) contains a unique MV-cycle \([S]_{\bL_\Rt}\) of \(\bL_\Rt\) as a
dense subset. The map \([S]\mapsto [S]_{\bL_\Rt}\) is not injective, so we will
need another ingredient to distinguish \([S]\) from all those cycles with image
\([S]_{\bL_\Rt}\). Suppose 
\([S]_{\bL_\Rt}\) contains a component of \(S_{\bL_\Rt,\nu_\gamma}^{\mu}\)
as a dense subset, then there is a unique irreducible component \([Q]\) of
\begin{align}
    S_{\bL_\Rt,\mu}^{\bG,\lambda}
    \defeq (\Loop_{\bar{v}}\bN_\Rt\Arc_{\bar{v}}\bL_\Rt)\pi^\mu\Arc_{\bar{v}}\bG/\Arc_{\bar{v}}\bG\cap
    \Gr_\bG^{\lambda}
\end{align}
such that the preimage of \([S]_{\bL_\Rt}\) in \([Q]\) contains a dense subset
of \(S\). The pair \(([Q],[S]_{\bL_\Rt})\) then uniquely determines \([S]\).
In other words, let \(\Irr[S_{\bL_\Rt,\mu}^{\bG,\lambda}]\)
    \nomenclature[\(Irr \)]{\(\Irr(\cdot)\)}{the set of irreducible components}
be the set of maximal dimensional
irreducible components of \(S_{\bL_\Rt,\mu}^{\bG,\lambda}\), then we have
canonical bijection of sets
\begin{align}
    \coprod_{\mu}\Irr[S_{\bL_\Rt,\mu}^{\bG,\lambda}]\x
    \Irr[S_{\bL_\Rt,\nu_\gamma}^{\mu}]\simeq
    \Irr[S_{\nu_\gamma}^{\lambda}].
\end{align}
These claims can be found in \cite{BG01} for example.
Therefore, for any fixed \(\mu\), the multiplicity space of \(V_{\Rt,\mu}\) in
\(V_\lambda\) may be identified with the vector space generated by
\(\Irr[S_{\bL_\Rt,\mu}^{\bG,\lambda}]\). In other words, we have
canonical isomorphism
\begin{align}
    V_\lambda\simeq\bigoplus_\mu
    \Qlb^{\oplus\Irr[S_{\bL_\Rt,\mu}^{\bG,\lambda}]}\otimes_{\Qlb}
    V_{\Rt,\mu},
\end{align}
and it is compatible with canonical isomorphisms of vector spaces
\begin{align}
    \RHc^\bullet(\Gr_\bG,\IC^\lambda)\simeq\RHc^\bullet(\Gr_{\bL_\Rt},(\IC^\lambda)_!^*)\simeq
    \RHc^\bullet(\Gr_\bT,(\IC^\lambda)_{!\,!}^{**}).
\end{align}

\begin{lemma}
    \label[lemma]{lem:MV_TF_reduce_to_rank_1_Levi}
    We have commutative diagram of isomorphic \(\Qlb\)-vector spaces
    \begin{equation}
        \begin{tikzcd}[column sep=huge]
            V_\lambda[\nu_\gamma] \ar[r, "\bFs_\Rt^\lambda"]\ar[d] & V_\lambda[s_\Rt( \nu_\gamma )] \ar[d]\\
            \bigoplus_\mu
            \Qlb^{\oplus\Irr[S_{\bL_\Rt,\mu}^{\bG,\lambda}]}\otimes_{\Qlb}V_{\Rt,\mu}[\nu_\gamma]
            \ar[r, "\bigoplus_\mu\Id\otimes\bFs_{\bL_\Rt,\Rt}^\mu"] & \bigoplus_\mu
            \Qlb^{\oplus\Irr[S_{\bL_\Rt,\mu}^{\bG,\lambda}]}\otimes_{\Qlb}V_{\Rt,\mu}[s_\Rt(\nu_\gamma)]
        \end{tikzcd}
    \end{equation}
    where the vertical maps are restrictions.
\end{lemma}
\begin{proof}
    Notice that \(S_{\bL_\Rt,\mu}^{\bG,\lambda}\) is insensitive to the choice
    between \(\bB_\Rt\) versus \(\bB_\Rt'\).
    Suppose \(Z\subset\SFM_\bB\) and \(Z'\subset\SFM_{\bB'}\) coincide in
    \(\SFM_\bG\), then they coincide in
    \(\SFM_{\bP_\Rt}\) hence so do their images in \(\SFM_{\bL_\Rt}\). This
    means that if the image of \([S]\) under the left vertical arrow is of the
    form \([Q]\otimes [S]_\mu\) where \([S]_\mu\in V_{\Rt,\mu}\), then that of
    \([S']\) under the right vertical arrow is of the form \([Q']\otimes
    [S']_\mu\) where the image of \([S]_\mu\) is \([S']_\mu\) under the map
    \(V_{\Rt,\mu}[\nu_\gamma]\to V_{\Rt,\mu}[s_\Rt(\nu_\gamma)]\). So it
    remains to show that \([Q]=[Q']\).

    The \([Q]\)-factor of \([S]\) may be
    computed from \(Z\) as follows: pick any sufficiently general point \(g\) in
    the image of \(Z\) in \(\Gr_{\bP_\Rt}\), then \(\Ad_g^{-1}(\gamma)\) is a
    well-defined \(\Arc_{\bar{v}}\bP_\Rt\)-orbit in \(\Gr_\bG^\lambda\), which
    is connected. It is thus contained in a unique component of
    \(S_{\bL_\Rt,\mu}^{\bG,\lambda}\), which is exactly \(Q\). Since \(g\)
    is also a sufficiently general point in the image of \(Z'\) (which coincides
    with that of \(Z\)), we have \([Q]=[Q']\) as desired.
\end{proof}

\begin{proof}[Proof of \Cref{lem:MV_transfer_factor_key}]
    Using \Cref{lem:MV_TF_reduce_to_rank_1_Levi}, we may assume \(\bG\) has
    semisimple rank \(1\) and \(\Rt\) is the unique positive root. Then the
    lemma is trivial because \(V_\lambda[\nu_\gamma]\) is one-dimensional.
\end{proof}

\begin{corollary}
    \label[corollary]{cor:MV_TF_independent_of_B_for_zero_coweight}
    Suppose \(\nu_\gamma\) is central. Let \(\bB'\) be any Borel subgroup
    containing \(\bT\). Let
    \(Z\) (resp.~\(Z'\)) be an irreducible component of
    \(\SFM_{\bB}\) (resp.~\(\SFM_{\bB'}\)), and let \([S]\)
    (resp.~\([S']\)) be the cohomology class of the corresponding MV-cycle in
    \(\RHc^0(\Gr_\bG,\IC^\lambda)\). If \(Z=Z'\) as
    subsets of \(\cM_{\bG,\bar{v}}^{\le\lambda}(\gamma)\), then \([S]=[S']\).
\end{corollary}
\begin{proof}
    Suppose \(\bB'=w(\bB)\) and let \(w=s_1\cdots s_m\) be a shortest (relative
    to \(\bB\)) expression. For each \(0\le i\le m\), let \(r_i=s_1\cdots s_i\)
    (so that \(r_0=1\)), and for \(i>0\) let
    \(\bB_i=r_i(\bB)=r_{i-1}s_ir_{i-1}^{-1}(\bB_{i-1})\). Then we have a
    sequence of Borel subgroups containing \(\bT\):
    \begin{align}
        \bB=\bB_0, \bB_1, \ldots, \bB_m=\bB',
    \end{align}
    such that \(\bB_i=s_i'(\bB_{i-1})\), and \(s_i'=r_{i-1}s_ir_{i-1}^{-1}\) is
    a simple reflection relative to \(\bB_{i-1}\). By the proof of
    \Cref{lem:MV_TF_reduce_to_rank_1_Levi}, each \(s_i'\) induces the identity
    map  on the subspace of \(\RHc^0(\Gr_\bG,\IC^\lambda)\) corresponding to
    coweight \(\nu_\gamma\), and so we are done.
\end{proof}

\begin{remark}
    We will prove a stronger version of
    \Cref{cor:MV_TF_independent_of_B_for_zero_coweight} in
    \Cref{prop:stable_twisted_Frobenius_is_the_crystal_action} by utilizing the
    concept of crystal bases. In fact, \Cref{lem:MV_transfer_factor_key} will be
    made more precise if we treat \([S]\) as a crystal basis element instead of
    a vector basis element: in that way there is a Weyl group action on
    those crystal basis elements and there is no need to lift Weyl group
    elements to \(\dual{\bG}\), and consequently the scalar of ambiguity \(c_S\)
    in the lemma can be eliminated.
\end{remark}

\subsection{}
Now we move to general \(G\) obtained by outer twist \(\OGT_G\) of \(\bG\). Let
\(\FRM\in\FM(G^\SC)\), and we want to consider
\(a\in\Stack*{\FRC_\FRM/Z_\FRM}(\cO_v)\) that is generically
regular semisimple and unramified. Base change to \(\bar{k}\), then as before we
have isomorphism
\begin{align}
    \cM_{G,v}(a)_{\bar{k}}\simeq \prod_{\bar{v}\colon k_v\to
    \bar{k}}\cM_{G,\bar{v}}(a),
\end{align}
where \(\bar{v}\) ranges over \(k\)-embeddings of \(k_v\) into \(\bar{k}\). In
particular, we have bijection of geometric irreducible components:
\begin{align}
    \Irr[\cM_{G,v}(a)_{\bar{k}}]\simeq \prod_{\bar{v}\colon k_v\to
    \bar{k}}\Irr[\cM_{G,\bar{v}}(a)].
\end{align}
The Frobenius \(\Frob_k\in\Gal(\bar{k}/k)\) acts on the \(\bar{k}\)-points on the
left-hand side, and it induces a \(\Frob_k\)-action on the right-hand side,
sending a \(\bar{k}\)-point in the \(\bar{v}\)-factor to one in the
\(\Frob_k(\bar{v})\)-factor. For convenience, let
\(\Frob_v=\Frob_k^{[k_v:k]}\) be the Frobenius of \(\bar{k}/k_v\), then
\(\Frob_v\) acts on each factor \(\cM_{G,\bar{v}}(a)\).

If we let \(\cM_{G,v}'(a)\) be the \(k_v\)-functor analogously defined as
\(\cM_{G,v}(a)\) except we treat \(\cO_v\)-points of \(\Stack{\FRM/G\x Z_\FRM}\)
with trivialization over \(F_v\) as \(k_v\)-points of \(\cM_{G,v}'(a)\) rather
than \(k\)-points, then \(\cM_{G,v}(a)\) is just the Weil restriction of
\(\cM_{G,v}'(a)\) from \(k_v\) to \(k\), and each \(\cM_{G,\bar{v}}(a)\) is the
base change of \(\cM_{G,v}'(a)\) via \(\bar{v}\colon k_v\to \bar{k}\).
Therefore, it is clear that the \(\Frob_k\)-action on \(\Irr[\cM_{G,v}(a)]\) is
induced by the \(\Frob_v\)-action on any one of
\(\cM_{G,\bar{v}}(a)\). Consequentially, without loss of generality we may
assume \(k_v=k\).

\subsection{}
Now we assume \(k_v=k\). Let
\(x_G\colon \breve{X}_v^\bullet\to \OGT_G\) be a fixed pinning. We may also
view it as a homomorphism \(\OGT_G^\bullet\colon\Gamma_v\to\Out(\bG)\) which factors
through \(\Gal(\bar{k}/k_v)\).
In later parts we will be using some notions in the
definition of transfer factor. The reader may refer to
\Cref{chap:Review_on_Transfer_Factors} for a full summary, as some portion of it
is formulated slightly differently than the original in \cite{LS87} to be more
coherent with the language in this book.

Recall that the definition of \(\cM_{G,v}(a)\) depends on a choice of
\(\gamma_a\in\Stack*{\FRM/Z_\FRM}(F_v)\) lying over \(a\), which necessarily
exists by \Cref{thm:A_2m_rationality_fix}. We view \(\gamma_a\) as an
element in \(G^\AD(F_v)\). In the situation relevant to fundamental lemma, it is
important to further assume that \(\gamma_a\) comes from a point \(\gamma\in
G(F_v)\), and the boundary divisor of \(a\) is \(-w_0(\lambda)\) where
\(\lambda\) is fixed by \(\Frob_v\). This assumption is not necessary for
a big portion of our discussion below, and so we will remind the reader when it
is needed later.

Since \(a\) is unramified, over
\(\breve{\cO}_{\bar{v}}\) it can be lifted to a point
\(x_a\in\FRT_\FRM(\breve{\cO}_{\bar{v}})\). The point \(x_a\) may also be
viewed as a point in \(\OGT_G\x \bar{\bT}_\bM\) lying over \(x_G\), or a lift of
\(\OGT_G^\bullet\) to a homomorphism \(\Gamma_v\to \bW\rtimes\Out(\bG)\) which
also factors through \(\Gal(\bar{k}/k_v)\).

Let \(h_{\bar{v}}\in
G^\SC(\breve{F}_{\bar{v}})\) such that \(h_{\bar{v}}^{-1}\gamma_a
h_{\bar{v}}=x_a\) (necessarily exists by Steinberg's theorem
on torsors over \(\breve{F}_{\bar{v}}\)). We will abuse notation and also use
\(x_a\) to denote the element \(h_{\bar{v}}^{-1}\gamma h_{\bar{v}}\in
T(\breve{F}_{\bar{v}})\). In contrast to \(\gamma\) versus \(\gamma_a\), the
distinction for \(x_a\) will not be very important and it will save some
notations. Then
\begin{align}
    \dot{w}_{x_a}=h_{\bar{v}}^{-1}\Frob_v(h_{\bar{v}})\in\Norm_{G}(T)(\breve{F}_{\bar{v}})
\end{align}
because \(a\) is generically regular semisimple, and we let
\(w_{x_a}\) be
the image of \(\dot{w}_{x_a}\) in \(W(\breve{F}_{\bar{v}})\cong\bW\). The
element \(w_{x_a}\) depends only on \(x_a\) but not
on \(\gamma_a\) nor \(h_{\bar{v}}\) (by contrast, \(\dot{w}_{x_a}\) \emph{does}
depend on such choices). Then \(x_a\) is fixed by
\(\Frob_v'=w_{x_a}\rtimes \Frob_v\),
    \nomenclature[\(sigma^prime_v \)]{\(\Frob_v'\)}{the image of \(\dot{\Frob}_v'\) in \(\bW\rtimes\Frob_v\)}
and \(\cM_{G,v}(a)\) may be regarded as a
\(k_v\)-structure on the \(\bar{k}\)-scheme
\(\cM_{\bG,\bar{v}}^{\le\lambda}(x_a)\) with
Frobenius \(\dot{\Frob}_v'\), where
\(\dot{\Frob}_v'=\dot{w}_{x_a}\rtimes \Frob_v\).
    \nomenclature[\(sigma'dot^prime_v \)]{\(\dot{\Frob}_v'\)}{a twisted Frobenius
    element in \(\Norm_{\bG}(\bT)\rtimes\Frob_v\) induced by a regular semisimple and unramified element of \(G\) or \(\FRM\)}
For convenience, we also define
\(\bar{\Frob}_v'=\bFn(\Frob_v')=\bFn(w_{x_a})\rtimes\Frob_v\)
    \nomenclature[\(n_bf \)]{\(\bFn\)}{the Tits section \(\bW\rtimes\Out(\bG)\to \Norm_{\bG}(\bT)\rtimes\Out(\bT)\)}
    \nomenclature[\(sigma'bar^prime_v \)]{\(\bar{\Frob}_v'\)}{the modified twisted Frobenius induced by \(\dot{\Frob}_v'\) and the Tits section}
to be the Tits
section of \(\Frob_v'\).

\subsection{}
Now we base change to \(x_G\) and so \(x_a\in\bar{\bT}_\bM\). Let
\begin{align}
    \bB'=\Frob_v'(\bB)
\end{align}
be another Borel containing \(\bT\) and \(\bU'\) its unipotent radical. Then
we may obtain functors \(\tilde{\KVU}_{\bG,\bar{v}}^{\le\lambda,\prime}(x_a)\),
\(\KVU_{\bG,\bar{v}}^{\le\lambda,\prime}(x_a)\), and so on by replacing \(\bU\)
with \(\bU'\), and \(\cM_{\bG,\bar{v}}^{\le\lambda}(x_a)\) is a
\(\CoCharG(\bT)\)-tiling of \(\KVU_{\bG,\bar{v}}^{\le\lambda,\prime}(x_a)\) as
well.

Since \(\dot{\Frob}_v'\) preserves \(\CoCharG(\bT)\), it acts on the stack
\(\Stack{\cM_{\bG,\bar{v}}^{\le\lambda}(x_a)/\CoCharG(\bT)}\), which induces
a \(\dot{\Frob}_v'\)-action
on the irreducible components of \(\Stack{\cM_{\bG,\bar{v}}^{\le\lambda}(x_a)/\CoCharG(\bT)}\).
The latter action factors through \(\Frob_v'\).
Similarly, the natural \(\bar{\Frob}_v'\)-action on \(\Loop_{\bar{v}}\bG\)
induces an isomorphism \(\Loop_{\bar{v}}\bU\to\Loop_{\bar{v}}\bU'\), which
descends to an isomorphism
\begin{align}
    \bar{\Frob}_v\colon
    \KVU_{\bG,\bar{v}}^{\le\lambda}(x_a)\Arc_{\bar{v}}{\bG}/\Arc_{\bar{v}}{\bG}\stackrel{\sim}{\longto}
    \KVU_{\bG,\bar{v}}^{\le\lambda,\prime}(x_a)\Arc_{\bar{v}}{\bG}/\Arc_{\bar{v}}{\bG}.
\end{align}
Therefore, \(\bar{\Frob}_v'\) induces a bijection
\begin{align}
    \bar{\Frob}_v'\colon
    \Irr(\KVU_{\bG,\bar{v}}^{\le\lambda}(x_a))\stackrel{\sim}{\longto}
    \Irr(\KVU_{\bG,\bar{v}}^{\le\lambda,\prime}(x_a)),
\end{align}
compatible with the action of \(\dot{\Frob}_v'\) on the irreducible components of
\(\Stack{\cM_{\bG,\bar{v}}^{\le\lambda}(x_a)/\CoCharG(\bT)}\).

\subsection{}
On the other hand, using geometric Satake (see \cites{MiVi07,Zh17}),
\eqref{eqn:diagram_connecting_GASF_to_MV} establishes an embedding of vector
spaces
\begin{align}
    \iota_\bU\colon
    \Qlb^{\oplus\Irr(\Stack{\cM_{\bG,\bar{v}}^{\le\lambda}(x_a)/\CoCharG(\bT)})}
    \longto \RHc^{\Pair{2\rho}{\nu_\gamma}}(\Gr_\bG,\IC^\lambda).
\end{align}
Similarly, by replacing \(\bU\) with \(\bU'\), we have identification
\begin{align}
    \iota_{\bU'}\colon
    \Qlb^{\oplus\Irr(\Stack{\cM_{\bG,\bar{v}}^{\le\lambda}(x_a)/\CoCharG(\bT)})}
    \longto \RHc^{\Pair{2\rho'}{\nu_\gamma}}(\Gr_\bG,\IC^\lambda).
\end{align}
Note that \(\Pair{2\rho}{\nu_\gamma}=\Pair{2\rho'}{\nu_\gamma}\).

The element \(\bar{\Frob}_v'\in\Arc_{\bar{v}}{\bG}\rtimes{\Frob_v}\) naturally acts
on \(\Gr_{\bG,\bar{v}}\), hence on its cohomologies. The complex \(\IC^\lambda\) is
supported on \(\Gr_{\bG,\bar{v}}^{\le\lambda}\), on which the
\(\Arc_{\bar{v}}{\bG}\)-action factors through a finite jet group.
By \Cref{lem:homotopy_lemma}, the action of
\(\bar{\Frob}_v'\) on \(\RHc^\bullet(\Gr_{\bG,\bar{v}},\IC^\lambda)\)
descends to a well-defined \(\Frob_v\)-action because jet groups of \(\bG\) are
connected. Thus, we have commutative diagram
\begin{equation}
    \begin{tikzcd}
        \Qlb^{\oplus\Irr(\Stack{\cM_{\bG,\bar{v}}^{\le\lambda}(x_a)/\CoCharG(\bT)})}
        \ar[r, "\iota_\bU"]\ar[d, "\dot{\Frob}_v'"'] & \RHc^{\bullet}(\Gr_{\bG,\bar{v}},\IC^\lambda) \ar[d, "\Frob_v"]\\
        \Qlb^{\oplus\Irr(\Stack{\cM_{\bG,\bar{v}}^{\le\lambda}(x_a)/\CoCharG(\bT)})}
        \ar[r, "\iota_{\bU'}"] & \RHc^{\bullet}(\Gr_{\bG,\bar{v}},\IC^\lambda)
    \end{tikzcd}
\end{equation}
If furthermore \(\nu_\gamma\) is central, then \(\iota_\bU\) is independent of
\(\bU\) by \Cref{cor:MV_TF_independent_of_B_for_zero_coweight}, so
\(\iota_\bU=\iota_{\bU'}\) becomes an embedding of Frobenius modules.
See \Cref{prop:stable_twisted_Frobenius_is_the_crystal_action} for the case
where \(\nu_\gamma\) is not necessarily central.

\subsection{}
Let \(\FRV_\lambda[\gamma]\)
    \nomenclature[\(V"frak_lambda_gamma \)]{\(\FRV_\lambda[\gamma]\)}{the space
    of \(\Qlb\)-valued functions on \(\Irr[\cM_{\bG,\bar{v}}^{\le\lambda}(x_a)]\), where \(x_a\in\bT(\breve{F}_{\bar{v}})\) is conjugate to \(\gamma\) over \(\breve{F}_{\bar{v}}\)}
be the space of \(\Qlb\)-valued functions on discrete set
\(\Irr[\cM_{\bG,\bar{v}}^{\le\lambda}(x_a)]\). It carries a natural
\(\CoCharG(\bT)\)-action,
and for any \(\kappa\in \dual{\bT}(\Qlb)=\Hom(\CoCharG(\bT),\Qlb^\x)\), we let
\(\FRV_\lambda[\gamma]_\kappa\)
    \nomenclature[\(V"frak_lambda_gamma_kappa \)]{\(\FRV_\lambda[\gamma]_\kappa\)}{the
    \(\kappa\)-isotypic subspace of \(\FRV_\lambda[\gamma]\), viewed as a \(\dot{\Frob}_v'\)-module}
be the \(\kappa\)-isotypic subspace, in other words,
\begin{align}
    \FRV_\lambda[\gamma]_\kappa
    =\Set*{f\in \FRV_\lambda[\gamma]\given
    f({\theta^{-1}\cdot Z})=\kappa(\theta)f(Z),\forall \theta\in\CoCharG(\bT)}.
\end{align}
When \(\kappa=1\), we also denote \(\FRV_\lambda[\gamma]_\kappa\) by
\(\FRV_\lambda[\gamma]_\hST\).
    \nomenclature[\(V"frak_lambda_gamma_st \)]{\(\FRV_\lambda[\gamma]_\hST\)}{\(\FRV_\lambda[\gamma]_\kappa\) when \(\kappa=1\)}
Note that as vector spaces, each
\(\FRV_\lambda[\gamma]_\kappa\) is isomorphic to
\(\Qlb^{\oplus\Irr(\Stack{\cM_{\bG,\bar{v}}^{\le\lambda}(x_a)/\CoCharG(\bT)})}\),
and when \(\kappa=1\) they are also isomorphic as \(\dot{\Frob}_v'\)-modules. In
general, \(\dot{\Frob}_v'\) acts on \(\FRV_\lambda[\gamma]_\kappa\) if
\(\kappa\) is fixed by \(\Frob_v'\).

\subsection{}
From now on, we assume that \(\kappa\) is fixed by \(\Frob_v'\).
The \(\dot{\Frob}_v'\)-action on \(\FRV_\lambda[\gamma]_\kappa\) is rather
complicated if \(\kappa\neq 1\), so we break it into two parts.
Write \(\dot{\Frob}_v'\) as
\begin{align}
    \dot{\Frob}_v'=\bFc \bar{\Frob}_v',
\end{align}
where \(\bFc\in \bT(\breve{F}_{\bar{v}})\) represents the
\(\Gal(\breve{F}_{\bar{v}}/F_v)\simeq\hat{\bbZ}\Frob_v'\)-cocycle
\begin{align}
    (\Frob_v')^i\longmapsto \prod_{j=0}^{i-1} (\Frob_v')^j(\bFc).
\end{align}
More generally, any element in \(\bT(\breve{F}_{\bar{v}})\) induces a
\(\hat{\bbZ}\Frob_v'\)-cocycle in the same way. By Tate--Nakayama duality,
\(\kappa\) can be evaluated at \(\bFc\). The element \(\bar{\Frob}_v'\)
still acts on \(\FRV_\lambda[\gamma]_\kappa\), and we denote the resulting
\(\bar{\Frob}_v'\)-module by \(\bar{\FRV}_\lambda[\gamma]_\kappa\).
    \nomenclature[\(V"frak_lambda_gamma_kappa_bar \)]{\(\bar{\FRV}_\lambda[\gamma]_\kappa\)}{same
    as \(\FRV_\lambda[\gamma]_\kappa\), but viewed as a \(\bar{\Frob}_v'\)-module instead}
So we have isomorphism
\begin{align}
    \FRV_\lambda[\gamma]_\kappa\otimes\kappa(\bFc)^{-1}\simeq
    \bar{\FRV}_\lambda[\gamma]_\kappa,
\end{align}
which is also equivariant with respect to the \(\dot{\Frob}_v'\)-action on the
left-hand side and \(\bar{\Frob}_v'\)-action on the right. Here
\(\otimes{\kappa(\bFc)^{-1}}\) means twisting the natural
\(\dot{\Frob}_v'\)-action by scalar \(\kappa(\bFc)^{-1}\).

\subsection{}
Choose an \(a\)-datum for the \(\Frob_v'\)-action on
\(\CoRoots\subset \CoCharG(\bT)\), then we have the \(1\)-cochain \(u_p\), which
in turn induces a cohomological class \(\lambda_T\in
\RH^1(F_v,I_{\gamma_a})\simeq\RH^1(\Frob_v',\bT)\), which is represented by
cocycle \(u_p(\Frob_v')\bFc^{-1}\)(see \Cref{sec:story_on_G_side}, especially
\Cref{sub:reformulation_of_lambda_T_in_transfer_factor}).

Suppose \(\kappa\) is now a part of endoscopic datum \((H,\kappa,\xi)\), and
\(\gamma\) is matched with stable conjugacy class \(\gamma_H\) in
\(\Stack*{\FRM_H/Z_\FRM^\kappa}(F_v)\), then
we have
\begin{align}
    \kappa(\bFc)^{-1} =
    \Delta_{\symup{I}}(\gamma_H,\gamma_G)\kappa(u_p(\Frob_v'))^{-1}
\end{align}
for any \(\gamma_G\) stably conjugate to \(\gamma\). There is some subtlety
here: it depends on \(\gamma\) being an element in \(G\), and \((H,\kappa,\xi)\)
may not be induced by an endoscopic group of \(G^\AD\). However, the
definition of \(\Delta_{\symup{I}}\) does not depend on
the admissible embedding \(\xi\) and only depends on the image of \(\kappa\) in
\(\dual{\bG}^\AD\), so \(\Delta_{\symup{I}}(\gamma_H,\gamma_G)\) is actually 
well-defined even for \(\gamma_a\).

Since the conjugate class is unramified, we may choose the \(a\)-datum in
\(\breve{\cO}_{\bar{v}}^\x\subset\breve{F}_{\bar{v}}^\x\) (for example, let
\(a_{\CoRt}=\Rt(\bFx)\) for any element \(\bFx\in
\Lie(\bT)(\breve{\cO}_{\bar{v}})^{\Frob_v'}\) that is regular at \(\bar{v}\)).
We shall call such \(a\)-data \notion{unramified}\index{\(a\)-datum!unramified}\index{unramified!\(a\)-datum}.
With this choice, \(u_p(\Frob_v')\in\bT(\breve{\cO}_{\bar{v}})\) and by Lang's theorem we may
find \(t_u\in\bT(\breve{\cO}_{\bar{v}})\) such that \(u_p(\Frob_v')=t_u\Frob_v'(t_u)^{-1}\),
therefore \(\kappa(u_p(\Frob_v'))=1\). We have thus proved the following result:

\begin{lemma}
    \label[lemma]{lem:Delta_I_in_MV}
    Suppose the \(a\)-datum is unramified, then we have
    isomorphism of \(\dot{\Frob}_v'\)-modules
    \begin{align}
        \FRV_\lambda[\gamma]_\kappa\otimes \Delta_{\symup{I}}(\gamma_H,\gamma_G)
        \stackrel{\sim}{\longto}\bar{\FRV}_\lambda[\gamma]_\kappa,
    \end{align}
    where \(\dot{\Frob}_v'\) acts on the right-hand side through
    \(\bar{\Frob}_v'\).
\end{lemma}

\subsection{}
The action of \(\bar{\Frob}_v'\) on \(\bar{\FRV}_\lambda[\gamma]_\kappa\) is
more subtle and is related to the theory of crystal bases, which is heavy in
combinatorial details. In fact, we have already seen such connection at work in
the proof of \Cref{lem:MV_transfer_factor_key}. To make it more digestible,
however, we will postpone further discussion until
\Cref{sec:Connection_with_Kashiwara_Crystals}.

On the other hand, observe that the whole process above can be done for any
representation \(V\in\Rep_{\LD{G}}^\Alg\) in place of
\(V_\lambda\). For simplicity, we assume that the monoid \(\FRM\)
corresponds to cocharacters in the support of \(V\) (in other words, the
cocharacter cone of \(\bA_\bM\) is freely generated by \(-w_0(\lambda)\) for
\(\lambda\) ranging over the support of \(V\)), that the sheaf \(\IC^\lambda\) is
replaced with the Satake sheaf corresponding to \(V\), and the space spanned by
irreducible components of
\(\Stack{\cM_{\bG,\bar{v}}^{\le\lambda}(x_a)/\CoCharG(\bT)}\) is replaced by
\begin{align}
    \bigoplus_\lambda
    \Qlb^{\oplus\Irr\Stack{\cM_{\bG,\bar{v}}^{\le\lambda}(x_a)/\CoCharG(\bT)}}\otimes
    \Hom_{\dual{\bG}}(V_\lambda,V).
\end{align}
Similarly, \(\FRV_\lambda[\gamma]_\kappa\) (resp.~\(\bar{\FRV}_\lambda[\gamma]_\kappa\)) is replaced by
\begin{align}
    \FRV[\gamma]_\kappa &\defeq\bigoplus_\lambda \FRV_\lambda[\gamma]_\kappa\otimes
    \Hom_{\dual{\bG}}(V_\lambda,V)\\
    \text{resp. }\quad
    \bar{\FRV}[\gamma]_\kappa &\defeq\bigoplus_\lambda \bar{\FRV}_\lambda[\gamma]_\kappa\otimes
    \Hom_{\dual{\bG}}(V_\lambda,V).
    \nomenclature[\(V"frak_gamma_kappa \)]{\(\FRV[\gamma]_\kappa,\bar{\FRV}[\gamma]_\kappa\)}{the
    generalization of \(\FRV_\lambda[\gamma]_\kappa\) or
\(\bar{\FRV}_\lambda[\gamma]_\kappa\) for a general \(\LD{G}\)-representation \(V\)}
\end{align}
Note that any direct summand involving \(\lambda\) is \(0\) if \(\nu_\gamma\le
\lambda\) does not hold.

\subsection{}
Let us now turn to endoscopic groups. By definition, we have an admissible
embedding
\begin{align}
    \xi\colon\LD{H}\longto \LD{G},
\end{align}
and we denote the image \(\xi(\Frob_v)\) of \(\Frob_v\)
by \(\Frob_H\).
We may restrict \(V\) to \(\LD{H}\) using \(\xi\). Further
restricting to \(\dual{\bH}\), we have canonical decomposition
\begin{align}
    \Res_{\dual{\bH}}^{\LD{G}}V\simeq
    \bigoplus_{\lambda_H}V_{\lambda_H}^H\otimes\Hom_{\dual{\bH}}(V_{\lambda_H}^H,V),
\end{align}
where \(V_{\lambda_H}^H=\RH^\bullet(\Gr_H,\IC_H^{\lambda_H})\) and \(\lambda_H\)
ranges over highest weights of \(\dual{\bH}\).

\subsection{}
Suppose \(a\) lies in the image of \(\Stack{\FRC_{\FRM,H}/Z_\FRM^\kappa}(\cO_v)\).
By the proof of
\Cref{lem:arc_lifting_to_endoscopic_monoid}, there are only finitely many lifts
of \(a\) to \(\Stack{\FRC_{\FRM,H}/Z_\FRM^\kappa}(\cO_{v})\), one for each
\(\lambda_{H}\) such that \(\nu_\gamma\le_\bH\lambda_{H}\). As
\(F_{v}\)-points these lifts are all equal, and we denote this common
point by \(a_H\). Using the Cartesian diagram
\begin{equation}
    \begin{tikzcd}
        \Stack*{\FRM_H/Z_\FRM^\kappa} \ar[r]\ar[d] & \Stack*{\Env(H^\SC)/Z_{\Env(H^\SC)}} \ar[d]\\
        \Stack*{\FRC_{\FRM,H}/Z_\FRM^\kappa} \ar[r] & \Stack*{\FRC_{\Env(H^\SC)}/Z_{\Env(H^\SC)}}
    \end{tikzcd}
\end{equation}
we see that since the image of \(a_H\) in \(\Stack*{\FRC_{\Env(H^\SC)}/Z_{\Env(H^\SC)}}\) always lifts
to a point in \(\Stack*{\Env(H^\SC)/Z_{\Env(H^\SC)}}(F_v)\) (due to
\(Z_{\Env(H^\SC)}\) being an induced torus, and
\Cref{thm:A_2m_rationality_fix}), we can always lift
\(a_H\) to a point in \(\Stack*{\FRM_H/Z_\FRM^\kappa}(\cO_v)\), which induces a
point in \(H/Z_G(F_v)\). It is now very crucial to assume that it further lifts to
\(\gamma_H\in H(F_v)\) with the same Newton point \(\nu_\gamma\). This
assumption will always be satisfied in any computations we will perform relevant
to fundamental lemma.

\subsection{}
Recall that by choosing a \(\chi\)-datum, we have an admissible embedding
\begin{align}
    \xi_T\colon
    \LD{\FRJ}_a\simeq\dual{\bT}\rtimes\Ggen{\Frob_v'}\longto
    \LD{G},
\end{align}
and we write \(\tilde{\Frob}_v'=\xi_T(\Frob_v')=\tilde{w}\rtimes\Frob_v\).
Clearly, \(\tilde{\Frob}_v'\) acts on \(V_\lambda\) and stabilizes
\(V_\lambda[\nu_\gamma]\).

The \(\chi\)-datum restricts to a \(\chi\)-datum of the \(\Frob_v'\)-action on
\(\Roots_{\bH}\subset\CharG(\bT)\) through canonical embedding
\(\bW_\bH\rtimes\Ggen{\Frob_v}\to\bW\rtimes\Ggen{\Frob_v}\) in
\Cref{lem:kappa_G_H_compatibility}. Thus, we also have admissible embedding
\begin{align}
    \xi_{T_H}\colon \LD{\FRJ}_{H,a_H}\longto \LD{H}.
\end{align}

The difference between \(\xi\circ\xi_{T_H}\) and \(\xi_T\) induces cohomology
class \(\bFa\in\RH^1(W_{F_v},\dual{\FRJ}_a)\) as seen in
\Cref{sub:def_of_Delta_III_2}. Here \(W_{F_v}\) is the Weil group of
\(F_v\), not the Weyl group of \(G\) restricted to \(F_v\). The class \(\bFa\)
may be paired with any element \(\gamma_G\in\FRJ_a(F_v)\simeq
\bT(\breve{F}_{\bar{v}})^{\Frob_v'}\) under Langlands correspondence for tori.

Since everything is unramified, we may choose an
\notion{unramified}\index{\(\chi\)-datum!unramified}\index{unramified!\(\chi\)-datum}
\(\chi\)-datum as follows (see \Cref{sec:cohomological_notations} for
notations): for each asymmetric orbit of \(\Frob_v'\) on
\(\Roots_G\), the \(\chi\)-datum is constantly trivial, while for each symmetric
orbit generated by root \(\Rt\), we let \(\chi_\Rt\) be the character on
\(F_{\Rt}^\x\) that equals \((-1)^{\val_{F_v}}\). With this choice of
\(\chi\)-datum, and if we assume that \(\gamma_G\) is \(\nu\)-regular semisimple,
the pairing is given by evaluation
\begin{align}
    \Pair{\bFa}{\gamma_G}=\nu_{\gamma_G}(\bFa(\Frob_v)),
\end{align}
where \(\bFa(\Frob_v)\in\dual{\FRJ}_a\) and \(\nu_{\gamma_G}\) is the image of
\(\gamma_G\) in \(\CoCharG(\bT)\simeq \bT(\breve{F}_{\bar{v}})/\bT(\breve{\cO}_{\bar{v}})\).
If \(\gamma_G\) is not \(\nu\)-regular semisimple, the pairing is still not hard to
compute, but we will skip the formula here.

\subsection{}
Here is where the existence of \(\gamma\) (rather than just \(\gamma_a\)) and
\(\gamma_H\) becomes important.
Let \(\gamma_G=\gamma\) be \(\nu\)-regular semisimple, then we have
\begin{align}
    \Pair{\bFa}{\gamma}=\Delta_{\symup{III}_2}(\gamma_H,\gamma).
\end{align}
When \(\gamma\) is not \(\nu\)-regular semisimple, we want to consider the
right-hand side above, but the left-hand
side is helpful for computations in special cases.

\subsection{}
With our choice of \(a\)-datum and \(\chi\)-datum, the last piece of transfer
factor, namely \(\Delta_{\symup{II}}\), has a simple description by unpacking
the definition:
\begin{align}
    \Delta_{\symup{II}}(\gamma_H,\gamma)=\prod_\Rt
    (-1)^{\val_{v}(1-\Rt(\gamma))},
\end{align}
where \(\Rt\) ranges over representatives of \(\Frob_v'\)-symmetric orbits of
roots of \(G\) not in \(H\).
However, roots being symmetric means that it is perpendicular to \(\nu_\gamma\),
and if \(\gamma\) is \(\nu\)-regular semisimple, we have
\(\Delta_{\symup{II}}(\gamma_H,\gamma)=1\), otherwise it is not hard to compute
by looking at root valuations.

\subsection{}
By replacing \(G\) with \(H\), and restrict \(V\) to \(\LD{H}\) via \(\xi\), we have
\begin{align}
    \bar{\FRV}^H[\gamma_H]_\hST^\xi
    =\FRV^H[\gamma_H]_\hST^\xi
    =\bigoplus_{\lambda_H}\FRV_{\lambda_H}^H[\gamma_H]_\hST\otimes\Hom_{\dual{\bH}}(V_{\lambda_H}^H,V),
    \nomenclature[\(V"frak_lambda_H_gamma_H_xi_st \)]{\(\FRV^H[\gamma_H]_\hST^\xi\)}{the
    analogue of \(\FRV[\gamma]_\hST\) for an endoscopic group \(H\) via an \(L\)-embedding \(\xi\)}
\end{align}
where the superscript \(\xi\) is to emphasize the dependence on \(\xi\).
Finally, we define yet another \(\Frob_v'\)-module \(\FRV[\gamma]_\kappa^\xi\)
which is supposed to ``match'' \(\FRV^H[\gamma_H]_\hST^\xi\), and this
matching is viewed as a skeletal shadow of the fundamental lemma (see
\Cref{thm:asymptotic_FL} below):
\begin{align}
    \FRV[\gamma]_\kappa^\xi\defeq
    \bigoplus_\lambda
    \bar{\FRV}_\lambda[\gamma]_\kappa\otimes
    \Hom_{\dual{\bG}}(V_\lambda,V)\otimes
    \Delta_{\symup{II}}\Delta_{\symup{III}_2}(\gamma_H,\gamma),
\end{align}
It depends on \(\xi\), but not on the specific
\(\chi\)-datum we chose above, because
\(\Delta_{\symup{II}}\Delta_{\symup{III}_2}\) does not.
When \(\gamma\) is \(\nu\)-regular semisimple, we have a slightly simper
expression
\begin{align}
    \FRV[\gamma]_\kappa^\xi\defeq
    \bigoplus_\lambda
    \bar{\FRV}_\lambda[\gamma]_\kappa\otimes
    \Hom_{\dual{\bG}}(V_\lambda,V)\otimes \Pair{\bFa}{\gamma}.
\end{align}

\subsection{}
We aim to prove the following statement:
\begin{theorem}[Asymptotic Fundamental Lemma]
    \label[theorem]{thm:asymptotic_FL}
    We have an isomorphism of Frobenius modules
    \begin{align}
        \FRV[\gamma]_\kappa^\xi\simeq\FRV^H[\gamma_H]_\hST^\xi
    \end{align}
    where \(\gamma_{H}\) is any fixed \(F_v\)-point of \(H\) lying over \(a_H\).
\end{theorem}
\begin{proof}
    This looks like a difficult result to prove purely combinatorially. We will
    first prove some important special cases in
    \Cref{sec:Connection_with_Kashiwara_Crystals}, and the general cases
    will be proved simultaneously with the fundamental lemma itself as their
    proofs intertwine with each other.
\end{proof}

\begin{remark}
    It would be interesting to extend the asymptotic fundamental
    lemma to ramified elements.
\end{remark}

\section{Connection with Kashiwara Crystals} 
\label{sec:Connection_with_Kashiwara_Crystals}

In this section we describe in detail the Frobenius action on the vector space
\(\bar{\FRV}_\lambda[\gamma]_\kappa\) or more generally
\(\bar{\FRV}[\gamma]_\kappa\) using Kashiwara crystals. See
\Cref{chapA:Review_of_Kashiwara_Crystals} for a basic review of crystals. I
thank Xinwen Zhu for drawing my attention to the topic of crystal basis.

\subsection{}
As before, we fix a closed point \(v\in X\) and without loss of generality may
assume \(k_v=k\). We fix a pinning at \(v\) so that there is a canonical
Borel subgroup \(\bB\subset\bG\) and its unipotent radical \(\bU\). Let
\(\lambda\) be an \(F_v\)-rational highest coweight of \(G\) relative to \(\bU\).
Let \(\gamma_\lambda\in \FRM^\x(F_v)\cap\FRM(\cO_v)\) be an element whose
boundary divisor is \(-w_0(\lambda)\) and its image in \(G^\AD\) coincides with
\(\gamma\). Let \(a\in\FRC_\FRM(\cO_v)\) be the image of
\(\gamma_\lambda\) and let \(x_a\in T_\FRM(\breve{F}_{\bar{v}})\) be a fixed
point over \(a\). Let \(\nu\in\CoCharG(\bT^\AD)\) be such that
\begin{align}
    x_a\in \pi_v^{(-w_0(\lambda),\nu_\AD)}T_\FRM(\breve{\cO}_{\bar{v}})
\end{align}
and \(\gamma\) is conjugate to an element in
\(\pi_v^{\nu_\AD}T^\AD(\breve{\cO}_{\bar{v}})\). Recall that the choice of \(x_a\)
determines a twisted Frobenius \(\Frob_v'\in
\Norm_\bG(\bT)(\breve{F}_{\bar{v}})\rtimes\Frob_v\) as well as a modified version
\(\bar{\Frob}_v'\in\Norm_{\bG^\SC}(\bT^\SC)(\bar{k})\rtimes\Frob_v\) defined using Tits
section. Since we have already studied the effect of the difference
\(\bar{\Frob}_v'(\Frob_v')^{-1}\), we are only interested in the action of
\(\bar{\Frob}_v'\). For simplicity, we denote
\(\Frob'\defeq\bar{\Frob}_v'=\dot{w}\rtimes\Frob_v\).
\nomenclature[\(sigma'_^prime \)]{\(\Frob'\)}{same as \(\bar{\Frob}_v'\) (not \(\dot{\Frob}_v'\) or \(\Frob_v'\))}

\subsection{}
We have seen in \Cref{sec:connection_with_mv_cycles} that the geometric
irreducible
components of \(\cM_{\bG,\bar{v}}^{\le\lambda}(x_a)\) may be identified with
\begin{align}
    \sB_\lambda[\nu]\x\CoCharG
\end{align}
where \(\sB_\lambda[\nu]\) is the set of MV-cycles of highest weight \(\lambda\)
and weight \(\nu\). We identify \(\sB_\lambda[\nu]\) with the crystal basis
elements of weight \(\nu\) in the highest-weight normal crystal \(\sB_\lambda\)
of \(\dual{\bG}\). Then we have isomorphism of vector spaces
\begin{align}
    \FRV_\lambda[\gamma]_\hST\simeq
    \bar{\FRV}_\lambda[\gamma]_\hST\simeq\Qlb^{\oplus\sB_\lambda[\nu]}.
\end{align}
Note that \textit{a priori} the last isomorphism depends on
\(\bB\), which we denote by \(\iota_\bB\).

\begin{proposition}
    \label[proposition]{prop:stable_twisted_Frobenius_is_the_crystal_action}
    When \(\kappa=1\), the \(\Frob'\)-action on
    \(\Qlb^{\oplus\sB_\lambda[\nu]}\) induced by
    \(\bar{\FRV}_\lambda[\gamma]_\hST\) is the same as the action induced by the
    natural \(\Frob'\)-action (see \Cref{sub:Appendix_actions_on_crystals}) on \(\sB_\lambda\) as a crystal.
\end{proposition}
\begin{proof}
    This is essentially an extension of
    \Cref{cor:MV_TF_independent_of_B_for_zero_coweight} to the case where
    \(\nu\) is not central. Let \(\Rt\) be a simple
    root and \(s_\Rt\) the corresponding simple reflection. Then as before we
    have the bijection of crystal basis elements
    \begin{align}
        \bFs_\Rt^\lambda\colon \sB_\lambda[\nu]\longto   \sB_\lambda[s_\Rt(\nu)]
    \end{align}
    by ``changing unipotent'' from \(\bU\) to \(s_\Rt(\bU)\). It fits into the
    commutative diagram
    \begin{equation}
        \begin{tikzcd}
            \Qlb^{\oplus\sB_\lambda[\nu]} \ar[r, "\bFs_\Rt^\lambda"]\ar[d, "\iota_\bB", swap] &
            \Qlb^{\oplus\sB_\lambda[s_\Rt(\nu)]} \ar[d, "\iota_{s_\Rt(\bB)}"]\\
            \bar{\FRV}_\lambda[\gamma]_\hST \ar[r, equal] & \bar{\FRV}_\lambda[\gamma]_\hST
        \end{tikzcd}
    \end{equation}
    Such diagram is
    compatible with restricting the crystal to the rank-\(1\) Levi subgroup
    \(\bL_\Rt\) generated by \(\Rt\). In other words, restricting
    \(\FRV_\lambda\) to \(\bL_\Rt\), then \(\bFs_\Rt^\lambda\) is the same as
    the one induced by changing unipotent \(\bU_\Rt\subset\bL_\Rt\) to its
    opposite. By direct computation on rank-\(1\) group, we see that
    \(\bFs_\Rt^\lambda\) is the same as the Weyl group (of \(\bL_\Rt\)) action
    on normal crystals.

    Using the same inductive process as in 
    \Cref{cor:MV_TF_independent_of_B_for_zero_coweight}, we write \(w=s_1\cdots
    s_m\) be a shortest expression of \(w\) relative to \(\bB\). Let
    \(r_i=s_1\cdots s_i\) and \(\bB_i=r_i(\bB)\) (where we
    use the convention \(r_0=1\) and \(\bB_0=\bB\)). Then each
    \(\bB_i\) is obtained from \(\bB_{i-1}\) by applying a \(\bB_{i-1}\)-simple
    reflection \(s_i'=r_{i-1}s_ir_{i-1}^{-1}\). Let \(\bFs_i'\) be the
    analogue of \(\bFs_\Rt\) defined for \(s_i'\) and relative to Borel subgroup
    \(\bB_{i-1}\). As actions on crystal basis, \(\bFs_i'\) is the same as
    \(\bFs_i\), the one defined by \(s_i\) relative to \(\bB\).
    We then have the following diagram
    \begin{equation}
        \begin{tikzcd}
            & &\bar{\FRV}_\lambda[\gamma]_\hST \ar[lld, "\iota_\bB",swap]\ar[ld,
            "\iota_{\bB_1}"]
            \ar[d,phantom,"\cdots"]
            \ar[rd,"\iota_{\bB_{m-1}}",swap]\ar[rrd,"\iota_{\bB_{m}}"]
            & & \\
            \Qlb^{\oplus\sB_\lambda[\nu]}\ar[r,"\bFs_1", swap]
            & \Qlb^{\oplus\sB_\lambda[s_1(\nu)]}\ar[r,"\bFs_2", swap] 
            & \cdots\ar[r,"\bFs_{m-1}", swap] 
            & \Qlb^{\oplus\sB_\lambda[s_{m-1}\cdots s_1(\nu)]}\ar[r,"\bFs_{m}", swap] 
            & \Qlb^{\oplus\sB_\lambda[w^{-1}(\nu)]}.
        \end{tikzcd}
    \end{equation}
    In other words, \(\iota_{w(\bB)}=\iota_{\bB_m}=w^{-1}\circ\iota_\bB\). On
    the other hand, we have seen in \Cref{sec:connection_with_mv_cycles} that
    \(\iota_{\bB}\) and \(\iota_{\bB_m}\) intertwines the \(\Frob'\)-action on
    \(\bar{\FRV}_\lambda[\gamma]_\hST\) with the \(\Frob_v\)-action on cohomology
    group \(\RH^\bullet(\Gr_{\bG,\bar{v}},\IC^\lambda)\), the latter of which is
    induced by the crystal action of \(\Frob_v\). This shows that
    \begin{align}
        \Frob_v\circ\iota_\bB=\iota_{\bB_m}\circ\Frob'=w^{-1}\circ\iota_\bB\circ\Frob',
    \end{align}
    where \(\Frob_v\) on the left-hand side and \(w^{-1}\) on the right-hand
    side mean their respective crystal actions. This implies that
    \begin{align}
        w\circ\Frob_v\circ\iota_\bB=\iota_\bB\circ\Frob',
    \end{align}
    as desired.
\end{proof}

\subsection{}
The result for general \(\kappa\) is not very difficult compared to the
\(\kappa=1\) case. To start, we note that if
\begin{align}
    (\FRx,\theta)\in \sB_\lambda[\nu]\x\CoCharG
\end{align}
is an irreducible component of \(\cM_{\bG,\bar{v}}^{\le\lambda}(x_a)\), then
by \Cref{prop:stable_twisted_Frobenius_is_the_crystal_action}, \(\Frob'\)
sends this component to
\begin{align}
    (\Frob'(\FRx), \Frob'(\theta)+\mu_\FRx)
\end{align}
for some \(\mu_\FRx\in\CoCharG\) independent of \(\theta\). This implies that the
Frobenius action on \(\bar{\FRV}_\lambda[\gamma]_\kappa\) sends the image of
\(\FRx\) (i.e., the function sending \((\FRx,0)\) to \(1\)) to
\(\kappa(\mu_\FRx)\Frob'(\FRx)\). The remaining task is then figuring out how to
describe \(\mu_\FRx\). To this end we define for each simple coroot
\(\CoRt\in\SimCoRts\) and each crystal basis element \(\FRx\in\sB_\lambda\) a
number
\begin{align}
    n_{\CoRt,\FRx}=\min\Set*{\epsilon_{\CoRt}(\FRx),\phi_{\CoRt}(\FRx)},
\end{align}
where \(\epsilon_{\CoRt}\) and \(\phi_{\CoRt}\) are the string length functions
associated with crystal operators \(e_{\CoRt}\) and \(f_{\CoRt}\) respectively.
Clearly, \(n_{\CoRt,\FRx}=n_{\Frob_v(\CoRt),\Frob_v(\FRx)}\).
For any (not necessarily simple) root \(\Rt\), we also define a number
\begin{align}
    d_\Rt^0(x_a)=\begin{cases}
        0 & \Pair{\Rt}{\nu}\neq 0,\\
        \val_{F_v}(1-\Rt(x_a)) & \Pair{\Rt}{\nu}=0.
    \end{cases}
\end{align}

\begin{proposition}
    \label[proposition]{prop:kappa_twisted_Frobenius_and_n_Rt_x}
    Let \(w=s_1\cdots s_m\) be a shortest expression of \(w\), and
    \(r_i=s_1\cdots s_i\) (with the convention \(r_0=1\)), then we have
    \begin{align}
        \mu_\FRx=\sum_{i=1}^m
        \bigl(n_{\CoRt_i,r_{i-1}^{-1}\Frob'(\FRx)}+d_{r_{i-1}(\Rt_i)}^0(x_a)\bigr)r_{i-1}(\CoRt_i).
    \end{align}
    In particular, if \(\gamma\) is \(\nu\)-regular semisimple, we have
    \begin{align}
        \mu_\FRx=\sum_{i=1}^m
        n_{\CoRt_i,r_{i-1}^{-1}\Frob'(\FRx)}r_{i-1}(\CoRt_i).
    \end{align}
\end{proposition}
\begin{proof}
    The proof is just a modification of that of
    \Cref{prop:stable_twisted_Frobenius_is_the_crystal_action}. Instead of
    analyzing \(\bFs_\Rt^\lambda\) only on the crystal basis, we take the lattice
    part into account. Such analysis is not difficult using diagram
    \eqref{eqn:Levi_induction_of_MV_big_diagram}. For groups of
    semisimple rank \(1\), the split form of universal monoid is just
    \(\Mat_2\cong\La{gl}_2\), hence the (geometric) multiplicative affine Springer
    fiber is essentially the same as the usual affine Springer fibers for 
    \(\La{gl}_2\), and they are already computed in \cite{GKM04}.

    Explicitly, by the same argument as in
    \Cref{prop:stable_twisted_Frobenius_is_the_crystal_action}, to analyze the
    action of \(\bFs_\Rt^\lambda\), it suffices to replace
    \(\bG\) (resp.~\(\bB\)) by \(\bL_\Rt\) (resp.~\(\bB_\Rt\)) and so \(\bG\) is
    now of semisimple rank \(1\) (and \(\lambda\) will be replaced by certain
    \(\bL_\Rt\)-highest coweight depending on \(\FRx\)). In addition, we may take the quotient of \(\bG\)
    by its connected center, so that \(\bG\) is either \(\SL_2\) or \(\PGL_2\).
    The universal monoid of \(\bG^\SC\) is now \(\Mat_2\), and \(x_a\) is a
    diagonal matrix with \emph{distinct} (by its regular semisimplicity) entries
    \(x_1,x_2\in\breve{\cO}_{\bar{v}}\cap\breve{F}_{\bar{v}}^\x\).
    By \cite{GKM04}*{Lemma~8.4}, we have
    \begin{align}
        \cM_{\bG,\bar{v}}^{\le\lambda}(x_a)(\bar{k})
        =\bigcup_{d=-\val_{F_v}(x_1-x_2)}^0\bT(\breve{F}_{\bar{v}})
        \begin{pmatrix}
            1 & \pi_v^{-d}\\
             & 1
        \end{pmatrix}
        \bG(\breve{\cO}_{\bar{v}})/\bG(\breve{\cO}_{\bar{v}}).
    \end{align}
    In this case, \(\FRx\) is the unique element in  \(\sB_\lambda[\nu]\)
    and corresponds to the closure of any
    \(\bT(\breve{\cO}_{\bar{v}})\)-orbit when \(d\) reaches its
    minimum. Similarly, we also have
    \begin{align}
        \cM_{\bG,\bar{v}}^{\le\lambda}(x_a)(\bar{k})
        =\bigcup_{d=-\val_{F_v}(x_1-x_2)}^0\bT(\breve{F}_{\bar{v}})
        \begin{pmatrix}
            1 & \\
            \pi_v^{-d} & 1
        \end{pmatrix}
        \bG(\breve{\cO}_{\bar{v}})/\bG(\breve{\cO}_{\bar{v}}),
    \end{align}
    and \(s_\Rt(\FRx)\) also corresponds to any
    \(\bT(\breve{\cO}_{\bar{v}})\)-oribt when \(d\) reaches its minimum. Using
    the Steinberg relation, we have for any \(d\),
    \begin{align}
        \begin{pmatrix}
            \pi_v^d & \\
             & \pi_v^{-d}
        \end{pmatrix}
        \begin{pmatrix}
            1 & \pi_v^{-d}\\
             & 1
        \end{pmatrix}
        \bG(\breve{\cO}_{\bar{v}})
        = 
        \begin{pmatrix}
            1 & \\
            \pi_v^{-d} & 1
        \end{pmatrix}
        \bG(\breve{\cO}_{\bar{v}}).
    \end{align}

    Therefore, if \((\FRx,0)\) and \((s_\Rt(\FRx),\mu_\Rt)\) represent (again,
    via \(\bB\) and \(s_\Rt(\bB)\) respectively) the same irreducible component
    of \(\cM_{\bG,\bar{v}}^{\le\lambda}(x_a)\), then we have
    \begin{align}
        \mu_\Rt=-\val_{F_v}(x_1-x_2)\CoRt =-\bigl(n_{\CoRt,\FRx}+d_{\Rt}^0(x_a)\bigr)\CoRt,
    \end{align}
    where the second equality is obtained by simply writing down the definition of
    discriminant valuations.
    By induction, if \((\FRx,0)\) and \((w^{-1}(\FRx),\mu)\) represent the same
    irreducible component \(S\) of \(\cM_{\bG,\bar{v}}^{\le\lambda}(x_a)\), via
    \(\bB\) and \(w(\bB)\) respectively, then we have
    \begin{align}
        \mu=-\sum_{i=1}^m
        \bigl(n_{\CoRt_i,r_{i-1}^{-1}(\FRx)}+d_{r_{i-1}(\Rt_i)}^0(x_a)\bigr)r_{i-1}(\CoRt_i).
    \end{align}
    Since \(\Frob'(S)\) is identified with \((\Frob_v(\FRx),0)\) with the choice
    of Borel being \(w(\bB)\), which is in turn identified with
    \((\Frob'(\FRx),\mu_\FRx)\) under \(\bB\). Substituting \(\FRx\) by
    \(\Frob'(\FRx)\) in the expression of \(\mu\) we obtain the desired result.
\end{proof}

\begin{corollary}
    \label[corollary]{cor:kappa_twisted_Frobenius_for_central_coweight}
    If \(\nu\) is central and \(\gamma\) is \(\nu\)-regular semisimple, then
    \begin{align}
        \mu_\FRx=\Frob_v(\rho_\FRx)-\Frob'(\rho_\FRx),
    \end{align}
    where
    \begin{align}
        \rho_\FRx=\sum_{\CoRt_i\in\SimCoRts}n_{\CoRt_i,\FRx}\CoWt_i.
    \end{align}
    In particular, if in addition \(G\) is split and adjoint, then
    \(\bar{\FRV}_\lambda[\gamma]_\kappa=\bar{\FRV}_\lambda[\gamma]_\hST\) is a
    trivial \(\Frob'\)-module.
\end{corollary}
\begin{proof}
    When \(\nu\) is central, then we necessarily have
    \(r_{i-1}^{-1}\Frob'(\FRx)=\Frob_v(\FRx)\) by the definition of Weyl group
    action on normal crystals (cf.~\Cref{sub:Appendix_actions_on_crystals}).
    \Cref{prop:kappa_twisted_Frobenius_and_n_Rt_x}
    then implies that
    \begin{align}
        \mu_\FRx=\sum_{i=1}^m n_{\CoRt_i,\Frob_v(\FRx)}r_{i-1}(\CoRt_i)
        =(1-w)\biggl(\sum_{\Rt_i\in\SimRts}n_{\CoRt_i,\Frob_v(\FRx)}\CoWt_i\biggr).
    \end{align}
    Then the corollary follows from the fact that
    \begin{align}
        \sum_{\Rt_i\in\SimRts}n_{\CoRt_i,\Frob_v(\FRx)}\CoWt_i
        =\sum_{\Rt_i\in\SimRts}n_{\Frob_v^{-1}(\CoRt_i),\FRx}\CoWt_i
        =\sum_{\Rt_i\in\SimRts}n_{\CoRt_i,\FRx}\Frob_v(\CoWt_i)
        =\Frob_v\biggl(\sum_{\Rt_i\in\SimRts}n_{\CoRt_i,\FRx}\CoWt_i\biggr)
        =\Frob_v(\rho_\FRx).
    \end{align}
    For the last claim, if \(G\) is also split and adjoint, then \(\Frob_v=1\),
    and \(\rho_\FRx\in\CoCharG\), therefore \(\kappa(\mu_\FRx)=1\) because
    \(\kappa\) is fixed by \(\Frob'\). Since in this case we necessarily have
    \(\nu=0\) and Weyl group acts trivially on \(\sB_\lambda[0]\), the claim
    follows.
\end{proof}

\begin{corollary}
    If \(\nu\in W\lambda\) and \(\gamma\) is \(\nu\)-regular semisimple, then
    \(\mu_\FRx=0\). Consequently,
    \(\bar{\FRV}_\lambda[\gamma]_\kappa=\bar{\FRV}_\lambda[\gamma]_\hST\) are
    trivial \(\Frob'\)-modules.
\end{corollary}
\begin{proof}
    For any \(\theta\in W\lambda\), the set \(\sB_\lambda[\theta]\) is a
    singleton \(\Set{\FRy}\), and we necessarily have either \(\epsilon_{\CoRt}(\FRy)=0\) 
    or \(\phi_{\CoRt}(\FRy)=0\) for any simple coroot \(\CoRt\), because otherwise
    \(\lambda+\CoRt'\) would be a weight of \(\sB_\lambda\) for some
    positive coroot \(\CoRt'\). This implies that \(n_{\CoRt,\FRy}=0\), and the
    corollary follows from
    \Cref{prop:stable_twisted_Frobenius_is_the_crystal_action,prop:kappa_twisted_Frobenius_and_n_Rt_x}.
\end{proof}

\begin{corollary}
    \label[corollary]{cor:special_case_of_unramified_GASF_irreducible_components_Frob_transfer}
    \Cref{thm:asymptotic_FL} holds when
    \(\gamma\) is \(\nu\)-regular semisimple and either
    \begin{enumerate}
        \item \(G\) is split adjoint and \(\nu=0\), or
        \item \(\nu\in W\lambda\) for some highest weight \(\lambda\) in the
            support of \(V\). In particular, it holds if every \(\lambda\) in
            the support of \(V\) is minuscule.
    \end{enumerate}
\end{corollary}
\begin{proof}
    It suffices to prove when \(V=V_\lambda\) and \(\lambda\) is fixed by
    \(\Frob_v\), because the general case follows from this case by inducing
    from some power of \(\Frob_v\).
    We shall fix the specific unramified \(\chi\)-datum as in the
    discussion before \Cref{thm:asymptotic_FL} and if the lift
    \(\gamma_G^\lambda\in G(F_v)\) exists, the result does not depend on such
    choice.

    If \(G\) is adjoint, the dual group \(\dual{\bG}\) is simply-connected, so
    the centralizer of \(\kappa\) in \(\dual{\bG}\) is connected. Consequently,
    if \(G\) is also split, then \(H\) must also be split. By
    \Cref{cor:kappa_twisted_Frobenius_for_central_coweight}, both
    \(\bar{\FRV}[\gamma]_\kappa\) and \(\FRV^H[\gamma_H]_\hST\) are trivial
    \(\Frob'\)-modules, and they have the same dimension by definition. Since we
    also have \(\Pair{\bFa}{\gamma}_\lambda=\nu(\bFa(\Frob_v))=1\) because
    \(\nu=0\), the theorem follows in this case.

    For general \(G\) with \(\nu\in W\lambda\), there is a unique
    \(\dual{\bH}\)-highest weight \(\lambda_H\) such that
    \(V^H_{\lambda_H}[\gamma_H]\) is non-trivial. Both
    \(\bar{\FRV}[\gamma]_\kappa\) and \(\FRV^H_{\lambda_H}[\gamma_H]_\hST\)  are
    \(1\)-dimensional and carry trivial \(\Frob'\)-actions.
    Recall we also have \(\tilde{\Frob}_v'=\xi_T(\Frob')\) induced by our choice
    of \(\chi\)-datum, and it acts on the \(1\)-dimensional weight space
    \(V_\lambda[\nu]\) via a character \(\phi\), and so we may identify
    \(V_\lambda[\nu]\) with \(\tilde{\Frob}_v'\)-action with
    \(\FRV_\lambda[\nu]_\kappa^\xi\otimes \phi\). Similarly, we have the
    character \(\phi_H\) for \(\tilde{\Frob}_H'=\xi_{T_H}(\Frob')\) on
    \(V_H^{\lambda_H}\). Lastly, the multiplicity space
    \(\Hom_{\dual{\bH}}(V_H^{\lambda_H},V)\) is \(1\)-dimensional and
    \(\Frob_H\) (hence \(\tilde{\Frob}_H'\)) acts via a character \(\psi\).
    By definition of \(\bFa\), we have
    \(\Pair{\bFa}{\gamma}_\lambda\phi=\phi_H\psi\), so we are left to prove that
    \(\phi=\phi_H\).

    The proof of this last equality is essentially just a toy case of
    Harish-Chandra descent for Levi subgroups. Indeed, if \(\nu\) is
    central, then so is \(\lambda\) and thus \(V_\lambda\) is a character of
    \(\LD{G}\), hence trivial on \(\Frob_v\). The construction of \(\xi_T\)
    shows that \(\tilde{\Frob}_v'\) is in fact contained in
    \(\dual{\bG}^\Der\rtimes\Frob_v\), which implies that \(\phi=1\) and
    similarly the same is true for \(\phi_H\). If \(\nu\) is not central, then
    \(\gamma\) is contained in the centralizer of some split torus inside
    \(G^\AD\), which induces a Levi subgroup \(L\subset G\). It is well-known that \(L\) is
    \(G(F_v)\)-conjugate to a standard Levi subgroup, and thus we have a
    canonical \(\dual{\bG}\)-conjugacy class of \(L\)-embedding \(\LD{L}\to
    \LD{G}\). It is also not hard to see that up to \(\dual{\bG}\)-conjugacy the
    composition
    \begin{align}
        \LD{T}\longto\LD{L}\longto\LD{G}
    \end{align}
    is exactly \(\xi_T\). Choose
    \(L\) to be minimal among the Levi subgroups containing \(\gamma\), then
    after potential conjugation by \(\dual{\bG}\) we may assume that
    \(\nu\) is a central character of \(\LD{L}\) and let \(V^L_\nu\) be the
    normalized \(\LD{L}\)-character (in other words, \(\Frob_v\) acts
    trivially). Restricting \(V_\lambda\) to \(\LD{L}\), then there
    exists a unique subrepresentation of \(\LD{L}\) in \(V_\lambda\)
    containing \(\nu\), which also has to be a character since \(\nu\) is
    extremal. This is exactly
    \begin{align}
        V_\nu^L\otimes\Hom_{\dual{L}}(V_\nu^L,V_\lambda).
    \end{align}
    Since \(\xi_{T_L}(\Frob_v)\) is contained in
    \(\dual{\bL}^\Der\rtimes\Frob_v\), its action on \(V_\nu^L\) is
    trivial, so we necessarily have (we abuse notation and regard character and
    \(1\)-dimensional module as the same thing)
    \begin{align}
        \phi=\Hom_{\dual{L}}(V_\nu^L,V_\lambda).
    \end{align}
    However, we can use \(\dual{\bG}\)-conjugation again and view \(\LD{L}\) as
    a \emph{standard} Levi subgroup of \(\LD{G}\), and \(V_\nu^L\) is
    replaced by some \(V_{\nu'}^L\) for a \(\Frob_v\)-fixed (the action induced
    by \(\LD{G}\)) weight \(\nu'\). This shows that \(\phi=\phi_H=1\) and
    finishes the proof.
\end{proof}

\subsection{}
Using the same Harish-Chandra descent argument as in
\Cref{cor:special_case_of_unramified_GASF_irreducible_components_Frob_transfer},
one can show that
\Cref{thm:asymptotic_FL} follows from
the special case where \(\nu=0\). Although we have not been able to find a
simple proof of this special case for general \(G\), we point out that it
becomes a purely combinatorial statement: indeed, one immediately sees
that both \(\Delta_{\symup{II}}(\gamma_H,\gamma)\) and
\(\Pair{\bFa}{\gamma}\) are trivial.

The multiplicity spaces
\(\Hom_{\dual{\bH}}(V_{\lambda_H}^H,V)\) are slightly trickier,
but it is not hard to see that if \(\Hom_{\dual{\bH}}(V_{\lambda_H}^H,V)\) is
non-trivial and \(V_{\lambda_H}^H\) contains the zero weight, then such
multiplicity space
does not depend on admissible embedding \(\xi\) (because the ambiguity lies in
\(Z_H\)) but only on the twists given
by \(\Out(\bG)\)-torsors and
\(\Out(\bH)\)-torsors, which are also combinatorial in nature. As a result,
it is independent of the field \(k\).

\subsection{}
Lastly, we prove four more special cases of
\Cref{thm:asymptotic_FL}, all of which
necessary to finish the proof of fundamental lemma.

\begin{proposition}
    \label[proposition]{prop:special_case_of_three_little_pigs}
    \Cref{thm:asymptotic_FL} holds when
    \(\nu=0\) and \(\gamma\) is \(\nu\)-regular semisimple in the following cases:
    \begin{enumerate}
        \item \(\dual{\bG}\) is simple of type \(\TypeB_r\) (\(r\ge 2\)) and
            \(V\) is the standard representation of \(\SO_{2r+1}\).
        \item \(\dual{\bG}\) is simple of types \(\TypeE_6\) (including the case
            where \(G\) is of type \(\prescript{2}{}{\TypeE}_6\)) or \(\TypeE_7\),
            and \(V\) is the adjoint representation.
    \end{enumerate}
\end{proposition}
\begin{proof}
    We prove this by hand. It is harmless to assume that \(\dual{\bG}\) is
    adjoint in all three cases, as we have seen that the statement is purely
    combinatorial for \(\nu=0\). We shall always call a root of \(\dual{\bG}\) a
    \emph{coroot}, but may use either \emph{weight} (in \(V\)) or
    \emph{coweight} (of \(G\)) for a character of \(\dual{\bT}\), and
    hopefully the context will clarify any potential confusion.

    Suppose \(\dual{\bG}\) is of type \(\TypeB_r\), and so it must be
    Frobenius-split. The
    weights in the standard representation are all in the same \(W\)-orbit,
    except for \(\nu=0\). The \(0\)-weight space is \(1\)-dimensional and
    spanned by the ``middle'' vector
    \begin{align}
        \FRx=(0,\ldots,0,1,0,\ldots,0),
    \end{align}
    where \(1\) sits in the \((r+1)\)-coordinate (assuming we use the
    anti-diagonal quadratic form to split \(\SO_{2r+1}\)). The element
    \(\rho_\FRx\) is simply the last fundamental coweight \(\CoWt_r\).
    Since \(2\CoWt_r\) is contained in the coroot lattice, the value of
    \(\kappa\) at
    \begin{align}
        \Frob_v(\rho_\FRx)-\Frob'(\rho_\FRx)=
        \CoWt_r-\Frob'(\CoWt_r)
    \end{align}
    must be either \(1\) or \(-1\). Call this sign \(\epsilon\). On the other
    hand, there is a unique \(\lambda_H\) such that
    \(\FRV_{\lambda}^H[\gamma_H]_\hST^\xi\neq 0\), and so it is a
    \(1\)-dimensional trivial \(\Frob'\)-module. Therefore, it amounts to showing
    that
    \begin{align}
        \Hom_{\dual{\bH}}(V_{\lambda_H}^H,V)=\Qlb\otimes\epsilon.
    \end{align}
    Note that \(\Frob'\) is now just an element \(w\in W\),
    but it is unnecessary to enumerate all possible \(w\)'s.
    Rather, by modifying \(w\) by an element of \(W_H\), and use the fact
    that \(\kappa\) is trivial on the coroots of \(H\), we may simply use the
    element \(w_H^G\) corresponding to \(\Frob_H\) instead. Since the
    \(\Frob'\)-action on \(\Hom_{\dual{\bH}}(V_{\lambda_H}^H,V)\) factors
    through \(\Frob_H\) as well, we reduce to proving that the
    \(\Frob_H\)-action on \(0\)-weight space \(V[0]\) is \(\epsilon\) (the
    sign characters from cohomological shift are trivial since both \(V\) and
    \(V_{\lambda_H}^H\) contains \(0\)-weight).
    As such, the result is clearly true if \(\dual{\bH}\) is also split, so we
    assume that \(\Frob_H\neq 1\).

    Since \(\dual{\bH}\) is the connected centralizer of a semisimple element
    in \(\dual{\bG}\), its Dynkin diagram is embedded into the \emph{completed}
    Dynkin diagram of \(\dual{\bG}\). The completed Dynkin diagram
    of \(\dual{\bG}\) has a unique automorphism of order \(2\), exchanging
    \(\CoRt_1\) with the lowest negative coroot \(\CoRt_0\) while keeping other
    simple coroots fixed. It is represented by the reflection
    \begin{align}
        w_H^G=s_{\CoRt_1+\cdots+\CoRt_r},
    \end{align}
    and the difference between \(\CoWt_r\) and
    \(\Frob_H(\CoWt_r)\) is \(\CoRt=\CoRt_1+\cdots+\CoRt_r\). Since \(\Frob_H\)
    stabilizes the positive coroots of \(H\) but inverts the coroot
    \(\CoRt\), \(\CoRt\) cannot be a coroot of \(H\), and so
    \(\epsilon=\kappa(\CoRt)=-1\). On the other hand, it is easy to compute that
    any lift of \(w_H^G\) to \(\SO_{2r+1}\) maps \(\FRx\) to \(-\FRx\) inside
    \(V\). This proves the case for \(\TypeB_r\).

    Suppose \(\dual{\bG}\) is of type \(\TypeE_6\) and \(V\) is the adjoint
    representation. There are \(6\) crystal basis elements
    \(\FRx_1,\ldots,\FRx_6\) with weight \(0\), and we have
    \(\rho_{\FRx_i}=\CoWt_i\). Similar to \(\TypeB_r\)-case, we need to prove
    that the \(\Frob_H\)-action on the Cartan algebra \(V[0]\simeq\La{t}_{\Qlb}\)
    is isomorphic to the module \(\bigoplus_{i=1}^6\Qlb\FRx_i\) with
    \(\Frob_H\)-action
    \begin{align}
        \FRx_i\longmapsto
        \kappa\bigl(\Frob_v(\CoWt_i)-\Frob_H(\CoWt_i)\bigr)\FRx_{\Frob_v(i)}.
    \end{align}
    The \(H\)-Frobenius \(\Frob_H\) is represented by an automorphism
    of the completed Dynkin diagram of \(\TypeE_6\), but not every automorphism
    can appear.

    If \(\dual{\bG}\) is \(\Frob_v\)-split, so that \(\Frob_v=1\), then
    \(\Frob_H\) is either trivial or a rotation of the completed Dynkin diagram
    of order \(3\). If \(\Frob_H=1\), the conclusion is trivial. If \(\Frob_H\)
    has order \(3\), then the characteristic polynomial (in variable \(x\)) of
    \(\Frob_H\) on \(\La{t}\) is \((x^3-1)^2\). On the other hand, since
    \(\CoWt_2\) is contained in the coroot lattice, we have
    \begin{align}
        \kappa\bigl(\Frob_v(\CoWt_2)-\Frob_H(\CoWt_2)\bigr)
        =\kappa(\CoWt_2)-\kappa(\Frob_H(\CoWt_2))
        =1.
    \end{align}
    Similarly, we have
    \begin{gather}
        \kappa\bigl(\Frob_v(\CoWt_2)-\Frob_H(\CoWt_2)\bigr)=1,\\
        \kappa\bigl(\Frob_v(\CoWt_1+\CoWt_3)-\Frob_H(\CoWt_1+\CoWt_3)\bigr)=
        \kappa\bigl(\Frob_v(\CoWt_1+\CoWt_6)-\Frob_H(\CoWt_1+\CoWt_6)\bigr)=1,\\
        \kappa\bigl(\Frob_v(\CoWt_6+\CoWt_5)-\Frob_H(\CoWt_6+\CoWt_5)\bigr)=1.
    \end{gather}
    So in this case we only need to prove that
    \(\kappa\bigl(\Frob_v(\CoWt_1)-\Frob_H(\CoWt_1)\bigr)\) is a primitive
     root of unity of order \(3\). Suppose it is not true, then
     \(\kappa\bigl(\Frob_v(\CoWt_i)-\Frob_H(\CoWt_i)\bigr)=1\) for every
     \(i\), which means that there exists a \(\Frob_H\)-fixed lift of \(\kappa\) to
     \(\dual{\bT}^\SC\), or in other words, \(H\) induces an endoscopic group of
     \(G^\AD\). But this is impossible because endoscopic groups of split
    adjoint group \(G^\AD\) are split, while \(H\) is not.

    If \(\Frob_v\) is non-trivial, then it has order \(2\), and
    \(\Frob_H\) must have order \(2\) as well. The characteristic polynomial of
    \(\Frob_H\) on \(\La{t}\) in this case is \((x^2-1)^2(x-1)^2\). On the other
    hand, since \(3\CoWt_i\) is contained in the coroot lattice for any \(i\),
    we have
    \begin{align}
        \label{eqn:kappa_on_Frob_rho_i_minus_Frob_H_rho_i}
        \kappa\bigl(\Frob_v(\CoWt_i)-\Frob_H(\CoWt_i)\bigr) =1,
    \end{align}
    because the left-hand side has order dividing \(2\).
    The conclusion is then readily seen since \(\Frob_H\) is a conjugate of
    \(\Frob_v\) as automorphisms of the completed Dynkin diagram.

    Lastly, we consider \(\TypeE_7\). The group \(G\) is again split and
    \(\Frob_v=1\). If \(\Frob_H=1\), then the conclusion is again trivial to
    see. So we consider when \(\Frob_H\) is the unique involution of the
    completed Dynkin diagram. The \(\Frob_H\)-action has characteristic
    polynomial \((x^2-1)^3(x-1)\). The coweights \(\CoWt_1\), \(\CoWt_3\),
    \(\CoWt_4\) and \(\CoWt_6\) are contained in the coroot lattice, so
    \eqref{eqn:kappa_on_Frob_rho_i_minus_Frob_H_rho_i} holds for \(i=1,3,4,6\).
    Moreover, any two of \(\CoWt_2\), \(\CoWt_5\) and \(\CoWt_7\) sum up to some
    element in the coroot lattice, so if
    \eqref{eqn:kappa_on_Frob_rho_i_minus_Frob_H_rho_i} holds for any of the
    three, it holds for all three, in which case \(H\) again induces a non-split
    endoscopic group of the split adjoint group \(G^\AD\), a contradiction.
    This finishes the proof.
\end{proof}

\begin{proposition}
    \label[proposition]{prop:special_case_of_big_bad_wolf}
    \Cref{thm:asymptotic_FL} holds when \(\dual{\bG}\) is simple of type
    \(\TypeE_7\), \(V\) is the second fundamental representation, and \(\gamma\)
    is \(\nu\)-regular semisimple.
\end{proposition}
\begin{proof}
    We may assume that \(G=G^\AD\) and \(\dual{\bG}\) is
    simply-connected. Then both \(G\) and \(H\) are necessarily split, and so
    \(\Frob_v=\Frob_H=1\). The \(L\)-embedding \(\xi\) is then canonically fixed.

    The weights of the second fundamental representation fall into two Weyl
    orbits, represented by \(\CoWt_2\) and \(\CoWt_7\) respectively. The orbit
    of \(\CoWt_2\) is already covered by
    \Cref{cor:special_case_of_unramified_GASF_irreducible_components_Frob_transfer},
    and so we may assume that \(\nu\in W\CoWt_7\). Replacing the endoscopic
    datum by a \(\dual{\bG}\)-conjugate, we may assume that \(\nu=\CoWt_7\).

    The crystal data for \(V\) can be explicitly computed. The
    \(\CoWt_7\)-weight space is \(6\)-dimensional, with basis denoted by
    \(\FRx_1,\ldots,\FRx_6\), and we have
    \begin{align}
        \epsilon_i(\FRx_j)
        &= \begin{cases}
            1 & 1\le i=j\le 6,\\
            0 & \text{otherwise}.
        \end{cases}
    \end{align}
    The stabilzer of \(\CoWt_7\) in \(W\) is exactly the subgroup generated by
    simple reflections \(s_1,\ldots,s_6\). As a result, if
    \(\Frob'=w\rtimes\Frob_v=w\in W\), then with the notations in
    \Cref{prop:kappa_twisted_Frobenius_and_n_Rt_x}, \(r_{j-1}^{-1}\) acts trivially
    on the set of \(\FRx_i\) for any \(1\le j\le m\) because \(r_{j-1}^{-1}\)
    fixes \(\CoWt_7\). By the same argument in
    \Cref{cor:kappa_twisted_Frobenius_for_central_coweight}, the Frobenius
    module structure on \(\bar{\FRV}[\gamma]_\kappa\) is given by
    \begin{align}
        \FRx_i&\longmapsto
        \kappa\bigl(\Frob_v(\CoWt_i)-\Frob'(\CoWt_i)\bigr)\FRx_i.
    \end{align}
    But since \(H\) is also split, \(\Frob'=w\) must lie in the Weyl group of
    \(H\), and so \(\kappa\bigl(\Frob_v(\CoWt_i)-\Frob'(\CoWt_i)\bigr)=1\), or
    in other words, \(\bar{\FRV}[\gamma]_\kappa\) is a trivial Frobenius
    module. Since \(\gamma\) is \(\nu\)-regular, we have
    \(\Delta_{\symup{II}}(\gamma_H,\gamma)=1\) as discussed before stating
    \Cref{thm:asymptotic_FL}, as well as \(\Pair{\bFa}{\gamma}_\lambda=1\)
    because \(\Frob_v=1\). Combining these together, we see that
    \(\FRV[\gamma]_\kappa^\xi\) is again a trivial
    Frobenius module of dimension \(6\).

    Finally, it is easy to see that the module structure on
    \(\FRV^H[\gamma_H]_{\hST}^\xi\) is also trivial: for
    \(\FRV_{\lambda_H}^H[\gamma_H]_{\hST}\) in
    the decomposition it is because \(\nu\) is fixed by \(w\), and for the multiplicity
    space \(\Hom_{\dual{\bH}}(V_{\lambda_H}^H,V)\) it is due to that the
    dual groups are both split. Therefore, \(\FRV[\gamma]_\kappa^\xi\cong
    \FRV^H[\gamma_H]_{\hST}^\xi\) as desired.
\end{proof}


\section{Fundamental Lemma in the Language of Monoids} 
\label{sec:matching_orbits}

In this section, we demonstrate how to translate the fundamental lemma from
group setting into the language of monoids.

\subsection{}
Given conjugacy class \(\gamma\in G(F_v)\), its image \(\gamma_\AD\) in
\(G^\AD(F_v)\) lifts to \(\Env(G^\SC)(F_v)\). For any monoid
\(\FRM\in\FM(G^\SC)\), we have natural map
\begin{align}
    \Stack*{\FRM/Z_\FRM}\longto \Stack*{\Env(G^\SC)/Z_{\Env(G^\SC)}},
\end{align}
which identifies their invertible part
\begin{align}
    \FRM^\x/Z_\FRM\simeq \Env(G^\SC)^\x/Z_{\Env(G^\SC)}\simeq G^\AD,
\end{align}
so we may view \(\gamma_\AD\) as elements in \(\FRM^\x/Z_\FRM(F_v)\) for any
\(\FRM\). Let \(\bM\) be the split form of \(\FRM\) and for simplicity
suppose \(\bM=\bM(\lambda_1,\ldots,\lambda_m)\), where \(\lambda_i\) are
dominant cocharacters of \(\bG\) (the general case does not add
anything essential, only technical nuisance). Suppose \(\nu\) is the Newton
point of \(\gamma\), then \(\gamma_\AD\) extends to an \(\cO_v\)-point in
\(\Stack*{\FRC_\FRM/Z_\FRM}\) if and only if \(\nu_\AD\in\WP(-w_0(\lambda_\AD))\)
for some \(\lambda_\AD\) in the integral cone generated by \(\lambda_{i,\AD}\).
In fact, in this case it extends to a point in \(\FRC_\FRM(\cO_v)\) since
\(F_v\) is local.

\subsection{}
Let \(H\) be an endoscopic group of \(G\) associated with endoscopic datum
\((\kappa,\OGT_\kappa)\) (ignoring the choice of \(L\)-embedding \(\xi\) for the
moment), and let \(\FRM_H\) be the endoscopic
monoid corresponding to \(\FRM\). Recall we also have a (not necessarily flat)
monoid \(\FRM_H'\) and a natural map
\begin{align}
    \Stack*{\FRM_H/Z_\FRM^\kappa}\longto\Stack*{\FRM_H'/Z_\FRM},
\end{align}
identifying their invertible part
\begin{align}
    \FRM_H^\x/Z_\FRM^\kappa\longto(\FRM_H')^\x/Z_\FRM.
\end{align}
We also have canonical embedding \(Z_G\to H\), and thus 
for any \(\gamma_H\in H(F_v)\), its image
\(\bar{\gamma}_H\in H/Z_G(F_v)\) can also be identified with a point in
\(\FRM_H^\x/Z_\FRM^\kappa(F_v)\).

Suppose \(\gamma_H\) matches \(\gamma\), then their images modulo \(Z_G\) match
as well, and it is compatible with the canonical map
\begin{align}
    \Stack*{\FRC_{\FRM,H}/Z_\FRM^\kappa}\longto\Stack*{\FRC_\FRM/Z_\FRM}.
\end{align}
By the construction of \(\FRM_H\),
\(\bar{\gamma}_H\) extends to a point in
\(\Stack*{\FRC_{\FRM,H}/Z_\FRM^\kappa}(\cO_v)\)
if and only if \(\gamma_\AD\) extends to a point in
\(\Stack*{\FRC_\FRM/Z_\FRM}(\cO_v)\). Moreover, if \(\gamma\) is strongly
regular semisimple, then its image in \(\Stack*{\FRC_\FRM/Z_\FRM}\) is
necessarily generically regular semisimple.

This finishes the matching of conjugacy classes, but one may notice that the
monoid \(\FRM\) is still arbitrary at this point, because the choice of monoid
is not tied to conjugacy classes but to the Satake functions.

\subsection{}
Suppose now we choose an algebraic representation \(V\) for \(\LD{G}\), whose
corresponding Satake function is \(\fS^V\). Let \(\FRM\) be the monoid whose
split form is \(\bM=\bM(-w_0(\lambda_1),\ldots,-w_0(\lambda_m))\) where
\(\Set*{\lambda_i}\) is the support of \(V\). Such monoid is canonically defined
because the Galois action on \(\Set*{\lambda_i}\) is canonical. Even though the
spherical Hecke algebra \(\cH_{G,0}\) of \(G\) is generated by those
\(\fS^\lambda=\fS^{V_\lambda}\) where \(\lambda\) is \(F_v\)-rational, we would
nevertheless want to consider the case of a general \(V\) because that allows
more flexibility in geometric constructions. For example, it would allow us to
include minuscule \(\lambda_i\)'s even when there is no \(F_v\)-rational
minuscule coweights (e.g., when \(G\) is the unramified unitary group of rank
\(3\)).

Note that the monoid does not see the multiplicity spaces at all, because the
latter will show up as local systems on the relevant affine Schubert
varieties determined by those \(\lambda_i\)'s, which we are free to
choose. In other words, if \(\FRM\) is given as above, it allows us to encode
information about any \(\fS^V\) as long as its support is contained in the set
of \(\lambda_i\)'s.

\subsection{}
The statement of fundamental lemma involves orbital integrals that are twisted
sums of the forms
\begin{align}
    \OI_\gamma(\fS^V,\dd t_v)=\int_{I_\gamma(F_v)\backslash
    G(F_v)}\fS^V(g_{v}^{-1}\gamma g_{v})\frac{\dd g_v}{\dd t_v}.
\end{align}
By \Cref{prop:KV_nonempty_criteria}, this integral is non-trivial only if
\(\gamma\) extends to a point
\begin{align}
    a\in\Stack*{\FRC_\FRM/Z_\FRM}(\cO_v).
\end{align}
In this case, if \(\gamma\) matches \(\gamma_H\), then \(\gamma_H\) necessarily extends
to a point
\begin{align}
    a_H\in\Stack*{\FRC_{\FRM,H}/Z_\FRM^\kappa}(\cO_v).
\end{align}
Note that since \(V\) is not necessarily irreducible, \(a\) is not necessarily
unique, and even it is, there can still be multiple \(a_H\)'s. However, if
viewed as \(F_v\)-points then both \(a\) and \(a_H\) are unique.

The projection \(G\to G^\AD\) induces inclusion \(\LD{(G^\AD)}\to\LD{G}\),
and the restriction of \(V\) induces function \(\fS_\AD^V\in\cH_{G^\AD,0}\).
Observing the definition of \(\OI_\gamma\), we see that there is no need to
treat the \emph{group} \(G\) and the \emph{\(G\)-space} \(G\) as the
same object. Indeed, we can define a modified orbital integral
\(\OI_{\gamma_\AD,G}\) as follows:
\begin{align}
    \OI_{\gamma_\AD,G}(\fS_\AD^V,\dd t_v)=\int_{I_\gamma(F_v)\backslash
    G(F_v)}\fS_\AD^V(g_{v}^{-1}\gamma_\AD g_{v})\frac{\dd g_v}{\dd t_v}.
\end{align}
The important thing to note is that \(g_v\) still ranges over \(G(F_v)\), not
\(G^\AD(F_v)\). Note also that \(\gamma_\AD\) is not necessarily strongly
regular semisimple (in \(G^\AD\)). It is not necessarily true that we have the
equality
\begin{align}
    \label{eqn:compare_OI_to_OIAD}
    \OI_{\gamma_\AD,G}(\fS_\AD^V,\dd t_v)=\OI_{\gamma}(\fS^V,\dd t_v),
\end{align}
however, by 
\Cref{prop:KV_nonempty_criteria}, this equality holds if and only if \(V\) satisfies
the following condition: for any \(\lambda\) in the support of \(V\), we
have \(\nu_\gamma\le_\bbQ\lambda\) if and only if
\(\nu_{\gamma,\AD}\le_\bbQ\lambda_\AD\). In particular, it is true if
the support lies in the same \(\bbQ\CoRoots\)-coset as \(\nu_\gamma\).

More generally, if the center \(Z_G\) of \(G\) is connected, we have by Lang's
theorem
\begin{align}
    \RH^1(X_v,Z_G)=0.
\end{align}
Therefore, if \(\nu_{\gamma,\AD}\le_\bbQ\lambda_\AD\), we can always
lift \(\gamma_\AD\) to a point \(\gamma_G^\lambda\in G(F_v)\) such that
\(\nu_{\gamma_G^\lambda}\le_\bbQ\lambda\). By breaking the support of \(V\)
into \(\bbQ\CoRoots\)-equivalent subsets, each represented by some \(\lambda\), we may write
\begin{align}
    \OI_{\gamma_\AD,G}(\fS_\AD^V,\dd
    t_v)=\sum_\lambda\OI_{\gamma_G^\lambda}(\fS^V,\dd t_v),
\end{align}
where the summation ranges over all representatives \(\lambda\). For
convenience, we may also assume that \(\gamma\) is one of those
\(\gamma_G^\lambda\)'s.

Choose an admissible \(L\)-embedding \(\xi\) of \(\LD{H}\) into \(\LD{G}\),
then we have the transferred function \(\fS_{H,\xi}^V\) by restricting \(V\) via
\(\xi\). We then have orbital integral \(\OI_{\gamma_H}(\fS_{H,\xi}^V,\dd t_v)\).
Similarly, using \(H\to \bar{H}=H/Z_G\), we have the modified version
\(\OI_{\bar{\gamma}_H,H}(\fS_{\bar{H},\xi}^V,\dd t_v)\). Note, however, that there
is in general no embedding from the \(L\)-group of \(\bar{H}\) to that of
\(G^\AD\) that also fixes \(\kappa\).
We necessarily have
\begin{align}
    \OI_{\gamma_H}(\fS_{H,\xi}^V,\dd t_v)=\OI_{\bar{\gamma}_H,H}(\fS_{\bar{H}, \xi}^V,\dd
    t_v)
\end{align}
whenever \eqref{eqn:compare_OI_to_OIAD} holds.

\subsection{}
Since \(\gamma_\AD\) lifts to an \(F_v\)-point in \(\Env(G^\SC)\), it will also
lift to an \(F_v\)-point in \(\Stack*{\FRM/Z_\FRM}\). Fix an arbitrary lift
\(\tilde{\gamma}\), then for any \(\gamma'\in
G(F_v)\) stably conjugate to \(\gamma\), we have a corresponding lift
\(\tilde{\gamma}'\) stably conjugate to \(\tilde{\gamma}\), and it is
straightforward to see that
\begin{align}
    \inv(\tilde{\gamma},\tilde{\gamma}')=\inv(\gamma,\gamma').
\end{align}
Conversely, any \(\tilde{\gamma}'\) stably conjugate to \(\tilde{\gamma}\)
naturally induces a unique \(\gamma'\) stably conjugate to \(\gamma\).
Therefore, the modified \(\kappa\)-orbital may be defined as
\begin{align}
    \OI_{a,\AD}^\kappa(\fS_\AD^V,\dd
    t_v)=\sum_{\gamma'}\Pair{\inv(\tilde{\gamma},\tilde{\gamma}')}{\kappa}\OI_{\gamma_\AD,G}(f_\AD^V,\dd
        t_v),
\end{align}
which now depends only on \(a\), \(x_a\in \OGT_G\x\bar{\bT}_\bM(F_v^\sep)\)
lying over \(a\) and the choice of \(\tilde{\gamma}\) (within its
stable conjugacy class). If the \(\Out(\bG)\)-torsor \(\OGT_G\) is pointed (by
\(x_G\)), we also require that \(x_a\) lies over \(x_G\), and if the endoscopic
datum is pointed as well, we will choose the same point
\(x_a=x_{a_H}\) for both \(G\) and \(H\), lying over \(\OGT_\kappa^\bullet\).

If \eqref{eqn:compare_OI_to_OIAD} holds, or if \(Z_G\) is connected, we may
reduce the fundamental lemma to the following equality:
\begin{align}
    \Delta_0(\gamma_H,\gamma)\OI_{a,\AD}^\kappa(\fS_\AD^{V},\dd
    t_v)=
    \SOI_{\bar{\gamma}_H,H}(\fS_{\bar{H}, \xi}^V,\dd t_v).
\end{align}
Note that although the above equality is independent of the choice of \(x_a\),
the transfer factor \(\Delta_0\) is. Moreover, by the definition of transfer
factor, the product \(\Delta_{\symup{I}}\Delta_{\symup{II}}\Delta_{\symup{IV}}\)
depends (aside from \(x_a\)) only on \(\gamma_\AD\), so among the factors of
\(\Delta_0\), only \(\Delta_{\symup{III}_2}\) depends on \(\gamma\) on the nose.

\chapter{Global Constructions}%
\label{chap:global_constructions}

In this chapter, we make several useful global constructions. The first is the
moduli of boundary divisors, which is an indispensable tool
for describing multiplicative Hitchin fibrations. Technically, the moduli of
boundary divisors in this chapter is slightly different from what we will use
for future chapters: here we formulate them as a space due to its
simplicity in language, but in \Cref{chap:multiplicative_hitchin_fibrations} we
will introduce a slightly modified version where the moduli functor is a
Deligne--Mumford stack instead of a space.

Another important construction is the
global affine Schubert schemes and its variants. Instead of using affine
Grassmannians, we use an alternative formulation using reductive monoids. We
will also state a natural factorization property that is especially useful in
\Cref{chap:deformation}.

\section{Global Affine Grassmannian}
\label{sec:review_BD_affine_grassmannian}
In this section, we review the standard theory of Beilinson--Drinfeld
Grassmannians. The reference of this section is \cite{Zh17} (especially its
\S~3), where proofs or further references can be found.
\subsection{}
Let \(G\) be a reductive group over a smooth curve \(X\)
    \nomenclature[\(X \)]{\(X\)}{a smooth curve over \(k\), for the most part assumed to be projective and geometrically connected}
induced by \(\bG\) and a \(\Out(\bG)\)-torsor
\(\OGT_G\). For any positive integer \(d\), we may define the so-called
Beilinson--Drinfeld Grassmannian over \(X_d=\Sym^dX\)
    \nomenclature[\(X_d \)]{\(X_d\)}{the \(d\)-th symmetric power of the curve \(X\)}
as the follows: any
\(S\)-point \(D\in X_d(S)\) may be interpreted as a finite flat \(S\)-family of
divisors of degree \(d\) in \(X\). The Beilinson--Drinfeld Grassmannian \(\Gr_{G,d}\)
    \nomenclature[\(Gr_G_d \)]{\(\Gr_{G,d},\Gr_{G,\cB}\)}{the Beilinson--Drinfeld
    Grassmannian over \(X_d\), or more generally any moduli stack \(\cB\) of flat families of divisors on \(X\)}
sends \(S\) to the groupoid
\begin{align}
    \Set*{\,(D,E,\phi)\,\given \,
        \begin{gathered}
            D\in X_d(S), E\in \Bun_G(S)\\
            \phi\colon E|_{X\x S-D}
            \stackrel{\sim}{\to} E_0|_{X\x S-D}\text{ is a trivialization}
        \end{gathered}
        \,}.
        \nomenclature[\(B{}un_G \)]{\(\Bun_G\)}{the classifying stack of \(G\)-bundles over \(X\)}
\end{align}
This is known to be an ind-scheme of ind-finite-type over \(X_d\). Since \(G\)
is reductive, it is also ind-projective. Similarly, we have the \(X_d\)-family
of jet groups and arc groups defined as follows: for \(D\in X_d(S)\), let
\(I_D\) be the ideal in \(\cO_{X\x S}\) defining divisor \(D\), then we have
infinitesimal neighborhoods \(D_n\) defined by \(I_D^n\) as well as the formal
completion \(\hat{X}_D\) of \(X\x S\) at \(D\). Then the jet groups and arc
group are defined as
\begin{align}
    \Arc_{X_d,n}{G}(S)&=\Set{(D,g)\given D\in X_d(S), g\in G(D_n)},\\
    \Arc_{X_d}{G}(S)&=\Set{(D,g)\given D\in X_d(S), g\in G(\hat{X}_D)}.
    \nomenclature[\(L{}Y"bbold_B^+ \)]{\(\Arc_{\cB}{Y}\)}{the global family of arc
    spaces of \(Y\) parametrized by a moduli \(\cB\) of divisors}
\end{align}
They are known to be representable by an affine scheme over \(X_d\), and the jet
groups are of finite type. In
addition, using the fact that \(G\) is smooth over \(X\) and the 
infinitesimal lifting criterion, we see that the jet group
schemes are smooth over \(X_d\), while the arc group is formally smooth over
\(X_d\).

\subsection{}
The definition of loop groups is trickier. First, one can show that the formal
scheme \(\hat{X}_D\) is ind-affine relative to \(S\). Without loss of
generality let \(S\) be affine. 
Taking its ring \(R_D\) of global
functions and let \(\hat{X}'_D=\Spec{R_D}\), then there is a canonical map
\(\hat{X}_D\to\hat{X}'_D\) through which the map \(\hat{X}_D\to X\x S\)
uniquely factors. Therefore, it makes sense to define
scheme \(\hat{X}^\bullet_D=\hat{X}'_D-D\). Then the loop group is defined
as
\begin{align}
    \Loop_{X_d}{G}(S)&=\Set{(D,g)\given D\in X_d(S), g\in G(\hat{X}^\bullet_D)}.
    \nomenclature[\(L{}Y"bbold_B_ \)]{\(\Loop_{\cB}{Y}\)}{the global family of loop
    spaces of \(Y\) parametrized by a moduli \(\cB\) of divisors}
\end{align}
This functor is represented by an ind-scheme over \(X_d\), and it is well-known
(see \cite{Zh17}*{\S~3})
that we have a natural isomorphism of \(k\)-spaces
\begin{align}
    \Gr_{G,d}\simeq \Stack{\Loop_{X_d}{G}/\Arc_{X_d}{G}}
\end{align}
using Beauville--Laszlo's descent theorem.

\subsection{}
If we have curves \(X^\lambda\)
parametrized by \(\lambda\) in a finite set \(\Lambda\) and a tuple of 
positive integers \(\ul{d}=(d_\lambda)_{\lambda\in\Lambda}\), then we can let
\(X=\coprod_\lambda X^\lambda\), and we
will have natural map
\begin{align}
    X^{\Lambda}_{\ul{d}}
    \defeq\prod_\lambda X^\lambda_{d_\lambda}
    \longto X_{\abs{\ul{d}}},
\end{align}
where \(\abs{\ul{d}}=\sum_\lambda d_\lambda\).
Therefore, we also have the affine Grassmannian, arc group, etc. over
\(X^{\Lambda}_{\ul{d}}\).

\section{Boundary Divisors}
\label{sec:boundary_divisors}

Let \(X\) be a smooth, geometrically connected, projective curve over \(k\).
Consider the mapping stack from \(X\) to the quotient stack
\(\Stack{\bbA^1/\Gm}\). It contains an open substack where the generic point of
\(X\) is mapped into the invertible locus \(\Stack{\Gm/\Gm}\), which is
representable by a scheme and in fact is
precisely the moduli space of effective divisors in \(X\).

In \cite{BNS16} (mainly its \S~3), the authors generalizes such results to the
case where \(\bbA^1\) is replaced by a split affine normal toric variety
\(\FRA\) and \(\Gm\) by \(\FRA^\x\), and the resulting moduli space parametrizes
divisors valued in cocharacters in the cone \(\sC\) of \(\FRA\) (often called
\notion{colored divisors}\index{divisor!colored} in the literature). Further case where \(\FRA\) is a
split \(L\)-monoid (cf.~\Cref{sub:different_monoids}) is also considered in
\textit{loc. cit.}

This section generalizes the results in \cite{BNS16} to any affine normal
reductive monoid scheme \(\FRM\) that is locally constant over \(X\) and whose unit
group \(\FRM^\x\) is quasi-split. The proofs in
\textit{loc. cit.} can be applied here with necessary modifications, but
it can be a bit lengthy due to many reduction steps involved. Therefore,
following an anonymous referee's suggestion, we will adopt a different approach
towards a more abstract and uniform proof.

Aside from a representability
statement, the focus of this section will still be the case where \(\FRM=\FRA\)
is a toric scheme over \(X\), and we will continue the discussion for the
general case in \Cref{sec:global_affine_schubert_scheme}.

\subsection{}
We denote \(G=\FRM^\x\) and let \(T\subset B\) be a Borel pair defined over
\(X\) that is part of a pinning of \(G\). Let \(\FRT=\FRT_\FRM\) be the closure
of \(T\) in \(\FRM\), and it is an affine normal toric variety with unit torus
\(T\). Let \(\CoCharG(T)\) be the \emph{\'etale sheaf} of cocharacters defined
as the internal-hom functor in the category of abelian sheaves
\begin{align}
    \CoCharG(T)\defeq \IHom_X(\Gm,T).
\end{align}
The toric scheme \(\FRT\) is then defined by a sheaf of saturated strictly
convex cones \(\sC\subset \CoCharG(A)\).
If \(G\) is split, then \(\CoCharG(T)\) is constant and we may treat it as a
discrete set, and \(\lambda\in\sC\) if and only if
\(\lambda\colon\Gm\to T\) extends to a monoidal homomorphism from \(\bbA^1\) to
\(\FRT\). It is clear that both \(\CoCharG(T)\) and \(\sC\) are representable by
(possibly countably infinite) \'etale covers of \(X\).

Similarly, we also have the sheaf of
characters \(\CharG(T)\) and the cone \(\sC^*\subset\CharG(T)\) dual to
\(\sC\), both of which are representable by \'etale covers of \(X\). and we have
\begin{align}
    \FRT=\RSpec_{X}\bigoplus_{\lambda\in\pi_0(\sC^*)}\cO_{\sC^*_\lambda},
\end{align}
where \(\sC^*_\lambda\) is the connected component of \(\sC^*\) associated with
\(\lambda\in\pi_0(\sC^*)\).

\subsection{}
Consider the \(k\)-mapping stack
\begin{align}
    \OGASch[+]_X\defeq\IHom_{X/k}(X,\Stack{\FRM/G}).
\end{align}
For a \(k\)-scheme \(S\), an \(S\)-point of \(\OGASch[+]_X\) consists of a pair
\((E,\phi)\) where \(E\) is an \(G\)-torsor over \(X\x S\) and \(\phi\) is a
section of the induced monoidal bundle \(\FRM_E\) over \(X\x S\).
We are interested in the open substack
\begin{align}
    \OGASch_X\subset\OGASch[+]_X
    \nomenclature[\(Q_X_wow \)]{\(\OGASch_X,\OGASch_{X,\FRM}\)}{the global affine
    Schubert scheme induced by a reductive monoid \(\FRM\) over \(X\) such that \(\FRM^\x\) is quasi-split}
\end{align}
consisting of points such that for any
geometric point \(s\in S\), the image
of the generic point of \(X\x \Set{s}\) under \(\phi\) is contained in the
open locus \(E\subset \FRM_E\).

When \(\FRM\) is not clear from the
context, we will use \(\OGASch_{X,\FRM}\), etc. to emphasize the monoid. When
\(\FRM\) is commutative, in other words toric, we will more oftenly use \(\BD_X\)
    \nomenclature[\(B"cal_X_wow \)]{\(\BD_X\)}{the scheme \(\OGASch_X\) when
    \(\FRM^\x\) is a torus}
(resp.~\(\BD[+]_X\)) in place of \(\OGASch_{X}\) (resp.~\(\OGASch[+]_X\)).
We make the following definitions, the terms in which will be justified in what
follows in this section:

\begin{definition}
    \label[definition]{def:global_affine_schubert_as_mapping_stack}
    The stack \(\OGASch_X\) is called the \notion{global affine Schubert scheme
    associated with monoid \(\FRM\)}\index{global!affine Schubert
    scheme}\index{affine!Schubert scheme, global}. When \(\FRM\) is toric, we also call
    \(\BD_X\) the \notion{moduli of boundary divisors}\index{moduli!boundary
    divisors}\index{divisor!boundary, moduli of}.
\end{definition}

\subsection{}
The structure of \(\OGASch_X\) may be described with the help of the \inotion{Cartan
decomposition} for the monoid \(\FRM\). We have already implicitly used such
decomposition before, notably in \Cref{lem:arc_space_of_monoid_and_Cartan_cells}
where \(\FRM\in\FM(G^\SC)\), and the general
case is not hard either.

Let \(v\in\abs{X}\) be a closed point, \(\cO_v\) its
local ring with residue field \(k_v\), and \(F_v\) the local fractional field.
Choose a local uniformizer \(\pi_v\), then we may identify \(\cO_v\cong
k_v\powser{\pi_v}\) and \(F_v\cong k_v\lauser{\pi_v}\).
Let \(\sC_+\subset\sC\) be the subcone consisting of \(B\)-dominant
cocharacters. Then the set \(\sC_{v+}\defeq\sC_+(\cO_v)\) is a submonoid inside
the cocharacter lattice of the maximal \(\cO_v\)-split subtorus of \(T\). Using
the Cartan decomposition of \(G(F_v)\), we have
\begin{align}
    \label{eqn:Cartan_decomposition_of_reductive_monoids}
    G(F_v)\cap \FRM(\cO_v)=\coprod_{\lambda\in\sC_{v+}}G(\cO_v)\pi_v^\lambda
    G(\cO_v).
\end{align}
We call \eqref{eqn:Cartan_decomposition_of_reductive_monoids} the
\notion{Cartan decomposition of \(\FRM\) at \(v\)}\index{Cartan decomposition!of
a monoid}.

\begin{definition}
    Given a smooth curve \(X\) over \(k\), a \notion{\(\sC_+\)-valued
    divisor}\index{divisor!\(\sC_+\)-valued}, \notion{\(\sC_+\)-colored
    divisor}\index{divisor!\(\sC_+\)-colored} or \notion{boundary
    divisor}\index{divisor!boundary} is a formal sum
    \begin{align}
        \lambda_X=\sum_{v\in\abs{X}} \lambda_v\cdot v,
    \end{align}
    where \(v\) ranges over the closed points of \(X\), and
    \(\lambda_v\in\sC_{v+}\) is non-zero for only finitely many \(v\).
\end{definition}

\begin{remark}
    We use the name boundary divisor because it records the intersection type
    of the curve with
    the boundary of the unit group inside the monoid. The dependence of this
    notion on cone
    \(\sC\) would not cause any confusion since the cone is usually
    clear from the context.
\end{remark}

\subsection{}
To give a sense of the connection between the boundary divisors and \(\OGASch_X\),
we have the following generalization of \cite{BNS16}*{Lemma~3.4}:

\begin{lemma}
    \label[lemma]{lem:BD_and_boundary_divisor_k_points}
    There is a natural morphism from the groupoid \(\OGASch_X(k)\) to the set of
    boundary divisors on \(X\). The fiber of this morphism at
    \(\lambda_X=\sum_v\lambda_v\cdot v\) is equivalent to the set of points in
    the product of affine Schubert varieties
    \begin{align}
        \prod_{v}\Gr_{G,v}^{\lambda_v}(k_v).
    \end{align}
    In particular, when \(\FRM\) is toric, the morphism is an equivalence.
\end{lemma}
\begin{proof}
    Let \((E,\phi)\in\OGASch_X(k)\). At some \(v\in\abs{X}\) such that
    \(\phi(v)\) is not contained in \(E\), we can choose a local trivialization
    of \(E\)  over the formal disc \(X_v\) (by Lang's theorem, since the residue
    field \(k_v\) is finite). Then \(\phi(X_v)\) is a point in
    \(G(F_v)\cap\FRM(\cO_v)\), and so
    \eqref{eqn:Cartan_decomposition_of_reductive_monoids} gives an associated
    cocharacter \(\lambda_v\in\sC_{v+}\), which does not depend on the
    trivialization of \(E\). Thus, we obtain a boundary divisor on \(X\).
    Changing the trivialization of \(E\) changes \(\phi(X_v)\) by a right
    translation of some \(g\in G(\cO_v)\). It is
    then not hard to see the remaining claims by using
    Beauville--Laszlo's descent theorem. We leave the details to the reader.
\end{proof}

\begin{remark}
    It is tempting to upgrade the above lemma to a functorial version.
    Unfortunately, unless \(\FRM\) is toric, the map on the \(k\)-points is
    misleading about the topology of \(\OGASch_X\). We will see in
    \Cref{prop:global_affine_schubert_is_representable} that \(\OGASch_X\) is in
    fact proper, while the affine Schubert cells are not. The correct
    formulation of the relation between \(\OGASch_X\) and affine Schubert
    varieties appears in \Cref{sec:global_affine_schubert_scheme} when
    \(\FRM\in\FM(G^\SC)\).
\end{remark}

\subsection{}
Let \(\FRE=\FRM-G\) be the complement divisor with reduced
structure. Since \(\FRM\) is normal, \(\FRE\) is necessarily a Cartier
divisor. Then we have a canonical map of stacks
\begin{align}
    \label{eqn:global_affine_schubert_to_Hilbert}
    \OGASch_X\longto H_X\defeq\coprod_{n\ge 0}\Hilb^n{X}
\end{align}
where the right-hand side is the Hilbert scheme of points of \(X\), sending
\((E,\phi)\) to the divisor \(\phi^*\FRE_E\). Equivalently, the effective
Cartier divisor \(\FRE\) defines a map between quotient stacks
\begin{align}
    \Stack{\FRM/G} \longto \Stack{\bbA^1/\Gm}.
\end{align}
Applying the \(\OGASch_X\)-construction to both stacks, we obtain the map
\eqref{eqn:global_affine_schubert_to_Hilbert}.
We have the following approximation lemma, which is essentially
\cite{Bo17}*{Corollaire~41}:

\begin{lemma}
    \label[lemma]{lem:approximation_in_arc_monoid}
    Let \(S\) be a \(k\)-scheme and \(D\in H_X(S)\). Let
    \(I_D\subset\cO_{X\x S}\) be the ideal cutting out \(D\) and
    \(\hat{X}_S\) the relative (to \(S\)) spectrum of the completion of
    \(\cO_{X\x S}\) with respect to \(I_D\). Denote
    \(\hat{X}_S^\bullet=\hat{X}_S-D\). Then there exists some integer
    \(N>0\) depending only on \(D\) such that: for any \(\gamma,\gamma'\in
    G(\hat{X}_S^\bullet)\cap \FRM(\hat{X}_S)\) with
    \(\gamma^*\FRE=\gamma^{\prime*}\FRE=D\), if \(\gamma'\equiv\gamma\bmod
    I_D^N\), then \(\gamma'=\gamma g\) for some \(g\in G(\hat{X}_{S})\).
\end{lemma}
\begin{proof}
    The statement is local in \(S\), so there is no harm to assume that both \(G\)
    and \(\FRM\) are split over \(\hat{X}_{S}\).
    Choose \((\rho,V)\) a faithful representation of \(\FRM\).
    View \(\gamma\) as an \(\hat{X}_{S}^\bullet\)-point of \(G\), then the lattice
    \(\rho(\gamma)V(\hat{X}_{S})\subset V(\hat{X}_{S}^\bullet)\) is bounded
    between lattices
    \begin{align}
        I_D^{N}V(\hat{X}_{S})
        \subset\rho(\gamma)V(\hat{X}_{S})\subset V(\hat{X}_{S})
    \end{align}
    for a sufficiently large integer \(N\) depending only on \(D\). Thus,
    if \(\gamma\equiv\gamma'\bmod I_D^N\), the lattices
    \(\rho(\gamma)V(\hat{X}_{S})\) and
    \(\rho(\gamma')V(\hat{X}_{S})\) coincide, hence
    \(g=\gamma^{-1}\gamma'\in G(\hat{X}_{S}^\bullet)\) is such that
    \(\rho(g)\in\GL(V)(\hat{X}_{S})\). Since \(\rho(G)\) is closed in
    \(\GL(V)\), we have \(\rho(G)(\hat{X}_S^\bullet)\cap
    \GL(V)(\hat{X}_S)=\rho(G)(\hat{X}_S)\), and the lemma follows.
\end{proof}

\subsection{}
Next, we will prove the representability of \(\OGASch_X\). More precisely, we
show the following:

\begin{proposition}
    \label[proposition]{prop:global_affine_schubert_is_representable}
    The stack \(\OGASch_X\) is representable by a scheme that is of
    finite type and proper when restricted to any connected component of \(H_X\).
\end{proposition}
\begin{proof}
    Let \(S\) be a \(k\)-scheme.
    It is well-known (see \cite{FaGoIl05}*{Proposition~7.3.3} for example) that
    \(H_X\) is a disjoint union of symmetric powers of \(X\) over
    \(k\) and the statement is Zariski-local, so it suffices to assume that
    \(S\) is a Noetherian affine local scheme. Denote the Beilinson--Drinfeld
    Grassmannian of \(G\) over \(H_X\) by \(\Gr_{G,H_X}\) and its pullback
    to \(S\) by \(\Gr_{G,S}\). We then have by definition a canonical map
    \begin{align}
        \label{eqn:global_affine_schubert_to_global_grassmannian_naive}
        \OGASch_X\longto \Gr_{G,H_X}
    \end{align}
    sending \((E,\phi)\) over \(D\in H_X(S)\) to \((E,\phi|_{X\x S-D})\). The
    right-hand side is known to be an ind-scheme relative to \(H_X\)
    (cf.~\Cref{sec:review_BD_affine_grassmannian}).

    We first show that \(\OGASch_X\) is a sheaf. Indeed, given
    \((E,\phi)\in\OGASch_X(S)\) lying over \(D\in H_X(S)\), and \(\iota\) an
    automorphism of \((E,\phi)\). Then over
    \(X\x S-D\), \(\iota\) is the identity by definition. Since \(D\) is cut
    out by non-zero divisors, \(\iota\) must then be the identity over all \(X\x
    S\).

    Next, we show that the map of functors
    \eqref{eqn:global_affine_schubert_to_global_grassmannian_naive} is
    injective. Indeed, by Beauville--Laszlo's descent theorem, a pair
    \((E,\phi)\) is completely determined by its restriction to \(\hat{X}_S\).
    Replace \(S\) by an \'etale cover, we may trivialize \(E\) over
    \(\hat{X}_S\) and then \(\phi\) induces a point in \(\FRM(\hat{X}_S)\).
    Further restriicting to the punctured disc
    \(\hat{X}_S^\bullet\), it becomes a point in \(G(\hat{X}_S^\bullet)\), whose
    image in \(\Gr_{G,H_X}(S)\) is precisely that
    of \(E\). Since \(\FRM(\hat{X}_S)\to
    G(\hat{X}_S^\bullet)\) is injective, \((E,\phi)\) is then completely
    determined by its image in \(\Gr_{G,H_X}(S)\).

    We claim that the image of
    \eqref{eqn:global_affine_schubert_to_global_grassmannian_naive}, after
    pulling back to \(S\) via \(D\), is
    contained in a closed subscheme of \(\Gr_{G,S}\). Indeed, for any
    \(S\)-scheme \(Y\) and \((E_Y,\phi_Y)\in \OGASch_X(Y)\), we may replace
    \(Y\) by an \'etale covering so that \(E_Y\) is trivial over \(\hat{X}_Y\).
    Then by \Cref{lem:approximation_in_arc_monoid}, we may find some large \(N\)
    depending only on \(D\) such that the image of \((E_Y,\phi_Y)\) is contained
    in that of \(\FRM(\hat{X}_{S,N}\x_S Y)\), where \(\hat{X}_{S,N}\) is the
    closed subscheme of \(\hat{X}_S\) cut out by the \(N\)-th power of the
    ideal of \(D\). Since the jet schemes (with respect to the divisor \(D\)) of
    \(\FRM\) are representable by schemes of finite type, the claim then
    follows.

    By \cite{StacksP}*{\href{https://stacks.math.columbia.edu/tag/03ZQ}{Tags 03ZQ},
    \href{https://stacks.math.columbia.edu/tag/0A4X}{0A4X}}, it remains to prove that
    \eqref{eqn:global_affine_schubert_to_global_grassmannian_naive} satisfies
    the existence part of valuative criteria for properness. Let \(R\) be a
    discrete valuation ring with fractional field \(F\). Let \((E_F,\phi_F)\in
    \OGASch_X(F)\) be a point lying over \(D_F\in H_X(F)\) and whose image in
    \(\Gr_{G,H_X}\) extends to an \(R\)-point. In other words, \(E_F\) extends
    to a \(G\)-torsor over \(X\x R\) and \(D_F\) extends to an \(R\)-family of
    divisors \(D\) on \(X\x R\), and \(\phi_F\) extends to a trivialization of
    \(E|_{X\x R-D}\). Since \(X\x R\) is regular and \(\FRM_E\) is affine,
    \(\phi_F\) then further extends to a section of \(\FRM_E\) over all \(X\x
    R\) by Hartogs's theorem. This finishes the proof.
\end{proof}

\subsection{}
From the proof of \Cref{prop:global_affine_schubert_is_representable}, we see
that \((E,\phi)\in\OGASch_X\) really only captures the behavior of \(\phi(X)\)
near the divisor \(D\), which is local in nature. In other
words, it should only be sensitive to the boundary \(\FRE\) of \(G\) in
\(\FRM\). The following result is a precise version of this idea originated
from \cite{BNS16}*{\S~3}, which will
give us some flexibility when studying \(\OGASch_X\):

\begin{proposition}
    \label[proposition]{prop:embedding_of_cone_implies_embedding_of_divisor_moduli}
    Suppose we have a homomorphism \(\FRM\to\FRM'\) of monoids compatible
    with homomorphism of quasi-split reductive groups \(G\to G'\) over \(X\),
    such that the induced map of cones \(\sC\to\sC'\) is a closed embedding,
    then the induced map \(\OGASch_{X,\FRM}\to\OGASch_{X,\FRM'}\) is a closed
    embedding.
\end{proposition}
\begin{proof}
    As both functors are representable by schemes locally of finite type over
    \(k\), it suffices to show the map is a monomorphism and proper. 
    Properness is clear by \Cref{prop:global_affine_schubert_is_representable}.
    So we only need to show injectivity as functors.

    The saturated cone \(\sC\) generates a subsheaf of lattices in \(\CoCharG(T)\),
    giving a subtorus \(T_1\subset T\). The cone
    \(\sC\subset\CoCharG(T_1)\) corresponds to an \(T_1\)-toric scheme \(\FRT_1\)
    which also embeds as a closed subscheme in \(\FRT\). Let \(G_1\subset
    G\) be the subgroup generated by \(T_1\) and coroots in \(T_1\). Since
    \(\sC\) is invariant under the Weyl group, the projection of \(G_1\) to
    \(G^\AD\) is a direct factor.
    By \Cref{thm:monoid_classification_over_k_bar} (and using
    quasi-splitness), we have a monoid homomorphism \(\FRM\to G/G_1\) that fits
    into the Cartesian diagram
    \begin{equation}
        \begin{tikzcd}
            \FRM_1\ar[r]\ar[d] & \FRM\ar[d]\\
            \Set{1}\ar[r] & G/G_1
        \end{tikzcd},
    \end{equation}
    where \(\FRM_1\) is the monoid induced by \(\sC\) with unit group \(G_1\).
    We then have an induced Cartesian diagram
    \begin{equation}
        \begin{tikzcd}
            \OGASch_{X,\FRM_1}\ar[r]\ar[d] & \OGASch_{X,\FRM}\ar[d]\\
            \Set{b_0}\ar[r] & \BD_{X,G/G_1}
        \end{tikzcd},
    \end{equation}
    where \(b_0\) is the \(k\)-point of \(\OGASch_{X,G/G_1}\) corresponding to the
    trivial \(G/G_1\)-bundle equipped with the natural trivialization.
    But since \(\OGASch_{X,G/G_1}\) is obviously equivalent to a \(k\)-point, the map
    \(\OGASch_{X,\FRM_1}\to\OGASch_{X,\FRM}\) is an isomorphism.

    Thus, by replacing \((\FRM,G)\) with \((\FRM_1,G_1)\),
    we reduce to the case where
    \(G\to G'\) has finite and central kernel \(Z\). Let \(\FRM_2\) be the
    closure of \(G/Z\) in \(\FRM'\), then it is easy to see that
    \(\OGASch_{X,\FRM}\) maps into \(\OGASch_{X,\FRM_2}\), which clearly further
    \emph{injects} into \(\OGASch_{X,\FRM'}\). So we may replace \((\FRM',G')\) by
    \((\FRM_2,G/Z)\). We note that since \(\sC\to\sC'\) is injective,
    \Cref{thm:monoid_classification_over_k_bar} also implies that the boundary
    of \(\FRM\) must be mapped into that of \(\FRM'\). Let \(H_X'\) be
    the Hilbert scheme defined by the boundary of \(\FRM'\), then the map
    \(\OGASch_{X,\FRM}\to\OGASch_{X,\FRM'}\) is relative to \(H_X'\). By injectivity
    of the map \eqref{eqn:global_affine_schubert_to_global_grassmannian_naive},
    it suffices to show the injectivity of
    \(\Gr_{G,H_X'}\to\Gr_{G',H_X'}\). But this is immediate since \(G\to G'\) is
    an isogeny, and we are done.
\end{proof}
\subsection{}
We will continue studying the scheme \(\OGASch_X\) in
\Cref{sec:global_affine_schubert_scheme} for a general very flat monoid
\(\FRM\in\FM(G^\SC)\). For the remaining of this section, we will assume that
\(\FRM\) is toric, and to emphasize this we use \(\FRA\) (resp.~\(A\)) in place
of \(\FRM\) (resp.~\(G\)). Recall that in this case we also use the notation
\(\BD_X\) instead of \(\OGASch_X\).

We want to give a more precise description of \(\BD_X\) that is a functorial
upgrade of the equivalence in \Cref{lem:BD_and_boundary_divisor_k_points},
similar to \cite{BNS16}*{\S~3}.
Since this description is technical when \(A\) is not split, we will first cover
the split case following \textit{loc. cit.} for the ease of reading.

\subsection{}
Assume for now that \(\FRA\) and \(A\) are split.
To start, we would like to have 
a notion of \notion{degree}\index{degree!of boundary divisor} to stratify space \(\BD_X\) into more accessible
objects, similar to that on the Hilbert scheme \(H_X\). Given a boundary
divisor \(\lambda_X\), we can define its degree to
be the following element in \(\sC\):
\begin{align}
    \deg(\lambda_X)\defeq \sum_{v\in\abs{X}}[k_v:k]\lambda_v.
\end{align}
This straightforward definition by itself is too coarse, so we have the
following definition in \cite{BNS16}.
\begin{definition}
    \label[definition]{def:multisets_and_degrees}
    A \notion{multiset}\index{multiset!in a cone, split case} in \(\sC\) is an element \(\ul{\mu}\) of the free
    abelian monoid generated by \(\sC-\Set{0}\)
    \begin{align}
        \ul{\mu}=\sum_{\lambda\in\sC}\ul{\mu}(\lambda)e^\lambda\in\bigoplus_{\lambda\in\sC-\Set{0}}\bbN
        e^\lambda,
    \end{align}
    where \(\ul{\mu}(0)=0\) by convention. The
    \notion{degree}\index{degree!of multiset, split case} of multiset
    \(\ul{\mu}\) is defined as
    \begin{align}
        \deg(\ul{\mu})\defeq
        \sum_{\lambda\in\sC}\ul{\mu}(\lambda)\lambda\in\sC.
    \end{align}
\end{definition}

There is a natural partial order on the set of multisets by refinement: we say
\(\ul{\mu}\) refines \(\ul{\mu}'\) or \(\ul{\mu}\vdash\ul{\mu}'\) if the
difference \(\ul{\mu}-\ul{\mu}'\), viewed as an element of the free abelian
group generated by \(\sC-\Set{0}\), can be written as a sum of elements
of the form \(e^{\lambda_1}+e^{\lambda_2}-e^{\lambda_1+\lambda_2}\). Clearly,
refinement does not have an effect on degree.

\begin{definition}
    A multiset \(\ul{\mu}\) 
    is called \notion{primitive}\index{multiset!primitive, split case} if there is no strict refinement
    in the sense that \(\ul{\mu}'\vdash\ul{\mu}\) and \(\ul{\mu}'\neq\ul{\mu}\).
    Equivalently, \(\ul{\mu}\) is primitive if 
    \(\ul{\mu}(\lambda)=0\) unless \(\lambda\) is primitive in
    \(\sC\) (i.e., \(\lambda\) is not the sum of two non-zero elements in
    \(\sC\)).
\end{definition}

We can associate to a boundary divisor \(\lambda_X\) a multiset
\begin{align}
    \ul{\lambda}_{X}=\bigoplus_{\lambda\in\sC}\left(\sum_{\lambda_v=\lambda}[k_v:k]\right)e^\lambda,
\end{align}
and one sees that \(\deg(\lambda_X)=\deg(\ul{\lambda}_{X})\). We have the
following results.
\begin{proposition}[\cite{BNS16}*{Proposition~3.5}]
    \label[proposition]{prop:boundary_divisor_space_strat}
    There exists a unique stratification of \(\BD_X\) indexed by multisets:
    \begin{align}
        \BD_X=\coprod_{\ul{\mu}}\BD_{\ul{\mu}},
    \end{align}
    such that
    \begin{enumerate}
        \item \(\BD_{\ul{\mu}}(k)\) consists of those boundary divisors whose
            associated multiset is \(\ul{\mu}\).
        \item \label{item:boundary_mult_free_identification} 
            \(\BD_{\ul{\mu}}\) is isomorphic to the multiplicity-free locus of
            the (necessarily finite) direct product
            \begin{align}
                X_{\ul{\mu}}\defeq\prod_{\lambda\in\sC-\Set{0}}X_{\ul{\mu}(\lambda)},
            \end{align}
            where \(X_{\ul{\mu}(\lambda)}\) is the \(\ul{\mu}(\lambda)\)-th
            symmetric power of \(X\). Here \((D_\lambda)_\lambda\in X_{\ul{\mu}}\)
            being multiplicity-free means that 
            the total divisor \(\sum_\lambda D_\lambda\)
            is a multiplicity-free divisor (not just \(D_\lambda\) individually).
        \item \(\BD_{\ul{\mu}'}\) lies in the closure of \(\BD_{\ul{\mu}}\) if
            and only if \(\ul{\mu}\vdash\ul{\mu}'\).
    \end{enumerate}
\end{proposition}

\begin{corollary}[\cite{BNS16}*{Corollary~3.6}]
    We have the following description of connected and irreducible components of
    \(\BD_X\):
    \begin{enumerate}
        \item Each connected component of \(\BD_X\) contains a unique closed
            stratum of the form \(\BD_{e^\lambda}\) for some
            \(\lambda\in\sC\). Consequently, there is a canonical
            bijection \(\pi_0(\BD_X)\simeq \sC\), and
            \(\BD_{\ul{\mu}}\) lies in the connected component \(\BD[,\lambda]_X\) 
            associated with \(\lambda\) if and only if
            \(\deg(\ul{\mu})=\lambda\).
        \item Each irreducible component of \(\BD_X\) is the closure of a
            stratum \(\BD_{\ul{\mu}}\) where \(\ul{\mu}\) is a primitive
            multiset.
    \end{enumerate}
\end{corollary}

For a primitive multiset \(\ul{\mu}\), let \(\lambda_1,\ldots,\lambda_m\) be the
coweights such that \(\ul{\mu}(\lambda_i)\neq 0\), then
the isomorphism in \Cref{prop:boundary_divisor_space_strat}
\eqref{item:boundary_mult_free_identification} can be extended to a
morphism
\begin{align}
    \label{eqn:sym_power_to_boundary_moduli}
    \prod_{i=1}^m X_{\ul{\mu}(\lambda_i)}\longto \bar{\BD_{\ul{\mu}}}.
\end{align}
\begin{corollary}[\cite{BNS16}*{Corollary~3.7}]
    \label[corollary]{cor:res_sing_of_boundary_moduli}
    The disjoint union of morphisms \eqref{eqn:sym_power_to_boundary_moduli} for
    which \(\ul{\mu}\) ranges over all primitive multisets is both a
    normalization and a resolution of singularity of \(\BD_X\). In particular,
    it is a finite map.
\end{corollary}

\subsection{}
Given \(\FRA\) with associated cone \(\sC\), there are two constructions to
associate a toric variety of standard type (cf.~\Cref{def:toric_standard_type})
to \(\FRA\). First, let
\(\lambda_1,\ldots,\lambda_e\) be the primitive elements of \(\sC\),
each of which gives a monoidal homomorphism \(\bbA^1\to\FRA\). The multiplication map
then defines a homomorphism
\begin{align}
    \tilde{\FRA}\defeq \bbA^e\longto\FRA.
\end{align}
On the other hand, let \(\beta_1,\ldots,\beta_{e'}\) be a set of generators of the
dual cone \(\sC^*\) \emph{as an \(\bbN\)-monoid}, and let \(\ul{\sC}\) be the
\(\bbN\)-dual of the free abelian monoid generated by \(\beta_i\), then we
have injective map of cones \(\sC\to \ul{\sC}\) and a monomorphism of
monoids
\begin{align}
    \FRA\longto
    \ul{\FRA}\defeq\Spec{k[e^{\beta_1},\ldots,e^{\beta_{e'}}]}=\bbA^{e'}.
\end{align}
Note that \(\tilde{\FRA}\) is canonical since \(\sC\) is strictly convex, 
while \(\ul{\FRA}\) depends on the choice of generators \(\beta_i\) and may be
replaced by any faithful \(k\)-algebraic representation of \(\FRA\).

Both \(\BD_{X,\tilde{\FRA}}\) and \(\BD_{X,\ul{\FRA}}\) are clearly smooth.
By \Cref{cor:res_sing_of_boundary_moduli}, the induced morphism
\(\BD_{X,\tilde{\FRA}}\to\BD_{X,\FRA}\) is a finite resolution of
singularity, and by
\Cref{prop:boundary_divisor_space_strat}, the morphism
\(\BD_{X,\FRA}\to\BD_{X,\ul{\FRA}}\) is a closed embedding.
Note that both \(\tilde{\FRA}\) and \(\ul{\FRA}\) may be defined for non-split
\(\FRA\) as well.

\subsection{}
Now we deal with the case where the torus is not necessarily split or constant.
Recall that the sheaf \(\CoCharG(A)\) can be represented by a countable 
union of finite \'etale covers of \(X\), and the cone \(\sC\) a subscheme
consisting of some connected components therein. There is a canonical component
in \(\sC\) corresponding to the zero cocharacter, still denoted by \(0\).
The relative
symmetric power \(\sC_{m/X}\defeq\Sym_X^m\sC\) is easily
seen smooth by looking \'etale-locally over \(X\). 
We also have the summation map by definition:
\begin{align}
    \Sigma\colon \sC_{m/X}\longto \sC.
\end{align}
Since \(\sC\) is a disjoint union of smooth projective curves, its scheme of
connected components \(\pi_0(\sC)\) is an \'etale scheme over \(k\) locally of
finite type. For \(\lambda\in\pi_0(\sC)\) a closed point and \(\sC^\lambda\) the
corresponding component in \(\sC\), we let
\(k_{\lambda}=\RH^0(\sC^\lambda,\cO_{\sC})\), then \(\lambda=\Spec{k_\lambda}\).
The following definitions generalize \Cref{def:multisets_and_degrees}:

\begin{definition}
    A \notion{multiset}\index{multiset!in a cone, non-split case} \(\ul{\mu}\) is an element in 
    the free abelian monoid generated by the closed points in
    \(\pi_0(\sC)-\Set{0}\). For
    \(0\neq\ul{\mu}=\sum_{\lambda\in\pi_0(\sC)}\ul{\mu}(\lambda)\lambda\) with
    \(\ul{\mu}(\lambda)\in\bbN\) (again, \(\ul{\mu}(0)=0\) by convention), we
    have
    \begin{align}
        \Sigma\colon \prod_{\lambda}\sC^\lambda_{\ul{\mu}(\lambda)/X}\longto \sC,
    \end{align}
    where \(\sC^\lambda_{\ul{\mu}(\lambda)/X}\) is relative symmetric power of
    \(\sC^\lambda\) over \(X\) with degree \(\ul{\mu}(\lambda)\) and the (finite) direct
    product is taken as the Cartesian product over \(X\). Any point
    \(\ul{\mu}'\) in such image is called \notion{a degree of
    \(\ul{\mu}\)}\index{degree!of multiset, non-split case}, regarded
    itself as a multiset.
\end{definition}

If \(\ul{\mu}=\ul{\mu}_1+\ul{\mu}_2\), and \(\ul{\mu}_2'\) is a degree of
\(\ul{\mu}_2\), then \(\ul{\mu}'=\ul{\mu}_1+\ul{\mu}_2'\) is again a multiset.
In this case \(\ul{\mu}\) can be seen as a refinement of \(\ul{\mu}'\).
More generally, a multiset \(\ul{\mu}\) is said to be a
refinement of \(\ul{\mu}'\) or \(\ul{\mu}\vdash\ul{\mu}'\) 
if \(\ul{\mu}'\) can be obtained
from \(\ul{\mu}\) after finite steps of taking degrees of its summands.
Note that there may be multiple possible degrees for a given multiset.
The notion of degree is less useful compared to the case of split \(\FRA\),
since they no longer correspond bijectively to connected components of
\(\BD_X\).

\begin{definition}
    A multiset \(\ul{\mu}\) is called
    \notion{primitive}\index{multiset!primitive, non-split case} if \(\ul{\mu}\neq 0\)
    and cannot be further refined.
\end{definition}

\subsection{}
Given a boundary divisor \(\lambda_X\), we may define an induced multiset
\(\ul{\lambda}_X\) as follows: let \(v_1,\ldots, v_m\) be the mutually distinct
points of \(X\) such that \(\lambda_{v_i}\neq 0\). Since \(\lambda_{v_i}\) can be
identified with a point in \(\sC(k_{v_i})\), it is contained in a unique
component \(\lambda_i\in\pi_0(\sC)\). Let \(d_i=[k_{v_i}:k]\) be the
degree of point \(v_i\) over \(k\).  
Then \(\ul{\lambda}_X\) is just the formal sum
\begin{align}
    \ul{\lambda}_X\defeq\sum_{i=1}^m d_i\lambda_i.
\end{align}

\subsection{}
Now we are ready to prove various properties of \(\BD_X\) parallel to the split
case.

\begin{proposition}
    \label[proposition]{prop:boundary_divisor_space_strat_nonsplit}
    There exists a unique stratification of \(\BD_X\) indexed by multisets:
    \begin{align}
        \BD_X=\coprod_{\ul{\mu}}\BD_{\ul{\mu}},
    \end{align}
    such that
    \begin{enumerate}
        \item \(\BD_{\ul{\mu}}(k)\) consists of those boundary divisors whose
            associated multiset is \(\ul{\mu}\).
        \item \label{item:boundary_mult_free_identification_nonsplit} 
            \(\BD_{\ul{\mu}}\) is isomorphic to the multiplicity-free locus of
            the (necessarily finite) direct product
            \begin{align}
                X_{\ul{\mu}}\defeq
                \prod_{\lambda\in\pi_0(\sC)}\sC^\lambda_{\ul{\mu}(\lambda)},
            \end{align}
            where \((D_\lambda)_\lambda\in
            X_{\ul{\mu}}\) being multiplicity-free means that 
            the total divisor \(\sum_\lambda D_\lambda\)
            is a multiplicity-free divisor on \(X\) after pushing forward to
            \(X\).
        \item \(\BD_{\ul{\mu}'}\) lies in the closure of \(\BD_{\ul{\mu}}\) if
            and only if \(\ul{\mu}\vdash\ul{\mu}'\).
    \end{enumerate}
\end{proposition}
\begin{proof}
    First, in the case of a simple multiset
    \(\ul{\mu}=\lambda\in\pi_0(\sC)-\Set{0}\), we show \(\BD_{\ul{\mu}}\simeq
    \sC^\lambda\). In the case \(\FRA=p_*\bbA^1\), where \(p\colon
    \sC^\lambda\to X\) is the natural map, the result is straightforward by
    noting \(\BD_{X,\FRA}=\BD_{\sC^\lambda,\bbA^1}\).
    In general, the inclusion \(\sC^\lambda\to\sC\) gives an embedding of
    cones corresponding to morphism of toric schemes \(p_*\bbA^1\to\FRA\)
    (with homomorphism of tori \(p_*\Gm\to A\)). We also clearly have a
    natural bijection on geometric points \(\BD_{\ul{\mu}}(\bar{k})\simeq
    \sC^\lambda(\bar{k})\) (cf.~\Cref{lem:BD_and_boundary_divisor_k_points}).
    Then we are done by applying
    \Cref{prop:embedding_of_cone_implies_embedding_of_divisor_moduli}.
    Let
    \begin{align}
        \theta_{\lambda}\colon\sC^\lambda\longto\BD_{\ul{\mu}}
    \end{align}
    be the isomorphism.

    The monoidal structure on \(\FRA\) induces a monoidal structure on
    \(\BD_X\): given \(\theta_1=(E_1,\phi_1)\) and \(\theta_2=(E_2,\phi_2)\) in
    \(\BD_X\), the product \(\theta_1\otimes\theta_2\) is the pair \((E_1\x^A
    E_2,\phi_1\otimes\phi_2)\), and this monoidal structure is clearly
    commutative. Given a multiset \(\ul{\mu}\), we can define a
    map (the product is the Cartesian product over \(k\))
    \begin{align}
        \label{eqn:sym_to_boundary_moduli_nonsplit}
        \iota_{\ul{\mu}}\colon
        X_{\ul{\mu}}=\prod_{\lambda\in\pi_0(\sC)}\sC^\lambda_{\ul{\mu}(\lambda)}&\longto\BD_X\\
        \prod_\lambda(x^\lambda_1,\ldots,x^\lambda_{\ul{\mu}(\lambda)})&\longmapsto
        \bigotimes_{\lambda}[\theta_\lambda(x^\lambda_1)
        \otimes\cdots\otimes\theta_\lambda(x^\lambda_{\ul{\mu}(\lambda)})].
    \end{align}
    Note that the ordering of the product has no effect on the resulting map, so
    it is well-defined. The map \eqref{eqn:sym_to_boundary_moduli_nonsplit} is
    proper because the source is a projective variety and the target is a
    separated scheme. Let \(\bar{\BD_{\ul{\mu}}}\) be the image of
    \(\iota_{\ul{\mu}}\), which is a reduced closed subscheme of \(\BD_X\)
    because the map is proper, and the source is reduced. Let
    \(X_{\ul{\mu}}^\circ\) be the multiplicity-free locus designated by
    part~\eqref{item:boundary_mult_free_identification_nonsplit}, then by
    checking the boundary divisors at
    \(\bar{k}\)-point level we can see that
    \begin{align}
        \iota_{\ul{\mu}}^{-1}(\iota_{\ul{\mu}}(X_{\ul{\mu}}))=X_{\ul{\mu}},
    \end{align}
    hence it has open image in \(\bar{\BD_{\ul{\mu}}}\), denoted by
    \(\BD_{\ul{\mu}}\). In fact, one may check that
    \(\bar{\BD_{\ul{\mu}}}-\BD_{\ul{\mu}}\) is just the union of the images of
    \(\iota_{\ul{\mu}'}\) such that \(\ul{\mu}\vdash \ul{\mu}'\).
    This stratification exhausts \(\BD_X\) by checking on \(\bar{k}\)-points
    which is straightforward because \(\BD_X(\bar{k})\) is in bijection with
    the set of  boundary divisors defined over \(\bar{k}\).

    It remains to prove that the restriction of \(\iota_{\ul{\mu}}\) to
    \(X_{\ul{\mu}}^\circ\) is an isomorphism. For this purpose we embed
    \(\BD_X\) into \(\BD_{X,\FRA'}\) for some \(\FRA'\) of standard type. In
    the latter case, the map \(X_{\ul{\mu}}^\circ\to \BD_{\ul{\mu},\FRA'}\) is an
    isomorphism by direct computation. Therefore, we have a section
    \(\BD_{\ul{\mu}}\to X_{\ul{\mu}}^\circ\) to \(\iota_{\ul{\mu}}\) by
    composing with embedding \(\BD_X\to \BD_{X,\FRA'}\). Since
    \(X_{\ul{\mu}}^\circ\)
    is integral and \(\iota_{\ul{\mu}}\) is a bijection on \(\bar{k}\)-points, 
    it has to be an isomorphism.
\end{proof}

\begin{corollary}
    \label[corollary]{cor:boundary_moduli_components_nonsplit}
    \begin{enumerate}
        \item Every connected component of \(\BD_X\) contains a (not necessarily
            unique) closed stratum \(\BD_{\ul{\mu}}\) where \(\ul{\mu}=\lambda\)
            is a simple multiset. In particular, there is an equivalence
            relation \(\sim\) on the set of closed points of \(\pi_0(\sC)\) such
            that \(\pi_0(\BD_X)\) (as a set) is in canonical bijection with
            \(\pi_0(\sC)/\sim\). Each \(\BD_{\ul{\mu}}\) lies in the connected
            component corresponding to \([\lambda]\) if and only if there is some
            \(\lambda'\sim\lambda\) such that \(\lambda'\) is a degree of
            \(\ul{\mu}\).
        \item Irreducible components of \(\BD_X\) are the closures of strata
        \(\BD_{\ul{\mu}}\) where \(\ul{\mu}\) are primitive multisets,
        i.e.~\(\ul{\mu}(\lambda)\neq 0\) only if \(\lambda\) cannot be further
        refined.
    \end{enumerate}
\end{corollary}
\begin{proof}
    For the first part, the minimal elements of the partial order \(\vdash\) are
    the simple ones, and the closure of \(\BD_{\ul{\mu}}\) contains and only
    contains those \(\BD_{\lambda}\) where \(\lambda\) is a degree of
    \(\ul{\mu}\).

    For the second part, the maximal elements of the partial order \(\vdash\)
    are the ones that cannot be refined any further, and those multisets are
    exactly the ones described by the statement.
\end{proof}

\begin{corollary}
    For each primitive multiset \(\ul{\mu}\), the map \(\iota_{\ul{\mu}}\) gives
    a normalization which is also a resolution of singularity of 
    \(\bar{\BD_{\ul{\mu}}}\).
\end{corollary}
\begin{proof}
    The map \(\iota_{\ul{\mu}}\) is birational and finite by
    \Cref{prop:boundary_divisor_space_strat_nonsplit}, thus must be a
    normalization map. Since the source \(X_{\ul{\mu}}\) is smooth, it is also a
    resolution of singularity.
\end{proof}

\begin{definition}
    \label[definition]{def:mult_free_boundary_divisor}
    A boundary divisor \(\lambda_X\in\BD_X(\bar{k})\) is called
    \notion{primitive}\index{divisor!boundary, primitive} or
    \notion{multiplicity-free}\index{divisor!boundary, multiplicity-free} if it is contained in
    \(\BD_{\ul{\mu}}\) for some primitive multiset \(\ul{\mu}\).
\end{definition}

\subsection{}
\label{sub:BD_grassmannian_by_boundary_divisor}
To facilitate the discussion in \Cref{sec:global_affine_schubert_scheme}, we
return to the quasi-split reductive group \(G\) and discuss a
particular instance of Beilinson--Drinfeld Grassmannian
(cf.~\Cref{sec:review_BD_affine_grassmannian}).
The map \(G\to G^\AD\) induces a
homomorphism \(\CoCharG(T)\to\CoCharG(T^\AD)\) as well as the dominant cones
therein. The dominant cone \(\CoCharG(T^\AD)_+\) is locally freely 
generated by the fundamental coweights, hence it corresponds to a \(T^\AD\)-toric
scheme \(\FRA^\AD\) of standard type over \(X\). Viewing \(\CoCharG(T^\AD)\) as
a countable \'etale cover of \(X\), then there is a canonical union of connected
components \(\sC^{\CoWt}\subset \CoCharG(T^\AD)_+\) corresponding to the set
of all fundamental coweights.
The monoid \(\FRA^\AD\) is then the pushforward
\begin{align}
    \FRA^\AD=p^{\CoWt}_*\bbA^1,
\end{align}
where \(p^{\CoWt}\colon \sC^{\CoWt}\to X\) is the natural map. The moduli of
\(\CoCharG(T^\AD)_+\)-valued boundary divisors \(\BD_{X,\FRA^\AD}\) can thus be
identified with union of symmetric powers
\begin{align}
    \BD_{X,\FRA^\AD}=\coprod_{d=0}^\infty \sC^{\CoWt}_d,
\end{align}
on which we have the affine Grassmannian, arc group, and so on (using either \(G\)
or \(G^\AD\)).
The dominant cone \(\CoCharG(T)_+\) for \(T\) is not necessarily strictly
convex, but any saturated strictly convex subcone \(\sC\subset\CoCharG(T)_+\) 
determines a \(T\)-toric scheme \(\FRA\),
and the homomorphism of cones \(\sC\to\CoCharG(T^\AD)_+\) induces a proper
map of boundary moduli
\begin{align}
    \BD_{X,\FRA}\longto\BD_{X,\FRA^\AD},
\end{align}
and thus one may pull back the affine Grassmannian, arc group, etc. from
\(\BD_{X,\FRA^\AD}\). Therefore, such notions make sense for
\(\BD_{X,\FRA}\) even if it is not a union of symmetric powers of curves (or
even smooth).

\section{Global Affine Schubert Scheme}%
\label{sec:global_affine_schubert_scheme}

Just like the Beilinson--Drinfeld Grassmannian is ``affine Grassmannians in a family'',
we also have the corresponding family of affine Schubert schemes. 
Usually, the affine Schubert varieties are defined as certain subschemes
of affine Grassmannian. However, things get tricky when multiple points are
involved due to the complication caused by the torsions in the fundamental group
of \(G\) (or equivalently, in \(\pi_0(\Bun_G)\)). See, for example,
\cite{Ya04}*{Definition~2.10} when the group is split.

Here we give a definition of affine Schubert scheme using reductive monoids.
This method allows us to directly define affine Schubert schemes without
referring to affine Grassmannian at all (although its properties are still
better described using the affine Grassmannian). In fact, we will give two
equivalent constructions: one already given in
\Cref{def:global_affine_schubert_as_mapping_stack}, and another
using arc schemes.

\subsection{}
We let curve \(X\) as in \Cref{sec:boundary_divisors} and let \(G\) be a
quasi-split reductive group over \(X\). Let \(\FRM\in\FM(G^\SC)\), and we have the moduli of boundary
divisors \(\BD_X\) associated with the abelianization \(\FRA_\FRM\). Here we no
longer assume that \(G=\FRM^\x\).
Recall from \Cref{sec:boundary_divisors} the functor
\(\OGASch_X=\OGASch_{X,\FRM}\), which we have shown in
\Cref{prop:global_affine_schubert_is_representable} is representable by a
disjoint union of projective schemes. The natural map
\begin{align}
    \Stack*{\FRM/\FRM^\x}\longto \Stack*{\FRA_\FRM/\FRA_\FRM^\x}
\end{align}
induces a map of schemes
\begin{align}
    \OGASch_X\longto\BD_X.
\end{align}

Let \(\FRM^\AD=\FRM\git Z^\SC\in\FM(G^\AD)\), and by
\Cref{prop:embedding_of_cone_implies_embedding_of_divisor_moduli}, we have a
closed embedding
\begin{align}
    \OGASch_X\longto\OGASch_{X,\FRM^\AD}.
\end{align}
The abelianization of \(\FRM^\AD\) is still \(\FRA_\FRM\), and so the above map
is defined over \(\BD_X\). Similar to
\eqref{eqn:global_affine_schubert_to_global_grassmannian_naive}, we
have a canonical closed embedding
\begin{align}
    \OGASch_{X,\FRM^\AD}\longto \Gr_{\FRM^{\AD\x},\BD_X},
\end{align}
where the right-hand side is the affine Grassmannian of \(\FRM^{\AD\x}\)
parametrized by \(\BD_X\) (see \Cref{sub:BD_grassmannian_by_boundary_divisor}
and note that \(\FRA_{\Env(G^\SC)}=\FRA^\AD\) therein).
We have an induced map by composition
\begin{align}
    \label{eqn:global_affine_schubert_to_global_grassmannian}
    \OGASch_{X}\longto \Gr_{G^\AD,\BD_X}.
\end{align}

\subsection{}
Let \(b=(\cL,\theta)\in\BD_X(S)\), where \(\cL\) is a \(\FRA_{\FRM}^\x\)-torsor
over \(X\x S\) and \(\theta\) is a section of \(\FRA_{\FRM,\cL}\) that is
generically contained in \(\FRA_{\FRM,\cL}^\x\) over every geometric point
\(s\in S\). The torus \(\FRA_{\FRM}^\x\) does not act on \(\FRM\),
but it does act on \(\FRM^\AD\). We have then
a Cartesian diagram
\begin{equation}
    \begin{tikzcd}
        \FRM^\AD_b \ar[r]\ar[d] & \FRM^\AD_\cL \ar[d]\\
        X\x S \ar[r, "\theta"] & \FRA_{\FRM,\cL}
    \end{tikzcd}
\end{equation}
where \(\FRM^\AD_b\) is a bundle of \(G^\AD\x G^\AD\)-spherical varieties over
\(X\x S\).

Recall that the product of all simple roots
cuts out the numerical boundary divisor \(\FRE_{\Env(G^\SC)}\) on
\(\FRA_{\Env(G^\SC)}\). It pulls back to a principal divisor on \(\FRA_\FRM\),
whose underlying topological space is contained in \(\FRA_\FRM-\FRA_\FRM^\x\). Further
pulling back through \(\theta\), and we have an effective Cartier divisor
\(\FRE_b\) on \(X\x S\), relative to \(S\). Let \(\hat{X}_{\FRE_b}\) be the
formal completion of \(X\x S\) at \(\FRE_b\). Let \(\Arc_b{G^\SC}\) and
\(\Arc_b{G^\AD}\) be the arc
groups defined using divisor \(\FRE_b\), and similarly \(\Arc_b{\FRM^\AD_b}\) the
arc space of \(\FRM^\AD_b\). Then we have arc spaces over \(\BD_X\):
\begin{align}
    \Arc_{\BD_X}{G}(S)&\defeq\Set*{(b,g)\given b\in\BD_X(S), g\in
    G(\hat{X}_{\FRE_b})},\\
    \Arc_{\BD_X}{G^\AD}(S)&\defeq\Set*{(b,g)\given b\in\BD_X(S), g\in
    G^\AD(\hat{X}_{\FRE_b})},\\
    \Arc_{\BD_X}(\FRM^\AD/\FRA_{\FRM})(S)&\defeq\Set*{(b,x)\given b\in\BD_X(S),
    x\in \FRM^\AD_b(\hat{X}_{\FRE_b})}.
\end{align}
The arc group \(\Arc_{\BD_X}{G}\) acts on \(\Arc_{\BD_X}(\FRM^\AD/\FRA_{\FRM})\)
both on the left and on the right by translation and these actions factor
through \(\Arc_{\BD_X}{G^\AD}\).
The affine Grassmannian \(\Gr_{G,\BD_X}\) may also be
defined using \(\Arc_{\BD_X}{G}\) and associated loop group.

\begin{proposition}
    \label[proposition]{prop:embed_global_affine_schubert_and_its_fibers}
    The map
    \eqref{eqn:global_affine_schubert_to_global_grassmannian} is
    a closed embedding and the image is stable under the left action of
    \(\Arc_{\BD_X}G\). Moreover, if \(b\in \BD_X(\bar{k})\) with the associated
    boundary divisor \(\sum_{\bar{v}\in X(\bar{k})}\lambda_{\bar{v}}\bar{v}\),
    the \(\bar{k}\)-points of its fiber in \(\OGASch_X\) is isomorphic to
    \begin{align}
        \prod_{\bar{v}\in
        X(\bar{k})}\Gr_{G_{\bar{v}}^\AD}^{\le-w_0(\lambda_{\bar{v},\AD})}(\bar{k}).
    \end{align}
\end{proposition}
\begin{proof}
    The map is proper because the source is a proper scheme and the image is
    contained in a subscheme of the target. To show it is a closed embedding, we
    need only to prove the injectivity. Let \(S\)  be a \(k\)-scheme, and
    \((E_1,\phi_1)\) and \((E_2,\phi_2)\) be two \(S\)-points of \(\OGASch_X\)
    with the same image in \(\Gr_{G^\AD,\BD_X}\). Let \(b\) be their common
    image in \(\BD_X\) with associated Cartier divisor \(\FRE_b\subset X\x S\).
    Since \(\OGASch_X\) embeds into \(\OGASch_{X,\FRM^\AD}\), we may replace
    \((E_i,\phi_i)\) by their images in the latter. Since \((E_i,\phi_i)\) have
    the same image in both \(\Gr_{G^\AD,\BD_X}\) and \(\BD_X\), both \(E_i\)
    induce the same \(\FRA_\FRM^\x\)-torsor and the same \(G^\AD\)-torsor. Since
    \(\FRM^{\AD\x}=\FRA_\FRM^\x\x G^\AD\), we see that \(E_1=E_2\). To see
    \(\phi_1=\phi_2\), we only need to show that they are isomorphic over \(X\x
    S-\FRE_b\), but this is none other than saying they have the same image in
    \(\Gr_{G^\AD,\BD_X}\).

    The action of \(\Arc_{\BD_X}G\) is defined as follows: given
    \((E,\phi)\in\OGASch_X(S)\) and \(g\in\Arc_{\BD_X}G^\SC(S)\) both over
    \(b\in\BD_X\), the trivialization \(\phi\) of \(E\) over the punctured disc
    \(\hat{X}_{\FRE_b}^\bullet\) may be modified by \(g\), which can then be glued
    back with \((E,\phi)|_{X\x S-\FRE_b}\) using Beauville--Laszlo's
    descent theorem to obtain another point in \(\OGASch_X\). This action is compatible
    with \(\Arc_{\BD_X}G^\SC\)'s action on \(\Gr_{G^\AD,\BD_X}\), and the latter
    factors through \(\Arc_{\BD_X}G^\AD\). This then induces an action of
    \(\Arc_{\BD_X}G\) on \(\OGASch_X\) which embeds into the affine
    Grassmannian.

    For the second claim, we let \(S=\bar{k}\). Over
    \(\hat{X}_{\FRE_b}\), the \(G\)-torsor \(E\) can be trivialized because the
    residue algebra is a finite product of \(\bar{k}\). Then \(\bar{k}\)-points
    in the fiber of \(b\) in \(\OGASch_X\) is represented by
    \(\hat{X}_{\FRE_b}\)-points of \(\FRM\) whose
    image in \(\FRA_\FRM\) lies in the product
    \begin{align}
        \prod_{\bar{v}\in
        X(\bar{k})}\pi_{\bar{v}}^{\lambda_{\bar{v}}}\FRA_\FRM^\x(\breve{\cO}_{\bar{v}}).
    \end{align}
    The claim then follows from \Cref{lem:arc_space_of_monoid_and_Cartan_cells}.
\end{proof}

\subsection{}
We now give an alternative definition of \(\OGASch_X\) akin to the loop group
construction of the affine Grassmannian.
For this, we will need the arc space of \(\FRM\), but
there is no direct way to
define the arc space of \(\FRM\) over \(\BD_X\) (we can do that only over base
\(\cB_X\) defined in \Cref{sec:Constructions}).

However, we can define arc space of \(\FRM\) locally over \(\BD_X\), namely, for
any \(b=(\cL,\theta)\in\BD_X(S)\), we may choose locally over \(S\) a
trivialization of \(\cL\) over \(\hat{X}_{\FRE_b}\), thus it lifts to a
\(Z_\FRM\)-torsor, and then we may define \(\Arc_b\FRM_b\) in a similar way as
\(\Arc_b\FRM^\AD_b\). In particular, if \(b\) is a \(\bar{k}\)-point, with
corresponding boundary divisor \(\sum_{\bar{v}\in
X(\bar{k})}\lambda_{\bar{v}}\bar{v}\), then we have
\begin{align}
    \Arc_b{\FRM_b}(\bar{k})\cong\prod_{\bar{v}\in
        X(\bar{k})}\Bigg[\sum_{\substack{\mu\in\CoCharG(T_{\bar{v}}^\AD)_+\\\mu\le-w_0(\lambda_{\bar{v},\AD})}}
    G^\SC(\breve{\cO}_{\bar{v}})\pi_{\bar{v}}^{(\lambda_{\bar{v}},\mu)}G^\SC(\breve{\cO}_{\bar{v}})\Bigg].
\end{align}
Observe that even though \(\Arc_{\BD_X}(\FRM/\FRA_\FRM)\) can only be locally
defined over \(\BD_X\) (we may fix local charts and define it as a presheaf),
we still have an injective map (as presheaves) of quotients
\begin{align}
    \Arc_{\BD_X}(\FRM/\FRA_\FRM)/\Arc_{\BD_X}G^\SC\longto
    \Arc_{\BD_X}(\FRM^\AD/\FRA_\FRM)/\Arc_{\BD_X}G^\AD,
\end{align}
and it does not depend on the choice of local trivialization of \(\cL\),
because the ambiguity caused by finite group \(Z^\SC\) is killed after taking
quotient. Since the right-hand side is in fact a quotient of sheaves, the
associated sheaf of the left-hand side is canonically defined (the presheaf
itself, however, depends on the choice of charts).

\begin{proposition}
    We have a canonical isomorphism of fpqc sheaves
    \begin{align}
        \OGASch_X\simeq\Arc_{\BD_X}(\FRM/\FRA_{\FRM})/\Arc_{\BD_X}{G^\SC}.
    \end{align}
\end{proposition}
\begin{proof}
    Given \((E,\phi)\in\OGASch_X(S)\) lying over \(b=(\cL,\theta)\in\BD_X\), we
    replace \(S\) by an \'etale cover so that the restriction of both \(E\) and
    \(\cL\) to formal disc \(\hat{X}_{\FRE_b}\) may be trivialized and so
    \(\phi\) induces a point \(\gamma\in\FRM_b(\hat{X}_{\FRE_b})\). The image of
    \(\gamma\) in \(\Arc_{\BD_X}(\FRM/\FRA_{\FRM})/\Arc_{\BD_X}{G^\SC}\) does
    not depend on the choice of the trivializations. Therefore we have a
    well-defined map from \(\OGASch_X\) to
    \(\Arc_{\BD_X}(\FRM/\FRA_{\FRM})/\Arc_{\BD_X}{G^\SC}\).

    On the other hand, suppose \(\Arc_{\BD_X}(\FRM/\FRA_\FRM)\) is defined over
    \(S\), then by gluing with the trivial torsor together with its trivial
    section, we obtain a point in \(\OGASch_X(S)\), and such map clearly factors
    through \(\Arc_{\BD_X}(\FRM/\FRA_{\FRM})/\Arc_{\BD_X}{G^\SC}\). This way we
    obtain the map in the opposite direction. It is then straightforward to see
    these two maps are inverses to each other.
\end{proof}

\subsection{}
The following two results are immediate.

\begin{lemma}
    If \(\phi\colon\FRM'\to \FRM\) is a map in \(\FM(G^\SC)\), then
    we have canonical isomorphism
    \begin{align}
        \phi^*\OGASch_{X,\FRM}\simeq\OGASch_{X,\FRM'}.
    \end{align}
\end{lemma}
\begin{proof}
    Straightforward from the definition.
\end{proof}

\begin{lemma}
    \label[lemma]{lem:arc_monoid_to_affine_schubert_is_torsor}
    Let \(\FRM\in\FM(G^\SC)\). Suppose \(\Arc_b\FRM_b\) makes sense for some
    \(b\in\BD_X(S)\), then the map
    \begin{align}
        \Arc_{b}\FRM_b\longto \OGASch_b\defeq\OGASch_X\x_{\BD_X} b
    \end{align}
    is formally smooth, and its fiber is a torsor under \(\Arc_b G^\SC\).
\end{lemma}
\begin{proof}
    The arc group \(\Arc_b{G^\SC}\) is formally smooth
    and its action by right translation on \(\Arc_b\FRM_b\) is
    free.
\end{proof}

\subsection{}
By replacing \(\FRM\) with the big-cell locus \(\FRM^\circ\), we obtain a
generalized notion of affine Schubert cells. More precisely, we define
\begin{align}
    \OGASch[\circ]_X=\Arc_{\BD_X}(\FRM^\circ/\FRA_{\FRM})/\Arc_{\BD_X}{G^\SC}
    \nomenclature[\(Q_X_wow_circ \)]{\(\OGASch[\circ]_X\)}{the big-cell locus of
    \(\OGASch_X\) for \(\FRM\in\FM(G^\SC)\)}
\end{align}
whose fiber at \(b\in\BD_X(\bar{k})\) is canonically isomorphic to product
\begin{align}
    \prod_{\bar{v}\in
    X(\bar{k})}\Gr_{G_{\bar{v}}^\AD}^{-w_0(\lambda_{\bar{v},\AD})}(\bar{k}).
\end{align}
It is clear that \(\OGASch[\circ]_X\) is open and fiberwise dense in \(\OGASch_X\)
over \(\BD_X\).

\begin{remark}
    If the reader compares the definition in \cite{Ya04} with the one we give
    here, one can see the reason why the traditional definition gets so
    complicated: when we define the universal monoid \(\Env(G^\SC)\) (cf.,
    \Cref{sec:review_of_very_flat_reductive_monoids}), we are not only using the
    representations of \(G^\SC\), but also the abelianization space
    \(\FRA_{\Env(G^\SC)}\cong\bbA^r\); on the other hand, the traditional
    definition only utilizes the representation part, thus certain ``integrality''
    condition is lost and has to be rebuilt using \textit{ad hoc} formulations.
\end{remark}

\subsection{}
We end this section by discussing two variants of affine Schubert schemes.
The first is the equivariant version, namely the Hecke stacks. The arc group
\(\Arc_{\BD_X}G^\AD\) acts on \(\OGASch_X\) by left translation, therefore we
want to consider the ``quotient stack''
\begin{align}
    \Stack*{\OGASch_X}\defeq \Stack*{\Arc_{\BD_X}G^\AD\backslash \OGASch_X}.
    \nomenclature[\(Q_X_wow_stack \)]{\(\Stack*{\OGASch_X},\Stack*{\OGASch_X}_{G}\)}{the
        global Hecke or \(G\)-Hecke stack associated with a monoid \(\FRM\in\FM(G^\SC)\)}
\end{align}
The most obvious definition for this stack is that its points are pairs
\((E,\phi)\) where \(E\) is an fpqc torsor under
\(\Arc_{\BD_X}G^\AD\) and \(\phi\) is an \(\Arc_{\BD_X}G^\AD\)-equivariant map
from \(E\) to \(\OGASch_X\), and the morphisms between pairs \((E,\phi)\) and
\((E',\phi')\) are defined in an obvious way. Since the arc group
\(\Arc_{\BD_X}G^\AD\) is affine of infinite type, this is indeed a stack by fpqc
descent, but not algebraic in any way because the automorphism groups are
infinite dimensional. However, for non-algebraic stacks in general, many
geometric tools are lacking in development, so we
need to find an alternative definition. It is possible because our stack is not
that far from being algebraic.

\subsection{}
If we fix any connected component of \(\BD_X\), the action of
\(\Arc_{\BD_X}G^\AD\) on \(\OGASch_X\) factors through jet group
\(\Arc_{\BD_X,N}G^\AD\) for some large \(N\) (by the same argument as
\Cref{lem:approximation_in_arc_monoid}). On the said component, we may
consider the usual quotient stack
\begin{align}
    \Stack*{\OGASch_X}_N\defeq \Stack*{\Arc_{\BD_X,N}G^\AD\backslash \OGASch_X},
    \nomenclature[\(Q_X_wow_stack_N \)]{\(\Stack*{\OGASch_X}_N,\Stack*{\OGASch_X}_{G,N}\)}{the
        truncated global Hecke or \(G\)-Hecke
    stack associated with a monoid \(\FRM\in\FM(G^\SC)\) at level \(N\)}
\end{align}
and define
\begin{align}
    \Stack*{\OGASch_X}'\defeq\varprojlim_N\Stack*{\OGASch_X}_N,
\end{align}
where the limit simply means taking the pro-object (or the limit in the
pro-category of the category of algebraic stacks). Note that the limit
definition makes sense over all \(\BD_X\).

Given an fpqc \(\Arc_{\BD_X}G^\AD\)-torsor \(E\), it naturally induces torsors
\(E_N\) of \(\Arc_{\BD_X,N}G^\AD\). Since jet groups are of finite type and
smooth, fpqc torsors can be trivialized \'etale-locally, hence an fpqc torsor
under a jet group is also an \'etale torsor. This way we have a well-defined map
by universal property of pro-objects
\begin{align}
    \label{eqn:map_between_two_defs_of_equiv_affine_Schubert}
    \Stack*{\OGASch_X}\longto \Stack*{\OGASch_X}'.
\end{align}
On the other hand, given a compatible system of \(E_N\) for all sufficiently
large \(N\), the limit \(E\) as a presheaf is automatically an fpqc sheaf
because each \(E_N\) is. Since the arc group is the limit of jet groups, \(E\)
also admits an action of \(\Arc_{\BD_X}G^\AD\), and this action is
simply transitive because it can be verified over strict Henselizations (where
\(E_N\) become trivial torsors). After that, one can easily verify that
\eqref{eqn:map_between_two_defs_of_equiv_affine_Schubert} is an equivalence.

\subsection{}
Replacing \(\Arc_{\BD_X}G^\AD\) by
\(\Arc_{\BD_X}G\), we may also define pro-algebraic stack
\begin{align}
    \Stack*{\OGASch_X}_G\defeq \Stack*{\Arc_{\BD_X}G\backslash \OGASch_X},
\end{align}
and similarly for the truncated version \(\Stack*{\OGASch_X}_{G,N}\).

\begin{definition}
    \label[definition]{def:equivariant_global_affine_Schubert}
    The pro-algebraic stack \(\Stack*{\OGASch_X}\)
    (resp.~\(\Stack*{\OGASch_X}_G\)) is called the \notion{global Hecke
    (resp.~\(G\)-Hecke) stack associated with \(\FRM\)}\index{Hecke stack!global \(G\)-}
    \index{Hecke stack!global}\index{global!\(G\)-Hecke stack}\index{global!Hecke stack}. Similarly, the
    algebraic stack
    \(\Stack*{\OGASch_X}_N\) (resp.~\(\Stack*{\OGASch_X}_{G,N}\)) is called the
    \notion{\(N\)-truncated global Hecke (resp.~\(G\)-Hecke)
    stack associated with \(\FRM\)}\index{Hecke stack!\(N\)-truncated global \(G\)-}
    \index{Hecke stack!\(N\)-truncated global}\index{global!\(G\)-Hecke stack, \(N\)-truncated}
    \index{global!Hecke stack, \(N\)-truncated}.
\end{definition}

\begin{remark}
    \begin{enumerate}
        \item For any fixed \(N\), the \(N\)-truncated version is only defined
            for finitely many components of \(\BD_X\). However, since the exact
            value of \(N\) is not very important, we may pretend that
            \(\Stack*{\OGASch_X}_N\) is defined over all \(\BD_X\) with \(N\)
            growing suitably with the component of \(\BD_X\).
        \item The stack \(\Stack*{\OGASch_X}_G\) is a globalized version of
            \(\Hk_{\bG,\bar{v}}^{\le\lambda}\) used in
            \Cref{sec:connection_with_mv_cycles}, where the latter is obtained
            from an unbounded local Hecke stack \(\Hk_{\bG,\bar{v}}\) by
            truncating at \(\lambda\). We will, however, keep using the notation
            \(\Stack*{\OGASch_X}_G\) to emphasize the fact that we do not need
            to manually specify \(\lambda\), because it is automatically bounded
            by its construction using the monoid.
    \end{enumerate}
\end{remark}

\subsection{}
\label{sub:partial_affine_Schubert}
The second variant of affine Schubert scheme is the partial version.
Let \(U\) be an algebraic stack over \(k\) and \(U\to\BD_X\) be a map, then we
have induced affine Schubert scheme over \(U\):
\begin{align}
    \OGASch_{X,U}=\OGASch_X\x_{\BD_X} U.
\end{align}
The support of boundary divisors \(\lambda_U\)
parametrized by \(U\) can be viewed as a purely codimension \(1\) subscheme of
\(U\x X\). The union of any subset of connected components of such subscheme can be seen
as a \(U\)-flat family of boundary subdivisors of \(\lambda_U\). Call this
subdivisor \(\lambda_U'\), then since the arc scheme construction is local, we
may also define a \(U\)-family of affine Schubert schemes
\begin{align}
    \label{eqn:partial_affine_Schubert_def}
    \OGASch[\lambda_U']_{X,U}\longto U,
\end{align}
whose fiber at \(u\in U(\bar{k})\) is just the product of the direct factors 
of \(\Gr_{G^\AD}^{\le-w_0(\lambda_{U,u})}\) supported on \(\lambda_U'\).
\begin{definition}
    \label[definition]{def:partial_affine_Schubert_def}
    The map \eqref{eqn:partial_affine_Schubert_def} is called the
    \notion{partial affine Schubert scheme supported on
    \(\lambda_U'\)}\index{affine!Schubert scheme, partial}.
\end{definition}
If the arc space of \(\FRM\) is defined over \(U\), then we have natural
projections fitting into commutative diagram
\begin{equation}
    \begin{tikzcd}
        \Arc_{U}(\FRM/\FRA_{\FRM}) \ar[r]\ar[d] & \Arc_{U,\lambda_U'}(\FRM/\FRA_{\FRM}) \ar[d]\\
        \OGASch_{X,U} \ar[r] & \OGASch[\lambda_U']_{X,U}
    \end{tikzcd}
\end{equation}
where \(\Arc_{U,\lambda_U'}\) means the arc space but only consider the discs
centered at subdivisor \(\lambda_U'\). The vertical maps of this diagram are
formally smooth.

\subsection{}
The idea behind partial affine Schubert schemes is the same as the factorization
property of Beilinson--Drinfeld Grassmannians. However, the stack \(U\) and
subdivisor \(\lambda_U'\) in \Cref{def:partial_affine_Schubert_def} is quite
arbitrary and not very convenient to use. Therefore, it would be beneficial if we
can reformulate it as some sort of factorization result like the Grassmannian.
This will be the content of the next section, and it will play an important role
particularly in \Cref{sec:deformation_of_mHiggs_bundles}.

\section{Factorizations} 
\label{sec:factorizations}

It is well-known that Beilinson--Drinfeld Grassmannians satisfy certain
factorization properties. It roughly says that if \(D\) is a divisor on \(X\)
which is the disjoint sum of two subdivisors \(D'\) and \(D''\), then the fiber
of the Beilinson--Drinfeld Grassmannian at divisor \(D\) is canonically
isomorphic to the direct product of the respective fibers at \(D'\) and at \(D''\).
Moreover, the same also holds if \(D\), \(D'\) and \(D''\) move in flat families
as long as \(D'\) and \(D''\) are always disjoint. This allows us to separate
local studies at different points, and facilitates various kinds of inductive
arguments in studying representation theory through specialization to
when \(D'\) and \(D''\) collide. In this section we study the analogue
for boundary divisors.

\subsection{}
We use the notations in \Cref{sec:boundary_divisors}. Let \(\FRA\) be an
\'etale-locally constant toric scheme over \(X\) with affine and normal fibers,
and \(A\) the corresponding torus. Given any finite set \(I\), we have the
\(I\)-fold direct product of \(\FRA\), denoted by \(\FRA^I\). The
multiplication
\begin{align}
    \FRA^I&\longto \FRA\\
    (a_i)_{i\in I}&\longmapsto \prod_{i\in I}a_i
\end{align}
induces a map of quotient stacks
\begin{align}
    \label{eqn:boundary_quotient_mult_map}
    \Stack{\FRA^I/A^I}\simeq \Stack{\FRA/A}^I\longto \Stack{\FRA/A}.
\end{align}
If \(I=\emptyset\), we use the convention that \(\FRA^I\) is a singleton and its
image is \(1\).
Since \(\FRA^I\) is again an \'etale-locally constant toric scheme with affine
and normal fibers corresponding to torus \(A^I\), we have the moduli of boundary
divisors \(\BD_{X,\FRA^I}\), which is an open substack of the mapping stack
\begin{align}
    \BD[+]_{X,\FRA^I}=\IHom_{X/k}(X,\Stack{\FRA^I/A^I})\simeq\IHom_{X/k}(X,\Stack{\FRA/A})^I=(\BD[+]_{X,\FRA})^I.
\end{align}
It also identifies \(\BD_{X,\FRA^I}\) with \((\BD_{X,\FRA})^I\).
The multiplication map \eqref{eqn:boundary_quotient_mult_map} induces a map of
mapping stacks
\begin{align}
    \BD[+]_{X,\FRA^I}\longto \BD[+]_{X,\FRA},
\end{align}
and the image of \(\BD_{X,\FRA^I}\) is contained in \(\BD_{X,\FRA}\). As a
result, we have the \notion{summation map}\index{map!summation, of boundary
divisors} for boundary divisors
\begin{align}
    \Sigma\colon\BD_{X,\FRA^I}\simeq (\BD_{X,\FRA})^I&\longto\BD_{X,\FRA}.
\end{align}
\begin{example}
    When \(\FRA=\bbA^1\) and \(I=\Set{1,2}\), the summation map is none other
    than the usual summation of effective Cartier divisors
    \((D_1,D_2)\mapsto D_1+D_2\). At the level of moduli spaces, it is the
    disjoint union of maps (over all \(d_1,d_2\in\bbN\))
    \begin{align}
        X_{d_1}\x X_{d_2}\longto X_{d_1+d_2}.
    \end{align}
\end{example}

\subsection{}
Let \(\FRM\in\FM(G^\SC)\)
and \(I\) be a non-empty finite set.
We may pull back \(\FRM\) through multiplication map
\begin{align}
    \FRA_\FRM^I\longto \FRA_\FRM,
\end{align}
and denote the resulting monoid \(\FRM^I\). To be clear, \(\FRM^I\) is
\emph{not} the \(I\)-fold direct product of \(\FRM\). We will also need the
monoid \(\FRM^\AD=\FRM\git Z^\SC\in\FM(G^\AD)\) and
\(\FRM^{\AD,I}=\FRM^I\git Z^\SC\). The abelianization of \(\FRM^\AD\) (resp.~\(\FRM^{\AD,I}\)) is still
\(\FRA_\FRM\) (resp.~\(\FRA_\FRM^I\)).

Since \(\FRA_\FRM^\x\) (resp.~\(\FRA_\FRM^{I,\x}\)) is canonically isomorphic to
\(Z_{\FRM^\AD}\) (resp.~\(Z_{\FRM^{\AD,I}}\)), we have Cartesian diagram
\begin{equation}
    \begin{tikzcd}
        \Stack*{\FRM^{\AD,I}/\FRA_\FRM^{\x,I}} \ar[r]\ar[d] &
        \Stack*{\FRM^\AD/\FRA_\FRM^\x} \ar[d]\\
        \Stack*{\FRA_\FRM^{I}/\FRA_\FRM^{I,\x}} \ar[r] &
        \Stack*{\FRA_\FRM/\FRA_\FRM^\x}
    \end{tikzcd}
\end{equation}
where the bottom horizontal map is the multiplication map.

Any element \(i\in I\) defines a projection
\begin{align}
    \pr_i\colon\Stack*{\FRA_\FRM^I/\FRA_\FRM^{I,\x}}\simeq
    \Stack*{\FRM_\FRM/\FRA_\FRM^\x}^I\longto
    \Stack*{\FRA_\FRM/\FRA_\FRM^\x}.
\end{align}
For simplicity, we denote \(\BD_{X,\FRA_\FRM}\) by \(\BD_X\) and
\(\BD_{X,\FRA_\FRM^I}\) by \(\BD[I]_X\).
There are \(\abs{I}+1\) different affine Schubert schemes
over \(\BD[I]_X\): one is the pullback of \(\OGASch_X\)
through the multiplication map, which we denote by \(\OGASch_I\), and the others
are the pullbacks through projection maps \(\pr_i\) (\(i\in I\)), which we
denote by \(\OGASch_i\). For any point \(b_I=(b_i)_{i\in I}\), each \(b_i\) induces
an effective Cartier divisor \(\FRE_{b_i}\) on \(X\), whose support is denoted
by \(\supp(b_i)\). Let \(\BD[I,\disj]_X\subset \BD[I]_X\)
    \nomenclature[\(B"cal_X_wow_I_disj \)]{\(\BD[I,\disj]_X\)}{the locus in the
    \(I\)-fold (\(I\) a finite set) product of \(\BD_X\) such that the support
    of each \(b_i\in\BD_X\) (\(i\in I\)) is pairwise disjoint}
be the open
substacks such that \(\supp(b_i)\) are pairwise disjoint. We have the following
factorization result:
\begin{proposition}
    \label[proposition]{prop:factorization_of_affine_Schubert}
    We have canonical isomorphism over \(\BD[I,\disj]_X\)
    \begin{align}
        \label{eqn:global_affine_schubert_scheme_factorization}
        f_I\colon \Bigl.\Bigl(\prod_{i\in
        I}\OGASch_i\Bigr)\Bigr|_{\BD[I,\disj]_X}\stackrel{\sim}{\longto}\OGASch_I|_{\BD[I,\disj]_X}.
    \end{align}
    Moreover, if \(\phi\colon J\rightarrow I\) is a surjective map of finite
    sets, and \(J_i=\phi^{-1}(i)\), then \(\BD[J,\disj]_X\subset\prod_{i\in
    I}\BD[J_i,\disj]_X\) and \(f_J=\prod_{i\in I}\pr_i^*(f_{J_i})\) over
    \(\BD[J,\disj]_X\).
\end{proposition}
\begin{proof}
    The proof is essentially the same as in the case of Beilinson--Drinfeld
    Grassmannian, except using the language of monoids. Indeed, fix an arbitrary
    total order on \(I\), and we let \(\FRM^{\AD,\otimes I}\) be the \(I\)-fold
    direct product of \(\FRM^\AD\) over \(X\) with the given order, then we have a
    natural map \(\FRM^{\AD,\otimes I}\to\FRM^{\AD,I}\) compatible with the multiplication
    map \(\FRM^{\AD,\otimes I}\to\FRM^\AD\). We also have map
    \begin{align}
        \Stack*{\FRM^{\AD,\otimes I}/\FRA_\FRM^{I,\x}}\longto
        \Stack*{\FRA_\FRM^I/\FRA_\FRM^{I,\x}}
    \end{align}
    which factors through \(\Stack{\FRM^{\AD,I}/\FRA_\FRM^{I,\x}}\). Taking the mapping
    stacks, then for any \(b_I\in\BD[I]_X(S)\), the map between
    respective fibers over \(b_I\) is
    \begin{align}
        \label{eqn:factorization_prod_to_orig_over_b}
        \prod_{i\in I}\FRM^\AD_{b_i}\longto\FRM^{\AD,I}_{b_I}.
    \end{align}
    The left-hand side of
    \eqref{eqn:global_affine_schubert_scheme_factorization} at \(b_I\) is by
    definition the quotient space
    \begin{align}
        \label{eqn:factorization_prod_of_affine_Schubert}
        \prod_{i\in I}\Arc_{b_i}\FRM_{b_i}/\Arc_{b_i}G^\SC,
    \end{align}
    while the right-hand side is
    \begin{align}
        \label{eqn:factorization_orig_affine_Schubert}
        \Arc_{b_I}\FRM^I_{b_I}/\Arc_{b_I}G^\SC.
    \end{align}
    Again, in \eqref{eqn:factorization_prod_of_affine_Schubert}
    (resp.~\eqref{eqn:factorization_orig_affine_Schubert}) the arc space
    \(\Arc_{b_i}\FRM_{b_i}\) (resp.~\(\Arc_{b_I}\FRM^I_{b_I}\)) is not
    well-defined, but the quotient is.

    If \(b_I\) is in addition contained \(\BD[I,\disj]_X\), then
    \(\prod_{i\in I}\Arc_{b_i}G^\SC\simeq \Arc_{b_I}G^\SC\), and
    the map \eqref{eqn:factorization_prod_to_orig_over_b} descends to a map from
    \eqref{eqn:factorization_prod_of_affine_Schubert} to
    \eqref{eqn:factorization_orig_affine_Schubert}. Checking fiber by fiber
    it is not hard to see that it is an isomorphism. The rest of the proposition
    is very straightforward and also very similar to the case of Beilinson--Drinfeld
    Grassmannian, so we leave it to the reader.
\end{proof}

\begin{remark}
    In general, there is no map from
    \eqref{eqn:factorization_prod_of_affine_Schubert} to
    \eqref{eqn:factorization_orig_affine_Schubert} due to overlapping supports
    between \(b_i\), and the fact that the right (or left) translation of
    \((G^\AD)^I\) on \(\FRM^{\AD,\otimes I}\) is not compatible with multiplication
    map. However, just like affine Grassmannians, we can use an alternative
    action of \((G^\AD)^I\) to form the ``convolution Schubert
    schemes''. Namely, we let every \(G^\AD\)-factor except for the last one to
    act between consecutive \(\FRM^\AD\)-factors, and the last \(G^\AD\)-factor act
    on the very right-end. This is the analogue of the convolution
    Grassmannian. We omit the details.
\end{remark}

\subsection{}
The factorization still works if we replace \(\OGASch_X\) by the Hecke stack
\(\Stack*{\OGASch_X}\) or the truncated version \(\Stack*{\OGASch_X}_N\).

\subsection{}
Consider the special case when \(I=\Set{1,2}\), and let \(U=(\BD_X\x\BD_X)^{\disj}\).
Then over \(U\) we have the partial affine Schubert scheme by letting
\(\lambda_U'\) be equal to the subdivisor given by \(b_1\), which is then
isomorphic to \(\OGASch_1\). Moreover, it is straightforward to see that
\((U,\lambda_U')\) is actually universal among all moduli of partial affine
Schubert schemes in the sense that any other pair \((V,\lambda_V')\) is a unique
pullback of \((U,\lambda_U')\).


\section{Perverse Sheaves}%
\label{sec:perverse_sheaves}

One important aspect about global affine Grassmannian is the category of 
equivariant perverse sheaves on it. In the single-point case, i.e., the
``usual'' affine Grassmannian case, it is well-known that
it is equivalent to the Tannakian category of the dual group \(\dual{\bG}\)
through geometric Satake equivalence. The affine Schubert varieties
characterize all the simple objects in the that category corresponding to
irreducible representations.
The results are similar in the case of Beilinson--Drinfeld Grassmannians, and we briefly
review it. See \cite{Zh17}*{\S\S~4, 5} for details.

\subsection{}
First of all, it is more convenient to start with de-symmetrized version of
Beilinson--Drinfeld
Grassmannians. In other words, instead of symmetric power \(X_d\), we
consider direct product \(X^d\) and \(\Gr_{G,X^d}\) on it. The factorization
property of Beilinson--Drinfeld Grassmannian states that over the multiplicity-free locus
\((X^d)^\disj\), we have canonical isomorphism
\begin{align}
    \Gr_{G,X^d}\x_{X^d} (X^d)^\disj\simeq
    \left(\prod_{i=1}^d\Gr_{G,X}\right)\x_{X^d}(X^d)^\disj.
\end{align}
We also have the
\(d\)-fold convolution affine Grassmannian
\begin{align}
    \Gr_{G,X}^{\tilde{\x}d}\defeq \Gr_{G,X}\tilde{\x}\cdots\tilde{\x}\Gr_{G,X}
\end{align}
whose \(S\)-points consist of tuples 
\begin{align}
    (x_1,\ldots,x_d; E_1, \ldots, E_d, \phi_1,\ldots,\phi_d),
\end{align}
where \(x_i\in X(S)\), \(E_i\in\Bun_G(S)\), and
\begin{align}
    \phi_i\colon E_i|_{X\x S-x_i}\stackrel{\sim}{\longto} E_{i-1}|_{X\x S-x_i}
\end{align}
is an isomorphism. Here again \(E_0\) denotes the trivial \(G\)-torsor. We have
the convolution map
\begin{align}
    m_d\colon \Gr_{G,X}^{\tilde{\x}d}&\longto \Gr_{G,X^d}\\
    (x_i,E_i,\phi_i)&\longmapsto
    (x_1,\ldots,x_d,E_d,\phi_1\circ\cdots\circ\phi_d),
\end{align}
which is also known to be an isomorphism over \((X^d)^\disj\).
We may factor \(m_d\) into smaller steps
\begin{align}
    \Gr_{G,X}^{\tilde{\x}d}\longto\Gr_{G,X^2}\tilde{\x}\Gr_{G,X}^{\tilde{\x}(d-2)}
    \longto\cdots\longto \Gr_{G,X^{d-1}}\tilde{\x}\Gr_{G,X}\longto \Gr_{G,X^d}.
\end{align}
Over any geometric point \((\bar{v}_i)\in X^d\), each step above is 
stratified semi-small in the sense of
\cite{MiVi07}*{p.~14} (to see this, one only need to combine factorization
property and the well-known fact
that over a single geometric point \(\bar{v}\), the convolution map
\(\Gr_{G,\bar{v}}\tilde{\x}\Gr_{G,\bar{v}}\to\Gr_{G,\bar{v}}\) is stratified semi-small). Over the
whole base \(X^d\), each step above is in fact small because it is an
isomorphism over the open dense subset \((X^d)^\disj\), hence making the
inequality in the definition of semi-smallness strict.

\subsection{}
Let \(F_i\) (\(1\le i \le d\)) be a \(\Arc{G}\)-equivariant perverse sheaf
on \(\Gr_{G,X}\). Then there is the notion of twisted external product
\(F_1\tilde{\boxtimes}\cdots\tilde{\boxtimes}F_d\)
on \(\Gr_{G,X}^{\tilde{\x}d}\) whose restriction to
\begin{align}
    \Gr_{G,X}^{\tilde{\x}d}\x_{X^d}(X^d)^\disj\simeq 
    \left(\prod_{i=1}^d\Gr_{G,X}\right)\x_{X^d}(X^d)^\disj
\end{align}
may be identified with external product \(F_1\boxtimes\cdots\boxtimes F_d\).
The classical result is that
\begin{align}
    m_{d!}(F_1\tilde{\boxtimes}\cdots\tilde{\boxtimes}F_d)=
    m_{d*}(F_1\tilde{\boxtimes}\cdots\tilde{\boxtimes}F_d)=j_{!*}(F_1\boxtimes\cdots\boxtimes
    F_d|_{(X^d)^\disj})
\end{align}
where \(j\) is the inclusion \(\Gr_{G,X}^{\tilde{\x}d}\x_{X^d}(X^d)^\disj\to
\Gr_{G,X}^{\tilde{\x}d}\). Suppose that
\(F_1\tilde{\boxtimes}\cdots\tilde{\boxtimes}F_d\) is fiberwise perverse over
\(X^d\), then we know that \(j_{!*}(F_1\boxtimes\cdots\boxtimes
F_d|_{(X^d)^\disj})\) is fiberwise perverse over \(X^d\).

\subsection{}
The scheme \(\Gr_{G,X}\) is a locally trivial fibration over \(X\), and so we also have
locally trivial fibrations of various affine Schubert varieties which can also be
defined using reductive monoids in \Cref{sec:global_affine_schubert_scheme}. 
Suppose \(F_i\) is the intersection sheaf corresponding to such a locally
trivial fibration of affine Schubert varieties, then
\(F_1\tilde{\boxtimes}\cdots\tilde{\boxtimes}F_d\) is fiberwise perverse over
\(X^d\) because \(\Gr_{G,X}^{\tilde{\x}d}\) is locally trivial over \(X^d\). Thus, we
have that \(j_{!*}(F_1\boxtimes\cdots\boxtimes F_d|_{(X^d)^\disj})\) is
fiberwise perverse over \(X^d\). Using the Cartesian diagram
\begin{equation}
    \begin{tikzcd}
        \Gr_{G,X^d} \ar[r, "q_d"]\ar[d] & \Gr_{G,d} \ar[d]\\
        X^d \ar[r] & X_d
    \end{tikzcd}
\end{equation}
we have that the sheaf
\begin{align}
    q_{d*}j_{!*}(F_1\boxtimes\cdots\boxtimes
    F_d|_{(X^d)^\disj})^{\FRS_d}=j_{!*}\left[(q_d|_{(X^d)^\disj})_*(F_1\boxtimes\cdots\boxtimes
    F_d|_{(X^d)^\disj})\right]^{\FRS_d}
\end{align}
is fiberwise perverse over \(X_d\), where \(\FRS_d\)
    \nomenclature[\(S"frak_d \)]{\(\FRS_d\)}{the symmetric group of \(d\) elements}
is the symmetric group of
\(d\) elements, and the inclusion \(X_d^\disj\to X_d\) is still denoted by
\(j\). Thus, we have the following result:
\begin{proposition}
    \label[proposition]{prop:fiberwise_perversity_on_global_Schubert}
    Let \(\FRM\in\FM(G^\SC)\). Then the intersection complex
    \begin{align}
        \IC_{\OGASch_X}=\IC(\OGASch_X,\Qlb)
    \end{align}
    is fiberwise \(\Arc_{\BD_X}{G}\)-equivariant and perverse over
    \(\BD_X\). More generally, the same holds for any
    \(\Arc_{\BD_X}{G}\)-equivariant perverse sheaf on \(\BD_X\).
\end{proposition}
\begin{proof}
    If \(\FRA_{\FRM}\) is of standard type (see \Cref{def:toric_standard_type}), the result
    follows directly from the discussion above and the description of \(\BD_X\)
    as a disjoint union of direct products of symmetric
    powers of smooth curves.

    If \(\FRA_{\FRM}\) is not of standard type, then we still have a
    finite birational cover of \(\BD_X\) such that each irreducible component is
    cover by a direct product of symmetric power of smooth curves.
    Therefore, we still have the same result, since for a finite birational map
    \(f\colon \sX\to\sY\) of Deligne--Mumford stacks locally of finite type, we
    have that \(f_*\IC_{\sX}=\IC_{\sY}\). The generalization to other
    equivariant perverse sheaves is straightforward.
\end{proof}

\begin{definition}
    Let \(\Sat_X\)
    \nomenclature[\(S{}at_X \)]{\(\Sat_X\)}{the Satake category of semisimple
    \(\Arc_{\BD_X}{G}\)-equivariant perverse sheaves on \(\Gr_{G,\BD_X}\)}
    be the category of semisimple
    \(\Arc_{\BD_X}{G}\)-equivariant perverse sheaves on \(\Gr_{G,\BD_X}\). We
    call \(\cF\in\Sat_X\) a \notion{Satake
    sheaf}\index{sheaf!Satake}\index{Satake!sheaf}.
\end{definition}

\chapter{Multiplicative Hitchin Fibrations}%
\label{chap:multiplicative_hitchin_fibrations}

In this section, we study multiplicative Hitchin fibrations (mH-fibrations).
A particularly important result of this chapter is a local model of singularity in
\Cref{sec:local_model_of_singularities}, the proof
of which will be postponed to \Cref{chap:deformation} due to its technical
complexity.
The remaining part of this chapter is organized in a similar fashion as in
\cite{Ng10}*{\S~4}, and many proofs will be similar with small modifications.

\section{Constructions}%
\label{sec:Constructions}

The first fruitful study of some prototype of mH-fibrations is perhaps
\cite{HM02} in mathematical physics, although they have completely different
motivations and their formulation has some severe limitations (for example, they only
allow the curve to be either \(\Gm\) or an elliptic curve).
In algebraic setting, there have also been some earlier explorations of mH-fibrations.
In Frenkel and Ng\^o's visionary paper \cite{FN11}, there is a primitive definition
of mH-fibration with a hint towards a general one. It is later
carried out and studied in various later papers, notably
\cites{Bo15,Bo17,BoCh18,Ch22}, all focusing on split groups.

Compared to Hitchin fibrations for Lie algebras, the
mH-fibrations in those earlier papers come in several flavors.
Each of these variations has its technical advantages
and weaknesses, but they do not differ in a very essential way. 
Although it looked like a unified construction existed in the dark, 
it has never been brought to light. Here we give such
a construction with the help of the moduli of boundary divisors.
In fact, we will see starting from \Cref{sec:the_case_of_endoscopic_groups} that
such generalization is crucial in geometrizing the fundamental lemma.

\subsection{}%
The construction of the usual Hitchin fibration starts with a smooth,
geometrically connected, projective curve \(X\) over \(k\) with genus \(g_X\) and 
reductive group \(G\) over \(X\) given by a \(\Out(\bG)\)-torsor \(\OGT_G\). 
We will use \(\breve{X}\) to denote the base change of \(X\) to \(\bar{k}\).
The scaling action of \(\Gm\) on \(\FRg=\Lie(G)\) commutes with the
adjoint action of \(G\) hence induces an action on \(\FRc=\FRg\git G\). We
then have maps
\begin{align}
    \Stack{\FRg/G\x\Gm}\longto \Stack{\FRc/\Gm}\longto\BG{\Gm}.
\end{align}
They induce maps of \(k\)-mapping stacks
\begin{align}\label{eqn:Hitchin_Lie_alg_mapping_stacks}
    \IHom(X,\Stack{\FRg/G\x\Gm})\longto
    \IHom(X,\Stack{\FRc/\Gm})\longto \PicS{X}.
\end{align}
Fix any \(k\)-point in \(\PicS{X}\), in other words, a line bundle \(\cL\) on \(X\), the
fiber of \eqref{eqn:Hitchin_Lie_alg_mapping_stacks} over \(\cL\) is then the
Hitchin fibration \(\cM_\cL^{\La{g}}\to\cA_\cL^{\La{g}}\).

\subsection{}
In multiplicative setting, we use reductive monoids as the analogue of
Lie algebra \(\FRg\), and the translation action of the 
central subgroup of the monoid as the analogue of the \(\Gm\)-action on \(\FRg\).

Let \(\FRM\) be a reductive monoid in \(\FM(G^\SC)\) with abelianization
\(\FRA_{\FRM}\).
The translation action of central group \(Z_{\FRM}\) on \(\FRM\) induces maps
\begin{align}
    \Stack*{\FRM/G\x Z_{\FRM}}\longto \Stack*{\FRC_{\FRM}/Z_{\FRM}}\longto
    \Stack*{\FRA_{\FRM}/Z_{\FRM}}\longto
    \BG{Z_{\FRM}},
\end{align}
which further induce maps of mapping stacks over \(k\)
\begin{equation}\label{eqn:mH_monoid_mapping_stacks}
    \begin{tikzcd}
        \IHom(X,\Stack*{\FRM/G\x
        Z_{\FRM}})\ar[r]\ar[d,phantom,"\rotatebox{90}{\(\defeq\)}" description]&
        \IHom(X,\Stack*{\FRC_{\FRM}/Z_{\FRM}})\ar[r]\ar[d,phantom,"\rotatebox{90}{\(\defeq\)}"
        description]& \IHom(X,\Stack*{\FRA_{\FRM}/Z_{\FRM}})
        \ar[r]\ar[d,phantom,"\rotatebox{90}{\(\defeq\)}" description] &
        \Bun_{Z_{\FRM}}\ar[dd]\\
        \cM_X^+ & \cA_X^+ & \cB_X^+\ar[d] & \\
                & & \BD[+]_{X,\FRA_{\FRM}} \ar[r]& \Bun_{\FRA_{\FRM}^\x}
    \end{tikzcd}.
\end{equation}
For all practical purposes,
we would only consider the open substack \(\cB_X\subset \cB_X^+\) being the
preimage of the moduli \(\BD_{X,\FRA_{\FRM}}\) of boundary divisors.
Note that \(\cB_X\to\BD_{X,\FRA_{\FRM}}\) 
is always a \(Z^\SC\)-gerbe over its image, hence a 
Deligne--Mumford stack, proper and locally of finite type. 

\begin{definition}
    We make several important definitions.
    \begin{enumerate}
        \item The map \(h_X\colon\cM_X\to\cA_X\),
            \nomenclature[\(h_X \)]{\(h_X\)}{the universal multiplicative Hitchin
            fibration associated with a monoid \(\FRM\in\FM(G^\SC)\) and curve \(X\)}
            being the preimage of
            \(\cB_X\), is called the \notion{universal
            multiplicative Hitchin fibration}\index{multiplicative!Hitchin fibration, universal}
            \index{mH-!fibration}\index{fibration!mH-}\index{fibration!universal multimplicative Hitchin}
            \index{universal!multiplicative Hitchin fibration}
            associated with monoid \(\FRM\).
            We will use mH-fibration for short and omit the monoid if it is
            clear from the context.
        \item The stack \(\cM_X\)
            \nomenclature[\(M_X"cal \)]{\(\cM_X\)}{the mH-total stack of \(h_X\)}
            is called the
            \notion{mH-total stack}\index{mH-!total stack} and
            \(\cA_X\)
            \nomenclature[\(A"cal_X \)]{\(\cA_X\)}{the mH-base of \(h_X\)}
            is called the \notion{mH-base}\index{mH-!base}. If \(S\) is a
            \(k\)-scheme, then an \(S\)-point of \(\cM_X\) is called an
            \notion{mHiggs bundle}\notion{mHiggs bundle} on \(S\).
        \item The stack \(\cB_X\)
            \nomenclature[\(B"cal_X \)]{\(\cB_X\)}{the moduli stack of boundary divisors induced by a monoid \(\FRM\in\FM(G^\SC)\)}
            is (still) called the \notion{moduli of
            boundary divisors}\index{divisor!boundary, moduli
            of}\index{moduli!boundary divisors}, and its points are (still) called
            \notion{boundary divisors}\index{divisor!boundary}. For any \(b\in\cB_X\), the pullback
            \(h_b\colon\cM_b\to\cA_b\) of \(h_X\) to \(b\) is called a
            \notion{restricted mH-fibration with boundary divisor
            \(b\)}\index{mH-!fibration, restricted}.
        \item We sometimes also use \(\cL\) (resp.~\(\lambda_b\)) to
            denote the image of \(b\) in \(\Bun_{Z_\FRM}\)
            (resp.~\(\BD_{X,\FRA_\FRM}\)).
        \item For any \(k\)-scheme \(S\) and any \(a\in\cA_X(S)\), we call the
            pullback \(\cM_a\)
            \nomenclature[\(M_a"cal \)]{\(\cM_a\)}{the mH-fiber over \(a\in\cA_X\)}
            of \(\cM_X\) the
            \notion{mH-fiber}\index{mH-!fiber} over \(a\).
    \end{enumerate}
\end{definition}

\subsection{}
Let \(\FRM'\to\FRM\) be a morphism in \(\FM(G^\SC)\), that is, an excellent map
of reductive monoids. Then recall we have the following Cartesian diagram
\begin{equation}
    \begin{tikzcd}
        \Stack{\FRM'/G\x Z_{\FRM'}}\ar[r]\ar[d] &
        \Stack{\FRC_{\FRM'}/Z_{\FRM'}}\ar[r]\ar[d] &
        \Stack{\FRA_{\FRM'}/Z_{\FRM'}}\ar[d]\\
        \Stack{\FRM/G\x Z_{\FRM}}\ar[r] & \Stack{\FRC_{\FRM}/Z_{\FRM}}\ar[r]
                                        & \Stack{\FRA_{\FRM}/Z_{\FRM}}
    \end{tikzcd},
\end{equation}
which induces Cartesian diagram of mH-fibrations
\begin{equation}
    \label{eqn:pullback_diagram_for_mH_fibrations}
    \begin{tikzcd}
        \cM_X'\ar[r]\ar[d] & \cA_X'\ar[r]\ar[d] & \cB_X'\ar[d]\\
        \cM_X\ar[r] & \cA_X\ar[r] & \cB_X
    \end{tikzcd}.
\end{equation}
Therefore, studying mH-fibrations in general can largely be reduced to studying
the mH-fibration associated with the universal monoid \(\Env(G^\SC)\) using
\eqref{eqn:pullback_diagram_for_mH_fibrations}.

As another generalization, we may replace \(Z_{\FRM}\) by a smooth group \(Z\) of
multiplicative type together with a finite unramified
homomorphism \(Z\to Z_{\FRM}\).
The resulting mapping stacks are obtained by pulling back the
ones above via map \(\Bun_Z\to\Bun_{Z_{\FRM}}\), and so it largely reduces to
studying
the homomorphism of commutative group stacks \(\Bun_Z\to\Bun_{Z_{\FRM}}\), which
is relatively straightforward.

\subsection{}
The map \(\cA_X\to\cB_X\) is easy to describe using the isomorphism
\(\FRC_{\FRM}\simeq \FRA_{\FRM}\x \FRC\). Over boundary divisor
\(\lambda_X\in\cB_X(\bar{k})\) lying over a \(Z_{\FRM}\)-torsor \(\cL\), the fiber is
the vector space \(\RH^0(X,\FRC_\cL)\), where \(\FRC_\cL=\FRC\x^{Z_{\FRM}}\cL\). When
\(\FRM=\Env(G^\SC)\), \(\cL\)
is a \(T^\SC\)-torsor. Let \(\Wt\in\CharG(T^\SC)\) be a connected component
consisting of fundamental weights, then \(\cL\) induces a line bundle on \(\Wt\)
denoted by \(\Wt(\cL)\), such that 
\begin{align}
    \FRC_\cL\cong\bigoplus_\Wt p_{\Wt*}\Wt(\cL),
\end{align}
where \(p_\Wt\colon\Wt\to X\) is the natural map. Let \(d_\Wt\) be the degree of
\(\Wt\) over \(X\), thus we see that if \(\deg{\Wt(\cL)}> d_\Wt(2g_X-2)\) for all
\(\Wt\), then \(\FRC_\cL\) has no higher cohomology and
\begin{align}
    \dim_{\bar{k}}\RH^0(X,\FRC_\cL)=\sum_\Wt\deg{\Wt(\cL)}-rg_X+r.
\end{align}
Therefore, \(\cA_X\to \cB_X\) is a
vector bundle (of varying rank) 
for all but finitely many connected components of \(\cB_X\).
For a general monoid \(\FRM\in\FM(G^\SC)\), one may use
\eqref{eqn:pullback_diagram_for_mH_fibrations} to draw the same conclusion, or
alternatively use the following fact:
the \(Z_{\FRM}\)-torsor \(\cL\) induces a \(T^\SC\)-torsor using any excellent
morphism \(\FRM\to \Env(G^\SC)\), thus the line bundles
\(\Wt(\cL)\) still make sense and so on.

\subsection{}
Let us see the construction more explicitly in some special cases.
First, we consider the case where the group is split and the
abelianization \(\FRA_{\FRM}\) is of standard type.
Suppose \(\FRA_{\FRM}\) is an affine space with coordinates given by characters
\(e^{\theta_1}, \ldots, e^{\theta_m}\). 
Let \(\cL\) be a \(Z_{\FRM}\)-torsor such that \(\theta_i(\cL)\) has degree larger
than \(2g_X-2\) for all \(1\le i\le m\). Recall \(\FRM\) corresponds to a
homomorphism \(\phi_{\FRM}\colon Z_{\FRM}\to T^\SC\), which induces a \(T^\SC\)-torsor
\(\cL_T\) from \(\cL\). 
Assume that \(\Wt_i(\cL_T)\) also has degree larger than \(2g_X-2\) for
all \(1\le i\le r\). Denote the fiber of \eqref{eqn:mH_monoid_mapping_stacks}
over \(\cL\) by
\begin{align}
    h_\cL^+\colon \cM_\cL^+\longto \cA_\cL^+.
\end{align}
Here \(\cA_\cL^+\) is the vector space
\begin{align}
    \cB_\cL^+\oplus \cC_\cL\defeq 
    \bigoplus_{i=1}^m\RH^0(X,\theta_i(\cL))\oplus\bigoplus_{i=1}^r\RH^0(X,\Wt_i(\cL_T)),
\end{align}
where \(\cB_\cL^+\) is the first \(m\) summands. The subspace \(\cB_\cL\) is the
open locus where the section to \(\theta_i(\cL)\) is non-zero for all \(i\). The
fibration \(h_\cL\) is defined as the pullback of \(h_\cL^+\) to \(\cB_\cL\).
This is the version of mH-fibration that looks closest to the Hitchin fibration
in the Lie algebra case.

\subsection{}
Next, as another special case, 
let us explain here the connection between our current construction and the one
in \cite{FN11}. We consider the case where \(G=G^\SC\) and is split.
Let \(\cM_d\) to be the classifying stack of tuples \((D,E,\phi)\) where
\(D\in X_d\), \(E\in\Bun_G\), and \(\phi\colon E|_{X-D}\to E|_{X-D}\) is an
automorphism of \(G\)-torsor over \(X-D\). For any point \(v\in D\), we have a
well-defined relative position
\(\lambda_v=\Inv_v(\phi)\in\CoCharG(T)_+\) by choosing any
trivialization of \(E\) over the formal disc \(X_v\). Since the 
relative position 
can be arbitrarily large, \(\cM_d\) must be of infinite type.
To obtain a finite-type object, one has to put a restraint on
relative positions. The simplest way is to fix a dominant cocharacter
\(\lambda\in\CoCharG(T)_+\), and let \(\cM_{d,\lambda}\) to be the (closed) 
substack of tuples \((D,E,\phi)\) such that \(\lambda_v\le d_v\lambda\) for all
\(v\in D\), where \(d_v\) is the multiplicity of \(v\) in \(D\).

The mH-fibration with total stack
\(\cM_{d,\lambda}\) is then constructed as follows: fix a tuple \((D,E,\phi)\),
and recall that \(\FRC=G\git G\) is an affine \(r\)-space whose coordinates
are given by the traces of the \(r\) fundamental representations of \(G\).
Taking the traces \(\chi_i\) of \(\phi\), the restraint on relative positions means
that the \(\chi_i(\phi)\) will have poles bounded by the divisor
\(\Pair{\Wt_i}{\lambda}D\). This means that there is a map
\begin{align}
    h_{d,\lambda}\colon \cM_{d,\lambda}\longto\cA_{d,\lambda},
\end{align}
where \(\cA_{d,\lambda}\) is the line bundle over \(X_d\) whose fiber over \(D\)
is \(\bigoplus_{i=1}^r\RH^0(X,\cO_X(\Pair{\Wt_i}{\lambda}D))\). Let
\(\alpha_{d,\lambda}\) be the map \(\cM_{d,\lambda}\to X_d\) and \(h_D\colon
\cM_D\to\cA_D\) to be the fiber over \(D\in X_d\). Then
\(h_D\) is the restricted mH-fibration with boundary divisor \(\lambda\cdot D\).

\subsection{}%
More generally (still assuming \(G=G^\SC\) and split), 
let \(\ul{\lambda}=(\lambda_1,\ldots,\lambda_m)\) be a tuple 
of dominant cocharacters (allowing repetitions). 
Fix a tuple of positive integers \(\ul{d}=(d_1,\ldots,d_m)\) and let
\begin{align}
    X_{\ul{d}}\defeq X_{d_1}\x\cdots\x X_{d_m}.
\end{align}
The sum \(d=\abs{\ul{d}}\defeq d_1+\cdots+d_m\) is called the \emph{total
degree} of \(\ul{d}\). For \(\ul{D}\in X_{\ul{d}}\), we define \(D\) to be the
sum of divisors \(\abs{\ul{D}}\defeq D_1+\ldots+D_m\), so that \(d=\deg{D}\). We
also define \(\CoCharG(T)_+\)-valued divisor on \(X\)
\begin{align}
    \lambda_D=\ul{\lambda}\cdot\ul{D}\defeq \sum_{i=1}^m \lambda_i\cdot
    D_i.
\end{align}
Let \(\cM_{\ul{d}}\) be the classifying stack of tuples \((\ul{D},E,\phi)\)
where \(\ul{D}\in X_{\ul{d}}\), \(E\in\Bun_G\), and \(\phi\colon E|_{X-D}\to
E|_{X-D}\) is an automorphism of \(G\)-torsor. We also define
\(\cM_{\ul{d},\ul{\lambda}}\) to be the closed substack such that 
\(\Inv(\phi)\le \lambda_D\). We can also define the mH-fibration by taking the
trace of \(\phi\)
\begin{align}
    h_{\ul{d},\ul{\lambda}}\colon \cM_{\ul{d},\ul{\lambda}}\longto
    \cA_{\ul{d},\ul{\lambda}},
\end{align}
where \(\cA_{\ul{d},\ul{\lambda}}\) is the vector bundle over \(X_{\ul{d}}\)
whose fiber over \(\ul{D}\) is
\(\bigoplus_{i=1}^r\RH^0(X,\cO_X(\Pair{\Wt_i}{\lambda_D}))\). Again, we have the
restricted mH-fibration \(h_{\ul{D}}\) over \(\lambda_D\).

\subsection{}
Consider \(\FRM=\FRM(\ul{\lambda})\) as in \Cref{sub:different_monoids} where
each \(\lambda_i\in\CoCharG(T)_+\) for \(1\le i\le m\). Then the
maximal torus in \(\FRM\) is \(\Gm^m\x T\) and \(Z_{\FRM}=\Gm^m\) maps to \(T\)
via \(\ul{\lambda}\). In this case, the maps
\eqref{eqn:mH_monoid_mapping_stacks} can be extended into a commutative diagram
\begin{equation}\label{eqn:compare_different_mHitchin}
    \begin{tikzcd}
        \cM_{\ul{d},\ul{\lambda}}\ar[r]\ar[d] &
        \cA_{\ul{d},\ul{\lambda}}\ar[d]\ar[r] & X_{\ul{d}}
        \ar[d, "{D_i\mapsto (\cO(D_i),\sigma_{D_i}))}"] & \\
        \cM_X\ar[r] & \cA_X\ar[r] & \cB_X\ar[r] & (\PicS{X})^m
    \end{tikzcd},
\end{equation}
where \(\sigma_{D_i}\) is the canonical section of \(\cO(D_i)\).
Here the stack \(\cB_X\) classifies \(m\)-tuples
pairs \((\cL_i,s_i)\) where \(\cL_i\) is a line bundle and \(s_i\) a non-zero section
therein, and \(X_{\ul{d}}\) embeds as an open and closed 
subspace of pairs with \(\deg(\cL_i)=d_i\).
The squares in \eqref{eqn:compare_different_mHitchin} are easily shown to
be Cartesian. If we fix for each \(i\) a line bundle 
\(\cL_i\in\PicS(X)\) with degree \(d_i\), they assemble into a \(Z_{\FRM}\)-bundle
\(\cL\) by taking direct product, then we have pullback diagram
\begin{equation}
    \begin{tikzcd}
        \cM_\cL\ar[r]\ar[d] & \cA_\cL\ar[r]\ar[d] & \cB_\cL\ar[d]\\
        \cM_{\ul{d},\ul{\lambda}}\ar[r] & \cA_{\ul{d},\ul{\lambda}}\ar[r] & X_{\ul{d}}
    \end{tikzcd},
\end{equation}
with non-empty fibers being \(\Gm^m\)-torsors.
When \(\ul{\lambda}=\lambda\) is a single
cocharacter and degree \(d\) is fixed,
we recover the construction in \cite{FN11}. This is the case where
the monoid is a so-called \(L\)-monoid.

\section{Symmetry of mH-fibrations}%
\label{sec:Symmetry_of_mH_Fibrations}

The \(Z_{\FRM}\)-action on \(\FRC_{\FRM}\) lifts to an action on \(\FRJ_{\FRM}\) and
is compatible with the group scheme structure, therefore we have group scheme
\begin{align}
    \Stack*{\FRJ_{\FRM}/Z_{\FRM}}\longto \Stack*{\FRC_{\FRM}/Z_{\FRM}},
\end{align}
which further induces a group scheme
\begin{align}
    \FRJ_X\to X\x\cA_X
    \nomenclature[\(J"frak_X \)]{\(\FRJ_X\)}{the universal regular centralizer group scheme over \(X\x\cA_X\)}
\end{align}
through pulling back along the evaluation map
\(X\x\cA_X\to\Stack*{\FRC_{\FRM}/Z_{\FRM}}\). We define the relative
\notion{Picard stack}\index{Picard!stack, global}
\begin{align}
    p_X\colon\cP_X\defeq \PicS\bigl(\FRJ_X/(X\x\cA_X)/\cA_X\bigr)\longto \cA_X
    \nomenclature[\(P"cal_X \)]{\(\cP_X\)}{the total stack of \(p_X\)}
    \nomenclature[\(p_X \)]{\(p_X\)}{the relative Picard stack of \(\FRJ_X\)-torsors on \(X\) over \(\cA_X\)}
\end{align}
classifying \(\FRJ_X\)-torsors on \(X\) relative to \(\cA_X\).
This is a relative algebraic stack over \(\cA_X\). It has a natural group stack
structure because \(\FRJ_X\) is commutative, and naturally acts on \(\cM_X\)
relative to \(\cA_X\) induced by the canonical homomorphism \(\chi_{\FRM}^*\FRJ\to I\).
For any point \(a\in\cA_X(S)\), we use \(\cP_a\)
        \nomenclature[\(P"cal_a \)]{\(\cP_a\)}{the fiber of \(\cP_X\) over \(a\in\cA_X\)}
to denote the fiber of
\(\cP_X\) over \(a\), and \(\FRJ_a\) the pullback of \(\FRJ_X\) to \(X\x S\) through \(a\).

\begin{proposition}
    The relative Picard stack \(p_X\colon\cP_X\to\cA_X\) of \(\FRJ_X\)-torsors
    is smooth.
\end{proposition}
\begin{proof}
    This is because the obstruction space of deforming \(\FRJ_X\)-torsors
    is \(\RH^2(X, \Lie(\FRJ_X))=0\) since \(X\) is a curve (see
    \cite{Ng10}*{Proposition~4.3.5} and \cite{Ch22}*{Proposition~4.2.2}).
\end{proof}

\subsection{}
Similar to the local situation, we can consider the open subset \(\cM_X^\reg\)
(resp.~\(\cM_X^\circ\)) of \(\cM_X\) consisting of tuples \((\cL,E,\phi)\) such
that the image of the Higgs field \(\phi\) is contained in
\(\Stack{\FRM_\cL^\reg/G}\) (resp.~\(\Stack{\FRM_\cL^\circ/G}\)).
According to \Cref{cor:M_reg_is_union_of_gerbs}, we have:
\begin{proposition}
    The action of \(p_X\colon\cP_X\to\cA_X\) on
    \(h_X^\reg\colon\cM_X^\reg\to\cA_X\) has trivial stabilizers.
\end{proposition}

\section{Cameral Curves}%
\label{sec:cameral_curves}

\begin{definition}
    The \notion{universal cameral curve}\index{cameral!curve, universal}
    \index{universal!cameral curve} \(\tilde{\pi}\colon \tilde{X}\to X\x\cA_X\)
    \nomenclature[\(pi'tilde \)]{\(\tilde{\pi}\)}{the universal cameral cover over \(X\x \cA_X\)}
    is the pullback of the cameral cover
    \begin{align}
        \Stack*{\pi}\colon\Stack*{\FRT_{\FRM}/Z_{\FRM}}\longto
        \Stack*{\FRC_{\FRM}/Z_{\FRM}}
    \end{align}
    to \(X\x\cA_X\) via the evaluation map.
\end{definition}
\begin{definition}
    The \notion{universal discriminant divisor}\index{divisor!discriminant, universal}
    \index{universal!discriminant divisor} \(\FRD_X\)
    \nomenclature[\(D"frak_X \)]{\(\FRD_X\)}{the universal discriminant divisor over \(X\x\cA_X\)}
    is the pullback of
    \(\Stack{\FRD_{\FRM}/Z_{\FRM}}\) to \(X\x\cA_X\) via the evaluation map. 
    If the numerical boundary divisor \(\FRE_{\FRM}\) makes sense for the monoid, 
    we define \(\FRE_X\)
    \nomenclature[\(E"frak_X \)]{\(\FRE_X\)}{the universal numerical boundary divisor over \(X\x\cA_X\)}
    also by pullback.
\end{definition}
Fix \(a\in\cA_X(S)\), we have the \notion{cameral curve}\index{cameral!curve} \(\pi_a\colon
\tilde{X}_a\to X\x S\)
that is the fiber of \(\tilde{\pi}\colon\tilde{X}\to X\x\cA_X\) over \(a\).
We also have the \notion{(extended) discriminant
divisor}\index{divisor!discriminant, extended} \(\FRD_a\)
    \nomenclature[\(D"frak_a \)]{\(\FRD_a\)}{the pullback of \(\FRD_X\) to \(X\) via \(a\in\cA_X\)}
and the
\notion{numerical boundary divisor}\index{divisor!numerical boundary} \(\FRE_a\)
    \nomenclature[\(E"frak_a \)]{\(\FRE_a\)}{the pullback of \(\FRE_X\) to \(X\) via \(a\in\cA_X\)}
by looking at the fibers of \(\FRD_X\) and
\(\FRE_X\) over \(a\), respectively. Note that despite the name, \(\FRD_a\) may
not be a (proper) divisor.
\begin{definition}
    \label[definition]{def:mH_fibration_reduced_locus}
    Let \notion{reduced locus}\index{locus!reduced} 
    \(\cA_X^\heartsuit\subset \cA_X\)
    \nomenclature[\(.heartsuit \)]{\((\cdot)^\heartsuit\)}{related to the locus where the cameral curve is reduced}
    be the open locus such that
    \(\FRD_a\) is either empty or an effective divisor.
\end{definition}

\begin{lemma}
    \label[lemma]{lem:discriminant_divisor_total_degree}
    Let \(a\in\cA_X^\heartsuit(\bar{k})\) and let \(b\) be its boundary divisor
    such that at each \(\bar{v}\in X(\bar{k})\) it gives a dominant cocharacter
    \(\lambda_{\bar{v}}\). Then
    \begin{align}
        \deg(\FRD_a)=\sum_{\bar{v}\in X(\bar{k})}\Pair{2\rho}{\lambda_{\bar{v}}}.
    \end{align}
\end{lemma}
\begin{proof}
    It suffices to prove for \(\FRM=\Env(G^\SC)\).
    The extended discriminant function is a map \(\FRC_\FRM\to \bbA^1\). It is
    \(Z_\FRM\)-equivariant if we let \(Z_\FRM\) act on \(\bbA^1\) by character
    \(2\rho\). Let \(\cL\) be the \(Z_\FRM\)-torsor under \(b\), then we
    have induced map
    \begin{align}
        \FRC_{\FRM,\cL}\longto \cO_{\breve{X}}(\FRD_{\FRM,\cL})
        =\cO_{\breve{X}}\Bigl(\sum_{\bar{v}}\Pair{2\rho}{\lambda_{\bar{v}}}\bar{v}\Bigr),
    \end{align}
    whose preimage of the zero section of \(\cO(\FRD_{\FRM,\cL})\) is the
    discriminant divisor. The point \(a\) is a map \(\breve{X}\to
    \FRC_{\FRM,\cL}\), and so its composition with the map above is a section of
    \(\cO(\FRD_{\FRM,\cL})\) whose zero divisor is exactly \(\FRD_a\).
\end{proof}

\begin{lemma}
    \label[lemma]{lem:reduced_cameral_curve}
    For \(a\in\cA_X^\heartsuit(\bar{k})\), 
    the curve \(\tilde{X}_a\) is reduced.
\end{lemma}
\begin{proof}
    Since the cameral cover is a Cohen--Macaulay morphism
    (cf.~\Cref{cor:cameral_cover_is_Cohen_Macaulay}), 
    \(\tilde{X}_a\) is Cohen--Macaulay. Since \(X_a\) is generically
    reduced being a finite flat  cover over \(X\), it must be reduced.
\end{proof}

\subsection{}%
Using the Galois description of regular centralizer \(\FRJ_{\FRM}\), we can reach
a similar description of \(\FRJ_X\) using cameral cover \(\tilde{\pi}\).
Indeed, for any \(a\in\cA_X(\bar{k})\), we have monomorphism of sheaf of commutative
groups
\begin{align}
    \FRJ_a\longto \FRJ_a^1=\pi_{a,*}(\tilde{X}_a\x_{\breve{X}} T)^W
\end{align}
with a finite cokernel of finite support relative to \(\breve{X}\). Sometimes it is
convenient to base-change so that \(G\) (and \(T\), etc.) becomes split.
After \cite{Ng10}*{\S~4.5.2}, let \(\OGT\colon X_\OGT\to \breve{X}\) be a finite 
Galois \'etale
cover with \(X_\OGT\) being connected and
over which \(\OGT_G\) becomes a trivial \(\Out(\bG)\)-torsor. Let
\(\Theta_\OGT\) be the Galois group. For any \(a\in\cA_X(\bar{k})\) whose image in
\(\Bun_{Z_{\FRM}}\) is \(\cL\), we have
Cartesian diagram
\begin{equation}
    \begin{tikzcd}
        \tilde{X}_{\OGT,a}\ar[r]\ar[d, "\pi_{\OGT,a}", swap] &
        X_\OGT\x_{\breve{X}}\FRT_{\FRM,\cL}\ar[d]\\
        \breve{X}\ar[r, "a"] & \FRC_{\FRM,\cL}
    \end{tikzcd}.
    \nomenclature[\(X'tilde_theta_a \)]{\(\tilde{X}_{\OGT,a}\)}{the fiber
    product \(\tilde{X}_a\x_X X_\OGT\) where \(X_\OGT\to X\) is a finite Galois cover over which \(\OGT_G\) becomes trivial}
    \nomenclature[\(pi_theta_a \)]{\(\pi_{\OGT,a}\)}{the splitting cameral cover \(\tilde{X}_{\OGT,a}\to X\)}
\end{equation}
Then we may describe \(\FRJ_a^1\) as
\begin{align}
    \label{eqn:global_J1_galois_desc}
    \FRJ_a^1=\pi_{\OGT,a,*}(\tilde{X}_{\OGT,a}\x\bT)^{\bW\rtimes \Theta_\OGT}.
\end{align}
Similarly, we have a global N\'eron model \(\FRJ_a^\flat\), and its Galois
description
\begin{align}
    \label{eqn:global_neron_galois_desc}
    \FRJ_a^\flat=\pi_{a,*}^\flat(\tilde{X}_a^\flat\x_{\breve{X}}
    T)^W=\pi_{\OGT,a,*}^\flat(\tilde{X}_{\OGT,a}^\flat\x\bT)^{\bW\rtimes \Theta_\OGT},
\end{align}
where \(\tilde{X}_a^\flat\) (resp.~\(\tilde{X}_{\OGT,a}^\flat\)) is the
normalization of \(\tilde{X}_a\) (resp.~\(\tilde{X}_{\OGT,a}\)).

\subsection{}
Using \(\tilde{X}_{\OGT,a}\), we can prove a connectivity result using the
same method as in \cite{Ng10}*{\S~4.6}. First we record a theorem.
\begin{theorem}[\cite{De96}*{Th\'eor\`eme~1.4}]
    \label[theorem]{thm:bertini_type_connectivity_theorem}
    Let \(M\) be an irreducible variety,
    \begin{align}
        m\colon M\to \bbP=\bbP^{n_1}\x\cdots\x\bbP^{n_r}
    \end{align}
    is a morphism from \(M\) into a product of several projective spaces.
    Let \(H_i\subset\bbP^{n_i}\) be a fixed linear subspace. Suppose for any
    subset \(I\subset \Set{1,\ldots,r}\), we have
    \begin{align}
        \dim(p_I(m(M)))>\sum_{i\in I}\codim_{\bbP^{n_i}}(H_i),
    \end{align}
    where \(p_I\) is the natural projection with multi-index \(I\). Suppose in
    addition there is an open subset \(V\subset \bbP\) containing the product
    \(H=H_1\x\cdots\x H_r\), and \(m^{-1}(V)\) is proper over \(V\). Then
    \(m^{-1}(H)\) is connected.
\end{theorem}

\begin{definition}
    For any integer \(N\ge 0\), a \(Z_{\FRM}\)-torsor \(\cL\) is called
    \inotion{very \((G,N)\)-ample} if
    \begin{align}
        \deg{\Wt(\cL)}>2d_\Wt(g_X-1)+2+N
    \end{align}
    for every connected component \(\Wt\) of
    fundamental weights of \(G^\SC\). If \(N=0\), we simply call \(\cL\)
    \inotion{very \(G\)-ample}. A point \(a\in\cA_X(\bar{k})\) or a
    boundary divisor \(b\in\cB_X(\bar{k})\) is very \((G,N)\)-ample
    (resp.~very \(G\)-ample) if its
    associated bundle \(\cL\) is. Let \(\cA_{\gg}\subset \cA_X\)
    \nomenclature[\(.gg \)]{\((\cdot)_\gg,(\cdot)_{\gg N}\)}{related to the locus
        where the boundary divisor is very \(G\)-ample (resp.~\((G,N)\)-ample)}
    and
    \(\cB_{\gg}\subset \cB_X\) be the respective loci of very \(G\)-ample
    points, and \(\cA_{\gg N}\), \(\cB_{\gg N}\) be the respective very
    \((G,N)\)-ample loci.
\end{definition}

\begin{remark}
    \label[remark]{rmk:G_ampleness_and_basechange}
    It is clear that \(\cB_{\gg}\) is both open and closed in \(\cB_X\).
    In some situation, a weaker condition on ampleness can be used, i.e.,
    one that requires \(\deg{\Wt(\cL)}>2d_\Wt(g_X-1)\) instead. We may call those
    points \emph{\(G\)-ample} but since obtaining a sharp numerical bound is
    not important, we will stick to very \(G\)-ampleness to save notations.
    Finally, note that if \(\OGT\colon X_\OGT\to X\) is a finite \'etale cover,
    then if \(\cL\) is \(G\)-ample or very \(G\)-ample, so is \(\OGT^*\cL\).
\end{remark}

\begin{proposition}
    \label[proposition]{prop:cameral_curve_connected}
    If \(a\in\cA_X^\heartsuit(\bar{k})\) is very \(G\)-ample,
    then both \(\tilde{X}_a\) and \(\tilde{X}_{\OGT,a}\) are reduced and 
    connected.
\end{proposition}
\begin{proof}
    Reducedness is already proved in \Cref{lem:reduced_cameral_curve}.
    For connectedness, 
    it suffices to prove for \(\tilde{X}_{\OGT,a}\) hence we may replace
    \(\breve{X}\) with \(X_\OGT\) and assume that \(G\) is split. Let
    \(b\in\cB_X\) and \(\cL\in\Bun_{Z_{\FRM}}\) be the image of \(a\). Consider the
    pullback diagram
    \begin{equation}
        \begin{tikzcd}
            \FRT_b\ar[r]\ar[d, "\pi_b"] & \FRT_{\FRM,\cL}\ar[d]\\
            \FRC_b\ar[r]\ar[d] & \FRC_{\FRM,\cL}\ar[d]\\
            \breve{X}\ar[r, "b"] & \FRA_{\FRM,\cL}
        \end{tikzcd}
    \end{equation}
    such that \(a\) is a section \(\breve{X}\to\FRC_b\). We claim that \(\FRT_b\) is an
    irreducible variety.

    Indeed, because \(b\) lies generically
    inside the open part \(\FRA_{\FRM,\cL}^\x\), there is an open dense subset
    \(U\subset \breve{X}\) such that the fibers of \(\FRT_b\to \breve{X}\) are
    \(T^\SC\)-torsors.  This means that \(\FRT_b\x_{\breve{X}} U\) is irreducible,
    hence so is its closure in \(\FRT_b\). The complement \(U'\) of the closure
    of \(\FRT_b\x_{\breve{X}} U\) is open, hence its image is open in \(\breve{X}\)
    by flatness. We know \(U'\) must be empty, since otherwise its image in
    \(\breve{X}\) would have
    non-trivial intersection with \(U\), which is impossible by the definition
    of \(U'\). Therefore, \(\FRT_b\) is irreducible. Similar to cameral curves,
    since \(a\in\cA_X^\heartsuit\), \(\FRT_b\) is Cohen--Macaulay
    and generically reduced, hence reduced. Thus, \(\FRT_b\) is an irreducible
    variety. 

    We can now apply
    \Cref{thm:bertini_type_connectivity_theorem} where \(M=\FRT_b\).
    Since \(G\) is split, \(\FRC_b\) is a direct sum of line bundles
    \(\Wt_i(\cL)\) for \(1\le i\le r\). Compactify \(\Wt_i(\cL)\) into a
    projective line bundle \(\bar{\Wt_i(\cL)}\), whose total space is a projective
    surface over \(\bar{k}\). Let \(Z_i=\bar{\Wt_i(\cL)}-\Wt_i(\cL)\) be the
    infinity divisor. Since by assumption \(\deg{\Wt_i(\cL)}>2g_X\), the line
    bundle \(\cO(1)\) on \(\bar{\Wt_i(\cL)}\) is very ample, hence induces a
    closed embedding
    \begin{align}
        \bar{\Wt_i(\cL)}\longto \bbP^{n_i}.
    \end{align}
    Let \(V\subset \bbP=\prod_{i=1}^r\bbP^{n_i}\) be the open subset
    \begin{align}
        \bbP-\bigcup_{i=1}^r\left(Z_i\x\prod_{j\neq i}\bar{\Wt_j(\cL)}\right),
    \end{align}
    then \(\prod_{i=1}^r\Wt_i(\cL)\) is a closed subscheme of  \(V\). Since
    \(\FRC_b\) is a closed subscheme of \(\prod_{i=1}^r\Wt_i(\cL)\), it is a
    closed subscheme of \(V\). Therefore, \(m\colon M\to V\) is proper since
    \(\pi_b\) is finite. The \(i\)-th component of \(a\) is a section
    \(a_i\in\RH^0(X,\Wt_i(\cL))\). Since
    \begin{align}
        (a_i,1)\in\RH^0(\bbP^{n_i},\cO(1))=\RH^0(X,\Wt_i(\cL))\oplus\RH^0(X,\cO_X),
    \end{align}
    it determines a hyperplane \(H_i\subset \bbP^{n_i}\)
    that has no intersection with \(Z_i\). Thus, \(H=\prod_{i=1}^r H_i\) is
    contained in \(V\).

    It is easy to see that for every \(I\subset\Set{1,\ldots,r}\),
    \(p_I(m(M))=p_I(\FRC_b)\) has dimension \(\abs{I}+1\) (being the dimension
    of the direct sum of line bundles \(\Wt_i(\cL)\) for \(i\in I\)), thus
    \begin{align}
        \dim(p_I(m(M)))>\sum_{i\in I}\codim_{\bbP^{n_i}}(H_i)=\abs{I}.
    \end{align}
    By \Cref{thm:bertini_type_connectivity_theorem} we are done as
    \(\tilde{X}_a=m^{-1}(H)\).
\end{proof}

\subsection{}
We consider open subset \(\cA_X^\diamondsuit\)
    \nomenclature[\(.diamondsuit \)]{\((\cdot)^\diamondsuit\)}{related to the
    locus in \(\cA_X\) where \(a(X)\) and \(\FRD_X\) intersects transversally}
of \(\cA_X\) consisting
of points \(a\) such that \(a(X)\) intersects with the discriminant divisor
transversally. We also consider another open subset
\(\cA_X^\sharp\subset\cA_X^\diamondsuit\)
    \nomenclature[\(.sharp \)]{\((\cdot)^\sharp\)}{related to the
    substack in \(\cA_X^\diamondsuit\) where \(a(X)\) does not intersect \(\FRD_X\) and \(\FRE_X\) simultaneously}
of points with the additional
condition that \(a(X)\) does
not intersect with \(\FRD_X\) and \(\FRC_X-\FRC_X^\x\) simultaneously.
Clearly both \(\cA_X^\diamondsuit\) and \(\cA_X^\sharp\) are subsets of
\(\cA_X^\heartsuit\).

\begin{proposition}
    \label[proposition]{prop:non_emptyness_of_A_sharp_diamond_heart}
    Suppose \(b\in\cB_X(\bar{k})\) is very \(G\)-ample, then
    \(\cA_b^\sharp\) is non-empty. As a result, both \(\cA_b^\diamondsuit\)
    and \(\cA_b^\heartsuit\) are also non-empty.
\end{proposition}
\begin{proof}
    The proof is completely parallel to that of \cite{Ng10}*{Proposition~4.7.1}.
    First we show that the discriminant divisor \(\FRD_b\subset\FRC_b\) is
    reduced. Indeed, as \(b\) lies generically (over \(X\)) inside
    \(\FRA_{\FRM,\cL}^\x\), there is an open dense subset \(U\) of \(X\) over
    which \(\FRD_b\) is reduced. Similar to the proof of reducedness of cameral
    curves (see also the proof of
    \Cref{lem:extended_discriminant_divisor_is_reduced}), we see that \(\FRD_b\)
    is Cohen--Macaulay and generically reduced,
    hence reduced.

    Next we show that for any \(x\in \breve{X}(\bar{k})\) with ideal
    \(\FRm_x\subset \cO_{\breve{X}}\), the map
    \begin{align}
        \label{eqn:cA_diamondsuit_surjectivity_at_one_point}
        \RH^0(X_{\bar{k}},\FRC_b)\longto
        \FRC_b\otimes_{\cO_{X_{\bar{k}}}}\cO_{X_{\bar{k}}}/\FRm_x^2
    \end{align}
    is surjective. Indeed, 
    let \(\OGT\colon X_\OGT\to \breve{X}\) be a connected finite Galois
    cover of \(\breve{X}\) with Galois group \(\Theta\). 
    Then \(\OGT^*\FRC_b\) is isomorphic to a direct sum of line
    bundles
    \begin{align}
        \OGT^*\FRC_b = \bigoplus_{i=1}^r \Wt_i(\OGT^*\cL).
    \end{align}
    Since \(\FRC_b\) is a direct summand of \(\OGT_*\OGT^*\FRC_b\) being the
    \(\Theta\)-fixed subbundle, it suffices to prove the surjectivity of map
    \begin{align}
        \RH^0(X_{\OGT},\OGT^*\FRC_b)\longto 
        \OGT^*\FRC_b\otimes_{\cO_{\breve{X}}}\cO_{\breve{X}}/\FRm_x^2,
    \end{align}
    and in turn it suffices to prove surjectivity after replacing
    \(\OGT^*\FRC_b\) by each \(\Wt_i(\OGT^*\cL)\). The numerical assumption on
    \(\Wt(\cL)\) ensures that \(\Wt_i(\OGT^*\cL)\) has degree greater than
    \(2g_{X_\OGT}\), hence the claim is implied by Riemann--Roch theorem.

    Now since \(\FRD_b\) is reduced, it has an open dense smooth locus
    \(\FRD_b^\diamondsuit=\FRD_b-\FRD_b^\Sg\) such that \(\FRD_b^\Sg\) has codimension
    \(2\) in \(\FRC_b\). Let
    \(\FRD_b^\sharp=\FRD_b^\diamondsuit-(\FRC_b-\FRC_b^\x)\),
    then \(\FRD_b-\FRD_b^\sharp\) still has codimension at most \(2\) in
    \(\FRC_b\) since \(\FRD_b\cap(\FRC_b-\FRC_b^\x)\) is so.
    Let \(Z_1\subset \FRD_b^\sharp\x\cA_b\) consisting
    of pairs \((c,a)\) such that \(a(X)\) passes \(c\), and has intersection
    multiplicity with \(\FRD_b\) at least \(2\) at \(c\). Fix any \(c\), then
    the subset of \(a\in\cA_b\) such that \((c,a)\in Z_1\) has codimension at
    least \(2r\) in \(\cA_b\) by the surjectivity of
    \eqref{eqn:cA_diamondsuit_surjectivity_at_one_point}. Hence,
    \begin{align}
        \dim{Z_1}\le \dim{\cA_b}-2r+\dim{\FRD_b}=\dim{\cA_b}-r-1\le \dim{\cA_b}-1.
    \end{align}
    Thus, the image of \(Z_1\) in \(\cA_b\) has codimension at least \(1\).
    Similarly, consider the subset \(Z_2\subset
    (\FRD_b-\FRD_b^\sharp)\x\cA_b\) of pairs
    \((c,a)\) such that \(a(X)\) passes \(c\). Then we also have that
    \begin{align}
        \dim{Z_2}\le\dim{\cA_b}-1.
    \end{align}
    Therefore, \(\cA_b-\bar{Z_1\cup Z_2}\subset \cA_b^\sharp\) is dense in
    \(\cA_b\) as desired.
\end{proof}

\begin{proposition}
    \label[proposition]{prop:A_diamondsuit_then_cameral_is_smooth}
    Let \(a\in \cA_X^\diamondsuit(\bar{k})\), then \(\tilde{X}_a\)
    is smooth.
\end{proposition}

We make some preparations before attempting to prove
\Cref{prop:A_diamondsuit_then_cameral_is_smooth}. 
Recall that in the absolute setting,
\(\bD_{\bM}\) is defined by the extended discriminant function
\begin{align}
    \Disc_+=e^{(2\rho,0)}\prod_{\alpha\in\Roots}(1-e^{(0,\alpha)}).
\end{align}
For each positive root \(\alpha\), we define a rational function
\begin{align}
    \Disc_\alpha= (1-e^{(0,\alpha)})(1-e^{(0,-\alpha)})
\end{align}
on \(\bar{\bT}_{\bM}\), which is a regular function on \(\bT_{\bM}\).
Let \(\bD_\alpha\) be the scheme-theoretic closure in \(\bar{\bT}_{\bM}\) of the
vanishing locus of \(\Disc_\alpha\).

\begin{lemma}\label[lemma]{lem:unique_root_for_d_equals_1}
    Let \(a\in\bC_{\bM}\) be a geometric point contained in the smooth locus
    \(\bD_{\bM}^\Sm\) of the discriminant divisor. Suppose in addition
    \(a\) is invertible, i.e., contained in \(\bC_{\bM}^\x\) and \(t\in\bT_{\bM}\) is
    a preimage of \(a\). Then there is a unique
    root \(\alpha\in\PosRts\) such that \(\Disc_\alpha(t)=0\) and
    \(\Disc_\beta(t)\neq 0\) for \(\beta\neq\alpha\).
\end{lemma}
\begin{proof}
    Since \(a\) is invertible, the condition \(a\in\bD_{\bM}^\Sm\) is not affected
    by changing monoid.
    Thus, it suffices to prove in the case \(\bM=\Env(\bG^\SC)\).

    Let \(a=\Spec{k(a)}\) and \(\cO(a)=k(a)\powser{\pi}\) be the ring of formal series
    with coefficients in \(k(a)\). 
    Since \(\bC_{\bM}\) is smooth, \(a\in\bD_{\bM}^\Sm\cap\bC_{\bM}^\x\) 
    implies that we can find a map 
    \begin{align}
        \tilde{a}\colon \Spec{\cO(a)}\longto \bC_{\bM}^\x
    \end{align}
    whose special point is \(a\), and \(\tilde{a}\) intersects with \(\bD_{\bM}\)
    transversally. This means that the \(\cO(a)\)-valuation of \(\bD_{\bM}\) at
    \(\tilde{a}\) is \(1\). By dimension formula of multiplicative affine Springer
    fibers, this forces the ramification index \(c(\tilde{a})=1\) (observe that
    the local dimension formula still hold if we replace \(\bar{k}\) by
    any algebraically closed field containing \(\bar{k}\)). This implies that
    \(\tilde{a}\) lifts to a point \(\tilde{t}\in\bT_{\bM}(k(a)\powser{\pi^{1/2}})\) 
    specializing to \(t\). 

    The \(\cO(a)\)-valuation of \(1-e^{(0,\alpha)}\) at \(\tilde{t}\)
    is contained in \(\bbN/2\) for
    all \(\alpha\in\Roots\). Since the valuation of \(1-e^{(0,\alpha)}\) equal to
    that of \(1-e^{(0,-\alpha)}\) (as \(\tilde{a}\) is invertible), we have our
    result.
\end{proof}

\begin{lemma}\label[lemma]{lem:reduce_to_sl2_prelim}
    Suppose \(a\in\bD_{\bM}^\Sm\cap\bC_{\bM}^\x\) and \(t\in\bT_{\bM}\) a lift of
    \(a\). Let \(\bW_{\bL_\alpha}\) be the Weyl group of the Levi subgroup
    \(\bL_\alpha\) generated by root \(\alpha\), where \(\alpha\) is the unique positive root
    such that \(\Disc_\alpha(t)=0\). Then the natural map
    \begin{align}
        \bC_\alpha\defeq\bar{\bT}_{\bM}\git \bW_{\bL_\alpha}\longto \bC_{\bM}
    \end{align}
    is \'etale over \(a\).
\end{lemma}
\begin{proof}
    It suffices to prove that it is \'etale at the image of \(t\) in
    \(\bC_\alpha\). Since the derived
    subgroup of \(\bG_+\) is simply-connected, the same is true for
    \(\bL_\alpha\), hence \(\bT_{\bM}\git\bW_{\bL_\alpha}\) 
    is smooth (treating \(\bG_+\) as a
    very flat reductive monoid containing \(\bL_\alpha^\Der\)). Then similar to
    the differential calculation done at the end of
    \Cref{chap:multiplicative_valuation_strata}, we see that the map
    \(\bC_\alpha\to\bC_{\bM}\) is smooth at the image of \(t\), because 
    the determinant of the 
    differential is given by a non-zero scaling of the product of
    \(\Disc_\beta\) for \(\beta\neq\alpha\).
\end{proof}

\begin{corollary}
    \label[corollary]{cor:etale_at_D_smooth}
    Suppose \(t\in\bar{\bT}_{\bM}\) is such that there exists a small enough
    \'etale neighborhood \(t\in U\) with non-empty intersection with at
    most one \(\bD_\alpha\), then the map
    \begin{align}
        \bC_\alpha\longto \bC_{\bM}
    \end{align}
    is \'etale at the image of \(t\).
\end{corollary}
\begin{proof}
    By \Cref{lem:reduce_to_sl2_prelim}, the map is \'etale on \(U\) outside
    the intersection of numerical boundary divisor and discriminant divisor, which
    have codimension at least \(2\). Since the map given is finite, we are done 
    since the branching locus has to have codimension \(1\).
\end{proof}

\begin{proof}[Proof of
    \Cref{prop:A_diamondsuit_then_cameral_is_smooth}]
    The statement is local, so we may look at the formal disc
    \(\breve{X}_{\bar{v}}\) at a point \(\bar{v}\in\breve{X}\) and the cameral cover
    over it. Let \(\cO=\bar{k}\powser{\pi}\) and \(F=\bar{k}\lauser{\pi}\).
    Since \(G\) is necessarily split over \(\breve{X}_{\bar{v}}\), we
    may assume without loss of generality that \(G=\bG\x\Spec{\cO}\) and
    \(\FRM=\bM\x\Spec{\cO}\) where \(\bM=\Env(\bG^\SC)\). Suppose the
    image of \(a\) in the abelianization \(\bA\) is contained in
    \(\pi^\lambda\bT^\AD(\cO)\).

    Suppose \(a\) intersects with \(\bD_{\bM}\) transversally. Let \(U\) be
    a small enough \'etale neighborhood of the special point of \(a\). 
    Then we can find a root
    \(\alpha\) determined by \Cref{lem:unique_root_for_d_equals_1}, and
    by \Cref{cor:etale_at_D_smooth}, the map \(\bC_\alpha\to \bC_{\bM}\)
    is \'etale over \(U\). Thus, we may lift \(a\) to a point
    \begin{align}
        a'\in\bC_\alpha(\cO),
    \end{align}
    and we only need to show that the preimage \(\tilde{X}_{a'}\) of \(a'\) in
    \(\bar{\bT}_{\bM}\) is smooth.

    Let \(\cO_2=\bar{k}\powser{\pi^{1/2}}\), \(F_2=\bar{k}\lauser{\pi^{1/2}}\),
    and \(\gamma\in\bar{\bT}_{\bM}(\cO_2)\cap\bT_{\bM}(F_2)\) be a point lifting \(a'\). 
    Let \(\nu_\gamma\) be the (dominant) Newton point of \(\gamma\). Let
    \(\gamma_0\bar{\bT}_{\bM}(\bar{k})\) be the special point of \(\gamma\),
    \(\gamma_1\in\bT_{\bM}(F_2)\) be the generic point.
    After \(\bW\)-conjugation, we may assume
    \(\gamma\in\pi^{(\lambda,\nu_\gamma)}\bT_{\bM}(\cO_2)\) (and \(\alpha\) may
    change to another root). The special point of \(\gamma\) is fixed by the
    reflection \(s_\alpha\) corresponding to \(\alpha\), so
    \(\alpha(\gamma)\in\cO_2^\x\) and \(\Pair{\alpha}{\nu_\gamma}=0\).
    Using the formula
    \begin{align}
        1=d_+(a)&=\Pair{2\rho}{\lambda-\nu_\gamma} 
        + \sum_{\Pair{\alpha}{\nu_\gamma}=0}\val_F(1-\alpha(\gamma))\\
                &=\Pair{2\rho}{-w_0(\lambda)-\nu_\gamma} 
        + \sum_{\Pair{\alpha}{\nu_\gamma}=0}\val_F(1-\alpha(\gamma)),
    \end{align}
    we see that there are two possibilities:
    \begin{enumerate}
        \item \(\nu_\gamma=-w_0(\lambda)\), and for any \(\beta\) with
            \(\Pair{\beta}{\nu_\gamma}=0\), \(1-\beta(\gamma)\) has non-zero
            valuation \(1/2\) if and only if \(\beta=\alpha\), or
        \item \(\nu_\gamma=-w_0(\lambda)-1/2\CoRt_i\) for some simple coroot
            \(\CoRt_i\), and \(\val_F(1-\beta(\gamma))=0\) for any root \(\beta\)
            with \(\Pair{\beta}{\nu_\gamma}=0\)
    \end{enumerate}

    The second possibility is easy to deal with: consider fundamental weight
    \(\Wt_i\) and let
    \begin{align}
        \mu=(-w_0(\Wt_i),-\Wt_i)\in\bar{k}[\bar{\bT}_{\bM}].
    \end{align}
    Then \(\Pair{\mu}{(\lambda,\nu_\gamma)}=1/2\) and so we have a morphism
    \begin{align}
        \mu\colon \bar{\bT}_{\bM}\longto \bbA^1
    \end{align}
    such that the composition with \(\gamma\) factors through a point
    \(\Spec{\cO_2}\to\bbA^1\) sending the coordinate of \(\bbA^1\) into
    \(\pi^{1/2}\cO_2^\x\). Since we have the factorization
    \begin{align}
        \Spec{\cO_2}\longto\tilde{X}_{a'}\longto \bbA^1,
    \end{align}
    we see that the ring of regular functions of \(\tilde{X}_{a'}\) is a subring
    of \(\cO_2\) that contains an element in \(\pi^{1/2}\cO_2^\x\), hence must
    be the whole ring \(\cO_2\), and thus \(\tilde{X}_{a'}\) is smooth.

    Now we deal with the first possibility. 
    If \(\mu\in\CharG(\bT_{\bM})\cap\bar{k}[\bar{\bT}_{\bM}]\) is any
    character such that \(\Pair{\mu}{\dual{\alpha}}=1\), then the map
    \begin{align}
        \bar{\bT}_{\bM}\longto \Spec{\bar{k}[e^\mu,e^{\mu-\alpha}]}\cong\bbA^2
    \end{align}
    is \(\bW_{\bL_\alpha}\) equivariant, hence induces commutative diagram
    \begin{equation}
        \begin{tikzcd}
            \bar{\bT}_{\bM} \ar[r]\ar[d] & \bbA^2 \ar[d]\\
            \bC_\alpha \ar[r] & \bC''\defeq\bbA^2\git \bW_{\bL_\alpha}\cong\bbA^2
        \end{tikzcd}
    \end{equation}
    Let \(a''\) be the image of \(a'\) in \(\bC''(\cO)\), and
    \(\tilde{X}_{a''}\) the corresponding \(\SL_2\)-cameral cover. Then we have
    commutative diagram
    \begin{equation}
        \label{eqn:temporary_cameral_double_prime}
        \begin{tikzcd}
            \tilde{X}_{a'} \ar[r]\ar[d] & \tilde{X}_{a''} \ar[d]\\
            a' \ar[r] & a''
        \end{tikzcd}
    \end{equation}

    If \(\lambda\neq 0\), since \(\dual{\alpha}\) is perpendicular to
    \(-w_0(\lambda)\neq 0\) which is dominant, it is not the highest coroot. By
    looking up the table of root systems, one can
    always find a fundamental weight \(\Wt\) such that
    \(\Pair{\Wt}{\dual{\alpha}}=1\). Then the character
    \(\mu=(-w_0(\Wt),-\Wt+\alpha)\in\bar{k}[\bar{\bT}_{\bM}]\) is such that
    \begin{align}
        \Pair{\mu}{(\lambda,\nu_\gamma)}=0\text{ and
        }\Pair{\mu}{\dual{\alpha}}=1.
    \end{align}
    By direct computation, we see that \(\tilde{X}_{a''}\) in
    \eqref{eqn:temporary_cameral_double_prime} is isomorphic to
    \(\Spec{\cO_2}\), hence smooth. Then 
    the map of \(1\)-dimensional schemes 
    \(\tilde{X}_{a'}\to\tilde{X}_{a''}\) is generically an isomorphism and finite,
    thus must be an isomorphism. This shows that \(\tilde{X}_{a'}\) is
    smooth. If \(\lambda=\nu_\gamma=0\), 
    then find an arbitrary weight \(\Wt\) with
    \(\Pair{\Wt}{\dual{\alpha}}=1\) and the dominant weight \(\Wt_+\) in
    \(\bW\)-orbit of \(\Wt\). Let \(\mu=(\Wt_+,\Wt)\), and we are done by
    repeating the calculation above.
\end{proof}

\begin{corollary}
    Let \(X_\OGT\to\breve{X}\) be a connected Galois \'etale covering making
    \(G\) split. Suppose \(a\in\cA_X^\diamondsuit\) is very \(G\)-ample,
    then \(\tilde{X}_{\OGT,a}\) is smooth and irreducible.
\end{corollary}
\begin{proof}
    Since \(\tilde{X}_a\) is smooth, so is \(\tilde{X}_{\OGT,a}\). It is
    also connected by \Cref{prop:cameral_curve_connected}, thus
    irreducible.
\end{proof}

Contrary to Lie algebra case (see \cite{Ng10}*{Lemme~4.7.3}), the converse to
\Cref{prop:A_diamondsuit_then_cameral_is_smooth}
is not true in general. For example, suppose \(\bG\) is of type
\(\TypeB_3\), and let
\begin{align}
    (\lambda,\nu_\gamma)&=
    (2\dual{\Wt}_2+2\dual{\Wt}_3+\dual{\Rt}_3,2\dual{\Wt}_2+2\dual{\Wt}_3),\\
    \mu&=(-w_0(\Wt_1),-\Wt_1+\Rt_1)=(\Wt_1,-\Wt_1+\Rt_1),\\
    \alpha &= \Rt_1,
\end{align}
where the labeling on the Dynkin diagram is the ``usual one'', i.e., the
vertex labeled with \(1\) is not directly joint with the one labeled with
\(3\), and \(\Rt_3\) is the short simple root.
Then \(\Pair{\mu}{(\lambda,\nu_\gamma)}=0\), \(\alpha\) is the only positive
root perpendicular to \(\nu_\gamma\), and \(\Pair{\mu}{\dual{\alpha}}=1\).
Thus, by carefully choosing an element \(t\in \bT_+(\cO_2)\) with
\(\val_F(1-\alpha(t))=1/2\), we can find \(a\in\bC_{\bM}(\cO)\) such that
\(\gamma=\pi^{(\lambda,\nu_\gamma)}t\) lies over \(a\), \(d_+(a)= 3\), and
\(\tilde{X}_a\) is smooth because in this case \(a'\) as in the
\Cref{prop:A_diamondsuit_then_cameral_is_smooth} exists and
\(\tilde{X}_{a''}\) hence \(\tilde{X}_{a'}\) is smooth by direct
computation. Since we do not use this in the remaining part of this book,
we leave the verification to the reader
(cf.~\Cref{chap:multiplicative_valuation_strata} on how to find such \(t\)).
Nevertheless, it is still easy to prove a partial converse as follows (we
will not use this result):

\begin{proposition}
    If \(\tilde{X}_a\) is smooth for some \(a\in\cA_X^\heartsuit(\bar{k})\),
    then for any points \(\bar{v}\in\breve{X}\), \(a(\bar{v})\) is contained in
    \(\bD_\alpha\) for at most one \(\alpha\).
\end{proposition}
\begin{proof}
    If \(\tilde{X}_a\) is smooth, then the local monodromy group
    \(\pi_a^\bullet(I_{v})\) is cyclic. But if \(a(\bar{v})\) is contained
    in two different \(\bD_\alpha\), then the local monodromy group would
    contain two different involutions. A contradiction.
\end{proof}

\section{Global N\'eron Model and \texorpdfstring{\(\delta\)}{δ}-Invariant}%
\label{sec:global_neron_model_and_delta_invariant}

Similar to the local case, we have seen in \eqref{eqn:global_neron_galois_desc}
that we have for any \(a\in\cA_X^\heartsuit(\bar{k})\) a N\'eron model
\(\FRJ_a^\flat\) of \(\FRJ_a\) described using normalization of cameral curves.
Another way of describing the global N\'eron model is to use
Beauville--Laszlo's descent theorem: let \(U=\breve{X}-\FRD_a\), then for any
closed point \(\bar{v}\in\FRD_a\), we have a local N\'eron model over the formal
disc \(\breve{X}_{\bar{v}}\) of the torus \(\FRJ_a|_{\breve{X}_{\bar{v}}^\bullet}\)
over the punctured disc. Gluing the local N\'eron models with the torus
\(\FRJ_a|_U\), we obtain a group scheme \(\FRJ_a^\flat\) which is precisely
the global N\'eron model.

Consider the Picard stack \(\cP_a^\flat\) of \(\FRJ_a^\flat\)-torsors. We have
a natural homomorphism of Picard stacks \(\cP_a\to\cP_a^\flat\).
    \nomenclature[\(P"cal_a_flat \)]{\(\cP_a^\flat\)}{the Picard stack of \(\FRJ_a^\flat\) torsors over \(X\)}
We will see
that \(\cP_a^\flat\) is an abelian stack in the following sense:
\begin{definition}
    [\cite{Ng10}*{D\'efinition~4.7.6}]
    A \(\bar{k}\)-abelian stack is the quotient of a \(\bar{k}\)-abelian variety
    by the trivial action of a diagonalizable group.
\end{definition}

\begin{proposition}
    \begin{enumerate}
        \item The homomorphism \(\cP_a(\bar{k})\to\cP_a^\flat(\bar{k})\) is
            essentially surjective.
        \item The neutral component \((\cP_a^\flat)_0\) of \(\cP_a^\flat\) is an
            abelian stack.
        \item The kernel \(\cR_a\)
            \nomenclature[\(R"cal_a \)]{\(\cR_a\)}{the kernel of \(\cP_a\to\cP_a^\flat\)}
            of \(\cP_a\to\cP_a^\flat\) is representable
            by the product of some affine algebraic groups of finite type
            \(\cR_{\bar{v}}(a)\) appearing in
            \Cref{lem:kernel_of_local_to_Neron} for finitely many
            \(\bar{v}\in\breve{X}\).
    \end{enumerate}
\end{proposition}
\begin{proof}
    For the first claim, consider short exact sequence
    \begin{align}
        1\longto \FRJ_a\longto\FRJ_a^\flat\longto \FRJ_a^\flat/\FRJ_a\longto 1,
    \end{align}
    and the induced cohomological exact sequence
    \begin{align}
        \RH^0(\breve{X},\FRJ_a^\flat/\FRJ_a)\longto \RH^1(\breve{X},\FRJ_a)\longto
        \RH^1(\breve{X},\FRJ_a^\flat)\longto \RH^1(\breve{X},\FRJ_a^\flat/\FRJ_a)=0,
    \end{align}
    where the last term vanishes because the sheaf \(\FRJ_a^\flat/\FRJ_a\) is
    supported on finite subset \(\breve{X}-U\). The first claim then follows.
    
    For the second claim, recall the Galois description
    \eqref{eqn:global_neron_galois_desc} of \(\FRJ_a^\flat\), and let
    \(\tilde{\cP}_a^\flat\) be the stack of torsors under group 
    \(\tilde{\FRJ}_a^\flat\defeq\pi_{\OGT,a,*}^\flat(\tilde{X}_{\OGT,a}^\flat\x\bT)\).
    We know that the neutral component of \(\tilde{\cP}_a^\flat\) is an abelian
    stack, because it is just the product of \(n_G\)-copies (\(n_G\) being the rank
    of \(G\)) of the usual Picard stack \(\PicS(\tilde{X}_{\OGT,a}^\flat)\)
    of line bundles on curve \(\tilde{X}_{\OGT,a}^\flat\). The norm
    homomorphism
    \begin{align}
        \tilde{\FRJ}_a^\flat&\longto
        \FRJ_a^\flat\\
        t&\longmapsto \prod_{w\in\bW\rtimes\Theta_\OGT}w(t)
    \end{align}
    induces homomorphism \(\tilde{\cP}_a^\flat\to\cP_a^\flat\) whose
    composition with the natural map \(\cP_a^\flat\to\tilde{\cP}_a^\flat\)
    is an isogeny of \(\cP_a^\flat\) to itself,
    as long as \(\Char(k)\) does not divide the order of
    \(\bW\rtimes\Theta_\OGT\). This proves the second claim.

    The third claim follows from a Beauville--Laszlo's descent theorem. Indeed, 
    by definition, the kernel \(\cR_a\) consists of pairs of a \(\FRJ_a\)-torsor
    together with a trivialization of its induced \(\FRJ_a^\flat\)-torsor. The
    local \(\cR_{\bar{v}}(a)\), on the other hand, consists of pairs of
    \(\FRJ_a|_{\breve{X}_{\bar{v}}}\)-torsors together with a trivialization of
    its induced \(\FRJ_a^\flat|_{\breve{X}_{\bar{v}}}\)-torsor, which also gives a
    trivialization of the said 
    \(\FRJ_a|_{\breve{X}_{\bar{v}}}\)-torsor over punctured disc
    \(\breve{X}_{\bar{v}}^\bullet\). Therefore, by the formal gluing theorem, the
    map of \(\bar{k}\)-functors
    \begin{align}
        \label{eqn:cR_a_local_to_global}
        \prod_{\bar{v}\in\breve{X}-U}\cR_{\bar{v}}(a)\longto \cR_a
    \end{align}
    obtained from gluing with the trivial torsor over \(U\) is an isomorphism.
    Hence, the claim.
\end{proof}

\begin{definition}
    \label[definition]{def:global_delta_def}
    Given \(a\in\cA_X^\heartsuit(\bar{k})\), the
    \notion{\(\delta\)-invariant}\index{global!\(\delta\)-invariant}\index{\(\delta\)-!invariant, global}
    associated with \(a\) is defined as
    \begin{align}
        \delta_a=\dim{\cR_a}.
        \nomenclature[\(delta_a \)]{\(\delta_a\)}{the global \(\delta\)-invariant of a point \(a\in\cA_X^\heartsuit(\bar{k})\)}
    \end{align}
\end{definition}
In view of isomorphism \eqref{eqn:cR_a_local_to_global}, we have that
\begin{align}
    \delta_a=\sum_{\bar{v}\in\breve{X}-U}\delta_{\bar{v}}(a).
\end{align}

\begin{corollary}
    \label[corollary]{cor:global_delta_formula_by_Neron}
    For \(a\in\cA_X^\heartsuit(\bar{k})\), we have formula
    \begin{align}
        \delta_a=\dim{\RH^0(\breve{X},\La{t}\otimes_{\cO_{\breve{X}}}(\pi_{a,*}^\flat\cO_{\tilde{X}_a^\flat}/\pi_{a,*}\cO_{\tilde{X}_a}))^W}
    \end{align}
\end{corollary}
\begin{proof}
    This is the result of the formula above and
    \Cref{lem:local_delta_formula_by_Neron}.
\end{proof}

\subsection{}
\label{sub:rigidification_of_global_Picard}
Fix a point \(\infty\in X(\bar{k})\), and consider open set
\(\cA_X^\infty\subset \cA_{X,\bar{k}}^\heartsuit\) 
    \nomenclature[\(A"cal_X_infty \)]{\(\cA_X^\infty\)}{the open subset of
    \(\cA_{X,\bar{k}}^\heartsuit\) where \(\FRD_X\) avoids a given geometric point \(\infty\in X(\bar{k})\)}
consisting of points \(a\) such that
\(a(\infty)\) is contained in \(\FRC_{\FRM}^\rss\). If \(\infty\) is defined over
a field extension \(k'/k\), then so is \(\cA_X^\infty\).

We may rigidify \(\cP_X\) over \(\cA_X^\infty\) as follows: let
\(\cP_X^\infty\)
    \nomenclature[\(P"cal_X_infty \)]{\(\cP_X^\infty\)}{a rigidified version of
    \(\cP_X\) over \(\cA_X^\infty\) by adding trivializations of \(\FRJ_a\)-torsors at \(\infty\)}
be the classifying stack over \(\cA_X^\infty\) of
\(\FRJ_a\)-torsors together with a trivialization at point \(\infty\). Then by
\cite{Ng10}*{Proposition~4.5.7}, \(\cP_X^\infty\to\cA_X^\infty\) is
representable by a smooth group scheme locally of finite type over
\(\cA_X^\infty\), and the forgetful
morphism \(\cP_X^\infty\to \cP_X\) induces a canonical isomorphism
\begin{align}
    \Stack{\cP_X^\infty/\FRJ_{X,\infty}}\longto \cP_X^\heartsuit,
\end{align}
where \(\FRJ_{X,\infty}\) is the group scheme over \(\cA_X^\infty\) being the
fiber of \(\FRJ_X\) at point \(\infty\).

For \(a\in\cA_X^\infty(\bar{k})\), let \(P_{a,0}\) be the neutral component of
\(\cP_a^\infty\). Chevalley structure theorem implies that
there is a canonical short exact sequence
\begin{align}
    1\longto R_a\longto P_{a,0}\longto A_a\longto 1,
\end{align}
where \(R_a\) is a smooth and connected affine algebraic group over \(\bar{k}\), 
and \(A_a\) is an abelian variety. Since \(R_a\) is smooth and connected, the
map \(R_a\to\cP_a^\flat\) is trivial. So the map \(P_{a,0}\to\cP_a^\flat\)
factors through \(A_a\) and is surjective. If 
\(\cP_a^\flat\) has finite automorphism groups
(cf.~\Cref{sec:automorphism_group}), then
\begin{align}
    \dim{A_a}=\dim{\cP_a^\flat}=\dim{\cP_a}-\delta_a.
\end{align}
In this situation, we call \(\delta_a\) the 
\notion{dimension of the affine part of \(\cP_a\)}\index{dimension!affine part of \(\cP_a\)}. Note that such notion is
only well-defined if the automorphism group is finite, otherwise since the
automorphism group here is also affine, it introduces ambiguity.

\section{Component Group}%
\label{sec:component_group}

In this section, we study the fiberwise component group \(\pi_0(\cP_X)\) defined 
as an \'etale sheaf on \(\cA_X\).

\subsection{}%
For a point \(a\in\cA_X^\infty(\bar{k})\), let \(U\subset \breve{X}\) be the
maximal open subset over which the cameral cover is \'etale. In particular,
\(\infty\in U\). Recall we have a pointed version of the group \(G\) after
fixing the point \(\infty\), given by a continuous homomorphism
\begin{align}
    \OGT_G^\bullet\colon \pi_1(X,\infty)\longto \Out(\bG).
\end{align}
If we also fix a point \(\tilde{\infty}\in \tilde{X}_a(\bar{k})\) lying
over \(\infty\), we can lift \(\OGT_G^\bullet\) into a commutative diagram
\begin{equation}
    \label{eqn:I_W_diagram}
    \begin{tikzcd}
        \pi_1(U,\infty) \ar[r, "\pi_{\tilde{a}}^\bullet"]\ar[d] 
                & \bW\rtimes\Out(\bG) \ar[d]\\
\pi_1(\breve{X},\infty) \ar[r, "\OGT_G^\bullet"] & \Out(\bG)
    \end{tikzcd},
\end{equation}
in which \(\tilde{a}=(a,\tilde{\infty})\). Let \(W_{\tilde{a}}\)
    \nomenclature[\(W_a'tilde \)]{\(W_{\tilde{a}}\)}{the monodromy subgroup in
    \(\bW\rtimes\Out(\bG)\) determined by a point \(\tilde{a}\in\tilde{\cA}_X(\bar{k})\)}
be the image of \(\pi_{\tilde{a}}^\bullet\) in \(\bW\rtimes\Out(\bG)\), and
\(I_{\tilde{a}}\)
    \nomenclature[\(I_a'tilde \)]{\(I_{\tilde{a}}\)}{the inertia subgroup of \(W_{\tilde{a}}\)}
the image of the kernel of \(\pi_1(U,\infty)\to
\pi_1(\breve{X},\infty)\) under \(\pi_{\tilde{a}}^\bullet\). By
commutativity of the diagram we have that \(I_{\tilde{a}}\subset \bW\).

\subsection{}
Recall that we may fix a connected Galois cover \(X_\OGT\to\breve{X}\) with Galois
group \(\Theta_\OGT\) over which the \(\Out(\bG)\)-torsor \(\OGT_G\) becomes
trivial. Choosing a point \(\infty_\OGT\) lying over \(\infty\), we have the
pointed variant \((X_\OGT,\infty_\OGT)\), and we may identify \(\Theta_\OGT\)
with the quotient \(\pi_1(\breve{X},\infty)/\pi_1(X_\OGT,\infty_\OGT)\). By
assumption of \(X_\OGT\), the map \(\OGT_G^\bullet\colon\pi_1(\breve{X},\infty)\to
\Out(\bG)\) factors through \(\Theta_\OGT\). The normalization
\(\tilde{X}_{\OGT,a}^\flat\) of cameral curve \(\tilde{X}_{\OGT,a}\)
maps to \(\breve{X}\), and it is an \'etale
\(\bW\rtimes\Theta_\OGT\)-cover over \(U\). By replacing
\((\tilde{X}_a,\tilde{\infty})\)
with \((\tilde{X}_{\OGT,a},\tilde{\infty}_\OGT)\), we may lift
\(\OGT_G^\bullet\) to 
\begin{align}
    \pi_{\tilde{a}_\OGT}^\bullet\colon\pi_1(U,\infty)\longto \bW\rtimes \Theta,
\end{align}
where \(\Theta\) is the image of \(\OGT_G^\bullet\) in \(\Out(\bG)\). 
Thus, we have inclusion
\begin{align}
    \label{eqn:W_a_tilde_in_W_rtimes_Theta}
    W_{\tilde{a}}\subset \bW\rtimes\Theta.
\end{align}

Let \(\FRJ_a^0\subset\FRJ_a\) be the open subgroup such that over any point
\(\bar{v}\in \breve{X}\) the fiber of \(\FRJ_a^0\) is the neutral component of the fiber of
\(\FRJ_a\). We have an induced homomorphism of Picard stacks \(\cP_a'\to\cP_a\).
    \nomenclature[\(P"cal_a_prime \)]{\(\cP_a'\)}{the analogue of \(\cP_a\) by
    replacing \(\FRJ_a\) with \(\FRJ_a^0\), or more generally any smooth group
    scheme \(\FRJ_a'\) that is a subsheaf of and generically isomorphic to
    \(\FRJ_a\)}

\begin{lemma}[\cite{Ng10}*{Lemme~4.10.2}]
    The homomorphism \(\cP_a'\to\cP_a\) is surjective with finite kernel. The
    same is true for induced homomorphism \(\pi_0(\cP_a')\to\pi_0(\cP_a)\).
\end{lemma}
\begin{proof}
    Since \(a\) is generically (over \(\breve{X}\)) regular semisimple, the sheaf
    \(\pi_0(\FRJ_a)=\FRJ_a/\FRJ_a^0\) has finite support on \(\breve{X}\). The short
    exact sequence
    \begin{align}
        1\longto\FRJ_a^0\longto \FRJ_a\longto\pi_0(\FRJ_a)\longto 1
    \end{align}
    induces cohomological long exact sequence
    \begin{align}
        \label{eqn:coh_long_exact_sequence_FRJ_a_and_FRJ_a_0}
        \RH^0(\breve{X},\pi_0(\FRJ_a))\longto\RH^1(\breve{X},\FRJ_a^0)\longto\RH^1(\breve{X},\FRJ_a)\longto
        \RH^1(\breve{X},\pi_0(\FRJ_a))=0.
    \end{align}
    Therefore, \(\cP_a'\to\cP_a\) is surjective with kernel being the image of
    \(\RH^0(\breve{X},\pi_0(\FRJ_a))\), which is necessarily finite. The induced
    map \(\pi_0(\cP_a')\to\pi_0(\cP_a)\) is then also surjective.
\end{proof}

\subsection{}%
Once we fix a point \(\infty\in X(\bar{k})\), there is a nice description of the
Cartier dual of the groups \(\pi_0(\cP_a')\) and
\(\pi_0(\cP_a)\). Let
\begin{align}
    \pi_0(\cP_a)^\odot&=\Spec(\Qlb[\pi_0(\cP_a)]),\\
    \pi_0(\cP_a')^\odot&=\Spec(\Qlb[\pi_0(\cP_a')]),
\end{align}
then the surjectivity of \(\pi_0(\cP_a')\to\pi_0(\cP_a)\) induces closed
embedding \(\pi_0(\cP_a)^\odot\subset\pi_0(\cP_a')^\odot\).

\begin{proposition}
    \label[proposition]{prop:global_pi_0_description}
    For any \(\tilde{a}=(a,\tilde{\infty})\), we have canonical 
    isomorphisms of diagonalizable groups
    \begin{align}
        \pi_0(\cP_a')^\odot&\simeq \dual{\bT}^{W_{\tilde{a}}},\\
        \pi_0(\cP_a)^\odot&\simeq \dual{\bT}(I_{\tilde{a}},W_{\tilde{a}}),
    \end{align}
    where \(\dual{\bT}(I_{\tilde{a}},W_{\tilde{a}})\subset
    \dual{\bT}^{W_{\tilde{a}}}\) is the subgroup of elements \(\kappa\)
    such that 
    \(I_{\tilde{a}}\subset\bW_{\bH}\) where \(\bW_{\bH}\) is the Weyl group of
    the neutral component \(\dual{\bH}\) 
    of the centralizer of \(\kappa\) in \(\dual{\bG}\).
\end{proposition}
\begin{proof}
    First, similar to the local case in
    \Cref{lem:pi_0_local_Picard_with_coinv_cochar},
    we have a canonical isomorphism
    \begin{align}
        \CoCharG(\bT)_{W_{\tilde{a}}}\longto \pi_0(\cP_a'),
    \end{align}
    see \cite{Ng06}*{Lemme~6.6 and Corollaire~6.7} and \cite{Ko85}*{Lemma~2.2}.
    Therefore, the first isomorphism follows by taking Cartier dual.

    Let \(U\subset \breve{X}\) be the regular semisimple locus of \(a\). Using
    \eqref{eqn:coh_long_exact_sequence_FRJ_a_and_FRJ_a_0}, we obtain exact
    sequence
    \begin{align}
        \RH^0(\breve{X},\pi_0(\FRJ_a))\longto\pi_0(\cP_a')\longto\pi_0(\cP_a)\longto
        0.
    \end{align}
    We may decompose \(\RH^0(\breve{X},\pi_0(\FRJ_a))\) as
    \begin{align}
        \RH^0(\breve{X},\pi_0(\FRJ_a))=\bigoplus_{\bar{v}\in\tilde{X}-U}\pi_0(\FRJ_{a,\bar{v}}),
    \end{align}
    where \(\FRJ_{a,\bar{v}}\) is the fiber of \(\FRJ_a\) at closed point
    \(\bar{v}\).
    Recall we also have local Picard groups \(\cP_{\bar{v}}(a_{\bar{v}})\)
    (resp.~\(\cP_{\bar{v}}^0(a_{\bar{v}})\)) associated with
    \(\FRJ_{a_{\bar{v}}}\) (resp.~\(\FRJ_{a_{\bar{v}}}^0\)) 
    induced by \(a\). Note here \(\FRJ_{a_{\bar{v}}}\)  is used to denote the
    restriction of \(\FRJ_a\) to the formal disc \(\breve{X}_{\bar{v}}\) (contrary
    to \(\FRJ_{a,\bar{v}}\) which is the special fiber). We then have exact
    sequence
    \begin{align}
        \pi_0(\FRJ_{a,\bar{v}})\longto \pi_0(\cP_{\bar{v}}^0(a_{\bar{v}}))\longto
        \pi_0(\cP_{\bar{v}}(a_{\bar{v}}))\longto 0,
    \end{align}
    compatible with the forgetful maps \(\cP_{\bar{v}}(a_{\bar{v}})\to\cP_a\) and
    \(\cP_{\bar{v}}^0(a_{\bar{v}})\to\cP_a'\).
    Taking Cartier dual, we have commutative diagram
    \begin{equation}
        \begin{tikzcd}
            0 \ar[r] & \pi_0(\cP_a)^\odot \ar[r]\ar[d] &
            \pi_0(\cP_a^0)^\odot\ar[r]\ar[d] & 
            \bigoplus_{\bar{v}\in\breve{X}-U}\pi_0(\FRJ_{a,\bar{v}})^\odot\ar[d,equal]\\
            0 \ar[r] &
            \bigoplus_{\bar{v}\in\breve{X}-U}\pi_0(\cP_{\bar{v}}(a_{\bar{v}}))^\odot\ar[r]
                     &
            \bigoplus_{\bar{v}\in\breve{X}-U}\pi_0(\cP_{\bar{v}}^0(a_{\bar{v}}))^\odot\ar[r] &
            \bigoplus_{\bar{v}\in\breve{X}-U}\pi_0(\FRJ_{a,\bar{v}})^\odot
        \end{tikzcd}.
    \end{equation}
    Thus \(\pi_0(\cP_a)^\odot\) is such subgroup of \(\pi_0(\cP_a')^\odot\) that its
    local restriction is contained in \(\pi_0(\cP_{\bar{v}}(a_{\bar{v}}))\),
    which is described in \Cref{prop:local_pi_0_description}.
    Since the group \(I_{\tilde{a}}\) is generated by the inertia groups at
    each \(\bar{v}\) (by a version of Riemann existence theorem; see e.g.,
    \cite{SGA1}*{Expos\'e~XIII, Corollaire~2.12}), such requirement is the same
    as that \(I_{\tilde{a}}\) is contained in \(\bW_{\bH}\) as desired since
    each inertia group is so.
\end{proof}

\begin{corollary}
    \(\pi_0(\cP_a)\) is finite if and only if
    \(\dual{\bT}^{W_{\tilde{a}}}\) (and equivalently,
    \(\bT^{W_{\tilde{a}}}\)) is.
\end{corollary}

\begin{definition}
    \label[definition]{def:anisotropic_locus}
    The \notion{anisotropic locus}\index{locus!anisotropic} \(\cA_X^{\ANI}\subset \cA_X^\heartsuit\)
    \nomenclature[\(.ani \)]{\((\cdot)^{\ANI}\)}{related to the anisotropic locus in \(\cA_X^\heartsuit\)}
    is the subset consisting of \(a\) such that \(\pi_0(\cP_a)\) is finite.
    Similarly, the \notion{elliptic locus}\index{locus!elliptic}
    \(\cA_X^{\ELL}\subset \cA_X^\heartsuit\)
    \nomenclature[\(.ell \)]{\((\cdot)^{\ELL}\)}{related to the elliptic locus in \(\cA_X^\heartsuit\)}
    is where \(\dual{\bT}^{W_{\tilde{a}}}\) (equivalently, \(\bT^{W_{\tilde{a}}}\)) is
    finite modulo the center.
\end{definition}

By definition, the anisotropic locus is contained in the elliptic one, and
equal to the latter if and only if \(Z_G\) does not contain a split
subtorus, or empty otherwise. Conceptually, the elliptic locus is a more natural
object to consider; however, for technical reasons, we will only use the
anisotropic locus for the most part.

It is not immediately obvious that \(\cA_X^{\ANI}\) or \(\cA_X^{\ELL}\) is an open subset of
\(\cA_X^\heartsuit\), or it is non-empty,
but we shall see in \Cref{sub:anisotropic_locus_is_open} and
\Cref{cor:non_emptyness_of_anisotropic_locus}
that both are true under mild assumptions.

\section{Automorphism Group}%
\label{sec:automorphism_group}

Let \((\cL,E,\phi)\in\cM_X^\heartsuit(\bar{k})\) with image
\(a\in\cA_X^\heartsuit\), \(b\in\cB_X\), and \(\cL\in\Bun_{Z_{\FRM}}\). 
Obviously, we have maps of automorphism groups 
\begin{align}
    \IAut(\cL,E,\phi)\longto \IAut(a)\longto \IAut(b)\longto\IAut_{Z_{\FRM}}(\cL).
\end{align}
We already know that \(\IAut(b)\) is finite and described using the kernel of
map \(Z_{\FRM}\to\FRA_{\FRM}^\x\). So the image
of \(\IAut(\cL,E,\phi)\) in \(\IAut_{Z_{\FRM}}(\cL)\) is finite, and we only need to
describe
\begin{align}
    \IAut(E,\phi)=\ker\left[\IAut(\cL,E,\phi)\to\IAut_{Z_{\FRM}}(\cL)\right],
\end{align}
or in other words, when the automorphism on \(\cL\) is identity.

Recall the universal centralizer \(I_{\FRM}\to\FRM\), whose fiber over
\(x\in\FRM\) is the centralizer \(G_x\subset G\) of \(x\) in \(G\). Since
\(I_{\FRM}\to\FRM\) is \(G\x Z_{\FRM}\)-equivariant, it descends to
\(\Stack{\FRM/G\x Z_{\FRM}}\). The pair
\((E,\phi)\) is a map \(\breve{X}\to \Stack{\FRM_\cL/G}\), and let \(I_{(E,\phi)}\)
be the pullback of \(I_{\FRM}\) along this map. It is easy to see that the sheaf
of automorphisms of \((E,\phi)\) is representable by \(I_{(E,\phi)}\). The group
\(I_{(E,\phi)}\) is not flat in general, but according to
\cite{BSL90}*{\S~7.1, Theorem~5}, there
exists a unique group scheme \(I_{(E,\phi)}^{\Sm}\)
    \nomenclature[\(I_E_phi_sm \)]{\(I_{(E,\phi)}^{\Sm}\)}{the smoothening of the
    centralizer group scheme \(I_{(E,\phi)}\) of a mHiggs bundle \((E,\phi)\) (with the boundary divisor being implicit)}
smooth over \(\breve{X}\) such
that for any \(\breve{X}\)-scheme \(S\) smooth over \(\breve{X}\), we have
\begin{align}
    \Hom_{\breve{X}}(S,I_{(E,\phi)})=\Hom_{\breve{X}}(S,I_{(E,\phi)}^{\Sm}).
\end{align}
The tautological map \(I_{(E,\phi)}^\Sm\to I_{(E,\phi)}\) is an isomorphism over
the open subset \(U=\breve{X}-\FRD_a\). Since \(\FRJ_a\) is smooth, the canonical map
\(\FRJ_a\to I_{(E,\phi)}\) induces the canonical map \(\FRJ_a\to I_{(E,\phi)}^\Sm\).
We also have the map
\begin{align}
    I_{(E,\phi)}^\Sm\longto \FRJ_a^\flat
\end{align}
by the universal property of N\'eron models. Both \(\FRJ_a\to I_{(E,\phi)}^\Sm\) and
\(I_{(E,\phi)}^\Sm\longto \FRJ_a^\flat\) are
isomorphisms over \(U\). At any point \(\bar{v}\in\breve{X}-U\), we have
inclusions
\begin{align}
    \FRJ_a(\breve{\cO}_{\bar{v}})\subset I_{(E,\phi)}^\Sm(\breve{\cO}_{\bar{v}})\subset
    \FRJ_a^\flat(\breve{\cO}_{\bar{v}}).
\end{align}
Using Beauville--Laszlo's descent theorem, we see that for any
\'etale cover \(S\to\breve{X}\), the maps of sets (in the category of
\(\breve{X}\)-schemes)
\begin{align}
    \FRJ_a(S)\longto I_{(E,\phi)}^\Sm(S)\longto \FRJ_a^\flat(S)
\end{align}
are both injective. This implies that we have inclusions
\begin{align}
    \RH^0(\breve{X},\FRJ_a)\subset
    \Aut(E,\phi)=\RH^0(\breve{X},I_{(E,\phi)})=\RH^0(\breve{X},I_{(E,\phi)}^\Sm)\subset
    \RH^0(\breve{X},\FRJ_a^\flat).
\end{align}
Using the Galois description of N\'eron model
\eqref{eqn:global_neron_galois_desc}, we have that
\(\RH^0(\breve{X},\FRJ_a^\flat)\subset\bT^{W_{\tilde{a}}}\) after fixing
\(\tilde{a}=(a,\tilde{\infty})\) over \(a\): indeed, let \(\Theta\) be
the image of \(\OGT_G^\bullet\) in \(\Out(\bG)\), then we have
\(W_{\tilde{a}}\subset \bW\rtimes\Theta\) by
\eqref{eqn:W_a_tilde_in_W_rtimes_Theta}.
Thus, we proved the analogue to \cite{Ng10}*{Corollaire~4.11.3}:
\begin{proposition}
    Fix \(\tilde{a}=(a,\tilde{\infty})\) over
    \(a\in\cA_X^\heartsuit(\bar{k})\), and any
    \((\cL,E,\phi)\in\cM_X^\heartsuit(\bar{k})\) lying over \(a\), the
    automorphism group \(\Aut(E,\phi)\) (after fixing \(\cL\)) can be canonically
    identified with a subgroup of \(\bT^{W_{\tilde{a}}}\).
\end{proposition}

If \(\Char(k)\) does not divide the order of \(\bW\rtimes\Theta\),
then \(\bT^{W_{\tilde{a}}}\) is unramified if it is finite.
Therefore, we have:
\begin{corollary}
    Assuming \(\Char(k)\) does not divide the order of
    \(\bW\rtimes\Theta\), then \(\cM_a^{\ANI}\) and \(\cP_a^{\ANI}\) are
    Deligne--Mumford stacks.
\end{corollary}
\begin{proof}
    Over \(a\in\cA_X^\ANI(\bar{k})\),
    \(\dual{\bT}^{W_{\tilde{a}}}\), hence also \(\bT^{W_{\tilde{a}}}\),
    are finite. This implies that \(\Aut(E,\phi)\) is finite for \((\cL,E,\phi)\)
    lying over \(a\). Since the image of \(\Aut(\cL,E,\phi)\) in
    \(\Aut_{Z_{\FRM}}(\cL)\) is finite, the group \(\Aut(\cL,E,\phi)\) is itself
    finite and the claim for \(\cM_a^{\ANI}\) follows. The claim for
    \(\cP_a^\ANI\) also follows since the automorphism groups (after fixing
    \(\cL\)) are just
    \(\RH^0(\breve{X},\FRJ_a)\) which is a subgroup of \(\bT^{W_{\tilde{a}}}\).
\end{proof}

In fact, for \(\cP_a\) we have a stronger result.
\begin{proposition}
    Assuming \(\Char(k)\) does not divide the order of
    \(\bW\rtimes\Theta\), then for any \(a\in \cA_{\gg}^\infty(\bar{k})\),
    we have
    \begin{align}
        \RH^0(\breve{X},\FRJ_a^1)&=\bT^{\bW\rtimes\Theta},\\
        \RH^0(\breve{X},\FRJ_a)&=\bZ_{\bG}^\Theta.
    \end{align}
    In particular, if \(Z_G\) does not contain a split torus over
    \(\breve{X}\) then \(\cP_a\) is a Deligne--Mumford stack.
\end{proposition}
\begin{proof}
    The curve \(\tilde{X}_{\OGT,a}\) is connected by
    \Cref{prop:cameral_curve_connected}. Then we have 
    \(\RH^0(\breve{X},\FRJ_a^1)=\bT^{\bW\rtimes\Theta}\), 
    and \(\RH^0(\breve{X},\FRJ_a)\) is a subgroup therein. If \(\Char(k)\) does
    not divide the order of \(\bW\rtimes\Theta\) and \(Z_G\) does not
    contain a split torus, they are both \'etale \(\bar{k}\)-groups, so
    \(\cP_a\) in this case is a Deligne--Mumford stack. 

    The description of
    \(\RH^0(\breve{X},\FRJ_a)\) is proved using the identification
    \(\FRJ_a=\FRJ_a'\) in \Cref{prop:Galois_reg_cent_monoid}. Since
    \(a(\breve{X})\) intersects with every irreducible component of discriminant
    divisor \(\FRD_{\FRM}\), the definition of \(\FRJ_a'\) implies that it is the
    subgroup of \(\bT^{\bW\rtimes\Theta}\) with elements lying in the kernel of
    every root, which is exactly \(\bZ_{\bG}^\Theta\).
\end{proof}

\section{Tate Module}%
\label{sec:tate_module}

Suppose the center of \(G\) does not contain a split torus over \(\breve{X}\).
Then over very \(G\)-ample locus \(\cB_\gg\),
\(\cP_\gg^\heartsuit\) is a Deligne--Mumford stack.
In this section, we fix a connected component \(\cU\) 
of \(\cA_\gg^\heartsuit\), and by abuse of notation we denote
\(\cP_{\gg}=\cP_X|_{\cU}\).

The Picard stack \(\cP_{\gg}\) is smooth over \(\cU\),
and let \(\cP_{\gg}^0\) be the open substack of fiberwise neutral component. Let
\(g\colon \cP_{\gg}^0\to\cU\) be the natural map. Let \(d\) be the
relative dimension of \(g\). Consider the sheaf of Tate modules
\begin{align}
    \TateM_{\Qlb}(\cP_{\gg}^0)=\RH^{2d-1}(g_!\Qlb).
    \nomenclature[\(Tate_P"cal \)]{\(\TateM_{\Qlb}(\cP)\)}{the sheaf of Tate
    modules of a Deligne--Mumford abelian stack \(\cP\) relative to a base stack}
\end{align}
Over the open subset \(\cU^\infty=\cA_X^\infty\cap\cU\), we may rigidify
\(\cP_{\gg}^0\) as the quotient of smooth group schemes
\(P^0\) by \(P^{-1}\) where \(P^{-1}\) is affine over \(\cU^\infty\). Thus, we
have a short exact sequence 
\begin{align}
    1\longto P^{-1}\longto P^0\longto
    \cP_{\gg}^0\longto 1,
\end{align}
and the induced short exact sequence of Tate modules over \(\Qlb\)
\begin{align}
    1\longto \TateM_{\Qlb}(P^{-1})\longto
    \TateM_{\Qlb}(P^0)\longto\TateM_{\Qlb}(\cP_{\gg}^0)\longto 1.
\end{align}
For any \(a\in\cU^\infty\), the Chevalley exact sequence
\begin{align}
    1\longto R_a\longto P_{a}^0\longto A_a\longto 1
\end{align}
also induces an exact sequence of Tate modules
\begin{align}
    1\longto \TateM_{\Qlb}(R_a)\longto
    \TateM_{\Qlb}(P_{a}^0)\longto\TateM_{\Qlb}(A_a)\longto 1.
\end{align}
Since \(P_{a}^{-1}\) is affine, the morphism
\(\TateM_{\Qlb}(P_{a}^0)\to\TateM_{\Qlb}(A_a)\) factors through
\(\TateM_{\Qlb}(\cP_{a}^0)\) (since \(\TateM_{\Qlb}(P_{a}^{-1})\) and
\(\TateM_{\Qlb}(A_a)\) have incompatible Frobenius weights over any sufficiently
large extension \(k'/k\) in \(\bar{k}\)). Therefore, we have an exact sequence
\begin{align}
    \label{eqn:canonical_decomp_Tate_module}
    1\longto \TateM_{\Qlb}(R_a)/\TateM_{\Qlb}(P_{a}^{-1})\longto
    \TateM_{\Qlb}(\cP_{a}^0)\longto \TateM_{\Qlb}(A_a)\longto 1.
\end{align}
It does not depend on the choice of rigidification, as if
\(\cP_{\gg}^0=\Stack{P^{0\prime}/P^{-1\prime}}\) is another rigidification, one may form
a third rigidification \(P^{0\prime\prime}=P^0\x_{\cP_{\gg}^0}P^{0\prime}\) over the two 
and identify exact sequences \eqref{eqn:canonical_decomp_Tate_module} in all
three cases. Since we only consider \(\Qlb\)-coefficients, an isogeny of
smooth commutative Deligne--Mumford group stacks induces isomorphism of
corresponding Tate modules. Thus, we have the canonical isomorphism
\begin{align}
    \TateM_{\Qlb}(R_a)/\TateM_{\Qlb}(P_{-1})\simeq
    \TateM_{\Qlb}(\cR_{a}),
\end{align}
and we may call
\(\TateM_{\Qlb}(\cR_{a})\) (resp.~\(\TateM_{\Qlb}(A_a)\)) the \notion{affine
part}\index{Tate module!affine part} (resp.~\notion{abelian part}\index{Tate module!abelian part}) of \(\TateM_{\Qlb}(\cP_{a}^0)\).

\begin{proposition}
    \label[proposition]{prop:tate_module_is_polarizable}
    Over \(\cU\), there exists an alternating bilinear form
    \begin{align}
        \psi\colon \TateM_{\Qlb}(\cP_{\gg}^0)\x\TateM_{\Qlb}(\cP_{\gg}^0)\longto
        \Qlb(-1)
    \end{align}
    such that at any \(a\in\cU\), the fiber \(\psi_a\) is identically
    zero on the affine part \(\TateM_{\Qlb}(\cR_a)\), and non-degenerate on the
    abelian part \(\TateM_{\Qlb}(A_a)\).
\end{proposition}
\begin{proof}
    The proof is identical to \cite{Ng10}*{\S~4.12} and we refer the reader to
    \textit{loc. cit}. The two key ingredients are
    one, the homomorphism in the Galois description of \(\FRJ_X\):
    \begin{align}
        \FRJ_X\longto \pi_{\OGT,*}(\tilde{X}_{\OGT}\x\bT),
    \end{align}
    and two, the Weil pairing theory, which is a general theory for any flat and
    proper family of reduced and connected curves. It is applied to the family
    of curves \(\tilde{X}_{\OGT}\), which we have shown to have reduced and
    connected fibers over \(\cU\) in
    \Cref{prop:cameral_curve_connected}.
\end{proof}

\section{Dimensions}%
\label{sec:dimensions}

In this section, we discuss the dimension estimates for both the mH-base and
mH-fibers. Most of this section resembles the Lie algebra case.

\subsection{}
Let \(\cL\in\Bun_{Z_{\FRM}}(\bar{k})\) be a very \(G\)-ample \(Z_{\FRM}\)-torsor.
As before, let \(\OGT\colon X_\OGT\to\breve{X}\) be a
connected finite Galois \'etale cover of \(\breve{X}\) with Galois group
\(\Theta\) and making \(G\) split. We also assume that \(\Char(k)\) does not
divide the order \(d_\OGT\) of \(\Theta\). Then for any \(b\in\cB_X(\bar{k})\)
lying over \(\cL\), \(\FRC_b\) is a direct summand
of \(\OGT^*\FRC_b\) being the \(\Theta\)-invariant subspace. We know
\(g_{X_\OGT}-1=d_\OGT(g_X-1)\). So if \(\cL\) is very \(G\)-ample, then for any
fundamental weight \(\Wt_i\) we have
\begin{align}
    \deg{\Wt_i(\OGT^*\cL)}>2g_{X_\OGT}-2.
\end{align}
By Riemann--Roch theorem, we have
\begin{align}
    \dim_{\bar{k}}{\RH^0(X_\OGT,\OGT^*\FRC_b)}
    &=\sum_{i=1}^r\deg{\Wt_i(\OGT^*\cL)}-r(g_{X_\OGT}-1),
\end{align}
and the first cohomology term vanishes. This means that \(\FRC_b\) has no first
cohomology term either. Taking \(\Theta\)-invariant, and note that \(\FRC_b\) is
the \(\Theta\)-trivial isotypic subbundle in \(\OGT_*\OGT^*\FRC_b\) which, as a
vector bundle with \(\Theta\)-action, is \(r\)-copies of regular ones (i.e., the
one whose fibers are regular \(\Theta\)-representations). Therefore, if the
boundary divisor of \(b\) can be written as 
\(\sum_{\bar{v}\in X(\bar{k})}\lambda_{\bar{v}}\cdot \bar{v}\), we have
\begin{align}
    \label{eqn:dimension_of_restricted_mH_base}
    \dim{\cA_b}=\sum_{\bar{v}\in
    X(\bar{k})}\Pair{\rho}{\lambda_{\bar{v}}}-r(g_X-1).
\end{align}
Note that due to \(\Out(\bG)\)-twisting the sum
\(\sum_{\bar{v}}\lambda_{\bar{v}}\) does not make sense, but its pairing with
\(\rho\) does since \(\rho\) is fixed by \(\Out(\bG)\). For convenience, we
denote this pairing by \(\Pair{\rho}{\lambda_b}\).

In case \(\cL\) is not very \(G\)-ample, we still have estimate by Riemann--Roch
theorem:
\begin{align}
    \label{eqn:dimension_estimate_of_restricted_mH_base_non_ample}
    \dim{\cA_b}\le\Pair{\rho}{\lambda_b}+r,
\end{align}
where the equality is reached only if \(\cL\) is trivial or \(g_X=0\).

\subsection{}
Suppose \(\FRM\in\FM_0(G^\SC)\) such that \(\OGT^*\FRM\cong\bM\x X_\OGT\) for
some \(\bM\in\FM_0(\bG^\SC)\). Then at each geometric point \(\bar{v}\in
X(\bar{k})\) we may identify \(\FRM\) with \(\bM\).
If \(\bA_{\bM}\cong\bbA^m\) is of standard type whose cone is freely generated
by cocharacters \(\theta_i\)
(\(1\le i\le m\)), and suppose at each \(\bar{v}\in X(\bar{k})\), we have
\begin{align}
    \lambda_{\bar{v}}=\sum_{i=1}^m c_{\bar{v},i}\theta_i,
\end{align}
then the local dimension of \(\cB_X\) at \(b\) is simply
\begin{align}
    \dim_b\cB_X=\sum_{\bar{v}\in X(\bar{k})}\sum_{i=1}^m c_{\bar{v},i},
\end{align}
which is locally constant. Combining with
\eqref{eqn:dimension_of_restricted_mH_base}, we have the dimension formula when
\(a\) is very \(G\)-ample:
\begin{align}
    \label{eqn:dimension_of_mH_base_at_a_point}
    \dim_a\cA_X= \dim_b\cB_X + \Pair{\rho}{\lambda_b} - r(g_X-1).
\end{align}
We also have dimension estimate in general regardless of ampleness
\begin{align}
    \label{eqn:dimension_of_mH_base_at_a_point_non_ample}
    \dim_a\cA_X\le  \dim_b\cB_X + \Pair{\rho}{\lambda_b} + r.
\end{align}

In the case where \(\bA_{\bM}\) is not an affine space, 
\eqref{eqn:dimension_of_mH_base_at_a_point} still holds since \(\cA_X\) is still
smooth over \(\cB_X\) at the very \(G\)-ample locus. However, the expression for
\(\dim_b\cB_X\) is less straightforward (although still possible to write down
using the normalization in \Cref{prop:boundary_divisor_space_strat_nonsplit})
because the numbers \(c_{\bar{v},i}\) are not uniquely determined.

\subsection{}
We would like to point out that in the above discussion, we effectively
achieved equality \eqref{eqn:dimension_of_mH_base_at_a_point} by proving that
the map \(\cA_X\to \cB_X\) is smooth at \(a\) if \(a\) is very \(G\)-ample
(through Riemann--Roch theorem). However, it is not strictly necessary: it is
possible that \(\cA_X\) is smooth at \(a\) over \(k\) but not over \(\cB_X\),
and if that is the case, \eqref{eqn:dimension_of_mH_base_at_a_point} should
still hold by computing the dimension of the tangent space of \(\cA_X\) at
\(a\). In other words, it seems possible (in some cases) to turn inequality
\eqref{eqn:dimension_of_mH_base_at_a_point_non_ample} into equality
\eqref{eqn:dimension_of_mH_base_at_a_point} even when \(G\)-ampleness is not
met. In order to do so, we will need to analyze the deformation of \(a\) in
\(\cA_X\) over \(k\), which is less straightforward compared to the deformation
relative to \(\cB_X\) as we have done. See \Cref{sec:the_case_of_mH_base} for
more details.

\subsection{}
Let \(\FRT_{\FRM}^\circ=\FRT_{\FRM}\cap\FRM^\circ\) and let
\(\FRC_{\FRM}^\circ\) be the image of \(\FRT_{\FRM}^\circ\).
It is an open subset of \(\FRC_{\FRM}\) since the cameral cover is flat, and its
complement has codimension at least \(2\) in \(\FRC_{\FRM}\), because it
contains both \(\FRC_\FRM^\rss\) and \(\FRC_\FRM^\x\). The torus \(T^\SC\) acts
freely (by translation) on the fibers of \(\FRT_{\FRM}^\circ\to\FRA_{\FRM}\) and has open orbits.
This implies that we have a description of relative tangent and cotangent bundles
\begin{align}
    \TanB_{\FRT_{\FRM}^\circ/\FRA_{\FRM}}&\cong \La{t}^\SC\x_X\FRT_{\FRM}^\circ,\\
    \CoTB_{\FRT_{\FRM}^\circ/\FRA_{\FRM}}&\cong (\La{t}^\SC)^*\x_X\FRT_{\FRM}^\circ.
    \nomenclature[\(T_Y_S \)]{\(\TanB_{Y/S}\)}{the relative tangent bundle of \(Y\) over \(S\)}
    \nomenclature[\(Omega_Y_S \)]{\(\CoTB_{Y/S}\)}{the relative cotangent bundle of \(Y\) over \(S\)}
\end{align}
This description is \(W\)-equivariant because \(T^\SC\) is commutative.
Since \(\FRC_{\FRM}^\circ\) is the \(W\)-invariant quotient of
\(\FRT_{\FRM}^\circ\), same is true for the total space of their tangent bundles,
in other words,
\begin{align}
    \TanB_{\FRT_{\FRM}^\circ/\FRA_{\FRM}}\git
    W\stackrel{\sim}{\longto}\TanB_{\FRC_{\FRM}^\circ/\FRA_{\FRM}}.
\end{align}
It implies that
\begin{align}
    (\pi_{\FRM*}\CoTB_{\FRT_{\FRM}^\circ/\FRA_{\FRM}})^W=\CoTB_{\FRC_{\FRM}^\circ/\FRA_{\FRM}}.
\end{align}
Since \(\CoTB_{\FRC_{\FRM}/\FRA_{\FRM}}=\FRC^*\x_X\FRC_{\FRM}\), we have that
\begin{align}
    \pi_{\FRM*}((\La{t}^\SC)^*\x_X\FRT_{\FRM}^\circ)^W=\FRC^*\x_X\FRC_{\FRM}^\circ.
\end{align}
The Killing form on \(\La{g}^\SC\) identifies \(\La{t}^\SC\) with its dual as
\(W\)-spaces. In addition, 
\(\FRC_{\FRM}-\FRC_{\FRM}^\circ\) has codimension \(2\), thus we have over
\(\FRC_{\FRM}\)
\begin{align}
    \pi_{\FRM*}(\La{t}^\SC\x_X\FRT_{\FRM})^W=\FRC^*\x_X\FRC_{\FRM}.
\end{align}
Let \(\cL\in\Bun_{Z_{\FRM}}\). Since \(W\) commutes with \(Z_{\FRM}\), the same
argument also applies to \(\cL\)-twisted cameral cover
\(\pi_{\FRM,\cL}\colon\FRT_{\FRM,\cL}\to\FRC_{\FRM,\cL}\). Since
\(\La{t}=\La{z}_G\x_X\La{t}^\SC\), we have
\begin{align}
    \pi_{\FRM,\cL,*}(\La{t}\x_X\FRT_{\FRM,\cL})^W=\La{z}_G\x_X\FRC_\cL^*\x_X\FRC_{\FRM,\cL}.
\end{align}
\begin{proposition}
    \label[proposition]{prop:desc_of_Lie_FRJ_a}
    For any \(a\in\cA_X(\bar{k})\) with associated \(Z_{\FRM}\)-torsor \(\cL\), we
    have canonical isomorphism
    \begin{align}
        \Lie(\FRJ_a)=\Lie(\FRJ_a^1)\simeq\La{z}_G\x\FRC_\cL^*.
    \end{align}
\end{proposition}
\begin{proof}
    The cameral cover is flat and finite. Using Cartesian diagram
    \begin{equation}
        \begin{tikzcd}
            \tilde{X}_a \ar[r]\ar[d, "\pi_a", swap] & \FRT_{\FRM,\cL} \ar[d,
            "\pi_{\FRM,\cL}"]\\
            X \ar[r, "a"] & \FRC_{\FRM,\cL}
        \end{tikzcd}
    \end{equation}
    and the Galois description of \(\FRJ_a^1\), we have the result by
    proper base change for coherent sheaves.
\end{proof}

\begin{remark}
    \label[remark]{rmk:duality_and_symplectic_on_mH_comments}
    There seems to be some kind of duality between the mH-base and mH-fibers,
    and \Cref{prop:desc_of_Lie_FRJ_a} can be seen as a primitive form of such
    duality. Such duality also hints to a symplectic structure on the
    mH-total-stack. It would be interesting to see developments in these
    directions.
\end{remark}

\begin{corollary}
    \label[corollary]{cor:dimension_of_P_a}
    For any \(a\in\cA_X^\heartsuit(\bar{k})\), we have
    \begin{align}
        \dim(\cP_a)&=\sum_{\bar{v}\in X(\bar{k})}\Pair{\rho}{\lambda_{\bar{v}}} +
        n_Gg_X-n_G,
    \end{align}
    where \(n_G\) is the rank of \(G\).
\end{corollary}
\begin{proof}
    Since
    \begin{align}
        \dim(\cP_a)=\dim_{\bar{k}}(\RH^1(\breve{X},\Lie(\FRJ_a)))
                -\dim_{\bar{k}}(\RH^0(\breve{X},\Lie(\FRJ_a))),
    \end{align}
    we have the desired equation using Riemann--Roch theorem.
\end{proof}

\section{Product Formula}%
\label{sec:product_formula}
Let \(a\in\cA_X^\heartsuit(k)\) and \(b\) (resp.~\(\cL\)) its image in
\(\cB_X\) (resp.~\(\Bun_{Z_{\FRM}}\)).  Let \(U_a=X-\FRD_a\). For each
closed point \(v\in \abs{X-U_a}\), we have the local multiplicative
affine Springer fiber \(\SP_{v}(a)\), defined by a choice of
\(\gamma_{v}\in\FRM_\cL^\x(F_{v})\) lying over \(a\) (which necessarily
exists by \Cref{sub:base_point_for_MASF_general}).

Let \(F=k(X)\) be the function field of \(X\), then \(a\) induces a point
\(a_F\in \Stack*{\FRC_\FRM/Z_\FRM}(F)\), which lifts to some point
\(\gamma_F\in\Stack*{\FRM/Z_\FRM}(F)\). Replacing \(U_a\) by some open dense
subset, we may also assume that \(\gamma_F\) extends to \(U_a\), denoted by
\(\gamma_{U_a}\). Accordingly, we shall choose \(\gamma_v\) at \(v\in
\abs{X-U_a}\) to be the restriction of \(\gamma_F\) to \(X_v^\bullet\).
The image of \(\gamma_{U_a}\) in 
\(\Stack*{\FRM_\cL/G}\) is then the trivial \(G\)-torsor \(E_{U_a}=E_0\) over
\(U_a\) together with a \(G\)-equivariant map \(\phi_{U_a}\colon E_{U_a}\to\FRM_\cL\).

Fix a \(k\)-algebra \(R\) and an \(R\)-point \(g_{v}\in\SP_{v}(a)(R)\).
This gives a \(G\)-torsor \(E_{v}\) 
over \(X_{v,R}\), a \(G\)-equivariant
map \(\phi_{v}\colon E_{v}\to \FRM_\cL\),
and a commutative diagram over \(X_{v,R}^\bullet\)
\begin{equation}
    \begin{tikzcd}
        E_{v}\ar[r, "\phi_{v}"]\ar[d, "\sim", "\iota_v"'] & \FRM_\cL\ar[d, equal]\\
        E_0\ar[r, "\gamma_{v}"] & \FRM_\cL
    \end{tikzcd}.
\end{equation}
Gluing \(E_{v}\) and \(E_{U_a}=E_0\) using \(\iota_{v}\),
we obtain a point \((E,\phi)\in\cM_a(R)\).  Therefore, we have
a morphism defined over \(k\)
\begin{align}
    \prod_{v\in X-U_a}\SP_{v}(a)&\longto \cM_a.
\end{align}
Similarly, by gluing with the trivial torsor over \(U_a\), we have
\begin{align}
    \label{eqn:P_a_local_to_global}
    \prod_{v\in \abs{X-U_a}}\cP_{v}(a)&\longto \cP_a.
\end{align}
The induced morphism 
\begin{align}
    \prod_{v\in \abs{X-U_a}}\SP_{v}(a)\x \cP_a&\longto \cM_a
\end{align}
is \(\prod_{v\in \abs{X-U_a}}\cP_{v}(a)\)-invariant (acting anti-diagonally), hence we have
a morphism
\begin{align}
    \label{eqn:product_formula_non_reduced}
    \prod_{v\in \abs{X-U_a}}\SP_{v}(a)\x^{\prod_{v\in \abs{X-U_a}}\cP_{v}(a)} \cP_a&\longto
    \cM_a
\end{align}
and its reduced version
\begin{align}
    \label{eqn:product_formula}
    \prod_{v\in \abs{X-U_a}}\SP_{v}^\Red(a)\x^{\prod_{v\in
    \abs{X-U_a}}\cP_{v}^\Red(a)} \cP_a&\longto
    \cM_a.
    \nomenclature[\(.red \)]{\((\cdot)^\Red\)}{the reduced structure of a space}
\end{align}
It is a equivalence of groupoids over geometric points and equivariant under
the Galois group of \(k\), and the same is true for any finite extension
\(k'/k\) if \(a\in\cA_X^\heartsuit(k')\) instead. Note that
\eqref{eqn:product_formula_non_reduced} is contingent on the fact that
\(\gamma_F|_{F_v}\) and \(\gamma_v\) are chosen to be \(G(F_v)\)-conjugate;
otherwise, it is still defined over \(\bar{k}\), because they are always
\(G(\breve{F}_{\bar{v}})\)-conjugate.

\begin{proposition}[Product Formula]
    \label[proposition]{prop:product_formula}
    For any finite extension \(k'/k\) and \(a\in\cA_X^\heartsuit(k')\),
    \eqref{eqn:product_formula} is a universal homeomorphism.
\end{proposition}
\begin{proof}
    We follow the strategy of \cite{BoCe22}. It suffices to assume \(k'=k\) by
    base change. Both sides of
    \eqref{eqn:product_formula} are algebraic
    stacks locally of finite type over \(k\).
    By \cite{BoCe22}*{Lemma~4.3.7}, it suffices to show that the map is an
    equivalence on \(R\)-points for all seminormal, strictly Henselian, local
    \(k\)-algebra \(R\).

    For any \(k\)-algebra \(R\), the map
    \eqref{eqn:product_formula_non_reduced} is fully faithful on
    \(R\)-points: indeed, let \(m=(\cL,E,\phi)\in\cM_a(R)\) where
    \(\cL\in\Bun_{Z_\FRM}(R)\), \(E\in\Bun_G(R)\) and \(\phi\colon
    E\to\FRM_\cL\) is \(G\)-equivariant. Suppose \(m\) is glued from
    \((m_v)\in\prod_{v\in \abs{X-U_a}}\cM_v(a)(R)\) and \(p\in\cP_a(R)\), and let
    \(\iota\in\Aut(m)\). The restriction of \(\iota\) to \(U_a\) is then given
    by a point \(j\in\FRJ_a(U_a)\). The restriction of \(\iota\) to the \(R\)-family of
    punctured discs \(\hat{X}_R^\bullet\) around \(X-U_a\) also gives a point
    \((p_v)\in\prod_{v\in \abs{X-U_a}}\cP_v(a)\). Twisting \(p\) by \((p_v)\),
    then \(m\) is still the image of \(((m_v),p)\), and then \(j\) extends to an
    automorphism of this new \(p\). This shows that there is a unique
    automorphism of \(((m_v),p)\) corresponding to \(\iota\) because both \(j\)
    and \((p_v)\) are uniquely determined by \(\iota\).

    By \cite{StacksP}*{\href{https://stacks.math.columbia.edu/tag/0EUQ}{Tag 0EUQ}}, any seminormal ring \(R\) is reduced,
    therefore \eqref{eqn:product_formula} is fully faithful on seminormal rings.
    To prove essential surjectivity, we note that \(m\) induces a
    tuple \((m_v)\) by restriction which depends on a choice of the points
    \((\gamma_v)\in\FRM_\cL^\x(\hat{X}_R^\bullet)\). By the construction process
    of \eqref{eqn:product_formula},
    we only need to show that any two such choices are conjugate under
    \(G(\hat{X}_R^\bullet)\). Indeed, the transporter between two choices of
    \((\gamma_v)\) is a \(\FRJ_a|_{\hat{X}_R^\bullet}\)-torsor, and so we need
    to show that
    \begin{align}
        \RH^1(\hat{X}_R^\bullet,\FRJ_a|_{\hat{X}_R^\bullet})=0.
    \end{align}
    This is implied by \cite{BoCe22}*{Theorem~3.2.4}, because
    \(\FRJ_a|_{\hat{X}_R^\bullet}\) becomes a split torus over the cameral
    cover.
\end{proof}

\begin{corollary}
    \label[corollary]{cor:mH_fibers_homeomorphic_to_proj_varieties}
    For any \(a\in\cA_X^\heartsuit(\bar{k})\), we have
    \begin{align}
        \dim{\cM_a}=\dim{\cP_a}.
    \end{align}
    If furthermore \(a\in\cA_X^\ANI(\bar{k})\), then \(\cM_a\) is homeomorphic to a
    projective \(\bar{k}\)-scheme.
\end{corollary}
\begin{proof}
    We shall see in \Cref{prop:properness_over_anisotropic_locus} that \(h_X\) is
    proper when restricted to the anisotropic locus \(\cA_X^\ANI\).
    Then the corollary follows from the product formula,
    \eqref{eqn:cR_a_local_to_global},
    \Cref{thm:GASF_dimension_formula,prop:projective_quotient_of_GASF}.
\end{proof}

\begin{corollary}
    \label[corollary]{cor:M_a_nonempty_k_fiber}
    For any \(a\in\cA_X^\heartsuit(\bar{k})\), the fiber \(\cM_a\) is non-empty.
    In fact, \(\cM_a^\reg\) is non-empty.
\end{corollary}
\begin{proof}
    A direct consequence of \eqref{eqn:product_formula}.
\end{proof}

\begin{remark}
    \Cref{prop:product_formula} is a stronger version of
    \cite{Ch22}*{Theorem~4.2.10}, since the existence of global Steinberg
    quasi-section is no longer required and the identification is made
    \(k\)-rational. On the other hand, because we do not know how to construct
    multiplicative affine Springer fibers in a family yet, this method so far
    only works when \(a\) is a \(k\)-point. Later we will upgrade it into a
    family in \Cref{prop:simultaneous_product_formula}.
\end{remark}

\subsection{}
There is another version of product formula which is also useful. Although it
can be proved using \eqref{eqn:product_formula_non_reduced}, there is an
independent formulation. Following
\cite{Ng06}, for any closed point \(v\in X\), we may consider the stack
\(\cM_{a,v}\) classifying pairs \((E_v,\phi_v)\) where \(E_v\) is a \(G\)-torsor
over the formal disc \(X_v\) and \(\phi_v\) is a \(G\)-equivariant map
\(E_v\to\FRM_\cL\). In other words, it is a sort of stacky version of
multiplicative
affine Springer fiber. We also let \(\cP_{a,v}\) be the classifying stack of
\(\FRJ_a\)-torsors over \(X_v\).

The multiplicative affine Springer fiber \(\SP_{v}(a)\)
naturally maps to \(\cM_{a,v}\) by forgetting the trivialization part. On the
other hand, \(\cM_a\) also maps to \(\cM_{a,v}\) by restricting to \(X_v\).
Similarly, we may define \(\cM_{a,\bar{v}}\) and \(\cP_{a,\bar{v}}\) for any
geometric point \(\bar{v}\in\breve{X}\).
Therefore, we have a commutative diagram
\begin{equation}
    \begin{tikzcd}
        \prod_{\bar{v}\in \abs{\breve{X}-U_a}}\SP_{\bar{v}}(a) \ar[rr]\ar[rd] && \cM_a
        \ar[ld]\\
                                                                      &
        \prod_{\bar{v}\in \abs{\breve{X}-U_a}}\cM_{a,\bar{v}}& 
    \end{tikzcd}
\end{equation}
where the horizontal arrow is defined over \(\bar{k}\) while the two diagonal
ones are defined over \(k\). So we have induced maps of \(2\)-categorical quotients
\begin{equation}
    \begin{tikzcd}
        \prod_{\bar{v}\in \abs{\breve{X}-U_a}}\Stack{\SP_{\bar{v}}(a)/\cP_{\bar{v}}(a)}
        \ar[rr]\ar[rd] && \Stack{\cM_a/\cP_a} \ar[ld]\\
                                                                      &
        \prod_{\bar{v}\in \abs{\breve{X}-U_a}}\Stack{\cM_{a,\bar{v}}/\cP_{a,\bar{v}}}& 
    \end{tikzcd}
\end{equation}
and all three maps are equivalences of groupoids after taking
\(\bar{k}\)-points. 
The map \eqref{eqn:P_a_local_to_global} is always
defined over \(k\), hence so are the diagonal maps above.
It implies that the horizontal map is defined over \(k\), and 
induces an equivalence of \(k'\)-points for any \(k'/k\)
regardless whether \eqref{eqn:product_formula} is defined over \(k\) or
not. We summarize it in the following proposition.
\begin{proposition}
    \label[proposition]{prop:product_formula_alt}
    For \(a\in\cA_X^\heartsuit(k)\), the \(2\)-stack \(\Stack{\cM_a/\cP_a}\) is
    a stack, and we have a natural map of stacks over \(k\)
    \begin{align}
        \prod_{v\in \abs{X-U_a}}\Stack{\SP_{v}(a)/\cP_{v}(a)}\longto
        \Stack{\cM_a/\cP_a}
    \end{align}
    that induces equivalence on \(k'\)-points for any field extension \(k'/k\).
\end{proposition}
\begin{proof}
    That \(\Stack{\cM_a/\cP_a}\) is a stack follows from the fact
    that the automorphism group of any object in \(\cP_a\) is canonically
    embedded in that of any object in \(\cM_a\)
    (cf.~\Cref{sec:automorphism_group}). The other claims follow from our
    discussions above.
\end{proof}

\section{Local Model of Singularities}%
\label{sec:local_model_of_singularities}

Unlike in the Lie algebra case, the total space \(\cM_X\)
of mH-fibration is no longer smooth. Instead, it admits a 
stratification induced by affine Schubert cells, which in turn translates to
representations via geometric Satake equivalence. Therefore, the existence of
singularities is in fact a feature of mH-fibrations.

The main result of this section was essentially predicted in
\cite{FN11}*{Conjecture~4.1}, and a weaker version was first proved in
\cite{Bo17}. Readers can find a more streamlined proof in \cite{Ch22}, due to
Zhiwei Yun. However, these previous results would turn out to be too weak for
studying endoscopy. The main reason behind is probably due to that the core
argument is \textit{ad~hoc} in nature, and in particular it does not try to
establish a tangent-obstruction theory for deformations of mHiggs bundles,
contrary to what is done in the Lie algebra case (see \cite{Ng10}*{\S~4.14}).
Our goal is then to establish such missing tangent-obstruction theory, after
which we will be able to prove a much stronger result.

\subsection{}
Let \(h_X\colon\cM_X\to\cA_X\) be the universal mH-fibration of a monoid
\(\FRM\in\FM(G^\SC)\).
Recall
that in \ref{sub:partial_affine_Schubert} we defined global affine Schubert
schemes on \(\BD_X\):
\begin{align}
    \OGASch_X\longto\BD_X,
\end{align}
and we can pull it back to \(\cB_X\) and denote it by \(\GASch_X\).
    \nomenclature[\(Q_X \)]{\(\GASch_X\)}{the pullback of \(\OGASch_X\) to \(\cB_X\)}
We also have
the Hecke stack \(\Stack*{\GASch_X}\),
    \nomenclature[\(Q_X_stack \)]{\(\Stack*{\GASch_X},\Stack*{\GASch_X}_N\)}{the pullback of
    \(\Stack*{\OGASch_X}\) (resp.~\(\Stack*{\OGASch_X}_N\)) to \(\cB_X\)}
as well as
the truncated version \(\Stack*{\GASch_X}_{N}\) (see
\Cref{def:equivariant_global_affine_Schubert}). Note that the number \(N\)
depends on the irreducible component of \(\cB_X\). We shall abuse the notation
and let \(N\) be a fixed, sufficiently large number since we can always replace
\(\cB_X\) by some of its connected components.
So we have an evaluation map
\begin{align}
    \ev_N\colon \cM_{X}&\longto \Stack*{\GASch_X}_{N},
    \nomenclature[\(ev_ev_N \)]{\(\ev,\ev_N\)}{the evalluation map \(\cM_X\to
    \Stack*{\GASch_X}\) (resp.~\(\cM_X\to\Stack*{\GASch_X}_N\))}
\end{align}
which factors through the limit
\begin{align}
    \ev\colon \cM_{X}&\longto \Stack*{\GASch_X}.
\end{align}
The \notion{local model of singularity}\index{local model of singularity} expects that the map \(\ev_{N}\)
(resp.~\(\ev\)) is smooth (resp.~formally smooth). Of course, it will be too
much to hope for unconditionally: it will be more reasonable at least
restrict to the reduced locus \(\cM_X^\heartsuit\) and perhaps also some \(G\)-ampleness
requirements. For example, we will be able to prove the
following result which is already a big improvement of \cite{Bo17} (also
compare with \cite{Ch22}):
\begin{theorem}[Local model of singularity, part 1]
    \label[theorem]{thm:local_singularity_model_weak}
    Let \(x=(\cL,E,\phi)\in \cM_X^\heartsuit(\bar{k})\) be a point and let
    \(a=h_X(x)\in\cA_X(\bar{k})\) be its image. If
    \begin{align}
        \RH^1(\breve{X},(\Lie(I_{(E,\phi)}^{\Sm})/\La{z}_G)^*)=0,
    \end{align}
    then \(\ev_{N}\) (resp.~\(\ev\)) is smooth
    (resp.~formally smooth) at \(x\). In particular, it is true when \(x\) is
    very \((G,\delta_a)\)-ample,
\end{theorem}
The second part of the theorem is contained in
\Cref{thm:local_singularity_model_main},
whose statement requires the factorization results in \Cref{sec:factorizations}.
Before introducing the other part, we need to understand why
\Cref{thm:local_singularity_model_weak} is insufficient and what the
idea is behind the enhancement.

\subsection{}
The main reason for a further improvement is that
\Cref{thm:local_singularity_model_weak} is still not quite strong enough for
fundamental lemma. In short, 
outside potential edge cases, the cohomological conditions in
\Cref{thm:local_singularity_model_weak} for \(G\) and for endoscopic group \(H\)
can be met in conjunction only when the endoscopic group \(H\) is ``large''
compared to \(G\) in the sense that its semisimple rank \(r_H\) is actually
equal to \(r\). This means that we do not even cover the case of \(H\) being an
elliptic torus. See \Cref{sub:starting_proof_of_local_singularity_model_main} for
more details.

\subsection{}
On the other hand, we observe that the conclusion of
\Cref{thm:local_singularity_model_weak} is unnecessarily strong. To illustrate,
suppose we have a potentially singular variety \(Y\) and a set of singular
varieties \(Y_1,\ldots,Y_m\) whose singularities are ``well understood'', and
suppose there is a map
\begin{align}
    f\colon Y\longto Y_1\x\cdots \x Y_m.
\end{align}
If we can show that \(f\) is smooth, then since the singularities of \(Y_i\) are
well understood, we also understand the singularity of \(Y\). This is what
is done in \Cref{thm:local_singularity_model_weak} because the target of the
map \(\ev\) is more or less a direct product of affine Schubert varieties (hence
the ``local model'' in its name). Now
suppose, for example, we know that the map \(f\) is such that the first
coordinate of its image is entirely contained in the smooth locus \(Y_1^\Sm\) of
\(Y_1\), in other words, \(f\) factors through the map
\begin{align}
    f_1\colon Y\longto Y_1^\Sm\x Y_2\x\cdots\x Y_m.
\end{align}
Then the factor \(Y_1\) has no contribution towards understanding the
singularity of \(Y\), because the target of \(f_1\) has the same
singularity as \(Y_2\x\cdots\x Y_m\). Therefore,
instead of showing that \(f_1\) is smooth, we only need to show the
smoothness of the composition
\begin{align}
    Y\longto Y_1^\Sm\x Y_2\x\cdots\x Y_m\longto Y_2\x\cdots\x Y_m,
\end{align}
which is a weaker statement.

\subsection{}
For mH-fibrations in particular, suppose we have an mHiggs bundle \((\cL,E,\phi)\)
with boundary divisor \(b\), and we write
\begin{align}
    \lambda_b=\sum_{i=1}^d\lambda_i\cdot\bar{v}_i.
\end{align}
Over \(b\), the fiber of \(\Stack*{\GASch_X}\) is isomorphic to
\begin{align}
    \prod_{i=1}^d\Stack*{\Arc_{\bar{v}_i} G^\AD\backslash \Gr_{G^\AD,\bar{v}_i}^{\le-w_0(\lambda_{i,\AD})}}.
\end{align}
If, for example, we further assume that at \(\bar{v}_1\) the mHiggs bundle is
contained in the big-cell locus, then the first factor of the right-hand side
has no contribution to the singularity and thus can be eliminated.
In fact, we can do even better: if we only look at the map
\begin{align}
    \cM_b\longto \prod_{i=2}^d\Stack*{\Arc_{\bar{v}_i} G^\AD\backslash
    \Gr_{G^\AD,\bar{v}_i}^{\le-w_0(\lambda_{i,\AD})}},
\end{align}
then there is no reason to keep \(\bar{v}_1\) (hence also \(b\) itself) fixed in
place, and adding such flexibility means more room for deformation, which in
turn means that we are more likely to obtain smoothness results.

This is precisely our plan to enhance \Cref{thm:local_singularity_model_weak}.
Doing so requires separating the boundary divisors into two disjoint
subdivisors, and so the factorization results in \Cref{sec:factorizations}
become useful. We will only provide the relevant statements in this section and
the proof will be postponed (see
\Cref{sec:deformation_of_mHiggs_bundles}) because the deformation argument is
quite involved.

\subsection{}
We use the notations from \Cref{sec:factorizations} and let \(I=\Set{1,2}\), so
that \(\BD[I]_X=\BD_X\x \BD_X\). We have three maps from \(\BD[I]_X\) to
\(\BD_X\) being the summation map \(\Sigma\) and two projections \(\pr_1\) and
\(\pr_2\). They induce three different affine Schubert schemes
\begin{align}
    \OGASch_I&=\Sigma^*\OGASch_X,\\
    \OGASch_i&=\pr_i^*\OGASch_X\quad (i=1,2).
\end{align}
By \Cref{prop:factorization_of_affine_Schubert}, we have canonical isomorphism
over the disjoint locus \(\BD[I,\disj]_X\)
\begin{align}
    f_I\colon \bigl.\OGASch_1\x_{\BD[I]_X}\OGASch_2
    \bigr|_{\BD[I,\disj]}\stackrel{\sim}{\longto}\bigl.\OGASch_I\bigr|_{\BD[I,\disj]}.
\end{align}
Similar results hold for the \(G\)-equivariant versions as well.

Let \(h_X^I\colon \cM_X^I\to \cA_X^I\) be the mH-fibration associated with
monoid \(\FRM^I\) (again, \(\FRM^I\) is \emph{not} the \(I\)-fold direct
product of \(\FRM\)), and we use \(h_X^{I,\disj}\), \(\cM_X^{I,\disj}\), etc. to
denote the restriction to \(\BD[I,\disj]_X\). We have the following commutative
diagram
\begin{equation}
    \begin{tikzcd}
        \cM_X^{I,\disj} \ar[r]\ar[rd] 
        & \Stack*{\GASch_I^{\disj}} \ar[r]\ar[d] 
        & \Stack*{\OGASch[\disj]_I} \ar[r, "\pr_1\circ f_I^{-1}"]\ar[d] 
        & \Stack*{\OGASch_1} \ar[r] \ar[ld]
        & \Stack*{\OGASch_X} \ar[lldd]\\
        & \cB_X^{I,\disj} \ar[r]\ar[rd]
        & \BD[I,\disj]_X \ar[d, "\pr_1"]
        &
        &\\
        &
        & \BD_X
    \end{tikzcd}
\end{equation}
The composition of the top row is denoted by
\begin{align}
    \ev_1^I\colon \cM_X^{I,\disj}\longto \Stack*{\OGASch_X},
\end{align}
which is a map over \(\BD_X\). Similarly, we have \(\ev_{1,N}^I\) for
sufficiently large integers \(N\).

\subsection{}
In general,
\(\ev_1^I\) may not be formally smooth for a rather stupid reason: the map from
\(\cM_X^{I,\disj}\) to \(\BD_X\) factors through projection \(\pr_1\), whose
fiber is \(\BD_X\) itself, but \(\BD_X\) may also have singularity (controlled
by the singularities of \(\FRA_\FRM\)) which has
no connection with the affine Schubert schemes. Therefore to avoid making
the statement unnecessarily complicated, we will now restrict to monoids such that
\(\FRA_\FRM\) is of standard type. Since these
will likely be the only monoids that appear in applications anyway, it is not a
serious drawback.
We will also need a technical notion that is related to characteristics:
\begin{definition}
    \label[definition]{def:goodness_of_boundary_divisors}
    Suppose \(\theta_1,\ldots,\theta_m\) are the free generators of the cone of
    \(\bA_\bM\). The cocharacter \(\theta_i\) is called
    \notion{good}\index{cocharacter!good} if
    \(\theta_{i,\AD}\not\in\Char(k)\CoCharG(\bT^\AD)\). A point
    \(b\in\BD_X(\bar{k})\) is \notion{good}\index{divisor!boundary, good} if at every \(\bar{v}\in
    X(\bar{k})\), the cocharacter of \(b\) at \(\bar{v}\) is a sum of good
    \(\theta_i\)'s (after choosing a local splitting \(\FRA_\FRM\cong\bA_\bM\)
    at \(\bar{v}\)). A connected component of \(\BD_X\) is \notion{good} if
    all its points are good.
\end{definition}

\begin{remark}
    \begin{enumerate}
        \item The notion of \(b\) being good clearly does not depend on the local
            identification \(\FRA_\FRM\cong\bA_\bM\). Moreover, if \(b\in\BD_X\) is
            good, then the connected component containing \(b\) is also good.
        \item When \(\FRM=\Env(G^\SC)\), every component of \(\BD_X\) is good
            because every fundamental coweight is good.
        \item For a fixed monoid \(\bM\), viewed as combinatorially constructed
            from root data and cones,
            there exists some \(N\) such that as long as \(\Char(k)>N\), every
            component of \(\BD_X\) is good for any twist \(\FRM\) of \(\bM\).
        \item Goodness is always satisfied if \(\Char(k)=0\) and \(\theta_i\) is
            not central.
    \end{enumerate}
\end{remark}

\begin{theorem}
    [Local model of singularity, part 2]
    \label[theorem]{thm:local_singularity_model_main}
    Let \(\FRM\in\FM(G^\SC)\) be such that \(\FRA_\FRM\) is of standard type
    with rank \(m\). Let
    \(x=(\cL^I,E,\phi)\in(\cM_X^{I,\disj}\cap\cM_X^{I,\heartsuit})(\bar{k})\) be
    a point where \(\cL^I\) is a \(Z_{\FRM^I}\)-torsor. Let \(a=h_X^I(x)\) and
    \(b_I\in\cB_X\) be the boundary divisor of \(x\). Let \(b_i\in\BD_X\) be the
    image of \(b_I\) in \(\BD_X\) through the \(i\)-th projection (\(i=1,2\)).
    Then \(\ev_{1,N}^I\) (resp.~\(\ev_{1}^I\)) is smooth
    (resp.~formally smooth) at \(x\) if the followings are satisfied:
    \begin{enumerate}
        \item \(b_2\) is good,
        \item \(\supp(\FRD_a)\cap\supp(b_2)=\emptyset\),
        \item \(\dim_{b_2}\BD_X > 2\abs{\bW}m(g_X-1)\).
    \end{enumerate}
\end{theorem}
\begin{proof}
    See \Cref{sec:deformation_of_mHiggs_bundles}.
\end{proof}

\begin{remark}
    \begin{enumerate}
        \item Although \Cref{thm:local_singularity_model_weak} and
            \Cref{thm:local_singularity_model_main} do not supersede each other,
            the latter is in most cases a lot stronger than the former,
            because \(2\abs{\bW}m(g_X-1)\) is a constant depending only on the
            curve \(X\) and monoid \(\FRM\), while the number \(\delta_a\) can
            become quite large.
        \item Even though \Cref{thm:local_singularity_model_weak} is weaker
            most of the time, it is sometimes more convenient to use due
            to being stated over the entire moduli of boundary divisors
            \(\cB_X\), and does not require \(\FRA_\FRM\) to be of standard
            type.
    \end{enumerate}
\end{remark}

\begin{corollary}
    \label[corollary]{cor:local_model_sing_sheaf_version}
    Let \(\cF\in\Sat_X\) be a Satake sheaf, then it descends to a perverse sheaf
    on \(\Stack*{\GASch_X}_{N}\). Moreover, over the locus where either
    \Cref{thm:local_singularity_model_weak} or
    \Cref{thm:local_singularity_model_main} holds, \(\ev_N^*\cF\) is
    perverse fiberwise over \(\cB_X\).
\end{corollary}
\begin{proof}
    Immediate from the two said theorems.
\end{proof}

\begin{definition}
    A \notion{Satake sheaf}\index{Satake!sheaf}\index{sheaf!Satake} on \(\cM_X\) is a complex such that over the locus
    where either \Cref{thm:local_singularity_model_weak} or
    \Cref{thm:local_singularity_model_main} holds it is of the form
    \(\ev_N^*\cF\) for some \(\cF\in\Sat_X\), but shifted and Tate-twisted so that it
    is pure of weight \(0\). We still denote the category of these complexes by
    \(\Sat_X\).
\end{definition}

\section{The Case of Endoscopic Groups}%
\label{sec:the_case_of_endoscopic_groups}

In this section, we generalize the framework of mH-fibrations to endoscopic
groups, for which one major difference is that the twist provided by the center
of the monoid needs to be replaced by a subgroup therein. We will also
the reason behind not fixing a \(Z_\FRM\)-torsor in previous sections.

\subsection{}
Let \((\kappa,\OGT_\kappa)\) be an endoscopic datum of \(G\) on \(X\), and
\(H\) is the endoscopic group. By \Cref{sec:endoscopic_groups_inv_theory}, there
is a canonical monoid \(\FRM_H\in\FM(H^\SC)\) associated with \(\FRM\) and a 
canonical map \(\nu_H\colon \FRC_{\FRM,H}\to\FRC_{\FRM}\). Let
\(Z_{\FRM,H}\subset\FRM_H^\x\) be the center, and we have mH-fibration
\begin{align}
    h_{H,X}\colon\cM_{H,X}\longto \cA_{H,X}
    \nomenclature[\(.H_X \)]{\((\cdot)_{H,X}\)}{global constructions related to \(\FRM_H\) and the curve \(X\)}
\end{align}
associated with \(\FRM_H\). However, there is no direct relation between \(h_X\)
and \(h_{H,X}\), because \(Z_{\FRM,H}\) does not map into \(Z_{\FRM}\). Instead,
we need to replace \(Z_{\FRM,H}\) with its subgroup \(Z_{\FRM}^\kappa\). Let
\begin{align}
    h_{H,X}^\kappa\colon\cM_{H,X}^\kappa\longto \cA_{H,X}^\kappa
    \nomenclature[\(.H_X_kappa \)]{\((\cdot)_{H,X}^\kappa\)}{global constructions
    related to \((\cdot)_{H,X}\) and endoscopic datum \(\kappa\)}
\end{align}
be the pullback of \(h_{H,X}\) through \(\Bun_{Z_{\FRM}^\kappa}\to
\Bun_{Z_{\FRM,H}}\). Let \(\cB_{H,X}^\kappa\)
be the same pullback of
\(\cB_{H,X}\). Every result in this chapter about \(h_{H,X}\) applies to
\(h_{H,X}^\kappa\) (as long as the statement makes sense) due to it being
defined via pullback from
\(\Bun_{Z_{\FRM,H}}\).

The canonical map \(\nu_H\colon\FRC_{\FRM,H}\to\FRC_{\FRM}\) induces the commutative
diagram
\begin{equation}
    \begin{tikzcd}
        \Stack{\FRC_{\FRM,H}/Z_{\FRM}^\kappa} \ar[r]\ar[d] & \Stack{\FRC_{\FRM}/Z_{\FRM}} \ar[d]\\
        \Stack{\FRA_{\FRM,H}/Z_{\FRM}^\kappa} \ar[r] & \Stack{\FRA_{\FRM}/Z_{\FRM}}
    \end{tikzcd}
\end{equation}
hence further induces the diagram
\begin{equation}
    \begin{tikzcd}
        \nu_{\cA}\colon\cA_{H,X}^\kappa \ar[r]\ar[d] & \cA_X \ar[d]\\
        \nu_{\cB}\colon\cB_{H,X}^\kappa \ar[r] & \cB_X
    \end{tikzcd}
    \nomenclature[\(nu_A"cal \)]{\(\nu_{\cA}\)}{the endoscopic transfer map \(\cA_{H,X}^\kappa\to\cA_X\)}
    \nomenclature[\(nu_B"cal \)]{\(\nu_{\cB}\)}{the endoscopic transfer map \(\cB_{H,X}^\kappa\to\cB_X\)}
\end{equation}

Recall that we also have the coarse endoscopic monoid \(\FRM_H^\star\), which
canonically embeds into \(\FRM_H\). By pulling back, we have induced maps
\(h_{H,X}^{\kappa\star}\), \(\nu_\cA^\star\), \(\nu_\cB^\star\), and so on.

\subsection{}
To further our discussion, we need the following lemma, which can be easily
deduced from Zariski's Main Theorem:
\begin{lemma} [\cite{Ng06}*{Lemme~7.3}]
    \label[lemma]{lem:lifting_section_to_finite_cover_on_source}
    Let \(S\) be a normal integral and separated \(k\)-scheme. Let \(v\colon
    \tilde{V}\to V\)
    be a finite morphism of \(S\)-schemes and \(h\colon S\to V\) a
    section of \(V\to S\). Suppose there is an open dense subset \(S'\subset S\)
    over which \(h\) lifts to a section \(h'\colon S'\to\tilde{V}\x_S S'\).
    Then \(h'\) extends uniquely to \(S\) such that \(v\circ h'=h\).
\end{lemma}

\begin{proposition}
    \label[proposition]{prop:endoscopic_transfer_finite_and_unramified}
    The map
    \begin{align}
        \nu_\cA^\heartsuit\colon
        \cA_{H,X}^{\kappa,G\hy\heartsuit}\longto\cA_X^\heartsuit
        \nomenclature[\(.G_heartsuit \)]{\((\cdot)^{G\hy\heartsuit}\)}{constructions related to
        the preimage of \((\cdot)^\heartsuit\) in \((\cdot)_{H,X}^\kappa\)}
    \end{align}
    is finite and unramified, and the same is true for
    \(\nu_\cA^{\star\heartsuit}\). Moreover, \(\nu_\cA^\heartsuit\) and
    \(\nu_\cA^{\star\heartsuit}\) have the same image.
\end{proposition}
\begin{proof}
    The \(\pi_0(\kappa)\)-torsor \(\OGT_\kappa\colon X_\kappa\to X\) induces 
    commutative diagram
    \begin{equation}
        \begin{tikzcd}
            \bar{\bT}_{\bM,\bH}\x X_\kappa\ar[r]\ar[d] & \bar{\bT}_{\bM}\x
            X_\kappa\ar[r, equal]\ar[d] & \bar{\bT}_{\bM}\x X_\kappa\ar[d]\\
            \FRC_{\FRM,H} \ar[r] & \FRC_{\FRM,H}' \ar[r] &
            \FRC_{\FRM}
        \end{tikzcd}
    \end{equation}
    where the left and middle vertical maps are obtained by taking
    \(\bW_{\bH}\rtimes\pi_0(\kappa)\)-quotient,
    and the right one is by taking \(\bW\rtimes\pi_0(\kappa)\)-quotient, and the two
    groups are connected by \Cref{lem:kappa_G_H_compatibility}.

    Let \(a\in\cA_X^\heartsuit(\bar{k})\), and let \(U_a\subset\breve{X}\) be the
    open subset whose image under \(a\) is contained in
    \(\Stack*{\FRC_{\FRM}^{\x,\rss}/Z_{\FRM}}\).
    We also fix a geometric point \(\infty\in U_a(\bar{k})\). 
    Over \(U_a\), the vertical
    maps in the above commutative diagram are respective torsors of groups
    \(\bW_{\bH}\rtimes\pi_0(\kappa)\) and \(\bW\rtimes\pi_0(\kappa)\). Then the
    restriction of \(a\) to \(U_a\) lifts to a section
    \(U_a\to\Stack*{\FRC_{\FRM,H}'/Z_{\FRM}}\x_X U_a\) if and only if the image of the monodromy
    \begin{align}
        \OGT_a^\bullet\colon \pi_1(U_a,\infty)\longto
        \bW\rtimes\pi_0(\kappa)
    \end{align}
    is conjugate to subgroup \(\bW_{\bH}\rtimes\pi_0(\kappa)\). Since the number
    of such subgroups is finite, so is the number of such lifts. Also over
    \(U_a\), the map
    \(\Stack*{\FRC_{\FRM,H}/Z_{\FRM}^\kappa}\to\Stack*{\FRC_{\FRM,H}'/Z_{\FRM}}\) 
    is an isomorphism. At each \(\bar{v}\in\breve{X}-U_a\), there are only
    finitely many ways to extend a
    \(\breve{X}_{\bar{v}}^\bullet\)-point of
    \(\Stack*{\FRC_{\FRM,H}/Z_{\FRM}^\kappa}\) to a \(\breve{X}_{\bar{v}}\)-point
    lying over a fixed \(\breve{X}_{\bar{v}}\)-point in
    \(\Stack*{\FRC_{\FRM,H}'/Z_{\FRM}}\), according to
    \Cref{lem:arc_lifting_to_endoscopic_quotient} (note that over
    \(\breve{X}_{\bar{v}}\) the monoid \(\FRM\) is necessarily split).
    Thus, \(\nu_\cA^\heartsuit\) is quasi-finite.
    The argument of \Cref{lem:arc_lifting_to_endoscopic_quotient} also shows
    that if \(a\) lifts to a point in \(\cA_{H,X}^\kappa\), it must lift to a
    point in \(\cA_{H,X}^{\kappa\star}\) as well, and this proves the last claim
    of the proposition.

    We then show that \(\nu_\cA^\heartsuit\) is proper using valuative criteria.
    Let \(R\) be a discrete
    valuation ring and \(S=\Spec{R}\). Let \(\eta\in S\) be the generic point.
    Now let \(a\) be an \(S\)-point instead of a \(\bar{k}\)-point of
    \(\cA_X^\heartsuit\), and \(U_a\subset X\x S\) is defined as above. Since
    \(U_a\) is normal and integral and the map of stacks
    \begin{align}
        \label{eqn:tmp_endoscopic_finite_unramifiedness_of_mHitchin}
        \Stack*{\FRC_{\FRM,H}^{\x,G\hy\rss}/Z_{\FRM}^\kappa}\longto \Stack*{\FRC_{\FRM}^{\x,\rss}/Z_{\FRM}}
    \end{align}
    is representable by schemes and finite \'etale (because
    \(\bT_{\bM,\bH}/\bZ_{\bM}^\kappa\to \bT_{\bM}/\bZ_{\bM}\) is an
    isomorphism), a lifting of \(a|_{U_a}\) over \(U_a\x_S\eta\) can be extended
    to a lifting
    over \(U_a\) by \Cref{lem:lifting_section_to_finite_cover_on_source}. Then
    since the complement of \(U_a\cup (X\x\eta)\) has codimension \(2\) in
    normal scheme \(X\x S\) and with no codimension-\(3\) subsets, any
    \(Z_{\FRM}^\kappa\)-torsor on \(U_a\cup
    (X\x\eta)\) can be uniquely extended over \(X\x S\), and since
    \(\FRC_{\FRM,H}\) is affine over \(X\), we can uniquely lift \(a\) over
    \(X\x S\). This proves properness.

    Finally, we show that \(\nu_\cA^\heartsuit\) is unramified. Let
    \(a_H,a_H'\in\cA_{H,X}^{\kappa,G\hy\heartsuit}(\bar{k}[\epsilon]/\epsilon^2)\) be
    two points with the same image
    \(a\in\cA_X^\heartsuit(\bar{k}[\epsilon]/\epsilon^2)\) as well as in
    \(\cA_{H,X}^{\kappa,G\hy\heartsuit}(\bar{k})\). The map
    \eqref{eqn:tmp_endoscopic_finite_unramifiedness_of_mHitchin} is \'etale, so
    \(a_H\) and \(a_H'\) are necessarily equal when restricted to
    \(U_a\cup\breve{X}\). Since
    \(\cO_{U_a}\) contains \(\cO_{\breve{X}}[\epsilon]/\epsilon^2\) as a
    subalgebra, \(a_H\) and \(a_H'\) must be equal over all
    \(\breve{X}\x_{\Spec{\bar{k}}}\Spec{\bar{k}[\epsilon]/\epsilon^2}\), and we are
    done.
\end{proof}

\subsection{}
Suppose \(b\in\cB_\gg(\bar{k})\) induces a very \(G\)-ample boundary divisor written as
\begin{align}
    \lambda_b=\sum_{\bar{v}\in X(\bar{k})}\lambda_{\bar{v}}\cdot\bar{v}.
\end{align}
Similarly, \(b_H\in\cB_{H,\gg}^\kappa(\bar{k})\) induces
very \(H\)-ample boundary divisor
\begin{align}
    \lambda_{H,b_H}=\sum_{\bar{v}\in X(\bar{k})}\lambda_{H,\bar{v}}\cdot\bar{v}.
\end{align}
Suppose \(b\) is the image of \(b_H\).
For each \(\bar{v}\in X(\bar{k})\), \(-w_{H,0}(\lambda_{H,\bar{v}})\) is
one of the \(\dual{\bH}\)-highest weights appearing in the decomposition
of irreducible \(\dual{\bG}\)-representation with highest-weight
\(-w_0(\lambda_{\bar{v}})\) into irreducible
\(\dual{\bH}\)-representations.
Let
\begin{align}
    r_H^G(b_H)=\Pair{\rho}{\lambda_b}
            -\Pair{\rho_H}{\lambda_{H,b_H}}.
\end{align}
This way we obtain a locally constant function
\begin{align}
    r_H^G\colon \cB_{H,X}^\kappa\longto \bbN,
\end{align}
and we use the same notation for its pullback to \(\cA_{H,X}^\kappa\),
\(\cM_{H,X}^\kappa\), etc.

By \eqref{eqn:dimension_of_restricted_mH_base}, we have
\begin{align}
    \label{eqn:difference_in_dimension_of_mH_bases_fixed_divisor}
    \dim{\cA_b}-\dim{\cA_{H,b_H}^\kappa}=r_H^G(b_H)-(r-r_H)(g_X-1),
\end{align}
where \(r_H\) is the semisimple rank of \(H\).
So the image of \(\cA_{H,b_H}^{\kappa,G\hy\heartsuit}\) in \(\cA_b^\heartsuit\) is a
closed subscheme of codimension \(r_H^G(b_H)-(r-r_H)(g_X-1)\).
In case \(b_H\) is not very \(H\)-ample (but \(b\) is still very \(G\)-ample),
we still have inequality according to
\eqref{eqn:dimension_estimate_of_restricted_mH_base_non_ample}:
\begin{align}
    \label{eqn:weak_estimate_difference_in_dimension_of_mH_bases_fixed_divisor}
    \dim{\cA_b}-\dim{\cA_{H,b_H}^\kappa}\ge r_H^G(b_H)-r(g_X-1)-r_H.
\end{align}

\subsection{}
The dimension formula
\eqref{eqn:difference_in_dimension_of_mH_bases_fixed_divisor}
on a fiberwise basis is not enough as we will need a similar formula for the
whole mH-base. However, things become trickier when \(b_H\) is not fixed.

Suppose \(\FRM\in\FM_0(G^\SC)\) and \(\FRA_{\FRM}\) is of standard type. If
\(\cL_\kappa\) is a \(Z_{\FRM}^\kappa\)-torsor, then it induces a
\(Z_{\FRM,H}\)-torsor \(\cL_H\). We say \(\cL_\kappa\) is very \(H\)-ample
if \(\cL_H\) is. It also induces a \(Z_{\FRM}\)-torsor which
we denote by \(\cL\). However, if \(\cL\) is very \(G\)-ample, there is no
\emph{guarantee} that \(\cL_\kappa\) is very \(H\)-ample, or vice versa.
Moreover, since \(Z_\FRM^\kappa\) is
usually smaller than \(Z_{\FRM,H}\) in rank, the stack \(\cB_{H,X}^\kappa\)
may not be smooth even if \(\cB_{H,X}\) is.
This forces us to use alternative
assumptions other than \(G\)-ampleness or \(H\)-ampleness.
The best way to do this is to use deformation theory, which we shall discuss in
more details in \Cref{chap:deformation} (see in particular
\Cref{sec:the_case_of_mH_base}).

At \(a_H\in\cA_{H,X}^\kappa(\bar{k})\) whose \(Z_\FRM^\kappa\)-torsor is
\(\cL_\kappa\), the action of \(Z_\FRM^\kappa\) on
\(\FRM_H\) induces tangent action
\begin{align}
    \La{z}_{\FRM}^\kappa\longto
    \TanB_{a_H}\FRC_{\FRM,H,\cL_\kappa}\simeq\FRC_{\FRM,H,\cL_\kappa}.
\end{align}
As we will see in \Cref{sec:the_case_of_mH_base}, \(\cA_{H,X}^\kappa\)
is smooth at \(a_H\) if the following cohomological condition holds:
\begin{align}
    \label{eqn:dim_of_endoscopic_strata_cohom_assumption}
    \RH^1(\breve{X}, \FRC_{\FRM,H,\cL_\kappa}/\La{z}_\FRM^\kappa)=0.
\end{align}
This condition can be met, for example, if \(\cL_\kappa\) is very \(H\)-ample
and \(\FRA_{\FRM,H,\cL_\kappa}\) also has vanishing first
cohomology. We shall assume
\eqref{eqn:dim_of_endoscopic_strata_cohom_assumption} holds for \(a_H\), and
then the dimension of \(\cA_{H,X}^\kappa\) at \(a_H\) is the expected dimension
\begin{align}
    \dim_{a_H}\cA_{H,X}^\kappa
    &=\RH^0(\breve{X}, \FRC_{\FRM,H,\cL_\kappa}/\La{z}_\FRM^\kappa)\\
    &=\deg\FRC_{\FRM,H,\cL_\kappa}-\deg\La{z}_\FRM^\kappa-(\rk\FRC_{\FRM,H,\cL_\kappa}-\rk\La{z}_\FRM^\kappa)(g_X-1).
\end{align}
Obviously \(\deg\La{z}_\FRM^\kappa=0\), and it is easy to see that
\begin{align}
    \deg\FRC_{\FRM,H,\cL_\kappa}&=\dim_b\cB_X+\Pair{\rho_H}{\lambda_{H,b_H}},\\
    \rk\FRC_{\FRM,H,\cL_\kappa}-\rk\La{z}_\FRM^\kappa&=\rk
    \FRC_{\FRM,\cL}-\rk\La{z}_\FRM=r.
\end{align}
Thus, we have
\begin{align}
    \dim_{a_H}\cA_{H,X}^\kappa
    =\dim_b\cB_X+\Pair{\rho_H}{\lambda_{H,b_H}}-r(g_X-1).
\end{align}

\subsection{}
If \(a\in\cA_X(\bar{k})\) is the image of \(a_H\), we have the following diagram
of tangent actions:
\begin{equation}
    \begin{tikzcd}
        \La{z}_\FRM^\kappa \ar[r]\ar[d] &
        \TanB_{a_H}\FRC_{\FRM,H,\cL_\kappa}\simeq \FRC_{\FRM,H,\cL_\kappa} \ar[d]\\
        \La{z}_\FRM \ar[r] & \TanB_{a}\FRC_{\FRM,\cL}\simeq\FRC_{\FRM,\cL}
    \end{tikzcd}.
\end{equation}
It is clear that the induced map of coherent sheaves on \(\breve{X}\)
\begin{align}
    \FRC_{\FRM,H,\cL_\kappa}/\La{z}_\FRM^\kappa\longto
    \FRC_{\FRM,\cL}/\La{z}_\FRM
\end{align}
is generically an isomorphism. This means that assuming
\eqref{eqn:dim_of_endoscopic_strata_cohom_assumption}, we also have
\begin{align}
    \RH^1(\breve{X},\FRC_{\FRM,\cL}/\La{z}_\FRM)=0.
\end{align}
By our general discussion in \Cref{sec:the_case_of_mH_base} again,
\(\cA_X\) is smooth at \(a\), and so \eqref{eqn:dimension_of_mH_base_at_a_point}
holds regardless of whether \(a\) is very \(G\)-ample or not. Consequently,
\begin{align}
    \label{eqn:difference_in_dimension_of_mH_bases}
    \dim_{a}\cA_X-\dim_{a_H}\cA_{H,X}^\kappa=r_H^G(a_H).
\end{align}
If \eqref{eqn:dim_of_endoscopic_strata_cohom_assumption} does not hold,
we still have estimate
\begin{align}
    \label{eqn:difference_estimate_in_dimension_of_mH_bases_non_ample}
    \dim_{a}\cA_X-\dim_{a_H}\cA_{H,X}^\kappa\ge r_H^G(a_H)-rg_X,
\end{align}
as long as we assume that \(a\) is very \(G\)-ample.

\subsection{}
Note that the discussion about \eqref{eqn:difference_in_dimension_of_mH_bases}
applies equally well to \(\cA_{H,X}^{\kappa\star}\). Since \(\nu_\cA^\heartsuit\) and
\(\nu_\cA^{\star\heartsuit}\) have the same image, and since the
cohomological condition \eqref{eqn:dim_of_endoscopic_strata_cohom_assumption} is
easier to satisfy for \(\cA_{H,X}^{\kappa\star}\), it makes the
coarse version of \eqref{eqn:difference_in_dimension_of_mH_bases} easier to
use. In fact, it is the sole reason why we consider coarse endoscopic monoid
\(\FRM_H^\star\) at all.

Replace endoscopic groups by Levi subgroups, then the above discussions apply
equally well to mH-fibrations associated with Levi monoids.

\subsection{}
Although the mH-total-stacks of \(G\) and \(H\) do not have direct connection,
there is a canonical map between their Picard stacks. Let
\(a_H\in\cA_{H,X}^{\kappa,G\hy\heartsuit}(\bar{k})\) and its image is
\(a\in\cA_X^\heartsuit(\bar{k})\). Recall that the map
\(\nu_H\colon\FRC_{\FRM,H}\to\FRC_{\FRM}\) induces a homomorphism
\(\nu_H^*\FRJ_{\FRM}\to\FRJ_{\FRM,H}\), therefore we have a homomorphism of
commutative group schemes \(\FRJ_a\to\FRJ_{H,a_H}\) over \(\breve{X}\). Since
\(a\in\cA_X^\heartsuit\), this is generically an isomorphism hence
also injective. Therefore, we have a surjective map
\begin{align}
    \cP_a\longto\cP_{H,a_H}.
\end{align}
Its kernel is the affine group 
\(\RH^0(\breve{X},\FRJ_{H,a_H}/\FRJ_a)\), which we denote by \(\cR_{H,a_H}^G\).
Using \Cref{cor:dimension_of_P_a}, we have
\begin{align}
    \dim{\cR_{H,a_H}^G} = r_H^G(a_H).
\end{align}
In particular, it is locally constant over
\(\cA_{H,X}^{\kappa,G\hy\heartsuit}(\bar{k})\). Let \(\FRJ_{H,a_H}^\flat\) be
the N\'eron model of \(\FRJ_{H,a_H}\), then the composition
\(\FRJ_a\to\FRJ_{H,a_H}^\flat\) is generically an isomorphism, so it is also the
N\'eron model of \(\FRJ_a\). Thus, we have an exact sequence
\begin{align}
    1\longto \cR_{H,a_H}^G\longto \cR_a\longto\cR_{H,a_H}\longto 1,
\end{align}
where \(\cR_a=\ker(\cP_a\to\cP_a^\flat)\) and
\(\cR_{H,a_H}=\ker(\cP_{H,a_H}\to\cP_{H,a_H}^\flat)\), and so
\begin{align}
    \label{eqn:difference_in_delta_invariants}
    \delta_a-\delta_{H,a_H}=r_H^G(a_H),
\end{align}
where \(\delta_a=\dim{\cR_a}\) and \(\delta_{H,a_H}=\dim{\cR_{H,a_H}}\).

\subsection{}
We have established local models of singularity
(\Cref{thm:local_singularity_model_weak,thm:local_singularity_model_main}) for
\(\cM_{H,X}\). \Cref{thm:local_singularity_model_weak} can be applied to
\(\cM_{H,X}^\kappa\) by a simple base change, but the deformation argument for
\Cref{thm:local_singularity_model_main} will be trickier to pull off due to the
difference in rank between \(Z_\FRM^\kappa\) and \(Z_{\FRM,H}\), even though
we expect the same idea should still work with some adjustment. Fortunately,
we do not have to do so: for endoscopic groups,
\Cref{thm:local_singularity_model_weak} will suit our needs perfectly well and
we will leave it at that.

\chapter{Deformations} 
\label{chap:deformation}

In this chapter we study deformation questions related to mHiggs bundles.
It is split into two parts: the first one involves studying deformation in the
mH-base, and the second is a local model
of singularities for the total stack of mH-fibrations. The first part, though
has some intrinsic value itself, mostly serves as a stepping stone towards the
goal of the second part and to make the argument easier to digest because it is
quite technical. The reader may compare to
\cite{Ng10}*{\S\S~4.13, 4.14} in the Lie algebra case.

\section{Deformation of Mapping Stacks} 
\label{sec:deformation_of_mapping_stacks}

We start with some generalities concerning deformations.
Consider variety \(M\) over \(\bar{k}\) with a smooth group \(G\) acting on
it. Consider the quotient
\begin{align}
    \chi\colon M\longto \Stack{M/G}.
\end{align}
We have distinguished triangle of cotangent complexes
\begin{align}
    \chi^*L_{\Stack{M/G}}\longto L_M\longto L_{M/\Stack{M/G}}\stackrel{+1}{\longto}.
    \nomenclature[\(L_Y_S \)]{\(L_{Y/S}\)}{the relative cotangent complex of \(Y\) over \(S\)}
\end{align}
Since \(\chi\) is a \(G\)-torsor, we have \(L_{M/\Stack{M/G}}\simeq
\cO_{M}\otimes\La{g}^*[0]\), and \(L_M\) is \(G\)-equivariant.
Descending to \(\Stack{M/G}\), we get
\begin{align}
    L_{\Stack{M/G}}\longto L_M/G \longto
    M\x^G\La{g}^*[0]\stackrel{+1}{\longto}.
\end{align}

\subsection{}
Let \(X\) be a \(\bar{k}\)-scheme of finite type and \(m=(E,\phi)\) be an
\(X\)-point of \(\Stack{M/G}\), then the derived pullback \(\LDF{m}^*L_{\Stack{M/G}}\)
    \nomenclature[\(L. \)]{\(\LDF{(\cdot)}\)}{the left-derived functor of a right-exact functor}
is isomorphic to the
cone of the map of complexes
\begin{align}
    \bigl(\LDF\phi^*L_{E\x^G M}\longto \ad(E)^*[0]\bigr)[-1].
\end{align}
Let \(\cT^\bullet=\RIHom_X(\LDF{m}^*L_{\Stack{M/G}},\cO_X)\),
    \nomenclature[\(R. \)]{\(\RDF{(\cdot)}\)}{the right-derived functor of a left-exact functor}
    \nomenclature[\(T"cal^bullet \)]{\(\cT^\bullet\)}{the derived tangent complex
    controlling the deformation of a point in certain mapping stack}
then it is isomorphic to
the cone of the map between complexes
\begin{align}
    \label{eqn:sT_as_cone}
    \ad(E)[0]\longto \RIHom_X(\LDF\phi^*L_{E\x^G M},\cO_X).
\end{align}
Since \(L_M\) is supported on degrees \((-\infty,0]\), so is
\(\LDF\phi^*L_{E\x^G M}\). Since \(\IHom_X(-,\cO_X)\) is left exact, we see
that \(\RIHom_X(\LDF\phi^*L_{E\x^G M},\cO_X)\) is supported on \([0,+\infty)\).
Therefore, by directly consulting the construction of mapping cones,
\(\cT^\bullet\) is isomorphic to the complex
\begin{align}
    \ad(E)[1]\longto  \RIHom_X(\LDF\phi^*L_{E\x^G M},\cO_X),
\end{align}
where the arrow is just the map (of coherent sheaves)
\begin{align}
    \ad(E)\longto
    \IHom_X(\LDF\phi^*L_{E\x^G M},\cO_X)
\end{align}
induced by \eqref{eqn:sT_as_cone}, which can also be seen as the derivative
of the \(G\)-action. In particular, \(\cT^\bullet\) is supported on degrees
\([-1,+\infty)\).

Similar to the Lie algebra case in \cite{Ng10}*{\S~4.14}, the deformation of
\(m\in\IHom(X,\Stack{M/G})\) is controlled by
\(\cT^\bullet\), such that the obstruction space is the hypercohomology group
\(\RH^1(X,\cT^\bullet),\) the tangent space is \(\RH^0(X,\cT^\bullet)\), and the
infinitesimal automorphism group is \(\RH^{-1}(X,\cT^\bullet)\).

\subsection{}
Now we compute these cohomology groups using \v{C}ech cohomology
when \(X\) is a curve. We shall assume that the generic points of \(X\) are
sent to the smooth locus of \(\Stack{M/G}\) under \(m\). Since over the smooth
locus of \(\Stack{M/G}\), the complex \(L_{\Stack{M/G}}\) is a bounded complex
with locally-free components, its derived pullback is given by the naive one,
and so we see that each quasi-coherent sheaf \(\RH^i(\cT^\bullet)\) is supported
on finitely many points on \(X\) if \(i\ge 1\).

Let \(\sX=\Set{X_i}_{i\in I}\) be a finite Zariski open affine covering of
\(X\). The forgetful functor from the category of sheaves of \(\cO_X\)-modules
to presheaves of \(\cO_X\)-modules is left exact, whose right-derived functors
\(\iH^q\) are given by (\(\cF\) is any sheaf of \(\cO_X\)-modules):
\begin{align}
    \iH^q(\cF)\colon U\subset X\longmapsto \RH^q(U,\cF).
\end{align}
For any finite subset \(\Set{i_0,i_1,\cdots, i_p}\) of \(I\), we let
\begin{align}
    X_{i_0,\ldots,i_p}=X_{i_0}\cap\cdots\cap X_{i_p}.
\end{align}
Then for each \(q\ge 0\) we have \v{C}ech double complex
\(\vC^\bullet(\sX,\iH^q(\cT^\bullet))\):
\begin{equation}
    \begin{tikzcd}
                 & \vdots       & \vdots       & \vdots       & \\
        0 \ar[r] & \prod_{ijk}\RH^q(X_{ijk},\cT^{-1})\ar[r]\ar[u] & \prod_{ijk}\RH^q(X_{ijk},\cT^{0})\ar[r]\ar[u] & \prod_{ijk}\RH^q(X_{ijk},\cT^{1})\ar[u]\ar[r] & \cdots \\
        0 \ar[r] & \prod_{ij}\RH^q(X_{ij},\cT^{-1})\ar[r]\ar[u] & \prod_{ij}\RH^q(X_{ij},\cT^{0})\ar[r]\ar[u] & \prod_{ij}\RH^q(X_{ij},\cT^{1})\ar[u]\ar[r] & \cdots \\
        0 \ar[r] & \prod_{i}\RH^q(X_{i},\cT^{-1})\ar[r]\ar[u] & \prod_{i}\RH^q(X_{i},\cT^{0})\ar[r]\ar[u] & \prod_{i}\RH^q(X_{i},\cT^{1})\ar[u]\ar[r] & \cdots \\
                 & 0\ar[u] & 0\ar[u] & 0\ar[u] &
    \end{tikzcd}.
\end{equation}
The general theory of \v{C}ech cohomologies
(\cite{StacksP}*{\href{https://stacks.math.columbia.edu/tag/01FP}{Tag 01FP}}) shows that there is a spectral sequence
\begin{align}
    \ESP_2^{p,q}=\RH^p(\Tot(\vC^\bullet(\sX,\iH^q(\cT^\bullet))))\Longrightarrow
    \RH^{p+q}(X,\cT^\bullet),
\end{align}
where \(\Tot\) means taking the total complex associated with a double complex.
Since \(X_{i_0,\ldots,x_p}\) is affine for any \(p\) and any
\(i_0,\ldots,i_p\), and \(\cT^i\) is quasi-coherent for any \(i\), we know that
the above double complex vanishes completely for all \(q>0\). Thus, the above
spectral sequence degenerates at \(\ESP_2\)-page, and so we have
\begin{align}
    \RH^p(X,\cT^\bullet)\simeq\RH^p(\Tot(\vC^\bullet(\sX,\iH^0(\cT^\bullet)))).
\end{align}

Now we use the spectral sequence of double complexes to compute
\begin{align}
    \RH^p(\Tot(\vC^\bullet(\sX,\iH^0(\cT^\bullet)))).
\end{align}
The \(0\)-th page is
\begin{equation}
    \begin{tikzcd}
                 & \cdots       & \cdots       & \cdots       & \\
        0 \ar[r] & \prod_{ijk}\RH^0(X_{ijk},\cT^{-1})\ar[r] & \prod_{ijk}\RH^0(X_{ijk},\cT^{0})\ar[r] & \prod_{ijk}\RH^0(X_{ijk},\cT^{1})\ar[r] & \cdots \\
        0 \ar[r] & \prod_{ij}\RH^0(X_{ij},\cT^{-1})\ar[r] & \prod_{ij}\RH^0(X_{ij},\cT^{0})\ar[r] & \prod_{ij}\RH^0(X_{ij},\cT^{1})\ar[r] & \cdots \\
        0 \ar[r] & \prod_{i}\RH^0(X_{i},\cT^{-1})\ar[r] & \prod_{i}\RH^0(X_{i},\cT^{0})\ar[r] & \prod_{i}\RH^0(X_{i},\cT^{1})\ar[r] & \cdots \\
    \end{tikzcd}
\end{equation}
Since all \(X_{i_0,\ldots,i_p}\) are affine, we may directly compute the
\(1\)-st page to be
\begin{equation}
    \begin{tikzcd}
        \vdots       & \vdots       & \vdots       &\\
        \prod_{ijk}\RH^0(X_{ijk},\RH^{-1}(\cT^\bullet))\ar[u] & \prod_{ijk}\RH^0(X_{ijk},\RH^{0}(\cT^\bullet))\ar[u] & \prod_{ijk}\RH^0(X_{ijk},\RH^{1}(\cT^\bullet))\ar[u] & \vdots\\
        \prod_{ij}\RH^0(X_{ij},\RH^{-1}(\cT^\bullet))\ar[u] & \prod_{ij}\RH^0(X_{ij},\RH^{0}(\cT^\bullet))\ar[u] & \prod_{ij}\RH^0(X_{ij},\RH^{1}(\cT^\bullet))\ar[u]  & \vdots \\
        \prod_{i}\RH^0(X_{i},\RH^{-1}(\cT^\bullet))\ar[u] & \prod_{i}\RH^0(X_{i},\RH^{0}(\cT^\bullet))\ar[u] & \prod_{i}\RH^0(X_{i},\RH^{1}(\cT^\bullet))\ar[u]  & \vdots \\
         0\ar[u] & 0\ar[u] & 0\ar[u] 
    \end{tikzcd}
\end{equation}
Observe that the \(j\)-th column above is exactly the \v{C}ech complex of
quasi-coherent sheaf \(\RH^j(\cT^\bullet)\) with respect to covering \(\sX\), so
its cohomologies are just \(\RH^i(X,\RH^j(\cT^\bullet))\). Since \(X\) is a
curve, \(\RH^i(X,\RH^j(\cT^\bullet))=0\) for all \(i\ge 2\), and since for any
\(j\ge 1\) the sheaf \(\RH^j(\cT^\bullet)\) is supported over finitely many
points, we also have \(\RH^1(\RH^j(\cT^\bullet))=0\) for \(j\ge 1\). Thus, the
\(2\)-nd page has only two non-zero rows:
\begin{equation}
    \begin{tikzcd}
        \RH^1(X,\RH^{-1}(\cT^\bullet)) & \RH^1(X,\RH^{0}(\cT^\bullet)) & 0 & 0 & \cdots\\
        \RH^0(X,\RH^{-1}(\cT^\bullet)) & \RH^0(X,\RH^{0}(\cT^\bullet)) &
        \RH^0(X,\RH^{1}(\cT^\bullet)) & \RH^0(X,\RH^2(\cT)) & \cdots
    \end{tikzcd}
\end{equation}
and all differentials on all pages hereafter are zero, because the
vertical component of any arrow will go upwards at least two rows. Thus, the
spectral sequence degenerates at the \(2\)-nd page, and we have canonical
isomorphisms
\begin{align}
    \RH^{i}(X,\cT^\bullet)\simeq \RH^0(X,\RH^{i}(\cT^\bullet)),\quad i=-1\text{
    or }i\ge 2,
\end{align}
as well as exact sequences
\begin{align}
    0\longto \RH^1(X,\RH^{-1}(\cT^\bullet))\longto \RH^0(X,\cT^\bullet)\longto
    \RH^0(X,\RH^0(\cT^\bullet))\longto 0,\\
    0\longto \RH^1(X,\RH^{0}(\cT^\bullet))\longto \RH^1(X,\cT^\bullet)\longto
    \RH^0(X,\RH^1(\cT^\bullet))\longto 0.
\end{align}

\subsection{}
Note that everything so far in this section still holds in the relative setting
over a base algebraic stack \(S\) where \(M\to S\) is a flat \(S\)-scheme with reduced
fibers, \(G\to S\) is a smooth algebraic group scheme over \(S\), and \(X\to S\)
is a flat \(S\)-family of smooth curves.

Similarly, we may also replace \(X\) with a formal disc \(X_v\) at a point
\(v\in X\). In this case \(\RH^1(X_v,\cF)=0\) for any quasi-coherent sheaf
\(\cF\). Let \(\cT_v^\bullet\) be the analogue of \(\cT^\bullet\) on \(X_v\),
then we have for all \(i\ge -1\)
\begin{align}
    \RH^i(X_v,\cT_v^\bullet)\simeq\RH^0(X_v,\RH^i(\cT_v^\bullet)).
\end{align}
Let \(\iota\colon X_v\to X\) be the natural map, then it is flat and the
(non-derived) functor \(\iota^*\) is exact. Thus, we have
\(\cT_v^\bullet\simeq\iota^*\cT^\bullet\), and the natural map
\begin{align}
    \RH^i(\cT^\bullet)\otimes_{\cO_X}\cO_{v}\longto \RH^i(\cT_v^\bullet)
\end{align}
is an isomorphism of \(\cO_v\)-modules. In particular, when \(i\ge 1\), since
\(\RH^i(\cT^\bullet)\) is finitely supported, we have injective map
\begin{align}
    \RH^0(X,\RH^i(\cT^\bullet))\longto\prod_{v}\RH^0(X_v,\RH^i(\cT_v^\bullet)),
\end{align}
where \(v\) ranges over the points \(v\in X\) such that \(m(v)\) is singular in
\(\Stack{M/G}\). Note that when \(i=1\),
the right-hand side is precisely the obstruction
space of deforming the \(\prod_v X_v\)-arc in \(\Stack{M/G}\)
induced by \(m\), in other words, the local obstruction space.

\subsection{}
We now look at the obstruction space \(\RH^1(X,\cT^\bullet)\). We have seen
above that the quotient
\(\RH^0(X,\RH^1(\cT^\bullet))\) is precisely the
space of local obstructions. As a result, if we can show that
\begin{align}
    \label{eqn:local_model_singularity_primitive}
    \RH^1(X,\RH^{0}(\cT^\bullet))=0,
\end{align}
we will be able to prove that the global obstruction to deforming \(m\) is
completely determined by its local obstructions. If, in addition, \(M\) is
smooth, then there is no local obstruction, so the mapping stack from \(X\) to
\(\Stack{M/G}\) is also smooth at
\(m\). On the other hand, if \(M\) is not smooth, then in general one should not
expect the mapping stack to be smooth,  but can hope to understand the
singularities using the arc space of \(M\). In other words, we hope to establish
a \inotion{local model of singularity} for the mapping stack
\(\IHom(X,\Stack{M/G})\) by proving cohomological statement
\eqref{eqn:local_model_singularity_primitive}.

To make it more concrete, consider two-step truncation of \(\cT^\bullet\), viewed
as a map of two quasi-coherent sheaves:
\begin{align}
    \cT^{\le 0}\colon\ad(E)\longto \IHom_X(\phi^*\CoTB_{E\x^G M}^1,\cO_X),
\end{align}
then we are interested in showing that
\begin{align}
    \label{eqn:needed_vanishing_statement_for_deformation}
    \RH^1(X,\coker(\cT^{\le 0}))=0.
\end{align}
This is a generalization of \cite{Ng10}*{\S~4.14} for additive Hitchin fibrations
where \(M=\La{g}_D\). In that case, \(M\) is smooth, so there is no local
obstruction, and if \(D\) is sufficiently ample, the global obstructions also
vanish using Serre duality.

\subsection{}
In the next two sections we will apply our general framework above to two cases.

The first case is where \(M=\FRC_\FRM\) and \(G=Z_\FRM\) for a very flat
reductive monoid \(\FRM\) whose abelianization is an affine space. In this case
\(M\) is smooth, so our goal will be showing that \(\cA_X\) is smooth under mild
assumptions. However, these assumptions, although mild, are not more so than
very \(G\)-ampleness, thus it can not be used to improve the dimension formulae
in \Cref{sec:dimensions}. Instead, it is mainly used to facilitate our
discussion in the second case.

The second case is where \(M=\FRM\) itself and \(G\) is
the reductive group we have been considering for mH-fibrations (which we also
have been using the same notation \(G\)). In this
case \(M\) is not smooth in general, and we are going to prove the local model of
singularity \Cref{thm:local_singularity_model_weak,thm:local_singularity_model_main}.


\section{The Case of mH-base} 
\label{sec:the_case_of_mH_base}

In this section, let \(\FRM\in\FM(G^\SC)\) be such that \(\FRA_\FRM\) is of
standard type. The invariant
quotient \(\FRC_\FRM\simeq \FRA_\FRM\x \FRC\) can be viewed as vector bundle
over \(X\), so its cotangent complex may be identified with the vector dual
\(\FRC_\FRM^*\). Here we abuse the notation to still use \(\FRC_\FRM^*\) for its
pullback from \(X\) to \(\FRC_\FRM\). We denote by
\(\Stack{\FRC_\FRM^*/Z_\FRM}\) the descent of \(\FRC_\FRM^*\) to
\(\Stack{\FRC_\FRM/Z_\FRM}\). Since \(Z_\FRM\) is commutative, the descent of
bundle \(\La{z}_\FRM\) from \(\FRC_\FRM\) to
\(\Stack{\FRC_\FRM/Z_\FRM}\) is still \(\La{z}_\FRM\) (again, understood as
pullback from \(X\)), and the cotangent
complex \(L_{\Stack{\FRC_\FRM/Z_\FRM}}\) of \(\Stack{\FRC_\FRM/Z_\FRM}\) is given by
\begin{align}
    \Stack{\FRC_\FRM^*/Z_\FRM}\longto \La{z}_\FRM^*,
\end{align}
placed at degrees \(0\) and \(1\) respectively. Let \(a\in\cA_X(\bar{k})\) be a
point in mH-base, and \(b\) (resp.~\(\cL\)) its image in
\(\cB_X\) (resp.~\(\Bun_{Z_\FRM}\)). Then
\(\LDF{a}^*L_{\Stack{\FRC_\FRM/Z_\FRM}}\) is
\begin{align}
    \FRC_{\FRM,\cL}^*\longto \La{z}_\FRM^*,
\end{align}
where the arrow is induced by the derivative of \(Z_\FRM\)-action at the
section of \(\FRC_{\FRM,\cL}\) given by \(a\). Taking \(\cO_{\breve{X}}\)-dual, we have
\begin{align}
    \cT^\bullet=\cT^{\le 0}=\La{z}_\FRM\longto \FRC_{\FRM,\cL},
\end{align}
where the map is the tangent action at \(a\).

\subsection{}
Write \(a=(b,c)\) where \(c\in\RH^0(\breve{X},\FRC_\cL)\). If we replace the stack
\(\Stack{\FRC_\FRM/Z_\FRM}\) with the morphism
\begin{align}
    \Stack{\FRC_\FRM/Z_\FRM}\to\Stack{\FRA_\FRM/Z_\FRM},
\end{align}
then we are looking at
deformation of \(a\) relative to \(b\in\cB_X\). Denote the analogue of
\(\cT^\bullet\) in this relative setting by \(\cT_b^\bullet\), then
\(\cT_b^\bullet\) is isomorphic to \(\FRC_\cL\) viewed as a locally free
sheaf. In this case the obstruction space is given by \(\RH^1(\breve{X},\FRC_\cL)\), which
vanishes if \(\cL\) is very \(G\)-ample as we have seen in \Cref{sec:dimensions}.

We have a short exact sequence of coherent sheaves
\begin{align}
    0\longto \FRC_\cL\longto
    \coker{\cT^\bullet}=\FRC_{\FRM,\cL}/\La{z}_\FRM\longto
    \FRA_{\FRM,\cL}/\La{z}_\FRM\longto 0.
\end{align}
The quotient \(\FRA_{\FRM,\cL}/\La{z}_\FRM\) is torsion and supported on the
boundary divisor \(b\), and its length is \(\dim_b\cB_X\).
It is reasonable to expect that the obstruction space
\(\RH^1(\breve{X},\FRC_{\FRM,\cL}/\La{z}_\FRM)\) is more likely to vanish than
\(\RH^1(\breve{X},\FRC_\cL)\) if \(\dim_b\cB_X\) is large. Of course, it 
depends on what \(b\) and \(c\) are, and it could well be the extreme case where
the above exact sequence splits and
\(\RH^1(\breve{X},\FRC_{\FRM,\cL}/\La{z}_\FRM)\) is equal to
\(\RH^1(\breve{X},\FRC_\cL)\).

It seems difficult to analyze the quotient
\(\FRC_{\FRM,\cL}/\La{z}_\FRM\) in a direct way. To solve this problem, we
need to use the fact that \(\FRC_\FRM\) is the invariant quotient of \(\FRT_\FRM\)
and utilize the cameral cover.

\subsection{}
Recall the big-cell locus \(\FRM^\circ\) and its intersection
\(\FRT_\FRM^\circ\) with \(\FRT_\FRM\). Let \(\iota\colon\FRT_\FRM^\circ\to
\FRT_\FRM\) be the inclusion.
Since \(\FRT_\FRM^\circ\) is smooth over \(X\), its cotangent bundle (relative
to \(X\)) is locally free, and fits into the canonical short exact sequence
\begin{align}
    0\longto \alpha_\FRM^*\FRA_\FRM^*\longto
    \CoTB_{\FRT_\FRM^\circ/X}\longto (\La{t}^\SC)^*\longto 0,
\end{align}
where \(\alpha_\FRM\) is the abelianization map of \(\FRM\).
Since \(\FRA_\FRM\) is of standard type, the numerical boundary divisor
\(\FRE_\FRM\) is a union of smooth Cartier divisors with normal crossings on \(\FRA_\FRM\).
Therefore, we may consider the logarithmic cotangent bundles and exact
sequence
\begin{align}
    0\longto \alpha_\FRM^*\FRA_\FRM^*(\log{\FRE_\FRM})\simeq\La{z}_\FRM^*\longto
    \CoTB_{\FRT_\FRM^\circ/X}(\log{\alpha_\FRM^*\FRE_\FRM})\simeq\La{t}_\FRM^*\longto
    (\La{t}^\SC)^*\longto 0,
\end{align}
and it splits canonically. It also realizes \(\iota_*\CoTB_{\FRT_\FRM^\circ/X}\)
as a subsheaf of the constant sheaf \(\La{t}_\FRM^*\) because \(\iota_*\) is
left-exact.

\subsection{}
Let \(\pi_a\colon \tilde{X}_a\to \breve{X}\) be the cameral cover and
\(\tilde{X}_a^\flat\) the normalization of \(\tilde{X}_a\). Let
\(\pi_a^\flat\colon\tilde{X}_a^\flat\to\breve{X}\) be the \(W\)-cover
and \(\tilde{a}^\flat\colon\tilde{X}_a^\flat\to\FRT_{\FRM,\cL}\) the natural
map induced by \(a\).
The differential action
\begin{align}
    (\tilde{a}^\flat)^*\CoTB_{\FRT_{\FRM,\cL}/X}\longto \La{z}_\FRM^*
\end{align}
factors through \(\La{t}_\FRM^*\),
and we have maps of coherent sheaves on \(\tilde{X}_a^\flat\)
\begin{align}
    \label{eqn:deformation_mH_base_pre_important_string}
    (\tilde{a}^\flat)^*\CoTB_{\FRT_{\FRM,\cL}/X}\longto
    (\tilde{a}^\flat)^*\iota_*\CoTB_{\FRT_{\FRM,\cL}^\circ/X}\longto
    \La{t}_\FRM^*\longto
    \La{z}_\FRM^*.
\end{align}
Pushing forward to \(\breve{X}\) and taking \(W\)-invariants, we get
\begin{align}
    \label{eqn:deformation_mH_base_important_string}
    (\pi_{a*}^\flat(\tilde{a}^\flat)^*\CoTB_{\FRT_{\FRM,\cL}/X})^W\longto
    (\pi_{a*}^\flat(\tilde{a}^\flat)^*\iota_*\CoTB_{\FRT_{\FRM,\cL}^\circ/X})^W\longto
    (\pi_{a*}^\flat\La{t}_\FRM^*)^W\longto
    \La{z}_\FRM^*.
\end{align}

\subsection{}
The usual cotangent sequence of map \(\FRT_\FRM\to\FRC_\FRM\) induces a sequence
on \(\breve{X}\) that is exact in the middle
\begin{align}
    \pi_{a*}^\flat(\pi_a^\flat)^*a^*\CoTB_{\FRC_{\FRM,\cL}/X}\longto
    \pi_{a*}^\flat(\tilde{a}^\flat)^*\CoTB_{\FRT_{\FRM,\cL}/X}\longto
    \pi_{a*}^\flat(\tilde{a}^\flat)^*\CoTB_{\FRT_{\FRM,\cL}/\FRC_{\FRM,\cL}}.
\end{align}
Note that the first sheaf above is locally free since \(\FRA_\FRM\) is of
standard type.
Assume \(a\in\cA_X^\heartsuit\), then the sheaf
\(\pi_{a*}^\flat(\tilde{a}^\flat)^*\CoTB_{\FRT_{\FRM,\cL}/\FRC_{\FRM,\cL}}\) is torsion, and
since \(\breve{X}\) is a smooth curve, we have injective map
\begin{align}
    \pi_{a*}^\flat(\pi_a^\flat)^*a^*\CoTB_{\FRC_{\FRM,\cL}/X}\longto
    (\pi_{a*}^\flat(\tilde{a}^\flat)^*\CoTB_{\FRT_{\FRM,\cL}/X})^\tfree,
\end{align}
where superscript \(\tfree\)
    \nomenclature[\(.tf \)]{\((\cdot)^\tfree\)}{the largest torsion-free quotient module/bundle}
means taking torsion-free quotient.
Since \(\pi_a^\flat\) is a finite flat \(W\)-cover,
we have
\begin{align}
    (\pi_{a*}^\flat(\pi_a^\flat)^*a^*\CoTB_{\FRC_{\FRM,\cL}/X})^W\simeq
    a^*\CoTB_{\FRC_{\FRM,\cL}/X}\simeq \FRC_{\FRM,\cL}^*.
\end{align}
This means that
\((\pi_{a*}^\flat(\tilde{a}^\flat)^*\CoTB_{\FRT_{\FRM,\cL}/X})^{W,\tfree}\) contains
\(\FRC_{\FRM,\cL}^*\) as a subsheaf with torsion quotient.
Note that for any coherent sheaf \(\cF\) on \(\breve{X}\) with a
\(\cO_{\breve{X}}\)-linear \(W\)-action, taking \(W\)-invariants commutes with
taking torsion-free quotient because \(\breve{X}\) is a smooth curve
and the order of \(W\) is invertible in \(k\).
Similarly, since the first map in
\eqref{eqn:deformation_mH_base_pre_important_string} is generically an
isomorphism, the sheaf
\((\pi_{a*}^\flat(\tilde{a}^\flat)^*\iota_*\CoTB_{\FRT_{\FRM,\cL}^\circ/X})^{W,\tfree}\)
also contains \(\FRC_{\FRM,\cL}^*\) with torsion quotient.

\subsection{}
By Serre duality, we have
\begin{align}
    \RH^1\bigl(\breve{X},\coker(\La{z}_\FRM\to
    \FRC_{\FRM,\cL})\bigr)
    \simeq\RH^0\bigl(\breve{X},\ker(\FRC_{\FRM,\cL}^*\to\La{z}_\FRM^*)\otimes_{\cO_{\breve{X}}}\CoTB_{\breve{X}}\bigr).
\end{align}
Thus, to show that these spaces vanish, it is sufficient by
our discussion above to show
\begin{align}
    \label{eqn:deformation_mH_base_simplified_goal}
    \RH^0\Bigl(\breve{X},\pi_{a*}^\flat\ker[(\tilde{a}^\flat)^*\iota_*\CoTB_{\FRT_{\FRM,\cL}^\circ/X}\to\La{z}_\FRM^*]^{W,\tfree}
    \otimes_{\cO_{\breve{X}}}\CoTB_{\breve{X}}\Bigr)=0.
\end{align}
For convenience, we denote
\(\cK=\ker[(\tilde{a}^\flat)^*\iota_*\CoTB_{\FRT_{\FRM,\cL}^\circ/X}\to\La{z}_\FRM^*]\),
then \eqref{eqn:deformation_mH_base_simplified_goal} is equivalent to
that if \(L\) is a locally free subsheaf of
\(\pi_{a*}^\flat\cK^{W,\tfree}\) of rank \(1\), then
\begin{align}
    \Hom_{\breve{X}}(\CoTB_{\breve{X}}^{-1},L)=0.
\end{align}
Given any locally free sheaf \(L\) of finite rank on \(\breve{X}\), we have
\(L\simeq (\pi_{a*}^\flat(\pi_a^\flat)^*L)^W\). So by adjunction
\begin{align}
    \Hom_{\breve{X}}(\CoTB_{\breve{X}}^{-1},L)
    &=\Hom_{\breve{X}}(\CoTB_{\breve{X}}^{-1},(\pi_{a*}^\flat(\pi_a^\flat)^*L)^W)\\
    &=\Hom_{\breve{X}}(\CoTB_{\breve{X}}^{-1},\pi_{a*}^\flat(\pi_a^\flat)^*L)^W\\
    &=\Hom_{\tilde{X}_a^\flat}^W((\pi_a^\flat)^*\CoTB_{\breve{X}}^{-1},(\pi_a^\flat)^*L).
\end{align}
Therefore, it suffices to prove the following statement: for any
\(W\)-equivariant locally free subsheaf \(L\) of \(\cK^\tfree\) of rank \(1\) we have
\begin{align}
    \label{eqn:deformation_mH_base_more_simplified_goal}
    \Hom_{\tilde{X}_a^\flat}^W((\pi_a^\flat)^*\CoTB_{\breve{X}}^{-1},L)=0.
\end{align}
We will prove \eqref{eqn:deformation_mH_base_more_simplified_goal} under some
mild but technical conditions. Before that, we need more preparations.

\subsection{}
First, we will do some standard reduction.
Let \(\OGT_G\colon X_\OGT\to X\) be a finite \'etale Galois cover over which
the monoid \(\FRM\) becomes split and let \(\Theta\) be the Galois group.
To prove \eqref{eqn:deformation_mH_base_more_simplified_goal} it suffices to
prove its analogue after base change to \(X_\OGT\). Therefore, by
replacing \(X\) by \(X_\OGT\), we may assume that \(\FRM\) is split and \(W\) is
constant.
Let \(\theta_1,\ldots,\theta_m\) be the free generators of the cone of
\(\FRA_\FRM\), in other words, \(\FRM=\bM(\theta_1,\ldots,\theta_m)\). We let
\(\beta_1,\cdots,\beta_m\) be the corresponding coordinates on \(\FRA_\FRM\),
also understood as characters of \(\FRA_\FRM^\x\). The
canonical map \(\FRA_\FRM\to \FRA_{\Env(G^\SC)}\) gives map on the coordinate
rings
\begin{align}
    \Rt_s\longmapsto \sum_{t=1}^m\beta_t^{c_{st}},
\end{align}
where \(c_{st}\in \bbN\).

\subsection{}
Next, we want to describe the sheaf \(\cK^\tfree\) in a more
explicit way. More precisely, we are interested in describing \(\cK^\tfree\) as
a subsheaf of constant sheaf \((\La{t}^\SC)^*=\ker[\La{t}_\FRM^*\to\La{z}_\FRM^*]\).

Let \(b=(\cL,\lambda_b)\in\cB_X(\bar{k})\) be the image of \(a\) with
\begin{align}
    \lambda_b=\sum_{i=1}^d\lambda_i\cdot v_i,
\end{align}
and
\begin{align}
    \lambda_i=\sum_{j=1}^m l_{ij}\theta_j.
\end{align}
Let \(D_i=(\pi_a^\flat)^{-1}(v_i)\), then it is an effective Cartier divisor
of degree \(\abs{\bW}\) on \(\tilde{X}_a^\flat\). Suppose
\(D_i=\sum_{j=1}^{e_i}d_{ij}\tilde{v}_{ij}\) where
\(\tilde{v}_{ij}\in\tilde{X}_a^\flat(\bar{k})\), then the map
\(\tilde{X}_a^\flat\to\Stack{\FRA_\FRM/Z_\FRM}\) induced by \(b\) may be
presented by
\begin{align}
    \left((\pi_a^\flat)^*\cL,
        \sum_{i=1}^d\sum_{j=1}^{e_i}\tilde{\lambda}_{ij}\cdot
    \tilde{v}_{ij}\right),
\end{align}
where
\begin{align}
    \tilde{\lambda}_{ij}=d_{ij}\lambda_i=\sum_{t=1}^m
    d_{ij}l_{it}\theta_t=\sum_{t=1}^m \tilde{d}_{ijt}\theta_t.
\end{align}

Let \(Y_{ij}=\Spec{\cO_{ij}}\) be the formal disc around \(\tilde{v}_{ij}\)
and \(\pi_{ij}\) be a fixed uniformizer of \(\cO_{ij}\). Then after
choosing a trivialization of \(\cL\) on \(Y_{ij}\) we may
identify \(\tilde{a}^\flat(Y_{ij})\) with a point in
\(\pi_{ij}^{(\tilde{\lambda}_{ij},\mu_{ij})}T_\FRM(\cO_{ij})\) such that
\(w(\mu_{ij})\le -w_0(\tilde{\lambda}_{ij,\AD})\) for any \(w\in W\). The
cocharacter \(\mu_{ij}\) does not depend on the choice of \(\pi_{ij}\) nor the
trivialization of \(\cL\) over \(Y_{ij}\).

The sheaf \((\La{t}^\SC)^*/\cK^\tfree\) is
supported on \(\tilde{v}_{ij}\), so we may restrict to each \(Y_{ij}\),
and we will only be focusing on those
such that \(\mu_{ij}\in -W(\tilde{\lambda}_{ij,\AD})\), in other words,
\(\tilde{v}_{ij}\) is mapped into the big-cell locus. It is admittedly less
ideal but good enough for our purposes.
Let \(\mu_{ij}=-w(\tilde{\lambda}_{ij,\AD})\) for some \(w\in W\),
and we know as a coherent sheaf on \(\tilde{X}_a^\flat\), \(\La{t}_\FRM^*\) may
be decomposed as
\begin{align}
    \La{t}_\FRM^*
    =\left[\bigoplus_{t=1}^m\cO_{\tilde{X}_a^\flat}\frac{\dd(\beta_t,0)}{(\beta_t,0)}\right]
    \oplus
    \left[\bigoplus_{s=1}^r \cO_{\tilde{X}_a^\flat}\frac{\dd
    (\Rt_s,w(\Rt_s))}{(\Rt_s,w(\Rt_s))}\right].
\end{align}
For convenience, we denote logarithmic differentials
\begin{align}
    \bFf_t=\frac{\dd(\beta_t,0)}{(\beta_t,0)},\quad
    \bFe_{w,s}=\frac{\dd (\Rt_s,w(\Rt_s))}{(\Rt_s,w(\Rt_s))}.
\end{align}
The kernel of \(\La{t}_\FRM^*\to\La{z}_\FRM^*\) is identified with
\((\La{t}^\SC)^*\) through projection
\begin{align}
    \bFe_{w,s}\longmapsto \bar{\bFe}_{w,s}\defeq\frac{\dd w(\Rt_s)}{w(\Rt_s)}.
\end{align}
The characters \((\Rt_s,w(\Rt_s))\) are invertible over
\(Y_{ij}\), so it is straightforward to see that over \(Y_{ij}\),
\begin{align}
    (\tilde{a}^\flat)^*\iota_*\CoTB_{\FRT_{\FRM,\cL}^\circ/X}
    =\left[\bigoplus_{t=1}^m\pi_{ij}^{\tilde{d}_{ijt}}\cO_{ij}\bFf_t\right]
    \oplus\left[\bigoplus_{\Rt\in\SimRts} \cO_{ij}\bFe_{w,s}\right]\subset \La{t}_\FRM^*.
\end{align}
For each \(\Rt_s\in\SimRts\), the map
\((\tilde{a}^\flat)^*\iota_*\CoTB_{\FRT_{\FRM,\cL}^\circ/X}\to\La{z}_\FRM^*\)
maps basis vector \(\bFe_{w,s}\) to the element
\begin{align}
    \sum_{t=1}^m c_{st}\frac{\dd\beta_t}{\beta_t},
\end{align}
and an element \(\sum_{s=1}^rx_s\bFe_{w,s}+\sum_{t=1}^m y_t\bFf_t\) (\(x_s,y_t\in\cO_{ij}\)) is
contained in \(\cK^\tfree\) if and only if \(y_t=-\sum_{s=1}^r c_{st}x_s\) and
\begin{align}
    \val_{\cO_{ij}}y_t=\val_{\cO_{ij}}\left(\sum_{s=1}^r c_{st}x_s\right)\ge
    \tilde{d}_{ijt}.
\end{align}
The image of such element in \((\La{t}^\SC)^*\) is
\(\sum_sx_s\bar{\bFe}_{w,s}\). Thus, \(\cK^\tfree\subset (\La{t}^\SC)^*\) is
spanned by elements \(\sum_sx_s\bar{\bFe}_{w,s}\) such that
\begin{align}
    \label{eqn:deformation_mH_base_cK_tfree_first_desc}
    \val_{\cO_{ij}}\left(\sum_{s=1}^rc_{st}x_s\right)\ge\tilde{d}_{ijt}.
\end{align}
Note that this condition does not depend on the choice of trivialization of
\(\cL\) over \(Y_{ij}\), but is only valid when
\(\mu_{ij}=-w(\tilde{\lambda}_{ij})\) for this specific \(w\).

\subsection{}
At this stage the condition
\eqref{eqn:deformation_mH_base_cK_tfree_first_desc} is not very convenient
to use (the matrix \((c_{ij})\) is not even a square matrix).
There is also a glaring issue in positive characteristic: if \(p=\Char(k)>0\) divides
every \(c_{ij}\), then
\eqref{eqn:deformation_mH_base_cK_tfree_first_desc} is always satisfied, and so
over \(Y_{ij}\) we actually have \(\cK^\tfree=(\La{t}^\SC)^*\), which defeats
the purpose of analyzing \(\cK^\tfree\). Therefore, we will need an additional
assumption on \(a\) which is best stated using the factorization results in
\Cref{sec:factorizations}.

We use the notations in \Cref{sec:factorizations} and let \(I=\Set{1,2}\). Let
\(h_X^I\colon\cM_X^I\to \cA_X^I\) be the mH-fibration associated with monoid
\(\FRM^I\). Suppose \(a_I\in\cA_X^I\) is a point mapping to \(a\) via the map
\(\cA_X^I\to\cA_X\) induced by the summation map \(\Sigma\colon\BD[I]_X\to
\BD_X\). Let \(b_I\) be the image of \(a_I\) in \(\BD[I]_X\) and 
\(b_i\) (\(i=1,2\)) the images of \(b_I\) in \(\BD_X\) under the \(i\)-th
projection map \(\pr_i\colon\BD[I]_X\to\BD_X\). Note that for any
\(a_I\in\Sigma^{-1}(a)\) we naturally have
equality of Cartier divisors \(\FRD_a=\FRD_{a_I}\), as well as the cameral cover
\(\tilde{X}_a=\tilde{X}_{a_I}\) and their normalizations.  We
assume that there exists some \(a_I\) mapping to \(a\) such that
\begin{enumerate}
    \item \(b_I\in\BD[I,\disj]_X\), in other words,
        \(\supp(b_1)\cap\supp(b_2)=\emptyset\),
    \item \(\supp(b_2)\cap\supp(\FRD_a)=\emptyset\), and
    \item \(b_2\) is good (cf.~\Cref{def:goodness_of_boundary_divisors}).
\end{enumerate}
Let \(\tilde{b}_i\in\BD_{\tilde{X}_a^\flat}\) be the pullback of \(b_i\) to
\(\tilde{X}_a^\flat\), then we have
\begin{align}
    \dim_{\tilde{b}_2}\BD_{\tilde{X}_a^\flat}
    =\sum_{t=1}^m\sum_{\tilde{v}_{ij}\in\supp(\tilde{b}_2)}\tilde{d}_{ijt}
    =\sum_{t=1}^m\Bigl(\abs{\bW}\sum_{v_i\in\supp(b_2)}d_{it}\Bigr)
    =\abs{\bW}\dim_{b_2}\BD_X,
\end{align}
and similarly for \(b_1\). For \(b_2\) specifically, by the goodness
assumption (and that \(\abs{\bW}\) is invertible in \(k\)), we have
\(\tilde{d}_{ijt}\neq 0\) only if \(\theta_t\) is good. For convenience, we will
also define for each \(t\)
\begin{align}
    \dim_{\tilde{b}_2}^t\BD_{\tilde{X}_a^\flat}
    \defeq\sum_{\tilde{v}_{ij}\in\supp(\tilde{b}_2)}\tilde{d}_{ijt}
    =\abs{\bW}\dim_{b_2}^t\BD_X,
\end{align}
and it is non-zero only if \(\theta_t\) is good. We have elementary inequality
\begin{align}
    \max_{1\le t\le
    m}\dim_{b_2}^t\BD_X\ge\frac{\dim_{b_2}\BD_X}{m}.
\end{align}
Furthermore, if \(\tilde{v}_{ij}\in\supp(\tilde{b}_2)\), then by our assumptions
we also have \(\mu_{ij}\in -W(\tilde{\lambda}_{ij,\AD})\). 

Since \(b_I\) is in the disjoint locus, the summation map \(\Sigma\) is an
isomorphism near \(a_I\), therefore everything about deformations so far carries
over identically to \(a_I\) by pulling back.
Without loss of generality, suppose \(\theta_1\) is good, \(\dim_{b_2}^t\BD_X\)
achieves maximum at \(t=1\), and \(c_{11}\) is invertible in \(k\). Then for any
fixed \(w\in W\) we may
find another basis of \(\La{t}_\FRM^*\) consisting of all \(\bFf_t\) and
\begin{align}
    \bFe_{w,s}'=\begin{cases}
        \bFe_{w,1} & s=1,\\
        \bFe_{w,s}-c_{11}^{-1}c_{s1}\bFe_{w,1} & s\neq 1.
    \end{cases}
\end{align}
If \(\sum_sx_s\bFe_{w,s}'+\sum_ty_t\bFf_t\) is contained in \(\cK^\tfree\), then it
is necessary that
\begin{align}
    \label{eqn:deformation_mH_base_cK_tfree_convenient_desc}
    \val_{\cO_{ij}}(c_{11}x_1)=\val_{\cO_{ij}}(x_1)\ge \tilde{d}_{ij1}.
\end{align}

\subsection{}
Now let \(L\) be a \(W\)-equivariant subsheaf of \(\cK^\tfree\) of rank \(1\).
Any root \(\Rt\in\Roots\) determines a constant subbundle
\begin{align}
    L_\Rt\defeq\cO_{\tilde{X}_a^\flat}\frac{\dd \Rt}{\Rt}\subset (\La{t}^\SC)^*,
\end{align}
and when \(\Rt=w(\Rt_1)\), the basis \(\bFe_{w,s}'\) induces a decomposition
\begin{align}
    (\La{t}^\SC)^*=L_{\Rt}\oplus
    \left[\bigoplus_{s=2}^r\cO_{\tilde{X}_a^\flat}\bar{\bFe}_{w,s}'\right].
\end{align}
Since any constant line in \((\La{t}^\SC)^*\) generate
the whole bundle under \(W\) and \(L\) is \(W\)-equivariant, the induced
projection map \(L\to L_{\Rt}\) must be non-zero, otherwise the inclusion map
\(L\to (\La{t}^\SC)^*\) would also be zero, which is a contradiction. Let
\(\cK_{\Rt}\) be the image of \(\cK^\tfree\) in \(L_{\Rt}\), then by
\eqref{eqn:deformation_mH_base_cK_tfree_convenient_desc}, we see that
\(\cK_{\Rt}\) is contained in the subsheaf
\begin{align}
    L_{\Rt}\Bigl(-\sum_{\tilde{v}_{ij}}\tilde{d}_{ij1}\tilde{v}_{ij}\Bigr)\subset
    L_{\Rt},
\end{align}
where \(\tilde{v}_{ij}\) ranges over those points in \(\supp(\tilde{b}_2)\) such
that \(\mu_{ij}=-w(\tilde{\lambda}_{ij,\AD})\) for a fixed \(w\). As a result, we
have an estimate of degree of \(L\):
\begin{align}
    \deg{L}\le\deg{\cK_{\Rt}}
    &\le -\sum_{\tilde{v}_{ij}}\tilde{d}_{ij1}\\
    &\le -\frac{1}{\abs{\bW}}\dim_{\tilde{b}_2}^1\BD_{\tilde{X}_a^\flat}\\
    &\le
    -\frac{1}{m\abs{\bW}}\dim_{\tilde{b}_2}\BD_{\tilde{X}_a^\flat}=-\frac{1}{m}\dim_{b_2}\BD_X,
\end{align}
We can now finally state the
precise result regarding \eqref{eqn:deformation_mH_base_more_simplified_goal}:

\begin{theorem}
    \label[theorem]{thm:deformation_mH_base_main}
    Suppose \(a\in\cA_X^\heartsuit(\bar{k})\) is such that there exists some
    \(a_I\in\Sigma^{-1}(a)\cap\cA_X^{I,\disj}\) such that \(b_2\) is good with
    \(\supp(\FRD_a)\cap \supp(b_2)=\emptyset\) and
    \begin{align}
        \dim_{b_2}\BD_X\ge 2m\abs{\bW}(g_X-1),
    \end{align}
    then \eqref{eqn:deformation_mH_base_more_simplified_goal} holds. As a
    consequence, \(\cA_X\) is smooth at \(a\).
\end{theorem}
\begin{proof}
    From our discussion so far, the line bundle \(L\) has smaller degree than
    \((\pi_a^\flat)^*\CoTB_{\breve{X}}^{-1}\). Since both are \(W\)-equivariant,
    the same is true over each irreducible component of \(\tilde{X}_a^\flat\),
    thus \eqref{eqn:deformation_mH_base_more_simplified_goal} holds.
\end{proof}


\section{Deformation of mHiggs bundles} 
\label{sec:deformation_of_mHiggs_bundles}

In this section, we study deformation of mHiggs bundles. As discussed in
\Cref{sec:local_model_of_singularities}, the main result of this section will be
a proof of \Cref{thm:local_singularity_model_main}. 
We will do this in two steps: the first step is a
proof of \Cref{thm:local_singularity_model_weak}, which uses certain
self-duality within monoid \(\FRM\) inspired by \cite{Br09}, and the second is
combining the first step and the results in \Cref{sec:the_case_of_mH_base} to
obtain a proof of \Cref{thm:local_singularity_model_main}.

For the second step,
on top of our usual assumption on \(\Char(k)\), we need to impose some
additional characteristic-related conditions
(see \Cref{thm:local_singularity_model_main}). It will not
affect our proof of the fundamental lemma because it can be proved using monoids
for which such conditions are void
(cf.~\Cref{sec:Reduction_to_minuscule_and_quasi_minuscule_coweights}).

\subsection{}
Let \(\FRM\in\FM(G^\SC)\) be a very flat monoid.
Recall we have the big-cell locus \(\FRM^\circ\subset
\FRM\) such that the restriction of abelianization map
\(\alpha_\FRM^\circ\colon\FRM^\circ\to\FRA_\FRM\) is smooth, and its fibers are
homogeneous spaces under \(G^\SC\x G^\SC\). The action of \(G^\SC\x G^\SC\)
induces injection map of vector bundles on \(\FRM^\circ\):
\begin{align}
    \CoTB_{\FRM^\circ/\FRA_\FRM}\longto (\La{g}^{\SC})^*\x (\La{g}^{\SC})^*,
\end{align}
and its (surjective) dual map
\begin{align}
    \La{g}^{\SC}\x \La{g}^{\SC}\longto\TanB_{\FRM^\circ/\FRA_\FRM}.
\end{align}
Choosing a non-degenerate \(G^\SC\)-invariant symmetric bilinear form \(Q\) on
\(\La{g}^\SC\), we may identify \(\La{g}^\SC\) with its dual. We fix the
\emph{anti-diagonal} form \((Q,-Q)\) on \(\La{g}^\SC\x\La{g}^\SC\) and use this
form to identify \(\La{g}^\SC\x\La{g}^\SC\) with its dual. Note that under such
identification, the diagonal subspace
\begin{align}
    \Delta\colon\La{g}^\SC\longto \La{g}^\SC\x\La{g}^\SC
\end{align}
and the  quotient
\(\La{g}^\SC\x\La{g}^\SC/\Delta(\La{g}^\SC)\) are \(G^\SC\)-equivariant duals.
In other words, we have a short exact sequence that is \(G^\SC\)-equivariantly
self-dual:
\begin{align}
    0\longto \La{g}^\SC\stackrel{\Delta}{\longto} \La{g}^\SC\x\La{g}^\SC\longto
    \La{g}^\SC\x\La{g}^\SC/\Delta(\La{g}^\SC)\longto 0.
\end{align}

The following result is inspired by \cite{Br09}*{Example~2.5}:
\begin{lemma}
    With the choice of \((Q,-Q)\), we have a \(G^\SC\)-equivariantly
    self-dual short exact sequence
    \begin{align}
        0\longto \CoTB_{\FRM^\circ/\FRA_\FRM}\longto \La{g}^{\SC}\x
        \La{g}^{\SC}\longto \TanB_{\FRM^\circ/\FRA_\FRM}\longto 0.
    \end{align}
\end{lemma}
\begin{proof}
    It suffices to consider \(\FRM=\Env(G^\SC)\) by universal property and show
    that the fibers of
    \begin{align}
        \ker\bigl(\La{g}^{\SC}\x \La{g}^{\SC}\longto
            \TanB_{\FRM^\circ/\FRA_\FRM}\bigr)
    \end{align}
    are maximal \((Q,-Q)\)-isotropic subspaces of \(\La{g}^{\SC}\x
    \La{g}^{\SC}\). Since \(Z_\FRM\) commutes with \(G^\SC\x G^\SC\)-action, it
    suffices to prove the statement for \(ge_{I,\SimRts}h\), where \(g,h\in
    G^\SC\), and \(e_{I,\SimRts}\) are the system of idempotents associated with
    subsets of simple roots \(I\subset \SimRts\)
    (cf.~\Cref{sub:desc_of_big_cell_idempotents}). Moreover, since \((Q,-Q)\)
    is \(G^\SC\x G^\SC\)-invariant, it suffices to prove the statement for
    \(e_{I,\SimRts}\).

    Indeed, Let \(P_I\) (resp.~\(P_I^-\)) be the standard parabolic subgroup of
    \(G^\SC\) containing \(B\) (resp.~\(B^-\)), \(U_I\) (resp.~\(U_I^-\)) be its
    unipotent radical, and \(L_I\) be the Levi factor containing \(T^\SC\), then
    the stabilizer of \(e_{I,\SimRts}\) in \(G^\SC\x G^\SC\) is the semidirect
    product \((U_I\x U_I^-)\rtimes\Delta(L_I)\), where \(\Delta(L_I)\) is the
    diagonal embedding of \(L_I\) in \(G^\SC\x G^\SC\). One then verifies by
    direct computation that \(\La{u}_I\oplus\La{u}_I^-\oplus\Delta(\La{l}_I)\)
    is \((Q,-Q)\)-isotropic and has half the dimension of
    \(\La{g}^\SC\x\La{g}^\SC\).
\end{proof}

\subsection{}
The monoid \(\FRM\) is usually not smooth, but since it is normal and
\(\FRM-\FRM^\circ\) has codimension at least \(2\), we still have identification
by Hartogs's theorem (here \(\iota\) is the inclusion map \(\FRM^\circ\to \FRM\)):
\begin{align}
    \iota_*(\La{g}^\SC\x\La{g}^\SC)=\La{g}^\SC\x \La{g}^\SC,
\end{align}
because locally over the curve \(X\), \(\La{g}^\SC\x \La{g}^\SC\) is a trivial
vector bundle. By functoriality, we have maps
\begin{align}
    \CoTB_{\FRM/\FRA_\FRM}\longto \iota_*\CoTB_{\FRM^\circ/\FRA_\FRM}\longto
    \La{g}^\SC\x \La{g}^\SC\longto \TanB_{\FRM/\FRA_\FRM},
\end{align}
where we already identified \(\La{g}^\SC\x \La{g}^\SC\) with its dual using
\((Q,-Q)\). Since \(\FRM\) is normal, and \(\TanB_{\FRM/\FRA_\FRM}\) is
reflexive, it may be identified with \(\iota_*\TanB_{\FRM^\circ/\FRA_\FRM}\).
Now since \(Z_\FRM\) commutes with \(G^\SC\x G^\SC\), everything descends to
quotient
\begin{align}
    \Stack{\FRM/Z_\FRM}\longto\Stack{\FRA_\FRM/Z_\FRM}.
\end{align}
To simplify notations, we denote \(\Stack{\FRM}=\Stack{\FRM/Z_\FRM}\) and
\(\Stack{\FRA_\FRM}=\Stack{\FRA_\FRM/Z_\FRM}\). So we have
\begin{align}
    \CoTB_{\Stack{\FRM}/\Stack{\FRA_\FRM}}\longto
    \iota_*\CoTB_{\Stack{\FRM^\circ}/\Stack{\FRA_\FRM}}\longto
    \La{g}^\SC\x \La{g}^\SC\longto \TanB_{\Stack{\FRM}/\Stack{\FRA_\FRM}}.
\end{align}
Let \(x=(\cL,E,\phi)\in \cM_X(\bar{k})^\heartsuit\) be an mHiggs bundle, viewed
as a \(\breve{X}\)-point of the stack \(\Stack{\FRM/G\x Z_\FRM}\), and let
\(b\in \cB_X(\bar{k})\subset\Stack{\FRA_\FRM}(\breve{X})\) be its boundary
divisor. Then
\eqref{eqn:needed_vanishing_statement_for_deformation} translates to the
following statement:
\begin{align}
    \RH^1\Bigl(\breve{X},\coker\Bigl[\ad(E)\stackrel{\DD_{\Ad}}{\longto}
    \bigl(\phi^*\CoTB_{E\x^G
\FRM_\cL/\FRA_{\FRM,\cL}}\bigr)^*\Bigr]\Bigr)=0.
\end{align}
Since the
image of \(\DD_{\Ad}\) stays the same if we replace \(\ad(E)\) by
\(\ad(E)^\SC\), it is reduced to proving that
\begin{align}
    \RH^1\Bigl(\breve{X},\coker\Bigl[\ad(E)^\SC\stackrel{\DD_{\Ad}}{\longto}
    \bigl(\phi^*\CoTB_{E\x^{G} \FRM_\cL/\FRA_{\FRM,\cL}}\bigr)^*\Bigr]\Bigr)=0.
\end{align}

\subsection{}
On the other hand, since \(\iota_*\) is left-exact, we have left-exact sequence
\begin{align}
    0\longto \iota_*\CoTB_{\FRM^\circ/\FRA_\FRM}\longto
    \La{g}^\SC\x \La{g}^\SC\longto
    \iota_*\TanB_{\FRM^\circ/\FRA_\FRM}\simeq\TanB_{\FRM/\FRA_\FRM},
\end{align}
which induces exact sequence
\begin{align}
    \phi^*\iota_*\CoTB_{E\x^G\FRM_\cL^\circ/\FRA_{\FRM,\cL}}\longto 
    \ad(E)^\SC\x \ad(E)^\SC\longto \phi^*\TanB_{E\x^G\FRM_\cL/\FRA_{\FRM,\cL}}.
\end{align}
Taking \(\cO_{\breve{X}}\)-dual and using bilinear form \((Q,-Q)\) on the middle
term, we have exact sequence
\begin{align}
    \bigl(\phi^*\TanB_{E\x^G\FRM_\cL/\FRA_{\FRM,\cL}}\bigr)^*\longto 
    \ad(E)^\SC\x \ad(E)^\SC\longto 
    \bigl(\phi^*\iota_*\CoTB_{E\x^G\FRM_\cL^\circ/\FRA_{\FRM,\cL}}\bigr)^*
\end{align}
Since \(\breve{X}\) is a smooth curve, both
\(\bigl(\phi^*\TanB_{E\x^G\FRM_\cL/\FRA_{\FRM,\cL}}\bigr)^*\) and
\(\bigl(\phi^*\iota_*\CoTB_{E\x^G\FRM_\cL^\circ/\FRA_{\FRM,\cL}}\bigr)^*\)
are locally free because they are torsion-free, so the above sequence is also
exact on the left. Let
\begin{align}
    K\subset
    \bigl(\phi^*\iota_*\CoTB_{E\x^G\FRM_\cL^\circ/\FRA_{\FRM,\cL}}\bigr)^*
\end{align}
be the cokernel of the first map, then it is again locally free. This implies
that for any \(v\in\breve{X}\), the fiber map
\begin{align}
    \bigl(\phi^*\TanB_{E\x^G\FRM_\cL/\FRA_{\FRM,\cL}}\bigr)_v^*\longto 
    (\ad(E)^\SC\x \ad(E)^\SC)_v
\end{align}
is injective. As a consequence, since generically over \(\breve{X}\) the fiber
of \(\bigl(\phi^*\TanB_{E\x^G\FRM_\cL/\FRA_{\FRM,\cL}}\bigr)\) is a
maximal \((Q,-Q)\)-isotropic subspace of \(\ad(E)^\SC\x \ad(E)^\SC\), the same
must be true over all \(\breve{X}\). Thus, it shows that the short exact sequence
\begin{align}
    0\longto 
    \bigl(\phi^*\TanB_{E\x^G\FRM_\cL/\FRA_{\FRM,\cL}}\bigr)^*\longto 
    \ad(E)^\SC\x \ad(E)^\SC\longto K\longto 0
\end{align}
is self-dual. Moreover, since the adjoint action of \(G^\SC\) is simply the
restriction of the \(G^\SC\x G^\SC\)-action to the diagonal, we see that its
derivative \(\DD_{\Ad}\) factors through \(K\).
Since we have inclusions of locally free sheaves on \(\breve{X}\)
\begin{align}
    K\subset \bigl(\phi^*\iota_*\CoTB_{E\x^G\FRM_\cL^\circ/\FRA_{\FRM,\cL}}\bigr)^*
    \subset
    \bigl(\phi^*\CoTB_{E\x^G\FRM_\cL/\FRA_{\FRM,\cL}}\bigr)^*,
\end{align}
to prove \eqref{eqn:needed_vanishing_statement_for_deformation} it suffices to
prove that
\begin{align}
    \RH^1\Bigl(\breve{X},\coker\bigl(\ad(E)^\SC\to K\bigr)\Bigr)=0.
\end{align}
Now consider the following self-dual diagram of locally free sheaves on
\(\breve{X}\):
\begin{equation}
    \begin{tikzcd}
         & & \ad(E)^\SC \ar[d, "\Delta"]\ar[rd, "\DD_{\Ad}'"] & &\\
        0\ar[r] & \bigl(\phi^*\TanB_{E\x^G\Stack{\FRM}/\Stack{\FRA_\FRM}}\bigr)^*\ar[r]\ar[rd] &
        \ad(E)^\SC\x\ad(E)^\SC\ar[r]\ar[d] & K\ar[r] & 0\\
         & & \ad(E)^\SC\x\ad(E)^\SC/\Delta(\ad(E)^\SC) & & 
    \end{tikzcd}
\end{equation}
The \(\cO_{\breve{X}}\)-dual of \(\coker(\DD_{\Ad}')\) is identified with the
kernel of the lower-left map. But that kernel is none other than the kernel
of \(\DD_{\Ad}'\), which is the same as the kernel of \(\DD_{\Ad}\) and equal to
\(\Lie(I_{(E,\phi)}^{\Sm})/\La{z}_G\)
(cf.~\Cref{sec:automorphism_group}). So we reduce
\eqref{eqn:needed_vanishing_statement_for_deformation} to
\begin{align}
    \label{eqn:I_sm_vanishing_statement_for_deformation}
    \RH^1(\breve{X},(\Lie(I_{(E,\phi)}^\Sm)/\La{z}_G)^*)=0.
\end{align}
The canonical inclusion map
\(I_{(E,\phi)}^\Sm\to \FRJ_a^\flat\) into the N\'eron model induces inclusion
\begin{align}
    \Lie(I_{(E,\phi)}^\Sm)/\La{z}_G\longto \Lie(\FRJ_a^\flat)/\La{z}_G,
\end{align}
hence we may also further reduce to
\begin{align}
    \RH^1(\breve{X},(\Lie(\FRJ_a^\flat)/\La{z}_G)^*)=0.
\end{align}

\subsection{}
Let \(\Ob(x)\) be the obstruction space of deforming \(x=(\cL,E,\phi)\)
relative to boundary divisor \(b\), and for any \(v\in\breve{X}\), let
\(\Ob_v(x)\) be the obstruction space of deforming the
\(\breve{X}_v\)-arc induced by \(x\). We have the following
global-local principle for the obstructions: 
\begin{proposition}
    \label[proposition]{prop:global_local_principle_obstructions}
    If \eqref{eqn:I_sm_vanishing_statement_for_deformation} holds, then we have
    canonical injective map
    \begin{align}
        \Ob(x)\longto\prod_v \Ob_v(x),
    \end{align}
    where \(v\) ranges over the points in \(\breve{X}\) that is sent to the
    singular locus of \(\Stack{\FRM/G\x Z_\FRM}\).
    In particular, it is true when \(x\) is very \((G,\delta_a)\)-ample.
\end{proposition}
\begin{proof}
    We only need to prove the second statement.
    By \Cref{prop:desc_of_Lie_FRJ_a}, we know that
    \begin{align}
        \Lie(\FRJ_a)/\La{z}_G\simeq \FRC_\cL^*.
    \end{align}
    The result then follows from the fact that \(\delta_a\) is the length of
    \(\Lie(\FRJ_a^\flat)/\Lie(\FRJ_a)\) and the definition of \(x\) being
    very \((G,\delta_a)\)-ample.
\end{proof}

\subsection{}
Next we turn to tangent spaces. As we have seen in
\Cref{sec:deformation_of_mapping_stacks}, the tangent space of
\(\cM_X\) at \((\cL,E,\phi)\) relative to \(\cB_X\) fits in the
short exact sequence
\begin{align}
    0\longto \RH^1(\breve{X},\ker(\DD_{\Ad}))\longto \TanB_{(\cL,E,\phi)}(\cM_X/\cB_X)\longto
    \RH^0(\breve{X},\coker(\DD_{\Ad}))\longto 0.
\end{align}
Restricting to a formal disc \(\breve{X}_v\to X\), we have commutative diagram
\begin{equation}
    \begin{tikzcd}
        0 \ar[r] & \RH^1(\breve{X},\ker(\DD_{\Ad}))\ar[r]\ar[d] &
        \TanB_{x}(\cM_X/\cB_X) \ar[r]\ar[d] &
        \RH^0(\breve{X},\coker(\DD_{\Ad})) \ar[r]\ar[d] & 0\\
          &       0\ar[r] &
        \TanB_{(E_v,\phi_v)}(\Arc_v\Stack{\FRM_\cL/G}/\Arc_v\FRA_{\FRM,\cL})
        \ar[r, "\sim"] &
        \RH^0(\breve{X}_v,\coker(\DD_{\Ad}))      \ar[r] & 0
    \end{tikzcd}
\end{equation}
Since \(\coker(\DD_{\Ad})\) is unaffected by the center of \(G\), we may replace
\(G\) by \(G^\AD\x G^\AD\), and since \(\La{g}^\AD\simeq\La{g}^\SC\), we have
natural map
\begin{align}
    \coker(\DD_{\Ad})\longto
    P\defeq \coker\Bigl(\ad(E)^\SC\x\ad(E)^\SC\to 
    \bigl(\phi^*\CoTB_{E\x^G\FRM_\cL/\FRA_{\FRM,\cL}}\bigr)^*\Bigr).
\end{align}
Since the map defining \(P\) is generically surjective, \(P\)
is finitely generated and torsion, and so
\begin{align}
    \RH^0(\breve{X},P)\simeq\prod_v\RH^0(\breve{X}_v,P),
\end{align}
where \(v\) ranges over the points in \(\breve{X}\) that are not sent to the
big-cell locus \(\Stack{\FRM^\circ/G\x Z_\FRM}\). Moreover, one easily sees that
\begin{align}
    \ker\bigl[\coker(\DD_{\Ad})\to P\bigr]^\tfree\simeq
    \coker(\DD_{\Ad}')^\tfree\simeq (\Lie(I_{(E,\phi)}^\Sm)/\La{z}_G)^*.
\end{align}
Thus, when \eqref{eqn:I_sm_vanishing_statement_for_deformation} holds, we have surjective map
\begin{align}
    \RH^0(\breve{X},\coker(\DD_{\Ad}))\longto \prod_v\RH^0(\breve{X}_v,P).
\end{align}
It induces surjective map of tangent spaces
\begin{align}
    \label{eqn:surjectivity_of_global_local_tangent}
    \TanB_{x}(\cM_X/\cB_X)\longto
    \TanB_{\ev(x)}\left(\Stack*{\GASch_X}/\cB_X\right).
\end{align}

\subsection{}
Now we are ready to prove \Cref{thm:local_singularity_model_weak}.
\begin{proof}
    [Proof of \Cref{thm:local_singularity_model_weak}]
    Given a small extension of Artinian \(\bar{k}\)-algebras \(R\to R'=R/I\)
    with residue fields \(\bar{k}\)  and an \(x_{R'}\in\cM_X(R')\) specializing
    to \(x\), let \(\bar{x}_{R'}\) be its image in \(\Stack*{\GASch_X}\).
    Suppose \(\bar{x}_R\) is an \(R\)-lifting of \(\bar{x}_{R'}\), then the local
    obstruction of deforming \(\bar{x}_{R'}\) vanishes. By
    \Cref{prop:global_local_principle_obstructions}, the global obstruction of
    deforming \(x_{R'}\) also vanishes. Here we use the fact that since
    \(G^\SC\) is smooth, the map \(\Stack{\FRM/G^\SC}\to\Stack{\FRM/G^\SC\x
    G^\SC}\) is smooth, hence the map between induced arc stacks is formally
    smooth. Then, by surjectivity of map
    \eqref{eqn:surjectivity_of_global_local_tangent}, there is a lifting of
    \(x_{R'}\) to \(R\) lying over \(\bar{x}_R\), thus we have the desired
    theorem.
\end{proof}

\subsection{}
\label{sub:starting_proof_of_local_singularity_model_main}
Finally, we turn to \Cref{thm:local_singularity_model_main}, whose proof is a
modification of that of \Cref{thm:local_singularity_model_weak}. Although what
follows seems to be very complicated, the idea is in fact very intuitive which we
will explain first.

Looking at the proof of
\Cref{thm:local_singularity_model_weak}, we can see exactly why it is lacking:
if endoscopic group \(H\) happens to be an elliptic torus and \(a\) is
the image of some \(a_H\) for \(H\), then
\(\Lie(\FRJ_a^\flat)\simeq\Lie(\FRJ_{H,a_H}^\flat)\) will have degree \(0\). It
is not a problem on \(H\)-side because the semisimple part
\(\Lie(\FRJ_a^\flat)/\La{z}_H\) is \(0\), but \(\Lie(\FRJ_a^\flat)/\La{z}_G\) is
usually not \(0\), so it cannot have vanishing first cohomology if
the genus \(g_X>1\).

As discussed in \Cref{sec:local_model_of_singularities}, the idea to fix this
issue is to let certain points (in this case, the subdivisor \(b_2\))
move, and similar to the case of mH-base in \Cref{sec:the_case_of_mH_base}, this
will ``add some ampleness'' into the sheaf \(\Lie(\FRJ_a^\flat)/\La{z}_G\)
and hopefully enough to kill the first cohomology. This ``additional ampleness''
can be analyzed by looking at the formal disc around \(b_2\).
By assuming that \(b_2\) and \(\FRD_a\) have disjoint supports, this local
analysis will look exactly like what we have done for mH-base due to the duality
between \(\Lie(\FRJ_a)/\La{z}_G\) and \(\FRC_\cL\) given by
\Cref{prop:desc_of_Lie_FRJ_a}. Our task is just to write down this idea in a
technical manner.

\subsection{}
Recall from
\Cref{sec:local_model_of_singularities} that we have the following commutative
diagram
\begin{equation}
    \begin{tikzcd}
        \cM_X^{I,\disj} \ar[r]\ar[rd] 
        & \Stack*{\GASch_I^{\disj}} \ar[r]\ar[d] 
        & \Stack*{\OGASch[\disj]_I} \ar[r, "\pr_1\circ f_I^{-1}"]\ar[d] 
        & \Stack*{\OGASch_1} \ar[r] \ar[ld]
        & \Stack*{\OGASch_X} \ar[lldd]\\
        & \cB_X^{I,\disj} \ar[r]\ar[rd]
        & \BD[I,\disj]_X \ar[d, "\pr_1"]
        &
        &\\
        &
        & \BD_X
    \end{tikzcd}
\end{equation}
The composition of the top row is the map \(\ev_1^I\), while the composition
\(\cM_X^{I,\disj}\to \BD_X\) of the lower-left maps  is induced by the map of
stacks
\begin{align}
    \Stack*{\FRM^I/G\x Z_{\FRM^I}}\longto\Stack*{\FRA_\FRM/\FRA_\FRM^\x},
\end{align}
which factors through stack \(\Stack*{\FRA_{\FRM^I}/Z_{\FRM^I}}\). As before, we
simplify the notations by letting
\begin{align}
    \Stack*{\FRM^I}&\defeq \Stack*{\FRM^I/Z_{\FRM^I}},\\
    \Stack*{\FRA_{\FRM^I}}&\defeq \Stack*{\FRA_{\FRM^I}/Z_{\FRM^I}},\\
    \Stack*{\FRA_\FRM}&\defeq \Stack*{\FRA_\FRM/\FRA_\FRM^\x}.
\end{align}
The last notation is a bit inconsistent with the first two, but since
\(Z_\FRM\to\FRA_\FRM^\x\) is an isogeny, it does not affect computations related
to cotangent complexes.
We then have distinguished triangle of cotangent complexes
\begin{align}
    \label{eqn:cotangent_triangle_factorized_for_local_model_singularity}
    \alpha_{\FRM^I}^*L_{\Stack{\FRA_{\FRM^I}}/\Stack*{\FRA_\FRM}}\longto
    L_{\Stack*{\FRM^I}/\Stack*{\FRA_\FRM}}\longto
    L_{\Stack*{\FRM^I}/\Stack{\FRA_{\FRM^I}}}\stackrel{+1}{\longto}
\end{align}
associated with maps of stacks
\begin{align}
    \Stack*{\FRM^I}\stackrel{\alpha_{\FRM^I}}{\longto}
    \Stack*{\FRA_{\FRM^I}}\stackrel{\pr_1}{\longto} \Stack*{\FRA_\FRM}.
\end{align}
Since \(\FRA_\FRM\) is of standard type, it may be regarded as a vector
bundle, thus we may identify the first term of
\eqref{eqn:cotangent_triangle_factorized_for_local_model_singularity} with the
complex
\begin{align}
    \alpha_{\FRM^I}^*\pr_2^*\FRA_\FRM^*\longto \La{z}_\FRM^*,
\end{align}
placed at degrees \(0\) and \(1\).

\subsection{}
Let \(x=(\cL,E,\phi)\in \cM_X^{I,\disj}(\bar{k})\) with image
\(a\in\cA_X^{I,\disj}\) and \(b_I\in\cB_X^{I,\disj}\). Let \(b_i\in\BD_X\)
(\(i=1,2\)) be the image of \(b_I\) through the \(i\)-th projection, and
\(\cL_i\) be the \(\FRA_\FRM^\x\)-torsor induced by \(\cL\). We assume
that \(b_2\) is good with \(\supp(\FRD_a)\cap\supp(b_2)=\emptyset\).

We have seen that the deformation of \(x\) in \(\cM_X^I\) relative to
\(\cB_X^I\) is controlled by the truncated complex
\begin{align}
    \ad(E)\stackrel{\DD_{\Ad}^I}{\longto}
    \bigl(\phi^*\CoTB_{E\x^G \FRM^I_\cL/\FRA_{\FRM^I,\cL}}\bigr)^*.
\end{align}
Similarly, the deformation of \(x\) relative to \(\BD_X\) is controlled by
\begin{align}
     \La{z}_\FRM\oplus\ad(E)\stackrel{\DD_{\Ad,1}^I}{\longto}
    \bigl(\phi^*\CoTB_{E\x^G \FRM^I_\cL/\FRA_{\FRM,\cL_1}}\bigr)^*.
\end{align}
We first consider the obstruction spaces. The triangle
\eqref{eqn:cotangent_triangle_factorized_for_local_model_singularity} induces
left-exact sequence of locally-free sheaves
\begin{align}
    0
    \longto \bigl(\phi^*\CoTB_{E\x^G \FRM^I_\cL/\FRA_{\FRM^I,\cL}}\bigr)^*
    \longto \bigl(\phi^*\CoTB_{E\x^G \FRM^I_\cL/\FRA_{\FRM,\cL_1}}\bigr)^*
    \longto \FRA_{\FRM,\cL_2}.
\end{align}
This sequence is in fact also right-exact because over \(b_2\), \(x\) is
contained in the regular locus by assumption, and so the abelianization map is smooth
over \(b_2\) (and the projection \(\pr_1\) is also smooth). As a result, we have
short exact sequence
\begin{align}
    0\longto\coker\bigl(\DD_{\Ad}^I\bigr)\longto \coker\bigl(\DD_{\Ad,1}^I\bigr)
    \longto \FRA_{\FRM,\cL_2}/\La{z}_\FRM\longto 0.
\end{align}
On the other hand, we can also consider the deformation of \(\cA_X^I\) over
\(\BD_X\), and obtain a similar short exact sequence
\begin{align}
    0\longto \FRC_\cL\longto
    (\FRA_{\FRM,\cL_2}\oplus\FRC_\cL)/\La{z}_\FRM\longto \FRA_{\FRM,\cL_2}/\La{z}_\FRM\longto 0.
\end{align}
They fit into a natural diagram
\begin{equation}
    \begin{tikzcd}
        0 \ar[r] & \coker\bigl(\DD_{\Ad}^I\bigr) \ar[r]\ar[d] & \coker\bigl(\DD_{\Ad,1}^I\bigr)\ar[r]\ar[d] & \FRA_{\FRM,\cL_2}/\La{z}_\FRM\ar[r]\ar[d, equal] & 0\\
        0 \ar[r] & \FRC_\cL\ar[r] & (\FRA_{\FRM,\cL_2}\oplus\FRC_\cL)/\La{z}_\FRM\ar[r] & \FRA_{\FRM,\cL_2}/\La{z}_\FRM\ar[r] & 0
    \end{tikzcd}
\end{equation}
Taking \(\cO_{\breve{X}}\)-duals, we have the following diagram with left-exact
rows
\begin{equation}
    \begin{tikzcd}
        0 \ar[r] & \bigl((\FRA_{\FRM,\cL_2}\oplus\FRC_\cL)/\La{z}_\FRM\bigr)^*
        \ar[r]\ar[d] & \FRC_\cL^*\simeq \Lie(\FRJ_a)/\La{z}_G\ar[r]\ar[d] &
        \IExt_{\breve{X}}^1(\FRA_{\FRM,\cL_2}/\La{z}_\FRM,\cO_{\breve{X}})\ar[d,
        equal]\\
        0 \ar[r] & \coker\bigl(\DD_{\Ad,1}^I\bigr)^* \ar[r] &
        \coker\bigl(\DD_{\Ad}^I\bigr)^*\simeq\Lie(I^\Sm_{(E,\phi)})/\La{z}_G
        \ar[r] &
        \IExt_{\breve{X}}^1(\FRA_{\FRM,\cL_2}/\La{z}_\FRM,\cO_{\breve{X}})
    \end{tikzcd}
\end{equation}
Since both \(\FRA_{\FRM,\cL_2}\) and \(\La{z}_\FRM\) are locally free, the torsion sheaf
\(\IExt_{\breve{X}}^1(\FRA_{\FRM,\cL_2}/\La{z}_\FRM,\cO_{\breve{X}})\) is
canonically isomorphic to \(\La{z}_\FRM^*/\FRA_{\FRM,\cL_2}^*\) and is
supported on \(\supp(b_2)\). The quotient
\(\Lie(\FRJ_a^\flat)/\Lie(\FRJ_a)\) is supported on \(\supp(\FRD_a)\), which
is disjoint from \(\supp(b_2)\) by assumption, so the right-most square above
extends to a map
\begin{align}
    \label{eqn:deformation_of_mHiggs_bundles_the_map_for_final_goal}
    \Lie(\FRJ_a^\flat)/\La{z}_G\longto \La{z}_\FRM^*/\FRA_{\FRM,\cL_2}^*
\end{align}
which in turn induces injective maps of locally free sheaves
\begin{align}
    \bigl((\FRA_{\FRM,\cL_2}\oplus\FRC_\cL)/\La{z}_\FRM\bigr)^*\longto
    \coker\bigl(\DD_{\Ad,1}^I\bigr)^* \longto 
    \ker\bigl[\Lie(\FRJ_a^\flat)/\La{z}_G\to
    \La{z}_\FRM^*/\FRA_{\FRM,\cL_2}^*\bigr].
\end{align}
Our goal is to prove that \(\RH^1\bigl(\breve{X},\coker(\DD_{\Ad,1}^I)\bigr)=0\), and
by Serre duality it suffices to prove that
\begin{align}
    \label{eqn:deformation_of_mHiggs_bundles_final_goal}
    \RH^0\Bigl(\breve{X},\ker\bigl[\Lie(\FRJ_a^\flat)/\La{z}_G\to
    \La{z}_\FRM^*/\FRA_{\FRM,\cL_2}^*\bigr]\otimes_{\cO_{\breve{X}}}\CoTB_{\breve{X}}\Bigr)=0.
\end{align}

\subsection{}
The proof of \eqref{eqn:deformation_of_mHiggs_bundles_final_goal}
is similar to that of \eqref{eqn:deformation_mH_base_simplified_goal} and our
strategy is to give a local description of the kernel of
\eqref{eqn:deformation_of_mHiggs_bundles_the_map_for_final_goal} near
\(b_2\), as a subsheaf of \(\FRC_\cL^*\).
Indeed, we have the following diagram
\begin{equation}
    \begin{tikzcd}
        0 \ar[r] & \FRA_{\FRM,\cL_2}^*
        \ar[r]\ar[d] & \FRA_{\FRM,\cL_2}^*\oplus\FRC_\cL^*\ar[r]\ar[d] &
        \FRC_\cL^*\ar[d]\ar[r] & 0\\
        0 \ar[r] & \La{z}_\FRM^* \ar[r, equal] & \La{z}_\FRM^*
        \ar[r] & 0 \ar[r] & 0
    \end{tikzcd}
\end{equation}
whose first row is the cotangent sequence
associated with maps
\begin{align}
    \FRC_{\FRM^I,\cL}\longto\FRA_{\FRM^I,\cL}\stackrel{\pr_1}{\longto}\FRA_{\FRM,\cL_1},
\end{align}
and the vertical maps are the differential actions.
The snake lemma then induces isomorphism
\begin{align}
    \bigl((\FRA_{\FRM,\cL_2}\oplus\FRC_\cL)/\La{z}_\FRM\bigr)^*
    \simeq\ker\Bigl[\FRA_{\FRM,\cL_2}^*\oplus\FRC_\cL^*\to\La{z}_\FRM^*\Bigr]
    \stackrel{\sim}{\longto}\ker\Bigl[\FRC_\cL^*\to
        \La{z}_\FRM^*/\FRA_{\FRM,\cL_2}^*\Bigr],
\end{align}
which is also isomorphic to the kernel of
\eqref{eqn:deformation_of_mHiggs_bundles_the_map_for_final_goal} near
\(\supp(b_2)\) because
\(\supp(\FRD_a)\cap\supp(b_2)=\emptyset\). In addition, let
\(\tilde{a}^\flat\colon \tilde{X}_a^\flat\to\FRT_{\FRM^I,\cL}\) be the natural
map, then  over \(b_2\) it maps to the big-cell locus and the cover
\(\pi_a^\flat\) is \'etale. Therefore, we also have
\begin{align}
    \FRA_{\FRM,\cL_1}^*\oplus\FRA_{\FRM,\cL_2}^*\oplus\FRC_\cL^*\simeq
    \CoTB_{\FRC_{\FRM^I,\cL}}\longto 
    \pi_{a*}^\flat\Bigl[(\tilde{a}^\flat)^*\CoTB_{\FRT_{\FRM^I,\cL}}\Bigr]^W,
\end{align}
which is an isomorphism near \(b_2\). Using the cotangent sequence associated
with maps
\begin{align}
    \FRT_{\FRM^I,\cL}\longto \FRA_{\FRM^I,\cL}\stackrel{\pr_1}{\longto} \FRA_{\FRM,\cL_1},
\end{align}
we finally have a map
\begin{align}
    \bigl((\FRA_{\FRM,\cL_2}\oplus\FRC_\cL)/\La{z}_\FRM\bigr)^*\longto
    \pi_{a*}^\flat\ker\bigl[(\tilde{a}^\flat)^*\CoTB_{\FRT_{\FRM^I,\cL}/\FRA_{\FRM^I,\cL}}\to\La{z}_\FRM^*\bigr]^{W},
\end{align}
and it is an isomorphism near \(b_2\).
This way, locally over \(b_2\), we identify the kernel of
\eqref{eqn:deformation_of_mHiggs_bundles_the_map_for_final_goal} with
the sheaf \(\pi_{a*}^\flat\cK^{W,\tfree}\) in \Cref{sec:the_case_of_mH_base}
as subsheaves of \(\FRC_\cL^*\).
Using the same arguments as in \Cref{sec:the_case_of_mH_base}, we see that if
\(b_2\) is good and
\begin{align}
    \dim_{b_2}\BD_X\ge 2m\abs{\bW}(g_X-1),
\end{align}
and \(L\) is any rank-one locally free subsheaf of the kernel of
\eqref{eqn:deformation_of_mHiggs_bundles_the_map_for_final_goal}, we have
\begin{align}
    \Hom_{\breve{X}}(\CoTB_{\breve{X}}^{-1},L)=0,
\end{align}
which implies \eqref{eqn:deformation_of_mHiggs_bundles_final_goal}.

\subsection{}
Similar to the proof of \Cref{thm:local_singularity_model_weak},
\eqref{eqn:deformation_of_mHiggs_bundles_final_goal} also implies that
the tangent map
\begin{align}
    \TanB_x(\cM_X^I/\BD_X)\longto \TanB_{\ev^I_1(x)}(\Stack*{\OGASch_X}_G/\BD_X)
\end{align}
is surjective. This finishes the proof of
\Cref{thm:local_singularity_model_main}.



\chapter{Stratifications}%
\label{chap:stratifications}

In this chapter we study several important stratifications on the mH-base
\(\cA_X\) in a similar fashion as in \cite{Ng10}*{\S~5}. The first one is
associated with \(\pi_0(\cP_a)\), whose Cartier dual group is the 
group of global endoscopic characters \(\kappa\). The endoscopic monoid
continues to play a key role here.

Another stratification is
associated with \(\delta\)-invariants \(\delta_a\), which is loosely speaking
the dimension of the ``affine part'' of \(\cP_a\). We will prove a codimension
estimate known as \notion{\(\delta\)-regularity}\index{\(\delta\)-!regularity}. For this part, the argument
we use here is parallel to \cite{Ng10}*{\S~5.7} and results in
\Cref{chap:multiplicative_valuation_strata} will be used. Unfortunately, this
means our result on \(\delta\)-regularity has the same drawback as in
\textit{loc. cit.}; however, one can see
\Cref{rmk:unconditional_delta_regularity} for a potential way to strengthen it.

The final section introduces inductive subsets. It further refines some
discussions about \(\delta\)-strata and will become important in
stating the support theorem in \Cref{sec:application_to_mH_fibrations}.

\section{Simultaneous Normalization of Cameral Curves}%
\label{sec:simultaneous_normalization_of_cameral_curves}

\subsection{}
Following \cite{Ng10}*{\S~5.1}, the existence of aforementioned stratifications
depends on a constructibility
result of certain subsets of \(\cA_X^\heartsuit\), whose proof relies on the theory of
simultaneous normalization of a family of curves. 

\begin{definition}
    Let \(f\colon Y\to S\) be a flat and proper map with reduced \(1\)-dimensional
    fibers. A \notion{simultaneous normalization}\index{simultaneous normalization!of curves} of \(Y\to S\) (or just \(Y\)
    if base \(S\) is clear from the context) is a proper birational map
    \(\xi\colon Y^\flat\to Y\) such that
    \begin{enumerate}
        \item There exists open subset \(U\subset Y\) over which \(\xi\) is an
            isomorphism and \(f(U)=S\).
        \item The composition \(f^\flat=f\circ\xi\colon Y^\flat\to S\) is smooth
            and proper.
    \end{enumerate}
\end{definition}

\begin{lemma}[\cite{Ng10}*{Proposition~5.1.2}]
    \label[lemma]{lem:pi_0_of_simultaneous_normalization}
    Let \(\xi\colon Y^\flat\to Y\) be a simultaneous normalization of \(f\colon
    Y\to S\), then
    \begin{enumerate}
        \item \(f_*(\xi_*\cO_{Y^\flat}/\cO_Y)\) is a locally free
            \(\cO_S\)-sheaf of finite type.
        \item There is a locally constant \'etale sheaf \(\pi_0(Y^\flat)\) on
            \(S\) whose fiber at a geometric point \(s\in S\) is the set of
            connected components of \(Y^\flat_s\).
    \end{enumerate}
\end{lemma}

Consider the functor \(\cA_X^\flat\) whose \(S\)-points are 
triples \((a,\tilde{X}_a^\flat,\xi)\) where
\(a\in\cA_X^\heartsuit(S)\), \(\tilde{X}_a^\flat\) is a smooth and proper
\(S\)-family of curves together with a map \(\tilde{X}_a^\flat\to X\x S\)
on which \(W\) acts, and \(\xi\colon \tilde{X}_a^\flat\to \tilde{X}_a\)
is a \(W\)-equivariant simultaneous normalization.
The forgetful functor \(\cA_X^\flat\to \cA_X^\heartsuit\) induces a bijection on
\(\bar{k}\)-points, because for any \(a\in\cA_X^\heartsuit(\bar{k})\) 
the normalization of \(\tilde{X}_a\) is unique.

The functor \(\cA_X^\flat\) has another description. Recall that for
\(a\in\cA_X^\heartsuit(S)\), its image \(\cL\in\Bun_{Z_{\FRM}}(S)\) is
a \(Z_{\FRM}\) torsor on \(X\x S\), which induces map \(\pi_\cL\colon
\FRT_{\FRM,\cL}\to \FRC_{\FRM,\cL}\). Let \(\cA_X^{\flat\prime}\) be the functor of triples
\((\cL,\tilde{X}_a^\flat,\gamma)\) where \(\cL\in\Bun_{Z_{\FRM}}(S)\),
\(\tilde{X}_a^\flat\) is a smooth and proper family of curves over \(S\)
with a finite flat map \(\pi_a^\flat\colon \tilde{X}_a^\flat\to X\x S\) on
which \(W\) acts, and \(\gamma\colon
\tilde{X}_a^\flat\to\FRT_{\FRM,\cL}\) is a \(W\)-equivariant map.
In addition, for any geometric point \(s\in S\) the map
\(\pi_{a,s}^\flat\) is generically a \(W\)-torsor, and the preimage of
\(\FRT_{\cL,s}^\rss\) is dense in \(\tilde{X}_{a,s}^\flat\).
Given \((a,\tilde{X}_a^\flat,\xi)\in\cA_X^\flat(S)\), one may define
\(\gamma\) as the composition of \(\xi\) and the embedding \(\tilde{X}_a\to
\FRT_{\FRM,\cL}\). Thus, we have a natural map \(\cA_X^\flat\to\cA_X^{\flat\prime}\).

\begin{lemma}
    The natural map \(\cA_X^\flat\to\cA_X^{\flat\prime}\) is an
    isomorphism.
\end{lemma}
\begin{proof}
    It suffices to define the inverse map. Let
    \((\cL,\tilde{X}_a^\flat,\gamma)\in\cA_X^{\flat\prime}(S)\), then we claim
    \begin{align}
        \pi_{a*}^\flat(\cO_{\tilde{X}_a^\flat})^W=\cO_{X\x S}.
    \end{align}
    Indeed, the left-hand side is a finite \(\cO_{X\x S}\)-algebra containing
    \(\cO_{X\x S}\) and over any
    geometric point \(s\in S\) they are generically isomorphic on \(X\x
    \Set{s}\). Since \(X\) is normal, it must be an isomorphism.

    Since \(\gamma\) is \(W\)-equivariant, it induces a map \(a\colon X\x
    S\to\FRC_{\FRM,\cL}\) by taking GIT-quotient by \(W\). Let \(\tilde{X}_a\)
    be the corresponding cameral curve, then \(\gamma\) factors through
    \(\tilde{X}_a\to\FRT_{\FRM,\cL}\) hence induces map
    \begin{align}
        \xi\colon \tilde{X}_a^\flat\longto\tilde{X}_a.
    \end{align}
    Clearly over any geometric point \(s\in S\), \(\xi_s\) is a normalization
    map. Thus, we obtain a point
    \((a,\tilde{X}_a^\flat,\xi)\in\cA_X^\flat(S)\). One can verify it is the
    inverse map as desired.
\end{proof}

 Consider
\((a,\tilde{X}_a^\flat,\xi)\in\cA_X^\flat(\bar{k})\), its automorphism
group is a subgroup of \(\ShSec(\breve{X},Z^\SC)\), so in general it
cannot be representable by a scheme. However, the \'etale group 
\(Z^\SC\) is the only ``obstruction'' of \(\cA_X^\flat\) being a sheaf.
Therefore, it is reasonable to
expect that it is a scheme relative to \(\cB_X\). In any case, \(\cB_X\) is
still a Deligne--Mumford stack, even a 
\(\IHom_X(X,Z^\SC)\)-gerbe, thus most topological and representability 
result about schemes and sheaves still apply.

It is not straightforward to prove the representability
of \(\cA_X^\flat\) in full,
but it is enough for our purpose to prove the following slightly weaker
result:
\begin{proposition}
    \label[proposition]{prop:representability_simultaneous_normalization}
    There exists two \(\cB_X\)-schemes \((\cA_X^{\flat})^\Red\to\cB_X\) and
    \(\cA_X''\to \cB_X\) such that we have fully faithful inclusions of functors
    \begin{align}
        (\cA_X^{\flat})^\Red\subset \cA_X^\flat\subset\cA_X'',
    \end{align}
    where the first induces a bijection on \(\bar{k}\)-points, and
    the composition is a closed embedding.
\end{proposition}
\begin{proof}
    We first consider functor \(\sH\) whose \(S\)-points are the
    isomorphism classes of the following data: a smooth and projective
    \(S\)-curve \(\tilde{X}_a^\flat\) together with a map
    \(\tilde{X}_a^\flat\to X\x S\) on which \(W\) acts, such that 
    over any geometric point \(s\in S\) it is generically a \(W\)-torsor.

    Suppose for the moment that \(\sH\) is representable by a
    quasi-projective \(k\)-scheme, then we have a universal family of smooth and
    proper curves
    \begin{align}
        \tilde{X}_\sH^\flat\longto \sH
    \end{align}
    together with a \(W\)-equivariant map \(\tilde{X}_\sH^\flat\to X\x\sH\).

    This induces a \(\FRT_{\FRM}\)-bundle \(\FRT_{\sH}\) over
    \(\tilde{X}_\sH^\flat\x\Bun_{Z_{\FRM}}\), and the fibers of forgetful map
    \(\cA_X^\flat\to\sH\x\Bun_{Z_{\FRM}}\) is just the \(W\)-equivariant sections of
    \(\FRT_{\sH}\) over the curve \(\tilde{X}_\sH^\flat\). Using the
    representation-theoretic description of \(\Env(\bG^\SC)\), we know that we
    may embed \(\FRT_{\FRM}\) into a vector bundle \(\FRE\) with \(W\)-action as a
    closed subbundle. Moreover, \(\FRE\) can be so chosen that it maps to
    \(\FRA_{\FRM}\) compatible with \(\FRT_{\FRM}\to\FRA_{\FRM}\) and the
    \(W\)-action. One can pull it back to get a \(W\)-vector bundle
    \(\FRE_\sH\) on \(\tilde{X}_\sH^\flat\x\Bun_{Z_{\FRM}}\).

    We define a functor \(\cA_X''\) same as \(\cA_X^{\flat\prime}\) except
    that \(\gamma\) is a \(W\)-equivariant map to \(\FRE\) instead of
    \(\FRT_{\FRM}\). By the definition of \(\FRE\), the map
    \(\cA_X''\to\Bun_{Z_{\FRM}}\) factors through \(\cB_X\), and it is clear
    that \(\cA_X''\) is a sheaf relative to \(\cB_X\). 
    It is also clear that \(\cA_X''\) is the \(W\)-fixed points of the
    section space of a vector bundle on \(\sH\x\Bun_{Z_{\FRM}}\), hence it is
    representable by Deligne--Mumford stack which is also a \(\cB_X\)-scheme. 
    Moreover, it contains \(\cA_X^\flat\) as a fully faithful
    subfunctor. Consider Cartesian diagram
    \begin{equation}
        \begin{tikzcd}
            \tilde{X}' \ar[r]\ar[d] & \FRT_\sH \ar[d]\\
            \tilde{X}_{\cA_X''}^\flat \ar[r, "\gamma"] & \FRE_\sH
        \end{tikzcd}
    \end{equation}
    where \(\tilde{X}_{\cA_X''}^\flat\) is the universal curve on
    \(\cA_X''\) and the bottom horizontal map is the universal map
    \(\gamma\) over \(\cA_X''\). The vertical maps in the above square are
    closed embeddings. Since
    \(\tilde{X}_{\cA_X''}\to\cA_X''\) is proper, the
    image of \(\tilde{X}'\) in \(\cA_X''\) is closed. By
    upper-semicontinuity of fiber dimensions for proper maps, the locus
    in \(\cA_X''\) where the fiber in \(\tilde{X}'\) 
    has dimension \(1\) is also a closed
    subset. Let \((\cA_X^\flat)^\Red\) be the reduced substack structure of this
    closed subset.

    We want to show that \(\tilde{X}'\to\tilde{X}_{\cA_X''}^\flat\)
    is an isomorphism over \((\cA_X^\flat)^\Red\). We already know that their
    fibers have the same dimension and are proper and reduced (in fact regular),
    but \(\tilde{X}_a^\flat\) is not necessarily irreducible.
    Nevertheless, since \(\gamma\) is \(W\)-equivariant and since
    \(\pi_a^\flat\) is generically a \(W\)-torsor, once a fiber of
    \(\tilde{X}'\) contains one irreducible component of
    \(\tilde{X}_a^\flat\), it has to contain all of them. This means that
    \((\cA_X^\flat)^\Red\) is a subfunctor of \(\cA_X^\flat\). It is also clear
    that \((\cA_X^\flat)^\Red\) has the same \(\bar{k}\)-points as
    \(\cA_X^\flat\). Thus, it remains to prove that \(\sH\) is representable a
    quasi-projective \(k\)-scheme. We leave it to a separate proposition below.
\end{proof}

\begin{proposition}
    The moduli functor \(\sH\) as in
    \Cref{prop:representability_simultaneous_normalization} is
    representable by a quasi-projective \(k\)-scheme.
\end{proposition}
\begin{proof}
    As commented in \cite{Ng10}*{\S~5.2}, the proof is similar to the
    representability of Hurwitz schemes. We include a proof for completeness, as
    it is not in \textit{loc.~cit.}

    Let \(S\) be a \(k\)-scheme and let \(f\colon\tilde{X}_a^\flat\to X\x S\) be
    an \(S\)-point of \(\sH\). For simplicity, we denote \(\tilde{X}_a^\flat\) by
    \(Y\) and \(X\x S\) by \(X_S\). The relative cotangent complex
    \(\CoTB_{Y/X_S}\) is the coherent sheaf \(\CoTB_{Y/S}/f^*\CoTB_{X_S/S}\)
    since \(Y\) and \(X_S\) are smooth curves over \(S\) and \(f\) is finite
    flat. It is a coherent sheaf of finite length, whose induced divisor on
    \(Y\) is the ramification divisor \(D_Y\) of \(f\). Since \(D_Y\) is stable
    under \(W\), it descends to a divisor \(B_Y\) on \(X_S\), finite and
    flat over \(S\). Therefore, we obtain a map
    \begin{align}
        f_B\colon\sH&\longto H_X\defeq\coprod_n \Hilb^n{X}\\
        Y&\longmapsto B_Y.
    \end{align}
    It then suffices to show that the map \(f_B\) is representable and
    \'etale, and we shall use the same criteria in \cite{Fu69}*{Theorem~6.5},
    which is a slightly modified version of a result due to Grothendieck
    reported in \cite{Mu95}. Explicitly, we need to show the followings:
    \begin{enumerate}
        \item \label{item:Hurwitz_rep_crit_sheaf}
            \(\sH\) is an fpqc sheaf over \(H_X\).
        \item \label{item:Hurwitz_rep_crit_ind_limit}
            \(\sH\) commutes with filtered inductive limit of rings over
            \(H_X\).
        \item \label{item:Hurwitz_rep_crit_formal_to_local}
            If \(A\) is a complete Noetherian local ring over \(H_X\) and
            \(\FRm\subset A\) is the maximal ideal, then
            \begin{align}
                \sH(A)\longto\varprojlim_i\sH(A/\FRm^i)
            \end{align}
            is bijective. In other words, formal deformation of \(f_B\) over
            \(H_X\) can be uniquely promoted to a local deformation.
        \item \label{item:Hurwitz_rep_crit_formal_sm_and_unramified}
            If \(R\) is an Artinian ring over \(H_X\) and \(I\subset R\) is a
            nilpotent ideal, then \(\sH(R)\to\sH(R/I)\) is bijective. In other
            words, \(f_B\) is formally \'etale.
        \item \label{item:Hurwitz_rep_crit_separated}
            \(f_B\) satisfies the uniqueness part of valuative criteria for
            all complete discrete valuation rings.
    \end{enumerate}

    Property~\eqref{item:Hurwitz_rep_crit_sheaf} is immediate from definition.

    Property~\eqref{item:Hurwitz_rep_crit_ind_limit}: let \(A_i\) be a filtered
    inductive system of rings over \(H_X\), and \(A\) is the limit of \(A_i\).
    Then \(Y_{A_i}\to\Spec{A_i}\) is finitely presented, quasi-compact and
    quasi-separated. Moreover, \(H_X\) is locally quasi-compact and
    quasi-separated. By \cite{EGAIV3}*{Th\'eor\`eme~8.8.2}, the natural map
    \begin{align}
        \varinjlim_i\sH(A_i)\longto \sH(A)
    \end{align}
    is bijective.

    Properties~\eqref{item:Hurwitz_rep_crit_formal_to_local}:
    is due to the fact that formal deformation of projective
    curves can always be uniquely algebraicized
    (\cite{EGAIII1}*{Th\'eor\`eme~5.4.1}).

    Properties~\eqref{item:Hurwitz_rep_crit_formal_sm_and_unramified} is proved
    using deformation theory. Let \(R'=R/I\) and \(f'\colon Y'\to X_{R'}\).
    Without loss of generality, we may assume that \(R\) is a small extension of
    \(R'\). Let \(f_0\) be the fiber of \(f'\) over the residue field \(k_R\) of
    \(R'\). By \cite{Il71}*{\S~III.2, Th\'eor\`eme~2.1.7}, the obstruction of flatly
    deforming map \(f'\) (without regarding \(W\)-action) lies in the group
    \begin{align}
        \Ext^2_{\cO_{Y'}}(\CoTB_{Y'/X_{R'}},I)
    \end{align}
    which vanishes by Serre duality. Therefore, all flat deformations of \(f'\)
    to \(R\) is a torsor under vector space
    \begin{align}
        \TanB_I&\defeq\Ext^1_{\cO_{Y'}}(\CoTB_{Y'/X_{R'}},I)\\
               &\cong\Ext^1_{\cO_{Y'}}(\CoTB_{Y'/X_{R'}},\cO_{Y'})\otimes
               I.
    \end{align}
    The group \(W_X=\ShSec(X,W)\) acts canonically on the set of deformations by
    composition, in other words, \(w\in W_X\) sends a deformation \(\iota\colon
    Y'\to Y\) to \(\iota\circ w^{-1}\). The tangent space \(\TanB_I\)
    also admits a canonical \(W_X\)-action, thus the set of formations induces a
    class in \(\RH^1(W_X,\TanB_I)\). Since \(\Char(k)\) does not divide the order
    of \(\bW\), the last cohomology group is trivial, being a vector space
    annihilated by the order of \(W_X\). Thus, we identify the set of
    deformations with \(\TanB_I\) as \(W_X\)-sets. Let \(f\in\TanB_I^{W_X}\)
    be a \(W_X\)-invariant deformation, the obstruction of extending
    \(W_X\)-action to \(f\) lies in the group
    \begin{align}
        \RH^2(W_X,\Aut_{f'}(f))=\RH^2(W_X,\Aut_{f_0}(f_0[\epsilon])\otimes
        I),
    \end{align}
    where \(f_0[\epsilon]\) is the trivial extension of \(f_0\) to
    \(k_R[\epsilon]/\epsilon^2\). This obstruction group also vanishes since it
    takes value in a vector space.

    These constructions can also be done locally
    over \(X\), hence we see that the set of deforming \(f'\) together with
    \(W\)-action is identified with the subspace of \(\TanB_I^{W_X}\) such that
    over any open subset \(U\subset X\), it is also fixed by \(\ShSec(U,W)\).
    Since \(f'\) is finite flat, we have by adjunction
    \begin{align}
        \TanB_I\cong\Ext_{\cO_{X_{R'}}}^1(\CoTB_{Y'/X_{R'}}^W,\cO_{Y'})\otimes
        I,
    \end{align}
    hence the deformation of \(Y'\) with \(W\)-action can be identified with
    vector space
    \begin{align}
        \Ext_{\cO_{X_{R'}}}^1(\CoTB_{Y'/X_{R'}}^W,\cO_{X_{R'}})\otimes I.
    \end{align}
    On the other hand, it is well-known (see \cite{FaGoIl05}*{Corollary~6.4.10}
    for example) that the deformation of
    \(B_{Y'}\) in \(H_X\) by \(I\) is identified with
    \begin{align}
        \Hom_{\cO_{X_{R'}}}(J,\cO_{B_{Y'}})\otimes I\simeq
        \Ext_{\cO_{X_{R'}}}^1(\cO_{B_{Y'}},\cO_{B_{Y'}})\otimes I,
    \end{align}
    where \(J\subset\cO_{X_{R'}}\) is the ideal of \(B_{Y'}\). Therefore, it
    boils down to showing that the derivative
    \begin{align}
        \DD_{f_0}(f_B)\colon\Ext_{\cO_{X_{R'}}}^1(\CoTB_{Y'/X_{R'}}^W,\cO_{X_{R'}})
        \longto \Ext_{\cO_{X_{R'}}}^1(\cO_{B_{Y'}},\cO_{B_{Y'}})
    \end{align}
    induced by tensoring the source by \((\CoTB_{Y'}^W)^{\vee}\) is an isomorphism,
    but this is simply the definition of \(B_{Y'}\) combined with Serre duality,
    hence we are done.

    Property~\eqref{item:Hurwitz_rep_crit_separated} is straightforward: if
    \(A\) is a complete discrete valuation ring and \(F\) its function field and
    \(Y\in \sH(A)\), then since \(Y\) is smooth over \(A\) its ring of functions
    can be characterized as the integral closure of \(\cO_{X_A}\) in
    \(\cO_{Y_F}\). Therefore, \(Y\) is completely determined by its fiber over
    \(F\). This finishes the proof.
\end{proof}

\section{Stratification by Monodromy}%
\label{sec:stratification_by_monodromy}

\subsection{}
Using the moduli stack \(\cA_X^\flat\), we may now study the stratifications on
\(\cA_X^\heartsuit\). Similar to \cite{Ng10}*{\S~5}, it is more convenient to
study an \'etale cover of \(\cA_X\) since it will simplify the description of
\(\pi_0(\cP_X)\) (see \Cref{sec:the_sheaf_pi_0_cp_x}) and the resulting
stratifications. Recall that after fixing
\(\infty\in X(\bar{k})\), we have the open subset \(\cA_X^\infty\subset\cA_X\)
consisting of points \(a\) such that the cameral cover \(\tilde{X}_a\) is
\'etale over \(\infty\). If \(\infty\in X(k)\), then \(\cA_X^\infty\) has a
natural \(k\)-structure. Consider the functor \(\tilde{\cA}_X\) whose
\(S\)-points are pairs \((a,\tilde{\infty})\) where \(a\in\cA_X^\infty(S)\)
and \(\tilde{\infty}\in \tilde{X}_a(S)\) lying over \(\infty_S\). It is
clear that \(\tilde{\cA}_X\) is an \'etale subset of \(\cA_X\), and has a
\(k\)-structure if \(\infty\in X(k)\).

We have the natural map \(\cA_X^\infty\to\cB_X\), and the point \(\infty\)
defines a locally trivial fiber bundle \(\FRC^\infty\) over \(\cB_X\) whose
fiber at \((\cL,\theta)\) is the fiber of \(\FRC_\cL\) at \(\infty\) (recall
that \(\FRC\) is the GIT quotient \(G^\SC\git G\)). In addition,
\(\theta\in\ShSec(X,\FRA_{\FRM,\cL})\) induces a cameral cover
\(\FRT_{(\cL,\theta)}\) of \(\FRC_\cL\), which further induces cameral cover
\(\FRT^\infty\to\FRC^\infty\). Consider Cartesian diagram 
\begin{equation}
    \begin{tikzcd}
        \tilde{\cA}_X \ar[r]\ar[d] & \FRT^\infty \ar[d]\\
        \cA_X^\infty \ar[r] & \FRC^\infty
    \end{tikzcd}.
    \nomenclature[\(A"cal_X'tilde\)]{\(\tilde{\cA}_X\)}{the \'etale cover of
    \(\cA_X^\infty\) by choosing a lift of \(\infty\) to a point in the cameral curve}
\end{equation}

\begin{lemma}
    Over the very \(G\)-ample locus \(\cB_{\gg}\),
    \(\tilde{\cA}_{\gg}\to\cB_{\gg}\) is smooth and induces bijection on
    irreducible components. 
\end{lemma}
\begin{proof}
    We already know \(\cA_X\to\cB_X\) is a vector bundle when restricted to
    \(\cB_{\gg}\), and \(\tilde{\cA}_{\gg}\) is \'etale over \(\cA_{\gg}\),
    hence the smoothness result.
    Also by very \(G\)-ampleness, the bottom arrow in the above Cartesian diagram is
    a surjective map of vector bundles. This means that the top horizontal map
    has connected fibers. It remains to show that the preimage of any
    geometric irreducible component of \(\cB_{\gg}\) is irreducible. Replacing
    \(\cB_\gg\) by its irreducible components, we need to show that
    \(\tilde{\cA}_\gg\) is irreducible.

    If \(\infty\) lies in the invertible locus, the conclusion is immediate
    since \(\FRT^\infty\), being a torsor under \(T^\SC\), is itself
    irreducible. If \(\infty\) does not lie in the invertible locus, a
    perturbation \(\infty_S\to S\) (\(S\) is some \(1\)-dimensional
    scheme over \(\bar{k}\)) of \(\infty\) into the invertible locus would induce a
    deformation of \(W_\infty\)-cover (\(W_\infty\) is the fiber of \(W\) at
    \(\infty\)) \(\tilde{\cA}_{\gg,S}\to\cA_\gg^{\infty_S}\). But then it is a
    \(S\)-flat
    family of \(W_{\infty_S}\)-\'etale covers of which the generic member has
    irreducible total space. This is impossible unless \(\tilde{\cA}_\gg\) is
    also irreducible, hence the lemma.
\end{proof}

Now we define
\(\tilde{\cA}_{\gg}^\flat=\cA_X^\flat\x_{\cA_X}\tilde{\cA}_{\gg}\). The
map \(\tilde{\cA}_{\gg}^\flat\to\tilde{\cA}_{\gg}\) is a bijection on
\(\bar{k}\)-points. With the help of the reduced substack
\((\cA_X^\flat)^\Red\), the image of any connected components of
\((\tilde{\cA}_{\gg}^\flat)_{\bar{k}}\) in
\((\tilde{\cA}_{\gg})_{\bar{k}}\) is a constructible subset relative to
\(\cB_{\gg,\bar{k}}\). 

Over any geometric connected component of \(\cB_{\gg}\),
one may decompose this constructible subset into a finite union of locally
closed subsets. Let \(\tilde{\cA}_1\) be one of such locally closed subsets
and \(\tilde{\cA}^\flat_1\) be its preimage in
\((\tilde{\cA}_{\gg}^\flat)_{\bar{k}}\).
By Zariski Main Theorem, there exists an open dense subset of
\(\tilde{\cA}_1\) over which the map
\((\tilde{\cA}^\flat_1)^\Red\to\tilde{\cA}_1\) is finite radical. 
Using Noetherian induction, we may refine those \(\tilde{\cA}_{\gg}\) into a locally
closed stratification
\begin{align}
    \label{eqn:A_tilde_stratification_prelim}
    (\tilde{\cA}_{\gg})_{\bar{k}}=\coprod_{\psi\in\Psi}\tilde{\cA}_{\psi},
\end{align}
such that the map \((\tilde{\cA}_\psi^\flat)^\Red\to\tilde{\cA}_\psi\) is
finite radical. We may further refine the stratification so that the closure of
any \(\tilde{\cA}_\psi\) is a finite union of strata. Hence, we have a
partial order on \(\Psi\) such that \(\psi\le\psi'\) if \(\tilde{\cA}_\psi\) is
contained in the closure of \(\tilde{\cA}_{\psi'}\). 
Since \(\tilde{\cA}_{\gg}\to\cB_{\gg}\) is smooth and a bijection on irreducible
components, over each geometric irreducible component \(B\) of \(\cB_X\) there
is a maximal element \(\psi_{B}\in\Psi\).

\subsection{}
Recall \(G\) comes from an \(\Out(\bG)\)-torsor \(\OGT_G\) on \(X\), and when we
fix a point \(\infty\), we may lift it to a pointed version \(\OGT_G^\bullet\),
and \(\Theta=\Theta_\OGT\) is the image of \(\pi_1(\breve{X},\infty)\) in \(\Out(\bG)\)
under \(\OGT_G^\bullet\). By assumption, \(\Theta\) is finite, and its
order is not divided by \(\Char(k)\).

Let \(\tilde{a}=(a,\tilde{\infty})\in\tilde{\cA}_X(\bar{k})\).
Recall the commutative diagram \eqref{eqn:I_W_diagram} which we reproduce here:
\begin{equation}
    \label{eqn:I_W_diagram2}
    \begin{tikzcd}
        \pi_1(U,\infty) \ar[r, "\pi_{\tilde{a}}^\bullet"]\ar[d] 
                & \bW\rtimes\Out(\bG) \ar[d]\\
\pi_1(\breve{X},\infty) \ar[r, "\OGT_G^\bullet"] & \Out(\bG)
    \end{tikzcd}.
\end{equation}
In this diagram \(U=\breve{X}-\FRD_a\).
Let \(W_{\tilde{a}}\) be the image of \(\pi_{\tilde{a}}^\bullet\) and
\(I_{\tilde{a}}\) the image of the kernel of
\(\pi_1(U,\infty)\to\pi_1(\breve{X},\infty)\) under
\(\pi_{\tilde{a}}^\bullet\). By construction, \(W_{\tilde{a}}\) is
contained in \(\bW\rtimes\Theta\), and \(I_{\tilde{a}}\) is a normal
subgroup of \(W_{\tilde{a}}\) and contained in \(W_{\tilde{a}}\cap\bW\).

An alternative approach to the above setup is as follows.
Let \(X_\OGT\to\breve{X}\) be a connected finite Galois \'etale cover with Galois
group \(\Theta'\) together with a point \(\infty_\OGT\) lying over \(\infty\)
such that \(\OGT_G\) becomes trivial on \(X_\OGT\). Such
requirement is the same as saying
\(\OGT_G^\bullet\colon\pi_1(\breve{X},\infty)\to\Out(\bG)\) factors through
\(\Theta'\), thus we may replace \(X_\OGT\) by a quotient cover so that
\(\Theta'=\Theta\). Let 
\begin{align}
    \tilde{X}_{\OGT,a}=\tilde{X}_a\x_{\breve{X}}X_\OGT,
\end{align}
then \(\bW\rtimes \Theta\) acts on \(\tilde{X}_{\OGT,a}\) as well as its
normalization \(\tilde{X}_{\OGT,a}^\flat\). Let \(C_{\tilde{a}}\subset
\tilde{X}_{\OGT,a}^\flat\) be the connected component containing
\(\tilde{\infty}_\OGT=(\tilde{\infty},\infty_\OGT)\). 
Let \(W_{\tilde{a}}\) be the subgroup of \(\bW\rtimes\Theta\) mapping
\(C_{\tilde{a}}\) to itself, and \(I_{\tilde{a}}\subset
W_{\tilde{a}}\) is the subgroup generated by elements with at least one
fixed point in \(C_{\tilde{a}}\). 

Since \(U\) is the maximal subset of
\(\breve{X}\) over which the cover \(\tilde{X}_{\OGT,a}\to \breve{X}\) is
a \(\bW\rtimes\Theta\)-torsor, \(W_{\tilde{a}}\) is exactly the image of
\(\pi_1(U,\infty)\) in \(\bW\rtimes\Theta\). For any closed point
\(\tilde{v}\in\tilde{X}_{\OGT,a}^\flat\) with image \(\bar{v}\in\breve{X}\), 
the elements in \(I_{\tilde{a}}\)
fixing \(\tilde{v}\) generates local inertia group of
\(\tilde{X}_{\OGT,a,\tilde{v}}\)
over \(\breve{X}_{v}\). Since \(I_{\tilde{a}}\) is by definition generated by
these elements, it is equal to the kernel of map
\(\pi_1(U,\infty)\to\pi_1(\breve{X},\infty)\), by a version of Riemann Existence
Theorem (see e.g., \cite{SGA1}*{Expos\'e~XIII, Corollaire~2.12}). Finally,
since
\(C_{\tilde{a}}\) maps onto \(X_\OGT\) and since \(\Theta\) acts freely on
\(X_\OGT\), the projection of \(I_{\tilde{a}}\) to \(\Theta\) is trivial.
Therefore, \(I_{\tilde{a}}=W_{\tilde{a}}\cap\bW\). Thus, we established
the equivalence of two formulations of pairs
\((I_{\tilde{a}},W_{\tilde{a}})\) associated with \(\tilde{a}\) (and
\(\OGT_G^\bullet\), which is independent of \(\tilde{a}\)).

\begin{proposition}
    The map \(\tilde{a}\mapsto (I_{\tilde{a}},W_{\tilde{a}})\) is constant on
    every stratum \(\tilde{\cA}_\psi\) in
    \eqref{eqn:A_tilde_stratification_prelim}. As a result, we have a
    well-defined map on the set of strata \(\Psi\) 
    \begin{align}
        \psi\longmapsto (I_\psi,W_\psi).
    \end{align}
\end{proposition}
\begin{proof}
    We have that \((\tilde{\cA}_\psi^\flat)^\Red\to\tilde{\cA}_\psi\) is a
    finite radical morphism. Over \((\tilde{\cA}_\psi^\flat)^\Red\), the
    normalizations of individual cameral curves form a smooth family
    \(\tilde{X}_{\psi}^\flat\to\tilde{\cA}_\psi^\flat\). In fact, using the
    cover \(X_\OGT\), we also have the smooth proper family
    \begin{align}
        \tilde{X}_{\OGT,\psi}^\flat\longto \tilde{\cA}_\psi^\flat
    \end{align}
    on which \(\bW\rtimes\Theta\) acts, as well as a section
    \(\tilde{\infty}_\OGT\). Using
    \Cref{lem:pi_0_of_simultaneous_normalization}, we have a locally
    constant sheaf \(\pi_0(\tilde{X}_{\OGT,\psi}^\flat/\tilde{\cA}_\psi^\flat)\)
    whose fibers are the connected components of the fibers of
    \(\tilde{X}_{\OGT,\psi}^\flat\), and \(\bW\rtimes\Theta\) acts transitively
    on the fibers. The existence of section \(\tilde{\infty}_\OGT\) means that
    \(\pi_0(\tilde{X}_{\OGT,\psi}^\flat/\tilde{\cA}_\psi^\flat)\) is constant.
    This means that the collection of \(C_{\tilde{a}}\) forms a smooth and
    proper family over \(\tilde{\cA}_\psi^\flat\) with connected fibers. Since
    \(\bW\rtimes\Theta\) is discrete, we see that \(W_{\tilde{a}}\) and
    \(I_{\tilde{a}}\) are constant.
\end{proof}

\begin{definition}
    We define a partial order on the set of pairs \((I_\psi,W_\psi)\) where
    \(W_\psi\) is a subgroup of \(\bW\rtimes\Theta\) and \(I_\psi\subset
    W_\psi\) is a normal subgroup contained in \(\bW\): \((I_\psi,W_\psi)\le
    (I_\psi',W_\psi')\) if \(W_\psi\subset W_\psi'\) and \(I_\psi\subset
    I_\psi'\).
\end{definition}

\begin{lemma}
    \label[lemma]{lem:I_W_increasing}
    The map \(\psi\mapsto (I_\psi,W_\psi)\) is an increasing map on \(\Psi\).
\end{lemma}
\begin{proof}
    Let \(S=\Spec{R}\) be the spectrum of a complete discrete valuation ring, with
    generic point \(\eta=\Spec{K}\) and special point \(s=\Spec{k(s)}\), 
    with \(k(s)\) being algebraically closed. Let \(\tilde{a}\colon
    S\to\tilde{\cA}\) be a morphism sending
    \(s\) to \(\tilde{\cA}_\psi\) and \(\eta\) to \(\tilde{\cA}_{\psi'}\). We
    want to show that \((I_{\psi},W_{\psi})\le (I_{\psi'},W_{\psi'})\).

    Since finite \'etale coverings can only be trivially deformed locally, we
    have canonical cospecialization maps
    \(\pi_1(U_s,\infty_s)\to\pi_1(U_\eta,\infty_\eta)\) and
    \(\pi_1(\breve{X}_s,\infty_s)\to\pi_1(\breve{X}_\eta,\infty_\eta)\) compatible
    with diagram \eqref{eqn:I_W_diagram2}. Then the result follows from the
    definition.
\end{proof}

\subsection{}
\label{sub:anisotropic_locus_is_open}
Using \Cref{lem:I_W_increasing}, if \((I_-,W_-)\) is a pair where \(W_-\)
is a subgroup of \(\bW\rtimes\Theta\) and \(I_-\subset W_-\) is a normal
subgroup also contained in \(\bW\), then the union of such \(\tilde{\cA}_\psi\)
that \((I_\psi,W_\psi)\le(I_-,W_-)\) is a closed subset of \(\tilde{\cA}_\psi\),
and the union \(\tilde{\cA}_{(I_-,W_-)}\) of those strata satisfying
\((I_\psi,W_\psi)=(I_-,W_-)\) is an open subset of that closed subset. Thus, we
have a locally closed stratification
\begin{align}
    (\tilde{\cA}_X)_{\bar{k}}=\coprod_{(I_-,W_-)}\tilde{\cA}_{(I_-,W_-)}.
    \nomenclature[\(.I_W\)]{\((\cdot)_{(I_-,W_-)}\)}{related to a stratum in
    \(\tilde{\cA}_X\) of monodromy type \((I_-,W_-)\)}
\end{align}
In particular, the union of those strata such that \(\bT^{W_\psi}\) is finite is
an open subset of \(\tilde{\cA}_X\). Combined with
\Cref{prop:global_pi_0_description}, we see that \(\cA_X^\ANI\), or more
generally \(\cA_X^{\ELL}\)
is an open subset of \(\cA_X^\heartsuit\). The stratifications on
\(\tilde{\cA}_X\) naturally induce stratifications on
\(\tilde{\cA}_X^\ANI=\cA_X^\ANI\x_{\cA_X}\tilde{\cA}_X\)

\begin{lemma}
    Let \(B\subset \cB_{\gg}\) be an irreducible component, and \(\psi_B\) be
    the maximal element in \(\Psi\) corresponding to \(B\). Then
    \begin{align}
        (I_{\psi_B},W_{\psi_B})=(\bW,\bW\rtimes\Theta).
    \end{align}
    In fact, for any \(b\in\cB_{\gg}(\bar{k})\), we can find \(\tilde{a}\)
    lying over \(b\) such that
    \((I_{\tilde{a}},W_{\tilde{a}})=(\bW,\bW\rtimes\Theta)\).
\end{lemma}
\begin{proof}
    For any \(b\in\cB_{\gg}(\bar{k})\), we know \(\cA_b^\diamondsuit\) is
    non-empty by
    \Cref{prop:non_emptyness_of_A_sharp_diamond_heart}.
    Let \(\tilde{a}\in\tilde{\cA}_{\gg}^\diamondsuit(\bar{k})\) be a point
    lying over \(b\). Then we know \(\tilde{X}_{\OGT,a}\) is smooth and
    irreducible, therefore \(W_{\tilde{a}}=\bW\rtimes\Theta\). We also know that
    \(\tilde{X}_{\OGT,a}\) intersects with \(\bD_\alpha\) for every positive
    root \(\alpha\).
    Therefore, \(I_{\tilde{a}}\) is a normal subgroup of \(\bW\)
    containing every reflection, thus must be \(\bW\) itself.
\end{proof}

\begin{corollary}
    \label[corollary]{cor:non_emptyness_of_anisotropic_locus}
    For any \(b\in\cB_\gg\), the substack \(\cA_b^{\ELL}\) is non-empty.
    If \(Z_G\) does not contain a split subtorus, then the same holds for
    \(\cA_b^{\ANI}\).
\end{corollary}
\begin{proof}
    It immediately follows from that
    \(\cA_b^{\ELL}\subset\cA_b^\diamondsuit\), and the fact that
    \(\cA_X^{\ELL}=\cA_X^{\ANI}\) if \(Z_G\) does not contain a split subtorus.
\end{proof}

\section{The Sheaf \texorpdfstring{\(\pi_0(\cP_X)\)}{π0(PX)}}%
\label{sec:the_sheaf_pi_0_cp_x}

We have already seen the descriptions of \(\pi_0(\cP_{\tilde{a}})\) and
\(\pi_0(\cP_{\tilde{a}}')\) for individual
\(\tilde{a}\in\tilde{\cA}_X^\heartsuit(\bar{k})\) in
\Cref{prop:global_pi_0_description}. By \cite{Ng06}*{\S~6.2}, there exists an
\'etale sheaf \(\pi_0(\cP_X)\) over \(\cA_X^\heartsuit\) such that its fiber
over a geometric point \(a\) is \(\pi_0(\cP_a)\). Here we want to describe the
restriction of sheaves \(\pi_0(\cP_X)\) and \(\pi_0(\cP_X')\) to each stratum
\(\tilde{\cA}_{(I_-,W_-)}\).

\begin{proposition}
    \label[proposition]{prop:pi_0_Picard_as_lattice_quotient}
    There exists canonical surjections of sheaves over \(\tilde{\cA}_X\)
    \begin{align}
        \CoCharG(\bT)\longto \pi_0(\cP_X')\longto \pi_0(\cP_X)
    \end{align}
    such that the fiber at any \(\tilde{a}\in\tilde{\cA}_X(\bar{k})\) are the
    surjections given in \Cref{prop:global_pi_0_description}.
\end{proposition}
\begin{proof}
    The section \(\tilde{\infty}\) of the cameral curve over \(\tilde{\cA}_X\)
    gives a fixed \(\tilde{\cA}_X\)-family of pinnings of
    \(G\) at \(\infty\). Hence, using the Galois description of \(\FRJ_X\), we
    may identify the fiber of \(\FRJ_X\) over \(\Set{\infty}\x\tilde{\cA}_X\) 
    with \(\bT\x\tilde{\cA}_X\). Such identification uniquely extends to the
    formal disc \(\breve{X}_\infty\x\tilde{\cA}_X\), because \(\CoCharG(\bT)\) is
    discrete. For each \(\mu\in\CoCharG\simeq\Gr_{\bT}^\Red\), we have an
    induced local \(\FRJ_{\breve{X}_\infty}\)-torsor \(E_\mu\) over
    \(\tilde{\cA}_X\). Gluing with the trivial \(\FRJ_X\)-torsor over the
    complement of \(\Set{\infty}\x\tilde{\cA}_X\), we obtain a well-defined map
    \(\CoCharG(\bT)\to\pi_0(\cP_X)\). The same argument works for \(\cP_X'\) as
    well since at \(\infty\), the fibers of \(\FRJ_X\) and \(\FRJ_{X}^0\) are
    the same. The surjectivity can be checked at stalk level using
    \Cref{prop:global_pi_0_description}.
\end{proof}

For any \'etale open subset \(U\to\tilde{\cA}_X\), the stratifications of
\(\tilde{\cA}_X\) induce stratifications of \(U\). In particular, we have strata
\(U_{(I_-,W_-)}\) for pairs \((I_-,W_-)\). For a given pair \((I_1,W_1)\), we
say that \(U\) is \notion{of type \((I_1,W_1)\)}\index{stratum!of monodromy type
\((I_-,W_-)\)}
if \(U_{(I_1,W_1)}\) is the unique non-empty closed stratum in
\(U\). Clearly, \'etale open subsets of this kind form a base of the small
\'etale site of \(\tilde{\cA}_X\). This base allows us to define \'etale
sheaves \(\Pi\) and \(\Pi'\) as follows: for \(U_1\) of type \((I_1,W_1)\), let
\begin{align}
    \ShSec(U_1,\Pi')&=(\dual{\bT}^{W_1})^*=\CoCharG(\bT)_{W_1},\\
    \ShSec(U_1,\Pi)&=\dual{\bT}(I_1,W_1)^*,
\end{align}
where \(\dual{\bT}(I_1,W_1)\) is as in
\Cref{prop:global_pi_0_description}. Suppose \(U_2\) is an \'etale
open subset of \(U_1\) of type \((I_2,W_2)\), since \(U_{(I_1,W_1)}\) is the
unique non-empty closed stratum in \(U_1\), we see that
\begin{align}
    (I_1,W_1)\le (I_2,W_2).
\end{align}
The following lemma is straightforward.
\begin{lemma}[\cite{Ng10}*{Lemme~5.5.3}]
    \label[lemma]{lem:T_I_W_decreasing}
    If \((I_1,W_1)\le (I_2,W_2)\), then
    \(\dual{\bT}(I_2,W_2)\subset\dual{\bT}(I_1,W_1)\).
\end{lemma}
Thus, we have canonical restriction maps by dualization
\begin{align}
    \ShSec(U_1,\Pi')&\longto \ShSec(U_2,\Pi'),\\
    \ShSec(U_1,\Pi)&\longto \ShSec(U_2,\Pi).
\end{align}
Thus, the sheaves \(\Pi\) and \(\Pi'\) are defined, and we immediately have
the following by \Cref{prop:global_pi_0_description} and
\Cref{prop:pi_0_Picard_as_lattice_quotient}.
\begin{corollary}
    \label[corollary]{cor:pi_0_Picard_description}
    We have canonical isomorphisms of sheaves over \(\tilde{\cA}_X\)
    \begin{align}
        \pi_0(\cP_X')&\simeq \Pi',\\
        \pi_0(\cP_X)&\simeq \Pi.
    \end{align}
\end{corollary}

\section{Stratification by \texorpdfstring{\(\delta\)}{δ}-invariant}%
\label{sec:stratification_by_delta_invariant}

Recall that for any \(a\in\cA_X^\heartsuit(\bar{k})\) we have a global
\(\delta\)-invariant \(\delta_a\) by \Cref{def:global_delta_def}.
\begin{lemma}
    The function
    \begin{align}
        \tilde{\cA}_X(\bar{k})&\longto \bbN\\
        a&\longmapsto \delta_a
    \end{align}
    is constant on every stratum \(\tilde{\cA}_\psi\) (\(\psi\in\Psi\)).
\end{lemma}
\begin{proof}
    Base change to \(S=(\tilde{\cA}_\psi^\flat)^\Red\), then the cameral curve
    \(\pi_\psi\colon\tilde{X}_\psi\to \breve{X}\x S\) 
    admits a simultaneous normalization
    \(\xi\colon\tilde{X}_\psi^\flat\to \tilde{X}_\psi\). Let \(p_S\) be the
    projection \(\breve{X}\x S\to S\). By
    \Cref{lem:pi_0_of_simultaneous_normalization}, the sheaf
    \(F=\pi_{\psi,*}(\xi_*\cO_{\tilde{X}_\psi^\flat}/\cO_{\tilde{X}_\psi})\) is a
    coherent sheaf with \(W\)-action such that \(p_{S,*}F\) is a locally free
    \(\cO_S\)-sheaf of finite type. Since \(\Char(k)\) does not
    divide the order of \(\bW\), the same is true for
    \begin{align}
        (F\otimes_{\cO_{\breve{X}\x S}}\La{t})^W.
    \end{align}
    Therefore, the result follows from
    \Cref{cor:global_delta_formula_by_Neron}.
\end{proof}

Recall also that we have a rigidification of \(\cP_X\) over
\(\cA_X^\infty\) hence also over \(\tilde{\cA}_X\) in
\ref{sub:rigidification_of_global_Picard}. Combined with the following lemma, we
see that \(\delta\)-invariant is upper semi-continuous.
\begin{lemma}[\cite{Ng10}*{Lemme~5.6.3}]
    Let \(P\to S\) is a smooth commutative group scheme of finite type. The
    function \(s\mapsto\tau_s\) sending a geometric point \(s\in S\) to the
    dimension of the abelian part of \(P_s\) is lower semi-continuous.
    Equivalently, the function \(s\mapsto\delta_s\) sending \(s\) to the
    dimension of the affine part of \(P_s\) is upper semi-continuous.
\end{lemma}

\begin{corollary}
    For any \(\delta\in\bbN\), let \(\tilde{\cA}_{\delta}\)
    (resp.~\(\tilde{\cA}_{\ge\delta}\), resp.~\(\tilde{\cA}_{\le\delta}\)) be the union of all
    strata \(\tilde{\cA}_\psi\) such that \(\delta_\psi=\delta\)
    (resp.~\(\delta_\psi\ge\delta\), resp.~\(\delta_\psi\le\delta\)). Then it is
    a locally closed (resp.~closed, resp.~open) subset of \(\tilde{\cA}_X\).
    In particular, we have a stratification by \(\delta\)-invariant
    \begin{align}
        \tilde{\cA}=\coprod_{\delta\in\bbN}\tilde{\cA}_{\delta}.
    \end{align}
\end{corollary}
The stratification by \(\delta\)-invariant induces a stratification on
\(\tilde{\cA}_X^\ANI\), such that for any
\(\tilde{a}\in\tilde{\cA}_\delta^\ANI(\bar{k})\), the dimension of the
affine part of \(\cP_a\) is \(\delta\). It is also clear from
\Cref{cor:global_delta_formula_by_Neron} that
\(\delta\)-invariant is independent of the choice of points \(\infty\) and
\(\tilde{\infty}\), thus the stratification by \(\delta\)-invariant descends to
a stratification on
\(\cA_X^\heartsuit\):\index{stratum!\(\delta\)-}\index{\(\delta\)-!stratification}
\begin{align}
    \cA_X^\heartsuit = \coprod_{\delta\in\bbN}\cA_\delta.
    \nomenclature[\(A"cal_delta\)]{\(\cA_\delta\)}{a stratum of \(\cA_X^\heartsuit\) with a given \(\delta\)-invariant}
\end{align}

\subsection{}
The important result about \(\delta\)-strata is its codimension in
\(\cA_X\). We want to prove that the codimension of \(\cA_\delta\) is at least
\(\delta\) for all \(\delta\), and we will use
the local-global argument in \cite{Ng10}*{\S~5.7} to prove the result we need.
The proof in
\textit{loc.~cit.} uses root valuation strata studied in \cite{GKM09}.
Its multiplicative counterpart has been studied in
\Cref{chap:multiplicative_valuation_strata} (albeit less thoroughly).

\begin{proposition}
    \label[proposition]{prop:delta_regularity_fix_divisor}
    For any \(\delta\in\bbN\), there exists an integer \(N=N(\delta)\) depending on \(G\)
    and \(\delta\) such that if \(b\in\cB_X(\bar{k})\) is very \((G,N)\)-ample,
    then we have
    \begin{align}
        \codim_{\cA_b}\cA_{b,\delta}\ge \delta.
    \end{align}
    In particular, for every irreducible component \(A\subset
    \cA_{X}^\heartsuit\) that is very \((G,N)\)-ample, we
    have \(\codim_{A}A_\delta\ge\delta\).
\end{proposition}
\begin{proof}
    The \(\delta=0\) case is trivial. Suppose \(\delta>0\) and let
    \(\delta_\bullet\) be a partition of \(\delta\) by positive integers
    \begin{align}
        \delta=\delta_1+\cdots+\delta_n.
    \end{align}
    Consider subset \(Z_{\delta_\bullet}\subset \cA_b^\heartsuit\x \breve{X}^n\)
    consisting of tuples \((a;\bar{v}_1,\ldots,\bar{v}_n)\) such that
    \(a\in\cA_b^\heartsuit(\bar{k})\), \(\bar{v}_i\in \breve{X}(\bar{k})\) such that
    \(\delta_{\bar{v}_i}(a)=\delta_i\). Refine \(Z_{\delta_\bullet}\) into a disjoint
    union of subsets \(Z_{[w_\bullet,\bar{\lambda}_\bullet/l_\bullet,r_\bullet]}\)
    consisting of points \((a;v_1,\ldots, v_n)\) such that the image of \(a\) in
    \(\FRC_{\FRM,b}(\breve{\cO}_{\bar{v}_i})\) lies in the stratum
    \(\FRC_{\FRM,b}(\breve{\cO}_{\bar{v}_i})_{[w_i,\bar{\lambda}_i/l_i,r_i]}\). Note that \(G\) is
    split over \(\breve{\cO}_{\bar{v}_i}\), so the valuation strata make
    sense. Also note that \(\bar{\lambda}_i/l_i\) is actually fixed since \(b\)
    is. Suppose this
    stratum is an \(N_i\)-admissible cylinder. Let \(N'=N_1+\cdots+N_n\) and
    suppose \(b\) is very \((G,N')\)-ample. Then the linear map
    \begin{align}
        \cA_b\longto \prod_{i=1}^n\FRC_{\FRM,b}(\breve{\cO}_{\bar{v}_i}/\pi_{\bar{v}_i}^{N_i}\breve{\cO}_{\bar{v}_i})
    \end{align}
    is surjective. Using a modified proof of \Cref{thm:codim_valuation_strata}
    by fixing the boundary divisor, we see that
    the codimension of
    \(Z_{[w_\bullet,\bar{\lambda}_\bullet/l_\bullet,r_\bullet]}\) in \(\cA_b\x
    \breve{X}^n\) is
    \begin{align}
        \sum_{i=1}^n(\delta_i+c_i+e_i).
    \end{align}
    Here there is no \(b_i\)-term because the boundary divisor \(b\) is fixed.
    At \(\bar{v}_i\), if \(\bar{\lambda}_i/l_i=0\), then since
    \(\delta_i\neq 0\), we must have either \(c_i>0\) or \(e_i>0\). If
    \(\bar{\lambda}_i/l_i\neq 0\), then \(\bar{v}_i\) has to be one of the finite
    points in the support of \(b\), so it cannot move freely in \(\breve{X}\).
    Thus, the codimension of
    \(Z_{[w_\bullet,\bar{\lambda}_\bullet/l_\bullet,r_\bullet]}\) in \(\cA_b\x
    \breve{X}^n\) is always at least
    \begin{align}
        \sum_{i=1}^n(\delta_i+1)=\delta+n,
    \end{align}
    hence its image in \(\cA_b\) is at least \(\delta\). Let \(N\) be the
    supremum of all \(N'\) for various
    \(Z_{[w_\bullet,\bar{\lambda}_\bullet/l_\bullet,r_\bullet]}\) (there are
    only finitely many for a fixed \(\delta\)) and we are done.
\end{proof}

\begin{remark}
    \label[remark]{rmk:unconditional_delta_regularity}
    The inequality in \Cref{prop:delta_regularity_fix_divisor} is called
    \notion{\(\delta\)-regularity}\index{\(\delta\)-!regularity}. We note that its dependence on very
    \((G,N(\delta_a))\)-ampleness cannot be eliminated. However, we think that
    over the entire mH-base (in other words, not fixing \(b\)),
    \(\delta\)-regularity should be unconditional on \(\delta_a\). For
    example, if one can show:
    \begin{enumerate}
        \item for any \(a\in\cA_X^\heartsuit\), there exists some mHiggs bundle
            \(x\) lying over \(a\) such that the smoothened centralizer (see
            \Cref{sec:automorphism_group}) \(I_{x}^\Sm=\FRJ_a^\flat\), and
        \item the collection of such \(x\) is a constructible (likely closed)
            substack of the total stack,
    \end{enumerate}
    then the deformation results in
    \Cref{sec:deformation_of_mHiggs_bundles} would imply that
    \(\delta\)-regularity holds as long as the local model of
    singularity holds. The same argument can also be used to show that
    \(\delta\)-regularity holds unconditionally in \cite{Ng10}, and likely
    can be extended to many Hitchin-type fibrations in general. The
    \(I_{x}^\Sm=\FRJ_a^\flat\) condition can be checked locally, and if \(a\)
    is unramified at \(\bar{v}\in X(\bar{k})\) for example, such a condition is
    satisfied when \(x_{\bar{v}}\) is semisimple. Nevertheless,
    this plan has not been carried out and we leave it to a future work.
\end{remark}

\subsection{}
Let \(A\subset\cA_X\) be an irreducible component. If \(A\) is very \(G\)-ample, 
then its image in \(\BD_X\)
is an irreducible component \(B\) of \(\BD_X\). The normalization of \(B\) is
isomorphic to certain direct product of symmetric power of curves and there is
an open dense subset \(B^\circ\subset B\)  of multiplicity-free locus (see
\Cref{def:mult_free_boundary_divisor}).
We are particularly interested in the \notion{\(\delta\)-critical} subsets of
\(A\) defined as follows: 

\begin{definition}
    \label[definition]{def:delta_critical_strata}
    An irreducible subset in \(\cA_X^\heartsuit\) is called
    \notion{\(\delta\)-critical}\index{\(\delta\)-!critical}\index{stratum!\(\delta\)-critical} if its codimension in \(\cA_X\) equals
    its minimal \(\delta\)-invariant. We denote the union
    of \(\delta\)-critical irreducible components in \(\cA_\delta\) by \(\cA_\delta^=\).
    \nomenclature[\(A"cal_delta^=\)]{\(\cA_\delta^=\)}{the components of
    \(\cA_\delta\) whose codimension in \(\cA_X\) is equal to \(\delta\)}
\end{definition}

\begin{corollary}
    \label[corollary]{cor:when_delta_regularity_equality_is_met}
    Suppose we have a fixed \(G\), \(\FRM\) and \(\delta\in\bbN\). Let
    \(N=N(\delta)\) be as in \Cref{prop:delta_regularity_fix_divisor}.
    Then for every irreducible component \(Z\subset\cA_{\delta}\)
    that is very \((G,N)\)-ample, the
    equality \(\codim_{\cA_X}Z=\delta\) is achieved if and only if either
    \(\delta=0\), or \(\delta>0\) and there exists some \(a\in Z\) such that
    \begin{enumerate}
        \item \(a\) has multiplicity-free boundary divisor, in other words, it
            lies over a point \(b\in B^\circ\), and
        \item for every \(\bar{v}\in \FRD_a\), one of the followings holds:
            \begin{enumerate}
                \item \(\lambda_{\bar{v}}=0\) and \(c_{\bar{v}}+e_{\bar{v}}=1\), or
                \item \(\lambda_{\bar{v}}\neq 0\) and \(c_{\bar{v}}=e_{\bar{v}}=0\). 
            \end{enumerate}
    \end{enumerate}
\end{corollary}
\begin{proof}
    Straightforward from the proof of \Cref{prop:delta_regularity_fix_divisor}.
\end{proof}

\subsection{}
\label{sub:delta_critical_local_delta}
The numerical conditions in \Cref{cor:when_delta_regularity_equality_is_met}
implies that when \(\delta>0\) is not too large, in view of product formula
\eqref{eqn:product_formula}, the multiplicative affine Springer fibers
\(\cM_{\bar{v}}(a)\) has to be sufficiently simple, which we
summarize below.

In case there is no boundary divisor at \(\bar{v}\), in other words,
\(\lambda_{\bar{v}}=0\), we have two possibilities. The first possibility is
\(c_{\bar{v}}=1\) and \(e_{\bar{v}}=0\). In this case the ramification happens
in a Levi subgroup generated by a single pair of roots and no other root has
any contribution to \(\delta\), and we must have \(d_{\bar{v}+}(a)=1\) and
\(\delta_{\bar{v}}(a)=0\). The second possibility is \(c_{\bar{v}}=0\) and
\(e_{\bar{v}}=1\). It implies that \(a\) is unramified at \(\bar{v}\) and the
contribution to \(\delta\) still only comes from a single pair of roots, so we
must have \(d_{\bar{v}+}(a)=2\) and \(\delta_{\bar{v}}(a)=1\). In both
possibilities the
local computations can be reduced to the case where \(\bG^\SC=\SL_2\).

On the other hand, when \(\lambda_{\bar{v}}\neq 0\), then
\(a\) must be unramified at \(\bar{v}\), and
suppose \(\nu_{\bar{v}}\) is the local Newton point of \(a\) at \(\bar{v}\) and
\(\nu_{\bar{v}}\) is chosen to be dominant, then since
\(c_{\bar{v}}=e_{\bar{v}}=0\), we have
\begin{align}
    \label{eqn:Newton_rss_equality}
    \delta_{\bar{v}}(a)&=\Pair{\rho}{\lambda_{\bar{v},\AD}-\nu_{\bar{v},\AD}}\\
                       &=\Pair{\rho}{-w_0(\lambda_{\bar{v},\AD})-\nu_{\bar{v},\AD}}.
\end{align}
In other words, \(a\) must be \(\nu\)-regular semisimple
(cf.,~\Cref{def:newton_rss}) when restricted to \(\breve{X}_{\bar{v}}\).

\section{Inductive Subsets}%
\label{sec:inductive_strata}

In the Lie algebra case, \(\delta\)-critical subsets play an important role
in Ng\^o's support theorem for \(\delta\)-regular abelian
fibrations (see \cite{Ng10}*{\S~7.2}). For Hitchin fibrations specifically,
only two types of \(\delta\)-critical subsets turn out to actually ``matter'':
the \(\delta=0\) stratum, and strata that come from endoscopy
(see \cite{Ng10}*{\S\S~6.2--6.3} and \Cref{sec:kappa_decomposition_and_endoscopic_embedding}).

In the multiplicative case, the story is similar, but there is a third kind of
\(\delta\)-critical subsets to consider due to the influence of
boundary divisors. We shall call these subsets
\notion{inductive}\index{subset!inductive}\index{stratum!inductive}, whose meaning
will be apparent soon.

\subsection{}
We first use split group \(G=\bG^\SC\) and base change to \(\bar{k}\) to illustrate
the idea. Let \(\lambda\in\CoCharG(\bT)_+\) be a dominant cocharacter and
\(d\ge 0\) a sufficiently large integer.
Then we may have a moduli space of boundary divisors \(\breve{X}_d\), viewed as the
system of divisors \(\sum_{\bar{v}\in X(\bar{k})}d_{\bar{v}}\lambda\cdot \bar{v}\) such that
\(\sum_{\bar{v}\in X(\bar{k})}d_{\bar{v}}=d\). We have finite \(\FRS_d\)-cover
\(\breve{X}^d\to
\breve{X}_d\), where \(\FRS_d\) is the symmetric group of \(d\) elements.
The fiber product \(\cA_{d,\lambda}\x_{\breve{X}_d}\breve{X}^d\to\breve{X}^d\) is a
vector bundle whose fiber at \(D=(\bar{v}_1,\ldots,\bar{v}_d)\in\breve{X}^d\) is
\begin{align}
    \cA_{d,\lambda,D}\defeq \bigoplus_{j=1}^r\RH^0
    (\breve{X},\cO_{\breve{X}}(\Pair{\Wt_j}{\lambda}D)),
\end{align}
viewed as a vector space.

For each \(1\le i\le d\), let \(\mu_i\in\CoCharG(\bT)_+\) be another
dominant cocharacter with \(\mu_i\le \lambda\). Then for each \(1\le j\le r\)
we have natural inclusion of divisors on \(\breve{X}\):
\begin{align}
    D_{\mu,j}\defeq\sum_{i=1}^d\Pair{\Wt_j}{\mu_i}\bar{v}_i\subset
    \Pair{\Wt_j}{\lambda}D.
\end{align}
Suppose \(D_{\mu,j}\) is very ample for all \(j\) (or
equivalently, \(\sum_i\mu_i\cdot \bar{v}_i\) is very \(G\)-ample), then
we have an inclusion of vector bundles
\begin{align}
    \cA'\subset \cA_{d,\lambda}\x_{\breve{X}_d}\breve{X}^d,
\end{align}
where \(\cA'\to\breve{X}^d\) is the vector bundle whose fiber over \(D\) is
\begin{align}
    \bigoplus_{j=1}^r\RH^0(\breve{X},\cO_{\breve{X}}(D_{\mu,j})).
\end{align}
Since \(\breve{X}^d\to\breve{X}_d\) is finite, the image of \(\cA'\) in
\(\cA_{d,\lambda}\) is closed, and it is \(\delta\)-critical in
\(\cA_{d,\lambda}\).

If we fix \(D\subset \breve{X}_d\), then the set-theoretic
image of \(\cA'_D\) in the fiber \(\cA_{d,\lambda,D}\) is a union of linear
subspaces, each of which is isomorphic to a \emph{restricted} mH-base associated with divisor
\(\sum_i\mu_i\cdot\bar{v}_i\) for some \(\mu_i\le \lambda\). This way we
have an inductive system on restricted mH-bases (indexed by boundary
divisors), and for this reason we call any non-empty open subset of the image of
\(\cA'\) an \notion{inductive subset}\index{subset!inductive}\index{stratum!inductive} in \(\cA_{d,\lambda}\).

\subsection{}
Now we formulate for arbitrary quasi-split group \(G\).
Let \(\FRM\in\FM(G^\SC)\) and let \(B\) be an
irreducible component of the very \(G\)-ample locus \(\cB_{\gg}\subset\cB_X\).
There is no harm in assuming that \(\FRA_\FRM\) is of standard type.
There is a unique open subset \(B^\circ\) being the multiplicity-free locus.

A point \(b\in B^\circ(\bar{k})\) may be written as a tuple \((\cL,\lambda_b)\)
where \(\cL\in\Bun_{Z_{\FRM}}\), and
\begin{align}
    \lambda_b=\sum_{i=1}^{d}\lambda_{i}\cdot \bar{v}_{i}
\end{align}
for some closed points \(\bar{v}_{i}\in\breve{X}(\bar{k})\) and 
dominant cocharacters \(\lambda_{i}\in\CoCharG(T^\AD_{\bar{v}_{i}})_+\). The
points \(\bar{v}_{i}\) are pairwise distinct, and if we identify
\(\FRM\) with its split model \(\bM\) at \(\bar{v}_i\), then \(\lambda_i\) must
be one of the minimal generators of the cone of \(\bA_\bM\) because \(b\) is
multiplicity-free.
Let \(\mu_{i}\in\CoCharG(T^\AD_{\bar{v}_{i}})_+\) be a dominant cocharacter with
\(\mu_{i}\le\lambda_{i}\).

Let \(\OGT\colon X_\OGT\to \breve{X}\) be a connected finite \'etale cover
over which \(G\) becomes split. Let \(\cL_\OGT\) be the pullback
of \(\cL\) to \(X_\OGT\). For each fundamental weight \(\Wt_j\), the
divisor
\begin{align}
    D_{j}\defeq \sum_{i}\Pair{\Wt_j}{\lambda_{i}-\mu_{i}}\cdot
    \OGT^*\bar{v}_{i}
\end{align}
defines an inclusion of vector bundles on \(X_\OGT\):
\begin{align}
    \bigoplus_{j=1}^r\Wt_j(\cL_\OGT)(-D_j)\subset
    \bigoplus_{j=1}^r\Wt_j(\cL_\OGT),
\end{align}
which descends to an inclusion of vector bundles on \(\breve{X}\)
\begin{align}
    \FRC_\cL''\subset \FRC_\cL.
\end{align}
Suppose each \(\Wt_j(\cL_\OGT)(-D_j)\) is still very ample over \(X_\OGT\), and
the degree of \(\OGT\) is invertible in \(k\), then
the section space
\begin{align}
    \cA_b''=\RH^0(\breve{X},\FRC_\cL'')
\end{align}
is a linear subspace of \(\cA_b\) with codimension
\(\delta=\sum_{i}\Pair{\rho}{\lambda_{i}-\mu_{i}}\), and is \(\delta\)-critical
in \(\cA_b\).

Let \(b\in B^\circ\) vary while keeping \(\mu_{i}\) locally constant,
and let \(\cA'\) be the resulting union with reduced structure.
Then it is not hard to see that \(\cA'\) is a locally-trivial fibration over
\(B^\circ\) whose fiber \(\cA_b'\) at \(b\) is a union of subspaces
looking like \(\cA_b''\) (in fact, it is just the union of \(\cA_b''\) and all
its conjugates under monodromy). It is
clearly \(\delta\)-critical in \(\cA_X\).

\begin{definition}
    \label[definition]{def:inductive_subset}
    A non-empty constructible subset of \(\cA'\) containing an open subset therein is called
    an \notion{inductive subset}\index{subset!inductive}\index{stratum!inductive} in \(\cA_X\). We denote by
    \(\cA_{\delta}^\equiv\)
    \nomenclature[\(A"cal_delta^equiv\)]{\(\cA_\delta^\equiv\)}{the components of \(\cA_\delta\) that are inductive}
    the union of irreducible components of
    \(\cA_\delta\) that are inductive.
\end{definition}

\subsection{}
We want to point out the key difference between being \(\delta\)-critical
and being inductive, assuming appropriate \(G\)-ampleness conditions are met. Let
\(a\) be a sufficiently general point in an inductive subset \(Z\), then for any
\(\bar{v}\in X(\bar{k})\) where \(\lambda_{\bar{v}}=0\), one must have
\(\delta_{\bar{v}}(a)=0\) and either \(d_{\bar{v}+}(a)=0\)
or \(d_{\bar{v}+}(a)=1\). On the other hand, if
\(Z\) is only \(\delta\)-critical, there is a third
possibility where \(d_{\bar{v}+}(a)=2\) and \(\delta_{\bar{v}}(a)=1\).
Therefore, we have strict inclusion
\begin{align}
    \cA_\delta^\equiv\subsetneq \cA_\delta^=.
\end{align}

\begin{proposition}
    \label[proposition]{prop:delta_critical_in_inductive}
    For any \(\delta\in\bbN\) and \(N=N(\delta)\) as in
    \Cref{prop:delta_regularity_fix_divisor}, let \(A'\subset\cA_{\gg N}\) be an
    irreducible locally closed \(\delta\)-critical subset. Then we can find an
    inductive subset \(A\) containing \(A'\) as a closed subset such that
    for any sufficiently general point \(a'\in A_b'\) (resp.~\(a\in A_b\)) with
    multiplicity-free boundary divisor \(b\), we have
    \begin{align}
        \delta_{\bar{v}}(a)=\delta_{\bar{v}}(a')
    \end{align}
    for any \(\bar{v}\in\supp(b)\).
\end{proposition}
\begin{proof}
    The proof is similar to that of
    \Cref{prop:delta_regularity_fix_divisor}: we may deform locally near every point
    \(\bar{v}\in\FRD_{a'}-\supp(b)\) so that \(\breve{X}\) intersects with
    \(\FRD_\FRM\) transversally near those points. With the assumption on
    \(G\)-ampleness, this extends globally to a point \(a\) that lies in an
    inductive subset \(A_b\).
\end{proof}

\chapter{Cohomologies}%
\label{chap:cohomologies}

In this chapter we study the general properties of cohomologies over the
anisotropic locus. The first two sections contain results similar to those in
\cite{Ng10}*{\S\S~6.1--6.3}. The main difference is that the global
transfer map will no longer be a closed embedding in general but only a finite
map.

Such a difference is completely explained
by representations of the \(L\)-groups, which will be manifested in the ensuing
sections
where we explain how to transfer a global Satake sheaf from \(G\) to
its endoscopic group \(H\). We then use it to give a statement of geometric
stabilization which we will prove at the end of this book, together with the
fundamental lemma.

The last part of this chapter studies top ordinary cohomologies, which is much
more complicated than the Lie algebra case. To this end, we will introduce a new
type of Hecke stack which allows us to upgrade the product formula in
\Cref{sec:product_formula} into a family.

\section{Properness over Anisotropic Locus}%
\label{sec:properness_over_anisotropic_locus}

So far we have studied the properties of individual mH-fibers \(\cM_a\) as
well as the Picard stack \(\cP_X\to\cA_X\). In this section we turn to the total
space \(\cM_X\) of mH-fibration. The first result is its finiteness properties.

\begin{proposition}
    \label[proposition]{prop:Mnatural_DM_finite_type_over_anisotropic_locus}
    The stack \(\cM_X\) is locally of finite type, and 
    \(h_X^\ANI\colon\cM_X^\ANI\to \cA_X^\ANI\) over the
    anisotropic locus is a relative Deligne--Mumford stack of finite
    type. The same is true for \(\cP_X\).
\end{proposition}
\begin{proof}
    The natural map \(\cM_X\to\Bun_G\x\Bun_{Z_{\FRM}}\) 
    is of finite type because the fiber over \((\cE,\cL)\in\Bun_G\x\Bun_{Z_{\FRM}}\)
    is the section space of the \'etale-locally trivial fiber bundle
    \((\cE\x\cL)\x^{G\x Z_\FRM}\FRM\) with affine fibers isomorphic to \(\FRM\).
    Since \(X\) is projective, such section space must be of finite
    type. Since \(\Bun_G\x\Bun_{Z_{\FRM}}\) is locally of finite type, so is
    \(\cM_X\).

    For the second claim, it suffices to prove for \(\FRM=\Env(G^\SC)\).
    The case when \(G\) is split is proved by
    \cite{Ch22}*{Proposition~4.3.3}. Although the statement in
    \textit{loc.~cit.} is about mH-fibration over a fixed \(Z_{\FRM}\)-torsor
    \(\cL\), the role of \(\cL\) is inconsequential since our claim is about
    relative finiteness of map \(h_X^\ANI\).

    When \(G\) is non-split, let \(\OGT\colon X_{\OGT}\to X\) be a connected
    finite Galois \'etale cover over which \(G\) becomes split,
    and \(\Theta\)
    be the Galois group. Then we have a map by pullback
    \begin{align}
        \OGT^*\colon \cM_X\longto \cM_{X_{\OGT}},
    \end{align}
    the latter being the total space of mH-fibration associated with split
    monoid \(\OGT^*\FRM\).
    Since \(\OGT\) is \'etale, \(G\to G_\OGT\defeq\OGT_*\OGT^*G\) is a closed embedding of
    \'etale-locally constant reductive group schemes over \(X\), hence \(G_\OGT/G\) is affine over
    \(X\) by \cite{Ric77}*{Theorem~A}. It is straightforward to see that a \(G\)-bundle may be identified
    with a \(G_\OGT\)-bundle \(E\) together with a section of the associated bundle
    \(E\x^{G_\OGT}G_\OGT/G\). Since \(G_\OGT/G\) is affine over \(X\) and \(X\)
    is projective, the section space is representable of finite type, hence so
    is the map
    \begin{align}
        \Bun_{G/X}\longto \Bun_{G/X_\OGT}.
    \end{align}
    This implies that \(\cM_X\to\cM_{X_\OGT}\)
    is of finite type. The proof for \(\cP_X\) is the same hence we are done.
\end{proof}

\begin{proposition}
    \label[proposition]{prop:properness_over_anisotropic_locus}
    The map \(h_X^\ANI\colon\cM_X^\ANI\to\cA_X^\ANI\) is proper.
\end{proposition}
\begin{proof}
    Since the map is a morphism of Deligne--Mumford stacks and is of finite
    type, we use valuative criteria as in \cite{CL10}*{\S\S~8--9} and 
    \cite{Ch22}*{Proposition~4.3.6}. Let
    \(R\) be a discrete valuation ring and \(K\) its fractional field. 

    For the existence part of valuative criteria, it is
    harmless to assume that the residue field \(k_R\) of \(R\) is algebraically
    closed, because we are allowed to take finite extensions of \(R\).
    Let \((E,\phi)\in\cM_X^\heartsuit(K)\) over \(a\in\cA_X^\heartsuit(R)\). Let
    \(\cL\) be the \(Z_{\FRM}\)-torsor corresponding to \(a\). Let \(R(X)\) be
    the local ring of the generic point of \(X\x\Spec{k_R}\) inside \(X\x\Spec{R}\)
    and \(K(X)\) its fractional field. Similar notation is used for any
    algebraic extension \(K'/K\) and \(R'/R\) where \(R'\) is the integral
    closure of \(R\) in \(K'\).

    Since \(a\in\cA_X^\heartsuit\), the restriction of \(a\) to
    \(R'(X)\) has image contained in \(\FRM^\rss_\cL\), thus \(\FRJ_a\) is a torus
    when restricted to \(R(X)\). Since \(\Stack{\FRM^\rss/G}\to\FRC_{\FRM}\) is a
    gerbe bound by \(\FRJ\), any trivialization of this gerbe over \(k_R(X)\)
    (necessarily exists as \(k_R\) is algebraically closed)
    gives a trivial \(G\)-torsor \(E_0\) together with a \(G\)-equivariant map
    \(\phi_0\colon E_0\to\FRM_\cL\) over \(k_R(X)\). Since
    \(\FRM^\rss\to\FRC_{\FRM}\) is smooth, we can extend \((E_0,\phi_0)\) to a
    pair \((E_1,\phi_1)\) over \(R(X)\) by formal lifting property of
    smoothness, where \(E_1\) is a trivial \(G\)-torsor over \(R(X)\).

    The transporter between \(\phi\) and \(\phi_1\) over \(K(X)\) is a
    \(\FRJ_a\)-torsor, which can be trivialized after passing to a finite
    extension \(K'/K\). This means that \((E,\phi)\) and \((E_1,\phi_1)\) can
    be glued into a pair \((E',\phi')\) over an open subset of \(X\x\Spec{R'}\)
    whose complement has codimension at least \(2\). Since any \(G\)-torsor can
    be extended over a subset of codimension at least \(2\), and since
    \(\FRM_\cL\) is affine over \(X\x\Spec{R'}\), the pair \((E',\phi')\) extends
    to a point in \(\cM_X(R')\) lying over \(a\). This proves the existence part
    of valuative criteria.

    Now for the uniqueness part. Suppose
    \((E,\phi),(E',\phi')\in\cM_X^\ANI(R)\) be such that their restriction
    to \(K\) are isomorphic. Let \(\iota_K\) be such isomorphism, then using
    codimension-\(2\) argument again it
    suffices to extend \(\iota_K\) to \(R'(X)\) for some finite extension
    \(R'/R\), because \(K\cap R'=R\). Hence, it is still harmless to assume
    \(k_R\) is algebraically closed. Therefore, \(E\) and \(E'\) are both trivial
    over \(R(X)\).

    Moreover, as in the existence part, we may pass to a finite
    extension \(K'/K\) and carefully choose trivializations
    so that both \(\phi\) and \(\phi'\) map the neutral point
    of \(E_{K'(X)}\cong E_{K'(X)}'\) to some \(\gamma\in\FRM_\cL^\rss(K'(X))\). This means
    that \(\iota_{K'(X)}\) may be represented by some element in
    \(G_\gamma(K'(X))\). Since \((E,\phi)\) and \((E',\phi')\) are contained in
    the anisotropic locus, the centralizer \(G_\gamma\) is an anisotropic torus
    over \(R'(X)\). Since \(R'(X)\) is a discrete valuation ring with fractional
    field \(K'(X)\), we have \(G_\gamma(R'(X))=G_\gamma(K'(X))\), and we are
    done.
\end{proof}

\section{\texorpdfstring{\(\kappa\)}{κ}-decomposition and Endoscopic Transfer}%
\label{sec:kappa_decomposition_and_endoscopic_embedding}
Let \(\tilde{h}_{X}^\ANI\colon \tilde{\cM}_{X}^\ANI\to
\tilde{\cA}_{X}^\ANI\) and
\(\tilde{p}_{X}^\ANI\colon\tilde{\cP}_{X}^\ANI\to\tilde{\cA}_{X}^\ANI\) be
the restriction of mH-fibration \(h_X\colon\cM_X\to\cA_X\) and Picard stack
\(p_X\colon\cP_X\to\cA_X\) to
\(\tilde{\cA}_{X}^\ANI\) respectively.
We know that both
\(\tilde{h}_X^\ANI\) and \(\tilde{p}_X^\ANI\) are morphisms of
Deligne--Mumford stacks with \(\tilde{h}_X^\ANI\) being proper and
\(\tilde{p}_X^\ANI\) being smooth.
Let \(\bcQ\in\Sat_X\) be a Satake sheaf
(e.g., the intersection complex on \(\tilde{\cM}_X^\ANI\)), then we know
that \(\tilde{h}_{X,*}^\ANI\bcQ\) is a pure complex, hence non-canonically
decomposes into a direct sum of shifted perverse sheaves \cite{Sun12b}:
\begin{align}
    \tilde{h}_{X,*}^\ANI\bcQ\cong\bigoplus_{n\in\bbZ}\PH^n(\tilde{h}_{X,*}^\ANI\bcQ)[-n],
\end{align}
where \(\PH^n(\tilde{h}_{X,*}^\ANI\bcQ)\) is a perverse sheaf of weight \(n\)
over \(\tilde{\cA}_{X}^\ANI\).

\subsection{}
The action of \(\tilde{\cP}_X^\ANI\) on \(\tilde{\cM}_X^\ANI\) relative to
\(\tilde{\cA}_X^\ANI\) induces an action on \(\tilde{h}_{X,*}^\ANI\bcQ\).
As we will see in \Cref{sub:kappa_decomposition_at_chain_level}, this action
factors through \(\pi_0(\tilde{\cP}_{X}^\ANI)\). Over
\(\tilde{\cA}_{X}^\ANI\), we have by
\Cref{prop:pi_0_Picard_as_lattice_quotient} a
canonical epimorphism
\begin{align}
    \CoCharG(\bT)\x\tilde{\cA}_X^\ANI\longto\pi_0(\tilde{\cP}_X^\ANI),
\end{align}
so we have an inflated action of \(\CoCharG(\bT)\) on
\(\tilde{h}_{X,*}^\ANI\bcQ\). For any \(\kappa\in\dual{\bT}\), we define
\(\PH^n(\tilde{h}_{X,*}^\ANI\bcQ)_\kappa\) to be the \(\kappa\)-isotypic
subspace of \(\PH^n(\tilde{h}_{X,*}^\ANI\bcQ)\), where \(\kappa\) is regarded
as a character \(\CoCharG(\bT)\to \Qlb^\x\). Therefore, we have a (necessarily
finite) decomposition
\begin{align}
    \PH^n(\tilde{h}_{X,*}^\ANI\bcQ)
    =\bigoplus_{\kappa\in\dual{\bT}}\PH^n(\tilde{h}_{X,*}^\ANI\bcQ)_\kappa.
    \nomenclature[\(.kappa \)]{\((\cdot)_\kappa\)}{\(\kappa\)-isotypic subspace}
\end{align}
When \(\kappa=1\), we write \(\PH^n(\tilde{h}_{X,*}^\ANI\bcQ)_{\hST}\)
    \nomenclature[\(.st \)]{\((\cdot)_\hST\)}{\((\cdot)_\kappa\) when \(\kappa=1\)}
instead of \(\PH^n(\tilde{h}_{X,*}^\ANI\bcQ)_1\).

Recall we have the stratification \eqref{eqn:A_tilde_stratification_prelim}
over \(G\)-very ample locus induced by simultaneous normalizations of cameral
curves:
\begin{align}
    \tilde{\cA}_\gg=\coprod_{\psi\in\Psi}\tilde{\cA}_\psi.
\end{align}
Let \(\tilde{\cA}_\gg^\ANI=\tilde{\cA}_\gg\cap\tilde{\cA}_X^\ANI\).
By \Cref{lem:I_W_increasing,lem:T_I_W_decreasing}, the union of all \(\psi\)
such that \(\kappa\in\dual{\bT}(I_\psi,W_\psi)\) for a fixed \(\kappa\) is a
closed subset of \(\tilde{\cA}_\gg^\ANI\). Call this subset
\(\tilde{\cA}_\kappa^\ANI\).
\begin{proposition}
    The support of \(\PH^n(\tilde{h}_{\gg*}^\ANI\bcQ)_\kappa\) is contained in
    \(\tilde{\cA}_\kappa^\ANI\).
\end{proposition}
\begin{proof}
    This is a direct consequence of the definition of
    \(\PH^n(\tilde{h}_{\gg*}^\ANI\bcQ)_\kappa\) and
    \Cref{cor:pi_0_Picard_description}.
\end{proof}

\subsection{}
Now we turn to the endoscopic side. Let \((\kappa,\vartheta_\kappa^\bullet)\) is a
pointed endoscopic datum with endoscopic group \(H\) over \(\breve{X}\). Let
\(\vartheta_\kappa\colon X_\kappa\to \breve{X}\) be the corresponding
\(\pi_0(\kappa)\)-torsor. Recall we have finite unramified map
\begin{align}
    \nu_{\cA}^\heartsuit\colon \cA_{H,X}^{\kappa,G\hy\heartsuit}\longto
    \cA_X^\heartsuit,
\end{align}
by \Cref{prop:endoscopic_transfer_finite_and_unramified}.

The pointed endoscopic datum is given by a continuous homomorphism
\(\pi_1(\breve{X},\infty)\to\pi_0(\kappa)\). There is a natural point
\(\infty_{\vartheta_\kappa}\) lying over \(\infty\). Given
\(a\in\cA_X^\infty(\bar{k})\), let
\begin{align}
    \tilde{X}_{\vartheta_\kappa,a}=\tilde{X}_a\x_{\breve{X}}X_\kappa.
\end{align}
Choosing a point \(\tilde{\infty}\in\tilde{X}_a\) is the same as choosing 
a point \(\tilde{\infty}_{\vartheta_\kappa}=(\tilde{\infty},
\infty_{\vartheta_\kappa})\), and let
\((a,\tilde{\infty}_{\vartheta_\kappa})\in\tilde{\cA}_\gg\). Suppose
\(a_H\in\cA_{H,X}^\kappa(\bar{k})\) and \(\nu_\cA(a_H)=a\), then we have a
finite map \(\tilde{X}_{\vartheta_\kappa,a_H}\to\tilde{X}_{\vartheta_\kappa,a}\)
by construction, and if \(a\in\cA_X^\heartsuit\), it
birationally identifies \(\tilde{X}_{\vartheta_\kappa,a_H}\) with the union
of some irreducible components of \(\tilde{X}_{\vartheta_\kappa,a}\). Thus, we
have a finite unramified map
\begin{align}
    \tilde{\nu}_{\cA}\colon\tilde{\cA}_{H,X}^\kappa\longto \tilde{\cA}_X,
\end{align}
and it is defined over \(k\) if \((\vartheta_\kappa,\kappa)\) is.

Unlike Lie algebra case, \(\tilde{\nu}_\cA\) is in general not a closed
embedding. This roughly corresponds to the fact that an irreducible representation
of \(\dual{\bG}\) is usually not irreducible when restricted to
\(\dual{\bH}\). More precisely, we have the following characterization:

\begin{proposition}
    \label[proposition]{prop:pointed_endoscopic_transfer_closed_embedding}
    The map \(\tilde{\cA}_{H,X}^\kappa\) is finite and unramified.
    Moreover, suppose \(\bA_\bM\) is of standard type generated by dominant
    cocharacters \(-w_0(\theta_1),\ldots,-w_0(\theta_m)\), and suppose the
    \(\dual{\bG}\)-representation \(V_{\theta_i}\) with highest-weight \(\theta_i\) decomposes as
    \begin{align}
        V_{\theta_i}=\bigoplus_{j=1}^{e_i}V_{\theta_{H,ij}}^H
        \otimes\Hom_{\dual{\bH}}(V_{\theta_{H,ij}}^H,V_{\theta_i})
    \end{align}
    where \(V_{\theta_{H,ij}}^H\) is the \(\dual{\bH}\)-representation with
    highest-weight \(\theta_{H,ij}\). Then for any
    \(a\in\tilde{\nu}_\cA(\tilde{\cA}_{H,X}^\kappa)(\bar{k})\), we have the canonical
    isomorphism:
    \begin{align}
        \label{eqn:fiber_of_geometric_transfer_map}
        \tilde{\nu}_\cA^{-1}(a)\stackrel{\sim}{\longto}\prod_{\bar{v}\in
        \supp(b)}\Set*{(d_{\bar{v}ij})\given\nu_{\bar{v}}\le_{H,\bbQ}\sum_{i=1}^m\sum_{j=1}^{e_i}d_{\bar{v}ij}\theta_{H,ij}},
    \end{align}
    where 
    \begin{enumerate}
        \item \(b\) is the boundary divisor of \(a\),
        \item \(\lambda_{\bar{v}}\) is the dominant cocharacter of \(G_{\bar{v}}\)
            given by \(b\) which we write as
            \(-w_0(\lambda_{\bar{v}})=\sum_{i}c_{\bar{v}i}\theta_i\),
        \item \(\nu_{\bar{v}}\) is the local Newton point, viewed as a dominant
            \(\bbQ\)-cocharacter of \(G_{\bar{v}}\) such that
            \(\nu_{\bar{v}}\le_\bbQ\lambda_{\bar{v}}\),
        \item for each \(i\), \((d_{\bar{v}ij})\) are partitions of \(c_{\bar{v}i}\) into
            natural numbers.
    \end{enumerate}
\end{proposition}
\begin{proof}
    The first claim is directly from
    \Cref{prop:endoscopic_transfer_finite_and_unramified}, whose proof also
    shows that \eqref{eqn:fiber_of_geometric_transfer_map} is surjective. So we
    only need to show that if \(a_H\) and \(a_H'\) are two points in
    \(\tilde{\nu}_\cA^{-1}(a)\) with the same \(H\)-boundary divisor, then they
    are the same. The proof of this fact is similar to
    \cite{Ng06}*{Proposition~6.3.2}.

    Indeed, as in \Cref{prop:endoscopic_transfer_finite_and_unramified}, let \(U_a\subset \breve{X}\)
    be the maximal open subset whose image under \(a\) is contained in
    \(\Stack*{\FRC_{\FRM}^{\x,\rss}/Z_\FRM}\), then the restrictions of \(a_H\)
    and \(a_H'\) to \(U_a\) are the same because it is uniquely determined by
    the smallest \(\bW_\bH\rtimes\pi_0(\kappa)\)-stable subset of
    \(\tilde{X}_{\OGT_\kappa,a}\) containing \(\tilde{\infty}_{\OGT_\kappa}\).
    If \(a_H\) and \(a_H'\) have the same boundary divisor, then they must also
    be the same over all \(\breve{X}\), and we are done.
\end{proof}

\begin{proposition}
    \label[proposition]{prop:kappa_locus_is_from_endoscopy}
    Over the very \(G\)-ample locus, the subset
    \(\tilde{\cA}_\kappa\subset\tilde{\cA}_\gg\) is the disjoint union of
    various closed subsets \(\tilde{\nu}_\cA(\tilde{\cA}_{H,X}^\kappa)\), where
    \(H\) is the endoscopic group corresponding to a continuous homomorphism
    \(\OGT_\kappa^\bullet\colon\pi_1(\breve{X},\infty)\to\pi_0(\kappa)\).
\end{proposition}
\begin{proof}
    Given a geometric point
    \(\tilde{a}=(a,\tilde{\infty})\in\tilde{\cA}_X(\bar{k})\), recall we have
    diagram \eqref{eqn:I_W_diagram} which we reproduce here:
    \begin{equation}
        \label{eqn:I_W_diagram3}
        \begin{tikzcd}
            \pi_1(U,\infty) \ar[r, "\pi_{\tilde{a}}^\bullet"]\ar[d] 
                & \bW\rtimes\Out(\bG) \ar[d]\\
            \pi_1(\breve{X},\infty) \ar[r, "\OGT_G^\bullet"] & \Out(\bG)
        \end{tikzcd}
    \end{equation}
    Here \(U=\breve{X}-\FRD_a\). Let \(W_{\tilde{a}}\) be the image of
    \(\pi_{\tilde{a}}^\bullet\) in \(\bW\rtimes\Out(\bG)\) and \(I_{\tilde{a}}\)
    the image of the kernel of \(\pi_1(U,\infty)\to\pi_1(\breve{X},\infty)\).

    If \(\tilde{a}\in\tilde{\cA}_\kappa\), then \(W_{\tilde{a}}\subset
    (\bW\rtimes\Out(\bG))_\kappa\) and \(I_{\tilde{a}}\subset\bW_{\bH}\). Here we
    canonically identify \((\bW\rtimes\Out(\bG))_\kappa\) with
    \(\bW_{\bH}\rtimes\pi_0(\kappa)\) by \cite{Ng06}*{Lemme~10.1}. Then
    \(\pi_{\tilde{a}}^\bullet\) induces a unique homomorphism
    \begin{align}
        \OGT_\kappa^\bullet\colon\pi_1(\breve{X},\infty)\longto \pi_0(\kappa),
    \end{align}
    and let \(H\) be the corresponding endoscopic group. Let
    \(C_{\tilde{a}}^\kappa\subset\tilde{X}_{\OGT_\kappa,a}\) be the union of
    irreducible components in the \(\bW_{\bH}\rtimes\pi_0(\kappa)\)-orbit of the
    unique component containing \(\tilde{\infty}_{\OGT_\kappa}\), then
    \(\bW_{\bH}\rtimes\pi_0(\kappa)\) acts transitively on fibers of 
    \(C_{\tilde{a}}^\kappa\to\breve{X}\). This means that \(\tilde{a}\) comes
    from a map \(\breve{X}\to\Stack*{\FRC_{\FRM,H}'/Z_{\FRM}}\) (recall that
    over \(\breve{X}\), 
    \(\FRC_{\FRM,H}'=(\bar{\bT}_{\bM}\x X_\kappa)/(\bW_{\bH}\rtimes\pi_0(\kappa))\)
    by definition). The map
    \begin{align}
        \Stack*{\FRC_{\FRM,H}/Z_{\FRM}^\kappa}\longto\Stack*{\FRC_{\FRM,H}'/Z_{\FRM}}
    \end{align}
    is an isomorphism over the intersection of invertible and \(G\)-regular
    semisimple loci. The point \(\tilde{a}\), viewed as a map
    \(\breve{X}\to\Stack*{\FRC_{\FRM}/Z_{\FRM}}\) that is generically contained in the
    invertible and regular semisimple locus, can then be lifted to a rational
    map from \(\breve{X}\) to \(\Stack*{\FRC_{\FRM,H}/Z_{\FRM}^\kappa}\). By
    \Cref{lem:arc_lifting_to_endoscopic_quotient}, it can be extended to a
    morphism \(\breve{X}\to\Stack*{\FRC_{\FRM,H}/Z_{\FRM}^\kappa}\). This shows that
    every \(\tilde{a}\in\tilde{\cA}_\kappa(\bar{k})\) comes from some (not
    necessarily unique)
    \(\tilde{a}_H\in\tilde{\cA}_{H,X}^\kappa(\bar{k})\).
    The argument also shows that \(\OGT_\kappa^\bullet\) hence \(H\) is uniquely
    determined by \(\tilde{a}\).

    Conversely, let \(H\) be an endoscopic group given by a homomorphism
    \(\OGT_\kappa^\bullet\) and
    \(\tilde{a}_H\in\tilde{\cA}_H^\kappa(\bar{k})\). Then the \(H\)-cameral
    cover \(\tilde{X}_{a_H}\to\breve{X}\) is \'etale over \(U\). Let
    \(U_H\subset\breve{X}\) be the largest subset over which this \(H\)-cameral
    cover is \'etale, then \(U\subset U_H\). So the homomorphism
    \begin{align}
        \pi_{\tilde{a}_H}^\bullet\colon\pi_1(U_H,\infty)\longto
        \bW_{\bH}\rtimes\Out(\bH)
    \end{align}
    induces homomorphism
    \begin{align}
        \pi_{\tilde{a}_H}^{\kappa,\bullet}\colon \pi_1(U,\infty)\longto
        \bW_{\bH}\rtimes\pi_0(\kappa)
    \end{align}
    lying over \(\OGT_\kappa^\bullet\). Let \(\tilde{a}\in \tilde{\cA}_\gg\) be
    the image of \(\tilde{a}_H\), then \(\pi_{\tilde{a}}^\bullet\) is the
    composition of map \(\pi_{\tilde{a}_H}^{\kappa,\bullet}\) with canonical map
    \(\bW_{\bH}\rtimes\pi_0(\kappa)\to\bW\rtimes\Out(\bG)\). Therefore, we have
    \(W_{\tilde{a}}\subset(\bW\rtimes\Out(\bG))_\kappa\) and
    \(I_{\tilde{a}}\subset \bW_{\bH}\), and so
    \(\tilde{a}\in\tilde{\cA}_\kappa(\bar{k})\).
\end{proof}

\subsection{}%
\label{sub:levi_strata}
Similar to \(\kappa\)-strata, one can describe the compliment
\(\tilde{\cA}-\tilde{\cA}^\ELL\) using proper Levi subgroups of \(\bG\)
containing \(\bT\). Suppose
\(\tilde{a}\in(\tilde{\cA}-\tilde{\cA}^\ELL)(\bar{k})\), then
\(\bT^{W_{\tilde{a}}}\) is not finite, hence contains a subtorus \(\bS\). Let
\(\bL\) be the centralizer of \(\bS\) in \(\bG\), then \(W_{\tilde{a}}\) is
contained in the centralizer of \(\bS\) in \(\bW\rtimes\Out(\bG)\), and
\(I_{\tilde{a}}\) is contained in \(\bW_{\bL}\). Let \(\tilde{\bL}\) be the
centralizer of \(\bS\) in \(\bG\rtimes\Out(\bG)\) as in
\Cref{sub:invariant_theory_for_levi}. Therefore, there are induced maps
\begin{align}
    \pi_1(\breve{X},\infty)\longto W_{\tilde{a}}/I_{\tilde{a}}\longto
    \pi_0(\tilde{\bL})\longto\Out(\bG).
\end{align}
Let \(\cA_X^L\) be the mH-base corresponding to
\(\Stack*{\FRC_{\FRM}^L/Z_{\FRM}^L}\), then
the above maps imply that \(\tilde{a}\) comes from a point in
\(\tilde{\cA}_X^L\), similar to endoscopic case.
Combining with
\eqref{eqn:difference_in_dimension_of_mH_bases_fixed_divisor}, we have the
following result:
\begin{proposition}
    \label[proposition]{prop:large_codim_isotropic_locus}
    Let \(b\in\cB_\gg\).
    The codimension of the complement \(\cA_b-\cA_b^\ELL\) of the
    anisotropic locus in \(\cA_b\) goes to \(\infty\) as
    \(\Pair{\rho}{\lambda_b}\to\infty\). If \(Z_G\) does not contain a split
    torus, then the same holds for \(\cA_b^\ANI\) in place of \(\cA_b^\ELL\).
\end{proposition}

\begin{remark}
    A slightly weaker version of \Cref{prop:large_codim_isotropic_locus} was
    proved for split groups in \cite{Ch22}*{Corollary~4.2.8}, however, the
    proof in \textit{loc. cit.} has a minor gap because it only considered
    \emph{standard} Levi subgroups but not their conjugates. Nevertheless,
    \Cref{prop:large_codim_isotropic_locus} can serve as a drop-in replacement
    for that result from Chi's paper because its proof does not require the study of
    multiplicative affine Springer fiber at all.
\end{remark}

\section{Transfer of Global Satake Sheaves}
\label{sec:transfer_of_Satake_sheaves}

In this section, we explain how to transfer a Satake sheaf on \(\cM_X\) to a Satake
sheaf on \(\cM_{H,X}^\kappa\). Recall that a Satake sheaf \(\bcQ\in\Sat_X\) on
\(\cM_X\) is defined as the pullback of a semisimple \(\Arc_{\cB_X}G\)-equivariant
perverse sheaf on the \(N\)-truncated (with \(N\gg 0\)) global Hecke stack
\(\Stack*{\GASch_X}_N\) via the evaluation map \(\ev_N\), appropriately shifted
and Tate-twisted so that it is pure of weight \(0\) after restricting to an
open substack of \(\cM_X\) where one of
\Cref{thm:local_singularity_model_weak,thm:local_singularity_model_main} holds.

\subsection{}
To avoid repeating various technical conditions, we list the ones needed in
this section all at once below:

\begin{enumerate}
    \item We shall restrict ourselves to (the preimage of) the open locus in
        \(\cA_X\) where:
        \begin{enumerate}
            \item it is very \(G\)-ample,
            \item one of
                \Cref{thm:local_singularity_model_weak,thm:local_singularity_model_main}
                holds.
        \end{enumerate}
    \item We further shrink \(\cA_X\) so that its preimage in
        \(\cA_{H,X}^\kappa\) satisfies:
        \begin{enumerate}
            \item it is very \(H\)-ample,
            \item \Cref{thm:local_singularity_model_weak} holds
                (\Cref{thm:local_singularity_model_main} is not stated for
                endoscopic mH-base \(\cA_{H,X}^\kappa\)).
        \end{enumerate}
    \item Since the question is local on the base \(\cA_X\), we may restrict to
        \(\tilde{\cA}_X\) and \(\tilde{\cA}_{H,X}^\kappa\).
    \item Without loss of generality, it suffices to consider Satake sheaf on
        \(\Stack*{\GASch_X}_N\) that is the intermediate extension of its
        restriction to the multiplicity-free locus \(\cB_X^\circ\). As a result,
        we may restrict ourselves to the multiplicity-free locus, so that the
        restricted mH-fibration \(\cM_b\to\cA_b\) is locally constant when \(b\)
        varies in \(\cB_X^\circ\).
    \item To simplify notations, we shall use \(\tilde{\cA}\) and
        \(\tilde{\cA}_H^\kappa\) for the resulting mH-base for \(G\) and \(H\)
        respectively. Let \(\tilde{\cM}\) and \(\tilde{\cM}_H^\kappa\) be
        the corresponding total stacks.
\end{enumerate}

\subsection{}
Recall we have a pointed endoscopic datum
\begin{align}
    \OGT_\kappa^\bullet\colon\pi_1(X,\infty)\longto\pi_0(\kappa),
\end{align}
which is defined over \(k_\infty\). Without loss of generality, we may assume
that \(k_\infty=k\), otherwise we may base change everything to \(k_\infty\) and
replace \(k\) by \(k_\infty\). Fix an algebraic closure \(\bar{k}\) of
\(k_\infty\).
Let \(F=F(X)\) be the function field of \(X\), and we fix a separable closure
\(F^\sep\) over \(\bar{k}\). Let \(\Gamma_F=\Gal(F^\sep/F)\) be the Galois group.
Let \(\breve{F}\) be the everywhere unramified closure of \(F\) in \(F^\sep\),
then the Galois group \(\Gamma=\Gal(\breve{F}/F)\) is canonically identified
with \(\pi_1(X,\infty)\) as a quotient of \(\Gamma_F\) (see
\cite{StacksP}*{\href{https://stacks.math.columbia.edu/tag/0BQM}{Tag 0BQM}} for example). We shall use \(\Gamma\) instead of
\(\Gamma_F\) to form the global \(L\)-group
\begin{align}
    \LD{G}\defeq
    \dual{\bG}\rtimes\Gamma\simeq\dual{\bG}\rtimes\pi_1(X,\infty),
\end{align}
and we have canonical map
\begin{align}
    \LD{G}\longto\dual{\bG}\rtimes\pi_0(\kappa).
\end{align}
Similarly, we have canonical map
\begin{align}
    \LD{H}\longto\dual{\bH}\rtimes\pi_0(\kappa).
\end{align}
This way, we have well-defined global \(L\)-groups of \(G\) and \(H\) over 
the entire \(\tilde{\cA}^\kappa\).
Since \(\Gamma\) has finite image in \(\pi_0(\kappa)\), we may replace
\(\breve{F}\) by a subfield \(F'\) that is finite Galois over \(F\), and replace
\(\Gamma\) by \(\Gal(F'/F)\).

\subsection{}
Recall we also have the canonical embedding of discrete groups
\begin{align}
    \label{eqn:compatibility_with_root_datum_and_W}
    \bW_\bH\rtimes\pi_0(\kappa)\longto \bW\rtimes\pi_0(\kappa).
\end{align}
For simplicity, we first consider the case where \(G^\Der\) is simply-connected,
then a classical result of Langlands
(cf.~\Cref{sec:Appendix_endoscopic_groups}) shows that
\eqref{eqn:compatibility_with_root_datum_and_W} lifts to an admissible
\(L\)-embedding
\begin{align}
    \xi_F\colon\LD{H_F}\to \LD{G_F},
\end{align}
and we shall choose one if there are many.

At any closed point \(v\in X\), the tensor product \(F'\otimes_F F_v\) splits
into a product of unramified extensions \(F_w'\) of \(F_v\). Any geometric point
\(\bar{v}\) over \(v\) (in other words, a \(k\)-embedding of \(k_v\) into
\(\bar{k}\)) determines an unramified closure \(\breve{F}_{\bar{v}}\) of
\(F_v\), and in turn any choice of \(w\) over \(v\) and embedding \(\bar{w}\colon
k_w\to\bar{k}\) over \(\bar{v}\) determines an unramified
closure of \(F_w'\) over \(\breve{F}_{\bar{v}}\). Then the \(L\)-embedding
\(\xi_F\) induces an embedding of local \(L\)-groups
\begin{align}
    \xi_{F_v}\colon \LD{H_{F_v}}\longto\LD{G_{F_v}},
\end{align}
which factors through its Frobenius form
\begin{align}
    \dual{\bH}\rtimes\Gal(\bar{k}/\bar{v}(k_v))\longto\dual{\bG}\rtimes\Gal(\bar{k}/\bar{v}(k_v)).
\end{align}
Changing the choice of \(\bar{w}\) changes the embedding of the Frobenius
element \(\Frob_{H,v}\) by
conjugation by some element in \(\dual{\bG}\). In fact, since the pinnings
are fixed and the \(L\)-embedding descends to
\eqref{eqn:compatibility_with_root_datum_and_W}, the conjugation must be given
by some element \(t\in\dual{\bT}\) such that \(t\Frob_{H,v}(t)^{-1}\in
Z_{\dual{\bH}}\).

\subsection{}
Let \(\bcQ\in\Sat_X\) be a Satake sheaf on \(\tilde{\cM}\).
To give such a \(\bcQ\) is equivalent to give a semisimple algebraic
representation of global \(L\)-group \(\LD{G}\). Indeed, given
representation \(V\), suppose it decomposes into irreducible
\(\dual{\bG}\)-representations as
\begin{align}
    V=\bigoplus_\lambda V_\lambda\otimes\Hom_{\dual{\bG}}(V_\lambda,V),
\end{align}
then for \(\Gamma\)-fixed \(\lambda\), the multiplicity space
\(\Hom_{\dual{\bG}}(V_\lambda,V)\) induces a local system on \(X\),
while \(V_\lambda\) gives a local Satake sheaf \(\IC^\lambda\) at \(\infty\). Since
\(\lambda\) is fixed by monodromy \(\Gamma\), it can be parallel transported to
any point \(v\in X\). If \(\lambda\) is not fixed by monodromy, we can still parallel
transport its \(\Gamma\)-orbit as well as the direct sum of the corresponding
multiplicity spaces. They satisfy obvious cocycle conditions that are also
compatible with local monodromy at any point \(v\in X\). We leave the
verification of the cocycle conditions of this part to the reader.

A more algebro-geometric formulation can be given using the language in
\Cref{sec:boundary_divisors}: the \(\Gamma\)-orbit of \(\lambda\) can be seen as
a connected component \(\sC^\lambda\) of the \'etale covering \(\CoCharG(T)\to
X\), hence we may canonically define \(\IC^\lambda\) at any point in
\(\sC^\lambda\) since \(\lambda\) is constant over \(\sC^\lambda\). The
multiplicity space of \(V_\lambda\) can then be realized as a local system on
\(\sC^\lambda\).

Let \(\bM\) be the monoid whose abelianization is generated by
\(-w_0(\theta_i)\) where \(\Set{\theta_i\given i=1,\ldots,m}\) is the support of
\(V\), and \(\FRM\) is the monoid induced by \(\bM\) and
\(\OGT_\kappa^\bullet\). Then \(V\) induces a canonical Satake sheaf \(\bcQ\) on
\(\tilde{\cM}\) as follows: for any \(b\in \BD_X(\bar{k})\) with
\(\lambda_b=\sum_{\bar{v}}\lambda_{\bar{v}}\bar{v}\), write
\begin{align}
    -w_0(\lambda_{\bar{v}})=\sum_{i=1}^m c_{\bar{v}i}\theta_i
\end{align}
by choosing an identification \(\FRM\cong\bM\), the stalk \(\bcQ_b\) is then
the outer tensor product of the local Satake sheaves corresponding to
\begin{align}
    \bigoplus_{i=1}^m\Sym^{c_{\bar{v}}}\left(V_{\theta_i}\otimes\Hom_{\dual{\bG}}(V_{\theta_i},V)\right).
\end{align}
This description does not depend on the identification \(\FRM\cong\bM\), because
the effect of a different identification on the cocharacters cancels out with
its effect on multiplicity spaces. Alternatively, the connected component of
\(\BD_X\) containing \(b\) can be written as a direct product of symmetric
products of some irreducible components in \(\CoCharG(T)\), over which the local
Satake sheaves are easily defined.

\subsection{}
Conversely, given Satake sheaf \(\bcQ\), the representation \(V\) can be
recovered as follows: let \(\sC\) be the irreducible component of
\(\CoCharG(T)\) corresponding to the generators of the cone of \(\FRA_\FRM\).
Then it canonically embeds into \(\BD_X\) as the space of simple boundary divisors.
Let \(\sC_\infty\) be the fiber of \(\sC\to X\) over \(\infty\), then \(\bcQ\)
induces a canonical local Satake sheaf \(\bcQ_\infty\) on the
affine Grassmannian over \(\sC_\infty\). The cohomology space of \(\bcQ_\infty\)
then recovers \(V\).

\subsection{}
The \(L\)-embedding \(\xi_F\) then induces a transfer of \(\bcQ\) to
a Satake sheaf \(\bcQ_\xi^\kappa\) on \(\tilde{\cM}_H^\kappa\) by restricting
\(V\). It is compatible with the local transfer
defined in \Cref{def:local_transfer_of_Satake_stuff}, because at each geometric place
\(\bar{v}\), the induced local \(L\)-embedding is well-defined up to a
conjugation by some \(t_{\bar{v}}\in\dual{\bT}\) with
\(t\Frob_{H,v}(t)^{-1}\in Z_{\dual{\bH}}\), which has no effect on the
resulting \(\Frob_{H,v}\)-action on \(V\).

\subsection{}
Next, we study the dependence of \(\bcQ^\kappa_\xi\) on \(\xi_F\). Suppose we
have another choice \(\xi_F'\), then it differs from \(\xi_F\) by a
\(1\)-cocycle \(\mu\) in \(Z_{\dual{\bH}}\).
This induces a cohomological class in \(\RH^1(W_F,\dual{\bT})\) (\(W_F\) is the
Weil group of \(F\)), still denoted by \(\mu\). In turn, we obtain a
quasi-character of \(T(\bA_F)\) (\(\bA_F\) is the ring of adeles of \(F\)) that
is trivial on \(T(F)\), hence also on \(Z(F)\).
Unfortunately, there seems to be no way to eliminate the dependence on \(\xi_F\)
unless the group \(G\) is somewhat special.

Fortunately, for arithmetic purposes, this dependence on \(\xi_F\) is
superfluous. Suppose we have a \(G\)-strongly regular semisimple conjugacy class
\(\gamma_H\in H(F)\) that matches \(\gamma\in G(F)\). Then their local conjugacy
classes also match each other. At place \(v\), fix any choice of
local \(L\)-embedding, whose conjugacy class is well-defined. Therefore,
\(\xi_F\) canonically induces local transfer factor \(\Delta_{v}(\gamma_{H,v},\gamma_v)\) as
well as the modified version without the \(\Delta_{\symup{III}_1}\)-factor
\(\Delta_{0,v}(\gamma_{H,v},\gamma_v)\). For any \(z_v\in Z(F_v)\), Langlands and
Shelstad in \cite{LS87}*{Lemma~4.4.A} have shown that
\begin{align}
    \Delta_{0,v}(z_v\gamma_{H,v},z_v\gamma_v)
    =\tau_v^G(z_v)\Delta_{0,v}(\gamma_{H,v},\gamma_v)
\end{align}
for some quasi-character \(\tau_v^G\) of \(Z(F_v)\) independent of
\(\gamma_{H,v}\) or \(\gamma_v\). Moreover, it is not hard to directly see
that when \(z_v\in Z_0(F_v)\) is in the \emph{connected} center,
\(\tau_v^G(z_v)\) is the same as the pairing of \(z_v\) with the
cohomological class in the definition of the term \(\Delta_{\symup{III}_2}\)
(the latter, \textit{a priori}, does depend on conjugacy classes). If \(z_v\)
comes from a point \(z\in Z_0(F)\), then we have product formula
\begin{align}
    \prod_v\tau_v^G(z_v)=1,
\end{align}
as long as the \(\chi\)-datum is chosen globally, because the left-hand side is
none other than the pairing of a global cohomological class in
\(\RH^1(W_F,\dual{\bZ}_0)\) (\(\dual{\bZ}_0\) being the dual group of \(Z_0\))
with an element in \(Z_0(F)\), the latter of which becomes trivial in the adelic
quotient. Therefore, we have
\begin{align}
    \Delta_0(z\gamma_H,z\gamma)=\Delta_0(\gamma_H,\gamma).
\end{align}
Similarly, changing \(\xi_F\) by \(\mu\) does not affect global transfer factor
either.

\subsection{}
Now we move on to the case where a global \(\xi_F\) does not exist. Let
\begin{align}
    1\longto Z_1\longto G_1\longto G\longto 1
\end{align}
be a \(z\)-extension, \(H_1\) be the endoscopic group of \(G_1\) induced by
\(H\). Then \(G_1^\Der\) is simply-connected and a global \(L\)-embedding
\(\xi_{F,1}\) exists for \(\LD{H}_1\). Consider any representation \(V^1\) of
\(\LD{G}_1\). It restricts to both a representation \(V\) of \(\LD{G}\) and
\(V^{H_1}\) of \(\LD{H}_1\), the latter of which further restricts to
\(\LD{H}\), denoted by \(V^H\). In any local place \(v\) of \(F\), the map
\(G_1(F_v)\to G(F_v)\) is surjective, and so every \(F_v\)-rational highest
coweight of \(G\) is an image of that of \(G_1\).
Thus, for the purpose of fundamental lemma, which is a local statement, it suffices
to consider \(\LD{G}\)-representations obtained from all such \(V^1\). Let \(\bcQ\) be the
Satake sheaf given by \(V\) and let \(\bcQ_\xi^\kappa\) be the Satake sheaf on
\(\tilde{\cM}^\kappa\) induced by \(V^H\).

The embedding \(\xi_{F,1}\) induces a cocycle of \(W_F\) in \(\dual{\bZ}_1\)
(the dual group of torus \(Z_1\)), which then corresponds to a quasi-character
\(\mu_1\) of \(Z_1(\bA_F)/Z_1(F)\). Let \(f_{1,v}\) (resp.~\(f_v\)) be the local Satake function on
\(G_1(F_v)\) (resp.~\(G(F_v)\)) at \(v\) induced by \(V^1\) (resp.~\(V\)), and let \(\tilde{f}_{1,v}\) be the
one on \(G_1(F_v)\) obtained by pulling back \(f_v\), then we have
\begin{align}
    \tilde{f}_{1,v}(g_{1,v})=\sum_{z_{1,v}\in
    Z_1(F_v)/Z_1(\cO_v)}f_{1,v}(z_{1,v}g_{1,v}),
\end{align}
where the right-hand side is necessarily a finite sum.
Let \(f_{H_1,\xi,v}\) (resp.~\(f_{H,\xi,v}\), resp.~\(\tilde{f}_{H_1,\xi,v}\))
be the analogue of \(f_{1,v}\) (resp.~\(f_v\), resp.~\(\tilde{f}_{1,v}\)).
Then we necessarily have for any
\(z_{1,v}\in Z_1(F_v)\) and \(h_{1,v}\in H_1(F_v)\) that
\begin{align}
    \tilde{f}_{H_1,\xi,v}(z_{1,v}h_{1,v})=\mu_1(z_{1,v})\tilde{f}_{H_1,\xi,v}(h_{1,v}).
\end{align}
It is shown
in \cite{LS87}*{\S~4} that we in fact have \(\mu_{1,v}=\tau_v^{G_1}\) as a
quasi-character of \(Z_1(F_v)\) for every \(v\). Therefore, proving the
fundamental lemma for \((\gamma_{H,v},\gamma_v)\in H(F_v)\x G(F_v)\) is the same as proving that
for any matching lift \((\gamma_{H_1,v},\gamma_{1,v})\in H_1(F_v)\x G_1(F_v)\).
In particular, the choice of \(G_1\) is inconsequential for proving the
fundamental lemma.

\begin{definition}
    \label[definition]{def:global_transfer_Satake}
    For any Satake sheaf \(\bcQ\in\Sat_X\),
    \nomenclature[\(Q"cal \)]{\(\bcQ\)}{a Satake sheaf on the mH-total-stack}
    suppose it comes from a
    \(G_1\)-representation as above, we define
    its \notion{transfer}\index{transfer!of Satake sheaves} \(\bcQ^\kappa_{H,\xi}\)
    \nomenclature[\(Q"cal_H_xi_kappa \)]{\(\bcQ^\kappa_{H,\xi},\bcQ^\kappa_{H}\)}{the endoscopic
    transfer of \(\bcQ\) to \(H\) (via \(L\)-embedding \(\xi\))}
    to \(H\) to be the complex
    \(\bcQ^\kappa_\xi\) defined using an arbitrarily fixed \(\xi_{F,1}\). If
    the choice of \(G_1\) and \(\xi_{F,1}\) is implicit in the context, we
    simply use \(\bcQ_H^\kappa\).
\end{definition}

\section{Geometric Stabilization}%
\label{sec:geometric_stabilization}

In \cite{Ng10}*{Th\'eor\`emes~6.4.1, 6.4.2},
Ng\^o established the geometric stabilization theorem for the
usual Hitchin fibrations. It is essentially the geometric side of
the stabilization process of the trace formula for Lie algebras over
function fields, and one can deduce from (a slightly weaker version of) it the
endoscopic fundamental lemma for Lie algebras. Our goal for
mH-fibrations is the same.

\subsection{}
We retain the notations from \Cref{sec:transfer_of_Satake_sheaves}, and let
\(\tilde{\cA}^\ANI\) be the anisotropic locus in \(\tilde{\cA}\). Let
\(\tilde{\cA}^\ANI_\kappa\) be the intersection of \(\tilde{\cA}^\ANI\)
with the union of the images \(\tilde{\cA}^\kappa_H\) for all possible \(H\)
associated with \(\kappa\) (but different \(\OGT_\kappa^\bullet\)). The reader
should refer to \Cref{sec:transfer_of_Satake_sheaves} about the technical
conditions satisfied by \(\tilde{\cA}^\ANI_\kappa\), and note that it is
very slightly smaller than the object with the same notation in sections before
\Cref{sec:transfer_of_Satake_sheaves}.

\begin{theorem}
    [Geometric Stabilization]
    \label[theorem]{thm:main_geometric_stabilization}
    Suppose \(Z_G\) contains no split subtorus. Let \(\bcQ\in\Sat_X\) be any
    Satake sheaf induced by some representation \(V^1\) of \(\LD{G}_1\) where
    \(G_1\) is a \(z\)-extension of \(G\). Then there exists an isomorphism
    between graded shifted perverse sheaves over \(\tilde{\cA}^\ANI_\kappa\)
    \begin{align}
        \bigoplus_n\PH^n(\tilde{h}_{*}^\ANI\bcQ)_\kappa\cong
        \bigoplus_n\bigoplus_{(\kappa,\OGT_{\kappa}^\bullet)}
        \tilde{\nu}_{\OGT_\kappa^\bullet,*}^\ANI\PH^n(\tilde{h}_{H,*}^{\kappa,\ANI}\bcQ_{H}^\kappa)_{\hST},
    \end{align}
    where \((\kappa,\OGT_{\kappa}^\bullet)\) ranges over all pointed endoscopic
    data for a fixed \(\kappa\), \(H\) is the corresponding endoscopic group,
    \(\tilde{\nu}_{\OGT_\kappa^\bullet}\) is the endoscopic transfer map
    \(\tilde{\nu}_\cA\) corresponding to \(H\), and \(\bcQ_{H}^\kappa\) is the
    transfer of complex \(\bcQ\) defined in
    \Cref{def:global_transfer_Satake}.
\end{theorem}

\subsection{}
Since \(\tilde{\nu}_{\OGT_\kappa^\bullet}\) for different
\(\OGT_\kappa^\bullet\) have disjoint image, we have a more refined version of
\Cref{thm:main_geometric_stabilization} in some specials cases:

\begin{theorem}
    \label[theorem]{thm:arithmetic_geometric_stabilization}
    With the assumptions in \Cref{thm:main_geometric_stabilization}, suppose in
    addition that \(Z_G\) is connected and satisfies Hasse principle\index{Hasse principle} over the
    function field \(F\) of \(X\), namely the group
    \begin{align}
        \ker^1(F,Z_G)\defeq\ker\left[\RH^1(F,Z_G)\longto\prod_{v\in\abs{X}}\RH^1(F_v,Z_G)\right]
    \end{align}
    is trivial. Suppose
    \((\kappa,\OGT_{\kappa,\xi}^\bullet)\) is defined over \(k\), then
    there exists a \(k\)-isomorphism of semisimplifications of graded shifted
    perverse sheaves over \(\tilde{\cA}^\ANI_\kappa\)
    \begin{align}
        \Big(\bigoplus_n\PH^n(\tilde{h}_{*}^\ANI\bcQ)_\kappa
        \Big)\Big|_{\tilde{\nu}_{\OGT_\kappa^\bullet}(\tilde{\cA}_{H,X}^{\kappa})\cap
        \tilde{\cA}^\ANI}
        \cong \tilde{\nu}_{\OGT_\kappa^\bullet,*}^\ANI\bigoplus_n
        \PH^n(\tilde{h}_{H,*}^{\kappa,\ANI}\bcQ_{H}^\kappa)_{\hST}.
    \end{align}
\end{theorem}

\begin{remark}
    When the conditions on \(Z_G\) in
    \Cref{thm:arithmetic_geometric_stabilization} are satisfied,
    \Cref{thm:main_geometric_stabilization} is a direct consequence of
    \Cref{thm:arithmetic_geometric_stabilization}. The Hasse principle condition
    is satisfied, for example, when the action of \(\Gal(\bar{F}/F)\) on \(Z_G\)
    factors through a cyclic quotient according to \cite{On63}*{Proposition~4.5.1}.
\end{remark}

\subsection{}
Both \Cref{thm:main_geometric_stabilization,thm:arithmetic_geometric_stabilization}
will be proved in \Cref{chap:counting_points} after some back-and-forth
between global and local arguments, similar to what is done in
\cite{Ng10}*{\S~8}.

\section{Top Ordinary Cohomology}%
\label{sec:top_ordinary_cohomology}

In this section we study the top ordinary cohomology of
\(\tilde{h}_{X,*}^\ANI\bcQ\). We hope to give a description similar to that
in \cite{Ng10}*{\S~6.5} which can be used in tandem with support
theorem. Interestingly, it turns out to be drastically different from the Lie
algebra case and much more complicated.

\subsection{}
First we consider a general map \(f\colon X\to Y\) of \(k\)-varieties (or more
generally Deligne--Mumford stacks locally of finite types over \(k\)) of relative
dimension \(d\). The complex \(f_!\Qlb\) will have cohomological amplitude
\([0,2d]\). The stalk of the sheaf of top cohomology \(\RDF^{2d}{f}_!\Qlb\) at
geometric point \(y\in Y\) is a \(\Qlb\)-vector space with a canonical basis in
bijection with irreducible components of fiber \(X_y\). If \(X\) is smooth, then
the constant sheaf \(\Qlb\) is a pure complex, and so is \(f_!\Qlb\) if in
addition \(f\) is proper.

On the other hand, if \(f\) is proper but \(X\) is not smooth, then
\(f_*\Qlb=f_!\Qlb\) is not necessarily pure, and one needs to replace \(\Qlb\)
with intersection complex \(\IC_X\) in order to restore purity. However, in
general \(f_*\IC_X\) may have larger cohomological amplitude than \(2d\).

Unlike the Lie algebra case, the total stack of mH-fibration is not smooth in
general; rather we have a local model of singularity established in
\Cref{thm:local_singularity_model_weak,thm:local_singularity_model_main}.
Therefore, in order to give a nice description of the top ordinary
cohomology, we need to establish cohomological amplitude of
\(\tilde{h}_{X,*}^\ANI\bcQ\). Without loss of generality, we may assume that
\(\bcQ\) is the intersection complex on the total stack, because the case of
general Satake sheaves can be easily deduced from this special case.

\subsection{}
The question is local in \(\tilde{\cA}_X^\ANI\), and by proper base change,
we may fix \(b\in\cB_X(\bar{k})\) and restrict to \(\tilde{\cA}_b^\ANI\).
We keep the assumption in \Cref{sec:transfer_of_Satake_sheaves} that the local
model of singularity as in
\Cref{thm:local_singularity_model_weak,thm:local_singularity_model_main} holds.
To simplify notations, let \(\cA=\tilde{\cA}_b^\ANI\) and let
\(\cM=\tilde{h}_X^{-1}(\cA)\). Here the boundary divisor
\(b\) is fixed, so
the local model is just a finite direct product
\begin{align}
    \SFQ\defeq \prod_{i=1}^m
    \Stack*{\Arc_{\bar{v}_i,N}G^\AD\big\backslash\Gr_{G^\AD,\bar{v}_i}^{\le-w_0(\lambda_i)}},
\end{align}
where \(\bar{v}_i\in X(\bar{k})\) are distinct points and \(N\) is a
sufficiently large integer depending on \(b\). Let \(\ev_N\colon\cM\to\SFQ\) be
the evaluation map as in \Cref{thm:local_singularity_model_weak}. Let \(e\)
be the relative dimension of \(\ev_N\).

When \Cref{thm:local_singularity_model_weak} or
\Cref{thm:local_singularity_model_main} holds, we have by
\Cref{cor:local_model_sing_sheaf_version} that
\begin{align}
    \label{eqn:IC_of_restr_mH_by_local_model_in_top_coh}
    \bcQ|_{\cM}=\ev_N^*\cF[e](e/2)
\end{align}
for some perverse sheaf \(\cF\) on \(\SFQ\). We also know that
\(\IC_{\SFQ}\) is the unique direct summand of \(\cF\) supported on the entire
\(\SFQ\), and if \(b\) is multiplicity-free, then \(\cF\) is exactly \(\IC_\SFQ\).

We know that simple equivariant perverse sheaves on \(\SFQ\) of weight \(0\) are
none other than the intersection complexes of substacks
\begin{align}
    \SFQ'\defeq \prod_{i=1}^m
    \Stack*{\Arc_{\bar{v}_i,N}G^\AD\big\backslash\Gr_{G^\AD,\bar{v}_i}^{\le-w_0(\lambda_i')}}
\end{align}
for some \(\lambda_i'\le\lambda_i\). Let
\(\lambda_b=\sum_{i=1}^m\lambda_i\cdot\bar{v}_i\) and
\(\lambda_b'=\sum_{i=1}^m\lambda_i'\cdot\bar{v}_i\). If
\Cref{thm:local_singularity_model_weak} holds, by smoothness of
\(\ev_N\), we can compute the codimension of \(\cM'=\ev_N^{-1}(\SFQ')\):
\begin{align}
    \label{eqn:codimension_of_smaller_Schubert_to_mH_total}
    2\delta'\defeq\codim_{\cM}(\cM')
    =\codim_{\SFQ}(\SFQ')
    =\sum_{i=1}^m\Pair{2\rho}{\lambda_i-\lambda_i'}
    =\Pair{2\rho}{\lambda_b-\lambda_b'}.
\end{align}
If instead \Cref{thm:local_singularity_model_main} holds, then \(\ev_N\) may not
be smooth but \(\ev_{1,N}^I\) is. Note that in
\Cref{thm:local_singularity_model_main} the boundary divisor \(b\) is not
fixed, but its subdivisor \(b_1\) is. When the complementary divisor \(b_2\)
varies, the fiber \(\cM_b\) stays locally constant, so
\eqref{eqn:codimension_of_smaller_Schubert_to_mH_total} still holds for any
\(b\) lying over \(b_1\).

Since the substack \(\cM'\) is the total stack of the restricted mH-fibration with boundary
divisor \(\lambda_b'\), we can use \eqref{eqn:dimension_of_restricted_mH_base} to compute
the codimension of \(\cA'=\tilde{h}_X^\ANI(\cM')\):
\begin{align}
    \codim_{\cA}(\cA')=
    \sum_{i=1}^m\Pair{\rho}{\lambda_i-\lambda_i'}
    =\OneHalf\codim_{\cM}(\cM')=\delta'.
\end{align}
The restriction of \(\tilde{h}_X^\ANI\) to \(\cM'\) is of
relative dimension \(d-\delta'\), where \(d=\dim{\cM}-\dim{\cA}\).

Since the collection of all \(\SFQ'\) for different \(\lambda_b'\le \lambda_b\)
induce a stratification of \(\SFQ\) by smooth substacks, the construction of
intermediate extension functor implies that the support of \(\RH^i(\IC_{\SFQ})\)
is contained in the union of those \(\SFQ'\) with dimension at most \(-i\), and
the equality is achieved if and only if \(-i=\dim{\SFQ}\). This implies that
\(\IC_{\cM}=\ev_N^*\IC_\SFQ[e](e/2)\) has its \(i\)-th cohomology supported on
those \(\cM'\) with dimension at most \(-i\), and the equality is achieved if
and only if \(-i=\dim{\cM}\). In other words, if \(\cM'\) is contained in the
support of \(\RH^i(\IC_{\cM})\), then \(2\delta'\ge i+\dim{\cM}\), with equality
achieved if and only if \(\delta'=0\). This shows that
\(\tilde{h}_{X,*}^\ANI(\RH^i(\IC_{\cM}))\) is supported on cohomological
degrees
\begin{align}
    [-\dim{\cM},-\dim{\cM}+2d],
\end{align}
and with the upper bound achieved if and only if \(i=-\dim{\cM}\). A standard
argument using spectral sequence then shows that
\(\tilde{h}_{X,*}^\ANI(\IC_\cM)\) has cohomological degree bounded above by
\(-\dim{\cM}+2d\), and the stalk of its top cohomology at any geometric point
\(a\in\cA\) has a basis in bijection with the irreducible components of \(\cM_a\).

The same argument can be applied to
\(\IC_{\cM'}=\ev_N^*\IC_{\SFQ'}[e-2\delta'](e/2-\delta')\),
and we can see that \(\tilde{h}_{X,*}^\ANI(\IC_{\cM'})\) has cohomological
degree bounded above by \(-\dim{\cM'}+2(d-\delta')=-\dim{\cM}+2d\), and its top
cohomology over \(a\) has
a basis in bijection with the irreducible components of \(\cM'_a\). This
allows us to cover cases where \(b\) is not multiplicity-free, as well as when
\(\bcQ\) is a general Satake sheaf.
Thus far we have shown that
\begin{proposition}
    \label[proposition]{prop:cohomological_amplitude_of_mH_fibrations}
    Let \(\cA\subset\tilde{\cA}_X^\ANI\) be any irreducible open substack
    over which there is a local model of singularity as in
    \Cref{thm:local_singularity_model_weak} or
    \Cref{thm:local_singularity_model_main},
    and \(\cM=\tilde{h}_X^{-1}(\cA)\).
    Then for any Satake sheaf \(\bcQ\in\Sat_X\), the complex
    \(\tilde{h}_{X,*}^\ANI\bcQ|_{\cA}\) is supported on
    cohomological degrees
    \begin{align}
        [-\dim{\cM},-\dim{\cM}+2d],
    \end{align}
    where \(d=\dim{\cM}-\dim{\cA}\). In addition, if \(\tilde{a}\in\cA\) lies over a
    multiplicity-free boundary divisor and \(\bcQ\) is the
    intersection complex \(\IC_\cM\), then
    \(\RH^{-\dim{\cM}+2d}(\tilde{h}_{X,*}^\ANI\bcQ)_{\tilde{a}}\) has a canonical
    basis being the irreducible components of \(\cM_{\tilde{a}}\).
\end{proposition}
 
\subsection{}
We retain the notations in
\Cref{prop:cohomological_amplitude_of_mH_fibrations}. We let
\(p\colon\cP\to\cA\) be the pullback of \(\cP_X\to\cA_X\) to \(\cA\), and let
\(h\) be that of \(\tilde{h}_X\) to \(\cA\). We also stick to the
case where \(\bcQ=\IC_{\cM}\) and the general case is straightforward.

As we will see in \Cref{sub:kappa_decomposition_at_chain_level}, the action of
\(\cP\) on
\(h_*\bcQ\) factors through sheaf of finite abelian
groups \(\pi_0(\cP)\), and
so does its action on \(\RDF^{-\dim{\cM}+2d}h_*\bcQ\). Let \(\tilde{a}\in\cA(\bar{k})\)
be a point lying over a multiplicity-free boundary divisor
\(\lambda_b=\sum_{i=1}^m\lambda_i\bar{v}_i\). By
\Cref{prop:cohomological_amplitude_of_mH_fibrations}, the
\(\cP_{\tilde{a}}\)-action on \((\RDF^{-\dim{\cM}+2d}h_*\bcQ)_{\tilde{a}}\) is just the action
induced by the \(\pi_0(\cP_{\tilde{a}})\)-action on the set of irreducible components of
\(\cM_{\tilde{a}}\).

Using product formula (\Cref{prop:product_formula}), we
see that \((\RDF^{-\dim{\cM}+2d}h_*\bcQ)_{\hST,\tilde{a}}\) has a canonical basis in
bijection with the (necessarily finite) direct product
\begin{align}
    \prod_{\bar{v}\in
    X(\bar{k})}\Irr\Stack*{\cM_{\bar{v}}(\tilde{a})/\cP_{\bar{v}}(\tilde{a})}.
\end{align}
If \(a\) is unramified at point \(\bar{v}_i\), then its Newton point \(\nu_i\)
is integral, and we know
\begin{align}
    \Cnt\Irr\Stack*{\cM_{\bar{v}_i}(\tilde{a})/\cP_{\bar{v}_i}(\tilde{a})}=m_{-w_0(\lambda_i)\nu_i}
\end{align}
according to \Cref{thm:local_irr_components_weight_mult}.
For the same reason, when \(\bar{v}\neq\bar{v}_i\)
for any \(i\), the regular locus \(\cM_{\bar{v}}(\tilde{a})^\reg\) is dense in
\(\cM_{\bar{v}}(\tilde{a})\) and is a \(\cP_{\bar{v}}(\tilde{a})\)-torsor. So in this case
\begin{align}
    \Cnt\Irr\Stack*{\cM_{\bar{v}}(\tilde{a})/\cP_{\bar{v}}(\tilde{a})}=1.
\end{align}

Suppose \(\tilde{a}\) is very \((G,N)\)-ample for some
\(N=N(\delta_{\tilde{a}})\) as in \Cref{prop:delta_regularity_fix_divisor}.
This ensures that the \(\delta\)-stratum in \(\cA\) containing \(\tilde{a}\) has
codimension at least \(\delta_{\tilde{a}}\) in \(\cA\). Furthermore, assume
\(\tilde{a}\) is a sufficiently general point in the
\(\delta_{\tilde{a}}\)-critical subset \(\cA_{\delta_{\tilde{a}}}^=\). Then we know by
\Cref{cor:when_delta_regularity_equality_is_met} that \(\tilde{a}\) is
unramified at each \(\bar{v}_i\) and its Newton point \(\nu_i\) at \(\bar{v}_i\)
is integral.
Thus, we have the following:
\begin{proposition}
    \label[proposition]{prop:rank_irr_stable_of_delta_critical_general}
    Let \(\tilde{a}\in\cA(\bar{k})\) lying over a multiplicity-free boundary
    divisor \(\lambda_b=\sum_{i=1}^m\lambda_i\bar{v}_i\) be such that it is
    unramified at every \(\bar{v}_i\). Then
    \((\RDF^{-\dim{\cM}+2d}h_*\IC_\cM)_{\hST,\tilde{a}}\) is a \(\Qlb\)-vector
    space of rank \(\sum_im_{-w_0(\lambda_i)\nu_i}\). In particular, it is true
    for \(\tilde{a}\) that is very \((G,N(\delta_{\tilde{a}}))\)-ample and is a
    sufficiently general point of \(\cA_{\delta_{\tilde{a}}}^=\).
\end{proposition}

\subsection{}
Suppose we have a \(\kappa\)-stratum
\(\cA_{(\kappa,\OGT_{\kappa}^\bullet)}\subset\cA\) corresponding to
endoscopic datum \((\kappa,\OGT_{\kappa}^\bullet)\) and endoscopic group
\(H\). Let \(\cA_H^\kappa\) be the preimage of
\(\cA_{(\kappa,\OGT_{\kappa}^\bullet)}\) in \(\tilde{\cA}_{H,X}^\kappa\).
Note that \(\cA_H^\kappa\) may still have multiple irreducible components of
various dimensions even if \(\cA_{(\kappa,\OGT_{\kappa}^\bullet)}\) is
irreducible.

Suppose that \(\cA_H^\kappa\) is very \(H\)-ample, then
\(\delta_{\tilde{a}_H}=0\) for any general \(\tilde{a}_H\) in each irreducible
component of \(\cA_H^\kappa\). Let \(\cU_H\subset\cA_H^\kappa\) be an
irreducible component and let \(\cU\) be its image in
\(\cA_{(\kappa,\OGT_{\kappa}^\bullet)}\).
By upper-semicontinuity, the \(\delta\)-invariant achieves minimal value
over an open dense subset of \(\cU\), and let
\(\delta_{\cU}\) be this value.

If \(\codim_{\cA}(\cU)=\delta_{\cU}\), then by definition \(\cU\) is
\(\delta\)-critical. In fact, we know that the image of any \(\delta_H\)-critical
stratum in \(\cU_H\) is also \(\delta\)-critical in \(\cA\), because the
difference \(\delta-\delta_H\) is constant throughout \(\cU\) by
\eqref{eqn:difference_in_delta_invariants}. For any \(\delta'>0\), if
we further assume that \(\cU_H\) is very \((H,N(\delta'))\)-ample (here
\(N(\delta')\) depends on both \(H\) and \(\delta'\)),
then by \Cref{prop:delta_regularity_fix_divisor}
applied to \(\cU_H\), we have that
\begin{align}
    \codim_{\cA}(\cU\cap\cA_{\delta_\cU+\delta'})\ge \delta_\cU+\delta'.
\end{align}

\subsection{}
Suppose now \(\tilde{a}_H\) is a sufficiently general \(\bar{k}\)-point of a
\(\delta\)-critical stratum \(\cU_{H,\delta'}\) in \(\cU_H\), and suppose it is
very \((H,N(\delta'))\)-ample. Let
\begin{align}
    \lambda_{H,b_H}=\sum_{i=1}^m\lambda_{H,i}\cdot\bar{v}_i
\end{align}
be the boundary divisor of \(\tilde{a}_H\). Let \(\tilde{a}\) be the image of
\(\tilde{a}_H\) in \(\cU\), then it is also a general point in a
\(\delta\)-critical stratum in \(\cA\). Let
\begin{align}
    \lambda_{b}=\sum_{i=1}^m\lambda_{i}\cdot\bar{v}_i
\end{align}
be the boundary divisor of \(\tilde{a}\).

There is no guarantee that
\(\tilde{a}\) is very \((G,N(\delta_{\tilde{a}}))\)-ample, so we
cannot directly apply \Cref{cor:when_delta_regularity_equality_is_met} to
\(\tilde{a}\). Nevertheless, we may apply it to \(\tilde{a}_H\)
since we assume \(\tilde{a}_H\) is \((H,N(\delta'))\)-very ample, and it implies that
\(\tilde{a}_H\) is unramified at every \(\bar{v}_i\). Since the local
ramification index \(c_{\bar{v}_i}(\tilde{a})\)
(resp.~\(c_{H,\bar{v}_i}(\tilde{a}_H)\)) depends only on the generic
fiber of the regular centralizer \(\FRJ_{\tilde{a}}\)
(resp.~\(\FRJ_{H,\tilde{a}_H}\)), and we know generically \(\FRJ_{\tilde{a}}\)
and \(\FRJ_{H,\tilde{a}_H}\) are canonically isomorphic, we have
\begin{align}
    c_{\bar{v}_i}(\tilde{a})=c_{H,\bar{v}_i}(\tilde{a}_H)=0.
\end{align}
In other words, \(\tilde{a}\) must be unramified at every \(\bar{v}_i\).
Using the same argument as
\Cref{prop:rank_irr_stable_of_delta_critical_general}, and replacing the stable
constituent by \(\kappa\)-isotypic constituent, we reach a similar description.
\begin{proposition}
    Suppose \(\tilde{a}_H\) is a general point of a \(\delta\)-critical stratum
    in \(\tilde{\cA}_{H,X}^{\kappa,\ANI}\) and is very
    \((H,N(\delta_{H,\tilde{a}_H}))\)-ample. Suppose \(\tilde{a}\in\cA\) is the
    image of \(\tilde{a}_H\), and suppose it has multiplicity-free boundary
    divisor. Then \((\RDF^{-\dim{\cM}+2d}h_*\IC_\cM)_{\kappa,\tilde{a}}\) has rank
    \(\sum_im_{-w_0(\lambda_i)\nu_i}\).
\end{proposition}
Since at each \(\bar{v}_i\) the Newton point \(\nu_i\) depends only on the
restriction of \(\tilde{a}\) to the generic point of \(\breve{X}_{\bar{v}_i}\), and
similarly for the Newton point \(\nu_{H,i}\) of
\(\tilde{a}_H\), we have that \(\nu_i=\nu_{H,i}\) viewed as \(\bbQ\)-cocharacters of
\(T(\breve{F}_{\bar{v}_i})\). If \(\tilde{a}\) has multiplicity-free
boundary divisor, then so is
\(\tilde{a}_H\) by construction of \(\FRM_H\). Using
\Cref{prop:rank_irr_stable_of_delta_critical_general} and suppose local model of
singularity exists in a neighborhood of \(\tilde{a}_H\) too, then we have
\begin{align}
    \dim_{\Qlb}(\RDF^{-\dim_{\tilde{a}_H}{\cM_H}+2d_{H}}h_{H,*}\bcQ_H^\kappa)_{\hST,\tilde{a}_H}
    =\sum_{i=1}^m
    m_{-w_{H,0}(\lambda_{H,i})\nu_i}\dim_{\Qlb}\Hom_{\dual{\bH}}(V^H_{-w_{H,0}(\lambda_{H,i})},V_{-w_0(\lambda_i)}),
\end{align}
where \(\bcQ_H^\kappa\) is the transfer of \(\bcQ\) (which we assume to be
\(\IC_\cM\) here) in
\Cref{sec:transfer_of_Satake_sheaves}, \(\dim_{\tilde{a}_H}{\cM_H}\) is the
dimension of the component of
\(\tilde{\cM}_{H,X}^\kappa\) containing \(h_H^{-1}(\tilde{a}_H)\), and
\(d_{H}\) is the relative dimension of mH-fibration at \(\tilde{a}_H\).
When \(\bcQ\) is a general Satake sheaf, we replace \(V_{-w_0(\lambda_i)}\) by
appropriate representations of \(\dual{\bG}\).

\subsection{}
Recall that \(-w_{H,0}(\lambda_{H,i})\) is an
\(\dual{\bH}\)-highest weight in the decomposition of the irreducible
\(\dual{\bG}\)-representation with highest-weight \(-w_0(\lambda_i)\), and
so we must have
\begin{align}
    \label{eqn:condition_for_extension_to_base_for_H}
    \nu_i\le_{H} -w_{H,0}(\lambda_{H,i})\le-w_0(\lambda_i).
\end{align}
With this restriction in mind, if the set of \(\tilde{a}_H\) is non-empty for a
fixed \(\tilde{a}\), then among those \(\tilde{a}_H\) mapping
to a fixed \(\tilde{a}\), they all restrict to the same map
\(\tilde{a}_H^\circ\) from
\(\breve{X}-\supp(b)\) to \(\Stack{\FRC_{\FRM,H}/Z_\FRM^\kappa}\). The set of
ways to extend \(\tilde{a}_H^\circ\) over \(\bar{v}_i\) is in natural bijection
with the set of \(-w_{H,0}(\lambda_{H,i})\) such that
\eqref{eqn:condition_for_extension_to_base_for_H} holds; see
\Cref{prop:pointed_endoscopic_transfer_closed_embedding} for the precise
statement. Therefore, we reach equality
\begin{align}
    \label{eqn:dimensional_stabilization_for_top_cohomology}
    \dim_{\Qlb}(\RDF^{-\dim{\cM}+2d}h_*\bcQ)_{\kappa,\tilde{a}}
    =\sum_{\tilde{a}_H\mapsto\tilde{a}}
    \dim_{\Qlb}(\RDF^{-\dim_{\tilde{a}_H}{\cM_H}+2d_{H}}h_{H,*}\bcQ_H^\kappa)_{\hST,\tilde{a}_H}.
\end{align}
Note that \Cref{eqn:dimensional_stabilization_for_top_cohomology} holds for
a general Satake sheaf, not just the intersection complex of \(\cM\), and it
does not require \(\tilde{a}\) to have multiplicity-free divisor either, due
to the representation-theoretic interpretation of \(\bcQ\) and \(\bcQ_H^\kappa\).

If furthermore both \(\tilde{a}\) and \(\tilde{a}_H\) are defined over some finite
extension \(k'/k\) inside \(\bar{k}\),
then product formula \Cref{prop:product_formula} holds
over \(k'\). Combining with some general facts about point-counting in
\Cref{chap:counting_points}, eventually (see
\Cref{prop:primal_stabilization_for_top_cohomology}) we are going to show that under some mild
conditions, \eqref{eqn:dimensional_stabilization_for_top_cohomology} can be
upgraded to an isomorphism of Frobenius modules
\begin{align}
    \label{eqn:primal_stabilization_for_top_cohomology_preview}
    (\RDF^{-\dim_{\tilde{a}}{\cM}+2d}h_*\bcQ)_{\kappa,\tilde{a}}
    \simeq
    \bigoplus_{\tilde{a}_H\mapsto\tilde{a}}(\RDF^{-\dim_{\tilde{a}_H}{\cM_H}+2d_{H}}h_{H,*}\bcQ_H^\kappa)_{\hST,\tilde{a}_H}.
\end{align}

\begin{remark}
    Isomorpohism \eqref{eqn:primal_stabilization_for_top_cohomology_preview} can
    be viewed as a sort of primal form of
    \Cref{thm:main_geometric_stabilization}.
\end{remark}

\subsection{}
Up until now in this section we made a lot of assumptions. For reader's
convenience we summarize the essential ones below:
\begin{enumerate}
    \item For both \(G\) and \(H\), we require the existence of respective local
        model of singularity as in either
        \Cref{thm:local_singularity_model_weak} or
        \Cref{thm:local_singularity_model_main}.
    \item For the stable constituent of the cohomology, some ampleness condition on
        the boundary divisor depending only on the group and
        \(\delta\)-invariant (the curve \(X\) is always fixed).
    \item For the \(\kappa\)-constituent, only ampleness condition on \(H\)-side is
        required, not for \(G\). We still require
        \(\delta\)-criticality on the \(G\)-side, but this condition will be
        automatic in practice, thanks to
        \eqref{eqn:difference_in_dimension_of_mH_bases}.
\end{enumerate}

\subsection{}

So far our description of the top cohomology \(\RDF^{-\dim{\cM}+2d}h_*\bcQ\) is
at the stalk level. Due to the jump in ranks, there does not appear to be an easy
description of the top cohomology as a sheaf even just for the stable
constituent. The reader can compare to the Lie algebra case where the stable top
cohomology is just the constant sheaf of rank \(1\) over the Hitchin base
(see \cite{Ng10}*{Proposition~6.5.1}).

Nevertheless, it is still possible to describe
\((\RDF^{-\dim{\cM}+2d}h_*\bcQ)_{\hST}\) over an open subset of each
\(\delta\)-critical stratum. For example, there exists an open dense subset
\(\cU\subset\cA\) over which \(\delta=0\) and \(\cM\) is a \(\cP\)-torsor,
\((\RDF^{-\dim{\cM}+2d}h_*\bcQ)_{\hST}\) is isomorphic to constant sheaf
\(\Qlb\) up to a Tate twist. We put a more precise statement in the following
lemma:
\begin{lemma}
    Let \(\cU\subset\cA\) be the open dense locus where the discriminant
    divisor and the boundary divisor do not collide. Then
    \((\RDF^{-\dim{\cM}+2d}h_*\bcQ)_{\hST}|_{\cU}\) is isomorphic to the
    constant sheaf \(\Qlb\) up to Tate twist.
\end{lemma}
\begin{proof}
    By \Cref{cor:central_lambda_and_d_1_implies_reg} and the product formula
    \Cref{prop:product_formula}, we have that for any
    \(\tilde{a}\in\cU\), \(\cM_{\tilde{a}}^\reg\) is dense in
    \(\cM_{\tilde{a}}\) and is a \(\cP_{\tilde{a}}\)-torsor. \'Etale-locally we
    may trivialize this torsor hence have an open embedding \(\cP\to\cM\). This
    embedding identifies \(\RDF^{-\dim{\cM}+2d}h_*\bcQ\) with
    \(\RDF^{2d}p_!\Qlb\) compatible with \(\pi_0(\cP)\)-actions up to Tate
    twist. If we only consider the stable constituent, then by homotopy lemma
    (see \Cref{lem:homotopy_lemma}), this description is
    independent of the choice of \(\cP\to\cM\) hence descends to the entire
    \(\cU\) as desired.
\end{proof}

\begin{remark}
    \begin{enumerate}
        \item Note that \(\cU\) already contains multiple \(\delta\)-critical
            loci: there is \(\cU\cap\cA_0\) (where \(\delta=0\)), but there can also be other strata
            with \(\delta>0\), and a general point therein is described in
            \Cref{sub:delta_critical_local_delta}.
        \item It is possible to prove that if \(\delta_a=0\), then \(\cM_a\) is
            a \(\cP_a\)-torsor. In other words, the open subset \(\cU\) in
            the above lemma may be enlarged to include the entire \(\delta=0\)
            stratum \(\cA_0\). It is an improvement over
            \Cref{cor:central_lambda_and_d_1_implies_reg}. The proof is
            postponed because it requires the support theorem in
            \Cref{chap:support_theorem}. See
            \Cref{cor:delta_equals_0_implies_torsor}. We are not aware of a
            purely local
            proof at the time of writing.
    \end{enumerate}
\end{remark}

\subsection{}
Now we consider a \(\delta\)-critical stratum outside \(\cU\) above.  Suppose
\(\cV\subset\cA\) is an irreducible \(\delta\)-critical locally closed subset
where the discriminant
divisor and the boundary divisor intersects at exactly one point.
Let \(\delta_\cV\) be the minimum of \(\delta\) on \(\cV\).

Let \(\tilde{a}\in\cV(\bar{k})\) be a sufficiently general point so that
\(\delta_{\tilde{a}}=\delta_\cV\), and \(\bar{v}\in X(\bar{k})\)
is the unique point where the discriminant and boundary divisors intersect. When
\(\tilde{a}\) varies in \(\cV\), \(\bar{v}\) moves in \(\breve{X}\), hence it may be seen
as a family of divisors in \(\breve{X}\x\cV\) over \(\cV\). Then \(\tilde{a}\) is
unramified at \(\bar{v}\) with local
Newton point \(\nu\), and \((\RDF^{-\dim{\cM}+2d}h_*\bcQ)_{\hST}\) has constant rank
given by weight multiplicity \(m_\cV=m_{-w_0(\lambda_{\bar{v}})\nu}\). By
\Cref{prop:delta_critical_in_inductive}, there exists an inductive subset
\(\cU'\) containing \(\cV\), such that
\begin{enumerate}
    \item there exists an open subset \(\cV'\subset\cU'\) containing \(\cV\)
        as a closed subset, and
    \item the Newton point at \(\bar{v}\) stays locally constant on \(\cV'\).
\end{enumerate}
Replacing \(\cV\) by \(\cV'\), 
we may assume that \(\cV\) is inductive. In particular, 
\((\RDF^{-\dim{\cM}+2d}h_*\bcQ)_{\hST}\) still has constant rank
\(m_\cV\).

We would like to show that
\((\RDF^{-\dim{\cM}+2d}h_*\bcQ)_{\hST}|_{\cV}\) is a local system with rank
\(m_{\cV}\). Since \(\cV\) is inductive and all local Newton points stay
locally constant, it is smooth (see \Cref{sec:inductive_strata}), and so
we can always extend a local system over any subset of
codimension \(2\). Thus, it suffices to only consider the subset of \(\cV\) such that
we either have \(\delta=\delta_\cV\) or \(\delta=\delta_\cV+1\). In such case,
we may assume that \(a\) is \(\nu\)-regular semisimple at \(\bar{v}\),
because the collection of points not satisfying this condition has codimension
at least \(2\) (see the proof of \Cref{prop:delta_regularity_fix_divisor})
hence can be safely deleted.

\subsection{}
If we look at product formula \eqref{eqn:product_formula}, it is easy to believe
that the jump in number of irreducible components (modulo \(\pi_0(\cP)\)-action)
is purely a local phenomenon at \(\bar{v}\). The problem is, however,
that \eqref{eqn:product_formula} is only formulated over one point \(\tilde{a}\) at a
time. Therefore, we must upgrade the formula over
Henselian bases. This is achieved using a new kind of
Hecke-type stack which we formulate in the next section, and afterward we shall
continue describing top cohomology using those Hecke stacks.

\section{MH-Hecke Stacks}
\label{sec:mH_Hecke_stacks}

Given mH-fibration \(h_X\colon\cM_X\to\cA_X\), let \(\cH_X\)
    \nomenclature[\(H"cal_X \)]{\(\cH_X\)}{the mH-Hecke stack associated with \(h_X\)}
be the stack whose
\(S\)-points are tuples \((\cL,E_1,\phi_1,E_2,\phi_2,\psi)\) where
\(\cL\in\Bun_{Z_{\FRM}}(S)\), \((\cL,E_i,\phi_i)\in\cM_X(S)\) are two points
mapping to the same point in \(\cA_X(S)\) whose boundary divisor is denoted by
\(\lambda_b\), and \(\psi\) is an isomorphism
\begin{align}
    \psi\colon (E_1,\phi_1)|_{X\x
    S-\lambda_b}\stackrel{\sim}{\longto}(E_2,\phi_2)|_{X\x S-\lambda_b}.
\end{align}
Note that since \(X\) is separated, the mHiggs field \(\phi_2\) (or \(\phi_1\),
but not both) is determined by other data in the tuple, but we still want to
keep both \(\phi_1\) and \(\phi_2\) to make the definition more symmetric. By
its definition \(\cH_X\) fits into the following diagram where the maps are the
obvious ones:
\begin{equation}
    \begin{tikzcd}
        & \cH_X \ar[d, "\olrto{h}"]\ar[ldd, "\olto{h}", swap, bend right]\ar[rdd,
        "\orto{h}", bend left] & \\
        & \cM_X\x_{\cA_X}\cM_X\ar[ld, "\pr_1", swap]\ar[rd, "\pr_2"] &\\
        \cM_X \ar[rd, "h_X", swap] &  & \cM_X \ar[ld, "h_X"] \\
        & \cA_X &
    \end{tikzcd}
\end{equation}

\begin{definition}
    The stack \(\cH_X\) is called the \notion{mH-Hecke stack}\index{mH-!Hecke stack}
    \index{Hecke stack!mH-} associated with
    mH-fibration \(h_X\colon\cM_X\to\cA_X\).
\end{definition}

Similarly, let \(\cH_{\Bun_G}\) be the stack classifying tuples
\((b,E_1,E_2,\psi)\), where \(b\in\cB_X\), \(E_i\in\Bun_G\) and \(\psi\) is
an isomorphism of \(E_1\) and \(E_2\) outside \(\lambda_b\),
then \(\cH_X\) fits into the larger diagram
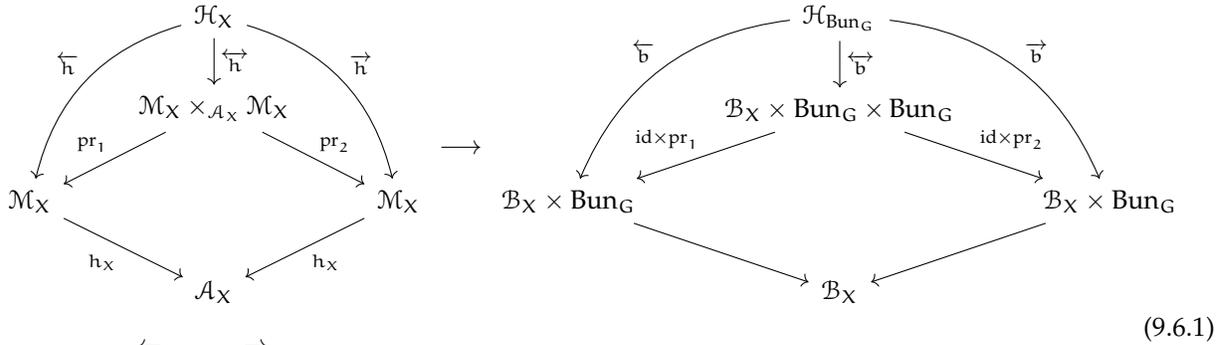
\begin{equation}
    \label{eqn:mH_Hecke_big_diagram}
    \begin{tikzcd}
        & \cH_X \ar[d, "\olrto{h}"]\ar[ldd, "\olto{h}", swap, bend right]\ar[rdd,
        "\orto{h}", bend left] & \\
        & \cM_X\x_{\cA_X}\cM_X\ar[ld, "\pr_1", swap]\ar[rd, "\pr_2"] &\\
        \cM_X \ar[rd, "h_X", swap] &  & \cM_X \ar[ld, "h_X"] \\
        & \cA_X &
    \end{tikzcd}
    \longto
    \begin{tikzcd}
        & \cH_{\Bun_G} \ar[d, "\olrto{b}"]\ar[ldd, "\olto{b}", swap, bend right]\ar[rdd,
        "\orto{b}", bend left] & \\
        & \cB_X\x\Bun_G\x\Bun_G\ar[ld, "\id\x\pr_1", swap]\ar[rd, "\id\x\pr_2"] &\\
        \cB_X\x\Bun_G \ar[rd] &  & \cB_X\x\Bun_G \ar[ld] \\
        & \cB_X &
    \end{tikzcd}
\end{equation}
It is clear that \(\olto{b}\) (resp.~\(\orto{b}\)) is a locally trivial
fibration of affine Grassmannians of \(G\) relative to \(\cB_X\), and that
\begin{align}
    \cH_X \longto \cH_{\Bun_G}\x_{\cB_X\x\Bun_G}\cM_X
\end{align}
induced by \(\olto{h}\) and \(\olto{b}\) (resp.~\(\orto{h}\) and \(\orto{b}\))
is a closed embedding of \(\cM_X\)-functors. Therefore, \(\cH_X\) is an
ind-algebraic stack of ind-finite type that is ind-proper over \(\cM_X\).

\begin{lemma}
    \label[lemma]{lem:fiber_of_h_left_or_h_right}
    Let \((\cL,E,\phi)\in\cM_X^\heartsuit(\bar{k})\) and let
    \(a\in\cA_X^\heartsuit\) be its image and
    \(\lambda_b\) the associated boundary divisor.
    Then the fiber of \(\olto{h}\) (resp.~\(\orto{h}\)) is isomorphic to the
    product of multiplicative affine Springer fibers at the support of
    \(\lambda_b\)
    \begin{align}
        \label{eqn:fiber_of_h_left_or_h_right}
        \cM_{\lambda_b}(a)\defeq\prod_{\bar{v}\in \lambda_b}\cM_{\bar{v}}(a).
    \end{align}
\end{lemma}
\begin{proof}
    The statements for \(\olto{h}\) and for \(\orto{h}\) are the same, so it
    suffices to prove for \(\orto{h}\). Given \((\cL,E,\phi)\), since
    \(\bar{k}\) is algebraically closed, we may choose
    and fix an isomorphism around the formal disc \(\breve{X}_{\lambda_b}\) around
    the support of \(\lambda_b\):
    \begin{align}
        \tau\colon (\cL,E,\phi)\longto (\cL_0,E_0,\gamma_{a,\lambda_b})
    \end{align}
    where \(\cL_0\) (resp.~\(E_0\)) is the trivial \(Z_{\FRM}\)-torsor
    (resp.~\(G\)-torsor), and
    \(\gamma_{a,\lambda_b}\in\FRM(\breve{X}_{\lambda_b})\) such that
    \(\chi_{\FRM}(\gamma_{a,\lambda_b})=a(\breve{X}_{\lambda_b})\).

    For any \(k\)-scheme \(S\), we have the map
    \begin{align}
        \orto{h}^{-1}(\cL,E,\phi)(S)&\longto \cM_{\lambda_b}(a)\\
        (\cL,E_1,\phi_1,\psi)&\longmapsto
        (\cL|_{\breve{X}_{\lambda_b}},E_1|_{\hat{X}_{\lambda_b}},
        \phi_1|_{\breve{X}_{\lambda_b}}, \beta),
    \end{align}
    where \(\beta\) is the composition of maps
    \begin{align}
        \beta\colon
        (\cL,E_1,\phi_1)|_{\breve{X}_{\lambda_b}^\bullet}
        \stackrel{\psi}{\longto}(\cL,E,\phi)|_{\breve{X}_{\lambda_b}^\bullet}
        \stackrel{\tau}{\longto}(\cL_0,E_0,\gamma_a)|_{\breve{X}_{\lambda_b}^\bullet},
    \end{align}
    and \(\breve{X}_{\lambda_b}^\bullet\) is the punctured disc. The map \(\beta\)
    is clearly injective: if \((\cL,E_1,\phi_1,\psi)\) and
    \((\cL,E_1',\phi_1',\psi')\) have isomorphic image under \(\beta\), they are
    isomorphic over both \(\breve{X}\x S-\lambda_b\) and \(\breve{X}_{\lambda_b}\),
    together with their gluing data. It implies they are isomorphic tuples in
    \(\orto{h}^{-1}(\cL,E,\phi)(S)\). On the other hand, \(\beta\) is also
    surjective, because any point in \(\cM_{\lambda_b}(a)\) can be glued with
    \((\cL,E,\phi)|_{\breve{X}\x S-\lambda_b}\) to obtain a point in
    \(\orto{h}^{-1}(\cL,E,\phi)\). This finishes the proof.
\end{proof}

\subsection{}
There are some useful variants of mH-Hecke stacks. Firstly, the mH-base can
be replaced by any algebraic \(\cA_X\)-stack \(\cU\to\cA_X\) and \(\cM_X\) by
its pullback to \(\cU\). Secondly, the boundary divisor \(\lambda_b\) can be
replaced by any finite flat family \(\cD\) of Cartier divisors in
\(X\x\cU\) over \(\cU\), so that the rational map \(\psi\) in the definition is
now an isomorphism outside \(\cD\). Let \(\cH_{\cD}\) be the
corresponding Hecke stack, together with maps \(\olto{h}_\cD\),
\(\orto{h}_\cD\), etc.
\begin{definition}
    Given a tuple \((h_X,\cU,\cD)\) as above, we call \(\cH_\cD\)
    \nomenclature[\(H"cal_D"cal \)]{\(\cH_\cD\)}{the generalized mH-Hecke stack
    associated with a family of divisors \(\cD\) on \(X\) parametrized by an \(\cA_X\)-stack \(\cU\to\cA_X\)}
    the
    \notion{generalized mH-Hecke stack}\index{Hecke stack!generalized mH-}
    \index{mH-!Hecke stack, generalized}associated with \((h_X,\cU,\cD)\).
\end{definition}

The representability of \(\cH_\cD\) can be seen using a similar diagram as
\eqref{eqn:mH_Hecke_big_diagram}, with \(\cB_X\) replaced by an appropriate
Hilbert scheme of \(X\).

\begin{example}
    Recall in \Cref{sec:factorizations} we have the space \(\BD[I]_X\)
    (\(I=\Set{1,2}\)), and let \(\cU=\cA_X^{I,\disj}\), so that the boundary
    divisor is a disjoint union of \(\lambda_1\) and \(\lambda_2\). We
    may let \(\cD=\lambda_1\), and the resulting Hecke stack \(\cH_\cD\)
    is called the \notion{partial mH-Hecke stack}\index{mH-!Hecke stack, partial}
    \index{Hecke stack!partial mH-} associated with \((h_X,\cU,\cD)\).
\end{example}

\begin{example}
    We can let \(\cU=\cA_X^\heartsuit\) and let \(\cD\) be the discriminant
    divisor \(\FRD_X\). In this case we denote the resulting Hecke stack by
    \(\cH_{\FRD}\),
    \nomenclature[\(H"cal_D"frak \)]{\(\cH_\FRD\)}{the \(\FRD\)-Hecke stack,
    i.e., \(\cH_\cD\) where \(\cU=\cA_X^\heartsuit\) and \(\cD=\FRD_X\)}
    and call it the \notion{\(\FRD\)-Hecke stack}\index{Hecke stack!\(\FRD\)-}.
\end{example}

If \(\cD'\to\cD\) is a \(\cU\)-morphism of divisor families,
then we have natural maps \(\cH_{\cD'}\to\cH_\cD\).
Using the same argument as in \Cref{lem:fiber_of_h_left_or_h_right},
we have the following result:
\begin{lemma}
    Let \((\cL,E,\phi)\in\cM_X^\heartsuit|_\cU(\bar{k})\) and \(a\in \cU\) be
    its image. Let \(\cD_a\) the associated Cartier divisor induced by
    \(\cD\). The fiber of \(\olto{h}_{\cD}\)
    (resp.~\(\orto{h}_{\cD}\)) over \((\cL,E,\phi)\) is isomorphic to
    \begin{align}
        \cM_{\cD_a}\defeq\prod_{\bar{v}\in\cD_a}\cM_{\bar{v}}(a).
    \end{align}
    Moreover, if \(\cD'\to\cD\) are two divisor families,
    then this isomorphism is compatible with the natural map
    \(\olto{h}_{\cD'}^{-1}(\cL,E,\phi)\to\olto{h}_{\cD}^{-1}(\cL,E,\phi)\)
    and the map induced by \((E,\phi)\):
    \begin{align}
        \cM_{\cD_a'}\longto\cM_{\cD_a},
    \end{align}
    and similarly for the \(\orto{h}\) side.
\end{lemma}

\subsection{}
The symmetry of mH-Hecke stacks can be described using regular centralizer just
like mH-fibrations. The action of Picard stack \(\cP_X\) can be pulled back to
\(\cM_X\) to give an action of \(\cM_X\x_{\cA_X}\cP_X\) on
\(\cM_X\x_{\cA_X}\cM_X\), relative to the first projection to \(\cM_X\).
This can also be achieved by pulling back regular centralizer \(\FRJ_X\to X\x
\cA_X\) to \(X\x \cM_X\), and form the relative Picard stack over \(\cM_X\).

Suppose for now  \(\FRM=\Env(G^\SC)\), then over \(\cB_X\) there is a finite
flat family of Cartier divisors in \(X\), namely the numerical boundary divisor
\(\FRE_X\). Since the regular centralizer is
affine smooth, we may use the same construction as in
\Cref{sec:review_BD_affine_grassmannian}, and define the
relative affine Grassmannian
\begin{align}
    \cP_{\cB_X}\defeq\Gr_{\FRJ_X,\cB_X}\longto \cA_X
\end{align}
whose fiber at \(a\in\cA_X\) is exactly the product of local Picard group
\begin{align}
    \cP_{\cB_X,a}=\prod_{\bar{v}\in\FRE_a}\cP_{\bar{v}}(a).
\end{align}

If \(\FRM\neq\Env(G^\SC)\), we have two subtlely different constructions:
we may simply pull back the construction for
\(\Env(G^\SC)\) to any \(\FRM\in\FM(G^\SC)\), and if \(\FRA_{\FRM}\) is of standard type,
we may also repeat the same construction using \(\FRE_\FRM\) in place of
\(\FRE_{\Env(G^\SC)}\). These two constructions are not the same if
the generating cocharacters of \(\FRA_\FRM\) contains a central one, in which
case the underlying topological space of \(\FRE_{\FRM}\) strictly contains
that of the pullback of \(\FRE_{\Env(G^\SC)}\). For our purposes in this
book, such difference is not very important, and for convenience, we will only
consider the case where \(\FRA_\FRM\) is of standard type and use
the second construction (the one involving \(\FRE_\FRM\)).

We have the forgetful map \(\cP_{\cB_X}\to\cP_X\) by forgetting the
trivialization of \(\FRJ_{a}\)-torsor over \(X-\FRE_{a}\), so \(\cP_{\cB_X}\)
naturally acts on \(\cM_X\). We claim that \(\cP_{\cB_X}\), after pulling back
to \(\cM_X\), acts on \(\olto{h}\colon\cH_X\to\cM_X\), making \(\olrto{h}\) a
\(\cP_{\cB_X}\)-equivariant map, and similarly for \(\orto{h}\). Indeed, suppose
we have tuple \((\cL,E_1,\phi_1,E_2,\phi_2,\psi)\in\cH_X(S)\) over \(a\in\cA_X(S)\)
and \((E_{\FRJ}, \tau)\in\cP_{\cB_X}(S)\) where \(E_{\FRJ}\) is a \(\FRJ_{a}\)-torsor over \(X\x
S\) and \(\tau\) is a trivialization of \(E_{\FRJ}\) outside \(\FRE_{a}\). The
action \(E_{\FRJ}\) on \(\cM_X\) sends \((\cL,E_2,\phi_2)\) to
\begin{align}
    \phi_2'\colon E_2'\defeq E_2\x_{\phi_2,\FRM_\cL}^{\FRJ_{a}}
    E_{\FRJ}\longto\FRM_\cL,
\end{align}
where the action of \(\FRJ_{a}\) on \(\phi_2\) is induced by the canonical map
\(\chi_{\FRM}^*\FRJ_{\FRM}\to I_{\FRM}\). The trivialization \(\tau\) induces
isomorphism
\begin{align}
    (E_2,\phi_2)|_{X\x S-\FRE_{a}}\stackrel{\sim}{\longto}
    (E_2',\phi_2')|_{X\x S-\FRE_{a}},
\end{align}
whose composition with \(\psi\) gives
\begin{align}
    \psi'\colon (E_1,\phi_1)|_{X\x S-\FRE_{a}}\stackrel{\sim}{\longto}(E_2',\phi_2')|_{X\x
    S-\FRE_{a}}.
\end{align}
This defines the \(\cP_{\cB_X}\)-action, and clearly it makes \(\olrto{h}\)
equivariant. The argument for \(\orto{h}\) is the same.

The story for generalized mH-Hecke stacks \(\cH_\cD\) associated with tuple
\((h_X,\cU,\cD)\) is the same, except one replaces \(\FRE_{\FRM}\)
by \(\cD\). For future convenience we denote the local Picard group in
this case by \(\cP_{\cD}\) in place of \(\cP_{\cB_X}\). We leave other
details to the reader.

\begin{lemma}
    Let \((\cL,E,\phi)\in\cM_X^\heartsuit(\bar{k})\) and
    \(a\in\cA_X^\heartsuit(\bar{k})\) be its image. Let
    \(\cH_{(\cL,E,\phi)}=\olto{h}^{-1}(\cL,E,\phi)\).
    The action of \(\cP_{\cB_X}\) on \(\cH_X\) induces a bijection between sets
    \begin{align}
        \Cnt\Irr\Stack*{\cH_{(\cL,E,\phi)}/\cP_{\cB_X,(\cL,E,\phi)}}
        \stackrel{\sim}{\longto}\Cnt\Irr\Stack*{\cM_a/\cP_a}.
    \end{align}
\end{lemma}
\begin{proof}
    By \Cref{lem:fiber_of_h_left_or_h_right}, the map \(\orto{h}\colon\cH_{(\cL,E,\phi)}\to \cM_a\)
    induces maps
    \begin{align}
        \Stack*{\cH_{(\cL,E,\phi)}/\cP_{\cB_X,(\cL,E,\phi)}}
        \stackrel{\sim}{\longto} \prod_{\bar{v}\in
        \lambda_b}\Stack*{\cM_{\bar{v}}(a)/\cP_{\bar{v}}(a)}\longto
        \Stack*{\cM_a/\cP_a}.
    \end{align}
    Note that when \(\bar{v}\not\in\FRD_a\), the stack
    \(\Stack*{\cM_{\bar{v}}(a)/\cP_{\bar{v}}(a)}\) is just a \(\bar{k}\)-point,
    and according to \Cref{thm:local_irr_components_weight_mult},
    if \(\bar{v}\in X\) does not support the boundary
    divisor, then \(\cM_{\bar{v}}^\reg(a)\) is dense in \(\cM_{\bar{v}}(a)\)
    and is a \(\cP_{\bar{v}}(a)\)-torsor. Combining these facts with
    product formula \Cref{prop:product_formula_alt} we obtain the
    lemma.
\end{proof}

\subsection{}
We may replace the mH-Hecke stack by \(\FRD\)-Hecke stack \(\cH_{\FRD}\),
and obtain a family of maps
\begin{equation}
    \begin{tikzcd}
        \cH_{\FRD}\x^{\cP_{\FRD}}\cP_X^\heartsuit \ar[r, "\olto{\Pi}_{\FRD}"]
        \ar[rd, "\olto{h}_{\FRD}\circ\pr_1", swap] &
        \cM_X^\heartsuit\x_{\cA_X^\heartsuit}\cM_X^\heartsuit \ar[d, "\pr_1"]\\
          & \cM_X^\heartsuit
    \end{tikzcd}
\end{equation}
where \(\olto{\Pi}_{\FRD}\) is the map
\begin{align}
    ((\cL,E_1,\phi_1,E_2,\phi_2,\psi),p\in\cP_X^\heartsuit)\longmapsto (\cL,E_1,\phi_1,p\cdot
    (E_2,\phi_2)).
\end{align}
Over any \((\cL,E,\phi)\) whose image is \(a\in\cA_X^\heartsuit(\bar{k})\), the
fiber of \(\olto{\Pi}_{\FRD}\) is clearly isomorphic to the (non-reduced)
product formula \eqref{eqn:product_formula_non_reduced}, so its
reduced version is isomorphic to \eqref{eqn:product_formula}.
If \((\cL,E,\phi)\in\cM_X^\reg(\bar{k})\), since
\(\cP_X\) acts on \(\cM_X^\reg\) freely, \(\cM_X^\reg\) is smooth over
\(\cA_X\), then we choose a section \(\cA_X\to\cM_X^\reg\) in a neighborhood of
\(a\), and pull back \(\olto{\Pi}_{\FRD}\). If a Steinberg
quasi-section exists on an irreducible component of \(\cA_X^\heartsuit\), we may
even do this over the entire said component. There is a symmetric
construction for \(\orto{h}_{\FRD}\) as well.

\subsection{}
The construction of map \(\olto{\Pi}_{\FRD}\) (and \(\orto{\Pi}_{\FRD}\)) works
for any generalized mH-Hecke stack \(\cH_\cD\), so we have morphism of stacks
\begin{equation}
    \begin{tikzcd}
        \cH_\cD\x^{\cP_{\cD}}\cP_\cU \ar[r, "\olto{\Pi}_\cD"]
        \ar[rd, "\olto{h}_\cD\circ\pr_1", swap] &
        \cM_\cU\x_{\cU}\cM_\cU\ar[d, "\pr_1"]\\
          & \cM_\cU
    \end{tikzcd}
\end{equation}
The quotient of \(\cH_\cD\) by \(\cP_{\cD}\) is an ind-algebraic
stack of ind-finite type over \(\cM_\cU\) whose geometric fibers are proper
algebraic stacks of finite type. In fact, since the reduced geometric fibers of
\(\olto{h}_\cD\) are schemes locally of finite type, we can, locally over
\(\cM_\cU\), find some affine open subset of \(\cH_\cD\) that maps surjectively
onto the quotient \(\Stack*{\cH_\cD^\Red/\cP_{\cD}^\Red}\). For
example, since \(\cH_\cD\) embeds into a locally constant affine Grassmannian
over \(\cM_\cU\), we can take a sufficiently large truncation in the affine
Grassmannian and take the preimage in \(\cH_\cD\). Hence, the stack
\begin{align}
    \label{eqn:general_simultaneous_product_formula_lhs_stack}
    \cH_\cD^\Red\x^{\cP_{\cD}^\Red}\cP_\cU
\end{align}
is algebraic over \(\cM_\cU\), locally of finite type (it is not necessarily of
finite type because \(\pi_0(\cP_\cU)\) may be infinite). When \(\cU\to\cA_X\)
has its image contained in \(\cA_X^\ANI\), then
\eqref{eqn:general_simultaneous_product_formula_lhs_stack} is a proper Deligne--Mumford
stack of finite type over \(\cU\).

\begin{proposition}[Simultaneous Product Formula]
    \label[proposition]{prop:simultaneous_product_formula}
    The induced map \(\olto{\Pi}_{\FRD}^\Red\) on respective reduced stacks
    is a universal homeomorphism over \(\cM_X^\heartsuit\).
\end{proposition}
\begin{proof}
    The same proof as in \Cref{prop:product_formula} shows that the fiber of
    \(\olto{\Pi}_{\FRD}^\Red\) over any field-valued point of \(\cM_X^\heartsuit\)
    is a universal homeomorphism and the result follows.
\end{proof}

\subsection{}
We return to the situation at the end of \Cref{sec:top_ordinary_cohomology}
and continue describing \((\RDF^{-\dim{\cM}+2d}h_*\bcQ)_\hST\).
Over \(\cV\), choose the divisor \(\cD=\bar{v}\), and consider the
generalized mH-Hecke stack \(\cH_{\bar{v}}\) associated with
\((h_X,\cV,\bar{v})\). This is a partial mH-Hecke stack with boundary
subdivisor \(\lambda_{\bar{v}}\). Recall the Newton point \(\nu\) at \(\bar{v}\)
is locally constant on \(\cV\), and the local \(\delta\)-invariant at
\(\bar{v}\) is exactly
\(\Pair{\rho}{\lambda_{\bar{v}}-\nu}\), which is constant over \(\cV\).

Since we want to prove that the stable top cohomology is a local system on
\(\cV\), we may replace \(\cV\) by any strict Henselian neighborhood therein.
Since \(\cM_{\cV}^\reg\) is smooth over \(\cV\) and \(\cV\) is strictly
Henselian, we have a section
\(\tau\) of \(\cV\) in \(\cM_{\cV}\). Such section induces a tuple
\((\cL,E,\phi)\in\cM_{\cV}(\cV)\). We
may trivialize \(E\) at \(\bar{v}\) over \(\cV\), and by smoothness lift it
to the formal disc \(\breve{X}_{\bar{v}}\) around \(\bar{v}\). Using the same
argument as in \Cref{lem:fiber_of_h_left_or_h_right}, we have that
\begin{align}
    \tau^*\cH_{\bar{v}}\cong\cV\x\cM_{\bar{v}}(a_0),
\end{align}
where \(a_0\) is the unique closed point in \(\cV\). Since by assumption
\(\bar{v}\) is the only point supporting both boundary and discriminant
divisors, \(\cM_{\bar{v}'}^\reg(a)\) is dense in \(\cM_{\bar{v}'}(a)\) for any
\(a\in\cV\) and \(\bar{v}'\neq\bar{v}\), and so the
\(\cP_{\cV}\)-equivariant map
\begin{align}
    \tau^*\left(\cH_{\bar{v}}\x^{\cP_{\bar{v}}}\cP_{\cV}\right)\longto
    \tau^*\left(\cM_{\cV}\x_{\cV}\cM_{\cV}\right)\cong\cM_{\cV}
\end{align}
has fiberwise dense image over \(\cV\). Since \(\cV\) is \(\delta\)-critical by
assumption, the boundary divisor must be multiplicity-free. This shows that
\((\RDF^{-\dim{\cM}+2d}h_*\bcQ)_\hST|_{\cV}\) is local system of rank
\(m_\cV\) over \(\cV\). The whole
argument clearly generalizes to more points than \(\bar{v}\), as well as when
\(\bcQ\) is a general Satake sheaf. Thus, we have the following result:
\begin{proposition}
    \label[proposition]{prop:stable_top_cohomology_local_system_main}
    Let \(\bcQ\in\Sat_X\).
    Suppose locally closed subset \(\cV\subset\tilde{\cA}_X^\ANI\) is such that:
    \begin{enumerate}
        \item There exists an open subset \(\cU\subset \tilde{\cA}_X^\ANI\)
            containing \(\cV\) over which the local model of singularity as in
            either \Cref{thm:local_singularity_model_weak} or
            \Cref{thm:local_singularity_model_main} exists.
        \item \(\cV\) is \(\delta\)-critical with multiplicity-free boundary
            divisor.
        \item For any \(a\in\cV\), \(a\) is unramified at all points supporting
            the boundary divisor.
        \item If locally over \(\cV\) we write boundary divisor  as
            \(\lambda_b=\sum_{i=1}^m\lambda_i\cdot\bar{v}_i\), then the Newton
            point \(\nu_i\) at each \(\bar{v}_i\) is locally constant over
            \(\cV\).
    \end{enumerate}
    Then \((\RDF^{-\dim{\cM}+2d}h_*\bcQ)_\hST|_{\cV}\) is a local system of 
    rank
    \begin{align}
        \prod_{i=1}^m
        m_{-w_0(\lambda_i)\nu_i}\dim_{\Qlb}\Hom_{\dual{\bG}}(V_{-w_0(\lambda_i)},V_{\bar{v}_i}),
    \end{align}
    where \(V_{\bar{v}_i}\) is the \(\dual{\bG}\)-representation at \(\bar{v}_i\) corresponding
    to \(\bcQ\).
\end{proposition}

\subsection{}
The same argument for \Cref{prop:stable_top_cohomology_local_system_main} can be
applied to \(\kappa\)-isotypic constituent as well:
\begin{proposition}
    \label[proposition]{prop:kappa_top_cohomology_local_system_main}
    For a fixed \(\kappa\), suppose locally closed substack
    \(\cV\subset\tilde{\cA}_\kappa^\ANI\) satisfies all the conditions in
    \Cref{prop:stable_top_cohomology_local_system_main}, then
    \((\RDF^{-\dim{\cM}+2d}h_*\bcQ)_\kappa|_{\cV}\) is a local system
    of rank
    \begin{align}
        \prod_{i=1}^m m_{-w_0(\lambda_i)\nu_i}\dim_{\Qlb}\Hom_{\dual{\bG}}(V_{-w_0(\lambda_i)},V_{\bar{v}_i}),
    \end{align}
    where the notations are as in
    \Cref{prop:stable_top_cohomology_local_system_main}.
\end{proposition}

\chapter{Support Theorem}%
\label{chap:support_theorem}

In this section we prove a slightly generalized version of the 
Support Theorem in \cite{Ng10}*{\S~7.2}. The method here follows the
outline in \cite{Ng10}*{\S\S~7.3--7.7}, with some modifications.

We then apply our abstract support theorem to mH-fibrations. One key input for
this application is the local model of singularity
(\Cref{thm:local_singularity_model_weak,thm:local_singularity_model_main}),
which will provide us with the bound
on cohomological amplitude that is only assumed in the abstract support theorem.

\section{Abelian Fibrations}%
\label{sec:abelian_fibrations}

Let \(f\colon M\to S\) be a proper map of varieties over a finite 
field \(k\). Let \(g\colon P\to S\) be a smooth commutative
group scheme over \(S\). Suppose
\(P\) acts on \(M\) relative to \(S\) and the stabilizers are affine. Let
\(P^0\subset P\) be the open group subscheme such that for any geometric point
\(s\in S\), \(P_s^0\) is the neutral component of \(P_s\). We then have the
canonical short exact sequence of Chevalley
\begin{align}
    \label{eqn:Chevalley_exact_sequence}
    1\longto R_s\longto P_s^0\longto A_s\longto 1,
\end{align}
where \(R_s\) is connected and affine and \(A_s\) is an abelian variety.
It induces a decomposition of Tate modules
\begin{align}
    \label{eqn:Chevalley_Tate_module}
    0\longto \TateM_{\Qlb}(R_s)\longto \TateM_{\Qlb}(P_s^0) \longto \TateM_{\Qlb}(A_s)\longto 0
\end{align}

We have an \(\bbN\)-valued function \(\delta(s)=\dim{R_s}\) defined for the 
topological points of \(S\), which is necessarily upper-semicontinuous (see
\cite{Ng10}*{\S~5.6.2}). Suppose \(\delta\) is constructible, then it induces a
locally closed stratification
\begin{align}
    S=\coprod_{\delta\in\bbN} S^{\delta},
\end{align}
so that if \(s\in S^{\delta}\), then \(\delta(s)=\delta\).

\begin{definition}
    \label[definition]{def:weak_abelian_fibration}
    We call \((f,g)\) a \notion{weak abelian fibration}\index{fibration!weak abelian}
    if the following conditions are satisfied:
    \begin{enumerate}
        \item \(f\) and \(g\) have the same relative dimension \(d\).
        \item For any geometric point \(s\in S\) and any \(m\in M\), its
            stabilizer in \(P_s\) is affine.
        \item The Tate module \(\TateM_{\Qlb}(P^0)\) is polarizable. In other words,
            there exists \'etale locally over \(S\) an alternating bilinear
            form on \(\TateM_{\Qlb}(P^0)\), such that for any \(s\) its
            restriction to \(\TateM_{\Qlb}(R_s)\) is zero, and it 
            induces a perfect
            pairing of \(\TateM_{\Qlb}(A_s)\) with itself.
    \end{enumerate}
\end{definition}

\begin{definition}
    \label[definition]{def:delta_regular_abelian_fibration}
    We call \((f,g)\) a \notion{\(\delta\)-regular abelian fibration}
    \index{\(\delta\)-!regularity}\index{fibration!\(\delta\)-regular abelian} if it is a
    weak abelian fibration, and for any \(\delta\in\bbN\), we have
    \(\codim_S(S^{\delta})\ge \delta\). (If \(S^{\delta}=\emptyset\), then the
    codimension is \(\infty\) by convention.) Equivalently, we have for any
    irreducible closed subset \(Z\subset S\), \(\codim_S(Z)\ge \delta_Z\) where
    \(\delta_Z\) is the minimum of \(\delta\) on \(Z\).
\end{definition}

\begin{remark}
    \label[remark]{rmk:abelian_fibrations_for_DM_stacks}
    Note that both
    \Cref{def:weak_abelian_fibration,def:delta_regular_abelian_fibration} make
    sense if we replace schemes with Deligne--Mumford stacks.
\end{remark}

\section{Goresky--MacPherson Inequality}%
\label{sec:goresky_macpherson_inequality}

Let \((f,g)\) be a weak abelian fibration.
Let \(\cF\in \Dbc(M,\Qlb)\) be a self-dual complex hence of pure weight \(0\), so
\(f_*\cF\in\Dbc(S,\Qlb)\) is also 
of pure weight \(0\) since \(f\) is proper. Thus, we
have (non-canonical) decomposition of Frobenius modules by \cite{BBD82}:
\begin{align}\label{eqn:BBD_decomp}
    f_*\cF\cong \bigoplus_{n\in\bbZ}\PH^n(f_*\cF)[-n].
\end{align}
Given an irreducible closed subset \(Z\) of \(S\), let \(\occ(Z)\subset\bbZ\) be
the set of numbers \(n\) such that \(Z\) appears in the set of supports
\(\supp{\PH^n(f_*\cF)}\). By Poincar\'e duality, \(\occ(Z)\) is symmetric about
\(0\). Suppose \(N\) is the largest number such that \(\RH^N(f_*\cF)\neq 0\).
Suppose \(Z\neq \emptyset\) and let \(0\le n\in\occ(Z)\). Then there exists an
open subset \(U\subset S\) such that \(U\cap Z\neq\emptyset\), and a local
system \(L\) on \(U\cap Z\), such that \(i_*L[\dim{Z}]\) (\(i\) being the map
\(U\cap Z\hookto Z\)) is a direct summand of
\(\PH^n(f_*\cF)|_U\), hence also a direct summand of \(f_*\cF[n]\). Taking the
usual cohomology, one has that \(i_*L\) is a direct summand of
\(\RH^{n-\dim{Z}}(f_*\cF)\). Therefore, \(n-\dim{Z}\le N\). Since \(n\ge 0\), we
obtain the Goresky--MacPherson inequality:
\begin{align}
    \codim_S(Z)\le \dim{S}+N-n\le \dim{S}+N.
\end{align}
In particular, if \(N\le -\dim{M}+2d=-\dim{S}+d\), then \(\codim_S(Z)\le d\).
Suppose further we have that \(n\) can be so chosen that \(n\ge (d-\delta_Z)\),
then we have an improved inequality \(\codim_S(Z)\le\delta_Z\). If equality
holds (e.g., it happens when \((f,g)\) is
\(\delta\)-regular), all the inequalities
just mentioned are equalities, and in particular \(N=-\dim{M}+2d\), and
the restriction of \(\PH^{n-\dim{Z}}(f_*\cF)\) to \(U\cap Z\) is a direct
summand of the top cohomology \(\RDF^{-\dim{M}+2d}f_*\cF\).

\section{Action by Cap Product}%
\label{sec:action_by_cap_product}

Let \((f,g)\) be a weak abelian fibration. Following
\cites{LaOl08I, LaOl08II}, we have the derived category \(\DCat(\Stack{M/P},\Qlb)\)
of quotient stack \(\Stack{M/P}\), which
we also define to be the \(P\)-equivariant derived category \(\DCat_P(M,\Qlb)\)
on \(M\) by pulling back through map \(q\colon M\to\Stack{M/P}\). 
We also have the full
subcategories of various boundedness and constructibility conditions\footnote{In
\cites{LaOl08I, LaOl08II}, they have two variants of subcategories for each
boundedness condition, but the distinction is not very important here.}.
Note that if a complex is \(P\)-equivariant, then it is also \(P^0\)
equivariant. For the rest of this section we are going to replace \(P\) with
\(P^0\) so that \(P\) is fiberwise connected.

Following \cite{Ng10}*{\S~7.4}, we define a complex in \(\Dbc(S,\Qlb)\)
\begin{align}
    \Lambda_P\defeq g_!\Qlb[2d](d),
\end{align}
which is concentrated in non-positive degrees and whose degree \(-1\) cohomology
is the sheaf \(\TateM_{\Qlb}(P)\) such that at any geometric point \(s\in S\) its
stalk is the Tate module \(\TateM_{\Qlb}(P_s)\). We have in fact canonical
isomorphisms in \(\Dbc(S,\Qlb)\)
\begin{align}\label{eqn:Lambda_P_decomp}
    \Lambda_P\simeq \bigoplus_{i\ge 0}\RH^{-i}(\Lambda_P)[i]\simeq \bigoplus_{i\ge
    0}\wedge^i\TateM_{\Qlb}(P)[i],
\end{align}
making it a graded algebra.

\subsection{}
More generally, suppose \(\cF\in\Dbc(M,\Qlb)\) is \(P\)-equivariant, in other
words, \(\cF\simeq q^*\bar{\cF}\) for some \(\bar{\cF}\) in
\(\Dbc(\Stack{M/P},\Qlb)\). Then the action morphism
\begin{align}
    a\colon P\x_S M\longto M
\end{align}
is smooth and of relative dimension \(d\). Therefore, we have by adjunction a
morphism of complexes
\begin{align}
    a_!a^*\cF[2d](d)\longto \cF.
\end{align}
Further pushing forward using \(f_!\), we obtain the morphism
\begin{align}
    (g\x_S f)_!a^*\cF[2d](d) \longto f_!\cF.
\end{align}
On the other hand, using the Cartesian diagram
\begin{equation}
    \begin{tikzcd}
        P\x_S M\ar[r, "a"]\ar[d,"p_2"] & M\ar[d, "q"]\\
        M\ar[r, "q"] & \Stack{M/P}
    \end{tikzcd},
\end{equation}
we see that 
\begin{align}
    \label{eqn:canonical_identification_for_equivariant_complex}
    a^*\cF\simeq a^*q^*\bar{\cF}\simeq p_2^*q^*\bar{\cF}\simeq
    \Qlb\boxtimes\cF.
\end{align}
So by K\"unneth formula, we have a morphism of cap products
\begin{align}\label{eqn:cap_product}
    \Lambda_P\otimes f_!\cF\longto f_!\cF.
\end{align}

\subsection{}

Now let \(\cF\in\Dbc(M,\Qlb)\) be a \(P\)-equivariant complex of pure weight
\(0\). Since \(f\) is proper, we have \(f_!\cF\simeq f_*\cF\) and it is also pure
of weight \(0\). 
For each \(n\in\bbZ\) we have a canonical isomorphism of Frobenius
modules
\begin{align}\label{eqn:support_decomp}
    \PH^n(f_*\cF)\simeq \bigoplus_{\alpha\in\PSP}K_\alpha^n,
\end{align}
where \(\PSP\) is the index set of supports of the perverse cohomologies of
\(f_*\cF\). For \(\alpha\in\PSP\), we denote by \(Z_\alpha\) corresponding
the irreducible closed subset in \(S_{\bar{k}}\)
and \(K_\alpha^n\) the perverse summand supported on
\(Z_\alpha\). Since \(\cF\) is bounded and \(f\) is of finite type, 
\(\PSP\) is necessarily finite.

By \eqref{eqn:Lambda_P_decomp} and \eqref{eqn:cap_product}, we have a morphism
\begin{align}
    \TateM_{\Qlb}(P)\otimes f_!\cF\longto f_!\cF[-1].
\end{align}
Composing with perverse truncation \(\Ptr^{\le n}(f_!\cF)\to f_!\cF\), we have
the induced map
\begin{align}
    \TateM_{\Qlb}(P)\otimes \Ptr^{\le n}(f_!\cF)\longto f_!\cF[-1].
\end{align}
Applying functor \({\PH^n}\), we have
\begin{align}
    \PH^n\left(\TateM_{\Qlb}(P)\otimes \Ptr^{\le n}(f_!\cF)\right)\longto \PH^{n-1}(f_!\cF).
\end{align}
Since tensoring with \(\TateM_{\Qlb}(P)\) is perverse right-exact, we know that
\begin{align}
    \TateM_{\Qlb}(P)\otimes \Ptr^{\le n-1}(f_!\cF)\in
    \prescript{\FRp}{}{\DCat}_{\mathup{c}}^{\le n-1}(S,\Qlb),
\end{align}
hence by tensoring with \(\TateM_{\Qlb}(P)\) and then taking the \(n\)-th
perverse cohomology, 
the exact triangle
\begin{align}\label{eqn:perverse_truncation_exact_triangle}
    \Ptr^{\le n-1}(f_!\cF)\longto \Ptr^{\le n}(f_!\cF) \longto
\PH^n(f_!\cF)[-n]\stackrel{+1}{\longto}
\end{align}
 induces an isomorphism
\begin{align}
    \PH^n\left(\TateM_{\Qlb}(P)\otimes \Ptr^{\le
    n}(f_!\cF)\right)\stackrel{\sim}{\longto}
    \PH^{0}\left(\TateM_{\Qlb}(P)\otimes \PH^{n}(f_!\cF)\right).
\end{align}
From this we have a map
\begin{align}
    \PH^{0}\left(\TateM_{\Qlb}(P)\otimes \PH^{n}(f_!\cF)\right)\longto
    \PH^{n-1}(f_!\cF).
\end{align}
Since \(\TateM_{\Qlb}(P)\otimes \PH^{n}(f_!\cF)\in 
\prescript{\FRp}{}{\DCat}_{\mathup{c}}^{\le 0}(S,\Qlb)\), it projects to its
\(0\)-th perverse cohomology, hence we have a canonical map
\begin{align}
    \TateM_{\Qlb}(P)\otimes \PH^{n}(f_!\cF)\longto \PH^{n-1}(f_!\cF).
\end{align}

The canonical decomposition \eqref{eqn:support_decomp} gives map
\begin{align}
    \bigoplus_{\alpha\in\PSP}\TateM_{\Qlb}(P)\otimes K_\alpha^n\longto
    \bigoplus_{\alpha\in\PSP}K_\alpha^{n-1},
\end{align}
and in particular a canonical map
\begin{align}\label{eqn:action_on_K_alpha}
    \TateM_{\Qlb}(P)\otimes K_\alpha^n\longto K_\alpha^{n-1}
\end{align}
for each \(\alpha\in\PSP\) and \(n\).

\section{Statement of Freeness}%
\label{sec:statement_of_freeness}

Following \cite{Ng10}*{\S\S~7.4.8--7.4.9}, for each \(\alpha\in\PSP\), we may find a
dense open subset \(V_\alpha\), such that \(K_\alpha^n\) can be expressed as
\(\cK_\alpha^n[\dim{V_\alpha}]\) for some local system
\(\cK_\alpha^n\) of weight \(n\) on \(V_\alpha\). One can also so choose
\(V_\alpha\) that there is a finite radical base change
\(V'_\alpha\to V_\alpha\) over which the Chevalley exact sequence 
exists:
\begin{align}
    1\longto R_\alpha\longto P|_{V'_\alpha}\longto A_\alpha\longto 1,
\end{align}
where \(R_\alpha\) is a smooth, fiberwise connected affine group scheme over
\(V'_\alpha\), and \(A_\alpha\) is an abelian scheme over \(V'_\alpha\).
It then induces short exact sequence of sheaves
\begin{align}
    0\longto \TateM_{\Qlb}(R_\alpha)\longto
    \TateM_{\Qlb}(P|_{V'_\alpha})\longto\TateM_{\Qlb}(A_\alpha)\longto 0,
\end{align}
which can be seen as a sequence of sheaves on \(V_\alpha\) since
\(V'_\alpha\to V_\alpha\) is a universal homeomorphism, and
\(\TateM_{\Qlb}(P|_{V'_\alpha})\) is identified with
\(\TateM_{\Qlb}(P|_{V_\alpha})\). One may further shrink
\(V_\alpha\) so that \(\TateM_{\Qlb}(R_\alpha)\) is a local system of weight
\(-2\) (coming from the multiplicative part) and \(\TateM_{\Qlb}(A_\alpha)\)
is a local system of weight \(-1\). Here the weight means the weight of
\(\Frob_{\bar{k}/k'}\) 
for some finite extension \(k'/k\) over which \(Z_\alpha\) is defined.
Finally, one shrinks \(V_\alpha\) further so that \(V_\alpha\cap
Z_{\alpha'}=\emptyset\) unless \(Z_{\alpha}\subset Z_{\alpha'}\).

Now we choose an open subset \(U_\alpha\subset S_{\bar{k}}\) such that
\(i_\alpha\colon V_\alpha\hookto U_\alpha\) is a closed embedding. Over
\(U_\alpha\), \eqref{eqn:action_on_K_alpha} becomes
\begin{align}
    \TateM_{\Qlb}(P)\otimes i_{\alpha *}\cK_\alpha^n[\dim{V_\alpha}]\longto i_{\alpha
    *}\cK_\alpha^{n-1}[\dim{V_\alpha}].
\end{align}
The projection formula gives
\begin{align}
    \TateM_{\Qlb}(P)\otimes i_{\alpha *}\cK_\alpha^n\simeq i_{\alpha
    *}\left(i_\alpha^*\TateM_{\Qlb}(P)\otimes\cK_\alpha^n\right).
\end{align}
Since \(i_\alpha\) is closed embedding, \(i_\alpha^* i_{\alpha *}\simeq \Id\),
hence we obtain canonical map over \(V_\alpha\)
\begin{align}
    \TateM_{\Qlb}(P|_{V'_\alpha})\otimes \cK_\alpha^n\longto \cK_\alpha^{n-1}.
\end{align}
Because \(\cK_\alpha^n\) is of weight \(n\) and \(\cK_\alpha^{n-1}\) is of
weight \(n-1\), the action of \(\TateM_{\Qlb}(P|_{V'_\alpha})\) factors through
\(\TateM_{\Qlb}(A_\alpha)\) since the affine part has weight \(-2\).
Therefore, we have a graded module structure on 
\(\cK_\alpha=\bigoplus_{n}\cK_\alpha^n[-n]\) over graded 
algebra \(\Lambda_{A_\alpha}\)
\begin{align}
    \Lambda_{A_\alpha}\otimes \cK_\alpha\longto \cK_\alpha.
\end{align}
We are going to prove the following result:
\begin{proposition}\label[proposition]{prop:freeness}
    Suppose \((f,g)\) is a weak abelian fibration and \(\cF\in\Dbc(M,\Qlb)\) is
    a self-dual \(P\)-equivariant complex. Then for any geometric point
    \(u_\alpha\in V_\alpha\), the stalk \(\cK_{\alpha,u_\alpha}\) is a free
    \(\Lambda_{A_\alpha,u_\alpha}\)-module.
\end{proposition}
\begin{proof}
    It will be proved in
    \Crefrange{sec:freeness_over_a_point}{sec:freeness_by_induction}.
\end{proof}

\begin{remark}\label[remark]{rmk:freeness_independent_of_point}
    According to \cite{Ng10}*{Lemme~7.4.11}, the freeness statement is independent of
    the geometric point \(u_\alpha\). In other words, if the freeness holds
    as stated at just one point \(u_\alpha\), then it holds for any point, and
    in addition one can find some graded local system \(E\) (i.e., a direct sum
    of shifted local systems) on
    \(V_\alpha\) such that \(\cK_{\alpha}\cong \Lambda_{A_\alpha}\otimes E\) as
    graded \(\Lambda_{A_\alpha}\)-modules.
\end{remark}

\section{Freeness over a Point}%
\label{sec:freeness_over_a_point}

We first prove a generalized version of \cite{Ng10}*{Proposition~7.5.1}, which will serve as
the base case for the inductive argument later on.

\begin{lemma}\label[lemma]{lem:one_point_freeness}
    Let \(M\) be a projective variety over algebraically closed field
    \(\bar{k}\) with an action of an abelian variety \(A\) over \(\bar{k}\).
    Suppose all stabilizers are finite. Then
    \begin{align}
        \bigoplus_n \RH_{\mathup{c}}^n(M,\cF)[-n]
    \end{align}
    is a free graded \(\Lambda_A\)-module for any \(A\)-equivariant complex
    \(\cF\in\Dbc(M,\Qlb)\).
\end{lemma}
\begin{proof}
    Denote by \(f\) (resp.~\(\bar{f}\)) the map \(M\to \Spec{\bar{k}}\)
    (resp.~\(\Stack{M/A}\to\Spec{\bar{k}}\)).
    Consider quotient map \(q\colon M\to \Stack{M/A}\). By definition,
    \(\cF\simeq q^*\bar{\cF}\) for some \(\bar{\cF}\in\Dbc(\Stack{M/A},\Qlb)\).
    Here  

    The map \(q\) is smooth and projective, and following
    \cite{Ng10}*{Proposition~7.5.1}, we
    have a non-canonical isomorphism
    \begin{align}\label{eqn:one_point_noncan_iso}
        q_*\Qlb\cong \bigoplus_n \RH_{\mathup{c}}^n(A,\Qlb)[-n],
    \end{align}
    where the right-hand side is viewed as constant sheaves on \(\Stack{M/A}\).
    Since \(\Stack{M/A}\) is a Deligne--Mumford stack, the cap product action
    constructed in \ref{sec:action_by_cap_product}, although not stated explicitly,
    can also be applied to \(S=\Stack{M/A}\) with \(P=A\x S\). Therefore, both
    sides of \eqref{eqn:one_point_noncan_iso} carry actions of \(\Lambda_A\)
    over \(\Stack{M/A}\).
    We then claim that we can choose the isomorphism to be compatible with the
    respective \(\Lambda_A\)-actions. Indeed, any choice of the isomorphism
    \eqref{eqn:one_point_noncan_iso} allows us to define a morphism in
    \(\Dbc(\Stack{M/A},\Qlb)\)
    \begin{align}
        \Qlb\longto q_*\Qlb[2d](d).
    \end{align}
    Applying tensor product, we have morphism
    \begin{align}
        \Lambda_A\longto \Lambda_A\otimes q_*\Qlb[2d](d).
    \end{align}
    The composition map \(\Lambda_A\longto q_*\Qlb[2d](d)\) with cap product induces
    isomorphism on cohomology groups because isomorphism can be checked
    stalkwise, and fibers of \(q\) are just trivial \(A\)-torsors (since
    \(\bar{k}\) is algebraically closed). This proves the claim by shifting
    and twisting back by \([-2d](-d)\).

    Now by projection formula, 
    \begin{align}
        q_*\cF\simeq q_*q^*\bar{\cF}\simeq
        q_*\Qlb\otimes\bar{\cF}\cong\left(\bigoplus_n
        \RH_{\mathup{c}}^n(A,\Qlb)[-n] \right)\otimes \bar{\cF},
    \end{align}
    with the last isomorphism compatible with \(\Lambda_A\)-action. Using
    projection formula again, we see that
    \begin{align}
        f_*\cF\cong \left(\bigoplus_n
        \RH_{\mathup{c}}^n(A,\Qlb)[-n] \right)\otimes \bar{f}_*\bar{\cF}
    \end{align}
    as complexes with \(\Lambda_A\)-action. One can then compute the
    cohomologies of \(f_*\cF\) 
    using the total complex of the tensor product, and because
    \(\bigoplus_n\RH_{\mathup{c}}^n(A,\Qlb)[-n]\)  is a direct sum of complexes with only
    one non-trivial term placed at different degrees, we have that
    \begin{align}
        \bigoplus_n\RH_{\mathup{c}}^n(M,\cF)[-n]\cong \Lambda_A[-2d](-d)
    \otimes\left(\bigoplus_n\RH_{\mathup{c}}^n(\Stack{M/A},\bar{\cF})[-n]\right).
    \end{align}
    This finishes the proof.
\end{proof}

Let \(P\) be a smooth connected commutative group scheme over a finite field
\(k\). Since \(k\) is perfect, Chevalley's exact sequence is defined over \(k\),
and let \(A\) be the abelian quotient of \(P\) in that sequence.
By \cite{Ng10}*{Proposition~7.5.3}, there is a
\emph{quasi-lifting} homomorphism \(a\colon A\to P\) such that the composition
with \(P\to A\) is the endomorphism of multiplication by some positive integer
\(N\) on \(A\). This quasi-lifting induces a canonical section of Tate modules
\begin{align}
    N^{-1}\TateM_{\Qlb}(a)\colon \TateM_{\Qlb}(A)\longto \TateM_{\Qlb}(P)
\end{align}
that is compatible with Galois action.
\begin{corollary}\label[corollary]{cor:one_point_freeness_P}
    Let \(M\) be a projective variety over algebraically closed field
    \(\bar{k}\) with an action of smooth connected commutative group scheme \(P\) 
    over \(\bar{k}\).
    Suppose all stabilizers are affine and \(P\) is defined over a finite field,
    so there is quasi-lifting \(a\colon A\to P\). Then 
    \begin{align}
        \bigoplus_n \RH_{\mathup{c}}^n(M,\cF)[-n]
    \end{align}
    is a free graded \(\Lambda_A\)-module for any \(P\)-equivariant complex
    \(\cF\in\Dbc(M,\Qlb)\), where the \(\Lambda_A\)-action is defined using the
    canonical section \(N^{-1}\TateM_{\Qlb}(a)\).
\end{corollary}
\begin{proof}
    Using the quasi-lifting \(a\), we have an action of \(A\) on \(M\) with
    finite stabilizers. If \(\cF\) is \(P\)-equivariant, then it is also
    \(A\)-equivariant since \(M\to\Stack{M/P}\) factors through \(\Stack{M/A}\).
    Then the claim follows from \Cref{lem:one_point_freeness} and scaling by
    \(N^{-1}\).
\end{proof}

\section{A Degenerate Spectral Sequence}%
\label{sec:a_degenerate_spectral_sequence}

Now we move from one-point base to strict Henselian bases. The setup and
argument here are completely parallel to those in \cite{Ng10}*{\S~7.6}, only with
constant sheaf \(\Qlb\) replaced by \(\cF\). Let \(S\) be a strict
Henselian scheme over \(\bar{k}\) with closed \(\bar{k}\)-point \(s\). 
Let \(\epsilon\colon S\to\Spec{\bar{k}}\) be the structure morphism. Retain
scheme \(f\colon M\to S\) with an action of group scheme \(g\colon P\to S\) as
before. Then \(\epsilon_*\Lambda_P\) is identified with stalk \(\Lambda_{P,s}\),
and we have by adjunction
\begin{align}
    \epsilon^*\Lambda_{P,s}\longto \Lambda_P.
\end{align}
Let \(\cF\) be a \(P\)-equivariant complex on \(M\).
Then the restriction of cap product \eqref{eqn:cap_product} gives map
\begin{align}
    \epsilon^*\Lambda_{P,s}\boxtimes f_!\cF\longto f_!\cF,
\end{align}
which defines an action of graded algebra \(\Lambda_{P,s}\) on 
\(f_!\cF\) (to simplify notation, we now drop \(\epsilon^*\) and treat
\(\Lambda_{P,s}\) as a constant sheaf on \(S\)). In particular, we have a map
\begin{align}
    \TateM_{\Qlb}(P_s)\boxtimes f_!\cF\longto f_!\cF[-1].
\end{align}
Since external tensor product is perverse \(t\)-exact, we have induced map
\begin{align}
    \TateM_{\Qlb}(P_s)\boxtimes \Ptr^{\le n}(f_!\cF)\longto \Ptr^{\le
    n}(f_!\cF[-1])=\Ptr^{\le n-1}(f_!\cF)[-1],
\end{align}
hence also a map
\begin{align}
    \TateM_{\Qlb}(P_s)\boxtimes \PH^{n}(f_!\cF)\longto \PH^{n-1}(f_!\cF).
\end{align}
This map
is compatible with the map constructed in \ref{sec:action_by_cap_product} restricted to
\(\TateM_{\Qlb}(P_s)\) as they are both canonical. Therefore, we obtain a graded
\(\Lambda_{P,s}\)-module structure on
\begin{align}
    \PH^\bullet(f_!\cF)=\bigoplus_n \PH^n(f_!\cF)[-n].
\end{align}

The canonical decomposition by supports \eqref{eqn:support_decomp} and 
properties of external tensor product allows us to
express the \(\TateM_{\Qlb}(P_s)\) action by a matrix indexed by
\((\alpha,\alpha')\in\PSP\x\PSP\), with entries in
\begin{align}
    \TateM_{\Qlb}(P_s)^*\otimes \Hom(K_\alpha^n,K_{\alpha'}^{n-1}),
\end{align}
where \(\TateM_{\Qlb}(P_s)^*\) means the \(\Qlb\)-dual vector space.
In this case, if \(\alpha\neq \alpha'\), then
\(\Hom(K_\alpha^n,K_{\alpha'}^{n-1})=0\) (because the base is strictly
Henselian), so we know the matrix is diagonal.

Again use the fact that tensoring with \(\Lambda_{P,s}\) commutes with
\(\Ptr^{\le n}\), we have an increasing
filtration \(\Ptr^{\le n}(f_!\cF)\) of \(f_!\cF\)
compatible with \(\Lambda_{P,s}\)-actions. This gives a spectral sequence, also
compatible with \(\Lambda_{P,s}\)-actions,
\begin{align}\label{eqn:key_spectral_sequence}
    \ESP_2^{m,n}=\RH^m(\PH^n(f_!\cF)_s)\Rightarrow
    \RH_{\mathup{c}}^{m+n}(M_s,\cF|_{M_s}).
\end{align}
Since \(f_!\cF\) is bounded, \eqref{eqn:key_spectral_sequence} is convergent, so we
have a \emph{decreasing} filtration \(F^m\RH_{\mathup{c}}^\bullet(M_s,\cF|_{M_s})\) 
of direct sum 
\begin{align}
    \RH_{\mathup{c}}^\bullet(M_s,\cF|_{M_s})=\bigoplus_n \RH_{\mathup{c}}^n(M_s,\cF|_{M_s})[-n],
\end{align}
such that
\begin{align}
    F^m\RH_{\mathup{c}}^\bullet(M_s,\cF|_{M_s})/F^{m+1}\RH_{\mathup{c}}^\bullet(M_s,\cF|_{M_s})
    =\bigoplus_n \ESP_\infty^{m,n}[-m-n].
\end{align}
The action of \(\Lambda_{P,s}\) is compatible with the filtration \(F^m\), hence
induces an action on the associated graded \(\Qlb\)-vector space
\begin{align}
    \bigoplus_{m,n} \ESP_\infty^{m,n}[-m-n].
\end{align}
This action is the same as the \(\Lambda_{P,s}\)-action induced from the
\(\ESP_2\)-page, which in turn comes from the action on \(\PH^\bullet(f_!\cF)\).

Now we assume \((f,g)\) is the strict Henselization at \(s\) 
of a weak
abelian fibration and \(\cF\) is self-dual. Then the non-canonical
decomposition \eqref{eqn:BBD_decomp} exists on \(S\) (as perverse cohomologies
commute with \'etale base change). So the spectral sequence
\eqref{eqn:key_spectral_sequence} degenerates at \(\ESP_2\)-page and
\(\ESP_\infty^{m,n}=\ESP_2^{m,n}\).

\section{Freeness by Induction}%
\label{sec:freeness_by_induction}

We can now finish the proof of \Cref{prop:freeness} using the same
inductive argument in \cite{Ng10}*{\S~7.7}. 
The proof is by induction on dimension of \(Z_\alpha\). Let \(\alpha_0\in\PSP\)
be the unique maximal element so that \(Z_{\alpha_0}=S_{\bar{k}}\).
Let \(V_{\alpha_0}\) be an open subset of \(S\) as in
\ref{sec:statement_of_freeness}, then \(\PH^n(f_*\cF)\) is a local system
\(\cK_{\alpha_0}^n[\dim{S}]\) when restricted to \(V_{\alpha_0}\). As in
\ref{sec:statement_of_freeness}, we have a short exact sequence of Tate modules
on \(V_{\alpha_0}\)
\begin{align}
    0\longto \TateM_{\Qlb}(R_{\alpha_0})\longto
    \TateM_{\Qlb}(P|_{V_{\alpha_0}})\longto\TateM_{\Qlb}(A_{\alpha_0})\longto
    0.
\end{align}
The action of \(\TateM_{\Qlb}(P|_{V_{\alpha_0}})\) factors through
\(\TateM_{\Qlb}(A_{\alpha_0})\) because of weights. By \cite{Ng10}*{Lemme~7.4.11}
(cf.~\Cref{rmk:freeness_independent_of_point}), we may choose a geometric
point \(u_{\alpha_0}\in V_{\alpha_0}\), defined over a finite field, so that we
can use \Cref{cor:one_point_freeness_P} to deduce that
\(\RH^\bullet(M_{u_{\alpha_0}},\cF|_{M_{u_{\alpha_0}}})\) is a free
\(\Lambda_{A_{\alpha_0,u_{\alpha_0}}}\)-module. The spectral sequence
\eqref{eqn:key_spectral_sequence} is especially simple in this case and only
consists of terms \(\cK_{\alpha_0,u_{\alpha_0}}^n[-n+\dim{S}]\) for \(n\in\bbZ\).
In other words, we have a (non-canonical)
isomorphism of \(\Lambda_{A_{\alpha_0,u_{\alpha_0}}}\)-modules
\begin{align}
    \RH^\bullet(M_{u_{\alpha_0}},\cF|_{M_{u_{\alpha_0}}})\cong \bigoplus_n
    \cK_{\alpha_0,u_{\alpha_0}}^n[-n+\dim{S}].
\end{align}
So we proved \Cref{prop:freeness} in the base case. 

Next, let \(\alpha\in\PSP\) and suppose \Cref{prop:freeness} is
proved for all \(\alpha'\in\PSP\) such that \(Z_\alpha\subset Z_{\alpha'}\).
Let \(u_\alpha\in V_\alpha\) be a geometric point defined over a finite field,
and \(S_\alpha\) be the strict Henselization of \(S\) at \(u_\alpha\).
As in \ref{sec:a_degenerate_spectral_sequence}, we have actions
\begin{align}
    \TateM_{\Qlb}(P_{u_\alpha})\boxtimes
    \PH^n(f_*\cF)|_{S_\alpha}\longto\PH^{n-1}(f_*\cF)|_{S_\alpha}.
\end{align}
These actions can be represented by a diagonal matrix using decomposition over
supports \(\PSP\). This means that we have a canonical isomorphism of
graded \(\Lambda_{P,u_\alpha}\)-modules
\begin{align}
    \bigoplus_n \PH^n(f_*\cF)|_{S_\alpha}[-n]\simeq \bigoplus_{\alpha,n}
    K_\alpha^n[-n].
\end{align}
Using a quasi-lifting \(A_{u_\alpha}\to P_{u_\alpha}\), one has induces 
\(\Lambda_{A,u_\alpha}\)-module structure that is compatible with Galois action
and, when restricted to \(V_\alpha\),
is the same as the \(\Lambda_{A,u_\alpha}\)-module structure on \(\cK_\alpha\)
defined in \ref{sec:statement_of_freeness} using factorization.

Using inductive hypothesis and the fact that \(\TateM_{\Qlb}(P)\) is polarizable
(one of the axioms for weak abelian fibrations), we have the
following result \cite{Ng10}*{Proposition~7.7.4}:
\begin{proposition}
    With notation above, we have that for any \(m\in\bbZ\),
    \begin{align}
        \bigoplus_n \RH^m(K_{\alpha',u_\alpha}^n)[-n]
    \end{align}
    is a free graded \(\Lambda_{A,u_\alpha}\)-module.
\end{proposition}
\begin{proof}
    We have by inductive hypothesis and \cite{Ng10}*{Lemme~7.4.11} an isomorphism of
    \(\Lambda_{A_{\alpha'}}\)-modules on \(V_{\alpha'}\)
    \begin{align}
        \cK_{\alpha'}\cong \Lambda_{A_{\alpha'}}\otimes E_{\alpha'}
    \end{align}
    for some graded local system \(E_{\alpha'}\) on \(V_{\alpha'}\). Restricting
    to \(V_{\alpha'}\cap S_\alpha\) and pick a geometric point
    \(\bar{y}_{\alpha'}\) over the generic point \(y_{\alpha'}\) of
    \(V_{\alpha'}\cap S_\alpha\). 
    Using polarizability assumption, we can show that the specialization map
    \(\TateM_{\Qlb}(A_{u_\alpha})\to \TateM_{\Qlb}(A_{\bar{y}_{\alpha'}})\) is
    injective, hence a direct summand as
    \(\Gal(\bar{y}_{\alpha'}/y_{\alpha'})\)-representations. So the stalk 
    \(\Lambda_{A,\bar{y}_{\alpha'}}\)  is a free
    \(\Lambda_{A,u_\alpha}\)-module, hence so is
    \(\cK_{\alpha',\bar{y}_{\alpha'}}\). Using \cite{Ng10}*{Lemme~7.4.11} again, we
    have that
    \begin{align}
        \cK_{\alpha'}|_{V_{\alpha'}\cap S_\alpha}\cong
        \Lambda_{A,u_\alpha}\boxtimes E'_{\alpha'}
    \end{align}
    as \(\Lambda_{A,u_\alpha}\)-modules for some graded local system
    \(E'_{\alpha'}\). Taking intermediate extension and then stalk at
    \(u_\alpha\), and since these operations commute with external tensor
    product, we have that \(K_{\alpha',u_\alpha}\) is the external tensor
    product of \(\Lambda_{A,u_\alpha}\) with another complex, hence the same is
    true for its \(m\)-th cohomology groups.
\end{proof}

We can now finish the inductive proof using a key property of \(\Qlb\)-algebra
\(\Lambda_{A,u_\alpha}\): its projective modules are also injective.
The spectral sequence \eqref{eqn:key_spectral_sequence} degenerates at
\(\ESP_2\)-page, and we have a decreasing filtration of
\begin{align}
    H\defeq \bigoplus_n \RH^n(M_{u_\alpha},\cF|_{M_{u_\alpha}})[-n]
\end{align}
whose \(m\)-th graded part is
\begin{align}
    \bigoplus_n
    \RH^m(\PH^n(f_*\cF)_{u_\alpha})[-m-n]=\RH^m\left(\bigoplus_n\PH^n(f_*\cF)_{u_\alpha}[-n]\right)[-m],
\end{align}
all compatible with \(\Lambda_{A,u_\alpha}\)-action. It can be further
decomposed using supports into
\begin{align}
    \RH^m\left(\bigoplus_{\alpha',n}K_{\alpha',u_\alpha}^n[-n]\right)[-m]
\end{align}
as \(\Lambda_{A,u_\alpha}\)-modules. 
If \(\alpha'\neq \alpha\), then \(K_{\alpha',u_\alpha}^n\neq 0\) only if
\(Z_{\alpha'}\supset Z_\alpha\), and if \(\alpha'=\alpha\),
\(\RH^m(K_{\alpha,u_\alpha}^n)\neq 0\) only if \(m=-\dim{Z_\alpha}\).
Let \(H''=F^{-\dim{Z_\alpha}}H\) and \(H'\) be the preimage of
\begin{align}
    \RH^{-\dim{Z_\alpha}}\left(\bigoplus_{\alpha'\neq\alpha,n}K_{\alpha',u_\alpha}^n[-n]\right)[\dim{Z_\alpha}]
\end{align}
in \(H''\).
Then we have a filtration
\begin{align}
    H'\subset H''\subset H,
\end{align}
in which \(H\) is a free \(\Lambda_{A,u_\alpha}\)-module by
\Cref{cor:one_point_freeness_P}, and so is \(H'\) and \(H/H''\)
by inductive hypothesis. Then since \(H\) and \(H/H''\) are free, so is \(H''\)
because \(\Lambda_{A,u_\alpha}\) is local (by Kaplansky's theorem, which applies
to non-commutative rings). 
Consider exact sequence
\begin{align}
    0\longto H'\longto H''\longto H''/H'\longto 0.
\end{align}
As \(H'\) is free hence also injective (using the special property of
exterior algebra \(\Lambda_{A,u_\alpha}\)), the sequence splits.
So \(H''/H'=\cK_{\alpha,u_\alpha}\) is projective hence free. This finishes
the proof of \Cref{prop:freeness}.

\section{Abstract Support Theorem}
\label{sec:statement_of_support_theorem}

Now we are ready to state and prove the main theorem of this chapter. The
common assumptions in this section are as follows:
\begin{enumerate}
    \item Let \(S\) be a Deligne--Mumford stack of finite
        type over \(k\).
    \item Let \(f\colon M\to S\) be a proper morphism of Deligne--Mumford stacks
        of relative dimension \(d\) together with an action of a smooth
        commutative Deligne--Mumford group stack \(g\colon P\to S\), such that
        \((f,g)\) is a weak abelian fibration
        (see \Cref{rmk:abelian_fibrations_for_DM_stacks}).
    \item The geometric fibers of \(f\) are homeomorphic to
        projective varieties.
    \item Let \(\cF\in\Dbc(M,\Qlb)\) be a self-dual
        \(P\)-equivariant complex, such that \(f_*\cF\) has cohomological
        degrees bounded above by \(-\dim{M}+2d\).
\end{enumerate}

\begin{theorem}
    \label[theorem]{thm:support_main}
    With the assumptions above, if \(K\) is a geometrically simple summand in
    the decomposition of some \(\PH^n(f_*\cF)\) with support \(Z\). Let
    \(\delta_Z\) be the minimal value of \(\delta\) on \(Z\), then we have
    \begin{align}
        \codim_S(Z)\le\delta_Z.
    \end{align}
    Moreover, if the equality holds, then there exists an open subset \(U\subset
    S_{\bar{k}}\) such that \(U\cap Z\neq\emptyset\), and a local system \(L\)
    on \(U\cap Z\) such that \(i_*L\) (where \(i\) is the closed embedding
    \(U\cap Z\to U\)) is a direct summand of the top cohomology
    \(\RDF^{-\dim{M}+2d}f_*\cF\) restricted to \(U\).
\end{theorem}
\begin{proof}
    Although \(f\) and \(g\) are not morphism of schemes, the argument for
    \Cref{prop:freeness} still applies. More precisely, everything except
    \Cref{lem:one_point_freeness} works verbatim for general
    Deligne--Mumford stacks,
    and \Cref{lem:one_point_freeness} only potentially breaks because \textit{a
    priori} \(\Stack{M_s/P_s}\) for a geometric point \(s\in S\) can be a
    \(2\)-stack instead of an algebraic stack. However, since we also assume
    that \(M_s\) is homeomorphic to a projective scheme (and \(P_s\) is always
    homeomorphic to a smooth group scheme), the argument in
    \Cref{lem:one_point_freeness} still works.

    By \Cref{prop:freeness}, \(K|_{U\cap Z}\) is a direct summand of a
    free \((\Lambda_P)|_{U\cap Z}\)-module, whose top cohomology is denoted by
    \(L\). Moreover, the intermediate extension of \(L\) to \(Z\) is also a
    simple perverse summand of \(f_*\cF\). Then the argument in
    \Cref{sec:goresky_macpherson_inequality} proves both claims of the theorem.
\end{proof}

\begin{corollary}
    \label[corollary]{cor:support_main_alt}
    If in \Cref{thm:support_main} \((f,g)\) is also \(\delta\)-regular, then the
    equality \(\codim_S(Z)=\delta_Z\) holds.
\end{corollary}
\begin{proof}
    Clear since the definition of \(\delta\)-regularity gives the inequality in
    the opposite direction.
\end{proof}

\subsection{}
\label{sub:kappa_decomposition_at_chain_level}
Let \(\pi_0(P)\) be the \'etale sheaf of finite abelian groups over \(S\) such that for
any geometric point \(s\in S\) its stalk is the group of connected components
\(\pi_0(P_s)\) (see~\cite{Ng06}*{\S~6.2}).
Suppose \(\pi_0(P)\) is a
quotient of some constant sheaf
\(\bX\) where \(\bX\) is a finite abelian group. Since \(\cF\) is
\(P\)-equivariant, \(P\) canonically acts on \(\cF\) by
\eqref{eqn:canonical_identification_for_equivariant_complex}. We have the
following homotopy lemma:
\begin{lemma}
    [Homotopy Lemma]
    \label[lemma]{lem:homotopy_lemma}
    Let \(S\) be a \(k\)-scheme and \(p\colon P\to S\) be a smooth group scheme
    with geometrically connected fibers. Let \(f\colon M\to S\) be an
    algebraic stack over \(S\) with a \(P\) action, and \(\cF\in\Dbc(M,\Qlb)\)
    be a \(P\)-equivariant sheaf. Then the set of global sections \(P(S)\)
    acts trivially on both \(\PH^n(f_*\cF)\) and \(\RH^n(f_*\cF)\).
\end{lemma}
\begin{proof}
    This is a simple generalization of \cite{LN08}*{Lemme~3.2.3} and the proof
    is the same. We record the slightly modified proof for completeness because
    unlike \textit{loc. cit.}, \(M\) is no longer a scheme and \(\cF\) is not
    the constant sheaf.

    Consider the diagram:
    \begin{equation}
        \begin{tikzcd}
            P\x_S M \ar[r, "\alpha"]\ar[rd, "\pr_P", swap] & P\x_S M \ar[r, "\pr_M"]\ar[d, "\pr_P"] & M \ar[d, "f"]\\
             & P \ar[r, "p"] & S
        \end{tikzcd}
    \end{equation}
    where \(\alpha\) is the map \((p,m)\mapsto (p,pm)\). By smooth base
    change (\cite{LaOl08II}*{Lemma~9.1.2, \S~12.1}), we have
    \begin{align}
        p^*f_*\cF\simeq \pr_{P*}\pr_M^*\cF,
    \end{align}
    and so \(\alpha\) together with the equivariant structure on \(\cF\) induces
    an endomorphism \([\alpha]\) of \(p^*f_*\cF\). By \(t\)-exactness of \(p^*\)
    (with respect to both the standard and perverse \(t\)-structures), we have
    \begin{align}
        p^*\bigl(\PH^n(f_*\cF)\bigr)=\PH^n(\pr_{P*}\pr_M^*\cF),\quad
        p^*\bigl(\RH^n(f_*\cF)\bigr)=\RH^n(\pr_{P*}\pr_M^*\cF),
    \end{align}
    and \([\alpha]\) induces an endomorphism of either cohomology sheaf.

    Let \(\cH\) be either \(\PH^n(f_*\cF)\) or \(\RH^n(f_*\cF)\).
    By \Cref{lem:fully_faithfulness_of_smooth_pullback} below, since \(p\) is
    smooth with geometrically connected fibers, the functor \(p^*\) is fully
    faithful when restricted to either the subcategory of perverse sheaves or
    the ordinary sheaves.
    This implies that
    \([\alpha]\) corresponds to a unique endomorphism \(\beta\) of
    \(\cH\) such that \(p^*\beta=[\alpha]\). Let \(g\colon S\to P\) be
    any section, then \(g^*[\alpha]\) is none other
    than the action of \(g\in P(S)\) on \(\cH\). Since \(p\circ g=\id_S\),
    \(g^*[\alpha]=g^*p^*\beta=\beta\). Letting \(g=1\), we see that \(\beta\) is
    the identity, and so is \([\alpha]\). The lemma follows.
\end{proof}

\begin{lemma}
    \label[lemma]{lem:fully_faithfulness_of_smooth_pullback}
    Let \(f\colon T\to S\) be a smooth surjective map of schemes over \(k\) with
    geometrically connected fibers. Then the restriction of the functor
    \(f^*\colon\Dbc(S,\Qlb)\to\Dbc(T,\Qlb)\) to either the full subcategory of
    perverse sheaves or that of the ordinary sheaves is fully faithful.
\end{lemma}
\begin{proof}
    The perverse version is \cite{BBD82}*{Proposition~4.2.5}, and the same
    proof works for the standard \(t\)-structure which we include here for
    completeness.

    Let \(\cF,\cG\in\Dbc(S,\Qlb)\). We will denote by \(f^0_*\) the
    \emph{non-derived} pushforward by \(f\), defined for any ordinary sheaf on
    \(T\). Note that the functor
    \(f^*\) is \(t\)-exact with respect to the standard \(t\)-structures on
    \(\Dbc(S,\Qlb)\) and \(\Dbc(T,\Qlb)\). By smoothness, we have a
    canonical isomorphism
    \begin{align}
        f^*\RIHom(\cF,\cG)\stackrel{\sim}{\longto} \RIHom(f^*\cF,f^*\cG).
    \end{align}
    Apply the functor \(f^0_*\RH^0\) and using
    \(t\)-exactness, we have
    \begin{align}
        f^0_*f^*\IHom(\cF,\cG)\stackrel{\sim}{\longto}f^0_*\IHom(f^*\cF,f^*\cG).
    \end{align}

    For any ordinary \(\Qlb\)-sheaf \(\cH\) on \(S\), we have by
    adjunction a map \(\cH\to f^0_*f^*\cH\). Since \(f\) is smooth, it admits
    a section \(e\) over an \'etale cover \(S'\to S\), and so the composition
    \begin{align}
        \cH(S')\longto f^0_*f^*\cH(S')\longto f^0_*e^0_*e^*f^*\cH(S')\simeq \cH(S')
    \end{align}
    is the identity. For any two sections \(u,v\in f^0_*f^*\cH(S')\), suppose they
    have the same image in \(\cH(S')\), then they must coincide over the
    Henselization \(T'_{e}\) of \(e(S')\) in \(T'=T\x_S S'\).
    For any two \'etale maps \(U\to T'\) and \(V\to T'\) with the same (necessarily
    open) image  in \(S'\), we have \(U\x_{T'} V\neq\emptyset\) by
    connectivity and smoothness of \(f\). As a result, the presheaf
    \begin{align}
        U\longmapsto \cH(f(U))
    \end{align}
    is already an \'etale sheaf, and is precisely the sheaf \(f^*\cH\) by
    definition. Therefore, we have
    \begin{align}
        u=v\in f^*\cH(T')=\cH(S')=\cH(f(T'_e))=f^*\cH(T'_e).
    \end{align}
    This shows that \(\cH\to
    f^0_*f^*\cH\) is an isomorphism. Consequently, we have
    \begin{align}
        \IHom(\cF,\cG)\simeq
        f^0_*f^*\IHom(\cF,\cG)\stackrel{\sim}{\longto}
        f^0_*\IHom(f^*\cF,f^*\cG).
    \end{align}

    Suppose \(\cF\) and \(\cG\) are both ordinary
    sheaves, then both \(\RIHom(\cF,\cG)\) and \(f_*\RIHom(f^*\cF,f^*\cG)\) are
    supported on degrees \(\ge 0\).
    Taking the hypercohomology functor \(\RH^0(S,-)\) and using a standard
    spectral sequence, we obtain
    \begin{multline}
        \Hom(\cF,\cG)=\RH^0\bigl(S,\RIHom(\cF,\cG)\bigr)\simeq\RH^0\bigl(S,\IHom(\cF,\cG)\bigr)
        \stackrel{\sim}{\longto}\RH^0\bigl(S,f^0_*\IHom(f^*\cF,f^*\cG)\bigr)\\
        \simeq\RH^0\bigl(S,f_*\RIHom(f^*\cF,f^*\cG)\bigr)=
        \RH^0\bigl(T,\IHom(f^*\cF,f^*\cG)\bigr)=\Hom(f^*\cF,f^*\cG).
    \end{multline}
    This finishes the proof.
\end{proof}

\begin{remark}
    \label[remark]{rmk:apply_homotopy_to_mH}
    In our applications, \(S\) will be a suitable locally closed substack of the mH-base
    \(\cA_X^\heartsuit\), and \(P\) will be the restriction of \(\cP_X\) to
    \(S\). Although now \(S\) is no longer a scheme, it is by construction
    still a finite \'etale gerbe over a scheme. Moreover, \(\cP_X\) is a group
    stack instead of a scheme, but we can still use the regidification process
    described in \Cref{sub:rigidification_of_global_Picard} to pullback the
    \(\cP_X|_S\)-action to an action of a smooth group scheme
    \(\cP_X^\infty|_S\) such that
    \(\pi_0(\cP_X^\infty|_S)\simeq\pi_0(\cP_X|_S)\). This way
    \Cref{lem:homotopy_lemma} still applies.
\end{remark}

By \Cref{lem:homotopy_lemma} and proper base change, the induced action of \(P\) on
\(\PH^\bullet(f_*\cF)\) factors through \(\pi_0(P)\),
hence \(\bX\) acts on \(\PH^n(f_*\cF)\)
for each \(n\). Let \(\PH^n(f_*\cF)_\kappa\) be the \(\kappa\)-isotypic direct
summand for any character \(\kappa\colon \bX\to \Qlb^\x\). Similarly, we have
\(\kappa\)-isotypic summand of ordinary cohomology \(\RH^n(f_*\cF)_\kappa\). By
\cite{LN08}*{Lemme~3.2.5}, there exists an integer \(N>0\) and a decomposition
\begin{align}
    f_*\cF=\bigoplus_{\kappa\in\bX^*}(f_*\cF)_\kappa,
\end{align}
such that for any \(\alpha\in\bX\), the restriction of
\(\bigl(\alpha-\kappa(\alpha)\id\bigr)^N\) on \((f_*\cF)_\kappa\) is zero, and we have
\begin{align}
    \PH^n(f_*\cF)_\kappa&=\PH^n\bigl((f_*\cF)_\kappa\bigr),\\
    \RH^n(f_*\cF)_\kappa&=\RH^n\bigl((f_*\cF)_\kappa\bigr).
\end{align}

\begin{proposition}
    \Cref{thm:support_main,cor:support_main_alt} hold if we replace
    \(\PH^n(f_*\cF)\) by \(\PH^n(f_*\cF)_\kappa\) and
    \(\RH^{-\dim{M}+2d}(f_*\cF)\) by \(\RH^{-\dim{M}+2d}(f_*\cF)_\kappa\).
\end{proposition}
\begin{proof}
    Apply the same proof to \((f_*\cF)_\kappa\).
\end{proof}

\section{Application to mH-fibrations} 
\label{sec:application_to_mH_fibrations}

In this section, we apply \Cref{thm:support_main} to mH-fibration
\(h_X\colon\cM_X\to\cA_X\), or more precisely a subset of the map
\(\tilde{h}_X^\ANI\colon\tilde{\cM}_X^\ANI\to\tilde{\cA}_X^\ANI\).
There is an action of Picard stack
\(\tilde{p}_X^\ANI\colon\tilde{\cP}_X^\ANI\to\tilde{\cA}_X^\ANI\) on
\(\tilde{h}_X^\ANI\), and both are relative Deligne--Mumford stacks of finite
type by \Cref{prop:Mnatural_DM_finite_type_over_anisotropic_locus}. By
\Cref{prop:properness_over_anisotropic_locus}, \(\tilde{h}_X^\ANI\) is also
proper. If we restrict to very \(G\)-ample locus \(\tilde{\cA}_\gg^\ANI\),
by \Cref{prop:tate_module_is_polarizable}, the Tate modules of
\(\tilde{\cP}_\gg^\ANI\) are polarizable. Therefore,
\((\tilde{h}_\gg^\ANI,\tilde{p}_\gg^\ANI)\) is a weak abelian fibration.

In general, \((\tilde{h}_\gg^\ANI,\tilde{p}_\gg^\ANI)\) may not be
\(\delta\)-regular, but we can still find an open subset over which
the fibration is \(\delta\)-regular.
This subset is dense in \(\tilde{\cA}_\gg^\ANI\)
because the \(\delta=0\)
stratum is by \Cref{prop:non_emptyness_of_A_sharp_diamond_heart}.
According to \Cref{prop:delta_regularity_fix_divisor}, for any fixed
\(\delta\in\bbN\), all but finitely many irreducible components of
\(\tilde{\cA}_{\delta}^\ANI\) is contained in this \(\delta\)-regular locus.

Let \(\tilde{\cA}^\dagger\subset\tilde{\cA}_\gg^\ANI\)
\nomenclature[\(.dagger \)]{\((\cdot)^\dagger\)}{the technical
    locus within the anisotropic locus where the local model of singularity holds,
    and the degree of \(G\)-ampleness is much larger than the
    \(\delta\)-invariant, so that the support theorem applies to the stable
    constituent of the cohomology}
be the largest open
subset satisfying the following conditions:
\begin{enumerate}
    \item For any \(\delta\in\bbN\), \(\tilde{\cA}^\dagger_\delta\) is either
        empty or very \((G,N(\delta))\)-ample (cf.,
        \Cref{prop:delta_regularity_fix_divisor}). In particular, the
        restriction of \((h_X,p_X)\) to \(\tilde{\cA}^\dagger\) is
        \(\delta\)-regular.
    \item The local model of singularity as in either
        \Cref{thm:local_singularity_model_weak} or
        \Cref{thm:local_singularity_model_main} holds.
\end{enumerate}
Let \(\bcQ\in\Sat_X\) be a Satake sheaf on \(\tilde{\cM}^\dagger\).
Then by \Cref{prop:cohomological_amplitude_of_mH_fibrations}, we have that over
any connected component \(\cU\) of \(\tilde{\cA}^\dagger\), the cohomological
degrees of \(\tilde{h}^\dagger_*\bcQ\) is bounded above by
\(-\dim{\cM_\cU}+2d\), where \(\cM_\cU\) is the preimage of \(\cU\) in
\(\tilde{\cM}^\dagger\), and \(d\) is the relative dimension of
\(\tilde{h}^\dagger\) over \(\cU\).
Finally, the fibers of \(\tilde{h}^\dagger\) are homeomorphic to projective
varieties by \Cref{cor:mH_fibers_homeomorphic_to_proj_varieties}.

\subsection{}
\label{sub:inductive_closed_subset}
Let \(Z\subset \tilde{\cA}^\dagger\) be a locally closed inductive
subset (cf.~\Cref{def:inductive_subset}).
Since we assume that
\(\tilde{\cA}_\delta^\dagger\) is very \((G,N(\delta))\)-ample, it is
by \Cref{prop:delta_critical_in_inductive} equivalent to saying that
there exists an open subset \(\cU\subset\tilde{\cA}^\dagger\) such that
\(\cU\cap Z\neq\emptyset\) has multiplicity-free boundary divisors
(\Cref{def:mult_free_boundary_divisor}), and for any \(a\in
(\cU\cap Z)(\bar{k})\) and any \(\bar{v}\in X(\bar{k})\), one of the following
cases holds:
\begin{enumerate}
    \item The boundary divisor is \(0\) at \(\bar{v}\), and
        the discriminant valuation \(d_{\bar{v}+}(a)\le 1\), or
    \item \(a\) is both unramified and \(\nu\)-regular semisimple at \(\bar{v}\).
\end{enumerate}

\subsection{}
We are now able to prove one of the main results of this book:
\begin{theorem}[Support Theorem for the Stable Constituent]
    \label[theorem]{thm:stable_support_of_mH_main}
    Let \(K\) be a simple perverse summand of
    \(\PH^\bullet(\tilde{h}^\dagger_*\bcQ)_{\hST}\), and \(Z\) is the support of
    \(K\), then \(Z\) is inductive. In particular, \(Z\) must be
    \(\delta\)-critical.
\end{theorem}
\begin{proof}
    Applying \Cref{cor:support_main_alt} to
    \((\tilde{h}^\dagger,\tilde{p}^\dagger)\), we see that \(Z\) must be
    \(\delta\)-critical. So 
    we only need to show that if \(Z'\subset\tilde{\cA}^\dagger\) is
    \(\delta\)-critical but not inductive, then
    \(Z'\) cannot be the support of any \(K\).

    Indeed, using \Cref{prop:delta_critical_in_inductive}, we may find an
    inductive subset \(Z\subset
    \tilde{\cA}^\dagger\) containing \(Z'\) such that there exists open subset
    \(\cU\subset\tilde{\cA}^\dagger\) with \(\cU\cap Z'\neq\emptyset\) and both
    the boundary divisor and all the local Newton points stay locally constant
    throughout \(\cU\cap Z\).
    By \Cref{prop:stable_top_cohomology_local_system_main}, the top cohomology
    of \((\tilde{h}^\dagger_*\bcQ)_{\hST}\) is a local system \(L\) on \(\cU\cap
    Z\). However, if \(Z'\) supports a perverse summand, then \(L\) contains a
    direct summand supported on proper subset \(\cU\cap Z'\), which is
    impossible. This finishes the proof.
\end{proof}

\subsection{}
An interesting corollary of \Cref{thm:stable_support_of_mH_main} is the
following local result. Since it will not be used in this book, we only give a
sketch and leave the details to the reader:

\begin{corollary}
    \label[corollary]{cor:delta_equals_0_implies_torsor}
    If \(\delta_{\bar{v}}(a)=0\), then \(\cM_{\bar{v}}(a)\) is a
    \(\cP_{\bar{v}}(a)\)-torsor.
\end{corollary}
\begin{proof}
    If \(a\) is unramified at \(\bar{v}\), then the conclusion is clear, so we
    may assume that it is ramified at \(\bar{v}\).
    We may find a point \(a\in\tilde{\cA}^\dagger\) such that for any
    \(\bar{v}'\neq\bar{v}\), we have either \(a\) is regular
    semisimple or \(\bar{v}'\) does not support the boundary divisor and
    \(d_{\bar{v}'+}(a)=c_{\bar{v}'}=1\). We may also assume that \(a\) has
    multiplicity-free boundary divisor. The subset of such \(a\) is a
    subset of the \(\delta=0\) stratum in \(\tilde{\cA}^\dagger\), and
    has codimension at least \(2\).

    Since in a neighborhood of such subset,
    the top ordinary cohomology of
    \((\tilde{h}^\dagger_*\IC_{\tilde{\cM}^\dagger})_\hST\) is a local
    system, which extends over codimension-\(2\) subsets, we have that the stalk
    of \((\tilde{h}^\dagger_*\IC_{\tilde{\cM}^\dagger})_\hST\) at
    \(a\) is \(1\)-dimensional as well. This means that
    \(\cP_{\bar{v}}(a)\) acts on the irreducible components of
    \(\cM_{\bar{v}}(a)\) simply transitively. Since \(\cM_{\bar{v}}(a)\) is
    \(0\)-dimensional and contains a \(\cP_{\bar{v}}(a)\)-torsor as an open
    subset, it cannot have other irreducible component and thus must be a
    \(\cP_{\bar{v}}(a)\)-torsor.
\end{proof}

\subsection{}
The case of \(\kappa\)-constituent is more complicated. The most notable issue
is that an endoscopic stratum may not be contained in
\(\tilde{\cA}^\dagger\), hence we cannot use
\Cref{prop:delta_regularity_fix_divisor} to deduce \(\delta\)-regularity.
Nevertheless, for \(\kappa\)-constituents, \(\delta\)-regularity can be replaced
by a weaker condition, namely that the intersection of any \(\delta\)-stratum
with the endoscopic strata has codimension at least \(\delta\).
This can be proved with the
help of \eqref{eqn:difference_in_dimension_of_mH_bases} and
\eqref{eqn:difference_in_delta_invariants}.
Some other ingredients such as \Cref{prop:delta_critical_in_inductive} also need
modification due to insufficient \(G\)-ampleness, and we will address them as
they come up.

Since we are not confined to \(\tilde{\cA}^\dagger\), we consider a
larger open subset \(\tilde{\cA}^\ddagger\)
\nomenclature[\(.ddagger \)]{\((\cdot)^\ddagger\)}{the technical
    locus within the anisotropic locus larger than \((\cdot)^\dagger\)
    where the local model of singularity holds but no requirement on
    \(G\)-ampleness versus \(\delta\)-invariant}
defined similar to
\(\tilde{\cA}^\dagger\) but with \((G,N(\delta))\)-ampleness condition removed
(in other words, we only require the local model of singularity to
hold). Let \((\kappa,\OGT_{\kappa}^\bullet)\) be an endoscopic datum with
endoscopic group \(H\). There is a canonical finite unramified map
\begin{align}
    \tilde{\nu}_\cA\colon \tilde{\cA}_{H,X}^\kappa\longto\tilde{\cA}_X.
\end{align}
By \eqref{eqn:difference_in_delta_invariants}, the difference
\(r_H^G=\delta-\delta_H\) is locally constant on \(\tilde{\cA}_{H,X}^\kappa\).

We would like to restrict to the locus \(\tilde{\cA}_H^{\kappa,\dagger}\),
\nomenclature[\(.dagger_H^kappa \)]{\((\cdot)_H^{\dagger,\kappa}\)}{the technical
    locus in \((\cdot)_{H,X}^\kappa\) analogous to \((\cdot)^\dagger\), with
    some extra conditions related to \((\cdot)^\ddagger\), in order to applying the
    support theorem on the \(\kappa\)-constituent of the cohomology}
which
is defined to be the open subset of
\(\tilde{\cA}_{H,X}^{\kappa,\ANI}\) satisfying the following conditions
similar to \(\tilde{\cA}^\dagger\) (with some additions):
\begin{enumerate}
    \item Its image under \(\tilde{\nu}_\cA\) is contained in
        \(\tilde{\cA}^\ddagger\) and
        \(\tilde{\cA}_H^{\kappa,\dagger}=\tilde{\nu}_\cA^{-1}\tilde{\nu}_\cA(\tilde{\cA}_H^{\kappa,\dagger})\).
    \item For any \(\delta_H\in\bbN\),
        \(\tilde{\cA}_{H,\delta_H}^{\kappa,\dagger}\) is either empty or very
        \((H,N(\delta_H))\)-ample
        (cf.,~\Cref{prop:delta_regularity_fix_divisor}), so that
        \((\tilde{h}_H^{\kappa,\dagger},\tilde{p}_H^{\kappa,\dagger})\) is
        \(\delta_H\)-regular.
    \item The local model of
        singularity as in \Cref{thm:local_singularity_model_weak} holds
        (but not \Cref{thm:local_singularity_model_main} because it is not
        formulated for \(\cM_{H,X}^\kappa\)).
    \item The codimension formula
        \eqref{eqn:difference_in_dimension_of_mH_bases} holds over the preimage
        of \(\tilde{\cA}^\ddagger\) in \(\tilde{\cA}_H^{\kappa,\dagger}\). In
        particular, it holds if cohomological condition
        \eqref{eqn:dim_of_endoscopic_strata_cohom_assumption} holds over
        \(\tilde{\cA}_H^{\kappa,\dagger}\).
\end{enumerate}

For any irreducible component \(\cU_H\) of \(\tilde{\cA}_H^{\kappa,\dagger}\),
let \(\cU=\tilde{\nu}_{\cA}(\cU_H)\). By
\eqref{eqn:difference_in_dimension_of_mH_bases}, the codimension of \(\cU\) in
\(\tilde{\cA}_X\) is \(r_H^G\) (necessarily a constant on \(\cU\)).
As a result, for any \(\delta\ge r_H^G\), we have
\begin{align}
    \codim_{\tilde{\cA}_X}(\cU_\delta)=r_H^G+\codim_{\tilde{\cA}_{H,X}}(\cU_{H,\delta-r_H^G})
    \ge r_H^G+(\delta-r_H^G)=\delta.
\end{align}
In other words, \(\delta\)-regularity holds for the intersection of any
\(\delta\)-strata with \(\cU\).
We will eventually show that \(\tilde{\cA}_H^{\kappa,\dagger}\) is non-empty in
\Cref{sec:construction_of_a_good_mH_base}.

\begin{theorem}
    [Support Theorem for the \(\kappa\)-Constituent]
    \label[theorem]{thm:support_endoscopic_main}
    Let \(\PH^\bullet(\tilde{h}^\ddagger_*\bcQ)_\kappa\) be the
    \(\kappa\)-isotypic summand in \(\PH^\bullet(\tilde{h}^\ddagger_*\bcQ)\) and
    \(K\) is a simple perverse summand in
    \(\PH^\bullet(\tilde{h}^\ddagger_*\bcQ)_\kappa\). Let \(Z\) be the support
    of \(K\). Then there exists a pointed endoscopic datum
    \((\kappa,\OGT_\kappa^\bullet)\) with endoscopic group \(H\) such that \(Z\)
    is contained in \(\tilde{\nu}_\cA(\tilde{\cA}_{H,X}^\kappa)\).
    Moreover, if we can find an open subset \(\cU\subset \tilde{\cA}^\ddagger\)
    such that \(\cU\cap Z\neq \emptyset\),
    then each irreducible component of
    \(\tilde{\nu}_\cA^{-1}(Z)\) that dominates \(Z\) is inductive in
    \(\tilde{\cA}_H^{\kappa,\dagger}\) and is a finite \'etale cover of \(Z\).
\end{theorem}
\begin{proof}
    The first claim is \Cref{prop:kappa_locus_is_from_endoscopy}, so we only
    need to prove the second claim. Similar to
    \Cref{thm:stable_support_of_mH_main}, \(Z\) must be \(\delta\)-critical.
    By \Cref{prop:pointed_endoscopic_transfer_closed_embedding} and shrinking
    \(\cU\) if necessary, we may assume that every irreducible
    component of \(\tilde{\nu}_\cA^{-1}(Z\cap\cU)\) is \(\delta_H\)-critical and
    is a finite \'etale cover of \(Z\cap \cU\).

    Let \(Z_H\) be one of the components in \(\tilde{\nu}_\cA^{-1}(Z\cap\cU)\),
    then by \Cref{prop:delta_critical_in_inductive}, we can find some
    irreducible locally closed subset \(\tilde{Z}_H\subset \tilde{\cA}_H^{\kappa,\dagger}\)
    containing \(Z_H\) whose closure is inductive, and the local Newton points
    stay locally constant. In particular, \(\tilde{Z}_H\) is smooth. Moreover, at any point
    \(\bar{v}\) supporting the boundary divisor, any \(a_H\in \tilde{Z}_H\) must be
    unramified and \(\nu\)-regular semisimple at \(\bar{v}\), which implies the
    same for \(\tilde{\nu}_\cA(a_H)\).

    Let \(\tilde{Z}\) be the
    isomorphic image of \(\tilde{Z}_H\) under \(\tilde{\nu}_\cA\), which is then smooth,
    locally closed in \(\tilde{\cA}^\ddagger\) and contains \(Z\cap \cU\). Apply
    \Cref{prop:kappa_top_cohomology_local_system_main}, and we see that the top
    cohomology of \((\tilde{h}_*^\ddagger\bcQ)_\kappa\) is a local system on
    \(\tilde{Z}\). This means that \(Z\cap \cU\) must be open in \(\tilde{Z}\), and
    the closure of \(Z_H\) is inductive in \(\tilde{\cA}_H^{\kappa,\dagger}\).
    This finishes the proof.
\end{proof}

\subsection{}

With \Cref{thm:stable_support_of_mH_main,thm:support_endoscopic_main}
in mind, we are now able to conditionally prove geometric stabilization over
\(\tilde{\cA}_H^{\kappa,\dagger}\). More precisely, we have the following
result:
\begin{proposition}
    \label[proposition]{prop:very_weak_stabilization}
    Suppose for any inductive subset \(Z_H\) of
    \(\tilde{\cA}_H^{\kappa,\dagger}\), there exists an open dense subset of
    \(\tilde{\nu}_\cA(Z_H)\) over which
    \eqref{eqn:primal_stabilization_for_top_cohomology_preview} holds. Then
    \Cref{thm:arithmetic_geometric_stabilization} hence also
    \Cref{thm:main_geometric_stabilization} holds over
    \(\tilde{\nu}_\cA(\tilde{\cA}_H^{\kappa,\dagger})\).
\end{proposition}

\begin{proof}
    By \Cref{thm:stable_support_of_mH_main}, the set of supports of
    \(\PH^\bullet(\tilde{h}_{H,*}^{\kappa,\dagger}\bcQ_H^\kappa)_\hST\) only
    contains inductive subsets in \(\tilde{\cA}_{H}^{\kappa,\dagger}\).
    Similarly, we may apply
    \Cref{thm:support_endoscopic_main} to \(\tilde{\cA}^\ddagger\), and it implies that the set of
    supports of \(\PH^\bullet(\tilde{h}^\ddagger_*\bcQ)_\kappa\) only contains
    subsets that are either disjoint from
    \(\tilde{\nu}_\cA(\tilde{\cA}_H^{\kappa,\dagger})\) or
    \(\delta\)-critical. For the latter kind, their preimages in
    \(\tilde{\cA}_{H}^{\kappa,\dagger}\) are inductive, and
    we shall call such subset \notion{\(\kappa\)-inductive}\index{subset!\(\kappa\)-inductive}
    \index{stratum!\(\kappa\)-inductive} for convenience.

    Let \(Z\) be an irreducible \(\kappa\)-inductive subset of
    \(\tilde{\cA}^\ddagger\). Suppose \(K_Z\subset
    \PH^\bullet(\tilde{h}^\ddagger_*\bcQ)_\kappa\) is the largest direct summand
    with support \(Z\), then by the proof of \Cref{thm:support_endoscopic_main}, we
    may find an open subset \(Z'\subset Z\) such that \(K_Z|_{Z'}\) is a local
    system. Shrinking \(Z'\) if necessary, we may also assume that the abelian part
    of the Tate module \(\Lambda_{A,Z}|_{Z'}\) forms a local system and
    \(K_Z|_{Z'}\) is a free
    graded \(\Lambda_{A,Z}|_{Z'}\)-module generated by the top degree.

    Similarly, if \(K_{H,Z}\subset
    \tilde{\nu}_{\cA,*}\PH^\bullet(\tilde{h}^\dagger_{H,*}\bcQ_H^\kappa)_\hST\) is
    the largest direct summand with support \(Z\), then further shrinking \(Z'\) if
    necessary, we may also assume that \(\tilde{\nu}_{\cA}^{-1}(Z')\) is just
    multiple disjoint isomorphic copies of \(Z'\) and \(K_{H,Z}|_{Z'}\) is a free
    graded \(\tilde{\nu}_{\cA,*}\Lambda_{H,A_H}|_{Z'}\)-module generated be the top
    cohomology. Moreover, the natural maps \(\cP_a\to\cP_{H,a_H}\) induces natural
    isomorphism \(\tilde{\nu}_\cA^*\Lambda_{A,Z}\simeq
    \Lambda_{H,A_H,\tilde{\nu}_\cA^{-1}(Z)}\) over \(\tilde{\nu}_\cA^{-1}(Z')\), so
    \(K_{H,Z}|_{Z'}\) is also a free graded \(\Lambda_{A,Z}|_{Z'}\)-module generated
    by the top degree.

    Assuming \eqref{eqn:primal_stabilization_for_top_cohomology_preview} holds
    generically over \(Z'\), then the stalks of the top cohomologies of
    \(\PH^\bullet(\tilde{h}^\ddagger_*\bcQ)_\kappa\) and
    \(\tilde{\nu}_{\cA,*}\PH^\bullet(\tilde{h}^\dagger_{H,*}\bcQ_H^\kappa)_\hST\)
    are isomorphic as Frobenius modules, which
    means that they are isomorphic as local systems after taking semisimplifications
    with respect to the Frobenius. We also know that both \(K_Z\) and \(K_{H,Z}\) are
    geometrically semisimple for any \(Z\). Therefore, we can deduce
    \Cref{thm:arithmetic_geometric_stabilization} hence also
    \Cref{thm:main_geometric_stabilization} over \(\tilde{\nu}_\cA(\tilde{\cA}_H^{\kappa,\dagger})\) using induction.

    Indeed, there is a partial order on (the closures of) \(\kappa\)-inductive subsets \(Z\) by
    inclusion. The maximal elements are just irreducible components of
    \(\tilde{\nu}_\cA(\tilde{\cA}_{H}^{\kappa,\dagger})\), and we let \(Z_0\)
    to be one of them. We have already established that the
    semisimplifications \(K_{Z_0}^\ssim\) and \(K_{H,Z_0}^\ssim\) are
    isomorphic over some open subset \(Z_0'\subset Z_0\). Although the functor of
    intermediate extension is not exact in general, it is still exact on
    geometrically semisimple local systems. Therefore, \(K_{Z_0}^\ssim\) and
    \(K_{H,Z_0}^\ssim\) are isomorphic. This proves the base case of induction.

    Over a sufficiently small open subset \(Z'\)
    in an arbitrary \(Z\), \(\PH^\bullet(\tilde{h}^\ddagger_*\bcQ)_\kappa\)
    decomposes into \(K_Z|_{Z'}\) and \(K_W|_{Z'}\) for all \(\kappa\)-inductive
    \(W\) containing \(Z\), and similarly for
    \(\tilde{\nu}_{\cA,*}\PH^\bullet(\tilde{h}^\dagger_{H,*}\bcQ_H^\kappa)_\hST\).
    By inductive hypothesis, all \(K_W^\ssim\) and \(K_{H,W}^\ssim\) are isomorphic,
    and in particular it is true for their top cohomologies. This implies that the
    top cohomologies of \(K_Z|_{Z'}^\ssim\) and \(K_{H,Z}|_{Z'}^\ssim\) are also
    isomorphic. Since \(K_Z|_{Z'}\) is generated by its top cohomology over
    \(\Lambda_{A,Z}|_{Z'}\), and similarly for \(K_{H,Z}|_{Z'}\), we see that
    \(K_Z|_{Z'}^\ssim\) and \(K_{H,Z}|_{Z'}^\ssim\) are isomorphic, hence so are
    \(K_Z^\ssim\) and \(K_{H,Z}^\ssim\). This finishes the proof of
    \Cref{thm:arithmetic_geometric_stabilization} over \(\tilde{\nu}_\cA(\tilde{\cA}_{H}^{\kappa,\dagger})\).
\end{proof}

\chapter{Counting Points} 
\label{chap:counting_points}

In this chapter we review general facts about counting points on
stacks defined over a finite field \(k\). We will then connect point counting on multiplicative affine
Springer fibers with both their cohomology and the orbital
integrals. Similarly, we will also apply the general results to mH-fibrations.
Most of the results are generalization of those in \cite{Ng10}*{\S~8}
with some small improvements.

\section{Generalities on Counting Points} 
\label{sec:generalities_on_counting_points}

We first review some general facts about counting points following
\cite{Ng10}*{\S~8.1}. In this book, we need more general coefficients than the
constant sheaf \(\Qlb\), and because there has been some significant development
on the algebro-geometric tools used for these purposes, we will base our
discussions on those frameworks as well.

\subsection{}
First we let \(f\colon\sX\to\Spec{k}\) be a \(k\)-variety, or more generally an
algebraic stack of finite type over a finite field \(k\), and
\(\cF\in\Dbc(\sX,\Qlb)\) a bounded constructible lisse-\'etale complex on
\(\sX\) with coefficients in \(\Qlb\). Then for any \(k\)-point \(x\in \sX(k)\)
and a fixed geometric point \(\bar{x}=\Spec{\bar{k}}\) over \(x\), we have a
continuous representation of \(\Gal(\bar{k}/k)\) on the stalk \(\cF_{\bar{x}}\).
Let Frobenius element \(\Frob_k\in\Gal(\bar{k}/k)\) be the one defined by
taking \(q=p^m\)-th root, where \(q\) is the cardinality of \(k\), then we may
define trace function as
\begin{align}
    \Tr_{\cF}\colon \sX(k)&\longto \Qlb\\
    x&\longmapsto \sum_{i\in\bbZ}(-1)^i\Tr\bigl(\Frob_k,\RH^i(\cF_{\bar{x}})\bigr).
    \nomenclature[\(Tr_F"cal \)]{\(\Tr_{\cF}\)}{the trace function induced by a
    bounded complex \(\cF\)}
\end{align}
If the complex \(\cF\) is not bounded, then the above construction still makes
sense, except that the trace function now takes values only as infinite series,
but one can impose convergence condition if so desired. For our purposes, we are
mostly interested in Deligne-Mumford stacks of finite types, therefore
boundedness condition is preserved by any six-functor operation, so we do not
worry about potential issues caused by unboundedness.

A classical result is Grothendieck--Lefschetz trace formula, which states that
(for \(\sX\) a \(k\)-scheme)
\begin{align}
    \label{eqn:gro_lef_trace_formula_abstract}
    \Cnt_\cF{\sX(k)}&\defeq \sum_{x\in\sX(k)}\Tr_{\cF}(x)=\Tr_{f_!\cF}(\Spec{k}).
    \nomenclature[\(\#_F"cal \)]{\(\Cnt_{\cF},\Cnt_\tau\)}{the weighted point count induced by a
    bounded complex \(\cF\) or a function \(\tau\)}
\end{align}
This is a special case of the relative version of the formula
\begin{align}
    \label{eqn:gro_lef_trace_formula_abstract_relative}
    \Cnt_\cF{\sX(k)}=\Cnt_{f_!\cF}{\sY(k)}
\end{align}
for a morphism \(f\colon\sX\to\sY\) between \(k\)-schemes of finite types.

A perhaps more common form of \eqref{eqn:gro_lef_trace_formula_abstract} is the
following: if \(k\) has \(q=p^m\) elements, let \(\Frob\) be the \(m\)-fold
iteration of the absolute Frobenius on \(\sX\), then \(\Frob\) is defined over
\(k\). It is a (non-trivial) fact that there is a functorial isomorphism
\(\iota\colon\cF\simeq\Frob^*\cF\) for any
\(\cF\), see \cite{StacksP}*{\href{https://stacks.math.columbia.edu/tag/03SL}{Tag 03SL}}.
The morphism \(\Frob\) is finite, hence proper, so we have adjunction map
\begin{align}
    \cF\longto \Frob_*\Frob^*\cF\simeq\Frob_!\Frob^*\cF.
\end{align}
Pushing forward using \(f_!\), we see that
\(\Frob\) induces a map on cohomological groups
\begin{align}
    \Frob\otimes_k\bar{k}\colon
    \RHc^i(\sX_{\bar{k}},\cF)\longto\RHc^i(\sX_{\bar{k}},\Frob^*\cF)
    \simeq\RHc^i(\sX_{\bar{k}},\cF),
\end{align}
where the second isomorphism is defined using the canonical identification
\(\iota\) above. It turns out that this map is the same map induced by
\(\Frob_k\)-action on \(f_!\cF\), and so we have
\begin{align}
    \label{eqn:gro_lef_trace_formula_basic_schemes}
    \Cnt_\cF{\sX(k)}
    =\sum_{i\in\bbZ}(-1)^i\Tr\bigl(\Frob\otimes_k\id_{\bar{k}},\RHc^i(\sX_{\bar{k}},\cF)\bigr).
\end{align}

\subsection{}
In \cite{Be93}, Behrend generalizes
\eqref{eqn:gro_lef_trace_formula_abstract_relative}
to the case where \(\sX\) is a smooth Deligne--Mumford stack of finite type over
\(k\), and \(\cF\) is the constant sheaf. In the stack case, the
\(\cF\)-weighted point counting is defined as
\begin{align}
    \Cnt_\cF{\sX(k)}\defeq
    \sum_{x\in\sX(k)/{\sim}}\frac{\Tr_\cF(x)}{\Cnt\Aut_{\sX}(x)(k)},
\end{align}
where \(x\) ranges over isomorphism classes of \(\sX(k)\). Note that since
\(\Aut_x\) is a \(k\)-group scheme of finite type and \(k\) is a finite field,
the above definition makes sense even for Artin stacks.
The result is further expanded to
smooth Artin stacks and arbitrary complexes in \cite{Be03}, assuming certain
technical boundedness or convergence condition is met.

In \cites{LaOl08I,LaOl08II}, Laszlo and Olsson developed six-functor formalism
for Artin stacks, and the smoothness assumption in Behrend's results is removed
under these frameworks by Sun in \cite{Sun12} because of the newly-available
duality. As noted by Sun in \cite{Sun12},
an important difference between \cite{Be03} and \cite{Sun12} is that the former
uses arithmetic Frobenius and ordinary cohomology (probably due to lack of
duality results at the time), while in the latter
the geometric version is used. Because in \cite{Sun12} the author considers a
very general setup, some complicated ``stratifiablility'' condition and a
convergence condition are used. However, for our purposes, we only consider
bounded complexes on Deligne--Mumford stacks, so those conditions are
automatically met. Therefore, to summarize, for any morphism
\(f\colon\sX\to\sY\) of Deligne--Mumford stacks of finite types over \(k\) and
any \(\cF\in\Dbc(\sX,\Qlb)\), the trace formula
\eqref{eqn:gro_lef_trace_formula_abstract_relative} always holds.

\subsection{}
In \cite{LN08}*{Appendice~A.3} and \cite{Ng10}*{\S~8.1}, an equivariant
version of \eqref{eqn:gro_lef_trace_formula_basic_schemes} is
developed with specific applications to affine Springer fibers and Hitchin
fibrations in mind. Here we reformulate those results using the general trace
formula \eqref{eqn:gro_lef_trace_formula_abstract_relative} for Deligne--Mumford
stacks.

Suppose \(\cM\) is a Deligne--Mumford stack of finite type over \(k\) together
with an action of a commutative Deligne--Mumford group stack \(\cP\)
of finite type over \(k\).
The groupoids \(\cM(\bar{k})\) and \(\cP(\bar{k})\) carry natural actions of
geometric Frobenius \(\Frob\), and the \(2\)-categorical quotient
\(\sX(\bar{k})\) is the \(2\)-category where:
\begin{enumerate}
    \item the objects are the objects \(m\in \cM(\bar{k})\);
    \item the \(1\)-morphisms \(m\to m'\) are pairs \((p,f)\)
        where \(p\in \cP(\bar{k})\) and \(f\in\Hom_{\cM(\bar{k})}(p(m),m')\);
    \item the \(2\)-morphisms \((p,f)\Rightarrow (p',f')\) are elements
        \(j\in\Hom_{\cP(\bar{k})}(p,p')\) such that the composition of
        \(j\colon p(m)\to p'(m)\) with \(f'\colon p'(m)\to m'\) is equal to
        \(f\colon p(m)\to m'\).
\end{enumerate}
According to \cite{Ng06}*{Lemme~4.7}, such \(2\)-category is equivalent to
a groupoid if and only if for any \(m\in \cM(\bar{k})\), the homomorphism induced
by \(\cP\) action
\begin{align}
    \Aut_{\cP(\bar{k})}(1_\cP)\longto \Aut_{\cM(\bar{k})}(m)
\end{align}
is injective. In view of product formula \Cref{prop:product_formula_alt},
this condition is always met over the reduced locus
(cf.~\Cref{def:mH_fibration_reduced_locus}) of mH-fibrations.
When this condition is met, then a \(1\)-morphism \((p,f)\) has only the
trivial \(2\)-auto\-mor\-phism, thus \(\sX(\bar{k})\) is equivalent to the groupoid
where:
\begin{enumerate}
    \item the objects are the objects \(m\in \cM(\bar{k})\);
    \item the morphisms \(m\to m'\) are isomorphism classes of pairs \((p,f)\).
\end{enumerate}
In this case, we omit \(f\) and simply write a morphism as \(p\colon m\to m'\),
or equivalently we have canonical isomorphism \(p(m)\simeq m'\).

We will now always assume that \(\sX(\bar{k})\) is equivalent to a groupoid. The
category \(\sX(k)\) of fixed points under the action of \(\Frob\) is
as follows:
\begin{enumerate}
    \item the objects are triples \((m,p,f)\) such that \(f\colon\Frob(m)\to p(m)\)
        is an isomorphism;
    \item a morphism \(h\colon (m,p,f)\to(m',p',f')\) is a pair \((h,\phi)\) where
        \(h\in\cP(\bar{k})\) and \(\phi\colon hm\to m'\) is an isomorphism, such
        that \(\Frob(h)ph^{-1}=p'\) and
        \(\Frob(\phi)\colon\Frob(hm)\to\Frob(m')\) is equal to the map
        \(f'^{-1}\circ (p'\phi)\circ(\Frob(h)f)\);
    \item the sets of \(2\)-morphisms are guaranteed either empty or singletons
        by assumption, making this category a groupoid.
\end{enumerate}
\begin{remark}
    \label[remark]{rmk:simplify_convoluted_condition_on_morphisms_in_sX_k}
    If we use \(f\) to identify \(\Frob(m)\) and \(p(m)\) and similarly for
    \(f'\), then the last condition on \((h,\phi)\) can be simplified to
    \(\Frob(\phi)=p'(\phi)\). If \(p=p'\), then the conditions further simplify
    to \(h\in\cP(k)/{\sim}\) and \(\Frob(\phi)=p(\phi)\). In this case,
    \(\phi\) and \(\phi'\) are equivalent if \(\phi'=\alpha\phi\) for some
    \(\alpha\in\Aut_{\cP}(1_\cP)\).
\end{remark}

Since \(\sX(k)\) has only finitely many isomorphism classes, and each object has
only finitely many automorphisms,
if \(\tau\) is a \(\Qlb\)-valued function on the set of isomorphism classes
\(\sX(k)/{\sim}\), the following sum makes sense
\begin{align}
    \Cnt_\tau{\sX(k)}=\sum_{x\in\sX(k)/{\sim}}\frac{\tau(x)}{\Cnt{\Aut_{\sX(k)}(x)}}.
\end{align}

\subsection{}
Let \(x=(m,p)\in\sX(k)\) as above, then the \(\Frob\)-conjugacy class of \(p\)
is determined by the isomorphism class of \(x\). Let \(P\) be the coarse space
of \(\cP\), then \(P\) is a commutative group scheme, and we have isomorphism
between the group of geometric connected components
\begin{align}
    \pi_0(\cP)\simeq\pi_0(P).
\end{align}
According to Lang's theorem, the group of \(\Frob\)-conjugacy classes of \(P(\bar{k})\)
can be identified with
\begin{align}
    P_\Frob\simeq\RH^1(k,P)\simeq\RH^1(k,\pi_0(P)).
\end{align}
This shows that \(\cP_\Frob\simeq \RH^1(k,\pi_0(\cP))\) as well, the former
denoting the group of the \(\Frob\)-conjugacy classes of the
isomorphism classes in \(\cP(\bar{k})\).
Therefore, given \(x\in\sX(k)\), we define \(\cl(x)\)
    \nomenclature[\(cl_x \)]{\(\cl(x)\)}{the Frobenius conjugacy class in
    \(\pi_0(\cP)\) of \(x\in\Stack*{\cM/\cP}(\bar{k})\)}
be the corresponding
\(\Frob\)-conjugacy class in \(\pi_0(\cP)\) also viewed as an element in
\(\RH^1(k,\pi_0(\cP))\). The class \(\cl(x)\) depends only on the isomorphism
class of \(x\). For any \(\Frob\)-invariant character
\(\kappa\colon\pi_0(\cP)\to\Qlb^\x\), we have a pairing
\begin{align}
    \Pair{\cl(x)}{\kappa}=\kappa(\cl(x))\in\Qlb^{\x}.
\end{align}
For any \(\Qlb\)-valued function \(\tau\) on the set of isomorphism classes \(\sX(k)/{\sim}\), we define
\((\tau,\kappa)\)-weighted point counting
\begin{align}
    \Cnt_{\tau}\sX(k)_\kappa
    =\sum_{x\in\sX(k)/{\sim}}\frac{\Pair{\cl(x)}{\kappa}\tau(x)}{\Cnt\Aut_{\sX(k)}(x)}.
\end{align}

\subsection{}
Let \(\cF\in\Dbc(\sX,\Qlb)\) now be \(\cP\)-equivariant, in other words, we are supplied with a fixed
isomorphism \(\psi\colon a^*\cF\simeq\Qlb\boxtimes\cF\) on \(\cP\x\cM\), where \(a\) is the
action map \(\cP\x\cM\to\cM\), and the usual cocycle and identity axioms of
equivariance is met for \(\phi\). In this way, the objects of \(\cP(\bar{k})\)
acts on the sheaf \(\cF|_{\cM_{\bar{k}}}\), hence also on any cohomology
groups.

According to \Cref{lem:homotopy_lemma,rmk:apply_homotopy_to_mH}, the induced
action on \(\RHc^\bullet(\cM_{\bar{k}},\cF)\) factors through \(\pi_0(\cP)\).
The group
\(\cP\) does not need to be smooth because we always have \(\cP^\Red\) to be
smooth. Thus, if \(\kappa\) is a
\(\Frob\)-invariant character of \(\pi_0(\cP)\), then \(\Frob\) acts on the
\(\kappa\)-isotypic subspace \(\RHc^\bullet(\cM_{\bar{k}},\cF)_\kappa\).

Finally, by equivariance, the
trace function \(\Tr_\cF\) descends to a function on groupoid \(\sX(k)\),
still denoted by \(\Tr_\cF\), as follows: for any pair \((m,p)\) such that
\(\Frob(m)\simeq p(m)\), the equivariance structure \(\psi\) identifies the
stalk of \(\cF\) at \(m\) with that at \(\Frob(m)\), hence defining an action of
\(\Frob\) on that stalk, and we can take the trace. Now we have the
following variant of Grothendieck--Lefschetz trace formula, which is a
generalization of \cite{LN08}*{Proposition~A.3.1} (see also
\cite{Ng10}*{Proposition~8.1.6}).

\begin{proposition}
    \label[proposition]{prop:M_mod_P_pt_counting_general}
    Let \(\cM\) be a Deligne--Mumford stack of finite type over \(k\) with an
    action of a commutative Deligne--Mumford group stack \(\cP\) of finite type
    over \(k\), and suppose \(\sX(\bar{k})=\Stack{\cM(\bar{k})/\cP(\bar{k})}\)
    is equivalent to a groupoid. Let \(\cF\) be a bounded constructible
    \(\cP\)-equivariant complex. Then for any \(\Frob\)-invariant character 
    \(\kappa\colon\pi_0(\cP)\to\Qlb^\x\) we have equality
    \begin{align}
        \Cnt{\cP_0(k)}\Cnt_\cF{\sX(k)}_\kappa=\sum_{i\in\bbZ}(-1)^i\Tr\bigl(\Frob,\RHc^i(\cM_{\bar{k}},\cF)_\kappa\bigr),
    \end{align}
    where \(\cP_0\) is the neutral component of \(\cP\).
\end{proposition}
\begin{proof}
    Using Fourier transform on the finite group \(\pi_0(\cP)\), we have that
    \begin{align}
        \sum_{i\in\bbZ}(-1)^i\Tr\bigl(\Frob,\RHc^i(\cM_{\bar{k}},\cF)_\kappa\bigr)
        =\frac{1}{\Cnt\pi_0(\cP)}\sum_{p\in\pi_0(\cP)}
        \kappa(p)\sum_{i\in\bbZ}(-1)^i\Tr\bigl(p^{-1}\circ\Frob,\RHc^i(\cM_{\bar{k}},\cF)\bigr).
    \end{align}
    Note that if a character \(\chi\) is not \(\Frob\)-invariant, then its
    isotypic space has no contribution to the right-hand side, which
    is then equal to
    \begin{align}
        \label{eqn:M_mod_P_pt_counting_tmp}
        \frac{1}{\Cnt\pi_0(\cP)_\Frob}\sum_{p\in\pi_0(\cP)_\Frob}
        \kappa(p)\sum_{i\in\bbZ}(-1)^i\Tr\bigl(\dot{p}^{-1}\circ\Frob,\RHc^i(\cM_{\bar{k}},\cF)\bigr),
    \end{align}
    where \(\dot{p}\) is any choice of representative of \(\Frob\)-conjugacy
    class \(p\) in \(\cP\).

    In the proof of \cite{LN08}*{Proposition~A.3.1}, another variant of
    \eqref{eqn:gro_lef_trace_formula_basic_schemes} by Deligne and Lusztig is
    used (see \cite{DeLu76}). The idea is that if we assume that \(\cM\)
    is a quasi-projective scheme, then \(\dot{p}^{-1}\circ\Frob\) is the Frobenius map
    of another \(k\)-model of \(\cM_{\bar{k}}\). Here \(\cM\)
    is a Deligne--Mumford stack, but we can still descend \(\cM_{\bar{k}}\) to a
    \(k\)-model \(\cM^{\dot{p}}\) such that \(\dot{p}^{-1}\circ\Frob\) is the induced Frobenius.
    Since \(\cP\) is commutative, the action
    \(\cP_{\bar{k}}\x\cM_{\bar{k}}\to \cM_{\bar{k}}\) is equivariant with
    respect to the usual \(\Frob\)-action on \(\cP_{\bar{k}}\) and the
    \(\dot{p}^{-1}\circ\Frob\)-action on \(\cM_{\bar{k}}\), hence descends to an
    action of \(\cP\) on \(\cM^{\dot{p}}\). Moreover, the identity map on
    \(\cM_{\bar{k}}\) induces an isomorphism of stacks over \(k\):
    \begin{align}
        \label{eqn:M_mod_P_isom_to_M_dotp_mod_P}
        \Stack{\cM/\cP}\simeq \Stack{\cM^{\dot{p}}/\cP}.
    \end{align}
    This shows that \(\cF\) may be canonically identified with a
    \(\cP\)-equivariant complex on \(\cM^{\dot{p}}\).
    Thus, using \eqref{eqn:gro_lef_trace_formula_basic_schemes} for
    Deligne--Mumford stacks, \eqref{eqn:M_mod_P_pt_counting_tmp} is equal to
    \begin{align}
        \frac{1}{\Cnt\pi_0(\cP)_\Frob}\sum_{p\in
        \pi_0(\cP)_\Frob}\sum_{\substack{m\in\cM(\bar{k})/{\sim}\\\Frob(m)=\dot{p}(m)}}
        \frac{\kappa(p)\Tr_\cF(m)}{\Cnt\Set*{\phi\in
        \Aut_{\cM(\bar{k})}(m)\given \Frob(\phi)=\dot{p}(\phi)}},
    \end{align}
    where \(\Tr_\cF\) is the trace function of \(\cF\) discussed in the
    paragraph preceding this proposition.
    Let \(P\) be the coarse space of \(\cP\). According to
    \Cref{rmk:simplify_convoluted_condition_on_morphisms_in_sX_k}, we have
    \begin{align}
        \Cnt{P(k)}\Cnt\Set*{\phi\in \Aut_{\cM(\bar{k})}(m)\given
        \Frob(\phi)=\dot{p}(\phi)}=\Cnt{\Aut_{\cP}(1_\cP)(k)}\Cnt{\Aut_{\sX(k)}(m,\dot{p})}.
    \end{align}
    Note that \(\pi_0(\cP)_\Frob\) and \(\pi_0(\cP)^\Frob=\pi_0(\cP)(k)\) have
    the same cardinality, and \(P(k)\to\pi_0(\cP)(k)\) is surjective by Lang's
    theorem. So we obtain
    \begin{align}
        \sum_{i\in\bbZ}(-1)^i\Tr\bigl(\Frob,\RHc^i(\cM_{\bar{k}},\cF)_\kappa\bigr)
        =\Cnt{\cP_0(k)}\sum_{p\in
        \pi_0(\cP)_\Frob}\sum_{\substack{m\in\cM(\bar{k})/{\sim}\\\Frob(m)=\dot{p}(m)}}
        \frac{\kappa(p)\Tr_\cF(m)}{\Cnt{\Aut_{\sX(k)}(m,\dot{p})}}.
    \end{align}
    Finally, note that the set of pairs \((m,\dot{p})\) in the summations above
    is in bijection with the isomorphism classes in \(\sX(k)\) by definition,
    thus we reach the desired equality
    \begin{align}
        \sum_{i\in\bbZ}(-1)^i\Tr\bigl(\Frob,\RHc^i(\cM_{\bar{k}},\cF)_\kappa\bigr)
        =\Cnt{\cP_0(k)}\sum_{(m,p)\in\sX(k)/{\sim}} \frac{\kappa(p)\Tr_\cF(m)}{\Cnt{\Aut_{\sX(k)}(m,p)}}.
    \end{align}
    This finishes the proof.
\end{proof}

\subsection{}
\label{sub:pt_counting_locally_finite_type_variant}
So far we only considered the case where \(\cM\) and \(\cP\) are of finite
types, but to apply our results to multiplicative affine Springer fibers, we need
another variant for locally finite type cases. Here we no longer need to
consider Deligne--Mumford stacks but only schemes. Let \(M\) be a \(k\)-scheme
with an action of a commutative \(k\)-group scheme \(P\), both locally of finite
types. We assume they satisfy the following assumptions:
\begin{enumerate}
    \item The \(\bar{k}\)-points of the
        group of connected components \(\pi_0(P)\) is a finitely generated
        abelian group.
    \item The stabilizer of any geometric point of \(M\) in \(P\) is
        of finite type over \(k\).
    \item We can find a torsion-free discrete subgroup \(\Lambda\subset P\) such
        that both \(P/\Lambda\) and \(M/\Lambda\) are of finite types.
\end{enumerate}
Note that in the above conditions since stabilizers of points in \(M\) are of
finite type, the action of \(\Lambda\) on \(M\) is necessarily free.
Note also that such \(\Lambda\) always exists if we only require \(P/\Lambda\)
to be of finite type. Indeed, let \(\Lambda_0\) be the largest free quotient
of \(P_{\bar{k}}^\Red\), then the kernel
\begin{align}
    P^{\tfree}=\ker(P_{\bar{k}}^\Red\longto\Lambda_0)
\end{align}
is of finite type. Since \(\Lambda_0\) is free, we may pick an arbitrary lifting
\(\gamma\colon\Lambda_0\to P_{\bar{k}}\), which is necessarily defined over some
finite extension \(k'/k\).

Since \(P^{\tfree}\) is of finite type, the group \(P^{\tfree}(k')\) is
finite. This means that when \(N\) is divisible by \(\Cnt{P^{\tfree}(k')}\), the
restriction of \(\gamma\) to \(\Lambda=N\Lambda_0\) is \(\Frob\)-equivariant
hence defined over \(k\), and \(P/\Lambda\) is clearly of finite type.
The condition that \(M/\Lambda\) is finite type is
clearly independent of the choice of \(\Lambda\), so the third condition above
is equivalent to saying \(M/\Lambda\) is of finite type for any specific
choice of \(\Lambda\).

From now on, since we only care about counting points, we
will use \(M\) and \(P\) to denote the respective sets of their
\(\bar{k}\)-points. Because the points of \(M\) and \(P\) have no
automorphisms, this notation will not cause any confusion.

\subsection{}
Consider quotient stack \(\sX=\Stack{M/P}\), whose groupoid of \(k\)-points
is as follows:
\begin{enumerate}
    \item the objects are pairs \((m,p)\) where \(m\in M\) and \(p\in
        P\) such that \(p(m)=\Frob(m)\);
    \item the morphisms \((m,p)\to (m',p')\) are \(h\in P\) such that
        \(hm=m'\) and \(p'=\Frob(h)ph^{-1}\).
\end{enumerate}
Given \(x=(m,p)\in\sX(k)\), the \(\Frob\)-conjugacy class \(\cl(x)\) of \(p\)
depends only on the isomorphism class of \(x\). Since \(p(m)=\Frob(m)\) and \(m\) is
defined over some finite extension \(k'/k\), the class \(\cl(x)\) is necessarily
torsion, hence it lies in the torsion subgroup of \(P_\Frob\), which is
identified with \(\RH^1(k,P)\). By \cite{Ng10}*{Lemme~8.1.12}, any character
\(\kappa\) of \(\RH^1(k,P)\) can be extended to a torsion character
\(\tilde{\kappa}\) of \(P_\Frob\).

Since we have equivalence of quotient stacks
\begin{align}
    \sX=\Stack{M/P}=\Stack{(M/\Lambda)/(P/\Lambda)},
\end{align}
we know that \(\sX(k)\) has only finitely many isomorphism classes,
and the automorphism group of each object is also finite. Therefore, for any
character \(\kappa\colon \RH^1(k,P)\to\Qlb^\x\), and any
\(\Qlb\)-valued function \(\tau\) on the isomorphism classes of \(\sX(k)\),
we may consider sum
\begin{align}
    \Cnt_\tau\sX(k)_\kappa
    =\sum_{x\in\sX(k)/{\sim}}\frac{\Pair{\cl(x)}{\kappa}\tau(x)}{\Cnt\Aut_{\sX(k)}(x)}.
\end{align}
Let \(\tilde{\kappa}\) be an extension of \(\kappa\) to \(P(\bar{k})_\Frob\),
which inflates to a \(\Frob\)-invariant character of \(P\). Since
\(\tilde{\kappa}\) has finite order, we can choose \(\Lambda\) so that the
restriction of \(\tilde{\kappa}\) to \(\Lambda\) is trivial. So
\(\tilde{\kappa}\) descends to a character of \(P/\Lambda\), and if
\(\Lambda'\subset\Lambda\) is a \(\Frob\)-stable sublattice of finite index,
\(\tilde{\kappa}\) induces a character of \(P/\Lambda'\) too.

Let \(\cF\) be a bounded locally constructible complex on \(M\) that is
\(P\)-equivariant, and still denote by \(\cF\) its descent to \(M/\Lambda\).
As before, \(\Tr_\cF\) descends to a function on \(\sX(k)\), also
denoted by \(\Tr_\cF\). Let \(\RHc^\bullet(M/\Lambda,\cF)_{\tilde{\kappa}}\) be the
\(\tilde{\kappa}\)-isotypic direct summand of \(\RHc^\bullet(M/\Lambda,\cF)\).
The following is a generalization of \cite{Ng10}*{Proposition~8.1.13}.
\begin{proposition}
    \label[proposition]{prop:pt_counting_locally_finite_type_variant}
    We have the following equality
    \begin{align}
        \Cnt{P_0(k)}\Cnt_\cF{\sX(k)}_\kappa
        =\sum_{i\in\bbZ}(-1)^i\Tr\bigl(\Frob,\RHc^i(M/\Lambda,\cF)_{\tilde{\kappa}}\bigr).
    \end{align}
    Moreover, if \(\Lambda'\subset\Lambda\) is a \(\Frob\)-stable sublattice of
    finite index, then we have for each \(i\in\bbZ\) a canonical isomorphism
    \begin{align}
        \RHc^i(M/\Lambda,\cF)_{\tilde{\kappa}}
        \stackrel{\sim}{\longto}\RHc^i(M/\Lambda',\cF)_{\tilde{\kappa}}.
    \end{align}
\end{proposition}
\begin{proof}
    Since \(\Lambda\) is free, we have \(P_0\simeq (P/\Lambda)_0\). Then the
    first formula is a consequence of \Cref{prop:M_mod_P_pt_counting_general}.
    If we replace \(k\) by a finite extension \(k'\) of degree \(e\), then for
    any \(x'\in \sX(k')\) we have a \(\Frob^e\)-conjugacy class in \(P\). Since
    \(\kappa\), \(\Lambda\), etc. are \(\Frob\)-invariant, they are also
    \(\Frob^e\)-invariant, so we also have
    \begin{align}
        \Cnt{P_0(k')}\Cnt_\cF{\sX(k')}_\kappa
        =\sum_{i\in\bbZ}(-1)^i\Tr\bigl(\Frob^e,\RHc^i(M/\Lambda,\cF)_{\tilde{\kappa}}\bigr).
    \end{align}
    Since \(e\) is arbitrary, it means that the map
    \begin{align}
        \bigoplus_{i}(-1)^i\RHc^i(M/\Lambda,\cF)_{\tilde{\kappa}}
        \longto\bigoplus_{i}(-1)^i\RHc^i(M/\Lambda',\cF)_{\tilde{\kappa}}
    \end{align}
    is an isomorphism of \(\Frob\)-modules which also respects grading \(i\).
    So we have the second claim.
\end{proof}

\begin{corollary}
    Let \(\sX=\Stack{M/P}\) (resp.~another such groupoid
    \(\sX'=\Stack{M'/P'}\)), and \(\kappa\) (resp.~\(\kappa'\)) be a \(\Frob\)-invariant
    character of \(P\) (resp.~\(P'\)) of finite order. Let \(\cF\) (resp.~\(\cF'\)) be a
    \(P\)-equivariant (resp.~\(P'\)-equivariant) bounded locally constructible
    complex on \(M\) (resp.~\(M'\)). Suppose there exists some \(N\) such
    that for all \(e>N\) and \(k'/k\) a finite extension of degree \(e\) we
    always have
    \begin{align}
        \Cnt{P_0(k')}\Cnt_\cF
        \sX(k')_\kappa=\Cnt{P'_0(k')}\Cnt_{\cF'}\sX'(k')_{\kappa'},
    \end{align}
    then we have
    \begin{align}
        \Cnt{P_0(k)}\Cnt_\cF
        \sX(k)_\kappa=\Cnt{P'_0(k)}\Cnt_{\cF'}\sX'(k)_{\kappa'}.
    \end{align}
\end{corollary}
\begin{proof}
    After choosing appropriate lattices \(\Lambda\subset P\) and
    \(\Lambda'\subset P'\), the assumptions imply that
    \(\RHc^\bullet(M/\Lambda,\cF)_{\tilde{\kappa}}\) and
    \(\RHc^\bullet(M'/\Lambda',\cF')_{\tilde{\kappa}'}\) are isomorphic
    \(\Frob\)-modules up to semisimplification. Taking the traces of \(\Frob\),
    and we get the result.
\end{proof}


\section{Counting Points on Multiplicative Affine Springer Fibers} 
\label{sec:counting_points_on_generalized_affine_springer_fibers}

In this section, we apply the general results from the previous section to
multiplicative affine Springer fibers. We use the notations from
\Cref{chap:generalized_affine_springer_fibers} which we now review.
Let \(X_v=\Spec{\cO_v}\) be the
formal disc around a closed point \(v\in \abs{X}\), and
\(X_v^\bullet=\Spec{F_v}\) is the punctured disc. Let \(\pi\) be a fixed choice
of uniformizer in \(\cO_v\) and \(k_v\) be the residue field. Let \(\bar{k}\) be
a fixed algebraic closure of \(k\). Denote
\(\breve{X}_v=\Spec{\breve{\cO}_v}\) where
\(\breve{\cO}_v=\cO_v\hat{\otimes}_k\bar{k}\) and similarly
\(\breve{X}_v^\bullet=\Spec{\breve{F}_v}\). For any \(k\)-embedding \(\bar{v}\colon
k_v\to\bar{k}\), we let \(\breve{X}_{\bar{v}}=\Spec{\breve{\cO}_{\bar{v}}}\) be
the component of \(\breve{X}_v\) containing \(\bar{v}\). If \(k_v=k\), then
\(\breve{X}_{\bar{v}}=\breve{X}_v\); in general, we have
\(\breve{X}_v=\prod_{\bar{v}}\breve{X}_{\bar{v}}\) where \(\bar{v}\) ranges over
all \(k\)-embeddings of \(k_v\) in \(\bar{k}\).

\subsection{}
Let \(\FRM\in\FM(G^\SC)\) be a very flat reductive monoid over \(X_v\).
Let \(\cL\) be a \(Z_\FRM\)-torsor and
\(a\in\FRC_{\FRM,\cL}(\cO_v)\cap\FRC_{\FRM,\cL}^\x(F_v)^\rss\), then there
exists some \(\gamma_\FRM\in\FRM_\cL^\x(F_v)\) with
\(\chi_{\FRM,\cL}(\gamma_\FRM)=a\). We have the multiplicative affine Springer fiber
\(\cM_v(a)\) (see \Cref{sec:definition_and_generalities}) defined using \(\gamma_\FRM\),
whose set of \(\bar{k}\)-points is
\begin{align}
    \prod_{\bar{v}\colon k_v\to \bar{k}}\Set*{g\in G(\breve{F}_{\bar{v}})/G(\breve{\cO}_{\bar{v}})\given \Ad_g^{-1}(\gamma_\AD)\in
    G^\AD(\breve{\cO}_{\bar{v}})\pi^{\lambda_v}G^\AD(\breve{\cO}_{\bar{v}})},
\end{align}
where \(\gamma_\AD\) is the image of \(\gamma_\FRM\) in \(G^\AD\), and
\(\lambda_v\) is the boundary divisor of \(a\), which may also be viewed as
an \(F_v\)-rational dominant cocharacter in \(\CoCharG(T^\AD)\) through map
\begin{align}
    \FRA_\FRM\to\FRA_{\Env(G^\SC)}.
\end{align}
We denote the Newton point of \(\gamma_\AD\) by \(\nu_v\). Since we only care
about point-counting, we will replace \(\cM_v(a)\) with its reduced subfunctor
and still use notation \(\cM_v(a)\) for simplicity. We know that the reduced
functor \(\cM_v(a)\) is represented by a \(k\)-scheme locally of finite type
(\Cref{thm:GASF_dimension_formula}).

\subsection{}
Following \cite{Ng10}, let \(\FRJ_a\) be the regular centralizer at \(a\), and
\(\FRJ_a'\) be a smooth group scheme over \(X_v\) with connected fibers and a
homomorphism \(\FRJ_a'\to\FRJ_a\) such that its restriction to \(X_v^\bullet\)
is an isomorphism. Note that \(\FRJ_a'\) is necessarily commutative.

Consider reduced commutative group scheme over \(k\)
\begin{align}
    \cP_v'(a)^\Red\defeq\Gr_{\FRJ_a'}^\Red=(\Loop{\FRJ_a'}/\Arc{\FRJ_a'})^\Red.
\end{align}
The second equality above is due to the fact that \(\FRJ_a'\) has connected
special fiber. This is a group scheme locally of finite type over \(k\), and we
use \(\cP_v'(a)\) for simplicity. Similarly, we have \(\cP_v(a)\) for
\(\FRJ_a\), and we have homomorphism
\begin{align}
    \cP_v(a)'\longto\cP_v(a).
\end{align}
The action of \(\cP_v(a)\) on \(\cM_v(a)\) induces an action of \(\cP_v'(a)\) on
\(\cM_v(a)\). By \Cref{prop:projective_quotient_of_GASF}, this action satisfies
conditions laid out in the beginning of
\Cref{sub:pt_counting_locally_finite_type_variant}, and so we may express
different \(\kappa\)-weighted point-countings using trace formula on
cohomologies. As in \Cref{sub:pt_counting_locally_finite_type_variant}, we also
use \(\cM_v(a)\), \(\cP_v(a)\), etc. to denote their respective sets of
\(\bar{k}\)-points.

\subsection{}
We now try to connect point-countings on stack \(\sX'=\Stack{\cM_v(a)/\cP_v'(a)}\)
with local orbital integrals. First, the following lemma relates the \(\kappa\)
in orbital integrals with \(\kappa\) in
\Cref{sec:generalities_on_counting_points}:
\begin{lemma}[\cite{Ng10}*{Lemme~8.2.4}]
    Assuming \(\FRJ_a'\) has connected special fiber, then we have canonical
    isomorphism
    \begin{align}
        \RH^1(F_v,\FRJ_a)\simeq\RH^1(k,\cP_v'(a)).
    \end{align}
\end{lemma}
By this lemma, any character \(\kappa\colon\RH^1(F_v,\FRJ_a)\to\Qlb^\x\) induces
a character of \(\RH^1(k,\cP_v'(a))\) and vice versa, and we still use
\(\kappa\) to denote the character of the latter. Let \(\cP_v'(a)_\Frob\) be the
group of \(\Frob\)-conjugacy classes of \(\cP_v'(a)\), and
\(\RH^1(k,\cP_v'(a))\) is the subgroup of torsion elements (due to continuity
requirement). Given \(x=(m,p)\in\sX'(k)\), we have associated class given by the
\(\Frob\)-conjugacy class of \(p\):
\begin{align}
    \cl(x)\in \RH^1(k,\cP_v'(a))\simeq\RH^1(F_v,\FRJ_a).
\end{align}
The element \(\gamma_\FRM\) is regular (as an \(F_v\)-point), so we have
canonical isomorphism \(\FRJ_a|_{F_v}\simeq I_{\gamma_\FRM}\subset G\).
Thus, \(\cl(x)\) may also
be regarded as an element of \(\RH^1(F_v, I_{\gamma_\FRM})\).

If \(g\in G(\breve{F}_v)\) is a representative of \(m\) and \(j\in
I_{\gamma_\FRM}(F_v)\) is a representative of \(p\), then \(p(m)=\Frob(m)\)
implies that \(\Frob(g)^{-1}j g\in G(\breve{\cO}_v)\), hence \(j\) is
\(\Frob\)-conjugate to \(1\) in \(G(\breve{F}_v)\). This shows that the image of
\(\cl(x)\) in \(\RH^1(F_v,G)\) is trivial. Note that here we are dependent on
a choice of \(\gamma_\FRM\), otherwise we cannot relate \(\RH^1(F_v,\FRJ_a)\)
with \(\RH^1(F_v,G)\).

We let \(\Gr^{\le\lambda_v}\) be the affine Schubert variety of the adjoint
group \(G^\AD\). The map of groupoids
\begin{align}
    \cM_v(a)(\bar{k})\longto
    \Stack*{\Arc{G^\AD}\backslash\Gr^{\le\lambda_v}}(\bar{k})
\end{align}
is \(\Frob\)-equivariant and \(\cP_v'(a)(\bar{k})\)-invariant, so any
function \(\tau\) on the set of isomorphism classes of
\begin{align}
    \Stack*{\Arc{G^\AD}\backslash\Gr^{\le\lambda_v}}(k)
\end{align}
induces a function on \(\sX'(k)\) by pullback.
Note that since \(G\) has connected fibers, the
said isomorphism classes are in bijection with double cosets
\begin{align}
    G(\cO_v)\pi^\mu G(\cO_v)
\end{align}
where \(\mu\in\CoCharG(T^\AD)\) is an \(F_v\)-rational dominant cocharacter with
\(\mu\le\lambda_v\). In particular, the function \(\tau\) can be the
trace function \(\Tr_\cF\) where \(\cF\) is
bounded and constructible \(\Arc{G^\AD}\)-equivariant complex on
\(\Gr^{\le\lambda_v}\).

\begin{proposition}
    \label[proposition]{prop:pt_counting_with_orbital_connected_fiber_case}
    Assuming \(\FRJ_a'\) has connected special fiber, and let
    \(\kappa\) be a character of \(\RH^1(F_v,\FRJ_a)\). Then for any
    function \(\tau\) on
    \(\Stack*{\Arc{G^\AD}\backslash\Gr^{\le\lambda_v}}(k)\), we have equality
    \begin{align}
        \Cnt_\tau\sX'(k)_\kappa=\vol\bigl(\FRJ_a'(\cO_v),\dd
        t_v\bigr)\OI_a^\kappa(\tau,\dd t_v),
    \end{align}
    where \(\dd t_v\) is any Haar measure on \(\FRJ_a(F_v)\).
\end{proposition}
\begin{proof}
    For any \(x\in\sX'(k)\), the class \(\cl(x)\) can be identified with an
    element in the kernel of the map
    \begin{align}
        \RH^1(F_v,I_{\gamma_\FRM})\longto \RH^1(F_v,G).
    \end{align}
    So we can decompose \(\sX'(k)\) into disjoint full subcategories
    \begin{align}
        \sX'(k)=\coprod_{\xi}\sX_\xi'(k),
    \end{align}
    where \(\xi\) ranges over the said kernel above, and \(\sX'(k)_\xi\)
    consists of objects \(x\) with \(\cl(x)=\xi\).

    For a fixed class \(\xi\), we pick a representative \(j_\xi\in
    I_{\gamma_\FRM}(\breve{F}_v)\). Let \(x=(m,p)\in\sX'_\xi(k)\), then we can
    find \(h\in\cP_v'(a)\) such that \(\Frob(h)ph^{-1}=j_\xi\). Replacing
    \((m,p)\) with an isomorphic object \((hm,j_\xi)\), we may always assume
    \(p=j_\xi\). In other words, we may restrict to an equivalent full
    subcategory of \(\sX_\xi'(k)\) consisting of objects \((m,j_\xi)\), and the
    morphisms from \((m,j_\xi)\) to \((m',j_\xi)\) are elements in
    \(h\in\cP_v'(a)(k)\) such that \(hm=m'\). Since \(\FRJ_a'\) has connected
    fibers, we have
    \begin{align}
        \cP_v'(a)(k)=\FRJ_a(F_v)/\FRJ_a'(\cO_v)\simeq
        I_{\gamma_\FRM}(F_v)/\FRJ_a'(\cO_v).
    \end{align}

    For an object \((m,j_\xi)\), pick a representative \(g\in G(\breve{F}_v)\),
    then \(j_\xi(m)=\Frob(m)\) implies that \(\Frob(g)^{-1}j_\xi g\in
    G(\breve{\cO}_v)\). Thus, we may pick \(g\) so that we actually have
    \(\Frob(g)^{-1}j_\xi g=1\), and if \(g\) and \(g'\) are two such choices,
    we must have \(g'=gg_0\) for some \(g_0\in G(\cO_v)\). So an object
    \((m,j_\xi)\) determines a unique element in quotient set
    \(G(\breve{F}_v)/G(\cO_v)\). Thus, the category \(\sX_\xi'(k)\) is
    equivalent to the following category \(O_\xi\):
    \begin{enumerate}
        \item the objects are elements \(g\in G(\breve{F}_v)/G(\cO_v)\), such
            that \(\Frob(g)^{-1}j_\xi g=1\) and
            \(\Ad_g^{-1}(\gamma_\FRM)\in\FRM_\cL(\breve{\cO}_v)\);
        \item a morphism \(g\to g_1\) is an element \(h\in \cP_v(a)(k)\) such
            that \(hg=g_1\), where the action of \(\cP_v'(a)(k)\) is induced by
            the isomorphism \(\FRJ_a\simeq I_{\gamma_\FRM}\) over
            \(X_v^\bullet\).
    \end{enumerate}

    Choose \(g_\xi\in G(\breve{F}_v)\) such that \(\Frob(g_\xi)^{-1}j_\xi
    g_\xi=1\), and let
    \begin{align}
        \gamma_\xi=\Ad_{g_\xi}^{-1}(\gamma_\FRM),
    \end{align}
    then we necessarily have \(\gamma_\xi\in\FRM_\cL^\x(F_v)\). One can verify that
    the \(G(F_v)\)-conjugacy class of \(\gamma_\xi\) does not depend on the
    choice of either \(j_\xi\) or \(g_\xi\), but only on the cohomology class
    \(\xi\in\RH^1(F_v,I_{\gamma_\FRM})\) (hence also implicitly depends on the
    choice of \(\gamma_\FRM\)).

    Let \(g'=g_\xi^{-1}g\), then \(O_\xi\) is equivalent to the following
    category \(O_\xi'\): the objects are \(g'\in G(F_v)/G(\cO_v)\), such that
    \(\Ad_{g'}^{-1}(\gamma_\xi)\in\FRM_\cL(\cO_v)\), and a morphism \(g'\to g_1'\)
    is an element \(h\in \cP_v'(a)(k)\) such that \(hg'=g_1'\). Here the
    action of \(\cP_v'(a)(k)\) is induced by the isomorphism \(\FRJ_a\simeq
    I_{\gamma_\xi}\) over \(X_v^\bullet\).

    Therefore, the isomorphism classes of \(\sX_\xi'(k)\) are in bijection with
    double cosets
    \begin{align}
        g'\in\dbq{I_{\gamma_\xi}(F_v)}{G(F_v)}{G(\cO_v)}
    \end{align}
    such that \(\Ad_{g'}^{-1}(\gamma_\xi)\in\FRM_\cL(\cO_v)\),  and automorphism
    group of \(g'\) is
    \begin{align}
        \bigl(I_{\gamma_\xi}(F_v)\cap g'G(\cO_v)(g')^{-1}\bigr)/\FRJ_a'(\cO_v).
    \end{align}
    As a result, for any Haar measure \(\dd t_v\) on \(\FRJ_a(F_v)\simeq
    I_{\gamma_\xi}(F_v)\) and any \(\Qlb\)-valued function \(\tau\) on
    \(\sX_\xi'(k)\), we have
    \begin{align}
        \Cnt_\tau\sX_\xi'(k)=\sum_{g'}\frac{\tau(g')\vol\bigl(\FRJ_a'(\cO_v),\dd
        t_v\bigr)}{\vol\bigl(I_{\gamma_\xi}(F_v)\cap g'G(\cO_v)(g')^{-1}\bigr)},
    \end{align}
    where \(g'\) ranges over the double cosets above. This implies that
    \begin{align}
        \Cnt_{\tau}{\sX_\xi'(k)}=\vol\bigl(\FRJ_a'(\cO_v),\dd
        t_v\bigr)\OI_{\gamma_\xi}(\tau,\dd t_v).
    \end{align}
    Summing over all classes \(\xi\), we have our result.
\end{proof}

\subsection{}
In general, \(\FRJ_a\) may have disconnected special fiber, and recall we
have open subgroup \(\FRJ_a^0\) of fiberwise neutral component.
We will connect the point-counting of
\(\sX(k)=\Stack*{\cM_v(a)/\cP_v(a)}(k)\) with that of
\(\sX^0(k)=\Stack*{\cM_v(a)/\cP_v^0(a)}(k)\), where \(\cP_v^0(a)\) is the
(reduced) affine Grassmannian of \(\FRJ_a^0\) (not to be confused with the
neutral component of \(\cP_v(a)\), which we denote by \(\cP_v(a)_0\)).

We have a homomorphism
\begin{align}
    \RH^1(F_v,\FRJ_a)\simeq\RH^1(k,\cP_v^0(a))\longto\RH^1(k,\cP_v(a)),
\end{align}
so a character \(\kappa\) of \(\RH^1(k,\cP_v(a))\) induces a character of
\(\RH^1(F_v,\FRJ_a)\), still denoted by \(\kappa\).
\begin{proposition}
    \label[proposition]{prop:pt_counting_with_orbital_general}
    Let \(\tau\) be a function on
    \(\Stack*{\Arc{G^\AD}\backslash\Gr^{\le\lambda_v}}(k)\).
    If we have a character \(\kappa\colon\RH^1(F_v,\FRJ_a)\to\Qlb^\x\) that is
    induced by a character of \(\RH^1(k,\cP_v(a))\), then we have equality
    \begin{align}
        \Cnt_\tau\sX(k)_\kappa=\vol(\FRJ_a^0(\cO_v),\dd
        t_v)\OI_a^\kappa(\tau,\dd t_v),
    \end{align}
    where \(\dd t_v\) is any Haar measure on \(\FRJ_a(F_v)\). If \(\kappa\) is
    not induced by a character of \(\RH^1(k,\cP_v(a))\), then
    \begin{align}
        \OI_a^\kappa(\tau,\dd t_v)=0.
    \end{align}
\end{proposition}
\begin{proof}
    We already know by
    \Cref{prop:pt_counting_with_orbital_connected_fiber_case} that
    \begin{align}
        \Cnt_\tau\sX^0(k)_\kappa
        =\vol(\FRJ_a^0(\cO_v),\dd t_v)\OI_a^\kappa(\tau,\dd t_v),
    \end{align}
    if \(\kappa\) is a character of \(\RH^1(k,\cP_v^0(a))\). So we need to
    prove that if \(\kappa\) is induced by a character of \(\RH^1(k,\cP_v(a))\),
    we have
    \begin{align}
        \Cnt_\tau\sX^0(k)_\kappa=\Cnt_\tau\sX(k)_\kappa,
    \end{align}
    and \(\Cnt_\tau\sX^0(k)_\kappa=0\) otherwise.

    Let \(\Pi\defeq\pi_0(\FRJ_{a,v})\) be the kernel of
    \(\cP_v^0(a)\to\cP_v(a)\). It is the group of connected components of the
    special fiber of \(\FRJ_a\). For each isomorphism class in \(\sX(k)\),
    we fix a representative \(x=(m,p)\).

    We consider a full subcategory
    \(\sX^0(k)_x\) of \(\sX^0(k)\) consisting of objects lying over \(x\) for
    each representative \(x\).
    If \(x\) and \(x'\) are representatives of two isomorphism classes in
    \(\sX(k)\), then no object in \(\sX^0(k)_x\) is isomorphic to any object in
    \(\sX^0(k)_{x'}\). On the other hand, since \(\cP_v^0(a)\to\cP_v(a)\) is
    surjective, any object in \(\sX^0(k)\) is
    isomorphic to an object of \(\sX^0(k)_x\) for some \(x\). Thus, we may
    replace \(\sX^0(k)\) with disjoint union of full subcategories
    \begin{align}
        \coprod_x\sX^0(k)_x.
    \end{align}

    For a given \(x=(m,p)\), the category \(\sX^0(k)_x\) is described as follows:
    \begin{enumerate}
        \item the objects are \((m,p_0)\) such that \(p_0\) maps to \(p\) (in
            particular, \(\sX^0(k)_x\) is non-empty);
        \item the morphism from \((m,p_0)\) to \((m,p_0')\) is an element
            \(h\in\cP_v^0(\FRJ_a)\) with \(hm=m\) and
            \(\Frob(h)p_0h^{-1}=p_0'\).
    \end{enumerate}
    In particular, if we fix \(x_0=(m,p_0)\in\sX^0(k)_x\), then any \((m,p_0')\)
    may be written as \((m,p_0p')\) for some \(p'\in\Pi\), and the automorphism
    of \((m,p_0)\) is the group
    \begin{align}
        H_0=\Set*{h\in\cP_v^0(a)\given hm=m, \Frob(h)h^{-1}\in\Pi}.
    \end{align}
    Therefore, \(\sX^0(k)_x\) is equivalent the categorical quotient
    \(\Stack{\Pi/H_0}\), where \(H_0\) acts on \(\Pi\) through map
    \begin{align}
        \alpha\colon H_0 &\longto \Pi\\
        h&\longmapsto \Frob(h)h^{-1}.
    \end{align}
    The isomorphism classes are represented by the \(\coker(\alpha)\), and
    automorphism groups are isomorphic to \(\ker(\alpha)\). It also implies that
    \(H_0\) is finite.

    Let \(\kappa\) be a character of \(\RH^1(k,\cP_v^0(a))\), the latter is
    identified with the torsion part of  \(\cP_v^0(a)_\Frob\). We still denote
    the restriction to \(\Pi_\Frob\) or \(\Pi\) by \(\kappa\). The restriction
    of \(\kappa\) to \(\alpha(H_0)\in \Pi\) is trivial by definition, so it
    induces a character on \(\coker(\alpha)\), again still denoted by
    \(\kappa\). The restriction of function \(\tau\) to \(\sX^0(k)_x\) is
    constant with value \(\tau(x)\). So we have equality of summations
    \begin{align}
        \sum_{x'\in\sX^0(k)_x/{\sim}}\frac{\Pair{\cl(x')\tau(x')}{\kappa}}{\Cnt\Aut_{\sX^0}(x')}
        =\Pair{\cl(x_0)}{\kappa}\tau(x)
        \sum_{z\in\coker(\alpha)}\frac{\Pair{z}{\kappa}}{\Cnt\ker(\alpha)}.
    \end{align}
    If \(\kappa\) is not induced by a character of \(\RH^1(k,\cP_v(a))\), that is,
    non-trivial on \(\Pi\), then the right-hand side is \(0\). Summing over
    \(x\), we have in this case
    \begin{align}
        \Cnt_\tau\sX^0(k)_\kappa=0.
    \end{align}
    If \(\kappa\) is trivial on \(\Pi\), then the same summation above is equal
    to
    \begin{align}
        \Pair{\cl(x_0)}{\kappa}\tau(x)
        \frac{\Cnt\coker(\alpha)}{\Cnt\ker(\alpha)}=
        \Pair{\cl(x_0)}{\kappa}\tau(x)
        \frac{\Cnt\Pi}{\Cnt{H_0}}.
    \end{align}
    If \(h_0\in H_0\) and let \(h\) be its image in \(\cP_v(a)\), then
    \(h\in\cP_v(a)(k)\), and we have short exact sequence
    \begin{align}
        1\longto \Pi\longto H_0\longto \Aut_{\sX(k)}(x)\longto 1.
    \end{align}
    As a result, we have \(\Cnt{H_0}=\Cnt\Pi\Cnt\Aut_{\sX(k)}(x)\), and summing
    over \(x\) of the sums above we have
    \begin{align}
        \Cnt_\tau\sX(k)_\kappa=\Cnt_\tau\sX^0(k)_\kappa.
    \end{align}
    This finishes the proof.
\end{proof}

\subsection{}
Finally, using the connection between point-counting on \(\sX(k)\) and Frobenius
trace on cohomologies, we can give a cohomological interpretation of orbital
integrals. According to \Cref{prop:projective_quotient_of_GASF} and
\Cref{sub:pt_counting_locally_finite_type_variant}, we may find a
\(\Frob\)-stable torsion-free subgroup \(\Lambda\subset\cP_v^0(a)\) such that
\(\Lambda\) acts freely on \(\cM_v(a)\), and both \(\cP_v^0(a)/\Lambda\) and
\(\cM_v(a)\) are of finite types over \(k\). Moreover, \(\cM_v(a)/\Lambda\) is
proper, so its cohomologies with compact support is canonically isomorphic to
ordinary cohomologies. Since the kernel of
\(\cP_v^0(a)\to\cP_v(a)\) is of finite type, \(\Lambda\) maps isomorphically to
its image in \(\cP_v(a)\), so we may also treat \(\Lambda\) as a subgroup of
\(\cP_v(a)\).
\begin{corollary}
    Let \(\cF\) be a bounded constructible \(\Arc{G^\AD}\)-equivariant complex
    on \(\Gr^{\le\lambda_v}\) and \(\kappa\) be a character of
    \(\RH^1(F_v,\FRJ_a)\). Let \(\FRJ_a^{\flat,0}\) be the open subgroup scheme
    of the N\'eron model of \(\FRJ_a\) of fiberwise neutral components. Then for
    any \(\Lambda\) as above, we have
    \begin{align}
        \sum_{i\in\bbZ}(-1)^i\Tr\bigl(\Frob,\RH^i(\Stack*{\cM_v(a)/\Lambda},\cF)_\kappa\bigr)
        =\vol(\FRJ_a^{\flat,0},\dd t_v)\OI_a^\kappa(\Tr_\cF,\dd t_v).
    \end{align}
\end{corollary}
\begin{proof}
    Apply \Cref{prop:pt_counting_locally_finite_type_variant} and
    \Cref{prop:pt_counting_with_orbital_connected_fiber_case} to \(\FRJ_a^0\),
    we have
    \begin{align}
        \sum_{i\in\bbZ}(-1)^i\Tr\bigl(\Frob,\RH^i(\Stack*{\cM_v(a)/\Lambda},\cF)_\kappa\bigr)
        =\Cnt{\cP_v^0(a)_0(k)}\vol(\FRJ_a^0,\dd t_v)\OI_a^\kappa(\Tr_\cF,\dd
        t_v),
    \end{align}
    where \(\cP_v^0(a)_0\) is the neutral component of \(\cP_v^0(a)\). The
    injective map \(\FRJ_a^0\to\FRJ_a^{\flat,0}\) induces short exact sequence
    \begin{align}
        1\longto\FRJ_a^{\flat,0}(\breve{\cO}_v)/\FRJ_a^0(\breve{\cO}_v)\longto
        \cP_v^0(a)\longto\cP_v^{\flat,0}(a)\longto 1.
    \end{align}
    The group \(\cP_v^{\flat,0}(a)\) is discrete, so the \(\bar{k}\)-point of
    the neutral component \(\cP_v^0(a)_0\) can be identified with
    \(\FRJ_a^{\flat,0}(\breve{\cO}_v)/\FRJ_a^0(\breve{\cO}_v)\), compatible with
    \(\Frob\)-action. Since both \(\FRJ_a^{\flat,0}\) and \(\FRJ_a^0\) have
    connected fibers, we have
    \begin{align}
        \cP_v^0(a)_0(k)\simeq\FRJ_a^{\flat,0}(\cO_v)/\FRJ_a^0(\cO_v).
    \end{align}
    This implies that
    \begin{align}
        \Cnt{\cP_v^0(a)_0(k)}\vol(\FRJ_a^0(\cO_v),\dd
        t_v)=\vol(\FRJ_a^{\flat,0}(\cO_v),\dd t_v)
    \end{align}
    for any Haar measure on \(\FRJ_a(F_v)\). Combining it with the equality
    above, we have our result.
\end{proof}


\section{Counting Points on mH-fibrations} 
\label{sec:counting_points_on_mh_fibrations}

In this section, we fix \(a\in\cA_X^\ANI(k)\) a point in the anisotropic locus
(cf.~\Cref{def:anisotropic_locus})
of the mH-base and consider point-counting
problems on \(\cM_a\). Let \(\FRJ_a'\to\FRJ_a\) be a morphism of smooth
commutative groups schemes on \(X\) that is an isomorphism over
\(U=X-\FRD_a\). Suppose \(\FRJ_a'\) has connected fibers. For example, we may
let \(\FRJ_a'=\FRJ_a^0\). Let \(\cP_a'\) be the Picard stack classifying
\(\FRJ_a'\)-torsors on \(X\), then the \(\cP_a\)-action on \(\cM_a\) induces an
action of \(\cP_a'\). Unlike \cite{Ng10}, we do not need to choose such
\(\FRJ_a'\) that \(\cP_a'\) is a scheme, because in
\Cref{prop:M_mod_P_pt_counting_general} we allow \(\cP\) to be a Deligne--Mumford
stack.

\subsection{}
By \Cref{prop:product_formula_alt}, we have
\(\Gal(\bar{k}/k)\)-equivariant equivalence of groupoids
\begin{align}
    \Stack{\cM_a/\cP_a}=\prod_{v\in \abs{X-U}}\Stack{\cM_v(a)/\cP_v(a)},
\end{align}
and similarly if we replace \(\cP_a\) (resp.~\(\cP_v(a)\)) by \(\cP_a'\)
(resp.~\(\cP_v'(a)\)). In particular, we have equivalence of groupoids of
\(k\)-points
\begin{align}
    \Stack{\cM_a/\cP_a'}(k)=\prod_{v\in \abs{X-U}}\Stack{\cM_v(a)/\cP_v'(a)}(k).
\end{align}
A constructible complex or a function on the left-hand side is called
\notion{factorizable}\index{function!factorizable}\index{sheaf!factorizable} if it is an outer tensor of the same on the
right-hand side. In this case, we write \(\cF=\boxtimes_{v\in\abs{X-U}}\cF_v\)
or \(\tau=\boxtimes_{v\in\abs{X-U}}\tau_v\) respectively.

Since \(a\) is anisotropic, the group \(\pi_0(\cP_a)\) is
finite with
\(\Frob\)-action. For any \(v\in\abs{X-U}\), the map \(\cP_v'(a)\to\cP_a'\)
induces map of connected components
\begin{align}
    \pi_0(\cP_v'(a))\longto \pi_0(\cP_a').
\end{align}
Any \(\Frob\)-invariant character
\begin{align}
    \kappa\colon\pi_0(\cP_a)_\Frob\longto\Qlb^\x
\end{align}
induces characters on \(\pi_0(\cP_a')\) and \(\pi_0(\cP_v'(a))\). On the other
hand, if \(x\in\Stack{\cM_a/\cP_a}(k)\) corresponds to a tuple of points
\(x_v\in\Stack{\cM_v(a)/\cP_v'(a)}(k)\), then by \Cref{prop:product_formula} we
have
\begin{align}
    \Pair{\cl(x)}{\kappa}=\prod_{v\in\abs{X-U}}\Pair{\cl(x_v)}{\kappa}.
\end{align}
Caution that the above equality is sensitive to the choice of the base points used to define
\(\cM_v(a)\), and we \emph{may need to shrink \(U\)} (see
the construction of \(\gamma_{U_a}\) at the beginning of
\Cref{sec:product_formula}), because otherwise
\Cref{prop:product_formula} may fail to hold over \(k\).
Thus, we have the following result:
\begin{proposition}
    \label[proposition]{prop:pt_counting_on_mH_fiber_using_prod_formula_alt}
    For any \(\Frob\)-invariant character \(\kappa\) of \(\pi_0(\cP_a)\) and any
    \(\Qlb\)-valued factorizable function \(\tau\) on
    \(\Stack{\cM_a/\cP_a'}(k)\), we have
    \begin{align}
        \Cnt_\tau\Stack{\cM_a/\cP_a'}(k)_\kappa=\prod_{v\in
        \abs{X-U}}\Cnt_{\tau_v}\Stack{\cM_v(a)/\cP_v'(a)}(k)_\kappa.
    \end{align}
\end{proposition}

\begin{corollary}
    \label[corollary]{cor:trace_formula_for_M_mod_P_global}
    For any \(\Frob\)-invariant character \(\kappa\) of \(\pi_0(\cP_a)\), and
    \(\cF\) a bounded, constructible and factorizable \(\cP_a\)-equivariant
    complex on \(\cM_a\), we have
    \begin{align}
        \sum_{n\in\bbZ}(-1)^n\Tr\bigl(\Frob,\RH^n(\cM_{a,\bar{k}},\cF)_\kappa\bigr)
        =\Cnt(\cP_a')_0(k)\prod_{v\in\abs{X-U}}\Cnt_{\cF_v}\Stack{\cM_v(a)/\cP_v'(a)}(k)_\kappa.
    \end{align}
\end{corollary}
\begin{proof}
    This follows from
    \Cref{prop:M_mod_P_pt_counting_general,prop:pt_counting_on_mH_fiber_using_prod_formula_alt}
    and that \(\cM_a\) is proper by
    \Cref{prop:properness_over_anisotropic_locus}.
\end{proof}

\subsection{}
Finally, we express \(\Cnt_\tau\Stack{\cM_v(a)/\cP_v(a)}(k)_\kappa\) as orbital
integrals. For that purpose, we choose at each point \(v\in\abs{X-U}\):
\begin{enumerate}
    \item a trivialization of the \(Z_\FRM\)-torsor \(\cL\) on the formal disc
        \(X_v\);
    \item a Haar measure \(\dd t_v\) on \(\FRJ_a(F_v)\).
\end{enumerate}
In this way we have equality by
\Cref{prop:pt_counting_with_orbital_connected_fiber_case}:
\begin{align}
    \Cnt_\tau\Stack{\cM_v(a)/\cP_v'(a)}(k)_\kappa=\vol(\FRJ_a'(\cO_v),\dd
    t_v)\OI_{a,v}^\kappa(\tau_v,\dd t_v).
\end{align}
Therefore, we have
\begin{align}
    \sum_{n\in\bbZ}(-1)^n\Tr\bigl(\Frob,\RH^n(\cM_{a,\bar{k}},\cF)_\kappa\bigr)
    =\Cnt(\cP_a')_0(k)\prod_{v\in\abs{X-U}}\vol(\FRJ_a'(\cO_v),\dd
    t_v)\OI_{a,v}^\kappa(\Tr_{\cF_v},\dd t_v).
\end{align}
If \(Z_G\) is connected, then by the discussion in
\Cref{sec:matching_orbits}, \(a\) can be treated as a conjugacy class in
\(G^\AD\), hence the trivialization requirement of \(\cL\) over \(X_v\) can be ignored.



\chapter{Proof of the Fundamental Lemma} 
\label{chap:proof_of_fundamental_lemma}

In this chapter we are going to prove \Cref{thm:FL_main},
\Cref{thm:asymptotic_FL}, and
\Cref{thm:main_geometric_stabilization,thm:arithmetic_geometric_stabilization}.
They will be proved together because their proofs intertwine with each other,
and similar to \cite{Ng10}, we need to go back and forth between the global
setting and the local setting. This process is more intricate in the
multiplicative case due to extra complication caused by boundary divisors.

\section{Preliminary Reductions} 
\label{sec:preliminary_reductions}

We start by doing some routine reductions. Let \(G\) be a quasi-split reductive
group over \(\cO_v\) with a fixed pinning. Suppose \(C\) is a split central subtorus in \(G\), then
it is a standard fact that fundamental lemma for \(G\) is equivalent to that for
\(G/C\). Therefore, we may (and shall) assume that the center of \(G\) contains no
split torus. Such groups are called \emph{elliptic}.\index{group!elliptic}

\subsection{}
We first reduce to the case where \(H\) is an \emph{elliptic} endoscopic
group, following a similar argument as in \cite{Ng10}*{\S~8.6.1}.
Let \(H\) be an endoscopic group of \(G\), and \(\gamma_H\in
H(F_v)\) and \(\gamma\in G(F_v)\) are matching (strongly regular
semisimple) conjugacy classes. Suppose \(C_H\) is a
central split subtorus in \(H\), then it naturally embeds into the centralizer
\(I_\gamma\). Let \(L\) be the centralizer of \(C_H\) in \(G\), which is a Levi
subgroup, and \(\gamma\in L(F_v)\). Let \(H_L\) be the endoscopic group of \(L\)
that is the descent of \(H\) (see \cite{LS90}*{\S~1.4}), then \(\gamma_H\in H_L(F_v)\).
The group \(H_L\) is naturally identified with a Levi subgroup of \(H\).
To emphasize different groups, we use notations \(\gamma_L\in L(F_v)\) and
\(\gamma_{H_L}\in H_L(F_v)\) instead.

We have a canonical transfer map between spherical Hecke algebras
\(\cH_{G,0}\to\cH_{L,0}\), which, via geometric Satake, may be defined as
follows: let \(\cF\) be any \(\Arc_v{G}\)-equivariant perverse sheaf on
\(\Gr_{G,v}\), and \(P\) be the unique parabolic subgroup of \(G\) containing
both \(L\) and the Borel \(B\) that is part of the pinning of \(G\). 
Then we have a canonical diagram
\begin{align}
    \Gr_{L,v}\stackrel{p}{\longot}\Gr_{P,v}\stackrel{i}{\longto}\Gr_{G,v}.
\end{align}
Then \(\cF_L=p_!i^*\cF\) is a \(\Arc_v{L}\)-equivariant perverse sheaf on
\(\Gr_{L,v}\) by the general theory of geometric Satake. The function \(f\) induced by
\(\cF\) on \(G(F_v)\) transfers to the function \(f_L\) induced by
\(\cF_L\) on \(L(F_v)\). More explicitly at the function level, we have
\begin{align}
    f_L(l_v)=\delta_P^{\OneHalf}(l_v)\int_{N(F_v)}f(l_vn_v)\dd n_v,
\end{align}
where \(N\) is the unipotent radical of \(P\) and the left Haar measure \(\dd n_v\)
is such that \(N(\cO_v)\) has volume \(1\), and
\begin{align}
    \delta_P\colon L(F_v)& \longto \Qlb^\x\\
    l_v&\longmapsto \abs*{\det\bigl(\Ad_{l_v}|\Lie(N)(F_v)\bigr)}_{F_v}.
\end{align}
We remind the reader that we always fix a half Tate twist, and so also a square
root of \(q\).

\begin{proposition}[Harish-Chandra descent]
    \label[proposition]{prop:HC_descent_of_kappa_OI}\index{Harish-Chandra descent}
    Suppose \(\gamma\in L(F_v)\) is strongly \(G\)-regular semisimple, then
    we have equality
    \begin{align}
        \OI_\gamma^\kappa(f)=q^{\deg(v)\sum_\Rt\val_{F_v}(\Rt(\gamma)-1)/2}\OI_{L,\gamma}^\kappa(f_L),
    \end{align}
    where \(\Rt\) ranges from roots of \(I_\gamma\) in \(G\) that are not in
    \(L\), and \(\deg(v)=[k_v:k]\).
\end{proposition}
\begin{proof}
    This result is well-known but we cannot locate a clean reference in the form
    stated, so we include a proof.
    Since \(H\) descends to an endoscopic
    group \(H_L\) of \(L\), the \(\kappa\)-twisting on the \(G\)-side is easily
    seen to be the same as that on the \(L\)-side by unpacking the definition.
    So we reduce to proving the same equality for the ordinary orbital integral. 

    By the definition of \(\OI_\gamma\), we have
    \begin{align}
        \OI_\gamma(f)=\int_{I_\gamma(F_v)\backslash
        G(F_v)}f\bigl(\Ad_{g_v}^{-1}(\gamma)\bigr)\frac{\dd g_v}{\dd t_v}.
    \end{align}
    By Iwasawa decomposition, we have \(G(F_v)=P(F_v)G(\cO_v)\), and so the
    above integral is the same as
    \begin{align}
        \int_{I_\gamma(F_v)\backslash P(F_v)}f\bigl(\Ad_{p_v}^{-1}(\gamma)\bigr)\frac{\dd p_v}{\dd t_v},
    \end{align}
    where \(\dd p_v\) is the left Haar measure on \(P(F_v)\) such that
    \(\vol(P(\cO_v))=1\). Using the semidirect product \(P=N\rtimes L\) and the
    fact that \(I_\gamma\subset L\), we further write the integral as
    \begin{align}
        \int_{N(F_v)\times (I_\gamma(F_v)\backslash
        L(F_v))}f\bigl(\Ad_{n_vl_v}^{-1}(\gamma)\bigr)\delta_P(l_v)^{-1}\dd n_v\frac{\dd l_v}{\dd t_v}.
    \end{align}
    To prove the proposition, it suffices to show that
    \begin{align}
        q^{\deg(v)\sum_\Rt\val_{F_v}(\Rt(\gamma)-1)/2}f_L\bigl(\Ad_{l_v}^{-1}(\gamma)\bigr)=\int_{N(F_v)}f\bigl(l_v^{-1}n_v^{-1}\gamma
        n_vl_v)\bigr)\delta_P(l_v)^{-1}\dd n_v
    \end{align}
    for every \(l_v\in L(F_v)\). Change the variable \(n_v\) by
    \(\Ad_{l_v}^{-1}(n_v)\), then the measure \(\dd n_v\) changes into
    \(\delta_P(l_v)\dd n_v\), and so the integral on the right-hand side becomes
    \begin{align}
        \int_{N(F_v)}f\bigl(n_v^{-1}\Ad_{l_v}^{-1}(\gamma)n_v\bigr)\dd
        n_v.
    \end{align}
    Let \(\gamma_1=\Ad_{l_v}^{-1}(\gamma)\), then
    \begin{align}
        n_v^{-1}\Ad_{l_v}^{-1}(\gamma)n_v=\gamma_1(\gamma_1^{-1}n_v^{-1}\gamma_1 n_v).
    \end{align}
    Since \(\gamma\) is strongly \(G\)-regular semisimple, so is \(\gamma_1\),
    and it implies that
    \begin{align}
        n_v\longmapsto \gamma_1^{-1}n_v^{-1}\gamma_1 n_v
    \end{align}
    is an automorphism of \(N(F_v)\) as a topological space. Looking at the
    induced tangent map, we see
    that this change of variable \(n_v\) will change \(\dd n_v\) at \(u\in N(F_v)\) into
    \begin{align}
        \abs*{\det\bigl(1-\Ad_{u^{-1}\gamma_1u}^{-1} |\Lie(N)(F_v)\bigr)}_{F_v}^{-1}\dd n_v
        =\abs*{\det\bigl(1-\Ad_{\gamma}^{-1}|\Lie(N)(F_v)\bigr)}_{F_v}^{-1}\dd n_v.
    \end{align}
    Therefore, we reduce to show that
    \begin{align}
        q^{\deg(v)\sum_\Rt\val_{F_v}(\Rt(\gamma)-1)/2}f_L\bigl(\Ad_{l_v}^{-1}(\gamma)\bigr)
        =\abs*{\det\bigl(1-\Ad_{\gamma}^{-1} |\Lie(N)(F_v)\bigr)}_{F_v}^{-1}
        \int_{N(F_v)}f\bigl(\Ad_{l_v}^{-1}(\gamma)n_v\bigr)\dd n_v.
    \end{align}

    On the other hand, we have by the transfer map
    \begin{align}
        f_L\bigl(\Ad_{l_v}^{-1}(\gamma)\bigr)
        =\delta_P^{\OneHalf}\bigl(\Ad_{l_v}^{-1}(\gamma)\bigr)\int_{N(F_v)}f\bigl(\Ad_{l_v}^{-1}(\gamma)n_v\bigr)\dd n_v
        =\delta_P^{\OneHalf}(\gamma)\int_{N(F_v)}f\bigl(\Ad_{l_v}^{-1}(\gamma)n_v\bigr)\dd n_v.
    \end{align}
    So we finally reduce to show that
    \begin{align}
        \label{eqn:HC_descent_of_kappa_OI_final_tmp}
        \abs*{\det\bigl(1-\Ad_{\gamma}^{-1} |\Lie(N)(F_v)\bigr)}_{F_v}\delta_P^{\OneHalf}(\gamma)
        =q^{-\deg(v)\sum_\Rt\val_{F_v}(\Rt(\gamma)-1)/2}.
    \end{align}
    Let \(\Roots_N\) be the set of roots in \(N\) of the maximal torus
    \(I_{\gamma}\), then the disjoint union \(\Roots_N\cup -\Roots_N\) is
    precisely the set of roots of \(I_\gamma\) in \(G\) that are not in \(L\).
    Then we see that
    \begin{align}
        \abs*{\det\bigl(1-\Ad_{\gamma}^{-1} |\Lie(N)(F_v)\bigr)}_{F_v}\delta_P^{\OneHalf}(\gamma)
        &=q^{-\deg(v)\sum_{\Rt\in\Roots_N}\val_{F_v}[(1-\Rt(\gamma)^{-1})\Rt(\gamma)^{\OneHalf}]}\\
        &=q^{-\deg(v)\sum_{\Rt\in\Roots_N}\val_{F_v}[\Rt(\gamma)^{\OneHalf}-\Rt(\gamma)^{-\OneHalf}]}.
    \end{align}
    The square root of \(q\) is always fixed, and so we may
    compare the squares of both sides of
    \eqref{eqn:HC_descent_of_kappa_OI_final_tmp}. But after squaring, both sides
    become the same quantity
    \begin{align}
        q^{-\deg(v)\sum_{\Rt\in\Roots_N}\val_{F_v}[2-\Rt(\gamma)-\Rt(\gamma)^{-1}]}.
    \end{align}
    This finishes the proof.
\end{proof}

By the same argument of \cite{Ha93}*{Lemma~9.2}, we have
\begin{align}
    \Delta_0(\gamma_H,\gamma)=q^{-\deg(v)\sum_{\Rt}\val_{F_v}(\Rt(\gamma)-1)}\Delta_0(\gamma_{H_L},\gamma_L),
\end{align}
where \(\Rt\) ranges from roots of \(I_\gamma\) in \(G\) that are neither in
\(L\) nor in \(H\). Combined with \Cref{prop:HC_descent_of_kappa_OI} (applied to
both \(\OI_\gamma^\kappa\) and \(\SOI_{\gamma_H}\)),
we see that the fundamental lemma holds for \((\gamma_H,\gamma)\) if and only if
it holds for \((\gamma_{H_L},\gamma_L)\). Therefore, by replacing \(G\) by \(L/C_H\) and \(H\) by
\(H_L/C_H\), we may assume that the center of \(H\) contains no split torus
either. In other words, \(H\) is also elliptic.

\subsection{}
The next reduction is due to \cite{Ha95} and rather subtle (though not
particularly difficult). Although the main body of that paper is written for
local number fields (mainly due to its method of simple trace formula),
the preliminary reductions in \S~3 of \textit{loc. cit.} uses only results
(Galois cohomology, Satake transform, etc.) valid for any local field. First,
according to Lemma~3.6 of \textit{loc. cit.},
the fundamental lemma for \(G\) can be reduced to that of \(G^\AD\), with one
caveat, however. A well-known fact is that groups like \(\PGL_2\) has no
non-trivial endoscopic group while \(\SL_2\) does have one being any elliptic
subtorus. The reason is that in the definition of \emph{standard} endoscopy, we
require \(\kappa\) to be fixed by admissible embedding \(\xi\), up to a
\emph{cohomologically trivial} \(1\)-cocycle in \(Z_{\dual{\bG}}\), which by duality
corresponds to the trivial quasi-character \(\omega=1\) of \(G(F_v)\). In order to allow
endoscopies like the ones for \(\SL_2\) to happen for \(G^\AD\), one has to
allow non-trivial \(\omega\), or equivalently, \(\kappa\) is fixed by \(\xi\)
only up to an \emph{arbitrary} cocycle in \(Z_{\dual{\bG}}\).

This allows us to only consider \(G\) such that \(G^\Der\) is simple, but it
also makes endoscopic theory more complicated. Fortunately such issue is
immediately solved by Lemma~3.7 of that same paper: for any unramified simple
adjoint group \(G^\AD\), and an endoscopic datum with potentially non-trivial
\(\omega\), there exists a group \(G_1\) such that:
\begin{enumerate}
    \item The center of \(G_1\) is connected, and so the derived group of
        \(\dual{\bG}_1\) is simply-connected.
    \item The derived subgroup of \(G_1\) is simply-connected.
    \item The image of \(G_1(F_v)\) in \(G^\AD(F_v)\) is exactly the kernel of
        \(\omega\).
\end{enumerate}
If we let \(G\) be the quotient of \(G_1\) by the maximal split subtorus \(Z_1\)
in its center, we can assume that \(G\) has connected and anisotropic center, but
\(G^\Der\) is no longer simply-connected. Nevertheless, the map \(G_1(F_v)\to
G(F_v)\) is surjective. Then it can be shown with some
elementary argument (see \textit{loc. cit.}) that fundamental lemma for
\(G^\AD\) and general \(\omega\) can be reduced to that for \(G\) and trivial
\(\omega\). Thus, we return to the category of standard endoscopy, and may
assume that \(G\) has simple derived subgroup.

\subsection{}
The group \(G_1\) and \(G\) are constructed in a case-by-case manner, because in fact there are
very few cases where \(\omega\neq 1\) can happen. Take a rank-\(1\) elliptic
torus as an endoscopic group of
\(\SL_2\) for example. It can be viewed as an endoscopic group of \(\PGL_2\)
with \(\omega\) having order \(2\). Following \cite{Ha95}*{Lemma~3.7}, let \(S\)
be the torus whose character lattice is
\begin{align}
    \CharG(S)=\bbZ^3/(2,1,1)\bbZ,
\end{align}
and let \(Z_{\SL_2}=\mu_2\to S\) be the embedding induced by
\begin{align}
    \CharG(S)&\longto \bbZ/2\bbZ\\
    (a,b,c)&\longmapsto a.
\end{align}
Let the Frobenius act on \(\CharG(S)\) by involution \((a,b,c)\mapsto
(a,c,b)\), and let
group \(G_1\) be the quotient of \(S\x\SL_2\) by the anti-diagonal
embedding of \(\mu_2\). Note that \(G_1\) is no longer split even if both
\(\PGL_2\) and \(\SL_2\) are. The kernel of \(G_1\to \PGL_2\) is
just the torus \(S\), and the quotient torus by its maximal split subtorus is
the \(1\)-dimensional elliptic torus whose character lattice is generated by
\((0,1,-1)\).

The \(L\)-group of \(G\) is simply
\begin{align}
    \LD{G}=\GL_2\rtimes\Gamma_F,
\end{align}
where the Frobenius acts trivially on subgroup \(\SL_2\) and inverts the
central torus, or equivalently it acts on the diagonal torus by
\((a,b)\mapsto (b^{-1},a^{-1})\). In this case the standard representation of
\(\GL_2\) is not stable under the Frobenius, but the direct sum of the standard
representation and its dual is. This is important because it gives us a
representation of \(\LD{G}\) whose \(\dual{\bG}\)-irreducible components are all
\emph{minuscule}
(cf.~\Cref{sec:Reduction_to_minuscule_and_quasi_minuscule_coweights}). On the
other hand, if we stick to \(G=\SL_2\), then there is no minuscule
representation for its dual group \(\PGL_2\). We shall see that we want to
exploit minuscule representations whenever it is possible, because combinatorial
computations like
\Cref{thm:asymptotic_FL} become
considerably easier.

\begin{remark}
    We would like to point out that even with non-trivial \(\omega\), it is
    ultimately unlikely a serious issue at conceptual level.
    In a more general framework such as \cite{KS99}, such
    additional twisting was already incorporated. Nevertheless, we will still
    use the reduction by \cite{Ha95}, as it saves quite a bit of technical
    nuisance.
\end{remark}

\subsection{}
If \(G^\Der\) is simple of either type \(\TypeE_8\), \(\TypeF_4\), or
\(\TypeG_2\), then the quasi-split group \(G^\Der=G^\SC=G^\AD\) is necessarily
split, and \(G\) is a direct product of a torus and \(G^\Der\). As a result, for
these three types, we can reduce to the case where \(G\) is split and
adjoint. These three types are thus especially nice and will be treated slightly
differently
(cf.~\Cref{sec:Reduction_to_minuscule_and_quasi_minuscule_coweights}), and so we
make the following definition:

\begin{definition}
    We say a reductive group \(G\) is \notion{unlucky}\index{group!unlucky} if it has no
    simple factor of types \(\TypeE_8\), \(\TypeF_4\), or \(\TypeG_2\).
    We say \(G\) is \notion{very lucky}\index{group!very lucky} if it consists only of simple
    factors of those three types.
\end{definition}

\section{Some Simple Cases}
\label{sec:some_simple_cases}

In this section, we compute some simplest cases of fundamental lemma, which are
similar to the Lie algebra case computed in \cite{Ng10}*{Lemme~8.5.7}.

\subsection{}
We first recall some basic facts about the resultant divisor \(\FRR_H^G\) from
\Cref{sub:absolute_resultant_divisor}.
Recall on \(\FRC_{\FRM,H}\) we have equality of principal divisors
\begin{align}
    \nu_H^*\FRD_\FRM=\FRD_{\FRM,H}+2\FRR_H^G,
\end{align}
where both \(\FRD_{\FRM,H}\) and \(\FRR_H^G\) are principal divisors. The
discriminant divisor \(\FRD_{\FRM,H}\) is reduced, but the resultant
\(\FRR_H^G\) may not be. However, we do know that \(\FRR_H^G\) is the sum of
divisors
\begin{align}
    \FRR_H^G=(\FRR_H^G)'+(\FRR_H^G)'',
\end{align}
where \((\FRR_H^G)'\) is reduced and intersects with the numerical boundary
divisor \(\FRE_{\FRM,H}\) properly, and \((\FRR_H^G)''\) is a
sum of multiples of irreducible components of \(\FRE_{\FRM,H}\).
The non-reducedness of \((\FRR_H^G)''\) is characterized by
\Cref{prop:characterizing_non_reducedness_of_resultant}.

\subsection{}
We now return to the setup for \Cref{prop:very_weak_stabilization} (we
have not showed that \(\tilde{\cA}_H^{\kappa,\dagger}\) exists yet,
which will be done in
\Cref{sec:construction_of_a_good_mH_base}). Pick a point
\(a\in\tilde{\nu}_\cA(\tilde{\cA}_{H}^{\kappa,\dagger})(\bar{k})\) and
without loss of generality we may assume it is defined over \(k\) (otherwise we
can base change to some finite extension \(k'/k\)
and replace \(k\) by \(k'\)). Then the endoscopic datum is also defined over
\(k\). By assumption, \(\tilde{\cA}_H^{\kappa,\dagger}\)
is very \(H\)-ample, therefore similar to
\Cref{prop:non_emptyness_of_A_sharp_diamond_heart}, we may find such \(a\) that
every \(a_H\in\tilde{\nu}_\cA^{-1}(a)\) intersects with
\(\FRD_{\FRM,H}+(\FRR_H^G)^\Red\) transversally, and
\(a_H^*\bigl(\FRD_{\FRM,H}+(\FRR_H^G)'\bigr)\) does not collide with the boundary divisor.
In particular, we have \(\delta_{H,a_H}=0\). Let \(\cU\) be the open
subset of \(\tilde{\nu}_\cA(\tilde{\cA}_{H}^{\kappa,\dagger})(\bar{k})\)
consisting of such \(a\). Without loss of generality, we may also assume
\(\cU\) to be irreducible by looking at each irreducible component.

\subsection{}
For any \(a\in \cU(k)\), we
look at the local situation at any \(\bar{v}\in X(\bar{k})\).
Note that \(a\) is unramified if and only if any of \(a_H\in\tilde{\nu}_\cA^{-1}(a)\)
is, and we always have equality of ramification indices
\(c_{\bar{v}}(a)=c_{H,\bar{v}}(a_H)\). Let
\(-w_0(\lambda_{\bar{v}})\) be the boundary cocharacter of \(a\) at \(\bar{v}\) and
\(-w_{H,0}(\lambda_{H,\bar{v}})\) that of \(a_H\), and \(\nu_{\bar{v}}\) be the local
Newton point. Let \(\nu_{\bar{v}}^+\) (resp.~\(\nu_{H,\bar{v}}^+\)) be the
\((G,B)\)-dominant (resp.~\((H,B_H)\)-dominant) element in the
\(W\)-orbit (resp.~\(W_H\)-orbit) of \(\nu_{\bar{v}}\). We then have the
following possibilities:
\begin{equation}
    \label{eqn:simple_GASF_cases}
    \begin{array}{|c|ccc|cccc|cc|}
        \hline
        \# & \lambda_{\bar{v}} & \lambda_{H,\bar{v}} & \nu_{\bar{v}} &
        d_{H,\bar{v}+} & d_{\bar{v}+} & c_{\bar{v}} & r_{H,\bar{v}}^G &
        \delta_{H,\bar{v}} & \delta_{\bar{v}}\\
        \hline
        1 & 0 & 0 & 0 &     0 & 0 & 0 & 0 &     0 & 0 \\
        2 & 0 & 0 & 0 &     1 & 1 & 1 & 0 &     0 & 0 \\
        3 & 0 & 0 & 0 &     0 & 2 & 0 & 1 &     0 & 1 \\
        4 & \lambda_{\bar{v}} & -w_{H,0}(\nu_{H,\bar{v}}^+) & \nu_{\bar{v}} &
        0 & \Pair{2\rho}{\lambda_{\bar{v}}-\nu_{\bar{v}}^+} & 0 &
        \Pair{\rho}{\lambda_{\bar{v}}-\nu_{\bar{v}}^+}  &     0 &
        \Pair{\rho}{\lambda_{\bar{v}}-\nu_{\bar{v}}^+} \\
        \hline
    \end{array}
\end{equation}
The first three cases are essentially the same as the ones in
\cite{Ng10}*{Lemme~8.5.7} (given in the same order), while the last case is
unique to the multiplicative setting. Note that in the final case \(a\) is
always \(\nu\)-regular semisimple at \(\bar{v}\) due to the multiplicity formula
(\Cref{prop:characterizing_non_reducedness_of_resultant}) for \((\FRR_H^G)''\).

\subsection{}
The natural homomorphism \(\FRJ_a\to\FRJ_{H,a_H}\) induces homomorphism
\(\cP_a\to\cP_{H,a_H}\) and similarly the local analogues
\(\cP_v(a)\to\cP_{H,v}(a_H)\) for any closed point \(v\in\abs{X}\). So we have
an action of \(\cP_v(a)\) on \(\cM_{H,v}(a_H)\) too. By local model of singularity
\Cref{thm:local_singularity_model_weak,thm:local_singularity_model_main}, the
respective functions induced by sheaves \(\bcQ\) and \(\bcQ_H^\kappa\)
(assuming \(\bcQ\) satisfies the condition in \Cref{def:global_transfer_Satake} and
\(\bcQ_H^\kappa\) is a choice of global transfer of \(\bcQ\)) are
factorizable, but not in a unique way since Tate twists can be adjusted among
local factors. We use the following normalization: if \(\bcQ\) is the
intersection complex of \(\tilde{\cM}^\ddagger\), then at a general 
point in \(\Stack{\cM_v(a)/\cP_v(a)}(\bar{k})\), the stalk of the local factor
of \(\bcQ\) at \(v\) is isomorphic to
\begin{align}
    (-1)^{\Pair{2\rho}{\lambda_{v}}}\Qlb[\Pair{2\rho}{\lambda_{\bar{v}}}](\Pair{\rho}{\lambda_{\bar{v}}}),
\end{align}
and the case where \(\bcQ\) is a general Satake sheaf is straightforward. The
product of all local sign characters above is \(1\) because the boundary divisor
necessarily extends to a \(Z_{\Env(G^\SC)}\simeq T^\SC\)-torsor.
The local factor of \(\bcQ_H^\kappa\) is the transfer of the corresponding
factor of \(\bcQ\). Note that this normalization does not depend on \(a\) or
\(a_H\).
By dimension formulae \eqref{eqn:dimension_of_mH_base_at_a_point} and
\Cref{cor:dimension_of_P_a,cor:mH_fibers_homeomorphic_to_proj_varieties}, we see
that this gives a factorization of \(\bcQ\) but only up to an additional
global shift and Tate twist. However, one may check using
\eqref{eqn:difference_in_dimension_of_mH_bases} that this additional global
correction for \(\bcQ\) is exactly the same as that for
the factorization of \(\bcQ_H^\kappa\). Therefore, it is harmless to
pretend the global correction does not exist when moving between local and
global situations. This way the weighted point-counts
\(\Cnt_{\bcQ_v}\Stack{\cM_v(a)/\cP_v(a)}(k)_\kappa\) and
\(\Cnt_{\bcQ_{H,v}^\kappa}\Stack{\cM_{H,v}(a_H)/\cP_v(a)}(k)_{\hST}\) both make
sense for any \(v\in\abs{X}\).

At each \(v\), the restriction \(a_v\) of \(a\)
to \(X_v\), viewed as a point in \(\Stack{\FRC_\FRM/Z_\FRM}(\cO_v)\) has
finitely many lifts \(a_{H,v}\in\Stack{\FRC_{\FRM,H}/Z_\FRM^\kappa}\). The
disjoint union
\begin{align}
    \coprod_{a_{H,v}\mapsto a_v}\Stack{\cM_{H,v}(a_H)/\cP_v(a)}
\end{align}
is still defined over \(k\). If in addition we have \(\nu_{v}\in
W\lambda_{v}\), then \(a_{H,v}\) is necessarily unique. We can now state the following
point-counting equality arising from the cases in \eqref{eqn:simple_GASF_cases}:

\begin{lemma}
    \label[lemma]{lem:local_simple_GASF_pt_counting}
    For any \(a\in \cU(k)\) and any closed
    point \(v\in \abs{X}\) such that \(\nu_{v}\in W\lambda_{v}\), we have equality
    \begin{align}
        \Cnt_{\bcQ_{H,v}^\kappa}\Stack{\cM_{H,v}(a_H)/\cP_v(a)}(k)_{\hST}
        =
        \Delta_0(\gamma_H,\gamma)\Cnt_{\bcQ_v}\Stack{\cM_v(a)/\cP_v(a)}(k)_\kappa 
    \end{align}
    Moreover, it is a non-zero rational number.
\end{lemma}
\begin{proof}
    In the first three cases of \eqref{eqn:simple_GASF_cases}, where
    \(\lambda_{\bar{v}}=0\), the argument is the same as
    \cite{Ng10}*{Lemme~8.5.7}, so we only give a
    sketch. The first two cases works verbatim as in the Lie algebra case, and
    \(\cM_v(a)\) is geometrically a lattice with simply transitive
    \(\cP_v(a)\)-action. In
    the third case, one reduces to groups of semisimple rank \(1\), and can
    easily show that geometrically \(\cM_{\bar{v}}(a)\) is a union of infinite
    chain of \(\bbP^1\) as in the case of \cite{Ng10}*{\S~8.3}. The existence of
    a regular \(k\)-rational point is provided by the Steinberg
    quasi-section, which always exists for groups of semisimple rank \(1\).

    The fourth and last case of \eqref{eqn:simple_GASF_cases} is slightly more
    complicated: \(\cM_v(a)\) is still geometrically a lattice and a
    \(\cP_v(a)\)-torsor, but one has to analyze the Frobenius action more
    carefully. This is done in
    \Cref{cor:special_case_of_unramified_GASF_irreducible_components_Frob_transfer},
    from which we deduce that
    \begin{align}
        \Cnt_{\Qlb}\Stack{\cM_{H,v}(a_H)/\cP_v(a)}(k)_{\hST}
        =
        \frac{\Delta_0(\gamma_H,\gamma)}{\Delta_{\symup{IV}}(\gamma_H,\gamma)}
        \Cnt_{\Qlb}\Stack{\cM_v(a)/\cP_v(a)}(k)_\kappa,
    \end{align}
    where both sides are non-zero.
    When restricted to \(\Stack{\cM_{H,v}(a_H)/\cP_v(a)}(k)_{\hST}\)
    (resp.~\(\Stack{\cM_v(a)/\cP_v(a)}(k)_\kappa\)), we have
    \begin{align}
        \bcQ_{H,v}^\kappa&=(-1)^{\Pair{2\rho}{\lambda_{v}}}\Qlb[\Pair{2\rho_H}{\lambda_{H,v}}](\Pair{\rho_H}{\lambda_{H,v}})\\
        \text{resp. }\bcQ_v&=(-1)^{\Pair{2\rho_H}{\lambda_{H,v}}}\Qlb[\Pair{2\rho}{\lambda_{v}}](\Pair{\rho}{\lambda_{v}}),
    \end{align}
    and so we have
    \begin{align}
        q^{\deg(v)\Pair{\rho_H}{\lambda_{H,v}}}
        \Cnt_{\bcQ_{H,v}^\kappa}\Stack{\cM_{H,v}(a_H)/\cP_v(a)}(k)_{\hST}
        =
        q^{\deg(v)\Pair{\rho}{\lambda_{v}}}
        \frac{\Delta_0(\gamma_H,\gamma)}{\Delta_{\symup{IV}}(\gamma_H,\gamma)}
        \Cnt_{\bcQ_v}\Stack{\cM_v(a)/\cP_v(a)}(k)_\kappa,
    \end{align}
    where \(\deg(v)\) is the degree of \(k_v\) over \(k\). Since we also have in this case that
    \begin{align}
        \Delta_{\symup{IV}}(\gamma_H,\gamma)=q^{\deg(v)\bigl(\Pair{\rho}{\lambda_{v}}-\Pair{\rho_H}{\lambda_{H,v}}\bigr)},
    \end{align}
    we are done.
\end{proof}

\begin{proposition}
    \label[proposition]{prop:primal_stabilization_for_top_cohomology}
    Suppose \(G\) and \(H\) are both elliptic.
    Let \(a\in\tilde{\cA}_X^\ANI(k)\) be a point lying in the image of
    \(\tilde{\cA}_{H,X}\) such that for any \(\bar{v}\in X(\bar{k})\), the
    restriction \(a_{\bar{v}}\) of \(a\) to \(\breve{X}_{\bar{v}}\) belongs to
    the cases in \eqref{eqn:simple_GASF_cases}.
    Let \(F\) be the function field of \(X\) and let \(\gamma_{F,\AD}\in
    G^\AD(F)\) be a point lying over \(a_F\). Suppose
    \Cref{thm:asymptotic_FL} holds for
    \(a\) at all \(v\) in the boundary divisor, and
    suppose that \(\gamma_{F,\AD}\) lifts to a point \(\gamma_F\in G(F)\), then
    have isomorphism of \(\Gal(\bar{k}/k)\)-modules:
    \begin{align}
        \label{eqn:primal_stabilization_for_top_cohomology}
        (\RDF^{-\dim_{a}{\cM_X}+2d}h_*\bcQ)_{\kappa,a}
        \simeq
        \bigoplus_{a_H\mapsto a}(\RDF^{-\dim_{a_H}{\cM_H}+2d_{H}}h_{H,*}\bcQ_H^\kappa)_{\hST,a_H},
    \end{align}
    where \(d=\dim\cM_a\) and \(d_H=\dim\cM_{H,a_H}\).
\end{proposition}
\begin{proof}
    Since \(\gamma_F\in G(F)\), it induces \(\gamma_v\in G(F_v)\) for each
    \(v\in\abs{X}\). Choose a global admissible \(L\)-embedding \(\xi_F\) which
    also induces local \(L\)-embeddings \(\xi_v\). Let \(\bcQ_v\) be the local
    Satake sheaf
    corresponding to \(\bcQ\) and \(\bcQ_{H,v}^\kappa\) be the local transfer of
    \(\bcQ_v\). When \(v\) is outside of boundary divisor, this means that
    \(\bcQ_v=\Qlb\) and \(\bcQ_{H,v}^\kappa=\Qlb\). We also choose a global
    \(a\)-datum (resp.~\(\chi\)-datum), which induces local \(a\)-data
    (resp.~\(\chi\)-data). This way, by \cite{LS87}*{Theorem~6.4.A}, we
    have the product formula for the transfer factor:
    \begin{align}
        \prod_v\Delta_{0}(\gamma_{H,v},\gamma_v)=1,
    \end{align}
    where
    \(\Delta_0=\Delta_{\symup{I}}\Delta_{\symup{II}}\Delta_{\symup{III}_2}\Delta_{\symup{IV}}\)
    is the modified transfer factor we defined in
    \Cref{sec:Appendix_Definition_of_Transfer_Factors} (but note that
    \(\Delta_0\) means something different in \cite{LS87}).

    At each \(v\), let \(\Lambda_v\) be a \(\Frob_v\)-stable lattice of maximal
    rank in \(\cP_v(a)\), and by replacing with a sublattice of finite index we
    can also assume that it lies in the kernel \(\kappa\). If \(v\) supports the
    boundary divisor, we then have by looking the Frobenius action on
    irreducible components:
    \begin{align}
        \Delta_0(\gamma_{H,v},\gamma_v)\otimes
        \RHc^{-2\Pair{\rho}{\lambda_v}+2\dim{\cM_v(a)}}\bigl(\Stack*{\cM_v(a)/\Lambda_v},
        \bcQ_v\bigr)_\kappa
        \simeq
        \FRV[\gamma_v]_\kappa^{\xi_v},
    \end{align}
    where the Tate twist on the cohomology group is cancelled out by the fact
    that
    \begin{align}
        \Delta_{\symup{IV}}(\gamma_{H,v},\gamma_v)
        =q^{\deg(v)[\Pair{\rho}{\lambda_v}-\dim{\cM_v(a)}]}.
    \end{align}
    If \Cref{thm:asymptotic_FL} holds at
    \(v\), then we have
    \begin{multline}
        \Delta_0(\gamma_{H,v},\gamma_v)\otimes
        \RHc^{-2\Pair{\rho}{\lambda_v}+2\dim{\cM_v(a)}}\bigl(\Stack*{\cM_v(a)/\Lambda_v},
        \bcQ_v\bigr)_\kappa\\
        \cong
        \bigoplus_{a_{H,v}\mapsto a_v}
        \RHc^{-2\Pair{\rho_H}{\lambda_{H,v}}+2\dim{\cM_{H,v}(a_H)}}
        \bigl(\Stack*{\cM_{H,v}(a_H)/\Lambda_v},
        \bcQ_{H,v}^\kappa\bigr)_\hST.
    \end{multline}
    If \(v\) does not support the boundary divisor, we can similarly prove
    \begin{align}
        \Delta_0(\gamma_{H,v},\gamma_v)\otimes
        \RHc^{2\dim{\cM_v(a)}}\bigl(\Stack*{\cM_v(a)/\Lambda_v},
        \Qlb\bigr)_\kappa
        \cong
        \RHc^{2\dim{\cM_{H,v}(a_H)}}
        \bigl(\Stack*{\cM_{H,v}(a_H)/\Lambda_v}, \Qlb\bigr)_\hST,
    \end{align}
    and the argument is just a simple adaptation of
    that for \Cref{lem:local_simple_GASF_pt_counting}.

    Using product formula, we have
    \begin{align}
        (\RDF^{-\dim_{a}{\cM_X}+2d}h_*\bcQ)_{\kappa,a}
        &\simeq
        \bigotimes_v\RHc^{-2\Pair{\rho}{\lambda_v}+2\dim{\cM_v(a)}}\bigl(\Stack*{\cM_v(a)/\Lambda_v},
        \bcQ_v\bigr)_\kappa
    \end{align}
    up to a Tate twist, and similarly on \(H\)-side. Note that the right-hand
    side above is necessarily a finite tensor product, and it is isomorphic to
    \begin{align}
        \biggl(\prod_v\Delta_0(\gamma_{H,v},\gamma_v)\biggr)
        \otimes\bigotimes_v\RHc^{-2\Pair{\rho}{\lambda_v}+2\dim{\cM_v(a)}}\bigl(\Stack*{\cM_v(a)/\Lambda_v},
        \bcQ_v\bigr)_\kappa
    \end{align}
    because \(\prod_v\Delta_0(\gamma_{H,v},\gamma_v)=1\). This proves
    \eqref{eqn:primal_stabilization_for_top_cohomology} up to a Tate twist
    \(\epsilon\).

    Finally, note that the left-hand side of
    \eqref{eqn:primal_stabilization_for_top_cohomology} has Frobenius weight
    \((\dim_{\tilde{a}}{\tilde{\cA}_X}-d)/2\), while the right-hand side has weight
    \((\dim_{\tilde{a}_H}{\tilde{\cA}_{H,X}^\kappa}-d_H)/2\) (for any fixed \(\tilde{a}_H\)). Using
    \Cref{cor:dimension_of_P_a,cor:mH_fibers_homeomorphic_to_proj_varieties} and
    \eqref{eqn:difference_in_dimension_of_mH_bases}, we have
    \begin{align}
        \dim_{\tilde{a}}{\tilde{\cA}_X}-d=\dim_{\tilde{a}_H}{\tilde{\cA}_{H,X}^\kappa}-d_H.
    \end{align}
    This implies \(\epsilon\) is trivial and thus finishes the proof.
\end{proof}


\section{Construction of a Good mH-base}
\label{sec:construction_of_a_good_mH_base}

So far we have not constructed the subset \(\tilde{\cA}_H^{\kappa,\dagger}\),
crucial for many results such as
\Cref{thm:support_endoscopic_main,prop:very_weak_stabilization}. We will do this
now, and the resulting
subset will have the additional property that every generic point therein will
have multiplicity-free boundary divisor.

\subsection{}
Without loss of generality, we shall assume \((H,\kappa,\OGT_\kappa^\bullet)\)
is defined over \(k\), and \(\OGT_\kappa\colon X_\OGT\to X\) is a
geometrically connected finite \'etale cover over which both \(G\) and \(H\)
become split.

Let \(\FRM\in\FM(G^\SC)\) be such that \(\FRA_\FRM\) is of standard type. Let
\begin{align}
    \bM=\bM(-w_0(\theta_1),\ldots,-w_0(\theta_m))
\end{align}
be the split form of \(\FRM\). Then
we have \(\FRM_H\in\FM(H^\SC)\) whose abelianization \(\FRA_{\FRM,H}\) is also of
standard type. Let
\begin{align}
    \bM_\bH=\bM_\bH(-w_{H,0}(\theta_{H,11}),\ldots,-w_{H,0}(\theta_{H,1e_1}),\ldots,-w_{H,0}(\theta_{H,m1}),\ldots,-w_{H,0}(\theta_{H,me_m}))
\end{align}
be the split form of \(\FRM_H\).

We naturally have central group
\(Z_{\FRM}^\kappa\in\FRM_H^\x\), and it surjects onto \(Z_\FRM\) with connected
kernel.  We have seen in \Cref{sec:endoscopic_groups_inv_theory} that
\(\bZ_\bM^\kappa\) maps dominantly to each \(\bbA^1\)-factor of
\(\bA_{\bM,\bH}\), and as such for any \(N\in\bbN\), we can always find
\(Z_\FRM^\kappa\)-torsor \(\cL\) such that every line bundle direct summand in
\(\OGT_\kappa^*\FRA_{\FRM,H,\cL}\)-bundle has degree greater than \(N\). In
particular, we may find irreducible component \(\cU_H\) of \(\tilde{\cA}_{H,X}^\ANI\)
such that \eqref{eqn:dim_of_endoscopic_strata_cohom_assumption}
holds. Moreover, for any fixed \(\delta_H\in\bbN\), we may also assume that
\(\cU_H\) is very \((H,\delta_H)\)-ample and very \((H,N(\delta_H))\)-ample. It
then implies that the local model of singularity
\Cref{thm:local_singularity_model_weak} holds over \(\cU_H\).

\subsection{}
It is not enough to obtain \(\tilde{\cA}_H^{\kappa,\dagger}\), since we
do not have either \(\cU_H=\tilde{\nu}_\cA^{-1}\tilde{\nu}_\cA(\cU_H)\) or
that \(\tilde{\nu}_\cA(\cU_H)\subset\tilde{\cA}^\ddagger\). Both can be fixed by
slightly shrinking \(\cU_H\), however. First, let us look at the intersection
\(\tilde{\nu}_\cA(\cU_H)\cap\tilde{\cA}^\ddagger\). Let \(b_H\in\cB_{H,X}^\kappa(\bar{k})\)
be in the image of \(\cU_H\), and suppose the cocharacter-valued
divisor associated with \(b_H\) is
\begin{align}
    \lambda_{H,b_H}=\sum_{s=1}^d-w_{H,0}(\lambda_{H,\bar{v}_s})\cdot\bar{v}_s,
\end{align}
where \(\lambda_{H,\bar{v}_s}\) is an \(\bbN\)-combination of
\(\theta_{H,ij}\).
Let \(b\in\cB_X\) be the image of \(b_H\), and similarly let
\begin{align}
    \lambda_{b}=\sum_{s=1}^d-w_0(\lambda_{\bar{v}_s})\cdot\bar{v}_s,
\end{align}
where \(\lambda_{\bar{v}_s}\) is an \(\bbN\)-combination of \(\theta_i\).

We first shrink \(\cU_H\) by putting the following restriction on
\(\lambda_{H,b_H}\): there are at least
\begin{align}
    d_\FRM\defeq m\abs{\bW}(2g_X-1)+1
\end{align}
points among \(\bar{v}_s\) such that \(\lambda_{H,\bar{v}_s}=\theta_{i}\) for
some \(i\). This is easily possible by choosing \(Z_\FRM^\kappa\)-torsors. We
then shrink \(\cU_H\) a little further by requiring \(a\in\tilde{\nu}_\cA(\cU_H)\)
to be regular semisimple (with respect to \(G\)) at no less than \(d_\FRM\) of those
\(\bar{v}_s\) where \(\lambda_{H,\bar{v}_s}=\theta_{i}\) for some \(i\).
With this modification, the image of
\(\cU_H\) in \(\tilde{\cA}_X^\ANI\) satisfies the conditions in
\Cref{thm:local_singularity_model_main} as long as each \(\theta_i\) is good in
the sense of \Cref{def:goodness_of_boundary_divisors} (which is very easy to
achieve by choosing \(\theta_i\) carefully; see
\Cref{sec:Reduction_to_minuscule_and_quasi_minuscule_coweights} for example). As
a result, \(\tilde{\nu}_\cA(\cU_H)\) is contained in \(\tilde{\cA}^\ddagger\).

\subsection{}
Finally, we further modify \(\cU_H\) so that the condition
\(\cU_H=\tilde{\nu}_\cA^{-1}\tilde{\nu}_\cA(\cU_H)\) holds.
Similar to the modification above, we put more restrictions on
\(\lambda_{H,b_H}\), so that in addition to the one mentioned above, we also
require that for each \(1\le i\le m\) and \(1\le j \le e_i\), there are at least
\(2g_X-1\) points among \(\bar{v}_s\) such that
\(\lambda_{H,\bar{v}_s}=\theta_{H,ij}\). It is possible as long as the boundary
divisor is not too degenerate. We shrink \(\cU_H\) to the locus consisting of
\(a_H\) such that for each \(i\) and \(j\), \(a_H\) has local Newton point
exactly \(\theta_{H,ij}\) at no less than \(2g_X-1\) points where
\(\lambda_{H,\bar{v}_s}=\theta_{H,ij}\), and is \(\nu\)-regular semisimple at
those points.

Suppose \(a_H'\in\tilde{\nu}_\cA^{-1}\tilde{\nu}_\cA(\cU_H)(\bar{k})\) and let
\(a_H\in\cU_H\) be such that \(a=\tilde{\nu}_\cA(a_H')=\tilde{\nu}_\cA(a_H)\) and
let \(\bar{v}\) is a point where \(\lambda_{H,\bar{v}}=\theta_{H,ij}\) for some
\(i,j\) and the local Newton point of \(a_H\) at \(\bar{v}\) is also equal to
\(\theta_{H,ij}\). Then we must have \(\lambda_{H,\bar{v}}'=\theta_{H,ij'}\)
(\(\lambda_{H,\bar{v}}'\) is the boundary cocharacter of \(a_H'\) at
\(\bar{v}\)) for some \(\theta_{H,ij}\le_\bH\theta_{H,ij'}\). So for \(a_H'\) it remains to show that
\eqref{eqn:difference_in_dimension_of_mH_bases} holds for
\(a_H'\). For this purpose, it suffices to prove the same identity for points in
a neighborhood of \(a_H'\) such that the local Newton point at \(\bar{v}\) is
\(\theta_{H,ij'}\) (so they are more general points than \(a_H'\) in
\(\tilde{\cA}_{H,X}^\kappa\)). Recall that \(\nu_\cA^\heartsuit\) and
\(\nu_\cA^{\star\heartsuit}\) have the same image, so it suffices to assume that
\(a_H'\) is a point in \(\tilde{\cA}_{H,X}^{\kappa\star}\).
But in this case, it is straightforward to see that
\eqref{eqn:dim_of_endoscopic_strata_cohom_assumption} (the
\(\FRM_H^\star\)-version) holds at \(a_H'\).

For any fixed \(\delta_{H,0}\in \bbN\), one can always choose even more ample
\(Z_\FRM^\kappa\)-torsors and add finitely many
\(\bar{v}_s\) to the restriction list so that
\Cref{thm:local_singularity_model_weak} as well as \(\delta_H\)-regularity
also hold at \(a_H'\) as long as \(\delta_H\) does not exceed \(\delta_{H,0}\).
We then let
\(\tilde{\cA}_H^{\kappa,\dagger}=\tilde{\nu}_\cA^{-1}\tilde{\nu}_\cA(\cU_H)\),
and delete any strata with \(\delta_{H}>\delta_{H,0}\) (depending on how large we
need \(\delta_{H,0}\) to be). Note that clearly every generic point in
\(\tilde{\cA}_H^{\kappa,\dagger}\) has multiplicity-free boundary divisor (and
so the same is true for its image in \(\tilde{\cA}^\ddagger\)),
and this finishes the construction of \(\tilde{\cA}_H^{\kappa,\dagger}\).

\section{Minuscule and Quasi-minuscule Coweights} 
\label{sec:Reduction_to_minuscule_and_quasi_minuscule_coweights}

In this section, we will focus on some special elements in the spherical Hecke
algebra \(\cH_{G,0}\) for which the fundamental lemma is somewhat easier. The
reason behind this is that we do not know how to prove
\Cref{thm:asymptotic_FL} purely
locally, thus has to rely on global geometry and induction.

\subsection{}
Consider split group \(\bG\). Let \(\lambda\in\CoCharG(\bT)\) be a dominant
cocharacter, and \(\lambda_\AD\in\CoCharG(\bT^\AD)\) be the image of \(\lambda\).
Let \(\bbM\) be the set of minimal elements of
\(\CoCharG(\bT^\AD)_+-\Set{0}\). In other words, if \(\mu\in\bbM\) and
\(\mu'\in\CoCharG(\bT^\AD)_+\) is such that
\(\mu-\mu'\) is an \(\bbN\)-combination of simple coroots, then either
\(\mu'=\mu\) or \(\mu'=0\). There exists a unique element \(\mu_0\in\bbM\) such
that it is not a minimal element of \(\CoCharG(\bT^\AD)_+\), and it is known
that \(\mu_0\) is exactly the unique short dominant coroot (see
\cite{NP01}*{Lemme~7.1}). The element \(\mu_0\) is called
\notion{quasi-minuscule} and elements in \(\bbM-\Set{\mu_0}\) are called
\notion{minuscule}. More generally, we have the following definition:
\begin{definition}
    An element \(\lambda\in\CoCharG(\bT)_+\) is called
    \notion{minuscule}\index{cocharacter!minuscule} (resp.~\notion{quasi-minuscule}\index{cocharacter!quasi-minuscule}) if \(\lambda_\AD\) is.
    A highest-weight representation is minuscule\index{representation!minuscule}
    or quasi-minuscule\index{representation!quasi-minuscule} if its
    highest weight is.
\end{definition}

\begin{lemma}
    Let \(V_\mu\) be the irreducible \(\dual{\bG}\)-representation of highest
    weight \(\mu\in\CoCharG(\bT)_+\), and \(\mu'\) is a weight of
    \(\dual{\bT}\) in \(V_\mu\). Then \(\mu'\in\bW\mu\) if \(\mu\) is minuscule
    and \(\mu'\in\bW\mu\cup\CoCharG(\bZ_\bG)\) if \(\mu\) is quasi-minuscule.
\end{lemma}
\begin{proof}
    This is \cite{NP01}*{Lemmes~6.1, 7.1}.
\end{proof}

Among all the quasi-minuscule elements, there is a canonical one, namely the
short dominant coroot, and every other quasi-minuscule coweight differs from
this one by a central cocharacter. Therefore, we shall always use this canonical
choice when
we need quasi-minuscule coweights. For minuscule ones there is no such
convenience, and we will choose them depending on the situation.

\subsection{}
A key feature of minuscule representation is that every weight is in the Weyl
orbit of the highest one and thus its weight spaces are \(1\)-dimensional. This
makes statement like
\Cref{thm:asymptotic_FL} easy to
prove, and we already did it in
\Cref{cor:special_case_of_unramified_GASF_irreducible_components_Frob_transfer}.
Below are some examples of minuscule and quasi-minuscule representations
(the indexing of the Dynkin diagram is based on \cite{Bou68}):
\begin{enumerate}
    \item When \(\dual{\bG}\) is of type \(\TypeA_r\), every fundamental representation of
        \(\dual{\bG}^\SC\) is minuscule.
    \item When \(\dual{\bG}\) is of types \(\TypeC_r\) or \(\TypeD_r\), the
        first fundamental representation \(V_1\), corresponding to the standard representation
        of \(\SYP_{2r}\) or \(\SO_{2r}\) respectively, is minuscule. For type
        \(\TypeD_r\), the last two fundamental representations \(V_{r-1}\) and
        \(V_r\) are also minuscule (the two half-spin ones).
    \item When \(\dual{\bG}\) is of type \(\TypeB_r\), the first fundamental representation
        \(V_1\) is still the standard representation of \(\SO_{2r+1}\), but only
        quasi-minuscule while the last fundamental representation \(V_r\) (the spin
        representation) is minuscule. In this case  the \(0\)-weight space of
        \(V_1\) is still \(1\)-dimensional, making it almost as easy to handle as a
        minuscule one.
    \item When \(\dual{\bG}\) is of type \(\TypeE_6\), the first
        fundamental representation \(V_1\) is again minuscule, and so is the last one
        \(V_6\) by duality. The second fundamental representation \(V_2\) is the adjoint
        representation and quasi-minuscule.
    \item When \(\dual{\bG}\) is of type \(\TypeE_7\), the first
        fundamental representation \(V_1\) is the adjoint representation and
        quasi-minuscule. The representation \(V_7\) is minuscule.
\end{enumerate}

There are more minuscule or quasi-minuscule representations we did not list
above as we will not use them. For very lucky groups, in other words, the
groups \(\TypeE_8\), \(\TypeF_4\) and \(\TypeG_2\), there is no need to restrict
ourselves to minuscule or quasi-minuscule representations.

\subsection{}
We summarize the representations that we will use in the table below (note that
for \(\TypeE_7\) we include \(V_2\), which is neither minuscule nor
quasi-minuscule):
\begin{equation}
    \label{eqn:simple_reps_for_each_type}
    \begin{array}{|c|c|c|}
        \hline
        \text{Type(s)} & \text{Representation(s)} & \text{Notes} \\
        \hline
        \hline
        \TypeA_r & V_1,\ldots,V_r & \text{minuscule} \\
        \hline
        \TypeB_r (r\ge 2) & V_1 & \text{quasi-minuscule} \\
                & V_r & \text{minuscule} \\
        \hline
        \TypeC_r (r\ge 2) & V_1 & \text{minuscule} \\
        \hline
        \TypeD_r (r\ge 3) & V_1, V_{r-1}, V_r & \text{minuscule} \\
        \hline
        \TypeE_6 & V_1, V_6 & \text{minuscule} \\
                & V_2 & \text{quasi-minuscule} \\
        \hline
        \TypeE_7 & V_7 & \text{minuscule} \\
                & V_1 & \text{quasi-minuscule} \\
                & V_2 & \text{\emph{not} (quasi-)minuscule} \\
        \hline
        \TypeE_8, \TypeF_4, \TypeG_2 & V & \text{any representation}\\
        \hline
    \end{array}
\end{equation}


\section{From Global to Local} 
\label{sec:from_global_to_local}

We are now ready to prove the fundamental lemma. According to
\Cref{sec:preliminary_reductions}, it suffices to consider a group \(G_v\) over
\(X_v\) of which the derived subgroup is simple, and the center is connected and
anisotropic. The local \(L\)-group then has simply-connected derived subgroup.
In addition, we use induction on the semisimple rank \(r\) of \(G\). When \(G\)
is a torus, the fundamental lemma is trivial.

\subsection{}
Pick a smooth projective and geometrically connected curve \(X\) (whose global
function field is denoted by \(F\)) over \(k\)
together with two distinct \(k\)-points \(v\) and \(\infty\) and a
\(\pi_0(\kappa)\)-torsor \(\OGT_\kappa\) with a trivialization at \(\infty\) (in
other words, a \(\pi_0(\kappa)\)-torsor \(\OGT_\kappa^\bullet\) pointed over
\(\infty\)),
such that we have an isomorphism between the completion of \(X\) at \(v\) and
\(X_v\) above and an isomorphism \(\OGT_\kappa|_{X_v}\cong \OGT_{\kappa,v}\).
Let \(G\) and \(H\) be the corresponding twists of \(\bG\) and
\(\bH\) respectively, then \(G|_{X_v}\) is isomorphic to \(G_v\) and similarly
for \(H_v\).

Since \(G_v\) and \(H_v\) are elliptic, the centers \(G\) and \(H\) contain no
split subtorus over \(X\). Since \(\infty\) is defined over \(k\), for any
geometric point \(\bar{\infty}\) over \(\infty\), \(\Aut(\bar{\infty}/\infty)\)
has trivial image in \(\pi_0(\kappa)\), so that
\(\pi_1(\breve{X},\bar{\infty})\) has the same image as
\(\pi_1(X,\bar{\infty})\). This means that the centers of \(G\) and \(H\)
contain no split subtorus over \(\breve{X}\) either.

Let \(G_1\) be a global \(z\)-extension of \(G\) and 
\(V^1\) be a representation of \(\LD{G}_1\) whose irreducible components,
when restricted to \(\dual{\bG}^\SC\), are listed in the table
\eqref{eqn:simple_reps_for_each_type} and have at least one copy of each.
This is possible if, for example, we choose \(G_v\) to the cases in \cite{Ha95}
as discussed
in \Cref{sec:preliminary_reductions}, and \(G_1\) is the group in \textit{loc.
cit.} with connected center and simply-connected derived subgroup such that
\(\ker(G_1\to G)\) is a split torus. Let \(V\) be the restriction to \(\LD{G}\)
and let \(\FRM\in\FM(G^\SC)\) be the monoid corresponding to \(V\) and \(\FRM_H\)
the corresponding endoscopic monoid.

Choose an arbitrary global \(L\)-embedding \(\xi_1\colon\LD{H}_1\to\LD{G}_1\), which
induces local embeddings \(\xi_v\) at every closed point \(v\in \abs{X}\). With
\(\xi_1\) we also have local and global transfers of representations, Satake
sheaves, and Satake functions. Let \(\bcQ\) be the Satake sheaf induced by \(V\)
and \(\bcQ_H^\kappa\) its transfer induced by \(\xi_1\).

For any \(\gamma\in G(F)\), we may lift it to some \(\gamma_1\in G_1(F)\). If
\(\gamma_{\bar{H}}\in H/Z_G(F)\simeq\Stack{\FRM_H^{\kappa,\x}/Z_\FRM^\kappa}(F)\) matches
\(\gamma_\AD\in G^\AD(F)=\Stack{\FRM^\x/Z_\FRM}(F)\), then there exists a unique
\(\gamma_{H,1}\) matching \(\gamma_1\), and \(\gamma_H\) matching \(\gamma\), and
we necessarily have equality of global transfer factors
\begin{align}
    \Delta(\gamma_H,\gamma)=\Delta(\gamma_{H,1},\gamma_1)=1,
\end{align}
and the same is true for the \(\Delta_0\)-version.

\subsection{}
Given any point \(a\in \tilde{\cA}^{\ddagger}(k)\), the restriction of
\(a\) to \(F\) (where \(F\) is the function field of \(X\)) induces a point
\(\gamma_{F,\AD}\in G^\AD(F)\), which then induces \(\gamma_{v,\AD}\in
G^\AD(F_v)\) for any \(v\in\abs{X}\). Since the center of \(G\) is connected,
\(\gamma_{v,\AD}\) lifts to a point \(\gamma_v\in G(\breve{F}_{\bar{v}})\).
Since the boundary divisor is a cocharacter of \(G\), we may choose \(\gamma_v\)
such that \(\Frob_v(\gamma_v)\gamma_v^{-1}\) is contained in
\(Z_G(\breve{\cO}_{\bar{v}})\). By Lang's theorem, it implies that
\(\gamma_v\) can be chosen in \(G(F_v)\). This means that \(\gamma_{F,\AD}\)
induces an obstruction class
\begin{align}
    o_a\in\ker^1(F,Z_G)\defeq\ker\left[\RH^1(F,Z_G)\longto\prod_v\RH^1(F_v,Z_G)\right],
\end{align}
where \(\ker^1(F,Z_G)\) is well-known to be a finite group.

To resolve the obstruction \(o_a\),
we can always choose \(\OGT_\kappa^\bullet\) so that the center \(Z_G\)
splits over a cyclic extension of
\(F\): indeed, we have a canonical homomorphism
\begin{align}
    \Out(\bG)\longto \Aut(\bZ_\bG),
\end{align}
and we may choose \(X\) and \(\OGT_\kappa^\bullet\), so that the latter has the same
cyclic image as \(\Frob_v\).
By \cite{On63}*{Proposition~4.5.1}, \(\ker^1(F,Z_G)\) is
trivial. Therefore, \(o_a=1\) and \(\gamma_{F,\AD}\) lifts to some \(\gamma_F\in
G(F)\). Moreover, \(a_{H,F}\) always lifts to some point in \(H/Z_G(F)\) and
since \(o_a=1\)
it also lifts to some \(\gamma_{H,F}\in H(F)\).

\subsection{}
We first consider unlucky \(G\). By
\Cref{cor:special_case_of_unramified_GASF_irreducible_components_Frob_transfer,prop:special_case_of_three_little_pigs,prop:special_case_of_big_bad_wolf},
and the existence of \(\gamma_F\) for any \(a_F\),
\eqref{eqn:primal_stabilization_for_top_cohomology} in
\Cref{prop:primal_stabilization_for_top_cohomology} holds for
\(a=\tilde{\nu}_\cA(a_H)\) where
\(a_H\) is any sufficiently general point of an
inductive subset of \(\tilde{\cA}_H^{\kappa,\dagger}\). By
\Cref{prop:very_weak_stabilization}, we then have that:
\begin{proposition}
    \label[proposition]{prop:arith_geom_stab_holds_over_good_locus}
    \Cref{thm:arithmetic_geometric_stabilization} holds over
    \(\tilde{\nu}_{\cA}(\tilde{\cA}_H^{\kappa,\dagger})\).
\end{proposition}
Let \(\gamma_v\in G(F_v)\)
and \(\gamma_{H,v}\in H(F_v)\) represent matching strongly regular semisimple
stable conjugacy classes. Let
\(a_v^\bullet\in\Stack*{\FRC_\FRM^\x/Z_\FRM}(F_v)\) be the image of
\(\gamma_v\) and suppose it extends to
\(a_v\in\Stack*{\FRC_\FRM^\x/Z_\FRM}(\cO_v)\) (otherwise the fundamental lemma
involving \(\gamma_v\) and \(\fS^{V_v}\) is trivial). 

Let \(a_{H,v,1},\ldots,a_{H,v,e}\)
be all the lifts of \(a_v\) to
\(\Stack*{\FRC_{\FRM,H}/Z_\FRM^\kappa}(\cO_v)\). Their restrictions to \(F_v\)
are the same and equals the image of \(\gamma_{H,v}\). By making the genus
\(g_X\) of \(X\) sufficiently large, we may always choose a
sufficiently ample \(Z_\FRM^\kappa\)-torsor \(\cL_{1}^\kappa\) extending the
local \(Z_\FRM^\kappa\)-torsor of \(a_{H,v,1}\) such that the intersection
\begin{align}
    \tilde{\cA}_{H,X,\cL_{1}^\kappa}^\kappa\cap\tilde{\cA}_H^{\kappa,\dagger}
\end{align}
is non-empty. This also implies that we can locally modify \(\cL_1^\kappa\) at
\(v\), so we have \(\cL_i^\kappa\) (\(i=2,\ldots,e\)) extending the torsor of
\(a_{H,v,i}\). Let \(N\) be a sufficiently large integer so that for each \(i\),
as long as \(a_{H,v}'\equiv a_{H,v,i}\) modulo \(\pi_v^{N+1}\), we will have
\(\cM_{H,v}(a_{H,v}')\cong\cM_{H,v}(a_{H,v,i})\) compatible with local Picard
action. Consider maps of \(k\)-varieties
\begin{align}
    \tilde{\cA}_{H,X,\cL_i^\kappa}^{\kappa,\dagger}\longto
    \Arc_{v,N}\FRC_{\FRM,H,\cL_i^\kappa}
\end{align}
by restricting from \(X\) to \(X_v\). Increase the
\(H\)-ampleness of \(\cL_i^\kappa\) if necessary, we may assume this map is
smooth with geometrically irreducible fibers.

\begin{lemma}
    \label[lemma]{lem:non_empty_of_large_k_prime_points}
    Let \(S\) be a geometrically connected scheme of finite
    type over finite field \(k\), then there exists some \(N\) such that for all
    finite extension \(k'/k\) with degree larger than \(N\), the set \(S(k')\)
    is non-empty.
\end{lemma}
\begin{proof}
    Without loss of generality, we may assume that \(S\) is a smooth affine
    \(k\)-variety of pure dimension \(d\). Since \(S\) is geometrically
    connected, we have
    \begin{align}
        \RHc^{2d}(S_{\bar{k}},\Qlb)\simeq \Qlb(-d).
    \end{align}
    All other cohomology groups have higher weights, hence the Frobenius
    eigenvalues are dominated by that on the top cohomology. This implies that
    as long as \([k':k]\) is sufficiently large, the trace of Frobenius on
    \(\RHc^\bullet(S_{\bar{k}},\Qlb)\) cannot be \(0\). By Grothendieck--Lefschetz trace
    formula, \(S\) has a \(k'\)-point, and we are done.
\end{proof}

\begin{remark}
    It is important that \(S\) is geometrically connected, while the dimension
    of \(S\) is irrelevant (this corrects a minor error in \cite{Ng10}*{8.6.6}):
    if \(S=\Spec{k'[x]}\)
    for a finite extension \(k'/k\), then \(S\otimes_k\bar{k}\) is a disjoint
    union of several lines, and \(S\) has no \(k''\)-point for any \(k''/k\)
    that is not also an extension of \(k'\).
\end{remark}

\subsection{}
\Cref{lem:non_empty_of_large_k_prime_points} shows that for all sufficiently large extensions \(k'/k\), we
can find \(a_{H,i}\in\tilde{\cA}_H^{\kappa,\dagger}(k')\) such that:
\begin{enumerate}
    \item \(a_{H,i}|_{X_v}\equiv a_{H,v,i}\) modulo \(\pi_v^{N+1}\).
    \item At any \(v'\neq v\) that supports
        the boundary divisor, \(a_{H,i}\) is \(H\)-regular semisimple and
        intersects with \((\FRR_H^G)^\Red\) transversally.
    \item At any \(v'\neq v\) outside 
        the boundary divisor, \(a_{H,i}\) intersects with \(\FRD_{\FRM,H}+\FRR_H^G\)
        transversally.
\end{enumerate}
The images of \(a_{H,i}\) are the same point \(a\in\tilde{\cA}^\ddagger(k')\),
and the above conditions imply that
\begin{align}
    \tilde{\nu}_\cA^{-1}(a)=\Set*{a_{H,i}}.
\end{align}
Combining \Cref{prop:arith_geom_stab_holds_over_good_locus} with
\Cref{cor:trace_formula_for_M_mod_P_global}, we have
\begin{align}
    \Cnt(\cP_a^0)_0(k')\prod_{v}\Cnt_{\bcQ_v}\Stack{\cM_v(a)/\cP_v^0(a)}(k')_\kappa
    &=\Cnt(\cP_a^0)_0(k')
    \sum_{i=1}^e\prod_{v}\Cnt_{\bcQ_{H,v}^\kappa}\Stack{\cM_{H,v}(a_{H,i})/\cP_v^0(a)}(k')_\hST,
\end{align}
the right-hand side of which equals
\begin{align}
    \Cnt(\cP_a^0)_0(k')\left[\sum_{i=1}^e\Cnt_{\bcQ_{H,v}^\kappa}\Stack{\cM_{H,v}(a_{H,i})/\cP_v^0(a)}(k')_\hST\right]
    \prod_{v'\neq
    v}\Cnt_{\bcQ_{H,v'}^\kappa}\Stack{\cM_{H,v'}(a_{H,1})/\cP_{v'}^0(a)}(k')_\hST.
\end{align}
Cancelling out the equalities at all points other than \(v\)
(cf.~\Cref{lem:local_simple_GASF_pt_counting}), we have
\begin{align}
    \Delta_0(\gamma_v,\gamma_{H,v})\Cnt_{\bcQ_v}\Stack{\cM_v(a)/\cP_v^0(a)}(k')_\kappa
    =
    \sum_{i=1}^e\Cnt_{\bcQ_{H,v}^\kappa}\Stack{\cM_{H,v}(a_{H,i})/\cP_v^0(a)}(k')_\hST.
\end{align}
Since \(k'\) is arbitrary as long as its degree is sufficiently large, the same
equality holds over \(k\) as well, which then translates to
\begin{align}
    \Delta_0(\gamma_v,\gamma_{H,v})\OI_{\gamma_v}^\kappa(\fS^{V_v},\dd
    t_v)=\SOI_{\gamma_{H,v}}(\fS_{H,\xi_v}^{V_v},\dd t_v).
\end{align}
This finishes the proof of \Cref{thm:FL_main} for functions \(\fS^{V_v}\). It
remains to examine what functions \(\fS^{V_v}\) are, and deduce
\Cref{thm:FL_main} for the entire \(\cH_{G_v,0}\).

\subsection{}
The function \(\fS^{V_v}\) can be described as follows: suppose
\(\lambda_1,\ldots,\lambda_m\) are the highest weights in the support of
\(V\), and suppose the boundary divisor at any geometric point \(\bar{v}\) over
\(v\) is
\begin{align}
    \lambda_{\bar{v}}=\sum_{i=1}^m c_i\lambda_i,
\end{align}
then \(\fS^{V_v}\) is the function corresponding to representation
\begin{align}
    \bigotimes_{i=1}^m\Sym^{c_i}\bigl(V_{\lambda_i}\otimes\Hom_{\dual{\bG}}(V_{\lambda_i},V)\bigr).
\end{align}
We may also replace \(V\) by \(V^{\oplus n}\) for arbitrary \(n\) (and \(\FRM\)
by the monoid whose abelianization is \(\FRA_\FRM^{n}\)). This way,
\(\fS^{V_v}\) can be any function corresponding to
\(\LD{G_v}\)-representations consisting solely of minuscule representations or
quasi-minuscule representations, as well as any tensor product or symmetric
tensor of them. Moreover, if we prove \Cref{thm:FL_main} for any two of
\(\fS^V\), \(\fS^{V'}\), and \(\fS^{V\oplus V'}\), we can also deduce the
remaining one. We claim that we may then deduce \Cref{thm:FL_main} for
\(\fS^{\wedge^n V_v}\): indeed, for any
\(\Qlb\)-vector space \(V\) and any \(A\in\GL(V)\), we have
generating functions in \(t\):
\begin{align}
    \det(1-tA)&=\sum_{n=0}^\infty (-1)^n\Tr(\wedge^n A)t^n,\\
    \frac{1}{\det(1-tA)}&=\sum_{n=0}^\infty \Tr(\Sym^n A)t^n.
\end{align}
Inverting the second expression and taking the Taylor expansion at \(t=0\), we
see that as a virtual representation, \(\wedge^n V\) may be written as a
\(\bbQ\)-combination of symmetric powers of \(V\), hence proving the claim.

\subsection{}
Examine type by type, we see that exterior products of representations in
\eqref{eqn:simple_reps_for_each_type} generate all fundamental
representations, hence the entire representation ring of \(\dual{\bG}\) (which
has simply-connected derived subgroup because \(Z_G\) is connected).
By assumption, the Galois group \(\Gamma_v\) acts on the highest weights of \(\dual{\bG}\)
through a cyclic quotient group, and so any virtual \(\LD{G}\)-representation
\(V\) can be obtained as a \(\bbQ\)-linear combination of representations of the form
\(\Sym^nV_1\otimes\cdots\otimes \Sym^nV_m\) where \(n\) is some integer and
\(\Set{V_i}\) is a \(\Gamma_v\)-orbit of irreducible representations appearing
in \eqref{eqn:simple_reps_for_each_type}. By our discussion in the previous
subsection, we are able to deduce the fundamental lemma for \(f^V\).
This finishes the proof of \Cref{thm:FL_main} for unlucky groups.

\subsection{}
Reverse the local argument and global argument, we can also prove
\Cref{thm:arithmetic_geometric_stabilization} hence also
\Cref{thm:main_geometric_stabilization} for \(G\):
indeed, \Cref{thm:FL_main} implies that for any finite extension \(k'/k\) and
any \(a\in
\tilde{\nu}_\cA(\tilde{\cA}_{H,X}^\kappa)\cap\tilde{\cA}_\kappa^\ANI(k')\),
we have
\begin{align}
    \Cnt(\cP_a^0)_0(k')\prod_{v}\Cnt_{\bcQ_v}\Stack{\cM_v(a)/\cP_v^0(a)}(k')_\kappa
    &=\Cnt(\cP_a^0)_0(k')
    \sum_{a_H\mapsto
    a}\prod_{v}\Cnt_{\bcQ_{H,v}^\kappa}\Stack{\cM_{H,v}(a_{H})/\cP_v^0(a)}(k')_\hST,
\end{align}
which by \Cref{cor:trace_formula_for_M_mod_P_global} implies
\begin{align}
    \sum_{n\in\bbZ}(-1)^n\Tr\bigl(\Frob_{k'},\RH^n(\cM_{a,\bar{k}},\bcQ)_\kappa\bigr)
    =
    \sum_{n\in\bbZ}(-1)^n\Tr\bigl(\Frob_{k'},\bigoplus_{a_H\mapsto
    a}\RH^n(\cM_{H,a_H,\bar{k}},\bcQ_H^\kappa)_\hST\bigr).
\end{align}
By pure perversity of \((\tilde{h}_{X,*}^\ANI\bcQ)_\kappa\) and
\(\tilde{h}_{H,X,*}^{\kappa,\ANI}\bcQ_{H}^\kappa\), this proves
\Cref{thm:arithmetic_geometric_stabilization} hence also
\Cref{thm:main_geometric_stabilization} for unlucky groups.

\subsection{}
The isomorphism in \Cref{thm:main_geometric_stabilization} induces an
isomorphism of the corresponding top-degree ordinary cohomology groups as
Frobenius modules. Since top-degree ordinary cohomology groups are generated by
irreducible components of the respective mH-fibers, we may use the product formula
again and 
prove \Cref{thm:asymptotic_FL} in full for unlucky groups.

\subsection{}
For very lucky groups, similar to above, we have that
\Cref{thm:asymptotic_FL} holds for
\(\nu\neq 0\) (and it already holds for \(\nu=0\) by
\Cref{cor:special_case_of_unramified_GASF_irreducible_components_Frob_transfer})
since \Cref{thm:FL_main} is proved for groups with smaller semisimple rank. Then
the conclusion of \Cref{prop:primal_stabilization_for_top_cohomology} holds in
general, which implies that \Cref{prop:very_weak_stabilization} is again true
over the image of \(\tilde{\cA}^{\kappa,\dagger}\), this time for any monoid
\(\FRM\). The rest of the proof proceeds the same way.

\subsection{}
Thus, we finished the proofs of
\Cref{thm:FL_main,thm:asymptotic_FL}
for all \(G_v\), as well as
\Cref{thm:arithmetic_geometric_stabilization,thm:main_geometric_stabilization}
for any \(G\) whose center \(Z_G\) is connected with trivial
\(\ker^1(F,Z_G)\).

\subsection{}
Finally, for an arbitrary \(G\), the only obstruction to
proving \Cref{thm:arithmetic_geometric_stabilization} is whether \(a_F\) can be
lifted to \(\gamma_F\in G(F)\), which is given by the local
obstruction groups \(\RH^1(\cO_v,\pi_0(Z_G))\) and the global obstruction group
\(\ker^1(F,Z_G)\). However, for any \(a\) in a fixed connected component of
\(\tilde{\cA}^\ANI\), both can be trivialized after base change to a
sufficiently large extension \(k'/k\). This means that the isomorphism in
\Cref{thm:arithmetic_geometric_stabilization} still holds for \(G\) after
base changing to \(k'\). This proves
\Cref{thm:main_geometric_stabilization} for arbitrary \(G\).



\appendix

\chapter{Review on Transfer Factors} 
\label{chap:Review_on_Transfer_Factors}

In this chapter, we review the definition of transfer factors in standard
endoscopy theory. The main reference is the original paper \cite{LS87} by
Langlands and Shelstad. The emphasis will be on unramified groups over function
fields since they are the most relevant to us.
We reassure the readers that although \cite{LS87} is formulated
for local number fields only, all results cited in this chapter apply
to local function fields as well, because the arguments involved are either
combinatorial or Galois theoretic in nature and no harmonic analysis is
involved. Sometimes the invariant theory of the adjoint action of reductive
group \(G\) on itself requires the characteristic to be not too small, but such
a condition is consistent with our general
setup in \Cref{chap:reductive_monoids_and_invariant_theory}. In fact, the
same considerations also apply to all discussions about
the transfer factor throughout this work.

\section{Review on \texorpdfstring{\(L\)}{L}-groups}

\subsection{}
Let \(F\) be a local or global field, and \(G\) a connected reductive group over
\(F\). Let \(F^\sep\) be a fixed separable closure of \(F\), and
\(\Gamma_F=\Gal(F^\sep/F)\),
    \nomenclature[\(Gamma_F \)]{\(\Gamma_F\)}{the Galois group \(\Gal(F^\sep/F)\)}
and suppose \(G\) is split over \(F^\sep\).
If \(F\) is local, let \(\abs{}_F\) be a fixed absolute value of \(F\). For
example, we may choose the normalized absolute value:\index{absolute value!normalized} if \(F=\bbR\),
then \(\abs{x}_\bbR^2=x^2\); if \(F=\bbC\), then
\(\abs{z}_\bbC=z\bar{z}\);
and in non-Archimedean case, \(\abs{\pi}_F=q^{-1}\) for any uniformizer \(\pi\)
of \(F\) and \(q\) is the order of the residue field of \(F\).

For any two Borel pairs \((T_1,B_1)\) and \((T_2,B_2)\) of \(G\) over
\(F^\sep\), any inner automorphism of \(G\) carrying \((T_1,B_1)\) to
\((T_2,B_2)\) induces the same isomorphism \(T_1\to T_2\). So we have canonical
isomorphisms of the associated based root datum
\begin{align}
  (\CharG(T_1),\SimRts(T_1,B_1),\CoCharG(T_1),\SimCoRts(T_1,B_1))\longto(\CharG(T_2),\SimRts(T_2,B_2),\CoCharG(T_2),\SimCoRts(T_2,B_2)).
\end{align}

Thus, we obtain a diagram of based root data indexed by Borel pairs of \(G\), on which
\(\Gamma_F\) acts as automorphisms. One can thus
form the limit of this diagram, represented by an abstract based root
datum \((\CharG,\SimRts,\CoCharG,\SimCoRts)\), on which \(\Gamma_F\)-acts.
\begin{definition}
    The tuple \((\CharG,\SimRts,\CoCharG,\SimCoRts)\) is the \notion{canonical
    based root datum}\index{root datum!canonical based} of \(G\), denoted by \(\Psi_0(G)\).
\end{definition}
The Weyl group attached to \(\Psi_0(G)\) is
denoted by\footnote{In this chapter we have to use \(\Omega\) instead of \(W\)
because the latter will be reserved for Weil groups.} \(\Omega_0(G)\),
    \nomenclature[\(Omega \)]{\(\Omega,\bOmega\)}{the Weyl group
    (when distinction is needed from the Weil group)}
and it
has a canonical set of simple reflections determined by \(\SimRts\).
Note that \(\Gamma_F\) acts on \(\Omega_0(G)\), and the semidirect product
\(\Omega_0(G)\rtimes\Gamma_F\) acts on \(\CharG\) and \(\CoCharG\) in a
natural way.
More generally, if \(\Gamma\) is any group acting on \(G\) over \(F^\sep\), then
\(\Gamma\) acts on \(\Psi_0(G)\) and \(\Omega_0(G)\)as well.

\subsection{}
Take the dual based root datum of \(\Psi_0(G)\), and because \(\bbC\) is
algebraically closed, we can obtain a reductive
group \(\dual{\bG}\) over \(\bbC\) with a pinning \(\dbSPL=(\dual{\bT},
\dual{\bB}, \dual{\bFx}_+=\Set{\dual{\bU}_{\dual{\alpha}}})\), on which
\(\Gamma_F\) acts. Fixing an isomorphism \(\bbC\cong\Qlb\) as abstract fields,
we may view \(\dual{\bG}\) as defined over
\(\Qlb\). We will use \(\bbC\) in this chapter to be more consistent with number
theory literature, but in other parts of this book we use \(\Qlb\) for
geometric purposes. The Weyl group \(\bOmega\) of \(\dual{\bT}\) in
\(\dual{\bG}\) may be identified with \(\Omega_0(G)\), compatible with
\(\Gamma_F\)-action.

\subsection{}
If \(G\) is quasi-split and defined by a \(\Out(\bG)\)-torsor \(\OGT_G^\bullet\)
over \(F\), pointed over \(F^\sep\), then \(\Psi_0\) is canonically identified
with the based root datum contained in the \(F\)-pinning
\(\bSPL=(\bT,\bB,\bFx_+=\Set{\bU_\Rt})\) associated with \(\OGT_G^\bullet\), and
\(\Omega_0\) with \(\bW\). This is the setup we have in
\Cref{sec:Quasi-split Forms}. For this chapter only we will use
\(\bOmega\) in place of \(\bW\).
Note, however, that \(\bSPL\) contains additional
information being the root vectors \(\bFx_+\).

\subsection{}
\label{sub:action_and_induced_L_action}
If a group \(\Gamma\) acts on a reductive group \(H\) over an
algebraically closed field, say \(\bbC\), the resulting \(\Gamma\)-action on
\(\Psi_0(H)\) induces a \(\Gamma\)-action on some pinning of \(H\) since we can
reconstruct \(H\) together with a canonical pinning from \(\Psi_0(H)\) over
\(\bbC\). Hence, we obtain
another action of \(\Gamma\) on \(H\).

These two actions coincide after
projecting to the outer automorphism group \(\Out(H)\) of \(H\), which may be
identified with the \emph{subgroup} of \(\Aut(H)\) fixing a given pinning (see
\cite{Sp79}*{Corollary~2.14}).
However, since the original \(\Gamma\)-action may not fix a pinning of \(H\) at
all, the two actions may not
coincide in \(\Aut(H)\). Therefore, we have the following definition:
\begin{definition}
    Let \(H\) be a reductive group over an algebraically closed field, and
    \(\Gamma\) a group acting on \(H\). Such action is called an
    \notion{\(L\)-action}\index{\(L\)-!action} if it stabilizes a pinning. If \(\Gamma\) is a
    topological group, we require the \(L\)-action to also be continuous with
    respect to the discrete topology on \(\Aut(H)\).
\end{definition}

\subsection{}
The Weil group \(W_F\)
    \nomenclature[\(W_F \)]{\(W_F\)}{the Weil group of a local or global field \(F\)}
acts on \(\dbSPL\)
via natural map \(W_F\to \Gamma_F\), and we can form the
\notion{\(L\)-group}\index{\(L\)-!group}
\begin{align}
    \LD{G}=\dual{\bG}\rtimes W_F.
\end{align}
The action of \(W_F\) on \(\dual{\bG}\) is an \(L\)-action by definition. Note
that when \(G\) is an inner twist of some quasi-split form \(G^{\symup{qs}}\), then
we have canonical isomorphism \(\LD{G}\simeq\LD{G^{\symup{qs}}}\). In other words,
\(\LD{G}\) is insensitive to whether \(G\) is quasi-split or not.

\subsection{}
Suppose we have a split extension
\begin{align}
    1\longto \dual{\bG}\longto \sG \longto W_F\longto 1,
\end{align}
then \(\sG\) is not necessarily an \(L\)-group since \(W_F\)-action induced by a
splitting of this extension may not be an \(L\)-action. Nonetheless, for a fixed
\(\Gamma_F\)-pinning of \(\dual{\bG}\), we may attach to this extension an
\(L\)-action. Let \(c\colon W_F\to \sG\) be any splitting of this extension,
then \(W_F\) acts on \(\dual{\Psi}_0\), the dual of \(\Psi_0\). The latter does
not depend on \(c\) because a different choice of \(c\) changes the action on
\(\dual{\bG}\) by an inner automorphism, hence has no effect on the action on
\(\dual{\Psi}_0\). This induces another action of \(W_F\) on any pinning (e.g.,
\(\dbSPL\)), as we have discussed in \Cref{sub:action_and_induced_L_action}.
Therefore, we obtain an \(L\)-action of \(W_F\) on \(\dual{\bG}\). We will call
this the \(L\)-action associated with \((\sG,\dbSPL)\).

\section{Cohomological Formulations}
\label{sec:cohomological_notations}

\subsection{}
Let \(\Lambda\) be a \(\bbZ\)-lattice, and \(R\subset \Lambda\) a finite subset such that \(-R=R\). 

\begin{definition}
    A \inotion{gauge} \(p\) of \(R\) is a map \(R\to\Set{\pm 1}\) such that
    \(p(-\alpha)=-p(\alpha)\) for all \(\alpha\in R\).
\end{definition}

For example, if \(R\subset \Lambda=\CharG(T)\) is the set of roots of a maximal
torus \(T\subset G\), and \(B\) a Borel containing \(T\), then \(B\) determines
a gauge \(p_B\) on \(R\) such that \(p_B(\alpha)=1\) if and only if \(\alpha\)
is a root of \(T\) in \(B\). As another example, if \(O\subset R\) is a subset
such that \(R=O\coprod -O\) is a disjoint union, then one can define gauge
\(p_O(\alpha)=1\) if and only if \(\alpha\in O\).

\subsection{}
Let \(\Sigma=\Gamma\x\Ggen{\epsilon}\) be a group acting on \(\Lambda\) and \(R\) such that
\(\epsilon\) acts as \(-1\). If \(\Sigma\) is a topological group, we require
the action to be continuous with respect to the discrete topology on
\(\Lambda\). In our applications we will have \(\epsilon^2=1\),
so we will assume this as well, even though it is not required in many of the
results below.
Then we define a product notation for any \(r\)-tuple
\(a=(a_1,\ldots, a_r)\in \Sigma^r\)
\begin{align}
    \prod_{\alpha:a}^p=\prod_{\alpha:a_1,\ldots,a_r}^p
    \defeq\prod_{\substack{\alpha\in R\\ p((a_1\cdots a_s)^{-1}\alpha)=(-1)^{s+1}\\1\le s\le r}}
\end{align}
Let \(F\) be \emph{any} field and let \(\Sigma\) act on \(F\) trivially, then it acts on
\(F^\x\otimes_\bbZ \Lambda\), and we denote \(c\otimes\lambda\) by \(c^\lambda\).

\begin{lemma}[\cite{LS87}*{Lemmas~2.1.B and 2.1.C}]
    \label[lemma]{lem:tp_is_2_cocycle}
    The \(2\)-cochain
    \begin{align}
        t_p(\sigma,\tau)=\prod_{\alpha:1,\sigma,\tau}^p(-1)^\alpha
    \end{align}
    is a \(2\)-cocycle \(\Sigma^2\to F^\x\otimes_\bbZ \Lambda\).
    Moreover, if \(q\) is another gauge, \(t_p/t_q\) is a coboundary.
\end{lemma}

\subsection{}
\label{sub:definition_of_s_pq}
If \(p\) is a gauge, then so is \(-p\), and we define for a pair of
gauge \((p,q)\) and an \(r\)-tuple \(a\in\Sigma^r\) another product
\begin{align}
    \prod_{\alpha:a}^{p,q}=\prod_{\alpha:a_1,\ldots,a_r}^{p,q}\defeq\prod_{\substack{\alpha\in R\\ p((a_1\cdots a_s)^{-1}\alpha)=(-1)^{s+1}\\q((a_1\cdots a_s)^{-1}\alpha)=1\\1\le s\le r}}
\end{align}
If we define \(1\)-cochain of \(\Gamma\) (not \(\Sigma\))
\begin{align}
    s_{p/q}(\sigma)=\prod_{\alpha:1,\sigma}^{p,q}(-1)^\alpha\prod_{\alpha:1,\sigma}^{-q,p}(-1)^\alpha,
\end{align}
then one can show that (see \cite{LS87}*{Lemma~2.4.A})
\begin{align}
    \partial s_{p/q}=t_p/t_q,
\end{align}
as cochains of \(\Gamma\). In other words, we have an explicit splitting of
\(t_p/t_q\) after restricting to \(\Gamma\).

\subsection{}
Now we review the notions of \(a\)-data and \(\chi\)-data introduced in
\cite{LS87}. The former makes sense for any field \(F\), while the latter only
makes sense for local or global fields because it uses Weil groups. We start
with \(a\)-data. Suppose we have an extension of the (trivial) action of
\(\Sigma\) on \(F\) to an action on \(F^\sep\) such that \(\epsilon\) still acts
trivially.

\begin{definition}
    An \inotion{\(a\)-datum} is a \(\Gamma\)-equivariant map
    \begin{align}
        a\colon R&\longto (F^{\sep})^\x\\
        \alpha&\longmapsto a_\alpha
    \end{align}
    such that \(a_{-\alpha}=-a_{\alpha}\) (i.e. \(a\) is \(\epsilon\)-antivariant).
    A \inotion{\(b\)-datum} is a \(\Gamma\)-equivariant map \(b\colon R\to
    F^{\sep\x}\)
    that is also \(\epsilon\)-equivariant (hence \(\Sigma\)-equivariant).
\end{definition}

\begin{remark}
    In \cite{LS87}, \(a\) is called \emph{a set of \(a\)-data}. We prefer to
    view it as a single object (\emph{an \(a\)-datum}), so that
    the plural form can refer to multiple choices of such \(a\). Similar change
    is made for \(\chi\)-data as well.
\end{remark}

\subsection{}
\label{sub:sum_a_data}
Suppose \(a\)-datum exists for the \(\Sigma\)-action on \(R\), and \(p\) is a
gauge of \(R\), then we form \(1\)-cochain of \(\Gamma\)
\begin{align}
    u_p(\sigma)\defeq \prod_{\alpha:1,\sigma}^pa_\alpha^\alpha\in F^{\sep\x}\otimes_\bbZ \Lambda.
\end{align}

\begin{lemma}[\cite{LS87}*{Lemma~2.2.A}]
    Viewing \(t_p\) as a \(2\)-cocycle of \(\Gamma\) with value in \(F^{\sep\x}\otimes_\bbZ \Lambda\supset
    F^\x\otimes_\bbZ \Lambda\), we have that
    \begin{align}
        \partial u_p=t_p.
    \end{align}
\end{lemma}

Similarly, for any \(b\)-datum, we can form \(1\)-cochain that is in fact a
cocycle of \(\Gamma\) (\cite{LS87}*{Lemma~2.2.B}):
\begin{align}
    v_p(\sigma)\defeq\prod_{\alpha:1,\sigma}^pb_\alpha^\alpha\in F^{\sep\x}\otimes_\bbZ \Lambda.
\end{align}

\subsection{}
Consider the special case where \(\Gamma=\Gamma_F\). For \(\alpha\in R\), let
\(F_\alpha\subset F^\sep\) be its splitting field, in other words, \(F_\alpha\)
is the subfield fixed by the subgroup \(\Gamma_\alpha\subset \Gamma_F\) fixing
\(\alpha\). We have \(\Gamma_\alpha=\Gal(F^\sep/F_\alpha)\).
Similarly, let \(F_{\pm\alpha}\) be the splitting field of \(\pm\alpha\) and
\(\Gamma_{\pm\alpha}=\Gal(F^\sep/F_{\pm\alpha})\).
Then \([F_\alpha:F_{\pm\alpha}]\) is \(1\) if \(\Gamma_F\)-orbit of \(\alpha\)
does not contain \(-\alpha\) and \(2\) otherwise.

Suppose we have \(a\)-data for \(\Gamma_F\)-action on \(R\). Then we always have
\(a_\alpha\in F_\alpha^\x\). If moreover \([F_\alpha:F_{\pm\alpha}]=2\), then
\(\sigma(a_\alpha)=-a_\alpha\) for the unique non-trivial element \(\sigma\in
\Gamma_{\pm\alpha}/\Gamma_{\alpha}\). This means \(a_\alpha^2\equiv-1\)
in \(F_\alpha^\x/\Nm_{F_\alpha/F_{\pm\alpha}}(F_\alpha^\x)\).

\subsection{}
Next we turn to \(\chi\)-data. Unlike
\(a\)-data, since arithmetic duality will be used, \(\chi\)-data can only be
formulated for \(\Gamma=\Gamma_F\) where \(F\) a local or global field, as we
shall assume so. Suppose the action of \(\Gamma\) on \(\Lambda\) is continuous
with respect to the profinite topology on \(\Gamma\) and discrete topology on
\(\Lambda\). We use \(C_F\)  to denote either the multiplicative group if \(F\)
is local or the idele \emph{class} group if \(F\) is global.
Let \(C^*\) denote the Pontryagin dual of a locally compact abelian
group \(C\).

Then for \(\alpha\in R\), similar to the \(a\)-data case above,
we have splitting fields \(F_\alpha=F_{-\alpha}\) and \(F_{\pm\alpha}\), 
Galois groups \(\Gamma_\alpha\), \(\Gamma_{\pm\alpha}\), as well as groups
\(C_{\alpha}\), \(C_{\pm\alpha}\), etc.
Since \(F\) is local or global, we also have Weil groups
\(W_F\), \(W_{\alpha}=W_{F_\alpha}\), and \(W_{\pm\alpha}=W_{F_{\pm\alpha}}\).
Then \(\Gamma\) acts on \(\coprod_{\alpha\in O}C_\alpha^*\) for any
\(\Gamma\)-stable subset \(O\subset R\): if
\(\chi_\alpha\in C_\alpha^*\), then
\(\sigma(\chi_\alpha)\defeq\chi_\alpha\circ\sigma^{-1}\in C_{\sigma(\alpha)}^*\)
for any \(\sigma\in\Gamma\).
\begin{definition}
    A \notion{\(\chi\)-datum}\index{\(\chi\)-datum} is defined to be a \(\Gamma\)-equivariant map
    \begin{align}
        \chi\colon R&\longto \coprod_{\alpha\in R}C_\alpha^*,
    \end{align}
    such that \(\chi_{-\alpha}=\chi_{\alpha}^{-1}\), and if
    \([F_\Rt:F_{\pm\Rt}]=2\) then \(\chi_\alpha\) is
    non-trivial on \(C_{\pm\alpha}\). A
    \notion{\(\zeta\)-datum}\index{\(\zeta\)-datum} is a map of the same
    definition as a \(\chi\)-datum except that \(\zeta_\alpha\) is trivial on
    \(C_{\pm\alpha}\) when \([F_\Rt:F_{\pm\Rt}]=2\).
\end{definition}

\subsection{}
Although \(\chi\)-data do not have to exist in general, they exist in
situations relevant to transfer factors.
Note that if \([F_\alpha:F_{\pm\alpha}]=2\), and let \(\sigma\in \Gamma_{\pm\alpha}\) be a
non-trivial representative of \(\Gamma_{\pm\alpha}/\Gamma_\alpha\), then
\(\sigma(\alpha)=\sigma^{-1}(\alpha)=-\alpha\), and
\(\chi_{\alpha}^{-1}=\chi_{-\alpha}=\chi_{\sigma^{-1}(\alpha)}=\chi_\alpha\circ\sigma\).
Thus, \(\chi_{\alpha}\) must be trivial on
\(\Nm_{C_\alpha/C_{\pm\alpha}}(C_\alpha)\), hence must be an extension of the
quadratic character of \(C_{\pm\alpha}\) associated with \(F_\alpha/F_{\pm\alpha}\).
In addition, we may regard \(\chi_\alpha\) as a character of \(W_{\alpha}\) via
Artin reciprocity.

\subsection{}
Recall that the \(2\)-cocycle \(t_p\) in \Cref{lem:tp_is_2_cocycle} may be
defined for \emph{any} field \(F\). Let the field be \(\bbC\), then
we obtain a \(2\)-cocycle of \(\Gamma\) in \(\bbC^\x\otimes_\bbZ \Lambda\).
This cocycle is in general cohomologically non-trivial, but becomes
cohomologically trivial if inflated to \(W_F\).
Note that the only difference between this cocycle and the one with \(F\)
being a local or global field is that \(\Gamma\) always acts trivially on
\(\bbC\), but non-trivially on \(F^\sep\). Therefore, \(a\)-data do not help
in splitting the cocycle, and we need \(\chi\)-data instead.

\subsection{}
Suppose \(O\) is a \(\Sigma\)-orbit in \(R\), and \(\alpha\in O\) a fixed element.
Since \([F_\alpha:F]\) is finite, we can choose a finite set of representatives
of \(W_{\pm\alpha}\backslash W_F\), denoted by \(w_1,\ldots,w_m\), whose images
\(\sigma_1,\ldots, \sigma_m\) is a set of representatives of
\(\Gamma_{\pm\alpha}\backslash \Gamma\). Then
\(O=\Set{\pm\sigma_i^{-1}\alpha\given 1\le i\le m}\). Define gauge \(p\) on
\(O\) by declaring \(p(\sigma_i^{-1}\alpha)=1\). We can then assemble \(p\) for
all orbits \(O\) in \(R\) to obtain a gauge \(p\) on \(R\).

Still fix \(\alpha\) and \(O\), we define contraction maps \(u_i\colon W_F\to
W_{\pm\alpha}\) for each \(1\le i\le m\) by letting \(u_i(w)=W_{\pm\alpha}\) to
be the element such that
\begin{align}
    w_iw=u_i(w)w_j
\end{align}
for appropriate \(1\le j\le m\). Similarly, choose any element \(v_0\in
W_\alpha\) and if \([F_\alpha:F_{\pm\alpha}]=2\)
also an element \(v_1\in  W_{\pm\alpha}-W_\alpha\). Define contraction
\(v\colon W_{\pm\alpha}\to W_\alpha\) by
\begin{align}
    v_0u=v(u)v_j
\end{align}
for \(j=0\) or \(1\) as appropriate. Note if we choose \(v_0\) in the center of
\(W_\alpha\), then \(v\) is identity when restricted to \(W_\alpha\).

\subsection{}
Define \(1\)-cochains of \(W_F\) in \(\bbC^\x\otimes_\bbZ \Lambda\) by
\begin{align}
    r_{O,p}(w)=\prod_{i=1}^m\chi_\alpha(v(u_i(w)))^{\sigma_i^{-1}\alpha},
\end{align}
and
\begin{align}
    r_p=\prod_{O\in R/\Sigma}r_{O,p}.
\end{align}
For any gauge \(q\) on \(R\), we let
\begin{align}
    r_q=s_{q/p}r_p,
\end{align}
where \(s_{q/p}\) is the cochain defined in \Cref{sub:definition_of_s_pq}.
The following lemma is \cite{LS87}*{Lemma~2.5.A} together with the paragraph of
discussion after it.
\begin{lemma}
    We have \(\partial r_q=t_q\) as cocycles of \(W_F\) for any gauge \(q\).
    Moreover, for a fixed \(\chi\)-datum, the choices of \(\alpha\in O\),
    \(w_i\) and \(v_j\) in constructing \(r_q\) only change 
    \(r_q\) by a coboundary.
\end{lemma}

\subsection{}
Similarly, we may replace \(\chi\)-data with \(\zeta\)-data and form
\(1\)-cochain \(c_p\). By \cite{LS87}*{Corollary~2.5.B}, it is in fact a
cocycle, whose cohomology class is independent of the choices of \(\alpha, w_i,
v_j\) (hence independent of \(p\) as well). We denote this class by \(c\).

\section{Story on \texorpdfstring{\(G\)}{G}-side}
\label{sec:story_on_G_side}

\subsection{}
Let \(G\) be quasi-split induced by pointed \(\Out(\bG)\)-torsor
\(\OGT_\bG^\bullet\). We have a canonical pinning \(\bSPL\) as before, as a part
of which \(\bU_\Rt\colon \Ga\to \bG\) is the root vector associated with simple
root \(\Rt\).
Tits defines a canonical section from the Weyl group \(\bOmega\) to \(\bG\) in
\cite{Ti66}:
\index{Tits section}
\begin{align}
    \bFn\colon \bOmega\rtimes \Gamma_F &\longto \Norm_\bG(\bT)\rtimes \Gamma_F\\
    w\rtimes \sigma &\longmapsto \bFn(w)\rtimes \sigma,
\end{align}
where \(\bFn(w)\) is such that if \(w=s_{\Rt_1}\cdots s_{\Rt_m}\) is a
reduced expression, then
\begin{align}
    \bFn(w)=\bFn(s_{\Rt_1})\cdots \bFn(s_{\Rt_m}),
\end{align}
where
\begin{align}
    \bFn(s_\Rt) \defeq \bU_\Rt(1)\bU_{-\Rt}(-1)\bU_{\Rt}(1).
\end{align}
In particular, \(\bFn(1)=1\).
One can show that \(\bFn(w)\) is independent of the reduced expression hence is
well-defined, and the restriction \(\bOmega\to \Norm_\bG(\bT)\) is
\(\Gamma_F\)-equivariant.
Therefore, for \(\theta\in\bOmega\rtimes\Gamma_F\), \(\bFn(\theta)\) acts on
\(\bT\) as \(\theta\), and
\begin{align}
    t(\theta_1,\theta_2)\defeq
    \bFn(\theta_1)\bFn(\theta_2)\bFn(\theta_1\theta_2)^{-1}
\end{align}
is a \(2\)-cocycle of \(\bOmega\rtimes\Gamma_F\) in \(\bT\).
\begin{lemma}[\cite{LS87}*{Lemma~2.1.A}]
    We have that for any \(\theta_1,\theta_2\in\bOmega\rtimes\Gamma_F\),
    \begin{align}
        t(\theta_1,\theta_2)=t_{p_\bB}(\theta_1,\theta_2)=\prod_{\dual{\alpha}:1,\theta_1,\theta_2}^{p_\bB}(-1)^{\dual{\alpha}},
    \end{align}
    where \(p_\bB\) is the gauge determined by \(\bB\) on coroots \(\CoRoots\) of
    \(\bT\) in \(\bG\).
\end{lemma}
Note that \textit{a priori}, the cocycle \(t\)  depends on root vectors
\(\bU_\Rt\). However, it is not the case because the right-hand side of the
lemma above does not depend on \(\bU_\Rt\).

\subsection{}
Let \(T\subset G\) be a maximal \(F\)-torus. Let \(h\in
G(F^\sep)\simeq\bG(F^\sep)\) be a
chosen transporter from \(\bT\) to \(T\), in other words, \(\Ad_h(\bT)=T\). Then
\(h^{-1}\sigma(h)\in
\Norm_\bG(\bT)\), whose image in \(\bOmega\) is denoted by \(\omega_T(\sigma)\).
Thus, if we denote by \(\sigma_T\) the action of \(\sigma\) on \(\bT\) by transporting that on
\(T\) to \(\bT\) using \(h\), then \(\sigma_T=\omega_T(\sigma)\rtimes\sigma\in
\bOmega\rtimes\Gamma_F\). Let \(\Gamma_T\) be the group generated by \(\sigma_T\).
Clearly \(\sigma_T\) depends only on the choice of \(B=\Ad_h(\bB)\), not \(h\) itself.

The action of \(\Sigma=\Gamma_F\x\Ggen{\epsilon}\) on \(\CoRoots(G,T)\subset \CoCharG(T)\)
admits \(a\)-data (for example, any regular element \(\bFx\in \Lie(T)(F)\) defines
an \(a\)-datum by letting \(a_{\CoRt}=\Rt(\bFx)\)), which transports to
\(a\)-data for the \(\Gamma_T\x\Ggen{\epsilon}\)-action on
\(\CoRoots\subset\CoCharG\). Let \(\Set{a_{\CoRt}}\) be an \(a\)-datum.
 We have gauge \(p=p_\bB\) on \(\CoRoots\), so we have
\begin{align}
    u_p(\sigma_T)=\prod_{\CoRt:1,\sigma_T}^pa_{\CoRt}^{\CoRt},
\end{align}
whose coboundary is
\begin{align}
    t_p(\sigma_T,\tau_T)=\bFn(\sigma_T)\bFn(\tau_T)\bFn(\sigma_T\tau_T)^{-1}.
\end{align}
Since \(a_{\CoRt}^{-1}\) is also an \(a\)-datum, we also have that
\begin{align}
    \partial u_p^{-1}=t_p.
\end{align}
Thus, the map
\begin{align}
    \Gamma_T&\longto \Norm_\bG(\bT)\rtimes \Gamma_F\\
    \sigma_T &\longmapsto u_p(\sigma_T)\bFn(\sigma_T)
\end{align}
is a homomorphism, hence induces \(1\)-cocycle \(\sigma_T\mapsto
u_p(\sigma_T)\bFn(\omega_T(\sigma))\eqdef \bFm(\sigma_T)\).

Transporting by \(\Ad_{h}\), one has a map
\begin{align}
    \Gamma_F&\longto \Norm_G(T)\rtimes \Gamma_F\\
    \sigma &\longmapsto hu_p(\sigma_T)\bFn(\sigma_T)h^{-1}=h\bFm(\sigma_T)\sigma(h)^{-1}\rtimes\sigma,
\end{align}
whose image lies in \(T\rtimes\Gamma_F\). One thus obtains a \(1\)-cocycle of
\(\Gamma_F\) in \(T\), whose cohomology class in \(\RH^1(F,T)\) depends only
possibly on \(B=\Ad_h(\bB)\), not \(h\), because it is straightforward to see
that changing \(h\) changes the cocycle by a coboundary in \(T\).
In fact, a long computation in \cite{LS87}*{(2.3.3)} shows that it does not
depend on \(B\) either. Call this cohomology class in \(\RH^1(F,T)\) by
\(\lambda_T\).

\subsection{}
Here we try to summarize the construction of \(\lambda_T\) using
\(a\)-data more concisely. To begin with, we have extension
\begin{equation}\label{eqn:extension-origin}
    \begin{tikzcd}
        1\ar[r]&\bT\ar[r]&\Norm_\bG(\bT)\rtimes\Gamma_F\ar[r]&\bOmega\rtimes\Gamma_F\ar[r]&1.
    \end{tikzcd}
\end{equation}
The Tits section \(\bFn\colon \bOmega\to
\Norm_\bG(\bT)\) induces a set-theoretic section of \eqref{eqn:extension-origin},
still denoted by \(\bFn\), which in turn gives a \(2\)-cocycle of
\(\bOmega\rtimes\Gamma_F\) in \(\bT\).

Given \(T\), we choose \(B\) containing \(T\), and \(h\in G(F^\sep)\) such that
\(\Ad_{h}\) maps \((\bT,\bB)\) to \((T,B)\). Then we obtain another splitting of
\(\bOmega\rtimes\Gamma_F\) via map \(\Gamma_F\to \Gamma_T\). Restricting the
extension \eqref{eqn:extension-origin} to \(\Gamma_T\), we have
\begin{equation}
    \begin{tikzcd}
        1\ar[r] & \bT\ar[r]\ar[d,equal] & \Ad_h^{-1}(T\rtimes\Gamma_F) \ar[r]\ar[d] & \Gamma_T\ar[r]\ar[d] & 1\\
        1\ar[r] & \bT\ar[r] & \Norm_\bG(\bT)\rtimes\Gamma_F\ar[r] & \bOmega\rtimes\Gamma_F\ar[r] & 1
    \end{tikzcd},
\end{equation}
Where the second square is Cartesian and the rows are exact. This extension of
\(\Gamma_T\) by \(\bT\) is
split, and a choice of an \(a\)-datum provides a splitting \(u_p\bFn\), whose
composition with (set-theoretic) projection to \(\Norm_\bG(\bT)\) gives
\(1\)-cocycle \(\bFm\).

Transporting using \(\Ad_{h}\), one has another splitting \(\Gamma_F\to
T\rtimes\Gamma_F\) defined by \(\sigma\mapsto h\bFm(\sigma_T)\sigma(h)^{-1}\),
whose difference from the natural splitting induces a class \(\lambda_T\in
\RH^1(F,T)\).

\subsection{}
\label{sub:reformulation_of_lambda_T_in_transfer_factor}
For our convenience, we will make a slight reformulation of the above process.
Using \(h\), we have a homomorphism
\begin{align}
    \Gamma_F&\longto \Norm_{\bG}(\bT)\rtimes\Gamma_F\\
    \sigma&\longmapsto h^{-1}\sigma(h)\rtimes \sigma,
\end{align}
hence a \(1\)-cocycle \(\sigma\mapsto h^{-1}\sigma(h)\). Write
\begin{align}
    h^{-1}\sigma(h)=\bFc(\sigma)\bFn(\omega_T(\sigma)),
\end{align}
where \(\bFc(\sigma)\in\bT\), then the \(1\)-cochain \(\bFc^{-1}\) splits
\(t_p\) just like \(u_p^{-1}\).
Their difference \(u_p\bFc^{-1}\) is thus a cocycle, whose class in
\(\RH^1(F,T)\) (after transporting by \(\Ad_h\)) is exactly \(\lambda_T\).

\subsection{}
The effect of all the choices in the
construction of \(\lambda_T\) as well as its various functorial properties
can be summarized as follows (see \cite{LS87}*{(2.3)} for details):
\begin{enumerate}
    \item \(\lambda_T\) does not depend on \(h\) or \(B\).
    \item A change in pinning \(\bSPL\) modifies \(\lambda_T\) by an
        element in the image of composition of maps
        \begin{align}
            \coker[G(F)\to G_\AD(F)]\to \RH^1(F,Z)\to\RH^1(F,T).
        \end{align}
    \item A change in \(a\)-data results in a \(b\)-datum by taking quotient of two
        \(a\)-data. Forming \(1\)-cocycle \(v_p\) using the said \(b\)-datum
        (cf.~\Cref{sub:sum_a_data}), then
        \(\lambda_T\) is modified by \(hv_ph^{-1}\).
    \item The construction of \(\lambda_T\) is compatible with conjugation of
        triples \((T,B,a)\).
    \item If \(F\) is a global field, then \(G\) restricts to a quasi-split
        group \(G_v\) over any place \(v\) of \(F\), together with pinning
        \(\bSPL\). Similarly, a maximal \(F\)-torus \(T\) restricts to a maximal
        \(F_v\)-torus \(T_v\). In this case, since the action of
        \(\Gamma_v=\Gamma_{F_v}\) on \(\bSPL\) factors through natural map
        \(\Gamma_v\to\Gamma_F\), an \(a\)-datum for \(\Gamma_F\) acting on
        \(\CoRoots(G,T)\) is also an \(a\)-datum for \(\Gamma_v\) on
        \(\CoRoots(G_v,T_v)\). Let \(\lambda_{T,v}\in \RH^1(F_v,T_v)\) be the
        class defined by such \(a\)-datum, then
        \(\lambda_{T,v}\) is the image of \(\lambda_T\) under natural map
        \begin{align}
            \RH^1(F,T)\longto \RH^1(F_v,T_v).
        \end{align}
\end{enumerate}

\section{Story on \texorpdfstring{\(\LD{G}\)}{LG}-side}
\label{sec:story_on_LG_side}

\subsection{}
On the dual side we have \(\LD{G}=\dual{\bG}\rtimes W_F\) instead of
\(G\rtimes\Gamma_F\). We do not need to assume \(G\) to be quasi-split because
the construction of \(\LD{G}\) is insensitive to it.
We fix a pinning \(\dbSPL\) of \(\dual{\bG}\) as before, on which \(W_F\)
acts through \(\Gamma_F\). Let \(T\) be a maximal \(F\)-torus of \(G\), then
\(\Gamma_F\) acts on the dual torus \(\dual{T}\) over \(\bbC\), hence we have
the \(L\)-group \(\LD{T}=\dual{T}\rtimes W_F\). However, there is no canonical
embedding of \(\dual{T}\) (resp.~\(\LD{T}\)) into \(\dual{\bG}\)
(resp.~\(\LD{G}\)).

Recall that any choice of Borel subgroup \(B\subset G\) (over
\(F^\sep\)) containing \(T\) induces an isomorphism between \(T\) and the
torus in canonical based root datum \(\Psi_0\). The pinning \(\dbSPL\)
identifies \(\dual{\bT}\) with the torus in the dual of \(\Psi_0\). Therefore,
any \(B\supset T\) induces an isomorphism between \(\dual{T}\) and
\(\dual{\bT}\).

\begin{definition}
    An embedding \(\xi_T\colon \LD{T}\to\LD{G}\) is called
    \notion{admissible}\index{admissible embedding} if
    \begin{enumerate}
        \item It induces an isomorphism \(\dual{T}\to \dual{\bT}\) that is the same
            as the one induced by some choice of Borel \(B\supset T\)
            and \(\dbSPL\),
        \item It is a morphism of extensions
            \begin{equation}
                \begin{tikzcd}
                    1 \ar[r] & \dual{T}\ar[r]\ar[d] & \LD{T} \ar[r]\ar[d] &
                    W_F\ar[r]\ar[d, equal] & 1\\
                    1 \ar[r] & \dual{\bG}\ar[r] & \LD{G} \ar[r] & W_F\ar[r] & 1
                \end{tikzcd}
            \end{equation}
    \end{enumerate}
\end{definition}
Note that in general the image of \(\xi_T\) is not \(\dual{\bT}\rtimes W_F\). It
is straightforward to see that the \(\dual{\bG}\)-conjugacy class of \(\xi_T\)
is independent of \(B\) or \(\dbSPL\).

\subsection{}
We are going to attach to each \(\chi\)-datum for the \(\Gamma_F\)-action on
\(\Roots(G,T)\) an admissible embedding \(\xi_T\colon \LD{T}\to\LD{G}\), whose
\(\dual{\bG}\) conjugacy class is canonical. In this section, we will use
\(\bOmega\) to denote the Weyl group of \(\dual{\bT}\) in \(\dual{\bG}\).

To start we choose a Borel \(B\) containing \(T\). For any
\(\sigma\in\Gamma_F\),
\(\sigma(B)\) is another Borel containing \(T\), so there exists a unique \(u\)
in the Weyl group of \(T\) such that \(u^{-1}Bu=\sigma(B)\). Identifying \(u\)
as an element \(\omega_T(\sigma)\in\bOmega\) through \(\Psi_0\), then we obtain
an embedding \(\Gamma_F\to \bOmega\rtimes \Gamma_F\), hence a \(1\)-cocycle
\(\omega_T\colon \Gamma_F\to \bOmega\). We then inflate \(\omega_T\) to \(W_F\).
Let \(W_T\subset \bOmega\rtimes W_F\) be the subgroup of elements
\(\omega_T(w)\rtimes w\) where \(w\in W_F\). This way we translate the
\(W_F\)-action on \(\CharG(T)\) to the \(W_T\)-action on
\(\CoCharG(\dual{\bT})\). The section defined by Tits also makes sense on the
dual side, so we have
\begin{align}
    \dual{\bFn}\colon \bOmega &\longto \Norm_{\dual{\bG}}(\dual{\bT})
\end{align}
that is \(\Gamma_F\)-equivariant, hence also \(W_F\)-equivariant. Still use
\(\dual{\bFn}\) to denote the map \(\bOmega\rtimes
W_F\to \Norm_{\dual{\bG}}(\dual{\bT})\rtimes W_F\). 

Let \(p=p_{\dual{\bB}}\) be the gauge on \(\CoRoots(\dual{\bG},\dual{\bT})\)
determined by \(\dual{\bB}\). Then the cocycle of \(\bOmega\rtimes
W_F\) with value in \(\dual{\bT}(\bbC)\)
\begin{align}
    t_p(w_1,w_2)=\dual{\bFn}(w_1)\dual{\bFn}(w_2)\dual{\bFn}(w_1w_2)^{-1}
\end{align}
is a coboundary when restricted to \(W_T\). A choice of a \(\chi\)-datum \(\Set{\chi_\alpha}\) of
\(\Gamma_F\)-action (hence \(W_F\)-action) on \(\Roots(G,T)\) transports to a
\(\chi\)-datum of \(W_T\)-action on \(\CoRoots(\dual{\bG},\dual{\bT})\) using
\(B\). Note that
\(\Set{\chi^{-1}_\alpha}\) is also a \(\chi\)-datum, and we use it to form
\(1\)-cochain \(r_p^{-1}\), so that \(\partial r_p^{-1}=t_p\).
Thus, we obtain homomorphism
\begin{align}
    \xi_T\colon \LD{T}&\longto \LD{G}\\
    t\rtimes w&\longmapsto t_{B}r_p(w)\dual{\bFn}(w),
\end{align}
where \(t\mapsto t_{B}\) is the map \(\dual{T}\to\dual{\bT}\) induced by the
choice of \(B\) (and \(\dbSPL\)).

\subsection{}
We now briefly compare the process on \(\LD{G}\)-side with the one on \(G\)-side
to get a better overview.

Here we embed \(\LD{T}\cong \dual{\bT}\rtimes W_T\) into
\(\Norm_{\dual{\bG}}(\dual{\bT})\rtimes W_F\subset \LD{G}\) using a
\(\chi\)-datum. On \(G\)-side, the splitting \(\Gamma_T\to
\Ad_h^{-1}(T\rtimes\Gamma_F)\) given by \(a\)-data effectively induces an
``admissible'' embedding \(T\rtimes\Gamma_F\to \Norm_\bG(\bT)\rtimes\Gamma_F\)
that maps \(T\) to \(\bT\) and \(B\) to \(\bB\).

On \(G\)-side, since there is already a natural inclusion of
\(T\rtimes\Gamma_F\) in \(\bG\rtimes\Gamma_F\) to begin with (because \(T\) is a
subgroup of \(G\)), the difference between the two induces cohomology class
\(\lambda_T\). On \(\LD{G}\)-side, however, there is no natural embedding of
\(\LD{T}\) in \(\LD{G}\) to compare with, so there is seemingly one step missing
so far. It turns out the missing comparison and the resulting cohomology class
is supplied by endoscopic group (see \Cref{sub:def_of_Delta_III_2}), and the two
``dual comparisons'' lead to two terms in the
transfer factor that is also somewhat dual to each other.

\subsection{}
The effect of various choices and functorial properties of \(\xi_T\) can be
summarized as follows (see \cite{LS87}*{(2.6)}):
\begin{enumerate}
    \item For fixed \(\dbSPL\), \(B\), and \(\chi\)-data, the embedding
        \(\xi_T\) is determined up to \(\dual{\bT}\)-conjugacy.
    \item Change of \(\dbSPL\) will change \(\xi_T\) by \(\Ad_{g}\) for
        some \(g\in\dual{\bG}^{\Gamma_F}\).
    \item Change of \(B\) into \(B'=vBv^{-1}\) where \(v\in\Norm_G(T)\) will change
        \(\xi_T\) in the following way: \(\Ad_{v}\) acts on \(T\) hence on \(\dual{T}\),
        and thus on \(\dual{\bT}\) using \(\xi_T\). Call this action \(\mu\). Let
        \(g\in\Norm_{\dual{\bG}}(\dual{T})\) acts on \(\dual{T}\) as \(\mu\), then
        \(\xi'\) obtained using \(B'\) is equal to \(\Ad_{g}^{-1}\circ \xi_T\).
    \item Change of \(\chi\)-data results in a \(\zeta\)-datum by taking quotient,
        and \(\xi_T\) is multiplied by the cocycle \(c\) obtained from that
        \(\zeta\)-datum. In other words, \(\xi_T'(t\rtimes w)=c(w)\xi_T(t\rtimes w)\),
    \item If \(\Ad_{g}\) transports \((T,\chi)\) to
        \((T',\chi')\), then \(\xi_{T'}\) is simply the composition of
        \(\xi_T\) with canonical map \(\LD{T'}\to \LD{T}\) induced by \(\Ad_{g}\),
    \item Finally, if \(F\) is global, there are two ways to pass from global to
        local an admissible embedding attached to a global \(\chi\)-datum: for
        any place \(v\) of \(F\), the embedding \(\xi_T\) directly induces
        \(\xi_{T,v}\colon \LD{T_v}\to\LD{G}\) via the natural map
        \(W_{v}=W_{F_v}\to W_F\); on the other hand, a global \(\chi\)-datum
        induces a local one, which in turn induces a local admissible embedding
        \(\xi_{T,v}'\). One can show that \(\xi_{T,v}=\xi_{T,v}'\).
\end{enumerate}

\section{Endoscopic Groups}
\label{sec:Appendix_endoscopic_groups}

In this section we review the definition of endoscopic group and the transfer
factor. For the purposes of this book, we will stick with unramified groups
because it considerably simplifies the definition of the terms
\(\Delta_{\symup{III}_1}\) and \(\Delta_{\symup{III}_2}\). The reader can refer
to \cite{LS87} for more general cases.

\subsection{}
We first give a full definition of endoscopic group, or rather
\notion{endoscopic datum}, which works for any (not
necessarily unramified) reductive group \(G\).
\begin{definition}
    An \notion{endoscopic datum}\index{endoscopic!datum} is a quadruple \((H,\sH,\kappa,\xi)\) where
    \begin{enumerate}
        \item \(\kappa\in \dual{\bG}\) is semisimple,
        \item \(H\) is quasi-split reductive over \(F\), with \(L\)-group
            \(\LD{H}=\dual{\bH}\rtimes W_F\) and
            a \(\Gamma_F\)-pinning \(\bSPL_H\),
        \item \(\sH\) is a split extension of \(W_F\) by \(\dual{\bH}\), whose associated
            \(L\)-action (see \Cref{sub:action_and_induced_L_action})
            is the same as the one given by \(\LD{H}\),
        \item \(\xi\colon \sH\to\LD{G}\) is an \(L\)-embedding, i.e., a morphism of
            extensions
            \begin{equation}
                \begin{tikzcd}
                    1\ar[r]&\dual{\bH}\ar[r]\ar[d] & \sH\ar[r]\ar[d] & W_F\ar[d, equal] \ar[r]& 1\\
                    1\ar[r]&\dual{\bG}\ar[r] & \LD{G}\ar[r] & W_F \ar[r] & 1
                \end{tikzcd},
            \end{equation}
            such that the isomorphic image of \(\dual{\bH}\) is equal to
            \(\Cent_{\dual{\bG}}(\kappa)_0\) (the connected centralizer of
            \(\kappa\)), and
            that \({\Ad_{\kappa}}\circ \xi=\bFb\xi\), where \(\bFb\) is a
            \(1\)-cocycle of \(W_F\) in \(Z_{\dual{\bG}}\), inflated to \(\sH\)
            that is (resp.~locally) cohomologically trivial if \(F\) is local
            (resp.~global).
    \end{enumerate}
\end{definition}

\begin{remark}
    \label[remark]{rmk:role_of_sH_in_endoscopic_datum}
    The presence of \(\sH\) is purely due to the technical complication that it
    is not always possible to find an \(L\)-embedding from \(\LD{H}\) to
    \(\LD{G}\). In fact, the
    definition of transfer factor is only first made for the case \(\sH=\LD{H}\),
    and other cases are reduced to this special case using \(z\)-extensions of
    \(G\), and the presence of \(\sH\) in the datum ensures that there is only
    one essential way of such reduction.  See
    the discussion in \cite{LS87}*{(4.4)} for more details.
\end{remark}

\begin{remark}
    \label[remark]{rmk:role_of_Z_cocycle_a_in_endoscopic_datum}
    The presence of cocycle \(\bFb\) has no consequence in the definition of
    transfer factor. In fact, since \(\bFb\) is a coboundary, one can find some
    \(z\in Z_{\dual{\bG}}\) such that \(\bFb(\sigma)=z^{-1}\sigma(z)\) for
    \(\sigma\in W_F\). Therefore, \((H,\sH,z\kappa,\xi)\) is again an endoscopic
    datum and \(\xi\) commutes with \(z\kappa\). The transfer factor obtained
    using this datum turns out to be exactly the same as the one induced by the
    original one.
\end{remark}

\subsection{}
When \(F\) is a local field and both \(G\) and \(H\) are unramified,
the following lemma (see also
\Cref{rmk:role_of_sH_in_endoscopic_datum,rmk:role_of_Z_cocycle_a_in_endoscopic_datum})
shows that the group \(\sH\) can be omitted from the endoscopic datum. This
reconciles with the definition given in \Cref{sec:endoscopic_groups_inv_theory}.

\begin{lemma}
    If \(F\) is a local field with residue field \(k\) and
    both \(G\) and \(H\) are unramified. Suppose the following conditions are
    satisfied:
    \begin{enumerate}
        \item There exists some element
            \(w_H\in\bOmega\) such that the action of \(w_H\rtimes\Frob_k\) on
            \(\dual{\bT}\) stabilizes the based root sub-datum of \(\dual{\bH}\)
            and such action coincides with the image of
            \(\Frob_k\in\Aut(\dual{\bH})\) in \(\Out(\dual{\bH})\).
        \item The image of \(\kappa\) in \(Z_{\dual{\bH}}/Z_{\dual{\bG}}\) is
            fixed by \(\Frob_k\), and its image in
            \(\RH^1(F,Z_{\dual{\bG}})=\RH^1(\Frob_k,Z_{\dual{\bG}})\) is \(0\).
    \end{enumerate}
    Then there exists an \(L\)-embedding \(\xi\colon\LD{H}\to\LD{G}\). Moreover,
    \((H,\LD{H},\kappa,\xi)\) is an endoscopic datum.
\end{lemma}
\begin{proof}
    This is just a different phrasing of \cite{Ha93}*{Lemma~6.1}. We include a
    proof for reader's convenience.
    Since \(G\) is unramified, the action of \(W_F\) on \(\dual{\bG}\) factors
    through the Frobenius element \(\Frob_k\), and same holds for \(\dual{\bH}\).
    Therefore, we may replace \(W_F\) with \(\bbZ\Frob_k\). Note that any
    extension of \(\bbZ\Frob_k\) is split by any set-theoretic lifting of
    \(\Frob_k\) to the extension group.

    By the first assumption, there exists a homomorphism
    \begin{align}
        \bOmega_{\bH}\rtimes\bbZ\Frob_k\longto\bOmega\rtimes\bbZ\Frob_k,
    \end{align}
    such that \(\Frob_k\) is sent to \(w_H\rtimes\Frob_k\) and is compatible
    with their respective actions on \(\dual{\bT}\). Choose a lifting
    \(\dot{w}_H\in\Norm_{\dual{\bG}}(\dual{\bT})\) of \(w_H\), then \(\Frob_k\)
    acts on \(\dual{\bT}\) through
    \(\dot{w}_H\rtimes\Frob_k\in\dual{\bG}\rtimes\bbZ\Frob_k\). Call this action
    \(\Frob_H'\). Clearly \(\Frob_H'\) stabilizes the Borel subgroup of
    \(\dual{\bH}\) that is part of the pinning \(\dbSPL_H\), but not necessarily
    its root vectors. However, since simple roots of \(\dual{\bH}\) are
    linearly independent and \(\bbC\) is algebraically closed, we can always
    find some \(t\in\dual{\bT}(\bbC)\) such that \(\Frob_H\defeq t\Frob_H'\)
    stabilizes the root vectors as well. Therefore, \(\Frob_k\mapsto \Frob_H\)
    induces an \(L\)-embedding of \(L\)-groups \(\LD{H}\to\LD{G}\).
    Combined with the second assumption of the lemma, we see that
    \((H,\LD{H},\kappa,\xi)\) is an endoscopic datum.
\end{proof}

\begin{remark}
    Note that the \(L\)-embedding \(\xi\) is not canonical, as
    the element \(\Frob_H\) can be modified with an element in
    \(Z_{\dual{\bH}}(\bbC)\). Therefore, the presence of a choice of
    \(\xi\) is still necessary in endoscopic data.
\end{remark}

\subsection{}
The global case is more complicated for general groups. However, there is
still similar results in the case where the center of \(\dual{\bG}\) is connected,
in other words, the derived subgroup of \(G\) is simply-connected.

\begin{proposition}
    Suppose \(F\) is a local or global field and \(G^\Der\) is simply-connected.
    Let \((H,\sH,\kappa,\xi_\sH)\) be an endoscopic datum for \(G\).
    Then there exists an \(L\)-embedding \(\xi\colon \LD{H}\to\LD{G}\) such
    that \(\xi_\sH|_{\dual{\bH}}=\xi|_{\dual{\bH}}\) and makes
    \((H,\LD{H},\kappa,\xi)\) also an endoscopic datum for \(G\).
\end{proposition}
\begin{proof}
    The existence of \(\xi\) is proved in \cite{La79}*{Proposition~1}. That
    \(\xi\) satisfies the commutativity requirement with \(\kappa\) can be seen
    as follows: for any \(\sigma\in W_F\), let \(h_\sigma\in \sH\) be any lift
    stabilizing \(\dual{\bT}\). Then both
    \(\xi_\sH(h_\sigma)\) and \(\xi(\sigma)\) must be contained in
    \(\Norm_{\dual{\bG}}(\dual{\bT})\rtimes W_F\), and their images in
    \(\bOmega\rtimes W_F\) must be the same, because they induce the same action
    on the canonical based root datum of \(\dual{\bH}\). On the other hand, the
    cocycle \(a\) in the definition of endoscopic datum depends only on the
    element \(\sigma\) not the lift \(h_\sigma\). This means that
    \(a(\sigma)=\Ad_\kappa(\xi(\sigma))\xi(\sigma)^{-1}=\Ad_{\kappa}(\xi_\sH(\sigma))\xi_\sH(\sigma)^{-1}\)
    is contained in \(Z_{\dual{\bG}}\) and is a coboundary when restricted to
    any place \(v\) of \(F\).
\end{proof}

\subsection{}
From now on we will assume both groups \(G\) and \(H\) are unramified and
\(\sH=\LD{H}\), and omit \(\sH\) from endoscopic data.
Given \(G\) and endoscopic datum \((H,\kappa,\xi)\), without loss of generality,
we may also assume that \(\xi\) maps the pinning \(\dbSPL_{H}\)
into \(\dbSPL\).

Assuming \(\Char(F)\) is either \(0\) or larger than twice of Coxeter number of
\(G\), one can construct a canonical map from the set of semisimple conjugacy classes
of \(H(F^\sep)\) to those of \(G(F^\sep)\). Indeed, 
let \(\bT\) be the maximal torus in \(F\)-pinning \(\bSPL\) of \(G\), then by
the result of Steinberg (and our assumption on \(\Char(F)\)) we have canonical
isomorphism of \(F^\sep\)-schemes
\begin{align}
    \bT\git \bOmega\stackrel{\sim}{\longto}G\git G,
\end{align}
whose \(F^\sep\)-points parametrizes semisimple conjugacy classes in
\(G(F^\sep)\). Using \(\xi\), we may identify the maximal torus in \(\bSPL_H\)
with \(\bT\), and so we have
\begin{align}
    \bT\git\bOmega_H\stackrel{\sim}{\longto}H\git H.
\end{align}
Furthermore, \(\xi\) induces an inclusion \(\bOmega_H\rtimes
W_F\to\bOmega\rtimes W_F\). Since \(\Gamma_F\)-actions on both \(G\) and \(H\)
factors through some finite quotient, and \(W_F\) surjects on to such
quotient, we have a canonical \(\Gamma_F\)-equivariant map
\begin{align}
    \nu_H\colon\FRC_H\defeq\bT\git\bOmega_H\longto \FRC_G\defeq\bT\git\bOmega,
\end{align}
as desired.

\begin{remark}
    We did not use the fact that \(G\) and \(H\) are unramified in the above
    construction. In fact, since inner twists have no effect on adjoint
    quotient, it works for non-quasi-split group \(G\) as well.
\end{remark}

\begin{remark}
    The map we just obtained is none other than the endoscopic transfer map
    defined in \Cref{sub:endo_transfer_map_in_group_form_srs}. It also has a
    monoidal variant (with the same notation) constructed in
    \Cref{sec:endoscopic_groups_inv_theory}. In literature such as
    \cite{LS87}, it is often denoted by \(\sA_{H/G}\).
\end{remark}

\subsection{}
The GIT quotients \(\FRC_G\) and \(\FRC_H\) are constructed using torus
\(\bT\), which is defined over \(F\) when viewed as either a maximal torus of
\(G\) or \(H\). However, the \(\Gamma_F\)-actions on \(\bT\) inherited from the
two groups are not the same, so the identity map of \(\bT\) is not
\(\Gamma_F\)-equivariant.

Since \(G\) is quasi-split, there is a way to remedy
this by replacing \(\bT\subset G\) with a \(G\)-conjugate (while fixing \(\bT\subset
H\)). See \cite{Ko82}*{Corollary~2.2} for details. In literature such choice
seems to be arbitrary most of
the time. However, if a pair of matching conjugacy classes in \(H\) and \(G\) is
given (by \(\nu_H\)), there is a natural choice given by the regular
centralizer, which we will now discuss. Since \(G^\Der\) is not
simply-connected, we will restrict to strongly regular semisimple classes
defined below.

\begin{definition}
    \begin{enumerate}
        \item An element in \(G(F^\sep)\) is called \notion{strongly regular
            semisimple}\index{orbit!strongly regular semisimple} if its centralizer is a torus. A conjugacy class
            \(\FRC_G(F^\sep)\) is called \notion{strongly regular semisimple} if
            any of its representative is.
        \item An element \(\gamma_H\in H(F^\sep)\) is called \notion{strongly
            \(G\)-regular semisimple}\index{orbit!strongly \(G\)-regular semisimple} if it is semisimple, and
            \(\nu_H([\gamma_H])\in\FRC_G\) is a strongly regular semisimple
            class. In this case, the class \([\gamma_H]\) is also called
            strongly \(G\)-regular semisimple.
    \end{enumerate}
\end{definition}

\begin{remark}
    A strongly \(G\)-regular semisimple element is necessarily strongly regular
    semisimple.
\end{remark}

\subsection{}
Recall that in \Cref{sec:fundamental_lemma_for_spherical_hecke_algebras}
we defined the regular centralizer scheme \(\FRJ\) over the strongly
regular semisimple locus \(\FRC_G^\srs\), which also admits a Galois
description. In addition, there is a canonical isomorphism over the strongly
\(G\)-regular semisimple locus
\begin{align}
    \nu_H^*\FRJ|_{\FRC_H^{G\hy\srs}}\longto \FRJ_H|_{\FRC_H^{G\hy\srs}}.
\end{align}
Let \(\gamma_H\in H(F)\) be strongly \(G\)-regular semisimple and
\(\gamma\in G(F)\) is such that \(\nu_H([\gamma_H])=[\gamma]\in\FRC_G(F)\).
Then we have canonical isomorphisms of \(F\)-tori
\begin{align}
    \Cent_{H}(\gamma_H)\simeq \FRJ_{H,a_H}\simeq \FRJ_a\simeq \Cent_G(\gamma),
\end{align}
where \(a_H=[\gamma_H]\) and \(a=[\gamma]\). Let \(T_H=\Cent_H(\gamma_H)\) and
\(T=\Cent_G(\gamma)\), both defined over \(F\), we have an embedding of
\(H\)-maximal torus \(T_H\) into \(G\) with image \(T\). Such embedding is
called \notion{admissible}\index{admissible embedding}.

\section{Definition of Transfer Factors} 
\label{sec:Appendix_Definition_of_Transfer_Factors}

In this section we define transfer factors, first in the local case and then for
the global case. As before, we will restrict to the case where \(G\) and \(H\)
are both unramified, and \(\sH=\LD{H}\) in endoscopic data. We mostly follow the
original definition in \cite{LS87}, but in some part we will make some
modifications because we use the language of regular centralizer, which is not
in \cite{LS87}. The end result is the same.

\subsection{}
Let \(F\) be local and \((H,\kappa,\xi)\) an endoscopic datum for \(G\)
such that \(\kappa\in\dual{\bT}\).
Let \(G^\SC\) be the simply-connected cover of \(G^\Der\), and superscript
\(\SC\) denotes lifts to \(G^\SC\) of objects related to \(G\).
The \(F\)-pinning \(\bSPL\) induces an \(F\)-pinning \(\bSPL^\SC\) on \(G^\SC\).

Let \(\gamma_H\in H(F)\) be strongly \(G\)-regular semisimple and
\(\gamma\in G(F)\) strongly regular semisimple, and
\(\nu_H([\gamma_H])=[\gamma]\). For simplicity, we let \(c_H=[\gamma_H]\) and
\(c=[\gamma]\). We do not use \(a_H\) and \(a\) here, contrary to the notations
in the main body of this book, in order to avoid any confusion from \(a\)-data.
To define the transfer factor, we need to choose a
point \(x_c\in\bT(F^\sep)\) lying over \(c\), which also induces an
isomorphism \(\FRJ_c\simeq\bT\). A different choice of \(x_c\) changes the
isomorphism by a \(\bOmega\)-conjugation.
To be more inline with our notations in
\Cref{sec:fundamental_lemma_for_spherical_hecke_algebras}, let \(\bG\)
be the split form of \(G\), and \(G\) is obtained by a pointed
\(\Out(\bG)\)-torsor \(\OGT_G^\bullet\), corresponding to a homomorphism
\begin{align}
    \OGT_G^\bullet\colon \Gamma_F\longto \Out(\bG),
\end{align}
and a choice of \(x_c\) lifts \(\OGT_G^\bullet\) to a homomorphism
\begin{align}
    \pi_c^\bullet\colon\Gamma_F\longto \bOmega\rtimes\Out(\bG),
\end{align}
whose induced action of \(\Gamma_F\) on \(\bT\) is just the translation of the
natural \(\Gamma_F\)-action on \(\Cent_G(\gamma)\) through a suitable
conjugation. Then it induces isomorphisms of \(F\)-tori
\begin{align}
    \Cent_G(\gamma_G)\simeq \FRJ_{c}\simeq \Spec F^\sep\wedge^{\Gamma_F,\pi_c^\bullet}\bT,
\end{align}
where the right-most term simply means the categorical quotient of
\(\Spec{F^\sep}\x\bT\) by the diagonal action of \(\Gamma_F\) induced by
\(\pi_c^\bullet\).
Tate--Nakayama duality implies that we have canonical isomorphism
\begin{align}
    \RH^1(F,\FRJ_c)^*\simeq \pi_0(\dual{\bT}^{\pi_c^\bullet(\Gamma_F)}),
\end{align}
where superscript \(*\) means Pontryagin dual.
Similarly, if we let \(G^\SC\) acting on \(G\) instead, we may define regular
centralizer \(\FRJ_c^\SC\) and \(x_c\) induces isomorphism
\begin{align}
    \RH^1(F,\FRJ_c^\SC)^*\simeq \pi_0(\dual{\bT}^{\AD,\pi_c^\bullet(\Gamma_F)}).
\end{align}
By definition, the image of \(\kappa\) in \(\dual{\bT}^\AD\) is fixed by
\(\Gamma_F\), so its image in
\(\pi_0(\dual{\bT}^{\AD,\pi_c^\bullet(\Gamma_F)})\) induces a character
\begin{align}
    \RH^1(F,\FRJ_c)\longto\bbC^\x.
\end{align}

\begin{remark}
    Due to our use of regular centralizers, the admissible embedding
    \(\Cent_H(\gamma_H)\to G\) is made canonical, and the image of \(\gamma_H\)
    is exactly \(\gamma\). This is different from the treatment in \cite{LS87}
    where \(\kappa\) is treated as an element in the dual
    torus of \(T=\Cent_G(\gamma)\), while here it is treated as an element of
    \(\dual{\bT}\). They can be identified only with a choice of \(x_c\), or
    equivalently a Borel subgroup containing \(T\). Therefore, there is no
    contradiction.
\end{remark}

\subsection{}
Similar construction can be carried out on the \(H\)-side as well. Now \(x_c\)
is also a point lying over \(c_H\), which induces homomorphism
\begin{align}
    \pi_{c_H}^\bullet\colon\Gamma_F\longto \bOmega_\bH\rtimes\Out(\bH),
\end{align}
whose induces \(\Gamma_F\)-action on \(\bT\) coincides with the one by
\(\pi_c^\bullet\). We will simply use \(\pi_c^\bullet\) to denote this action
and ignore \(\pi_{c_H}^\bullet\). This is also the same action denoted by
\(\Gamma_T\) in \Cref{sec:story_on_G_side} where \(T=\Cent_G(\gamma)\).

\subsection{}
We are now ready to define the transfer factor \(\Delta(\gamma_H,\gamma_G)\) for
any \(\gamma_G\in G(F)\) lying over \(c\) (in other words, \(F\)-conjugacy
classes that is stably conjugate to \(\gamma\)). Choose an \(a\)-datum for
the \(\Gamma_F\)-action on coroots \(\CoRoots(G,\Cent_G(\gamma))\), or equivalently
the \(\pi_c^\bullet(\Gamma_F)\)-action on \(\CoRoots\).
The first term in the transfer factor is
\begin{align}
    \Delta_{\symup{I}}(\gamma_H,\gamma_G)=\Pair{\lambda_{\FRJ_c}}{\kappa},
    \nomenclature[\(Delta_I" \)]{\(\Delta_{\symup{I}}\)}{the term in the transfer
    factor by pairing \(\kappa\) with the obstruction to extending the centralixer maximal torus to a pinning}
\end{align}
where \(\lambda_{\FRJ_c}\) is computed using identification
\(\FRJ_c\simeq\Cent_G(\gamma)\), \(\bSPL^\SC\) and the chosen \(a\)-datum, and
the pairing is Tate--Nakayama duality discussed above. Note that
\(\Delta_{\symup{I}}(\gamma_H,\gamma_G)\) is independent of \(\gamma_G\), but
only \(c\) and \(x_c\).

\subsection{}
The second term is
\begin{align}
    \Delta_{\symup{II}}(\gamma_H,\gamma_G)=\prod_{[\Rt]\in
    [\Roots-\Roots_{\bH}]/\pi_c^\bullet(\Gamma_F)}\chi_\alpha\left(\frac{\Rt(x_c)-1}{a_\Rt} \right),
    \nomenclature[\(Delta_II" \)]{\(\Delta_{\symup{II}}\)}{the compensational
    term in the transfer factor for the choices of \(a\)-data and \(\chi\)-data}
\end{align}
which can be easily verified to be well-defined. Again,
\(\Delta_{\symup{II}}(\gamma_H,\gamma_G)\) is independent of \(\gamma_G\). It
is also independent of \(x_c\) because \(\pi_c^\bullet\) and \(x_c\) change
together by a \(\bOmega_\bH\)-conjugate, so a different choice only changes how
the product above is indexed.

\subsection{}
The third term of the transfer factor is split into two terms
\(\Delta_{\symup{III}_1}\) and \(\Delta_{\symup{III}_2}\).
Since \(\gamma_G\) and \(\gamma\) are stably conjugate, there exists \(h\in
G^\SC(F^\sep)\) such that \(h\gamma_Gh^{-1}=\gamma\), and the
cohomology class of cocycle \(v\colon \sigma\mapsto h\sigma(h)^{-1}\) in
\(\RH^1(F,\FRJ_c^\SC)\) is independent of \(h\). We use \(\inv(\gamma,\gamma_G)\)
for this class, and use \(x_c\) to identify it with a class in
\(\RH^1(F,\bT^\SC)\) (with \(\Gamma_F\) acting through \(\pi_c^\bullet\)). Then
the first part of the third term is
\begin{align}
    \Delta_{\symup{III}_1}(\gamma_H,\gamma_G)
    =\Pair{\inv(\gamma_H,\gamma_G)}{\kappa}^{-1}.
    \nomenclature[\(Delta_III1 \)]{\(\Delta_{\symup{III}_1}\)}{the
    \(\kappa\)-twisting term in \(\OI_\gamma^\kappa\)}
\end{align}
This is the only term dependent on \(\gamma_G\). It is also dependent on
\(x_c\), and the effect of such choice cancels out with its effect on
\(\Delta_{\symup{I}}\).

\begin{remark}
    Note that both the definitions of \(\Delta_{\symup{I}}\) and
    \(\Delta_{\symup{III}_1}\) only depend on the image of \(\kappa\) in
    \(\dual{\bT}^\AD\), which is \(\Gamma_F\)-invariant, so it is harmless
    to assume \(\kappa\) itself is also \(\Gamma_F\)-invariant
    (cf.~\Cref{rmk:role_of_Z_cocycle_a_in_endoscopic_datum}).
\end{remark}

\subsection{}
\label{sub:def_of_Delta_III_2}
The second part of the third term is the only part \(\sH=\LD{H}\) will be
used. Choose a \(\chi\)-datum for the \(\pi_c^\bullet(\Gamma_F)\)-action on
the set of roots \(\Roots\), which also induces a \(\chi\)-datum for the
\(\pi_{c_H}^\bullet(\Gamma_F)\)-action on \(\Roots_H\).
The \(\chi\)-datum (together with the point \(x_c\)) gives us admissible embeddings
\begin{align}
    \xi_{T_H}&\colon \LD{\FRJ_{H,c_H}}\longto \LD{H},\\
    \xi_T&\colon \LD{\FRJ_c}\longto \LD{G},
\end{align}
both identifying \(\dual{\FRJ}_c\) with \(\dual{\bT}\) which is compatible with
canonical identification \(\FRJ_{H,c_H}\simeq\FRJ_c\).
Composing with \(L\)-embedding \(\xi\), we have two embeddings of \(\FRJ_c\)
into \(\LD{G}\), namely \(\xi_T\) and \(\xi\circ\xi_{T_H}\). Suppose
\begin{align}
    \xi\circ\xi_{T_H}= \mu\xi_T,
\end{align}
then \(\mu\) is a \(1\)-cocycle of \(W_F\) in \(\dual{\FRJ}_c\simeq\dual{\bT}\).
Let \(\bFa\in\RH^1(F,\dual{\FRJ}_c)\) be the class of \(\mu\), then we define
\begin{align}
    \Delta_{\symup{III}_2}(\gamma_H,\gamma_G)=\Pair{\bFa}{\gamma},
    \nomenclature[\(Delta_III2 \)]{\(\Delta_{\symup{III}_2}\)}{the
    ``dual'' term in the transfer factor related to the discrepency between two
    \(L\)-embeddings of the \(L\)-group of the centralizer maximal torus}
\end{align}
where the pairing is Langlands reciprocity for tori:
\begin{align}
    \RH^1(W_F,\dual{\FRJ}_c)\simeq\Hom_{\symup{cont}}(\FRJ_c(F),\bbC^\x).
\end{align}
This term is independent of \(\gamma_G\) or \(x_c\).

\subsection{}
The final term of transfer factor is the cohomological shift given by the
discriminant valuation. Let
\begin{align}
    d(\gamma_G)&=\sum_{\Rt\in \Roots}\val_F(\Rt(x_c)-1),\\
    d_H(\gamma_H)&=\sum_{\Rt\in \Roots_\bH}\val_F(\Rt(x_c)-1),
\end{align}
which is clearly dependent only on \(c\) and \(c_H\), not \(x_c\). Then
\begin{align}
    \Delta_{\symup{IV}}(\gamma_H,\gamma_G)=q^{(-d(\gamma_G)+d_H(\gamma_H))/2}.
    \nomenclature[\(Delta_IV \)]{\(\Delta_{\symup{IV}}\)}{the
    term in the transfer factor induced by cohomological shifts}
\end{align}

\subsection{}
Finally, we can define the \inotion{transfer factor}
\begin{align}
    \Delta(\gamma_H,\gamma_G)
    =(\Delta_{\symup{I}}\Delta_{\symup{II}}\Delta_{\symup{III}_1}\Delta_{\symup{III}_2}\Delta_{\symup{IV}})(\gamma_H,\gamma_G).
    \nomenclature[\(Delta_() \)]{\(\Delta(\gamma_H,\gamma_G)\)}{the full transfer
    factor \(\Delta_{\symup{I}}\Delta_{\symup{II}}\Delta_{\symup{III}_1}\Delta_{\symup{III}_2}\Delta_{\symup{IV}}\)}
\end{align}
For convenience, we also define
\begin{align}
    \Delta_0(\gamma_H,\gamma_G)
    =\Delta_0(\gamma_H,\gamma)
    =(\Delta_{\symup{I}}\Delta_{\symup{II}}\Delta_{\symup{III}_2}\Delta_{\symup{IV}})(\gamma_H,\gamma_G).
    \nomenclature[\(Delta_0() \)]{\(\Delta_0(\gamma_H,\gamma)\)}{the modified
    transfer factor on a matching pair \((\gamma_H,\gamma)\), being \(\Delta\) without the \(\Delta_{\symup{III}}\)-term}
\end{align}
\begin{theorem}[\cite{LS87}*{Theorem~3.7.A}]
    The transfer factor \(\Delta(\gamma_H,\gamma_G)\) is independent of choice of
    \(x_c\), \(a\)-data, or \(\chi\)-data.
\end{theorem}

\begin{remark}
    \begin{enumerate}
        \item On the other hand, \(\Delta_0(\gamma_H,\gamma)\) \emph{does} depend on the
            choice of \(x_c\), but no more else than \(\Delta(\gamma_H,\gamma_G)\).
        \item Caution that \(\Delta_0\) defined here is not the same as the
            object with the same notation in \cite{LS87}. It is simply \(\Delta\) without the
            \(\Delta_{\symup{III}_1}\)-factor.
    \end{enumerate}
\end{remark}


\chapter{Review on Kashiwara Crystals} 
\label{chapA:Review_of_Kashiwara_Crystals}

In this chapter, we review several standard facts about Kashiwara crystals.
Roughly speaking, they are objects that capture the combinatorial
essence of any representation theory that resembles the classical
highest-weight theory of compact reductive Lie groups, like rational
irreducible representations of reductive groups over an algebraically closed
field of characteristic \(0\) (which will be the focus of this section),
highest-weight modules of Kac--Moody algebras,
or possibly something more general.
We shall use the book \cite{BuSc17} as our main reference due to its recency and
comprehensiveness, as well as accessibility because it establishes a purely
combinatorial approach thus avoiding the quantum groups in original constructions
independently by Kashiwara (\cites{Ka90,Ka91}) and Lusztig (\cite{Lu90a}).

\section{Definition of Crystals} 
\label{sec:Definition_of_Crystal}

Before we lay out the formal definitions, we note that \cite{BuSc17} use a
slightly different notational convention from Kashiwara's papers. We fix a
reductive group \(G\) over an algebraically closed field \(K\) of characteristic
\(0\). Since the content is entirely combinatorial, the field \(K\) itself is
not important, and we use \(\bbC\) here (in the main body of this book we use
\(\Qlb\) instead). We also fix a Borel pair \((T,B)\) in \(G\) and let
\((\CharG,\SimRts,\CoCharG,\SimCoRts)\) be the associated based root datum.

\begin{definition}
    A \notion{Kashiwara crystal}\index{Kashiwara crystal}\index{crystal!Kashiwara} (or \notion{crystal} for short)
    associated with based root datum \((\CharG,\SimRts,\CoCharG,\SimCoRts)\) is a
    non-empty set \(\sB\) such that \(0\not\in\sB\) together with the following
    maps:
    \begin{itemize}
        \item A \notion{weight map}
            \begin{align}
                \CRWT\colon \sB\longto \CharG,
            \end{align}
        \item For each simple root \(\Rt_i\in \SimRts\), a pair of
            \notion{crystal operators}\index{crystal!operator}
            \begin{align}
                e_i,f_i\colon \sB\longto \sB\sqcup\Set{0}
                \nomenclature[\(e_i_f_i \)]{\(e_i,f_i\)}{the crystal operators
                induced by the simple root \(\Rt_i\)}
            \end{align}
            together with \notion{string lengths}\index{crystal!string length}
            \begin{align}
                \epsilon_i,\phi_i\colon\sB\longto\bbZ\sqcup\Set{-\infty},
                \nomenclature[\(epsilon_i_phi_i \)]{\(\epsilon_i,\phi_i\)}{the
                string length functions in a crystal induced by the simple root \(\Rt_i\)}
            \end{align}
    \end{itemize}
    satisfying the following axioms:
    \begin{enumerate}
        \item If \(x,y\in\sB\), then \(e_i(x)=y\) if and only if \(f_i(y)=x\),
            and in this case we also have
            \begin{gather}
                \CRWT(y)=\CRWT(x)+\Rt_i,\\
                \epsilon_i(y)=\epsilon_i(x)-1,\quad \phi_i(y)=\phi_i(x)+1.
            \end{gather}
        \item For any \(i\) and any \(x\in\sB\), we have symmetry between string
            lengths
            \begin{align}
                \phi_i(x)=\Pair{\CRWT(x)}{\CoRt_i}+\epsilon_i(x).
            \end{align}
            In particular, \(\phi_i(x)=-\infty\) if and only if
            \(\epsilon_i(x)=-\infty\).
        \item If \(\epsilon_i(x)=\phi_i(x)=-\infty\), then \(e_i(x)=f_i(x)=0\).
    \end{enumerate}
\end{definition}

\begin{definition}
    A crystal \(\sB\) is called \notion{of finite type} if \(\epsilon_i\) and
    \(\phi_i\) do not take \(-\infty\) as value. It is called \notion{upper
        (resp.~lower) seminormal} if
    \begin{align}
        \epsilon_i(x)=\max\Set*{n\in \bbN\given e_i^n(x)\neq 0}\quad\text{resp.}\quad
        \phi_i(x)=\max\Set*{n\in \bbN\given f_i^n(x)\neq 0}
    \end{align}
    for all \(i\) and all \(x\). It is called
    \notion{seminormal}\index{crystal!seminormal} if it is both
    upper and lower seminormal.
\end{definition}

We will mostly focus on seminormal crystals, in which case the crystal elements
in many ways resemble (particularly carefully chosen) vector basis of
highest-weight representations, and the crystal operators look a lot like
corresponding Lie algebra elements from \(\La{sl}_2\)-triples. In fact, such
similarities can be made precise in the important case of those seminormal
crystals associated with highest-weight representations of \(G\). We will
discuss this later.

\subsection{}
To any crystal \(\sB\), we may associate a (labeled) directed graph, called the
\notion{crystal graph}\index{crystal!graph} of \(\sB\) as follows: the vertices are elements in
\(\sB\), and if \(f_i(x)=y\) for some \(i\) and \(x,y\in\sB\), then there is an
edge from \(x\) to \(y\) with label \(i\). If \(\sB\) is seminormal, then it is
completely determined by its crystal graph (provided we also label the vertices
by their weights) because the string lengths
\(\epsilon_i\) and \(\phi_i\) are determined by the corresponding crystal
operators.

For each connected component of the crystal graph of \(\sB\), there is a natural
crystal structure on the set of its vertices inherited from \(\sB\) in an
obvious way. We call such crystal \(\sB'\) a \notion{connected component} of
\(\sB\).

\begin{definition}
    A \notion{full subcrystal} of crystal \(\sB\) is a union of connected
    components of \(\sB\).
\end{definition}

\begin{definition}
    An element \(x\in\sB\) is called \notion{dominant} if \(\CRWT(x)\) is, and
    a \notion{highest-weight element} if \(e_i(x)=0\) for all \(i\). The
    crystal \(\sB\) is called a \notion{highest-weight
    crystal}\index{crystal!highest-weight} if it is
    connected and contains a unique highest-weight element.
\end{definition}

\begin{remark}
    If \(\sB\) is seminormal, highest-weight elements are automatically dominant.
\end{remark}

\begin{example}
    \label{exa:two_crystals_of_SL_3}
    Here are two examples of seminormal highest-weight crystals of
    \(\SL_3\), where the highest weight equals the highest root
    \(\Rt_1+\Rt_2=\Wt_1+\Wt_2\):
    \begin{equation}
        \begin{tikzpicture}
            \node[circle, draw] at (AII cs:x=0.15,y=0.15) (H1) {};
            \node[circle, draw] at (AII cs:x=-0.15,y=-0.15) (H2) {};
            \node[circle, draw] at (AII cs:x=1,y=0) (E1) {};
            \node[circle, draw] at (AII cs:x=0,y=1) (E2) {};
            \node[circle, draw] at (AII cs:x=1,y=1) (E3) {};
            \node[circle, draw] at (AII cs:x=-1,y=0) (F1) {};
            \node[circle, draw] at (AII cs:x=0,y=-1) (F2) {};
            \node[circle, draw] at (AII cs:x=-1,y=-1) (F3) {};

            \draw[-Stealth] (E3) -- (E2) node[midway, above]{\(1\)};
            \draw[-Stealth] (E1) -- (H1) node[midway, above]{\(1\)};
            \draw[-Stealth] (H1) -- (F1) node[midway, above]{\(1\)};
            \draw[-Stealth] (F2) -- (F3) node[midway, above]{\(1\)};
            \draw[-Stealth,dashed] (E3) -- (E1) node[midway, above right]{\(2\)};
            \draw[-Stealth,dashed] (E2) -- (H2) node[midway, above right]{\(2\)};
            \draw[-Stealth,dashed] (H2) -- (F2) node[midway, above right]{\(2\)};
            \draw[-Stealth,dashed] (F1) -- (F3) node[midway, above right]{\(2\)};
        \end{tikzpicture}
        \qquad\qquad\qquad\qquad
        \begin{tikzpicture}
            \node[circle, draw] at (AII cs:x=0,y=0) (H) {};
            \node[circle, draw] at (AII cs:x=1,y=0) (E1) {};
            \node[circle, draw] at (AII cs:x=0,y=1) (E2) {};
            \node[circle, draw] at (AII cs:x=1,y=1) (E3) {};
            \node[circle, draw] at (AII cs:x=-1,y=0) (F1) {};
            \node[circle, draw] at (AII cs:x=0,y=-1) (F2) {};
            \node[circle, draw] at (AII cs:x=-1,y=-1) (F3) {};

            \draw[-Stealth] (E3) -- (E2) node[midway, above]{\(1\)};
            \draw[-Stealth] (E1) -- (H) node[midway, above]{\(1\)};
            \draw[-Stealth] (H) -- (F1) node[midway, above]{\(1\)};
            \draw[-Stealth] (F2) -- (F3) node[midway, above]{\(1\)};
            \draw[-Stealth,dashed] (E3) -- (E1) node[midway, above right]{\(2\)};
            \draw[-Stealth,dashed] (E2) -- (H) node[midway, above right]{\(2\)};
            \draw[-Stealth,dashed] (H) -- (F2) node[midway, above right]{\(2\)};
            \draw[-Stealth,dashed] (F1) -- (F3) node[midway, above right]{\(2\)};
        \end{tikzpicture}
    \end{equation}
    The crystal on the left turns out to be the unique crystal associated with
    the adjoint representation of \(\SL_3\). The crystal on the right, however,
    may be viewed as a degeneration of the left one by collapsing the two
    weight-\(0\) elements in the center, and it does not
    correspond to any highest-weight representation of \(\SL_3\).
\end{example}

\subsection{}
\Cref{exa:two_crystals_of_SL_3} shows that the category of seminormal crystals
is still too large for us since it contains objects that are not related to
representations. To fix this, we need to further refine the notion into
\notion{normal crystals}. Such refinement is subtle, and quite technical if we
want the combinatorial characterization. Therefore, instead we will digress from
\cite{BuSc17} and adopt a categorical approach.

\begin{definition}
    Let \(\sB_1,\sB_2\) be two crystals of \(G\). The \notion{tensor product}
    \(\sB_1\otimes\sB_2\) is the crystal with underlying set \(\sB_1\x\sB_2\),
    whose elements are denoted by \(x_1\otimes x_2\), such that
    \begin{itemize}
        \item The weight function is defined by
            \begin{align}
                \CRWT(x_1\otimes x_2)=\CRWT(x_1)+\CRWT(x_2).
            \end{align}
        \item The crystal operators are
            \begin{align}
                e_i(x_1\otimes x_2)&=\begin{cases}
                    e_i(x_1)\otimes x_2 & \phi_i(x_2)<\epsilon_i(x_1),\\
                    x_1\otimes e_i(x_2) & \phi_i(x_2)\ge\epsilon_i(x_1).
                \end{cases},\\
                f_i(x_1\otimes x_2)&=\begin{cases}
                    f_i(x_1)\otimes x_2 & \phi_i(x_2)\le\epsilon_i(x_1),\\
                    x_1\otimes f_i(x_2) & \phi_i(x_2)>\epsilon_i(x_1).
                \end{cases},
            \end{align}
            where \(x\otimes 0=0\otimes x=0\) by convention.
        \item The string lengths are
            \begin{align}
                \epsilon_i(x_1\otimes x_2)&=\max\Set*{\epsilon_i(x_2),
                \epsilon_i(x_1)-\Pair{\CRWT(x_2)}{\CoRt_i}},\\
                \phi_i(x_1\otimes x_2)&=\max\Set*{\phi_i(x_1),
                \phi_i(x_2)+\Pair{\CRWT(x_1)}{\CoRt_i}}.
            \end{align}
    \end{itemize}
\end{definition}
One can check that the tensor product of two crystals is indeed a
crystal, and that the tensor product is associative using the identification map
\((x\otimes y)\otimes z\mapsto x\otimes(y\otimes z)\) on the underlying
sets. Moreover, the tensor product of two seminormal crystals is again
seminormal.

\subsection{}
Similar to representations, we may define the \notion{character} of crystal
\(\sB\) as the sum
\begin{align}
    \chi_\sB\defeq \sum_{x\in\sB}\CRWT(x)\in \bbC[T].
\end{align}

\section{Normal Crystals}
The miraculous result about crystals is that for any semisimple and
simply-connected \(G\) there is a unique (up to equivalence) groupoid
\(\NC=\NC_G\) of \(G\)-crystals such that
\begin{enumerate}
    \item For each highest-weight representation \(V_\lambda\) of \(G\), there is a unique
        isomorphism class of irreducible crystals \(\sB_\lambda\in \NC\)
        \nomenclature[\(B"scr_lambda \)]{\(\sB_\lambda\)}{the unique normal
        crystal associated with the irreducible representation of highest weight \(\lambda\)}
        with the same character as the representation.
    \item Any object in \(\NC\) is isomorphic to a finite disjoint union of
        \(\sB_\lambda\) for various \(\lambda\)'s.
    \item \(\NC\) is closed under tensor product. In other words, for any
        \(\sB,\sB'\in\NC\), the tensor product \(\sB\otimes\sB'\) is isomorphic
        to an object in \(\NC\).
\end{enumerate}
This is a deep result. See \cite{Jo95} for a proof. Note that for general \(G\),
we may restrict \(G\)-crystals to obtain a \(G^\SC\)-crystal induced by map
\(\CharG\to\CharG^\SC\).

\begin{definition}
    A \(G\)-crystal \(\sB\) is called \notion{normal}\index{crystal!normal} if it is isomorphic to an
    object in \(\NC_{G^\SC}\) after restricting to \(G^\SC\).
\end{definition}

\subsection{}
\label{sub:Appendix_actions_on_crystals}
For a normal crystal \(\sB\) there is a canonical Weyl group action
on the elements of \(\sB\) as follows: for any simple reflection \(s_i\) and
\(x\in\sB\), we define
\begin{align}
    s_i(x)=\begin{cases}
        f_i^{\Pair{\CRWT(x)}{\CoRt_i}}(x) & \Pair{\CRWT(x)}{\CoRt_i}\ge 0\\
        e_i^{-\Pair{\CRWT(x)}{\CoRt_i}}(x) & \Pair{\CRWT(x)}{\CoRt_i}<0
    \end{cases}
\end{align}
It turns out that this action of \(s_i\) is an involution, and they satisfy the
braid relations. As a result, it defines a \(W\)-action on \(\sB\). This implies
that for any \(x\in\sB\), there is a unique dominant element in the
\(W\)-orbit of \(x\). 

If, in addition, the normal crystal \(\sB\) consists of irreducible ones with
distinct highest weights stable under \(\Out(G)\), then we have a natural
action of \(\Out(G)\) on \(\sB\) by permuting vertices and crystal
operators in an obvious way. Clearly, it further extends to an action of
\(W\rtimes\Out(G)\).

\subsection{}
Let \(L_I\) be the standard Levi subgroup of \(G\) corresponding to a subset
\(I\subset\SimRts\) of simple roots, then we have a canonical functor
\begin{align}
    \NC_G\longto \NC_{L_I}
\end{align}
by simply deleting the edges in the crystal graphs corresponding to roots not in
\(I\). However, unlike representations, there is no straightforward way to
define such ``restriction'' functor for arbitrary group homomorphisms (or even
subgroups). For example, our
\Cref{thm:asymptotic_FL} can be viewed
as some sort of ``endoscopic \(L\)-restriction'', but it is not a purely
combinatorial result and its interaction with tensor product seems very
complicated.

\subsection{}
A particularly relevant model of \(\NC_G\) is by using MV-cycles in the affine
Grassmannian \(\Gr_{\dual{\bG}}\) of the group \(\dual{\bG}\) dual to
\(G\). The topic of MV-cycles as crystal basis elements is intensely studied in
recent years, and it is beyond us to give a comprehensive review. The reader can refer to
\cite{BG01} for the original proof.



\cleardoublepage
\phantomsection
\addcontentsline{toc}{chapter}{Bibliography}
\bibliography{MHFL_arxiv}

\begin{bibdiv}
\begin{biblist}

\bib{BoCe22}{article}{
      author={Bouthier, Alexis},
      author={{\v C}esnavi{\v c}ius, K{\k e}stutis},
       title={Torsors on loop groups and the {H}itchin fibration},
        date={2022},
        ISSN={0012-9593,1873-2151},
     journal={Ann. Sci. \'Ec. Norm. Sup\'er. (4)},
      volume={55},
      number={3},
       pages={791\ndash 864},
      review={\MR{4553656}},
}

\bib{BBD82}{incollection}{
      author={Be\u{i}linson, A.~A.},
      author={Bernstein, J.},
      author={Deligne, P.},
       title={Faisceaux pervers},
        date={1982},
   booktitle={Analysis and topology on singular spaces, {I} ({L}uminy, 1981)},
      series={Ast\'{e}risque},
      volume={100},
   publisher={Soc. Math. France, Paris},
       pages={5\ndash 171},
      review={\MR{751966}},
}

\bib{BoCh18}{article}{
      author={Bouthier, Alexis},
      author={Chi, Jingren},
       title={Correction to ``{D}imension des fibres de {S}pringer affines pour les groupes'' [{MR}3376144]},
        date={2018},
        ISSN={1083-4362},
     journal={Transform. Groups},
      volume={23},
      number={4},
       pages={1217\ndash 1222},
         url={https://doi.org/10.1007/s00031-018-9496-3},
      review={\MR{3869433}},
}

\bib{Be03}{article}{
      author={Behrend, Kai~A.},
       title={Derived {$l$}-adic categories for algebraic stacks},
        date={2003},
        ISSN={0065-9266},
     journal={Mem. Amer. Math. Soc.},
      volume={163},
      number={774},
       pages={viii+93},
         url={https://doi.org/10.1090/memo/0774},
      review={\MR{1963494}},
}

\bib{Be93}{article}{
      author={Behrend, Kai~A.},
       title={The {L}efschetz trace formula for algebraic stacks},
        date={1993},
        ISSN={0020-9910},
     journal={Invent. Math.},
      volume={112},
      number={1},
       pages={127\ndash 149},
         url={https://doi.org/10.1007/BF01232427},
      review={\MR{1207479}},
}

\bib{Be96}{article}{
      author={Bezrukavnikov, Roman},
       title={The dimension of the fixed point set on affine flag manifolds},
        date={1996},
        ISSN={1073-2780},
     journal={Math. Res. Lett.},
      volume={3},
      number={2},
       pages={185\ndash 189},
         url={https://doi.org/10.4310/MRL.1996.v3.n2.a5},
      review={\MR{1386839}},
}

\bib{BG01}{article}{
      author={Braverman, Alexander},
      author={Gaitsgory, Dennis},
       title={Crystals via the affine {G}rassmannian},
        date={2001},
        ISSN={0012-7094,1547-7398},
     journal={Duke Math. J.},
      volume={107},
      number={3},
       pages={561\ndash 575},
         url={https://doi.org/10.1215/S0012-7094-01-10736-9},
      review={\MR{1828302}},
}

\bib{BSL90}{book}{
      author={Bosch, Siegfried},
      author={L\"{u}tkebohmert, Werner},
      author={Raynaud, Michel},
       title={N\'{e}ron models},
      series={Ergebnisse der Mathematik und ihrer Grenzgebiete (3) [Results in Mathematics and Related Areas (3)]},
   publisher={Springer-Verlag, Berlin},
        date={1990},
      volume={21},
        ISBN={3-540-50587-3},
         url={https://doi.org/10.1007/978-3-642-51438-8},
      review={\MR{1045822}},
}

\bib{BNS16}{article}{
      author={Bouthier, A.},
      author={Ng\^{o}, B.~C.},
      author={Sakellaridis, Y.},
       title={On the formal arc space of a reductive monoid},
        date={2016},
        ISSN={0002-9327},
     journal={Amer. J. Math.},
      volume={138},
      number={1},
       pages={81\ndash 108},
         url={https://doi.org/10.1353/ajm.2016.0004},
      review={\MR{3462881}},
}

\bib{Bo15}{article}{
      author={Bouthier, Alexis},
       title={Dimension des fibres de {S}pringer affines pour les groupes},
        date={2015},
        ISSN={1083-4362},
     journal={Transform. Groups},
      volume={20},
      number={3},
       pages={615\ndash 663},
         url={https://doi.org/10.1007/s00031-015-9326-9},
      review={\MR{3376144}},
}

\bib{Bo17}{article}{
      author={Bouthier, Alexis},
       title={La fibration de {H}itchin-{F}renkel-{N}g\^{o} et son complexe d'intersection},
        date={2017},
        ISSN={0012-9593},
     journal={Ann. Sci. \'{E}c. Norm. Sup\'{e}r. (4)},
      volume={50},
      number={1},
       pages={85\ndash 129},
         url={https://doi.org/10.24033/asens.2316},
      review={\MR{3621427}},
}

\bib{Bou68}{book}{
      author={Bourbaki, N.},
       title={\'{E}l\'{e}ments de math\'{e}matique. {F}asc. {XXXIV}. {G}roupes et alg\`ebres de {L}ie. {C}hapitre {IV}: {G}roupes de {C}oxeter et syst\`emes de {T}its. {C}hapitre {V}: {G}roupes engendr\'{e}s par des r\'{e}flexions. {C}hapitre {VI}: syst\`emes de racines},
      series={Actualit\'{e}s Scientifiques et Industrielles [Current Scientific and Industrial Topics], No. 1337},
   publisher={Hermann, Paris},
        date={1968},
      review={\MR{0240238}},
}

\bib{BiRa94}{article}{
      author={Biswas, I.},
      author={Ramanan, S.},
       title={An infinitesimal study of the moduli of {H}itchin pairs},
        date={1994},
        ISSN={0024-6107,1469-7750},
     journal={J. London Math. Soc. (2)},
      volume={49},
      number={2},
       pages={219\ndash 231},
         url={https://doi.org/10.1112/jlms/49.2.219},
      review={\MR{1260109}},
}

\bib{Br09}{article}{
      author={Brion, Michel},
       title={Vanishing theorems for {D}olbeault cohomology of log homogeneous varieties},
        date={2009},
        ISSN={0040-8735},
     journal={Tohoku Math. J. (2)},
      volume={61},
      number={3},
       pages={365\ndash 392},
         url={https://doi.org/10.2748/tmj/1255700200},
      review={\MR{2568260}},
}

\bib{BuSc17}{book}{
      author={Bump, Daniel},
      author={Schilling, Anne},
       title={Crystal bases},
   publisher={World Scientific Publishing Co. Pte. Ltd., Hackensack, NJ},
        date={2017},
        ISBN={978-981-4733-44-1},
         url={https://doi.org/10.1142/9876},
        note={Representations and combinatorics},
      review={\MR{3642318}},
}

\bib{BoTi65}{article}{
      author={Borel, Armand},
      author={Tits, Jacques},
       title={Groupes r\'eductifs},
        date={1965},
        ISSN={0073-8301,1618-1913},
     journal={Inst. Hautes \'Etudes Sci. Publ. Math.},
      number={27},
       pages={55\ndash 150},
         url={http://www.numdam.org/item?id=PMIHES_1965__27__55_0},
      review={\MR{207712}},
}

\bib{BZSV24}{misc}{
      author={Ben-Zvi, David},
      author={Sakellaridis, Yiannis},
      author={Venkatesh, Akshay},
       title={Relative langlands duality},
        date={2024},
         url={https://arxiv.org/abs/2409.04677},
}

\bib{CaCeHa19}{article}{
      author={Casselman, William},
      author={Cely, Jorge~E.},
      author={Hales, Thomas},
       title={The spherical {H}ecke algebra, partition functions, and motivic integration},
        date={2019},
        ISSN={0002-9947,1088-6850},
     journal={Trans. Amer. Math. Soc.},
      volume={371},
      number={9},
       pages={6169\ndash 6212},
         url={https://doi.org/10.1090/tran/7465},
      review={\MR{3937321}},
}

\bib{Ch22}{article}{
      author={Chi, Jingren},
       title={Geometry of {K}ottwitz-{V}iehmann varieties},
        date={2022},
        ISSN={1474-7480,1475-3030},
     journal={J. Inst. Math. Jussieu},
      volume={21},
      number={1},
       pages={1\ndash 65},
         url={https://doi.org/10.1017/S1474748019000604},
      review={\MR{4366333}},
}

\bib{CL10}{article}{
      author={Chaudouard, Pierre-Henri},
      author={Laumon, G\'{e}rard},
       title={Le lemme fondamental pond\'{e}r\'{e}. {I}. {C}onstructions g\'{e}om\'{e}triques},
        date={2010},
        ISSN={0010-437X},
     journal={Compos. Math.},
      volume={146},
      number={6},
       pages={1416\ndash 1506},
         url={https://doi.org/10.1112/S0010437X10004756},
      review={\MR{2735371}},
}

\bib{De96}{article}{
      author={Debarre, Olivier},
       title={Th\'{e}or\`emes de connexit\'{e} pour les produits d'espaces projectifs et les grassmanniennes},
        date={1996},
        ISSN={0002-9327},
     journal={Amer. J. Math.},
      volume={118},
      number={6},
       pages={1347\ndash 1367},
         url={http://muse.jhu.edu/journals/american_journal_of_mathematics/v118/118.6debarre.pdf},
      review={\MR{1420927}},
}

\bib{DoGa02}{article}{
      author={Donagi, R.~Y.},
      author={Gaitsgory, D.},
       title={The gerbe of {H}iggs bundles},
        date={2002},
        ISSN={1083-4362},
     journal={Transform. Groups},
      volume={7},
      number={2},
       pages={109\ndash 153},
         url={https://doi.org/10.1007/s00031-002-0008-z},
      review={\MR{1903115}},
}

\bib{DeGa70}{book}{
      author={Demazure, Michel},
      author={Gabriel, Pierre},
       title={Groupes alg\'ebriques. {T}ome {I}: {G}\'eom\'etrie alg\'ebrique, g\'en\'eralit\'es, groupes commutatifs},
   publisher={Masson \& Cie, \'Editeurs, Paris; North-Holland Publishing Co., Amsterdam},
        date={1970},
        note={Avec un appendice {\it Corps de classes local}\ par Michiel Hazewinkel},
      review={\MR{302656}},
}

\bib{DeLu76}{article}{
      author={Deligne, P.},
      author={Lusztig, G.},
       title={Representations of reductive groups over finite fields},
        date={1976},
        ISSN={0003-486X},
     journal={Ann. of Math. (2)},
      volume={103},
      number={1},
       pages={103\ndash 161},
         url={https://doi.org/10.2307/1971021},
      review={\MR{393266}},
}

\bib{FaGoIl05}{book}{
      author={Fantechi, Barbara},
      author={G\"ottsche, Lothar},
      author={Illusie, Luc},
      author={Kleiman, Steven~L.},
      author={Nitsure, Nitin},
      author={Vistoli, Angelo},
       title={Fundamental algebraic geometry},
      series={Mathematical Surveys and Monographs},
   publisher={American Mathematical Society, Providence, RI},
        date={2005},
      volume={123},
        ISBN={0-8218-3541-6},
         url={https://doi.org/10.1090/surv/123},
        note={Grothendieck's FGA explained},
      review={\MR{2222646}},
}

\bib{FN11}{article}{
      author={Frenkel, Edward},
      author={Ng\^{o}, Bao~Ch\^{a}u},
       title={Geometrization of trace formulas},
        date={2011},
        ISSN={1664-3607},
     journal={Bull. Math. Sci.},
      volume={1},
      number={1},
       pages={129\ndash 199},
         url={https://doi.org/10.1007/s13373-011-0009-0},
      review={\MR{2823791}},
}

\bib{Fu69}{article}{
      author={Fulton, William},
       title={Hurwitz schemes and irreducibility of moduli of algebraic curves},
        date={1969},
        ISSN={0003-486X},
     journal={Ann. of Math. (2)},
      volume={90},
       pages={542\ndash 575},
         url={https://doi.org/10.2307/1970748},
      review={\MR{260752}},
}

\bib{GKM04}{article}{
      author={Goresky, Mark},
      author={Kottwitz, Robert},
      author={Macpherson, Robert},
       title={Homology of affine {S}pringer fibers in the unramified case},
        date={2004},
        ISSN={0012-7094},
     journal={Duke Math. J.},
      volume={121},
      number={3},
       pages={509\ndash 561},
         url={https://doi.org/10.1215/S0012-7094-04-12135-9},
      review={\MR{2040285}},
}

\bib{GKM09}{article}{
      author={Goresky, Mark},
      author={Kottwitz, Robert},
      author={MacPherson, Robert},
       title={Codimensions of root valuation strata},
        date={2009},
        ISSN={1558-8599},
     journal={Pure Appl. Math. Q.},
      volume={5},
      number={4, Special Issue: In honor of John Tate. Part 1},
       pages={1253\ndash 1310},
         url={https://doi.org/10.4310/PAMQ.2009.v5.n4.a4},
      review={\MR{2560317}},
}

\bib{EGAIII1}{article}{
      author={Grothendieck, Alexander},
       title={\'{E}l\'{e}ments de g\'{e}om\'{e}trie alg\'{e}brique. {III}. \'{E}tude cohomologique des faisceaux coh\'{e}rents. {I}},
        date={1961},
        ISSN={0073-8301},
     journal={Inst. Hautes \'{E}tudes Sci. Publ. Math.},
      number={11},
       pages={167},
         url={http://www.numdam.org/item?id=PMIHES_1961__11__167_0},
      review={\MR{217085}},
}

\bib{EGAIV3}{article}{
      author={Grothendieck, A.},
       title={\'{E}l\'{e}ments de g\'{e}om\'{e}trie alg\'{e}brique. {IV}. \'{E}tude locale des sch\'{e}mas et des morphismes de sch\'{e}mas. {III}},
        date={1966},
        ISSN={0073-8301},
     journal={Inst. Hautes \'{E}tudes Sci. Publ. Math.},
      number={28},
       pages={255},
         url={http://www.numdam.org/item?id=PMIHES_1966__28__255_0},
      review={\MR{217086}},
}

\bib{SGA1}{book}{
      author={Grothendieck, Alexander},
       title={Rev\^{e}tements \'{e}tales et groupe fondamental},
      series={Lecture Notes in Mathematics, Vol. 224},
   publisher={Springer-Verlag, Berlin-New York},
        date={1971},
        note={S\'{e}minaire de G\'{e}om\'{e}trie Alg\'{e}brique du Bois Marie 1960--1961 (SGA 1), Dirig\'{e} par Alexandre Grothendieck. Augment\'{e} de deux expos\'{e}s de M. Raynaud},
      review={\MR{0354651}},
}

\bib{Ha93}{incollection}{
      author={Hales, Thomas~C.},
       title={A simple definition of transfer factors for unramified groups},
        date={1993},
   booktitle={Representation theory of groups and algebras},
      series={Contemp. Math.},
      volume={145},
   publisher={Amer. Math. Soc., Providence, RI},
       pages={109\ndash 134},
         url={https://doi.org/10.1090/conm/145/1216184},
      review={\MR{1216184}},
}

\bib{Ha95}{article}{
      author={Hales, Thomas~C.},
       title={On the fundamental lemma for standard endoscopy: reduction to unit elements},
        date={1995},
        ISSN={0008-414X,1496-4279},
     journal={Canad. J. Math.},
      volume={47},
      number={5},
       pages={974\ndash 994},
         url={https://doi.org/10.4153/CJM-1995-051-5},
      review={\MR{1350645}},
}

\bib{HM02}{article}{
      author={Hurtubise, J.~C.},
      author={Markman, E.},
       title={Elliptic {S}klyanin integrable systems for arbitrary reductive groups},
        date={2002},
        ISSN={1095-0761},
     journal={Adv. Theor. Math. Phys.},
      volume={6},
      number={5},
       pages={873\ndash 978 (2003)},
         url={https://doi.org/10.4310/ATMP.2002.v6.n5.a4},
      review={\MR{1974589}},
}

\bib{Hu95}{book}{
      author={Humphreys, James~E.},
       title={Conjugacy classes in semisimple algebraic groups},
      series={Mathematical Surveys and Monographs},
   publisher={American Mathematical Society, Providence, RI},
        date={1995},
      volume={43},
        ISBN={0-8218-0333-6},
         url={https://doi.org/10.1090/surv/043},
      review={\MR{1343976}},
}

\bib{Il71}{book}{
      author={Illusie, Luc},
       title={Complexe cotangent et d\'eformations. {I}},
      series={Lecture Notes in Mathematics},
   publisher={Springer-Verlag, Berlin-New York},
        date={1971},
      volume={Vol. 239},
      review={\MR{491680}},
}

\bib{Jo95}{book}{
      author={Joseph, Anthony},
       title={Quantum groups and their primitive ideals},
      series={Ergebnisse der Mathematik und ihrer Grenzgebiete (3) [Results in Mathematics and Related Areas (3)]},
   publisher={Springer-Verlag, Berlin},
        date={1995},
      volume={29},
        ISBN={3-540-57057-8},
         url={https://doi.org/10.1007/978-3-642-78400-2},
      review={\MR{1315966}},
}

\bib{Ka90}{article}{
      author={Kashiwara, Masaki},
       title={Crystalizing the {$q$}-analogue of universal enveloping algebras},
        date={1990},
        ISSN={0010-3616,1432-0916},
     journal={Comm. Math. Phys.},
      volume={133},
      number={2},
       pages={249\ndash 260},
         url={http://projecteuclid.org/euclid.cmp/1104201397},
      review={\MR{1090425}},
}

\bib{Ka91}{article}{
      author={Kashiwara, M.},
       title={On crystal bases of the {$Q$}-analogue of universal enveloping algebras},
        date={1991},
        ISSN={0012-7094,1547-7398},
     journal={Duke Math. J.},
      volume={63},
      number={2},
       pages={465\ndash 516},
         url={https://doi.org/10.1215/S0012-7094-91-06321-0},
      review={\MR{1115118}},
}

\bib{KaLu88}{article}{
      author={Kazhdan, D.},
      author={Lusztig, G.},
       title={Fixed point varieties on affine flag manifolds},
        date={1988},
        ISSN={0021-2172},
     journal={Israel J. Math.},
      volume={62},
      number={2},
       pages={129\ndash 168},
         url={https://doi.org/10.1007/BF02787119},
      review={\MR{947819}},
}

\bib{Ko82}{article}{
      author={Kottwitz, Robert~E.},
       title={Rational conjugacy classes in reductive groups},
        date={1982},
        ISSN={0012-7094},
     journal={Duke Math. J.},
      volume={49},
      number={4},
       pages={785\ndash 806},
         url={http://projecteuclid.org/euclid.dmj/1077315531},
      review={\MR{683003}},
}

\bib{Ko85}{article}{
      author={Kottwitz, Robert~E.},
       title={Isocrystals with additional structure},
        date={1985},
        ISSN={0010-437X},
     journal={Compositio Math.},
      volume={56},
      number={2},
       pages={201\ndash 220},
         url={http://www.numdam.org/item?id=CM_1985__56_2_201_0},
      review={\MR{809866}},
}

\bib{Ko97}{article}{
      author={Kottwitz, Robert~E.},
       title={Isocrystals with additional structure. {II}},
        date={1997},
        ISSN={0010-437X},
     journal={Compositio Math.},
      volume={109},
      number={3},
       pages={255\ndash 339},
         url={https://doi.org/10.1023/A:1000102604688},
      review={\MR{1485921}},
}

\bib{Ko99}{article}{
      author={Kottwitz, Robert~E.},
       title={Transfer factors for {L}ie algebras},
        date={1999},
     journal={Represent. Theory},
      volume={3},
       pages={127\ndash 138},
         url={https://doi.org/10.1090/S1088-4165-99-00068-0},
      review={\MR{1703328}},
}

\bib{KS99}{article}{
      author={Kottwitz, Robert~E.},
      author={Shelstad, Diana},
       title={Foundations of twisted endoscopy},
        date={1999},
        ISSN={0303-1179},
     journal={Ast\'{e}risque},
      number={255},
       pages={vi+190},
      review={\MR{1687096}},
}

\bib{KoVi12}{article}{
      author={Kottwitz, Robert},
      author={Viehmann, Eva},
       title={Generalized affine {S}pringer fibres},
        date={2012},
        ISSN={1474-7480},
     journal={J. Inst. Math. Jussieu},
      volume={11},
      number={3},
       pages={569\ndash 609},
         url={https://doi.org/10.1017/S147474801100020X},
      review={\MR{2931318}},
}

\bib{La79}{article}{
      author={Langlands, R.~P.},
       title={Stable conjugacy: definitions and lemmas},
        date={1979},
        ISSN={0008-414X,1496-4279},
     journal={Canadian J. Math.},
      volume={31},
      number={4},
       pages={700\ndash 725},
         url={https://doi.org/10.4153/CJM-1979-069-2},
      review={\MR{540901}},
}

\bib{LN08}{article}{
      author={Laumon, G\'{e}rard},
      author={Ng\^{o}, Bao~Ch\^{a}u},
       title={Le lemme fondamental pour les groupes unitaires},
        date={2008},
        ISSN={0003-486X},
     journal={Ann. of Math. (2)},
      volume={168},
      number={2},
       pages={477\ndash 573},
         url={https://doi.org/10.4007/annals.2008.168.477},
      review={\MR{2434884}},
}

\bib{LaOl08I}{article}{
      author={Laszlo, Yves},
      author={Olsson, Martin},
       title={The six operations for sheaves on {A}rtin stacks. {I}. {F}inite coefficients},
        date={2008},
        ISSN={0073-8301},
     journal={Publ. Math. Inst. Hautes \'{E}tudes Sci.},
      number={107},
       pages={109\ndash 168},
         url={https://doi.org/10.1007/s10240-008-0011-6},
      review={\MR{2434692}},
}

\bib{LaOl08II}{article}{
      author={Laszlo, Yves},
      author={Olsson, Martin},
       title={The six operations for sheaves on {A}rtin stacks. {II}. {A}dic coefficients},
        date={2008},
        ISSN={0073-8301},
     journal={Publ. Math. Inst. Hautes \'{E}tudes Sci.},
      number={107},
       pages={169\ndash 210},
         url={https://doi.org/10.1007/s10240-008-0012-5},
      review={\MR{2434693}},
}

\bib{LS87}{article}{
      author={Langlands, R.~P.},
      author={Shelstad, D.},
       title={On the definition of transfer factors},
        date={1987},
        ISSN={0025-5831},
     journal={Math. Ann.},
      volume={278},
      number={1-4},
       pages={219\ndash 271},
         url={https://doi.org/10.1007/BF01458070},
      review={\MR{909227}},
}

\bib{LS90}{incollection}{
      author={Langlands, R.},
      author={Shelstad, D.},
       title={Descent for transfer factors},
        date={1990},
   booktitle={The {G}rothendieck {F}estschrift, {V}ol. {II}},
      series={Progr. Math.},
      volume={87},
   publisher={Birkh\"{a}user Boston, Boston, MA},
       pages={485\ndash 563},
      review={\MR{1106907}},
}

\bib{Lu90a}{article}{
      author={Lusztig, G.},
       title={Canonical bases arising from quantized enveloping algebras},
        date={1990},
        ISSN={0894-0347,1088-6834},
     journal={J. Amer. Math. Soc.},
      volume={3},
      number={2},
       pages={447\ndash 498},
         url={https://doi.org/10.2307/1990961},
      review={\MR{1035415}},
}

\bib{Mu95}{incollection}{
      author={Murre, J.~P.},
       title={Representation of unramified functors. {A}pplications (according to unpublished results of {A}. {G}rothendieck)},
        date={1995},
   booktitle={S\'{e}minaire {B}ourbaki, {V}ol. 9},
   publisher={Soc. Math. France, Paris},
       pages={Exp. No. 294, 243\ndash 261},
      review={\MR{1608802}},
}

\bib{MiVi07}{article}{
      author={Mirkovi\'{c}, I.},
      author={Vilonen, K.},
       title={Geometric {L}anglands duality and representations of algebraic groups over commutative rings},
        date={2007},
        ISSN={0003-486X},
     journal={Ann. of Math. (2)},
      volume={166},
      number={1},
       pages={95\ndash 143},
         url={https://doi.org/10.4007/annals.2007.166.95},
      review={\MR{2342692}},
}

\bib{Ng06}{article}{
      author={Ng\^{o}, Bao~Ch\^{a}u},
       title={Fibration de {H}itchin et endoscopie},
        date={2006},
        ISSN={0020-9910},
     journal={Invent. Math.},
      volume={164},
      number={2},
       pages={399\ndash 453},
         url={https://doi.org/10.1007/s00222-005-0483-7},
      review={\MR{2218781}},
}

\bib{Ng10}{article}{
      author={Ng\^{o}, Bao~Ch\^{a}u},
       title={Le lemme fondamental pour les alg\`ebres de {L}ie},
        date={2010},
        ISSN={0073-8301},
     journal={Publ. Math. Inst. Hautes \'{E}tudes Sci.},
      number={111},
       pages={1\ndash 169},
         url={https://doi.org/10.1007/s10240-010-0026-7},
      review={\MR{2653248}},
}

\bib{NP01}{article}{
      author={Ng\^{o}, B.~C.},
      author={Polo, P.},
       title={R\'{e}solutions de {D}emazure affines et formule de {C}asselman-{S}halika g\'{e}om\'{e}trique},
        date={2001},
        ISSN={1056-3911},
     journal={J. Algebraic Geom.},
      volume={10},
      number={3},
       pages={515\ndash 547},
      review={\MR{1832331}},
}

\bib{On63}{article}{
      author={Ono, Takashi},
       title={On the {T}amagawa number of algebraic tori},
        date={1963},
        ISSN={0003-486X},
     journal={Ann. of Math. (2)},
      volume={78},
       pages={47\ndash 73},
         url={https://doi.org/10.2307/1970502},
      review={\MR{156851}},
}

\bib{Re05}{book}{
      author={Renner, Lex~E.},
       title={Linear algebraic monoids},
      series={Encyclopaedia of Mathematical Sciences},
   publisher={Springer-Verlag, Berlin},
        date={2005},
      volume={134},
        ISBN={3-540-24241-4},
        note={Invariant Theory and Algebraic Transformation Groups, V},
      review={\MR{2134980}},
}

\bib{Re88}{article}{
      author={Renner, Lex~E.},
       title={Conjugacy classes of semisimple elements, and irreducible representations of algebraic monoids},
        date={1988},
        ISSN={0092-7872},
     journal={Comm. Algebra},
      volume={16},
      number={9},
       pages={1933\ndash 1943},
         url={https://doi.org/10.1080/00927878808823667},
      review={\MR{952217}},
}

\bib{Ric77}{article}{
      author={Richardson, R.~W.},
       title={Affine coset spaces of reductive algebraic groups},
        date={1977},
        ISSN={0024-6093,1469-2120},
     journal={Bull. London Math. Soc.},
      volume={9},
      number={1},
       pages={38\ndash 41},
         url={https://doi.org/10.1112/blms/9.1.38},
      review={\MR{437549}},
}

\bib{Ri01}{article}{
      author={Rittatore, Alvaro},
       title={Very flat reductive monoids},
        date={2001},
        ISSN={0797-1443},
     journal={Publ. Mat. Urug.},
      volume={9},
       pages={93\ndash 121 (2002)},
      review={\MR{1935829}},
}

\bib{She79}{article}{
      author={Shelstad, Diana},
       title={Orbital integrals and a family of groups attached to a real reductive group},
        date={1979},
        ISSN={0012-9593},
     journal={Ann. Sci. \'{E}cole Norm. Sup. (4)},
      volume={12},
      number={1},
       pages={1\ndash 31},
         url={http://www.numdam.org/item?id=ASENS_1979_4_12_1_1_0},
      review={\MR{532374}},
}

\bib{Sl80}{book}{
      author={Slodowy, Peter},
       title={Simple singularities and simple algebraic groups},
      series={Lecture Notes in Mathematics},
   publisher={Springer, Berlin},
        date={1980},
      volume={815},
        ISBN={3-540-10026-1},
      review={\MR{584445}},
}

\bib{Sp79}{incollection}{
      author={Springer, T.~A.},
       title={Reductive groups},
        date={1979},
   booktitle={Automorphic forms, representations and {$L$}-functions ({P}roc. {S}ympos. {P}ure {M}ath., {O}regon {S}tate {U}niv., {C}orvallis, {O}re., 1977), {P}art 1},
      series={Proc. Sympos. Pure Math.},
      volume={XXXIII},
   publisher={Amer. Math. Soc., Providence, RI},
       pages={3\ndash 27},
      review={\MR{546587}},
}

\bib{StacksP}{misc}{
      author={{Stacks project authors}, The},
       title={The stacks project},
        date={2022},
         url={\url{https://stacks.math.columbia.edu}},
}

\bib{St65}{article}{
      author={Steinberg, Robert},
       title={Regular elements of semisimple algebraic groups},
        date={1965},
        ISSN={0073-8301},
     journal={Inst. Hautes \'{E}tudes Sci. Publ. Math.},
      number={25},
       pages={49\ndash 80},
         url={http://www.numdam.org/item?id=PMIHES_1965__25__49_0},
      review={\MR{180554}},
}

\bib{St74}{book}{
      author={Steinberg, Robert},
       title={Conjugacy classes in algebraic groups},
      series={Lecture Notes in Mathematics, Vol. 366},
   publisher={Springer-Verlag, Berlin-New York},
        date={1974},
        note={Notes by Vinay V. Deodhar},
      review={\MR{0352279}},
}

\bib{Sun12b}{article}{
      author={Sun, Shenghao},
       title={Decomposition theorem for perverse sheaves on {A}rtin stacks over finite fields},
        date={2012},
        ISSN={0012-7094,1547-7398},
     journal={Duke Math. J.},
      volume={161},
      number={12},
       pages={2297\ndash 2310},
         url={https://doi.org/10.1215/00127094-1723657},
      review={\MR{2972459}},
}

\bib{Sun12}{article}{
      author={Sun, Shenghao},
       title={{$L$}-series of {A}rtin stacks over finite fields},
        date={2012},
        ISSN={1937-0652},
     journal={Algebra Number Theory},
      volume={6},
      number={1},
       pages={47\ndash 122},
         url={https://doi.org/10.2140/ant.2012.6.47},
      review={\MR{2950161}},
}

\bib{SV17}{article}{
      author={Sakellaridis, Yiannis},
      author={Venkatesh, Akshay},
       title={Periods and harmonic analysis on spherical varieties},
        date={2017},
        ISSN={0303-1179},
     journal={Ast\'{e}risque},
      number={396},
       pages={viii+360},
      review={\MR{3764130}},
}

\bib{Ti11}{book}{
      author={Timashev, Dmitry~A.},
       title={Homogeneous spaces and equivariant embeddings},
      series={Encyclopaedia of Mathematical Sciences},
   publisher={Springer, Heidelberg},
        date={2011},
      volume={138},
        ISBN={978-3-642-18398-0},
         url={https://doi.org/10.1007/978-3-642-18399-7},
        note={Invariant Theory and Algebraic Transformation Groups, 8},
      review={\MR{2797018}},
}

\bib{Ti66}{article}{
      author={Tits, J.},
       title={Normalisateurs de tores. {I}. {G}roupes de {C}oxeter \'{e}tendus},
        date={1966},
        ISSN={0021-8693},
     journal={J. Algebra},
      volume={4},
       pages={96\ndash 116},
         url={https://doi.org/10.1016/0021-8693(66)90053-6},
      review={\MR{206117}},
}

\bib{Ya04}{article}{
      author={Varshavsky, Yakov},
       title={Moduli spaces of principal {$F$}-bundles},
        date={2004},
        ISSN={1022-1824},
     journal={Selecta Math. (N.S.)},
      volume={10},
      number={1},
       pages={131\ndash 166},
         url={https://doi.org/10.1007/s00029-004-0343-0},
      review={\MR{2061225}},
}

\bib{Vi95}{incollection}{
      author={Vinberg, E.~B.},
       title={On reductive algebraic semigroups},
        date={1995},
   booktitle={Lie groups and {L}ie algebras: {E}. {B}. {D}ynkin's {S}eminar},
      series={Amer. Math. Soc. Transl. Ser. 2},
      volume={169},
   publisher={Amer. Math. Soc., Providence, RI},
       pages={145\ndash 182},
         url={https://doi.org/10.1090/trans2/169/10},
      review={\MR{1364458}},
}

\bib{Wal06}{article}{
      author={Waldspurger, J.-L.},
       title={Endoscopie et changement de caract\'eristique},
        date={2006},
        ISSN={1474-7480,1475-3030},
     journal={J. Inst. Math. Jussieu},
      volume={5},
      number={3},
       pages={423\ndash 525},
         url={https://doi.org/10.1017/S1474748006000041},
      review={\MR{2241929}},
}

\bib{Wal97}{article}{
      author={Waldspurger, J.-L.},
       title={Le lemme fondamental implique le transfert},
        date={1997},
        ISSN={0010-437X},
     journal={Compositio Math.},
      volume={105},
      number={2},
       pages={153\ndash 236},
         url={https://doi.org/10.1023/A:1000103112268},
      review={\MR{1440722}},
}

\bib{Zh17}{incollection}{
      author={Zhu, Xinwen},
       title={An introduction to affine {G}rassmannians and the geometric {S}atake equivalence},
        date={2017},
   booktitle={Geometry of moduli spaces and representation theory},
      series={IAS/Park City Math. Ser.},
      volume={24},
   publisher={Amer. Math. Soc., Providence, RI},
       pages={59\ndash 154},
      review={\MR{3752460}},
}

\end{biblist}
\end{bibdiv}


\cleardoublepage
\phantomsection
\addcontentsline{toc}{chapter}{List of Notations}
\printnomenclature


\cleardoublepage
\phantomsection
\addcontentsline{toc}{chapter}{Index}
\printindex

\end{document}